**W. B. VASANTHA KANDASAMY**

**FLORENTIN SMARANDACHE**

# ANALYSIS OF SOCIAL ASPECTS OF MIGRANT LABOURERS LIVING WITH HIV/AIDS USING FUZZY THEORY AND NEUTROSOPHIC COGNITIVE MAPS

## With Specific Reference To Rural Tamilnadu In India

# ANALYSIS OF SOCIAL ASPECTS OF MIGRANT LABOURERS LIVING WITH HIV/AIDS USING FUZZY THEORY AND NEUTROSOPHIC COGNITIVE MAPS

## With Specific Reference To Rural Tamilnadu In India

**W. B. Vasantha Kandasamy**

Department of Mathematics
Indian Institute of Technology, Madras
Chennai – 600036, India
e-mail: **vasantha@iitm.ac.in**
web: **http://mat.iitm.ac.in/~wbv**

**Florentin Smarandache**
Department of Mathematics
University of Gallup, NM
New Mexico, USA
e-mail: **smarand@gallup.unm.edu**

Translation of the Tamil interviews by
**Meena Kandasamy**

*2004*

# Contents

# Preface

Neutrosophic logic grew as an alternative to the existing topics and it represents a mathematical model of uncertainty, vagueness, ambiguity, imprecision, undefined-ness, unknown, incompleteness, inconsistency, redundancy and contradiction. Despite various attempts to reorient logic, there has remained an essential need for an alternative system that could infuse into itself a representation of the real world. Out of this need arose the system of neutrosophy and its connected logic, neutrosophic logic. This new logic, which allows also the concept of indeterminacy to play a role in any real-world problem, was introduced first by one of the authors Florentin Smarandache.

In this book for the first time we have ventured into the total analysis of migrant labourers in rural Tamil Nadu who are victims of HIV/AIDS using FCM, BAM and Neutrosophic Cognitive Maps. As in our study and analysis we felt several of the factors related with the psycho, socio, economic problems of these HIV/AIDS patients from rural Tamil Nadu (a southernmost state in India) remain indeterminate apart from the data being an unsupervised one. At the outset, we first emphasize that the study and analysis (and there by the conclusions and suggestions) pertain only to migrant labourers from rural Tamil Nadu who are poor and uneducated and who are HIV/AIDS infected. So this study cannot be extended to urban people or rich/ sophisticated/ educated populations of Tamil Nadu.

We have purposefully chosen this study for we saw majority of the people who come to take treatment for HIV/AIDS as inpatients of the Government Hospital for Thoracic Medicine at Tambaram (called the "Tambaram Sanatorium") hospital are from rural areas with no education, they are poor or become poor due to this disease and have caught this disease due to migration. Further these people are not empowered with trade unions or welfare associations, for these labourers are used till they serve the purpose and once they are ill there is no one to take care of them. At least 75% of them, when they came to know that they had HIV/AIDS, their first reaction was to think of committing suicide. These men at large are shaken to the core in their hearts with fear of stigma and poverty and they suffer from deep depression when they are just in the prime of life. Thus our sample under analysis is these migrant poor uneducated labourers from rural Tamil Nadu. This book is an extended mathematical analysis of the book *Love.Life.Lust.Loss: 101 Real Life Stories of Migration and AIDS—A Fuzzy Analysis*.

This book has seven chapters. The first chapter is introductory in nature and it speaks about the migrant labourers. In chapter two we use Fuzzy Cognitive Maps to analyze the socio-economic problems of HIV/AIDS infected migrant labourers in rural areas of Tamil Nadu. In chapter three we analyze the role played by the government helping these migrant labourers with HIV/AIDS and factors of migration and their vulnerability in catching HIV/AIDS. For the first time Neutrosophic Cognitive Maps are used in the study of migrant labourers who have become HIV/AIDS victims. This study is done in Chapter IV. In chapter V we use Neutrosophic Relational Maps and we define some new neutrosophic tools like Combined Disjoint Block FRM, Combined Overlap NRM and linked NRM. We adopt these new techniques in the study and analysis of this problem. Chapter VI gives a very brief sketch of the life history of these 60 HIV/AIDS infected migrant labourers so that people from different social and cultural backgrounds follow our analysis. The last chapter gives suggestions and conclusions based on our study. This book has 8 appendices. The Appendix 1 is the questionnaire used in interviewing these migrant labourers. The questionnaire was formatted with the help of Dr. K. Kandasamy. Appendix 2 gives the table of these 60 migrant labourers. Appendix 3 to 8 gives the C programs used in finding the resultant vector for any given state vector using the dynamical system. Appendix 9 deals with Neutrosophic logic. We have included a list of further reading.

We deeply acknowledge the TNSACS for its financial support and encouragement. We personally acknowledge the project director Mr. Deenabandu, I.A.S of TNSACS for his constant encouragement. We are indebted to Ms.Kama Kandasamy who had taken all pains to format, verify the equations and design the books. We deeply acknowledge the services of Dr.K.Kandasamy who had done a lot of background work and for formatting such a scientific linguistic questionnaire without which we would not have been in a position to get this data. We thank our research scholars Mr.Ganesan, Mr.Kanagamuthu and Mr.Narayanamurthy for taking the interviews. Our deep thanks are due to Meena Kandasamy who translated the Tamil interviews of these migrant labourers.

**We dedicate this book to the peoples and nations that are working to rehabilitate and support the poor men and women living with HIV/AIDS in the Third-World countries.**

W.B. Vasantha Kandasamy
Florentin Smarandache

Chapter One

# INTRODUCTION

Migration acquires great significance in the study of peoples and populations, for it not only involves the merely mathematical spatial redistribution of people, but also because it has enormous impact on livelihood, life-styles, employment, socioeconomic and political stability; or in other words, it influences the entire society. When used in the geographical context, the term 'migration' refers to the 'permanent or semipermanent change in the residence of an individual person or group of people.' Traditional studies classify migrants under two broad heads: voluntary migrants and involuntary migrants.

In this Indian scenario, voluntary migration (where the migrants move of their own choice) often takes place with the view to secure a livelihood. On the other hand, involuntary migration might take place due to natural disasters or war and the consequent persecution or in some of the cases because of environmental degradation. In between these two diverse classifications, there exist a great mass of people, whose migration is a direct result of their poverty, and whose migration is often, as much voluntary as involuntary.

In the same manner, researchers also divide the cause for migration into push factors (factors that push people to search for new areas to live: poverty, unemployment, agricultural failure) and pull factors (factors that pull people by exerting an attraction for a new location), taking into their foremost view the economic factors. Most migrations are a combination of both push and pull factors; the best illustration of this is the rural-urban migration. In this book, we skip doing an exhaustive review of the literature on the said topic. We like to point out that most of the mathematical research into migration has dealt with it on purely statistical lines and studied it on economic terms, while ignoring that migration is the result of an inter-relationship of an aggregate of several factors.

As a part of our research study of people with HIV/AIDS we have selected to study those affected as a result of migration: migrant labourers, long distance truckers, military servicemen, commercial sex workers. The term migrant labourers in this book means not only labourers who migrate from one place to another place and settle temporarily or permanently for livelihood but also those workers who travel from one place to another for completion of work (technically the bridge population: linking high-risk and low-risk communities and geographic areas) i.e. labourers like truck/lorry/cab drivers, also long distance government or private bus drivers and those accompanying them like the cleaners etc. This category of migrant labourers also includes construction labourers who stay for a month or so in a different town and carry out their job. Above all, this includes the adults or children who leave home due to frustration or employment i.e. to be precise—runaways from home. The HIV transmission typically occurs as a result of unprotected sex in destination areas and between destination and home areas. Thus, while we use the term migrant workers/labourers, we aim to include the mobile workers just as much.

Migrant workers are more vulnerable to HIV/AIDS than the local population because of their poverty, lack of power, lack of health awareness and unstable life-style. Sometimes, wrong and mythical notions about sex might lead to their high-risk behaviour. The main reason for taking up this study is that in the state of Tamil Nadu, migrant labourers (including truck/lorry drivers) constitute the major portion of those affected with HIV/AIDS. This is especially true in districts like Namakkal. Thus we are justified in studying the psychological, social and economic status of these HIV/AIDS infected migrant labourers and from this analysis we can draw conclusions on preventive measures to be taken to stop the spread of HIV as well as take steps to counsel and rehabilitate the affected. We observe from our study that the large-scale patronization of Commercial Sex Workers (CSWs) by heterosexual (and often married) men is to a great extent responsible for the spread of the disease. Further it is to be noted that most of the infected are illiterate or with very low educational qualifications (secondary/high school dropouts), low levels of awareness about HIV/AIDS and are addicted to habits like smoking, drinking and visiting CSWs that too mainly from rural areas. Further most of them were persons with low sense of personal accomplishment and they took pride in their flimsy life-styles.

One of our most important and vital understanding of AIDS came through when we worked towards a feminization of understanding and when we employed fuzzy methods of psychosocial analysis. We could comprehend the risk-behaviours better and correlate it with the raw data/information/feelings that we derived from the interviews when we looked into gender-stereotypes and socialization patterns.

Such conventional beliefs about male (hetero)sexuality have to be shed and impetus has to be given to shattering stereotyped patterns of behaviour. For, this embedded concept of 'masculinity' often defeats the very premise of the prevention programs: for, masculinity is predominantly (mis)associated with having multiple sex partners, and a complete void between emotional/passionate involvement and sexual activity.

In fact campaigns and earlier researches have also targeted the promiscuous heterosexual men and has sought to downplay and distinguish/segregate the concepts of masculinity and risk-taking (that includes multiple sex partners, unprotected sex with CSWs). This perceived and heightened sense of masculinity is garnered from the words the AIDS affected migrant workers use to describe the reasons for visiting commercial sex workers: to enjoy, to be jolly, etc. It also is an underlying reason as to why they pick up habits like alcoholism, addiction to tobacco and drugs. Doctors associated with our study also linked this prevalence of high-risk behaviour to popular myths present among sections of migrants.

For instance, many truck drivers had said to doctors that they visited CSWs because during the many days of travelling their body's heat became greater and because they believed sexual intercourse would lower the heat, they visited CSWs. They said they wanted this kind of 'release.' In our study we observed that majority of the HIV infections occurred only in Mumbai, followed by Andhra Pradesh. This is because of the situation of both Kamatipura, the largest red-light area in the world, and Cotton Green, the largest truck terminus in Asia both at Mumbai. While Kamatipura is notoriously famous for its brothels, Cotton Green deserves some mention: everyday more than 5,000 trucks come from different states of India. The inflow and outflow of the truck traffic is intense, majority of the trucks visit the terminus once in four to eight days and this perhaps sheds some more light on the risk-behaviour and increased chances of casual sex by migrant workers (in this case the truck drivers) in Mumbai.

Another major fact that we could unearth during this research on migrant workers with HIV/AIDS was: several of them, (56 to be precise) were unaware of how HIV/AIDS spreads. This was very shocking for it shows to us that the ad and awareness campaigns in the media have not had their sufficient impact and secondly in the case of women whose glaring ignorance about this disease has rendered them helpless in averting it mainly in rural areas. While knowledge of AIDS has been widespread, the knowledge of its method of spread is not so well known, sadly, even among the high-risk groups. This is not to say that the campaigns have failed in their objectives, this is to merely suggest that methods have to be developed to have greater outreach and awareness among the vulnerable categories of the people.

This is especially important in the cases of truck-drivers and lorry-drivers, because this mobile population spreads the disease at alarming levels. We selected this sample from around 200 people whom we had interviewed as a part of our research project. Other categories like NGOs, people with AIDS, and members of the public had been interviewed. From that we have 60 men living with HIV/AIDS (who had acquired the disease directly or indirectly due to migration) for our sample study on migration workers and their families. The interviews were taken from the people with HIV/AIDS taking treatment from the Government Hospital for Thoracic Medicine at Tambaram (we refer to this as the Tambaram Sanatorium). Some of the interviews are exhaustive, and give us great insight into the minds of people with AIDS and at the same time some of the patients exhibit a reluctance to talk about the disease, were we have been left with the thought of indeterminancy.

Out of 60 HIV/AIDS patients we have taken for this analysis, they were, at some point of their lives, migrant workers Further, we could note that in case of workers of the unorganized sector, medical testing for AIDS remains remote. Moreover, in intervention programs among migrant labourers, NGOs have faced several difficulties, one of the outstanding reasons being that the migrant workers are unorganized labourers and do not have unions or associations through which successful awareness programmes can be carried. As a result, they are an exploited working class that does not enjoy legal safeguards or claim workers' rights. Further the peer-group influence among such migrants causes further havoc to the fabric of society. Our data indicates that 85% of the migrant workers (51 of 60) were led into visiting CSWs because of their friends.

In this study we do not analyze the gender difference in treatment of HIV/AIDS affected though one is well aware of the established statistic in India that heterosexual intercourse is in more than 85% of the cases the cause of HIV transmission. Difference in status between men and women, supported by social and cultural systems make it difficult for women to take preventive and safe-sex measures. Further, the patriarchal Indian family set-up ensures tragically that men do not divulge their HIV seropositivity to their wives until at the very last stages. As a result women come for treatment with full-blown AIDS, suffering from all types of opportunistic infections and with nobody to look after them.

Thus women who contract AIDS (from their migrant husbands) become innocent and hapless victims, and, our research confirmed the oft-reported fact that women in the dying stages of AIDS are pathetic as the society shuns them and there is no one to look after them, whereas women always stay with their husbands and take care of them when the husbands are affected with AIDS. Another psychological aspect, which we could notice, was that couples infected with HIV/AIDS, where both the partners were alive appeared less depressed and less perplexed than individuals/surviving members with AIDS. We were also taken aback by the number of cases where the migrant workers (who had visited CSWs) did not adopt safe-sex practices even though the CSWs had themselves requested them to use safe-sex methods and had offered condoms. This proves that though NGO intervention has given CSWs responsibility and alacrity in their functioning, the men from rural areas visiting them do not adopt safe-sex practices.

Further, by forcing even the CSWs to have unprotected sex, they expose them to greater risk. Some of the reasons, which were attributed to not using condoms, were: intoxication at time of intercourse, non-availability of condoms, no satisfaction, ignorance/indifference, it is against raw masculinity and inhibition to buy condoms. We declare that throughout this book, by CSWs we mean Female Sex Workers (FSWs). Further, we also like to put forth the economic data about the sample we have chosen for this study. The 60 men in our selection were not highly affluent or surrounded by the plushness of rich-urban life. They did not have educational qualifications or family backgrounds to boast of, and most of them were from the poor, semi-poor, and lower middle class. As a result of our sample remaining in this manner, we cannot predict or give information about the distribution/prevalence of AIDS among varied economic strata. This study is also only restricted to the remote villagers with illiteracy prevailing in them.

Further we felt it essential to study only this sample first; for the rich and middle class are educated, have awareness and exposure as well as access to treatment. An important aspect of our study probed into the hitherto unresearched areas of religion and its role in the AIDS epidemic. We observed in our study the attitudes of the people with HIV/AIDS towards religion. We found out that a substantial number, of men had converted from their native religion to Christianity. This is an unprobed area and there is a complete absence of information with regard to the vital role of religious faith.

Finally as several interlinking of the relations and its effect on the nodes in some cases may remain to be an indeterminate and to the best of our knowledge that study of indeterminate relations is carried out mainly and solely by the neutrosophic models. We in this book have for the first time adopted the neutrosophic models in the analysis of the socio-economic problems of the HIV/AIDS affected patients who are migrant labourers. In our opinion neutrosophic models in the study of this problems would be better than the fuzzy models for we can also give the answer as indeterminate for some relations, where as fuzzy models can say the existence or non existence of any relation only. Thus we conclude in the study of the problem of HIV/AIDS in the migrant labourers the implementation of neutrosophic models is much more powerful than the fuzzy model. This book uses the notion of neutrosophic cognitive maps, combined neutrosophic cognitive maps, neutrosophic relational maps and further define some new techniques for first time called combined disjoint block FCM, combined overlapping block FCM. The corresponding notions in the case FRMs have been defined for the first time and the model is implemented in the study of the problem. Also the analogous neutrosophic models, is introduced for the first time in this book and the analysis of the study HIV/AIDS migrant labourers. To adopt the technique of neutrosophy, a neutrosophic questionnaire was used to interview the patients.

The book is organized into 7 chapters. Chapter I is introductory in nature describing the migrancy problem in India especially in Tamil Nadu. In chapter II fuzzy matrices are used for the first time in the analysis of age-group and its related diseases among migrant HIV/AIDS affected labourers. In this chapter we just recall the notion of fuzzy cognitive maps and study the model using FCMs. Also the new notion of combined disjoint block FCM and combined overlap block FCM are introduced in this chapter and these newly defined models are used in the study of the problem.

In chapter III yet another fuzzy theory technique namely the bi-directional associative memories is used in the analysis of the cause for vulnerability to HIV/AIDS and factors for migration. A C-program is given in Appendix 6, which is used to simplify working of the problem. In Chapter IV we first just recall some basic notions of neutrosophic theory, neutrosophic cognitive maps that are utilized in the analysis of the problem. We use the combined neutrosophic cognitive maps for the study. In this chapter we have introduced for the first time the notion of combined disjoint block NCM and combined overlap block NCM. They are adopted for this problem and conclusions are based on the analysis.

In chapter V we analyse the problem of socio economic conditions of migrant labourers relative to the vulnerability of being affected by HIV/AIDS, using both fuzzy theory and neutrosophy theory. We use the fuzzy models like relational maps, combined fuzzy relational maps and also use neutrosophic relational maps and combined neutrosophic relational maps. We in this chapter have introduced the new notions of combined disjoint block FRM, combined overlap block FRM, combined disjoint block NRM, combined overlap block NRM. These models are adopted in the analysis of the problem.

In chapter VI we briefly recall the interviews of the 60 migrants HIV/AIDS affected patients. This is mainly done for any reader to have a view of the social, economic and the psychological conditions of the rural migrant men who are mostly from the economically poor strata from Tamil Nadu.

For unless one knows more about this it will not be easy for a reader to understand the situation other than the native of Tamil Nadu; so we have recalled these interviews to make this book more contextual. This final chapter gives suggestion and conclusions based on our study.

This book has 8 appendixes; the first one is the neutrosophic questionnaire. Second the table of the 60 migrant labourers. The appendixes 3-8 are the C-programs, which are utilized to make the calculations simple and easy.

It is important to mention that the relevant mathematical notions used in each of these chapters are described in the respective chapters there by making each chapter a self-contained one. The conclusions are based on several of the state vectors and many experts opinions are utilized, to make the book a comprehensive one we have given one or two experts for illustrating and describing the model. The study is completely based on the interviews with the patients and the conclusions are valid mainly for the rural uneducated poor migrants of rural Tamil Nadu. By no means are these conclusions true for the urban rich people of Tamil Nadu, for even a proper statistics of how many urban people suffer with HIV/AIDS cannot be got by any one as they keep it as a secret due to social stigma.

Chapter Two

# ANALYSIS OF THE FEELINGS OF HIV/AIDS AFFECTED MIGRANT LABOURERS USING FUZZY MARICES AND FCM

This chapter has five sections. In the first section we introduce the fuzzy matrix model and employ it on the sample of 60 HIV/AIDS patients who are migrant labourers. Conclusions are given based on the model analysis.

In section two we analyse the socio economic problem of the HIV/AIDS affected migrant labourers using FCM. By studying the hidden pattern one can easily derive the results. In section three we study this model using combined FCMs simplified the mathematical part of it, so that a social scientist can also follow the book.

In section four we introduce the notion of combined disjoint block  FCM and use it to analyse this problem. In section five the same problem is analysed using the new model viz. combined overlap block FCM model.

So, the demand made on the mathematical knowledge of the reader is very meager. Also graphs are given to make the understanding simple and easy. The tools used in section 1 and section 2 are made very simple for any common reader to follow.

By this we have to a large extent suppressed some of the laborious technical working but give the results in an explicit way. The third and the fourth section however happens to be little difficult but our simple descriptions from the basic level will certainly make every one understand and follow it, if not at the first reading, but in the subsequent readings.

The data that we deal with, happens to be an unsupervised one, for we do not and cannot categorically put a number for a person visiting the CSWs, for it might be due to uncontrollable feeling, no one to monitor them in a different place, because of being away from the family, it might be easy money, or leisure, or bad company, or bad habits like alcohol or drugs, family dispute, luring by the CSWs when the migrants are in a drunken state and so on. So we cannot give it any statistical form/data representation for the psychological study for the cause may remain an indeterminate for which this chapter cannot give any analysis. Thus, for this study it has become inevitable for us to use fuzzy theory, as these models function mainly based on the working of the neurons of the brain.

The dynamical system approximately works like a human brain, so only we feel these methods are best-suited for this study.

## 2.1 Estimation of the Overall Age Group of the HIV/AIDS Patients from migrant labourers in which they are Maximum Infected, Using Fuzzy Matrices

In this section we give an algebraic approach to the problem of studying the psychological aspects and the age group of HIV/AIDS patients using matrices. We determine the peak age group viz. the maximum number of HIV/AIDS patients with diseases like TB, STD, skin ailments, stomach problems, fever and cold and their habits like smoking, alcoholism, bad company and CSWs. Also we estimate the overall age group. The study is significant for it can be adopted to other states like Andhra Pradesh, Madhya Pradesh, Kerala, etc. with suitable modifications. Thus our study is valid not only to any other state but also to every district/village/region/town in Tamil Nadu or in any HIV/AIDS infected nation.

The analysis of our study is mainly focused on the in/out patients of the Tambaram Sanatorium who have contracted the disease by migration that is  who had migrated (temporarily or permanently). Migration does not mean a permanent change of place, even if the period of mobility is for a week, it is taken as infected or contracted the disease by migration. The raw data obtained from the 60 HIV/AIDS patients from the in/out patients of the Tambaram Sanatorium is taken as the raw data for our study.

This raw data is converted into time-dependent matrices. By time dependent matrices we mean the matrices which are dependent on the age of the HIV/AIDS patients. So when we say time we only mean the age groups. After obtaining the time dependent matrices (age-dependent matrices) using the technique of usual mean/average and standard deviation we identify both age group in which they are maximum affected by the disease like STD, TB etc. and also the habits like smoking, alcohol, bad company, CSWs etc. The identification of the age group when they are maximum affected plays a vital role in improving the awareness about HIV/AIDS, and preparing well-tailored intervention programs. Our study of the related symptoms and the relevant habits prevalent in those age groups is also effective, for the government can take proper steps and precautionary methods to have an effective outreach among that age group and spread awareness preventing those already infected from further infecting others.

The maximum affected indicators or the standards used depend on the type of data collected by us. The majority of such standards are based on their education, habits, profession, migration (staying away from family) and bad company and so on. Study of HIV/AIDS male patients, have been analyzed taking mainly into account that migration is the main cause for them to visit CSWs and thereby acquire/contract the disease. The raw data (who are HIV/AIDS patients) under investigation is classified under the six broad heads: CSWs, Other Women, Smoking, Alcohol, Negative peer-group influence (Bad Company), and Quacks. These six broad heads forms the columns of the matrices. The time periods are taken as the age groups varying from > 20 to < 46. Estimation of the age group, which is most likely to be vulnerable to these bad habits and hence face the risk of becoming HIV/AIDS patients, is carried out in five stages.

In the first stage we get a matrix representation of the raw data. Entries corresponding to the intersection of rows and columns are values corresponding to a live network. This initial M × N matrix is not uniform i.e. the number of years in each interval of time period may not be the same; for instance, one period may be 21-30 and other 31-35 so on.. So, in the second stage, in order to obtain an unbiased uniform effect on each and every data so collected, we transform this initial matrix into an Average Time-dependent Data (ATD) matrix ($a_{ij}$). To make the calculations simpler, in the third stage, we use the simple average techniques and convert the above Average Time-dependent Data (ATD) matrix with entries $e_{ij}$ where $e_{ij} \in \{-1, 0, 1\}$. We name this matrix as the Refined Time-dependent Data (RTD) matrix.

$e_{ij}$ are calculated using the formula

$$\text{if } a_{ij} \leq (\mu_j - \alpha * \sigma_j) \text{ then } e_{ij} = -1$$
$$\text{else if } a_{ij} \in (\mu_j - \alpha * \sigma_j, \mu_j + \alpha * \sigma_j) \text{ then } e_{ij} = 0$$
$$\text{else if } a_{ij} > (\mu_j + \alpha * \sigma_j) \text{ then } e_{ij} = 1$$

where $a_{ij}$ are entries of the ATD matrix and $\mu_j$ and $\sigma_j$ are the mean and the standard deviation of the $j^{th}$ column of the ATD matrix respectively. $\alpha$ is a parameter taken from the interval [0, 1].

The value of the $e_{ij}$ corresponding to each entry is determined in a special way. At the next stage using the RTD matrix we get the Combined Effect Time-dependent Data matrix (CETD matrix), which gives the cumulative effect of all these entries. In the final stage we obtain the row sums of the CETD matrix. A program in C language is written which easily estimates all these five stages. The C program shall give details of the mathematical working and method of calculation. Now we give the description of the problem and the proposed solution to the problem.

PROPOSED SOLUTION TO THE PROBLEM

Using the raw data available, that is the data collected by us; we analyze it via matrices and predict mathematically the maximum age group in which migrant labourers are affected by HIV/AIDS among migrant labourers in rural areas of Tamil Nadu. Thus to be more precise our chief problem is the prediction of the age group which is most affected and most vulnerable to HIV/AIDS among migrant labourers in rural areas of Tamil Nadu. We have established that the predicted age group with maximum HIV/AIDS patients coincides with the real data. This makes the mathematical research carried out by us extremely valid. This age group only concerns the HIV/AIDS patients who have migrated and have acquired the disease due to migration: only of the men. Our analysis not only predicts the maximum patients in the age group but also suggests that for certain age group the awareness program is not so essential and instead the government can pitch upon that age group which may be given a rigorous awareness program. The advertisement and awareness program should be very suggestive to make them aware of the disease and the consequence of their high-risk behavior which will ultimately land them with AIDS.

The Initial Raw data matrix of the 55 migrant labourers affected with HIV/AIDS relating their age group and diseases are given in the following:

Initial Raw Data Matrix

| Age Group | CSWs | Other Women | Smoking | Alcohol | Bad Company | Quacks |
|---|---|---|---|---|---|---|
| 21-30 | 22 | 10 | 21 | 20 | 18 | 12 |
| 31-35 | 17 | 4 | 14 | 15 | 12 | 8 |
| 36-46 | 16 | 8 | 16 | 16 | 14 | 9 |

ATD Matrix

| Age Group | CSWs | Other Women | Smoking | Alcohol | Bad Company | Quacks |
|---|---|---|---|---|---|---|
| 21-30 | 2.2 | 1 | 2.1 | 2 | 1.8 | 1.2 |
| 31-35 | 3.4 | 0.8 | 2.8 | 3 | 2.4 | 1.6 |
| 36-46 | 1.46 | 0.73 | 1.46 | 1.46 | 1.27 | 0.82 |

The Average and the Standard Deviation

| Average | 2.35 | 0.84 | 2.12 | 2.15 | 1.82 | 1.21 |
|---|---|---|---|---|---|---|
| Standard Deviation | 0.80 | 0.11 | 0.55 | 0.64 | 0.46 | 0.32 |

The RTD matrix or fuzzy matrix for $\alpha = 0.3$

$$\begin{bmatrix} 0 & 1 & 0 & 0 & 0 & 0 \\ 1 & -1 & 1 & 1 & 1 & 1 \\ -1 & -1 & -1 & -1 & -1 & -1 \end{bmatrix}$$

The row sum matrix

$$\begin{bmatrix} 1 \\ 4 \\ -6 \end{bmatrix}$$

The graph depicting the maximum age group of the migrant labourers who are affected by HIV/AIDS for $\alpha = 0.3$

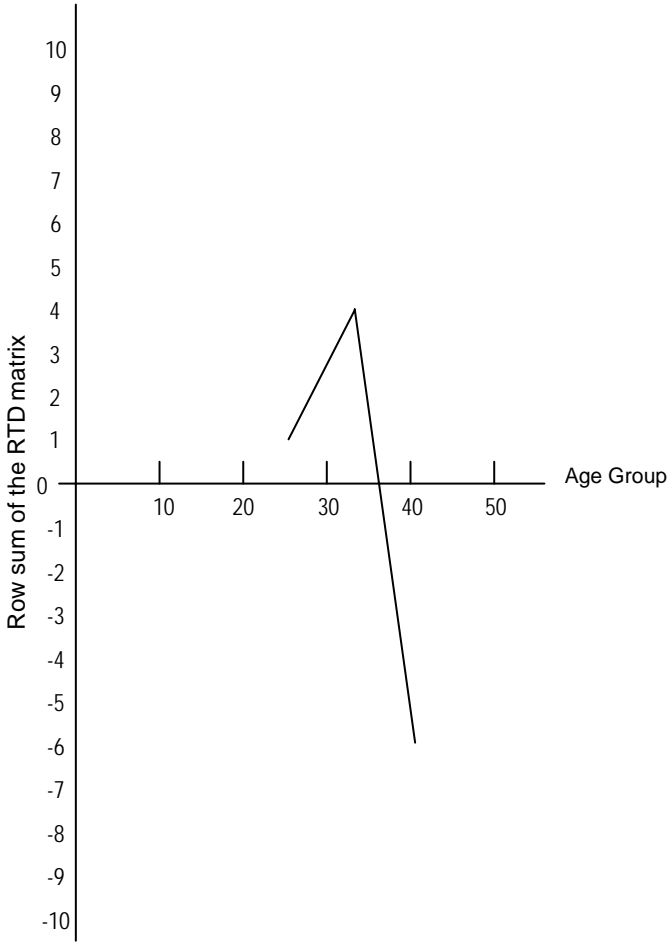

The RTD or Fuzzy Matrix for $\alpha = 0.7$

$$\begin{bmatrix} 0 & 1 & 0 & 0 & 0 & 0 \\ 1 & 0 & 1 & 1 & 1 & 1 \\ -1 & -1 & -1 & -1 & -1 & -1 \end{bmatrix}$$

The row sum matrix

$$\begin{bmatrix} 1 \\ 5 \\ -6 \end{bmatrix}$$

The graph depicting the maximum age group of the migrant labourers who are affected by HIV/AIDS for α = 0.7

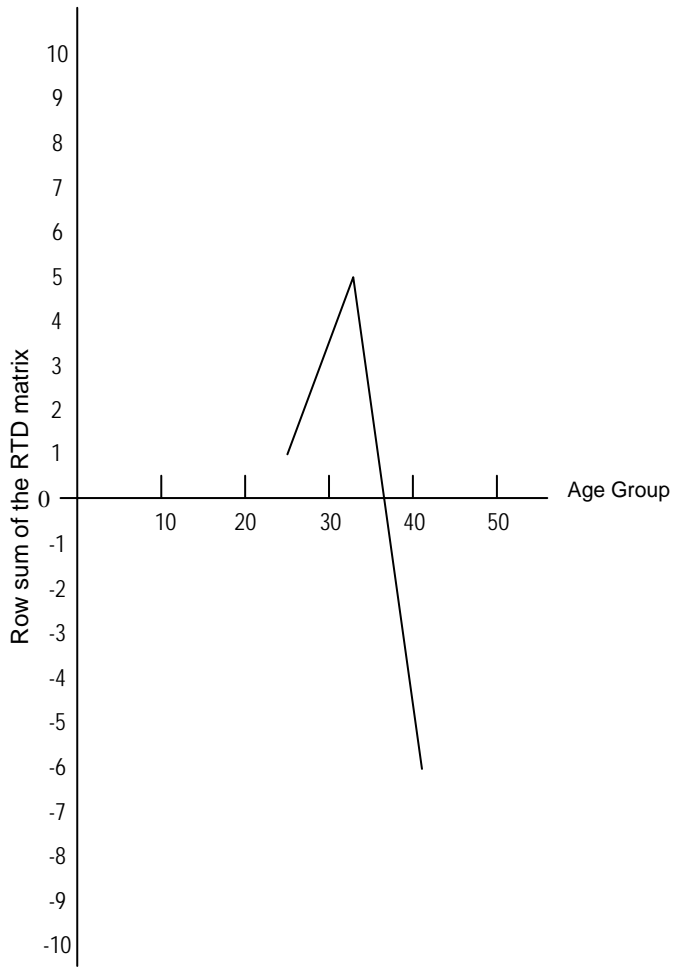

The RTD or Fuzzy Matrix for α = 1

$$\begin{bmatrix} 0 & 1 & 0 & 0 & 0 & 0 \\ 1 & 0 & 1 & 1 & 1 & 1 \\ -1 & -1 & -1 & -1 & -1 & -1 \end{bmatrix}$$

The row sum matrix

$$\begin{bmatrix} 1 \\ 5 \\ -6 \end{bmatrix}$$

The graph depicting the maximum age group of the migrant labourers who are affected by HIV/AIDS for α = 1

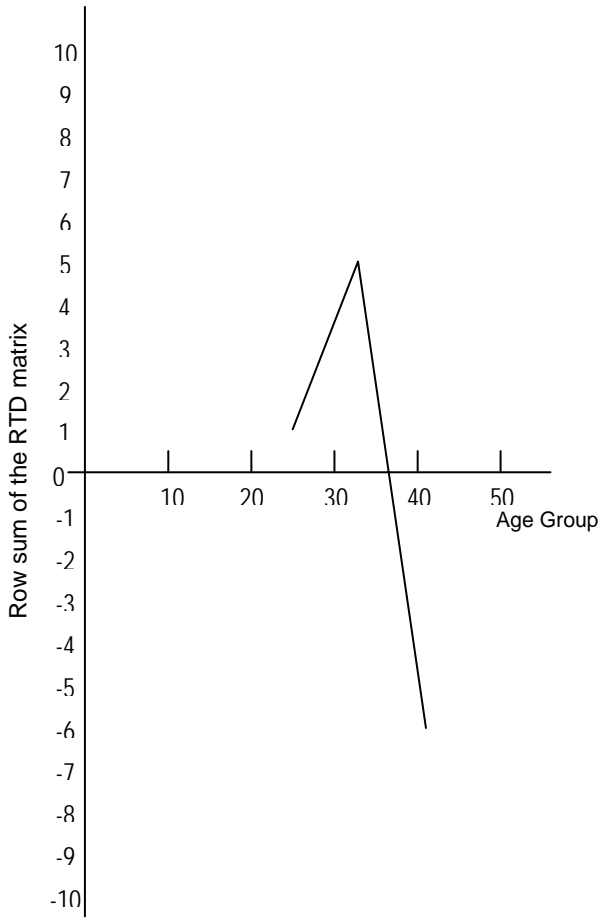

The only mathematical inference is that the most vulnerable age group to have HIV/AIDS is 31-35 among the migrant labourers of rural area in Tamil Nadu who are uneducated. To ones surprise they are regular visitors (or have casually visited) of the commercial sex workers.

Their chances of contracting AIDS from other women is zero and goes in the negative, whereas it is very important to note that the age group 21-30 visit other women also. Here by other women, we mean women who are not CSWs, but those who share a casual sex partnership with these men, which is extramarital or premarital.

So one is partly forced to think whether married men seek CSWs whereas unmarried men seek other women. This is the very important note which is evident from this model that after the age of 35 the migration is itself meagre, so the HIV/AIDS patients among migrants gives a negative value viz. -18, and we feel that this needs more analysis and attention for after an age the uneducated rural men do not migrate for livelihood but live in the same village but their children take up the role of migration.

The CETD Matrix is given in the following, which is formulated using the three RTD matrices or fuzzy matrices for the values of $\alpha$ = 0.3, $\alpha$ = 0.7 and $\alpha$ = 1.

The CETD Matrix

$$\begin{bmatrix} 0 & 3 & 0 & 0 & 0 & 0 \\ 3 & -1 & 3 & 3 & 3 & 3 \\ -3 & -3 & -3 & -3 & -3 & -3 \end{bmatrix}$$

Row sum of CETD Matrix

$$\begin{bmatrix} 3 \\ 14 \\ -18 \end{bmatrix}$$

Graph depicting the concentrated age group of HIV/AIDS patients among the migrant labourers obtained from the CETD matrix is given in the following page:

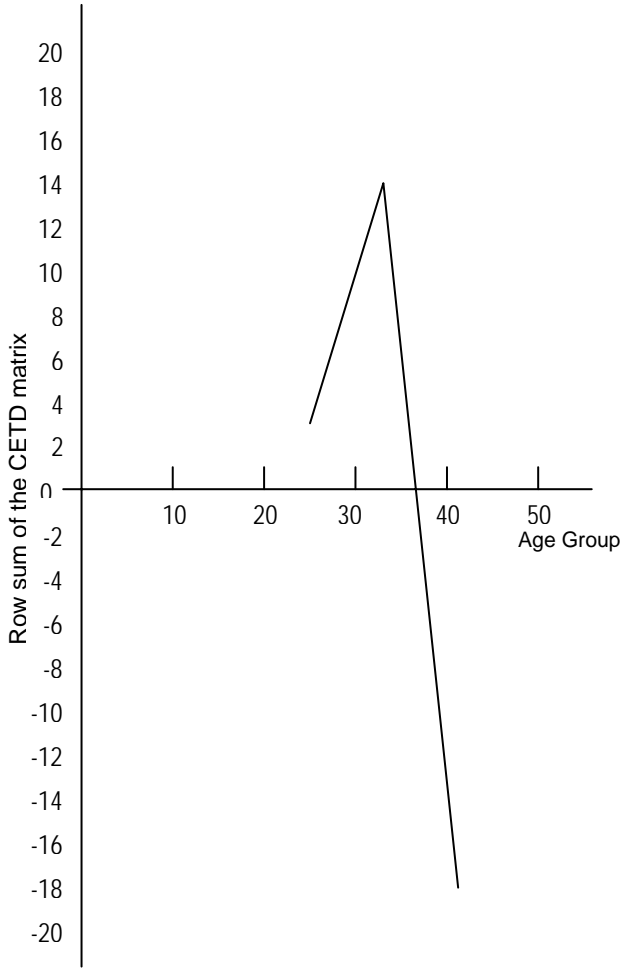

Now we increase the number of rows that is the grouping of the age (making the study more sensitive) and see whether the decision arrived is more sensitive to the earlier one we have discussed.

Thus we give the raw data of the matrix for these 58 HIV/AIDS male patients:

Initial Raw Data Matrix

| Age Group | CSWs | Other Women | Smoking | Alcohol | Bad Company | Quacks |
|---|---|---|---|---|---|---|
| 20-25 | 7 | 3 | 5 | 9 | 10 | 6 |
| 26-30 | 16 | 7 | 16 | 11 | 8 | 6 |
| 31-35 | 17 | 4 | 14 | 15 | 12 | 8 |
| 36-47 | 18 | 8 | 18 | 16 | 14 | 9 |

ATD Matrix

| Age Group | CSWs | Other Women | Smoking | Alcohol | Bad Company | Quacks |
|---|---|---|---|---|---|---|
| 20-25 | 1.16 | 0.5 | .83 | 1.5 | 1.67 | 1 |
| 26-30 | 3.20 | 1.4 | 3.2 | 2.2 | 1.6 | 1.2 |
| 31-35 | 3.40 | 0.8 | 2.8 | 3 | 2.4 | 1.6 |
| 36-47 | 1.50 | 0.67 | 1.5 | 1.3 | 1.17 | 0.75 |

Average and Standard Deviation

| Average | 2.32 | 0.84 | 2.08 | 2 | 1.71 | 1.14 |
|---|---|---|---|---|---|---|
| Standard Deviation | 0.99 | 0.34 | 0.96 | 0.67 | 0.44 | 0.31 |

The RTD Matrix or fuzzy matrix for $\alpha = 0.2$

$$\begin{bmatrix} -1 & -1 & -1 & -1 & 0 & -1 \\ 1 & 1 & 1 & 0 & -1 & 0 \\ 1 & 0 & 1 & 1 & 1 & 1 \\ -1 & -1 & -1 & -1 & -1 & -1 \end{bmatrix}$$

The row sum of the RTD matrix

$$\begin{bmatrix} -5 \\ 2 \\ 5 \\ -6 \end{bmatrix}$$

The graph for the parametric value α = 0.2 is given below:

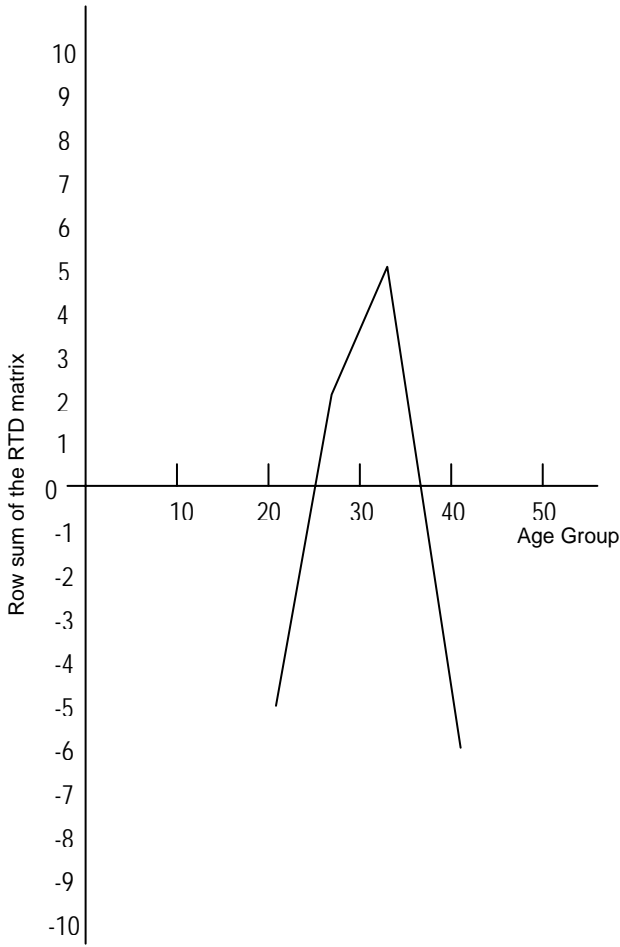

The age-group 31-35 is most vulnerable just followed by 26-30 least affected or relatively unaffected in the age group 36-47.

The RTD Matrix for α = 0.7

$$\begin{bmatrix} -1 & -1 & -1 & -1 & 0 & 0 \\ 1 & 1 & 1 & 0 & 0 & 0 \\ 1 & 0 & 1 & 1 & 1 & 1 \\ -1 & 0 & 0 & -1 & -1 & -1 \end{bmatrix}$$

Row sum of RTD matrix is given by

$$\begin{bmatrix} -4 \\ 3 \\ 5 \\ -4 \end{bmatrix}$$

The graph for the parametric value $\alpha = 0.7$ is given below:

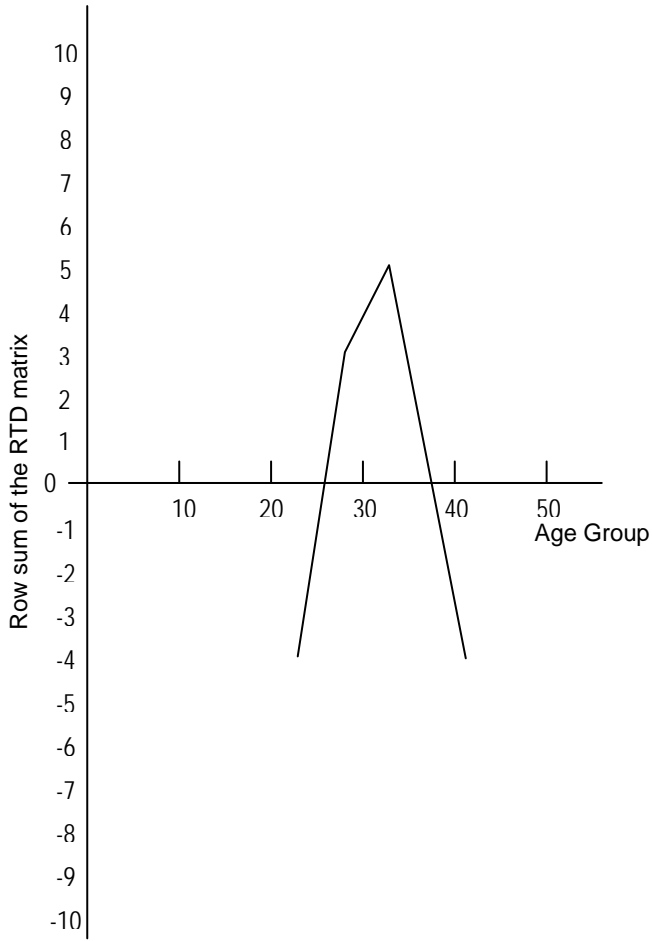

Next we proceed on to present the RTD matrix and the row sum of the matrix, along with the graph.

The  RTD Matrix for $\alpha = 1$         Row sum of the RTD matrix

$$\begin{bmatrix} -1 & 0 & -1 & 0 & 0 & 0 \\ 0 & 1 & 1 & 0 & 0 & 0 \\ 1 & 0 & 0 & 1 & 1 & 1 \\ 0 & 0 & 0 & -1 & -1 & -1 \end{bmatrix} \qquad \begin{bmatrix} -2 \\ 2 \\ 4 \\ -3 \end{bmatrix}$$

The graph of the fuzzy matrix for $\alpha = 1$ is given below:

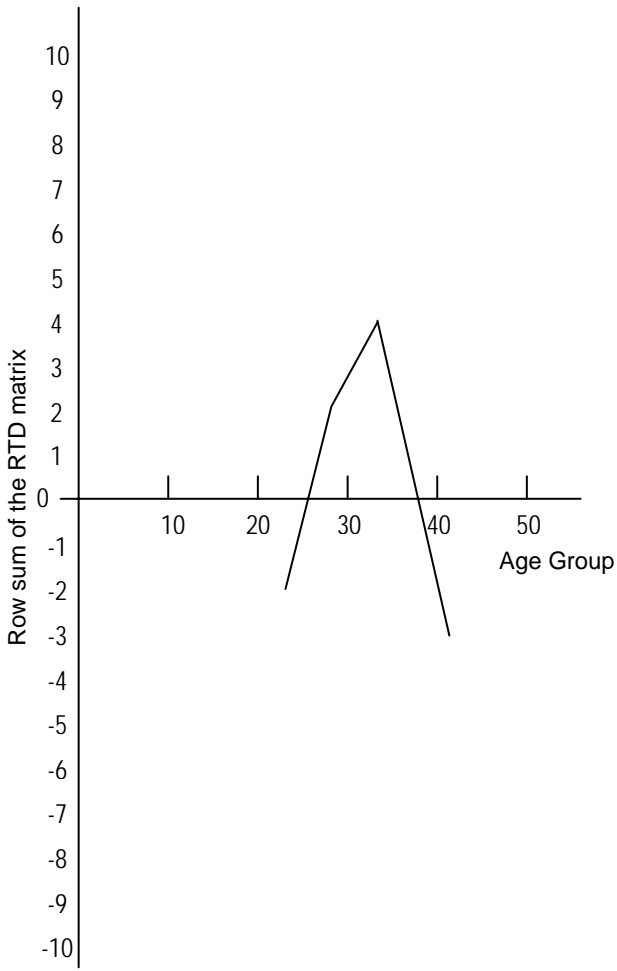

The age-group 31-35 is very highly affected age group followed by 26-30. The deviation of 20-25 and 36-47 are equally placed which implies people do not become victims of HIV/AIDS in the age group 20-25 and 36-47. Now we give the

Combined Time Dependent Matrix.

$$\begin{bmatrix} -3 & -2 & -3 & -2 & 0 & -1 \\ 2 & 3 & 3 & 0 & -1 & 0 \\ 3 & 0 & 2 & 3 & 3 & 3 \\ -2 & -1 & -1 & -3 & -3 & -3 \end{bmatrix}$$

Row sum of the CETD matrix

$$\begin{bmatrix} -11 \\ 7 \\ 14 \\ -13 \end{bmatrix}$$

The migrant labourers with HIV/AIDS is very low in the age group 20-25 there by it implies that migration does not start when the men are below 20 years, considering the fact that it takes at least four to five years for the disease to reach recognizable and advanced stages.

Having worked with four sets of age group we felt it would be appropriate to try by further increasing the number of age groups i.e. the time period to test whether more refinement may give more sensitive prediction. So we now work with 6 sets of age groups, 20-23, 24-30, 31-34, 35-37, 38-40 and 41-47 and analyze them using the same type of matrix models.

It is pertinent to mention here that our first model with 3 age groups could only suggest one class of age groups who are most vulnerable to the disease whereas the second model with four age groups suggests two sets of age groups which are vulnerable to the disease with gradations.

So now we will work with by redefining the age group to six intervals and draw conclusions from the model.

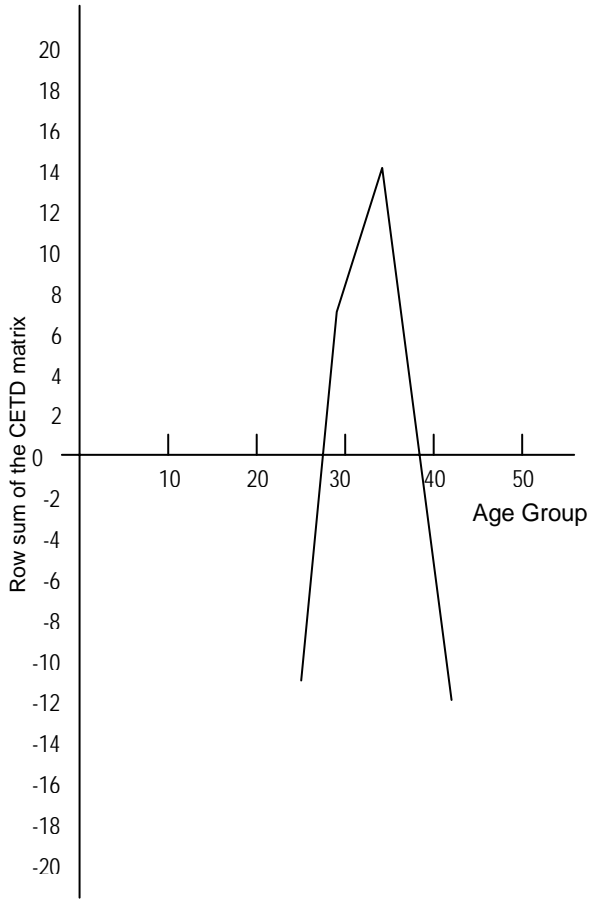

## Initial Raw Data Matrix

| Age Group | CSWs | Other Women | Smoking | Alcohol | Bad Company | Quacks |
|-----------|------|-------------|---------|---------|-------------|--------|
| 20-23 | 3 | 2 | 2 | 3 | 3 | 2 |
| 24-30 | 20 | 6 | 15 | 16 | 14 | 10 |
| 31-34 | 15 | 4 | 13 | 14 | 9 | 6 |
| 35-37 | 8 | 2 | 6 | 6 | 3 | 2 |
| 38-40 | 6 | 1 | 4 | 4 | 2 | 1 |
| 41-47 | 8 | 1 | 5 | 6 | 3 | 2 |

ATD Matrix

| Age Group | CSWs | Other Women | Smoking | Alcohol | Bad Company | Quacks |
|-----------|------|-------------|---------|---------|-------------|--------|
| 20-23 | 0.75 | 0.5 | 0.5 | 0.75 | 0.75 | 0.5 |
| 24-30 | 3.33 | 1 | 2.5 | 2.67 | 2.33 | 1.67 |
| 31-34 | 3.75 | 1 | 3.25 | 3.5 | 2.25 | 1.5 |
| 35-37 | 2.67 | 0.67 | 2 | 2 | 1 | 0.67 |
| 38-40 | 2 | 0.33 | 1.33 | 1.33 | 0.67 | 0.28 |
| 41-47 | 1.14 | 0.14 | 0.71 | 0.86 | 0.43 | 0.28 |

Average and Standard Deviation of the data

| Average | 2.27 | 0.61 | 1.72 | 1.85 | 1.24 | 0.83 |
|---------|------|------|------|------|------|------|
| Standard Deviation | 1.09 | 0.32 | 0.97 | 0.99 | 0.76 | 0.55 |

The RTD matrix for $= 0.5$

$$\begin{bmatrix} -1 & 0 & -1 & -1 & -1 & -1 \\ 1 & 1 & 1 & 1 & 1 & 1 \\ 1 & 1 & 1 & 1 & 1 & 1 \\ 0 & 0 & 0 & 0 & 0 & 0 \\ 0 & -1 & 0 & -1 & -1 & -1 \\ -1 & -1 & -1 & -1 & -1 & -1 \end{bmatrix}$$

Row sum of the RTD or fuzzy matrix

$$\begin{bmatrix} -5 \\ 6 \\ 6 \\ 0 \\ -4 \\ -6 \end{bmatrix}$$

Graph of the row sum of the RTD matrix is given for is $\alpha = 0.5$

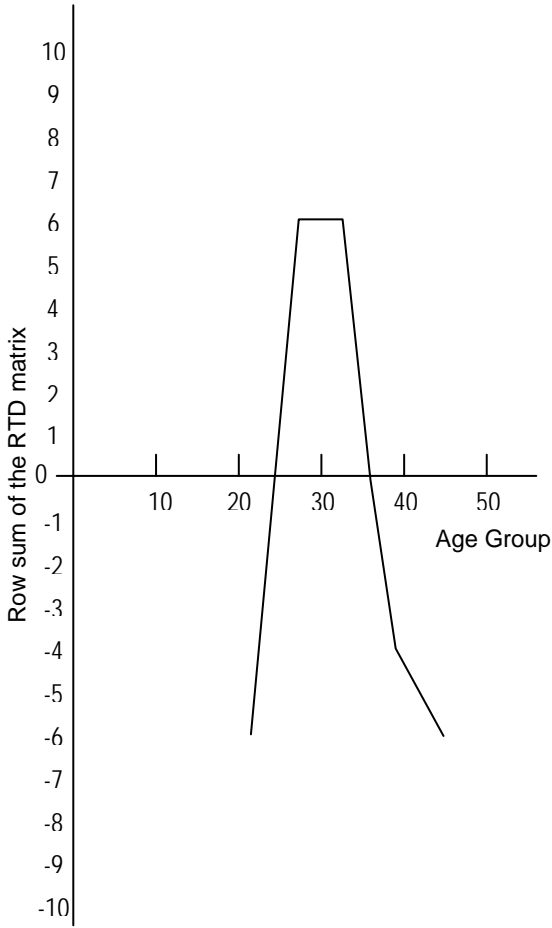

The maximum vulnerable age group being 24-30 and 31-34. The RTD matrix for the parameter $\alpha = 0.2$ is given below:

$$\begin{bmatrix} -1 & -1 & -1 & -1 & -1 & -1 \\ 1 & 1 & 1 & 1 & 1 & 1 \\ 1 & 1 & 1 & 1 & 1 & 1 \\ 1 & 0 & 1 & 0 & -1 & -1 \\ 0 & -1 & -1 & -1 & -1 & -1 \\ -1 & -1 & -1 & -1 & -1 & -1 \end{bmatrix}$$

The row sum of the RTD matrix

$$\begin{bmatrix} -6 \\ 6 \\ 6 \\ 0 \\ -5 \\ -6 \end{bmatrix}$$

We see that there are nil or very negligible HIV/AIDS patients in the age group 20-23, as its value is - 6. The graph for the parameter $\alpha = 0.2$ of the RTD matrix is given below:

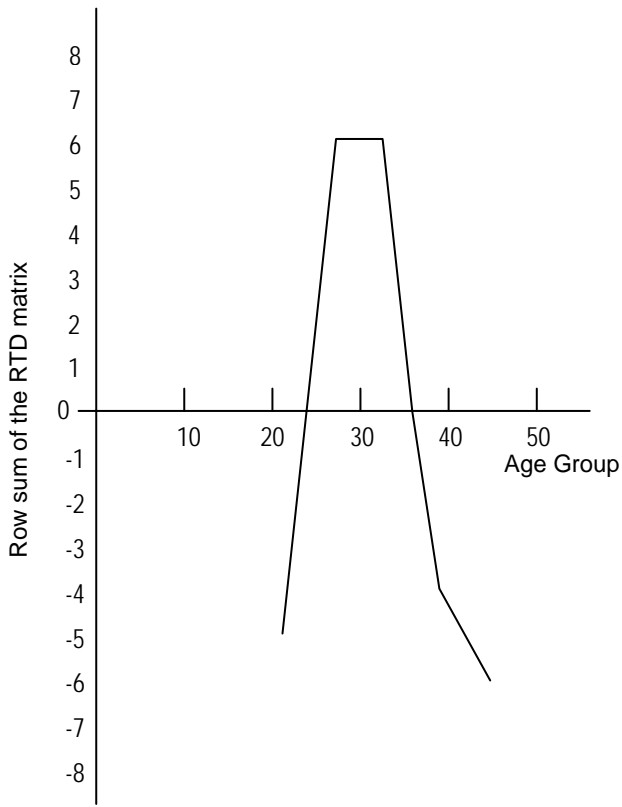

For α = 1, the RTD matrix is          Row sum of RTD matrix

$$\begin{bmatrix} -1 & 0 & -1 & -1 & 0 & 0 \\ 0 & 1 & 0 & 0 & 1 & 1 \\ 1 & 1 & 1 & 1 & 1 & 1 \\ 0 & 0 & 0 & 0 & 0 & 0 \\ 0 & 0 & 0 & 0 & 0 & 0 \\ -1 & -1 & -1 & -1 & -1 & -1 \end{bmatrix} \qquad \begin{bmatrix} -3 \\ 3 \\ 6 \\ 0 \\ 0 \\ -6 \end{bmatrix}$$

Graph of the RTD matrix for α = 1

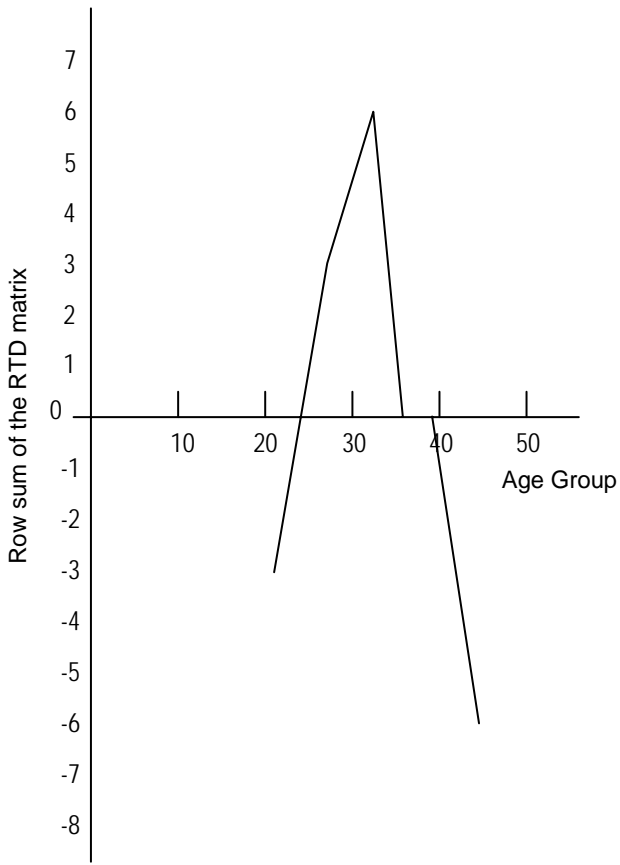

CETD matrix                    Row sum of CETD Matrix

$$\begin{bmatrix} -3 & -1 & -3 & -3 & -2 & -2 \\ 2 & 3 & 2 & 2 & 3 & 3 \\ 3 & 3 & 3 & 3 & 3 & 3 \\ 1 & 0 & 1 & 0 & -1 & -1 \\ 0 & -2 & -1 & -2 & -2 & -2 \\ -3 & -3 & -3 & -3 & -3 & -3 \end{bmatrix} \qquad \begin{bmatrix} -14 \\ 15 \\ 18 \\ 0 \\ -9 \\ -18 \end{bmatrix}$$

Graph of CETD Matrix

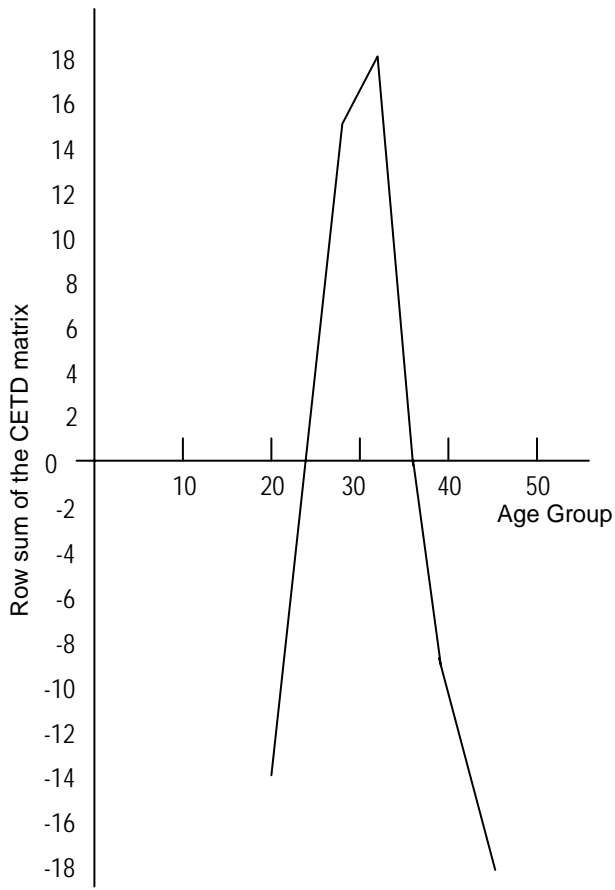

## 2.2 ANALYSIS OF SOCIO ECONOMIC POSITION OF THE HIV/AIDS PATIENTS USING FCM

Fuzzy cognitive maps (FCMs) are more applicable when the data in the first place is an unsupervised one. The FCMs work on the opinion of experts. FCMs model the world as a collection of classes and causal relations between classes.

FCMs are fuzzy signed directed graphs with feedback. The directed edge $e_{ij}$ from causal concept $C_i$ to concept $C_j$ measures how much $C_i$ causes $C_j$. The time varying concept function $C_i(t)$ measures the non-negative occurrence of some fuzzy event, perhaps the strength of a political sentiment, historical trend or military objective.

FCMs are used to model several types of problems varying from gastric-appetite behavior, popular political developments etc. FCMs are also used to model in robotics like plant control.

The edges $e_{ij}$ take values in the fuzzy causal interval $[-1, 1]$. $e_{ij} = 0$ indicates no causality, $e_{ji} > 0$ indicates causal increase, $C_j$ increases as $C_i$ increases (or $C_j$ decreases as $C_i$ decreases). $e_{ji} < 0$ indicates causal decrease or negative causality. $C_j$ decreases as $C_i$ increases (or $C_j$ increases as $C_i$ decreases). Simple FCMs have edge values in $\{-1, 0, 1\}$.

Then if causality occurs, it occurs to a maximal positive or negative degree.

Simple FCMs provide a quick first approximation to an expert stand or printed causal knowledge.

We illustrate this by the following, which gives a simple FCM of a socio-economic model.

This socio-economic model is constructed with Population, Crime, Economic condition, Poverty and Unemployment as nodes or concepts.

Here the simple trivalent directed graph is given by the following figure which is the expert's opinion.

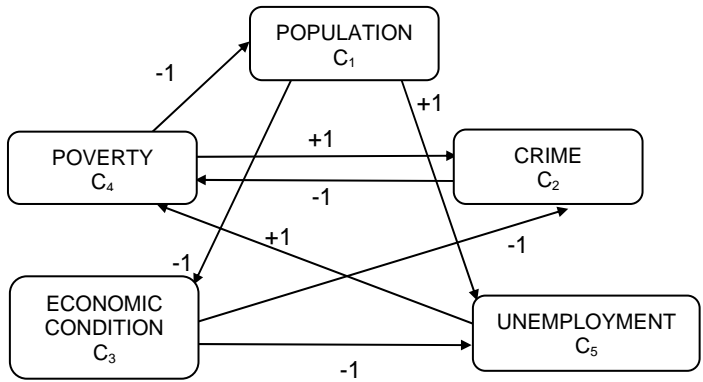

**FIGURE: 2.2.1**

Causal feedback loops abound in FCMs in thick tangles. Feedback precludes the graph-search techniques used in artificial-intelligence expert systems. FCMs feedback allows experts to freely draw causal pictures of their problems and allows causal adaptation laws, infer causal links from simple data. FCM feedback forces us to abandon graph search, forward and especially backward chaining. Instead we view the FCM as a dynamical system and take its equilibrium behavior as a forward-evolved inference. Synchronous FCMs behave as Temporal Associative Memories (TAM). We can always, in case of a model, add two or more FCMs to produce a new FCM. The strong law of large numbers ensures in some sense that knowledge reliability increases with expert sample size. We reason with FCMs. We pass state vectors C repeatedly through the FCM connection matrix E, thresholding or non-linearly transforming the result after each pass. Independent of the FCMs size, it quickly settles down to a temporal associative memory limit cycle or fixed point which is the hidden pattern of the system for that state vector C. The limit cycle or fixed-point inference summarizes the joint effects of all the interacting fuzzy knowledge.

**DEFINITION 2.2.1:** *An FCM is a directed graph with concepts like policies, events etc. as nodes and causalities as edges. It represents causal relationship between concepts.*

The following example illustrates how FCMs can be used to model social problems.

***Example 2.2.1****:* In Tamil Nadu (a southern state in India) in the last decade several new engineering colleges have been approved and started. The resultant increase in the production of engineering graduates in these years is disproportionate with the need of engineering graduates.

This has resulted in thousands of unemployed and underemployed graduate engineers. Using an expert's opinion we study the effect of such unemployed people on the society. An expert spells out the five concepts relating to the unemployed graduated engineers as

$E_1$  –  Frustration
$E_2$  –  Unemployment
$E_3$  –  Increase of educated criminals
$E_4$  –  Under employment
$E_5$  –  Taking up drugs etc.

The directed graph where $E_1$, …, $E_5$ are taken as the nodes and causalities as edges as given by an expert is given in the following figure:

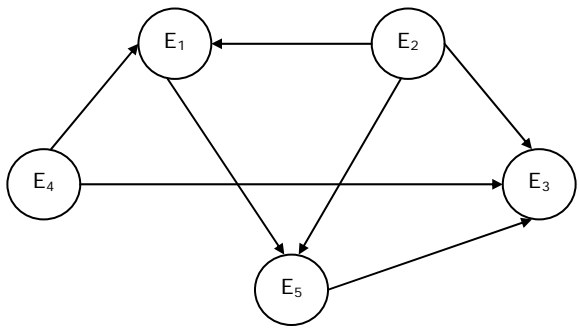

FIGURE: 2.2.2

According to this expert, increase in unemployment increases frustration. Increase in unemployment, increases the educated criminals. Frustration increases the graduates to take up to evils like drugs etc.

Unemployment also leads to the increase in number of persons who take up to drugs, drinks etc. to forget their worries and unoccupied time. Under-employment forces them to do criminal acts like white-collar crimes for want of more money and so on. Thus one cannot actually get data for this but can use the expert's opinion for this unsupervised data to obtain some idea about the real plight of the situation.

This is just an illustration to show how FCM is described by a directed graph.

{If increase (or decrease) in one concept leads to increase (or decrease) in another, then we give the value 1. If there exists no relation between two concepts the value 0 is given. If increase (or decrease) in one concept decreases (or increases) another, then we give the value –1. Thus FCMs are described in this way.}

**DEFINITION 2.2.2:** *When the nodes of the FCM are fuzzy sets then they are called as fuzzy nodes.*

**DEFINITION 2.2.3:** *FCMs with edge weights or causalities from the set {–1, 0, 1}, are called simple FCMs.*

**DEFINITION 2.2.4:** *Consider the nodes/ concepts $C_1$, …, $C_n$ of the FCM. Suppose the directed graph is drawn using edge weight $e_{ij}$ $\in$ {0, 1, –1}. The matrix E be defined by $E = (e_{ij})$ where $e_{ij}$ is the weight of the directed edge $C_iC_j$ . E is called the adjacency matrix of the FCM, also known as the connection matrix of the FCM.*

It is important to note that all matrices associated with an FCM are always square matrices with diagonal entries as zero.

**DEFINITION 2.2.5:** *Let $C_1$, $C_2$, … , $C_n$ be the nodes of an FCM. A = ($a_1$, $a_2$, … , $a_n$) where $a_i$ $\in$ {0, 1}. A is called the instantaneous state vector and it denotes the on-off position of the node at an instant.*

$a_i = 0$ *if $a_i$ is off and*
$a_i = 1$ *if $a_i$ is on*
*for i = 1, 2, …, n.*

**DEFINITION 2.2.6:** *Let $C_1$, $C_2$, …, $C_n$ be the nodes of an FCM. Let be the edges of the FCM (i≠ j). Then the edges form a directed cycle. An FCM is said to be cyclic if it possesses a directed cycle. An FCM is said to be acyclic if it does not possess any directed cycle.*

**DEFINITION 2.2.7:** *An FCM with cycles is said to have a feedback.*

**DEFINITION 2.2.8:** *When there is a feedback in an FCM, i.e., when the causal relations flow through a cycle in a revolutionary way, the FCM is called a dynamical system.*

**DEFINITION 2.2.9:** *Let* $\overline{C_1C_2}, \overline{C_2C_3}, ..., \overline{C_{n-1}C_n}$ *be a cycle. When $C_i$ is switched on and if the causality flows through the edges of a cycle and if it again causes $C_i$, we say that the dynamical system goes round and round. This is true for any node $C_i$, for i = 1,2,..., n. The equilibrium state for this dynamical system is called the hidden pattern.*

**DEFINITION 2.2.10:** *If the equilibrium state of a dynamical system is a unique state vector, then it is called a fixed point.*

**Example 2.2.2:** Consider a FCM with $C_1$, $C_2$, …, $C_n$ as nodes. For example let us start the dynamical system by switching on $C_1$. Let us assume that the FCM settles down with $C_1$ and $C_n$ on i.e. the state vector remains as (1, 0, 0, …, 0, 1) this state vector (1, 0, 0, …, 0, 1) is called the fixed point.

**DEFINITION 2.2.11:** *If the FCM settles down with a state vector repeating in the form*

$$A_1 \rightarrow A_2 \rightarrow ... \rightarrow A_i \rightarrow A_1$$

then this equilibrium is called a limit cycle.

Methods of finding the hidden pattern are discussed in the following section.

**DEFINITION 2.2.12:** *Finite number of FCMs can be combined together to produce the joint effect of all the FCMs. Let $E_1$, $E_2$,..., $E_p$ be the adjacency matrices of the FCMs with nodes $C_1$, $C_2$, …, $C_n$ then the combined FCM is got by adding all the adjacency matrices $E_1$, $E_2$, …, $E_p$.*

We denote the combined FCM adjacency matrix by $E = E_1 + E_2 + … + E_p$.

**NOTATION:** Suppose A = $(a_1, … , a_n)$ is a vector which is passed into a dynamical system E. Then AE = $(a'_1, … , a'_n)$ after thresholding and updating the vector suppose we get $(b_1, \ldots, b_n)$ we denote that by $(a'_1, a'_2, \ldots, a'_n) \hookrightarrow (b_1, b_2, \ldots, b_n)$.

Thus the symbol '$\hookrightarrow$' means the resultant vector has been thresholded and updated.

FCMs have several advantages as well as some disadvantages. The main advantage of this method: it is simple. It functions on expert's opinion. When the data happens to be an unsupervised one the FCM comes handy. This is the only known fuzzy technique that gives the hidden pattern of the situation. As we have a very well known theory, which states that the strength of the data depends on, the number of experts' opinion we can use combined FCMs with several experts' opinions.

At the same time the disadvantage of the combined FCM is when the weightages are 1 and –1 for the same $C_i$ $C_j$, we have the sum adding to zero thus at all times the connection matrices $E_1, \ldots, E_k$ may not be conformable for addition.

Combined conflicting opinions tend to cancel out and assisted by the strong law of large numbers, a consensus emerges as the sample opinion approximates the underlying population opinion. This problem will be easily overcome if the FCM entries are only 0 and 1.

We have just briefly recalled the definitions. For more about FCMs please refer Bart Kosko [58, 62].

The study of the HIV/AIDS affected patients who have contracted the disease by migrancy involves reasons like socially free, socially irresponsible behaviour, above all men with no education or ethics or values with no motivation, no care for the family, no mind to think about their children or wife, with all addictive habits are mainly vulnerable to HIV/AIDS. This was explicitly recorded when we met 60 men who had acquired the disease due to migration.

Most of them say with pride that they were infected by CSWs in the Mumbai city. Almost all of them barring one who is a cook who got infected from his wife, proudly disclose that the number of times they have visited the CSWs is "innumerable".

This throws light on their absence of social responsibility, patronage of commercial sex workers and the underlying machismo. They acknowledge without the least bit of hesitancy or reserve that they have/had addiction to alcohol and smoking. There are others who say that they had visited CSWs and also taken precautionary methods to prevent the infection, and they also say that because they were drunk sometimes they might have forgotten to wear condoms etc. Most of them quite freely talk about the CSWs whom they visited.

It is important to mention at least 4% of them acknowledged that they are reaping the consequences of their deeds. Some of them were stoic. When we asked what the government should do for them most of them demanded that the government should 'fastly' discover some medicine to cure them all, except only one individual who was of the opinion that discovery of any medicine to cure the HIV/AIDS patients will result in complete degeneration of moral values, hence government should not find any cure for HIV/AIDS.

Another factor to be worth mentioning is that only one out of sixty who had contracted the disease from his wife said that whether horoscope was found to be matching or not one should test for HIV/AIDS before marriage; for no fault of his he had acquired this disease. Since the feeling of each and every individual is so sensitive and intricate to be given any mathematical expression other than fuzzy mathematics; as it is the only tool that can give complete expression to such feelings, we are justified in approaching this problem using fuzzy mathematics and neutrosophic theory.

Here we approach the problem via attributes using Fuzzy Cognitive Maps (FCMs) that are basically matrices which predict the feelings of all the attributes under consideration. It is still unfortunate and unbearable to state that of the 60 infected men most of them are very reluctant to carry out the HIV/AIDS test for their wife for they are in the first place unconcerned and secondly they fear disclosing it to their wife for fear of their reaction. It is important to mention here that some of them answered several questions and some of the answers were inderterminate which will analyzed using neutrosophic cognitive maps in chapter IV.

Before we proceed on to apply the Fuzzy Cognitive Map model to this problem we define or precisely recall what we mean by the 12 attributes given by the expert. So that when we work with it the reader will have no confusion in noting each attribute and in analyzing them via the directed graph and its connection matrices.

## $A_1$ - Easy Availability of Money

Most of these men worked for eight or ten hours and earn a bulk amount from Rs.60 to Rs.200 per day. Because of the regular and the easy availability of this money, these men lavished it on drinks, smoking, and in visiting the CSWs. Since they lacked other social commitments and obligations, these daily wagers spent their money in all inappropriate ways, did not save or give considerable amount to their family members.

## A₂- Lack of Education

The term lack of education when taken as a attribute needs to be analyzed or given a different interpretation for when we interviewed the 60 HIV/AIDS patients we saw that 15 have not even attended the primary school.

Ten of them had studied standards 2, 3, 4 or 5. 14 of them have studied up to 6, 7, or 8 or 9. A few either completed 10[th] or failed 10[th]. Some completed 11[th] or failed 11[th]. But categorically when we say no education, we also mainly mean absence of moral education: all them had little or no moral values, social conscience, or social responsibility or an iota of self analysis to arrive at a conclusion whether they were correct or wrong in carrying out such a life-style. Miserably with no proper knowledge of HIV/AIDS became victims of it.

## A₃ - Visiting CSWs

When we studied the sample of 60 HIV/AIDS male patients we found that of the 60, 59 of them accepted their constant visit to CSW not only once or twice but regularly. Of these '55' of the men said they have visited CSWs innumerable number of times.

Here few things are worth mentioning: Most of them had not acquired the disease from Tamil Nadu though the HIV/AIDS patients are currently taking treatment in Tamil Nadu. One of the most fascinating spot for them was Mumbai and they paid Rs.50 to Rs.100 per visit to the CSWs there. Some of them had frequently visited Andhra Pradesh also. In Andhra Pradesh the CSWs were very cheap, ranging only from Rs.20 to Rs.50 and most of them were only aged between 18-23.

Another shocking information we derived from them was that most of the Andhra Pradesh based CSWs had STD/VD even at such early ages, several of these migrant patients have also taken medicine for the same after contracting it from them. Another notable point is only when they are in other states they do not have any fear or tension of being noted or observed by friends or relatives: their risk-behaviour heightens during outstation travel. Yet another information given by them is that the CSWs from the Chennai city are a bit costly ranging from Rs.200 to Rs.1000 or sometimes even more, so they cannot afford for the same. They also said that commercial sex workers in Chennai was delocalized and that there were a lot of brokers, middlemen etc. to fix them at a cheap rate. They (CSWs) in Chennai city carried out their trade with the help of police and brokers very successfully in all places which is well known to the interested.

## A₄ - Profession

Since almost all the 60 male patients had not entered college their profession varied from truck or lorry drivers, weavers, contract labourers, most of them were agriculture labourers who have taken up these profession mainly due to failure of agriculture or due to poor yield. Here also though they were weavers these persons were mostly businessmen who sold stitched goods or un-stitched goods from one state to other. One was a cook, one person was never employed but enjoyed life in cities like Bombay, Hyderabad and Nagari. One was a milkman (a resident of Chennai city) contracted the disease in Namakkal, when he went to buy a cow. Thus these people most of them had forcibly taken some other profession due to agriculture failure or no proper yield. It is not wrong if we say that mainly their change of profession made them as HIV/AIDS victims by giving them easy money and separation from family, here family means not only for married men who live away from their wives but also men who live away from their parents. Finally some (2 or 3) of them were runaways and they were platform dwellers.

## A₅ - Wrong/Bad Company

This term needs lot of explanation as an attribute for it carried lot of effect on the individual becoming an HIV/AIDS victim. In the first place wrong company does not mean all are bad: we just use this broad category to stress on negative peer-group influence that leads these youth astray. Youngsters join together to do deeds which have no moral or material significance, like first trying to smoke, trying to take drugs, trying to consume alcohol, seeing pornographic movies , pickpockets finally visiting CSWs. Thus bad company not only induces others in that gang of friends to act in this indiscrete manner.

It was really surprising for us when several of the drivers aged from 30 to 40 years humbly acknowledged that it was the bad company that landed them to visit CSWs. They proudly say this sort of friendship was to "enjoy life", "just in jolly mood" and so on.

## $A_6$ - Addiction to habit-forming substances and visiting CSWs

By this term we mean smoking, alcoholism, drug addiction, ganja and visiting CSWs. At the first stage when enquired we were made to understand that most of these HIV/AIDS patients had taken up smoking for the sake of male identity, as a display of aggressiveness. Then it became a habit and they could not come out of it.

It is pertinent to mention here that most of them had acquired this habit due to the bad company. The alcohol or 'drink' was attributed to their physical labour; they said that they drank after a day's long labour for relief from bodyache. Only few were ganja or drug addicts. Here it has become pertinent to mention that in Tamil Nadu (unlike in certain states like Manipur etc. where the spread of HIV/AIDS is due to Injecting Drug Users (IDUs)) the spread of HIV/AIDS is mainly through their heterosexual partners.

Some of them, nearly 50% acknowledged that only under intoxication they had visited CSWs so they were out of their mind to use any safe sex methods in spite of the CSWs insisting them to use. Thus to some extent it is not wrong if we say that these patients had gone off the track due to alcohol. Some said that their profession damaged their health and raised a lot of the heat of their bodies--so they had to particularly have intercourse in order to 'cool' their bodies.

The shocking observation was that they are proud to say that they have all such habits, it is still surprising to see that not only they take pride in doing so, but that even their wives are very proud to state that their husbands have all these addictive and bad habits. It is very sad and unless this attitude is changed it may be impossible for the government or NGOs to prevent the spread of HIV/AIDS whatever precautionary steps they may take.

## $A_7$ - Absence of Social Responsibility

It has become a fashion of the day for men to behave in an irresponsible and careless way and their family or wife takes it not as a mistake or a sin but on the other hand accepts it. Parents are to some extent to be blamed for the same as they do not bring up (especially) their male children properly: instead, male children are pampered, are not taken to task for making mistakes etc.. This has twofold effect on society (1) they are irresponsible to the extent of contracting the disease and on the other hand, they are even more irresponsible for they pass it on wilfully to their life partners and (there by) to their children. They do not even think for a minute that for no sin their wife and children become victims of HIV/AIDS.

(2) The female children are so rigidly brought up by our society, especially in rural areas of Tamil Nadu, they are not allowed to question even the improper behaviour of the males, they are taught the 'virtues' of being silent and all-suffering and as a result that they not only accept the irresponsible nature of their husbands but take even the disease in a very passive way.

None of them ever had the mind to question his character or even give derogatory remarks on their husband's bad activity. On the other hand, even when the wife was infected through her husband and whether the husband was dead or living with HIV/AIDS, the in-laws had the courage and arrogance to say that their son was good it was only the 'loose' character of their daughter-in-law which had infected him.

A few of them had even beaten up their daughter-in-laws. Even if we take into account the fact that it is also equally possible for the women to have had extramarital and unprotected sex, it does not really account for the husband's condemnable habit of visiting CSWs. Thus all these complacent attitudes has made men more socially irresponsible.

## $A_8$ - Socially Free

This attribute as a fuzzy node needs some more explanation. When we say 'socially free' we mean this only for the males in the patriarchal Indian society, who have greater freedom within the society.

For in the first stage when a male child is brought up he is free as ever in contrast with a female child who is from the time of her birth not free and is enslaved to some patriarchal control. None of the male child's behaviour is questioned, he is pampered beyond limits, he is petted, called the 'darling' of the family, and right from a very young age, the parents put up even with wrong behaviour among male children. Secondly the male child is given good food and professional education in contrast with the female child who may or may not be even given education depending on their social and economic conditions. Then when the male is married a huge wealth is demanded from the female partner as if he is incapable of supporting the family and the parents of female child have to bear the money/dowry. Thus all these makes the male of society very irresponsible and socially free for whatever way he leads his life his marriage is to compensate the money he needs.

This pattern of selective social freedom (given only to males, and at the same time denied to the females) makes the men in our society even more vulnerable to HIV. Unless this pattern of the society is changed it may not be possible for us to implement any awareness program or stop the spread of the disease. It is high time for the male members of the society to stop for the minute to think that if their female partners were to behave in an identical irresponsible manner as they do what would be their situation, unless self-realization comes to each and every male citizen of our nation, we would be reduced to the state of another South Africa.

## $A_9$ - Economic Status

This term when we accept as a fuzzy node has a very different meaning from the term "Economic status" in general. This will be explained shortly. This does not imply the economic conditions under which the family is living namely as poor, very poor, middle class, lower middle class, upper middle class, rich, very rich, or the hereditary richness or economic status of a family or community. But here by the term the economic status of a person we mean the individual economic status: what he has earned for the day and how lavishly it is spent on him for somke, drink, food and CSWs without a second thought. For an individual may be in the age group from 20 years to 48 years we say has status as individual depends on the money he earns per day. As per our observation these 60 men whom we have interviewed did not have any craze to dress well (i.e. even meagre clean clothes) or any interest in jewellery (like ring or watch or a chain) or savings but what was their economic status: then it is in their capability to spend money on smoking, drinks, CSWs and in bad company.

When they are irresponsible they psychologically suffer to maintain their male identity by aggressive attitudes like smoking, alcohol and CSWs. Thus by economic status we do not view it collectively, but only as individuals. And this includes their capacity to spend on CSWs, smoking, drinks, drugs etc.

## A₁₀ - More leisure

When we say more leisure as a fuzzy node we clearly mean the following: These persons have come from their hometown or villages to carry out some work, be as a truck driver or a lorry driver or as a contract labour or as a business errand to self or deliver the goods.

So once their work is over the rest of the time they are only free or they have more leisure. To while away their leisure time they visit the CSWs.

We wish to point out that because they cannot spend/while away their time, and because of their loneliness in a new place, they visit the CSWs. This is especially true of those who have migrated to other cities, and stay alone and have a job there. After the working hours, they are brazenly free.

Thus they have "more leisure" which is only and mainly utilized for visiting CSWs.

## A₁₁ - Machismo/ Exaggerated Masculinity

This term is very sensitive but at the same time very essential for us to describe it in the context of HIV/ AIDS patients.

So as a fuzzy node when we say "Machismo" (meaning, exaggerated masculinity or alternatively male ego) we mean exaggerated, emphasized masculinity, macho behaviour and the following:

(i)     they take pride in stating that they have all addictive habits like CSWs, drink and smoking.

(ii)    They are "socially free".

(iii)   They are socially irresponsible.

(iv)    They have no conscience of how they spend the money.

(v)     They have no binding on the family members, spend money only on themselves.

The main reason for all this is due to their upbringing and the existing patriarchy.

## $A_{12}$ - No awareness of the disease

It is a pity that the infected people are unaware of the disease, its mode of spread, cure etc..

When we say they are unaware of the disease we as a fuzzy node mainly mean this:

1. They are still not aware of the fact that this disease does not have a cure. That is why of the 60 persons 59 of them said that the government should find medicines soon to cure them.

2. They are not aware of the relation between HIV/AIDS and their main source of infection is CSWs.

3. They (most of them) are aware of STD / VD for they are at one stage or other victims of these diseases. Infact most of them had suffered this disease very frequently.

4. As all diseases like cancer, diabetes are incurable, so is HIV/AIDS. But this attitude has to change and people must understand that HIV/AIDS is now a livable chronic disease very frequently.

   But it is unfortunate that the doctors when they treat these patients for STD/VD do not make them know or convince them that they may end up with HIV/AIDS if they continue to have these habits. The doctors can greatly intervene at this stage, and counsel the patients about risk-behaviour patterns, advice them the importance of safe-sex methods and also try to slowly wean them away from this practice of visiting CSWs. Unless the doctors/ pharmacists make them fear about the dreadfulness of HIV/AIDS and tell them its methods of spreading etc. it is not possible to prevent the further spread of HIV/AIDS.

5. The society at large is frightened of this disease more as a social stigma but are not fully aware of it, at large in

the rural areas where the majority of them are uneducated.

6. Some of them acknowledge that they were somewhat aware of HIV/AIDS and even the CSWs advised them for safe methods yet because of their macho behaviour and their male ego, they did not listen to the CSWs. This is doubly risky for it not only puts these men at risk, but also the CSWs at high risk.

7. Some women said that awareness must not only be made through television or posters but should be made through street dramas or songs or radio as most of them were uneducated so they cannot read the posters or some of them said that TV was not in their reach.

Now taking these 12 attributes as fuzzy nodes we proceed to give the directed graphs using the experts opinion.

Here it is important to mention that we do not say we have exhausted all the attributes related to these patients or it is not necessary for one to take all 12 attributes the expert can choose to take only a few out of them or more attributes. Thus it is purely in the hands of the investigator to take the related attributes as per his/her wishes.

Further we by no means guarantee that the directed graph has any motive of the author for they were impersonally taken from the experts where the experts are taken as NGOs, doctors, public and the patients themselves.

Here the term 'expert' need not strictly mean an expert: it means someone who gives his views to us and has an understanding and expertise on the said disease and the related issues.

However using the fuzzy directed graph we obtain the connection matrix. Using these connection matrix we determine the hidden pattern.

Here we make use of only simple FCMs. The directed graph given by the first expert with these 12 attributes is as follows:

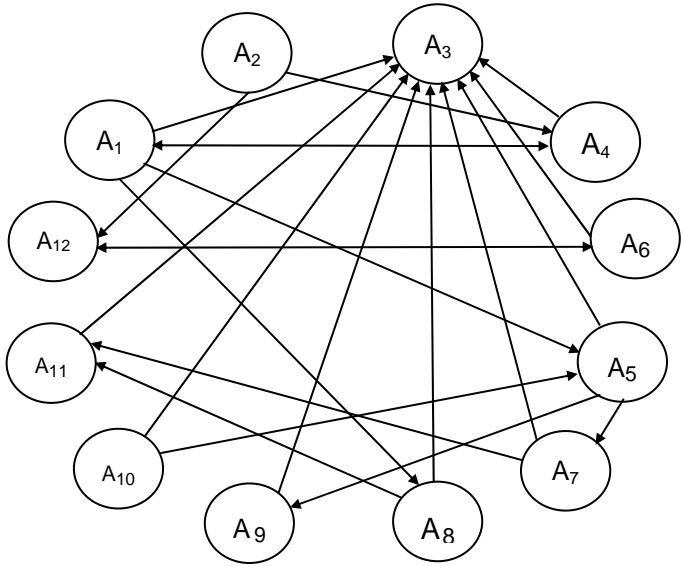

DIRECTED GRAPH OF THE FCM – FIGURE 2.2.3

The connection matrix D of the directed graph

$$
\begin{array}{c}
\quad\ \ A_1\ A_2\ A_3\ A_4\ A_5\ A_6\ A_7\ A_8\ A_9\ A_{10}\ A_{11}\ A_{12} \\
\begin{array}{c}
A_1 \\ A_2 \\ A_3 \\ A_4 \\ A_5 \\ A_6 \\ A_7 \\ A_8 \\ A_9 \\ A_{10} \\ A_{11} \\ A_{12}
\end{array}
\left[
\begin{array}{cccccccccccc}
0 & 0 & 1 & 1 & 1 & 0 & 0 & 1 & 0 & 0 & 0 & 0 \\
0 & 0 & 0 & 1 & 0 & 0 & 0 & 0 & 0 & 0 & 0 & 1 \\
0 & 0 & 0 & 1 & 0 & 0 & 0 & 0 & 0 & 0 & 0 & 0 \\
1 & 1 & 1 & 0 & 0 & 0 & 0 & 0 & 0 & 0 & 0 & 0 \\
0 & 0 & 1 & 0 & 0 & 0 & 0 & 0 & 1 & 0 & 0 & 0 \\
0 & 0 & 1 & 0 & 0 & 0 & 0 & 0 & 0 & 0 & 0 & 1 \\
0 & 0 & 1 & 0 & 1 & 0 & 0 & 0 & 0 & 0 & 1 & 0 \\
0 & 0 & 1 & 0 & 0 & 0 & 0 & 0 & 0 & 0 & 1 & 0 \\
0 & 0 & 1 & 0 & 1 & 0 & 0 & 0 & 0 & 0 & 0 & 0 \\
0 & 0 & 1 & 0 & 0 & 1 & 0 & 0 & 0 & 0 & 0 & 0 \\
0 & 0 & 1 & 0 & 0 & 0 & 0 & 0 & 0 & 0 & 0 & 0 \\
0 & 0 & 1 & 0 & 0 & 1 & 0 & 0 & 0 & 0 & 0 & 0
\end{array}
\right]
\end{array}
$$

Now using this connection matrix which we denote by D we find the effect on the state vectors. Let R = (0 0 0 1 0 0 0 0 0 0 0 0) i.e.

only the attribute profession is in the on state, we study the effect of the state vector R using the experts opinion

$$RD \quad \hookrightarrow \quad (1\ 1\ 1\ 1\ 0\ 0\ 0\ 0\ 0\ 0\ 0) \quad = \quad R_1$$

$$R_1D \quad \hookrightarrow \quad (1\ 1\ 1\ 1\ 1\ 0\ 0\ 1\ 0\ 0\ 1) \quad = \quad R_2$$

$$R_2D \quad \hookrightarrow \quad (1\ 1\ 1\ 1\ 1\ 1\ 0\ 1\ 1\ 0\ 1) \quad = \quad R_3$$

$$R_3D \quad \hookrightarrow \quad (1\ 1\ 1\ 1\ 1\ 1\ 0\ 1\ 1\ 0\ 1) \quad = \quad R_3$$

where "$\hookrightarrow$" denotes that the resultant vector that has been updated and thresholded. Thus the 'profession' especially mentioned as the attribute makes the on state all the nodes except $A_7$ and $A_{10}$. Thus their profession which is described in $A_4$ from the hidden pattern gives a fixed point. Thus in a way in case of migrant labourers the profession has a great impact on contracting HIV/AIDS.

Now suppose let us take the state vector $R_1 = (1\ 0\ 0\ 0\ 0\ 0\ 1\ 0\ 0\ 1\ 0\ 0)$ the effect of $R_1$ in which the vectors easy money, no social responsibility, and more leisure i.e. $A_1$, $A_7$ and $A_{10}$ are in the on state and all other states are off. The effect of $R_1$ on the dynamical system D is given by

$$R_1D \quad \hookrightarrow \quad (1\ 0\ 1\ 1\ 1\ 0\ 1\ 1\ 0\ 1\ 1\ 0) \quad = \quad R_2$$

$$R_2D \quad \hookrightarrow \quad (1\ 1\ 1\ 1\ 1\ 0\ 1\ 1\ 1\ 1\ 1\ 0) \quad = \quad R_3$$

$$R_3D \quad \hookrightarrow \quad (1\ 1\ 1\ 1\ 1\ 0\ 1\ 1\ 1\ 1\ 1\ 1) \quad = \quad R_4$$

$$R_4D \quad \hookrightarrow \quad (1\ 1\ 1\ 1\ 1\ 1\ 1\ 1\ 1\ 1\ 1\ 1) \quad = \quad R_5$$

$$R_5D \quad \hookrightarrow \quad (1\ 1\ 1\ 1\ 1\ 1\ 1\ 1\ 1\ 1\ 1\ 1) \quad = \quad R_5.$$

Thus when 3 nodes are on all, the nodes become on which is the hidden pattern. Thus when easy money, absence of social responsibility and more leisure is in the on state all the attributes become on. The expert was specially interested to study when Exaggerated Masculinity, i.e. $A_{11}$ alone is in the on state and all other states are off, the effect on the system D is:

$$\text{Let } X \quad = \quad (0\ 0\ 0\ 0\ 0\ 0\ 0\ 0\ 0\ 0\ 1\ 0)$$

$$XD \quad \hookrightarrow \quad (0\ 0\ 1\ 0\ 0\ 0\ 0\ 0\ 0\ 0\ 1\ 0) \quad = \quad X_1$$

$$X_1D \quad \hookrightarrow \quad (0\ 0\ 1\ 0\ 0\ 0\ 0\ 0\ 0\ 0\ 1\ 0) \quad = \quad X_1$$

thus if a person has only Exaggerated Masculinity, the person visits the CSW, nothing can be said or it has no relation to other attributes. Thus only the 'macho' person with positive attitude is not harmful according to this expert. Now let us consider the extreme case when X' = (0 0 0 0 0 0 0 0 0 0 0 1) i.e. only the person is not aware of anything about the disease is only in the on state and all other states are off. Now we study the effect of X' on D

$$X'D \quad \hookrightarrow \quad (0\ 0\ 1\ 0\ 0\ 1\ 0\ 0\ 0\ 0\ 0\ 1) \quad = X'_1$$
$$X'_1D \quad \hookrightarrow \quad (0\ 0\ 1\ 0\ 0\ 1\ 0\ 0\ 0\ 0\ 0\ 1) \quad = X'_1$$

i.e. when a person is unaware of the disease he visits CSW and has addictive habits according to this expert some times one may not be in a position to comprehend what the expert says on certain attributes, to overcome this problem usually we take the combined effect of all the experts opinion which gives a consolidated vision about the problem and its solution.

When we approached the second expert who is an NGO who has worked over a decade to spread awareness and who was himself a HIV/AIDS infected persons the following was given as the fuzzy directed graph. He did not take the twelve attributes as given by the first expert but he said he would not include the attribute 'no education' for, according to him both the educated and the uneducated behave in the same way except in case of uneducated they come to government hospitals for treatment and the educated get treatment in a different, secret and anonymous way. He says we are not certain about the percentage of educated persons infected with HIV/AIDS, they are cautious and don't come out.

Further this expert coupled the bad company and bad habits as a single attribute. Also according to this expert, no social responsibility and socially free are put in the same node for he feels a socially free person is socially irresponsible and so on. The concept of easy money and profession are put under one attribute for mainly the daily labourers are the persons who think the money to be an easily earned one, which reflects on their profession. According to him the economic status is relative, so he did not wish to deal with it as he is of the opinion the very rich are also easily susceptable to the disease but they take treatment staying in their homes or in big hospitals were confidentiality of the disease is well maintained in the hospital as well as by the doctors.

In case of upper middle class or middle class also where ever the disease is prevalent it is kept secret. In our opinion, a large percentage of those coming to the Tambaram Sanatorium for the treatment of the disease are poor uneducated rural people. This is because, unlike the rich for whom anonymity can be afforded, it is not the same way with those who are poor.

Thus this expert justifies his non-inclusion of the Economic Status as one of the attributes, for in his opinion the disease is prevalent in all economic status. It is only more well-known outside in the case of lower middle class or middle class. Thus this expert was contented and he justified for working with only these seven attributes the directed graph given by the second expert is as follows.

$A'_1$ - Easy Availability of Money
$A'_2$ - Wrongful company, Addictive habits
$A'_3$ - Visits CSWs
$A'_4$ - Socially irresponsible/free
$A'_5$ - Macho Behaviour
$A'_6$ - More leisure
$A'_7$ - No awareness about disease.

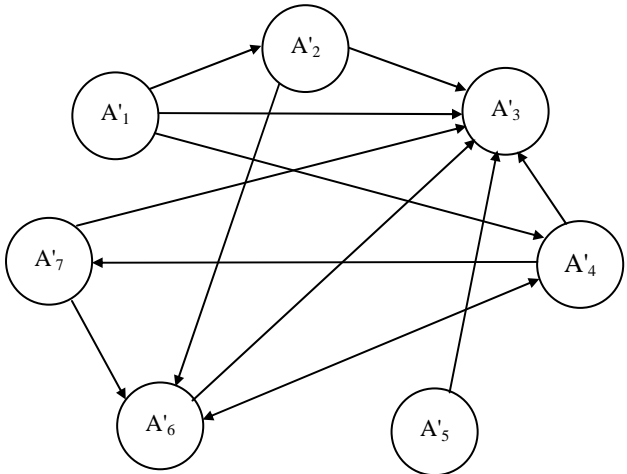

FIGURE 2.2.4

The related fuzzy relational matrix is as follows:

$$
\begin{array}{c}
& A'_1\; A'_2\; A'_3\; A'_4\; A'_5\; A'_6\; A'_7 \\
\begin{array}{c}
A'_1 \\ A'_2 \\ A'_3 \\ A'_4 \\ A'_5 \\ A'_6 \\ A'_7
\end{array}
\begin{bmatrix}
0 & 1 & 1 & 1 & 0 & 0 & 0 \\
0 & 0 & 1 & 0 & 1 & 1 & 0 \\
0 & 0 & 0 & 0 & 0 & 0 & 0 \\
0 & 0 & 1 & 0 & 0 & 1 & 0 \\
0 & 1 & 1 & 0 & 0 & 0 & 0 \\
0 & 1 & 1 & 1 & 0 & 0 & 0 \\
0 & 0 & 1 & 1 & 0 & 1 & 0
\end{bmatrix}
\end{array}
$$

Suppose A denotes the connection matrix of the directed graph. Now to find the stability of the dynamical system or to be more precise the hidden pattern of the system which may be a fixed point or a limit cycle .

Consider the state vector / initial vector

$X_1$ = (1 0 0 0 0 0 0) i.e. the only node 'easy availability of money which is dependent on the profession' i.e. $A'_1$ alone is in the on state, all other state vectors are in the off state. Now passing $X_1$ into the connection matrix A we get

$$X_1A \quad \hookrightarrow \quad (1\ 1\ 1\ 1\ 0\ 0\ 0)$$

where ' $\hookrightarrow$ ' denotes the resultant of the vector Xi after passing through A that has been updated and thresholded.

Let

$$X_1A \quad \hookrightarrow \quad (1\ 1\ 1\ 1\ 0\ 0\ 0) \quad = \quad X_2$$
$$X_2A \quad \hookrightarrow \quad (1\ 1\ 1\ 1\ 1\ 1\ 0) = \quad X_3$$

we see $X_3A = X_3$. Thus $X_i$, on passing through the dynamical system gives a fixed point so mathematically without any doubt we can say easy money based on profession (ie coolie, construction labourer, truck driver etc. easy money related to the individual who do not care for the family) leads to wrongful company, socially irresponsible and free, visiting of CSWs, male ego, more leisure. Only the node awareness of the disease remains in the off state, that is resulting in a fixed point.

Now let us take the node 'macho behaviour' to be in the on state i.e. $Y = (0\ 0\ 0\ 0\ 1\ 0\ 0)$ i.e. all other states are in the off state, passing Y in the connection matrix $YA \hookrightarrow (0\ 1\ 1\ 0\ 1\ 0\ 0)$ of course after updating and thresholding the resultant vector YA, let

$$YA \quad \hookrightarrow \quad (0\ 1\ 1\ 0\ 1\ 0\ 0) \quad = \quad Y_1.$$

$$Y_1A \quad \hookrightarrow \quad (0\ 1\ 1\ 0\ 1\ 1\ 0) \quad = \quad Y_2.$$

$$Y_2A \quad \hookrightarrow \quad (0\ 1\ 1\ 1\ 1\ 1\ 0) \quad = \quad Y_3.$$

$Y_3$ which results in a fixed point. Hence we see that when a person displays a 'macho' and male-chauvinist behaviour he would opt for bad company and bad habits and would be socially irresponsible, visit the CSWs thus all the nodes are on except $A'_7$ and $A'_7$: 'Easy money' and 'No awareness about the disease.'

Now suppose we take the node $A'_7$ i.e. no awareness about the disease IS in the on state we will now find the hidden pattern with only $A'_7$ in the on state, and all other nodes to be in the off state. Let $Z = (0\ 0\ 0\ 0\ 0\ 0\ 1)$, passing Z into the connection matrix A we get $ZA \hookrightarrow (0\ 0\ 1\ 1\ 0\ 1\ 1)$. Let $Z_1 = (0\ 0\ 1\ 1\ 0\ 1\ 1)$, now

$$Z_1A \quad \hookrightarrow \quad (0\ 1\ 1\ 1\ 0\ 1\ 1) \quad = \quad Z_2 \text{ say}$$

$$Z_2A \quad \hookrightarrow \quad (0\ 1\ 1\ 1\ 1\ 1\ 1) \quad = \quad Z_3 \text{ say}$$

$$Z_3A \quad \hookrightarrow \quad (0\ 1\ 1\ 1\ 1\ 1\ 1) \quad = \quad Z_3.$$

Thus we see if the person is not aware of the disease, it implies that the person is generally a victim of wrongful company, that he is socially irresponsible, that he visits CSWs, has leisure and is full of macho behaviour and male chauvinism. However, only the concept of 'easy money and profession' has no relation with this. Thus we can in all cases obtain an hidden pattern which is never possible by any other model. Thus this model is well-suited to give the impact of each attribute in an HIV/AIDS patient and the effect or inter-relation between these attributes which no other mathematical model has given.

The third expert is a HIV/AIDS patient himself who is an inpatient but he says he has become "reformed" and "refined" after the disease that is, "he is willing to do service and spread awareness about HIV/AIDS to people as a social service." However he was not willing to give his name or identity in this book. He was very sincere in giving his opinion about patients as he has spent over eight years as an in and out patient of the hospital and has done social work for many other very invalid patients but he is not a member of any NGO.

As the attributes given by him was different from the other two experts justification of them in certain cases differ which we have listed.

The third experts opinion with conceptual nodes $P_1$, ..., $P_{15}$ and the explanation or the definition of the nodes as given by him and the inter relation as mapped by him are given.

$P_1$    -    No binding with family
$P_2$    -    Male chauvinism
$P_3$    -    Women as inferior objects
$P_4$    -    Bad company and bad habits
$P_5$    -    Socially irresponsible
$P_6$    -    Uncontrollable sexual feelings
$P_7$    -    Food habits: gluttony
$P_8$    -    More leisure
$P_9$    -    No other work for the brain
$P_{10}$    -    Only physically active
$P_{11}$    -    Visits Commercial Sex Workers (CSWs)
$P_{12}$    -    Enjoys life -- jolly mood
$P_{13}$    -    Unreachable by friends or relatives
$P_{14}$    -    Pride in visiting "countless" CSWs
$P_{15}$    -    Failure of agriculture.

He had taken fifteen attributes as fuzzy nodes of any HIV/AIDS patient i.e. a man who has acquired the disease. He gave what he precisely means by each of the attributes.

Further a NGO or a doctor expert cannot be so good an expert as a HIV/AIDS patients who is not only living with the disease for the past eight years but who is also very closely acquainted with only HIV/AIDS patients for the past 8 years. He says not only he has shared his views with HIV/AIDS patients but also spent several hours talking about it with the kith and kin of many of the patients both male and female and he was very happy and excited while giving his opinion and discussing with us.

He is a 7th standard failed man in his late thirties hailing from a village, and an agriculturist labourer by profession. Since at the first instance we were not very convinced of the nodes/attributes he has given we asked the reason for taking those as attributes and we in order to make the reader know the reason for taking those nodes we have given the explanations as given by him. Also as this is the unsupervised data and once we seek the opinion of any expert we have no right to delete or add any new nodes, we had strictly followed and recorded every thing without any deletion or addition. Further we felt an opinion by an HIV/AIDS patient may be much more closer to the solution than by any other expert.

Now we proceed on to give the explanation of each of the nodes $P_1$, $P_2$, …, $P_{15}$ as given by that expert:

$P_1$ – Most of the HIV/AIDS patients do not have binding relationship with family that is why they are often tempted to visit CSW, drink, smoke and waste money on bad company.

$P_2$ – Male chauvinism and macho behaviour.

$P_3$ – Most of the men with male ego and who have no binding in family often think of women as an inferior object and an object of only pleasure / sex.

$P_4$ – Bad company and bad habits.

$P_5$ – Socially irresponsible.

$P_6$ – Several of these persons go to CSW due to uncontrollable emotions that is why some of HIV/AIDS patients had acknowledged to him if ever he had self control he would not have erred. That is why most of them do not use safe sex methods even if the CSW advises them. Further even in the interview some men wanted medicines to control their urge to have sex.

$P_7$ – Gluttony. This man says that ravenous appetite for food leads to the same kind of urge for sex.

$P_8$ – More leisure.

$P_9$ & $P_{10}$ – According to him that when men/women after a age do not have a goal or an ideal in life then he terms them as person for whom only the body functions and they have no work for the brain.

$P_{11}$ – CSW.

$P_{12}$ –  Most when asked why they visited CSWs the answer will be to 'enjoy life' in a "jolly mood"/"pastime" and so on.

$P_{13}$ –  When mainly they are in a different state they are more free from fear of being noticed so without any sense of shame / sin they visit CSW.

$P_{14}$ –  Most of them reply that they have visited countless CSWs, this is their attitude.

$P_{15}$ –  This point he says with his personal experience that most of the agricultural labour have taken up the profession of truck or lorry drivers for livelihood. Thus one of the major reasons is failure of agriculture/poor yield or not up to anticipation. He says in his village alone over 70 families have migrated some to Bombay and some to Madras. Many men have taken the profession of a drivers.

With these as nodes we give the directed graph of this expert.

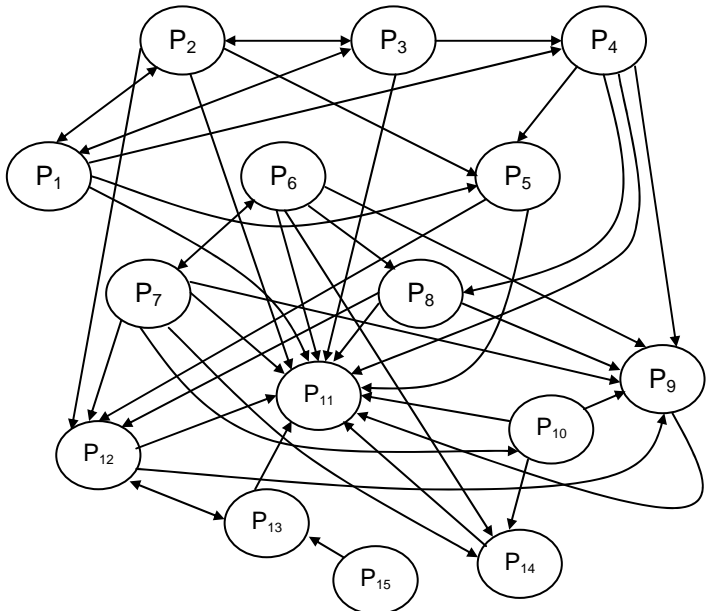

FIGURE 2.2.5

P is the related matrix of the directed graph given in figure 2.2.5.

$$\begin{bmatrix} 0 & 1 & 1 & 1 & 1 & 0 & 0 & 0 & 0 & 0 & 1 & 0 & 0 & 0 & 0 \\ 1 & 0 & 1 & 0 & 1 & 0 & 0 & 0 & 0 & 0 & 1 & 1 & 0 & 0 & 0 \\ 1 & 1 & 0 & 1 & 0 & 0 & 0 & 0 & 0 & 0 & 1 & 0 & 0 & 0 & 0 \\ 0 & 0 & 0 & 0 & 1 & 0 & 0 & 1 & 1 & 0 & 1 & 0 & 0 & 0 & 0 \\ 0 & 0 & 0 & 0 & 0 & 0 & 0 & 0 & 0 & 0 & 1 & 1 & 0 & 0 & 0 \\ 0 & 0 & 0 & 0 & 0 & 0 & 1 & 1 & 1 & 0 & 1 & 0 & 0 & 1 & 0 \\ 0 & 0 & 0 & 0 & 0 & 1 & 0 & 0 & 1 & 1 & 1 & 1 & 0 & 1 & 0 \\ 0 & 0 & 0 & 0 & 0 & 0 & 0 & 0 & 1 & 0 & 1 & 1 & 0 & 0 & 0 \\ 0 & 0 & 0 & 0 & 0 & 0 & 0 & 0 & 0 & 0 & 1 & 0 & 0 & 0 & 0 \\ 0 & 0 & 0 & 0 & 0 & 0 & 0 & 0 & 1 & 0 & 1 & 0 & 0 & 1 & 0 \\ 0 & 0 & 0 & 0 & 0 & 0 & 0 & 0 & 0 & 0 & 0 & 0 & 0 & 0 & 0 \\ 0 & 0 & 0 & 0 & 0 & 0 & 0 & 0 & 1 & 0 & 1 & 0 & 1 & 0 & 0 \\ 0 & 0 & 0 & 0 & 0 & 0 & 0 & 0 & 0 & 0 & 1 & 1 & 0 & 0 & 0 \\ 0 & 0 & 0 & 0 & 0 & 0 & 0 & 0 & 0 & 0 & 1 & 0 & 0 & 0 & 0 \\ 0 & 0 & 0 & 0 & 0 & 0 & 0 & 0 & 0 & 0 & 0 & 0 & 1 & 0 & 0 \end{bmatrix}$$

Let $T_1 = (0\ 0\ 0\ 0\ 0\ 0\ 1\ 0\ 0\ 0\ 0\ 0\ 0\ 0\ 0)$ be the state vector in which only $P_7$ in the on the state, food habits – glutton. The effect of $T_1 = (0\ 0\ 0\ 0\ 0\ 0\ 1\ 0\ 0\ 0\ 0\ 0\ 0\ 0\ 0)$ on the dynamical system P;

$T_1P \hookrightarrow \quad (0\ 0\ 0\ 0\ 0\ 1\ 1\ 0\ 1\ 1\ 1\ 1\ 0\ 1\ 0) \quad = T_2$

$T_2P \hookrightarrow \quad (0\ 0\ 0\ 0\ 0\ 1\ 1\ 1\ 1\ 1\ 1\ 1\ 1\ 1\ 0) \quad = T_3$

$T_3P \hookrightarrow \quad (0\ 0\ 0\ 0\ 0\ 1\ 1\ 1\ 1\ 1\ 1\ 1\ 1\ 1\ 0) \quad = T_4.$

Thus if one is a glutton, he gets all the other nodes to be 'on', except $P_1$, $P_2$, $P_3$, $P_4$, $P_5$ and $P_{15}$. Since the working becomes very labourious, we have provided a pseudo code for the C Program (Appendix 4) The main advantage of using the above pseudo code program is that it can calculate the hidden pattern viz. the fixed point and the limit cycle for any number of concepts.

Also it works for all 2n - 2 different combinations of the concepts involved in the analysis where just mere observation helps in the prediction of the system behaviour.

Now we use this program and give the fixed point or the limit cycle of the given state vector. Suppose we have S = (0 1 0 0 0 0 0 0 0 0 0 0 0 0) i.e. only male ego is in the on state and all other nodes are in the off state. To find the effect of S on the system P consider

$$SP \hookrightarrow (1\ 1\ 1\ 0\ 1\ 0\ 0\ 0\ 0\ 0\ 1\ 1\ 0\ 0\ 0) = S_1$$

$$S_1P \hookrightarrow (1\ 1\ 1\ 1\ 1\ 0\ 0\ 0\ 0\ 0\ 1\ 1\ 0\ 0\ 0) = S_2$$

$$S_2P \hookrightarrow (1\ 1\ 1\ 1\ 1\ 0\ 0\ 0\ 1\ 0\ 1\ 1\ 1\ 0\ 0) = S_3$$

$$S_3P \hookrightarrow (1\ 1\ 1\ 1\ 1\ 0\ 0\ 1\ 1\ 0\ 1\ 1\ 1\ 0\ 0) = S_3.$$

Thus we see if a person has male chauvinism, all the concepts $P_1$, $P_2$, $P_3$, $P_4$, $P_5$, $P_8$, $P_9$, $P_{11}$, $P_{12}$ and $P_{13}$ comes to the on state. Likewise any vector can be used and conclusions can be drawn. We give here the conclusions drawn from several experts which has been calculated. The three expert's opinion is given only for the reader to follow the procedure.

Thus having given the C programming one can use any number of attributes as per the expert's need and give the conclusions. This study reveals that -

(1) No social responsibility and socially free attitude leads to a further increase of the macho behaviour and male ego forcing the migrants to take up wrongful company and habits which in turn leads to treating women as inferior, hence no binding in the family. So easily, without any sense of shame they visit the CSW, consequently contracting the disease. This does not stop here for they consequently infect their wives and (indirectly their future) children thus the number of HIV/AIDS cases will raise at an alarming number.

Not only near relatives are infected but these men also are the root cause of infecting the CSW and this process will certainly lead to an exponential increase in HIV/AIDS patients. Another important factor to be observed is that migration of these persons has forcefully occurred due to poor agricultural yield or failure of agriculture which is directly attributed to the globalization and modernization. For most of these labourers have easily taken up the post of drivers for they have no education or do not have the knowledge of any other trade. Thus it goes without saying that one of the causes for increase in HIV/AIDS infected among migrants is due to their type of job and the full freedom they enjoy when they are away from the family. It is above all the irresponsibility of the government to provide these agricultural coolis with alternative job opportunities but on the other by importing harvesting machines etc., and using it had still made these agriculture labourer a victim of migration. This also shows if they have concern or any form of binding with the family members under no natural circumstances they will stoop to betray them. For the hard money they earn would be saved and spent on the family members. The general conclusions are given in Chapter 7 of this book.

## 2.3 COMBINED FCM TO STUDY THE HIV/AIDS AFFECTED MIGRANTS LABOURERS SOCIO-ECONOMIC PROBLEM

The concept of combine fuzzy cognitive maps (CFCMs) was just defined in section 2.2. For more about CFCM please refer Kosko [58, 62]. We analyze the problems of HIV/AIDS affected migrant labourer using CFCM.

Now we analyze the same problem using combined Fuzzy cognitive Maps (CFCM). Now we seek experts opinion about the seven attributes $A'_1$, $A'_2$, …, $A'_7$ given in page 54 of chapter II of this book.

The directed graph given by him is

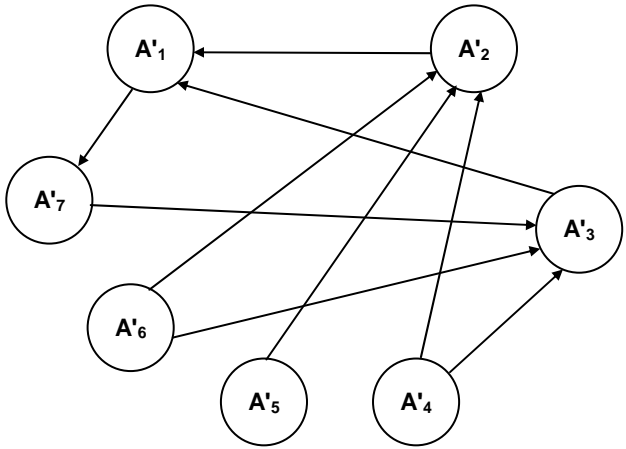

FIGURE 2.3.1

The related fuzzy relational matrix B is as follows:

$$\begin{array}{c} \begin{array}{ccccccc} A_1' & A_2' & A_3' & A_4' & A_5' & A_6' & A_7' \end{array} \\ \begin{array}{c} A_1' \\ A_2' \\ A_3' \\ A_4' \\ A_5' \\ A_6' \\ A_7' \end{array} \begin{bmatrix} 0 & 0 & 0 & 0 & 0 & 0 & 1 \\ 1 & 0 & 0 & 0 & 0 & 0 & 0 \\ 1 & 0 & 0 & 0 & 0 & 0 & 0 \\ 0 & 1 & 1 & 0 & 0 & 0 & 0 \\ 0 & 1 & 0 & 0 & 0 & 0 & 0 \\ 0 & 1 & 1 & 0 & 0 & 0 & 0 \\ 0 & 0 & 1 & 0 & 0 & 0 & 0 \end{bmatrix} \end{array}$$

Suppose we consider the state vector X = (1 0 0 0 0 0 0) i.e., easy availability of money to be in the on state and all other nodes are in the off state. The effect of X on the dynamical system B is

$$XB \quad \hookrightarrow \quad (1\ 0\ 0\ 0\ 0\ 0\ 1) \qquad = \quad X_1 \text{ (say)}$$

$$X_1B \quad \hookrightarrow \quad (1\ 0\ 1\ 0\ 0\ 0\ 1) \qquad = \quad X_2 \text{ (say)}$$

$$X_2B \quad \hookrightarrow \quad (1\ 0\ 1\ 0\ 0\ 0\ 1) \qquad = \quad X_3 \ = X_2$$

($X_2$ which a fixed point of the dynamical system.) i.e., easy money forces one to visit CSWs and he is also unaware of the disease and how it is communicated so he not only visits the CSWs but has unprotected sex with them there by infecting himself. As our motivation is to study the combined effect of the system we now for the same set of seven attributes seek the opinion of another expert. The directed graph given by the second expert is as follows:

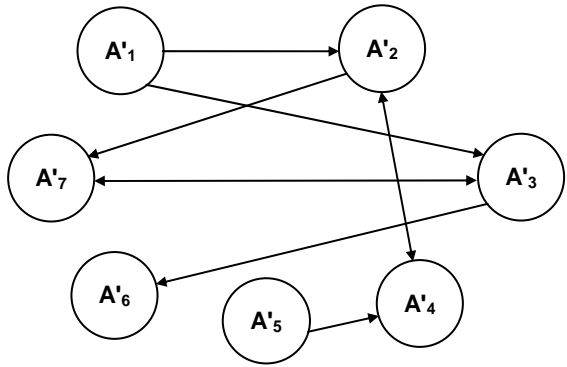

FIGURE 2.3.2

According to this expert when they have wrongful company and additive habits they will certainly become blind to awareness; for he feels as most of the migrant labourer who are HIV/AIDS patients acknowledge with full sense that they were in a complete drunken state when they visited CSWs and their statement they were not able to recollect whether they used *nirodh* (condom) is sufficient to explain when have additive habits they fail to be aware of and to use protected sexual methods. Now we proceed on to give the related connection matrix C

$$
\begin{array}{c c}
 & \begin{array}{ccccccc} A_1' & A_2' & A_3' & A_4' & A_5' & A_6' & A_7' \end{array} \\
C = \begin{array}{c} A_1' \\ A_2' \\ A_3' \\ A_4' \\ A_5' \\ A_6' \\ A_7' \end{array} & \begin{bmatrix} 0 & 1 & 1 & 0 & 0 & 0 & 0 \\ 0 & 0 & 0 & 1 & 0 & 0 & 1 \\ 0 & 0 & 0 & 0 & 0 & 1 & 1 \\ 0 & 1 & 0 & 0 & 0 & 0 & 0 \\ 0 & 0 & 0 & 1 & 0 & 0 & 0 \\ 0 & 0 & 0 & 0 & 0 & 0 & 0 \\ 0 & 0 & 1 & 0 & 0 & 0 & 0 \end{bmatrix}
\end{array}
$$

Suppose for this dynamical system C we consider the on state of the node visits CSWs to be in the on state and all other nodes are in the off state that is Y = (0 0 1 0 0 0 0); to analyze the effect of Y on C

$$YC \quad \hookrightarrow \quad (0\,0\,1\,0\,0\,1\,1) \quad = \quad Y_1 \text{ (say)}$$

$$Y_1C \quad \hookrightarrow \quad (0\,0\,1\,0\,0\,1\,1) \quad = \quad Y_2 \;=\; Y_1 \text{ (say)}$$

The hidden pattern is a fixed point which shows these migrant labourer who visit CSWs have more leisure and are not aware of how HIV/AIDS spreads.

Suppose we consider the on state of the vector easy availability of money to be in the on state the effect of this on the dynamical system C. Let P = (1 0 0 0 0 0 0) The effect of P on C is given by

$$PC \quad \hookrightarrow \quad (1\,1\,1\,0\,0\,0\,0) \quad = \quad P_1 \text{ (say)}$$

$$P_1C \quad \hookrightarrow \quad (1\,1\,1\,1\,0\,1\,1) \quad = \quad P_2 \text{ (say)}$$

$$P_2C \quad \hookrightarrow \quad (1\,1\,1\,1\,0\,1\,1) \quad = \quad P_3 \text{ (a fixed point)}$$

Thus if they earn easy victims of wrongful company and addictive habits, viots CSWs they are socially irresponsible, with more leisure and not knowing about the how the disease spreads.

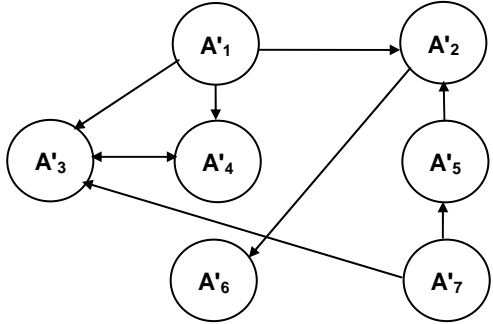

FIGURE 2.3.3

Now we proceed on to take the third experts opinion. The directed graph given by him and the related matrix D is as follows:

$$D = \begin{array}{c} \\ A_1^{'} \\ A_2^{'} \\ A_3^{'} \\ A_4^{'} \\ A_5^{'} \\ A_6^{'} \\ A_7^{'} \end{array} \begin{array}{c} A_1^{'} \ A_2^{'} \ A_3^{'} \ A_4^{'} \ A_5^{'} \ A_6^{'} \ A_7^{'} \\ \begin{bmatrix} 0 & 1 & 1 & 1 & 0 & 0 & 0 \\ 0 & 0 & 0 & 0 & 0 & 1 & 0 \\ 0 & 0 & 0 & 1 & 0 & 0 & 0 \\ 0 & 0 & 1 & 0 & 0 & 0 & 0 \\ 0 & 1 & 0 & 0 & 0 & 0 & 0 \\ 0 & 0 & 0 & 0 & 0 & 0 & 0 \\ 0 & 0 & 1 & 0 & 1 & 0 & 0 \end{bmatrix} \end{array}$$

Suppose we consider the on state of the node socially irresponsible /free and all other nodes are in the off state.

The effect of $Z = (0\ 0\ 0\ 1\ 0\ 0\ 0)$ on the dynamical system D is given by

$$ZD \quad \hookrightarrow \quad (0\ 0\ 1\ 1\ 0\ 0\ 0) \quad = \quad Z_1 \text{ say}$$

$$Z_1 D \quad \hookrightarrow \quad (0\ 0\ 1\ 1\ 0\ 0\ 0) \quad = \quad Z_2 \ = \ Z_1$$

Thus the hidden pattern is a fixed point, leading to the conclusion socially irresponsible persons visit CSWs for they do not have binding on the family or even the responsibility that they would face health problem. Thus a socially free / irresponsible person is more prone to visit CSWs.

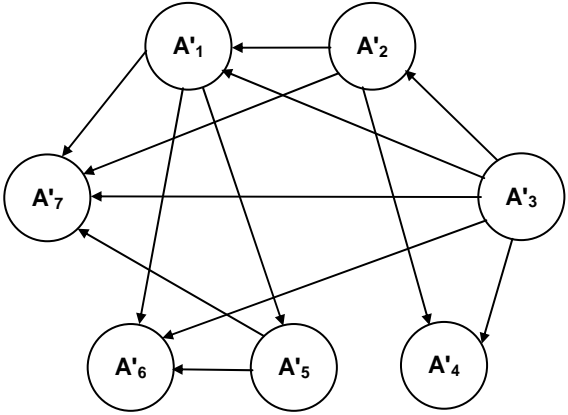

FIGURE 2.3.4

Now we take yet another experts opinion so that we can make conclusions on the combined FCMs. The directed graph given by the expert for figure 2.3.4.

The related matrix E is given below

$$E = \begin{array}{c} \\ A_1^{'} \\ A_2^{'} \\ A_3^{'} \\ A_4^{'} \\ A_5^{'} \\ A_6^{'} \\ A_7^{'} \end{array} \begin{array}{c} A_1^{'} \ A_2^{'} \ A_3^{'} \ A_4^{'} \ A_5^{'} \ A_6^{'} \ A_7^{'} \\ \begin{bmatrix} 0 & 0 & 0 & 0 & 1 & 1 & 1 \\ 1 & 0 & 0 & 1 & 0 & 0 & 0 \\ 1 & 1 & 0 & 1 & 0 & 1 & 1 \\ 0 & 0 & 0 & 0 & 0 & 0 & 0 \\ 0 & 0 & 0 & 0 & 0 & 1 & 1 \\ 0 & 0 & 0 & 0 & 0 & 0 & 0 \\ 0 & 1 & 0 & 0 & 0 & 0 & 0 \end{bmatrix} \end{array}$$

Suppose the node $A_2^{'}$ i.e., wrongful company and addictive habits are in the on state. The effect of U = (0 1 0 0 0 0 0) on the dynamical system E is given by

$$UE \quad \hookmapsto \quad (1\ 1\ 0\ 1\ 0\ 0\ 0) \quad = \quad U_1 \text{ (say)}$$

$$U_1 E \quad \hookmapsto \quad (1\ 1\ 0\ 1\ 1\ 1\ 1) \quad = \quad U_2 \text{ (say)}$$

$$U_2 E \quad \hookmapsto \quad (1\ 1\ 0\ 1\ 1\ 1\ 1) \quad = \quad U_3 = U_2.$$

The hidden pattern is a fixed point. Thus according to this expert when a person is in a wrongful company and has additive habits it need not always imply that he will visit CSWs but for he has easy money, he is socially free /responsible, may have macho behaviour, he has more leisure and has no awareness of about the disease; as our data is an unsupervised one we have no means to change or modify the opinion of any expert we have to give his opinion as it is for otherwise the data would become biased.

Now we formulate the Combined Fuzzy Cognitive Maps using the opinion of the 5 experts including the expert opinions given in chapter two, page 54.

Let S denote the combined connection matrix by $S = A + B + C + D + E$.

$$S = \begin{array}{c} \\ A_1^{'} \\ A_2^{'} \\ A_3^{'} \\ A_4^{'} \\ A_5^{'} \\ A_6^{'} \\ A_7^{'} \end{array} \begin{array}{c} \begin{array}{ccccccc} A_1^{'} & A_2^{'} & A_3^{'} & A_4^{'} & A_5^{'} & A_6^{'} & A_7^{'} \end{array} \\ \left[ \begin{array}{ccccccc} 0 & 3 & 3 & 3 & 1 & 0 & 2 \\ 2 & 0 & 1 & 1 & 1 & 2 & 1 \\ 2 & 1 & 0 & 2 & 0 & 3 & 2 \\ 0 & 2 & 3 & 0 & 0 & 1 & 0 \\ 0 & 3 & 1 & 1 & 0 & 1 & 2 \\ 0 & 2 & 2 & 1 & 0 & 0 & 0 \\ 0 & 1 & 4 & 1 & 1 & 0 & 0 \end{array} \right] \end{array}$$

We threshold in a different way if a quality addsup to 2 or less than two we put 0 and if it is 3 or greater 3 than we put as 1.

Now we consider the effect of only easy availability of money to be in the on state.

The effect of $P = (1\ 0\ 0\ 0\ 0\ 0\ 0)$ on the combined dynamical system S.

$$PS \quad = \quad (0\ 3\ 3\ 3\ 1\ 0\ 2) \quad \hookmapsto \quad (0\ 1\ 1\ 1\ 0\ 0\ 0) \quad = \quad P_1$$

$$P_1 S \quad = \quad (4\ 3\ 4\ 3\ 1\ 5\ 3) \quad \hookmapsto \quad (1\ 1\ 1\ 1\ 0\ 1\ 1) \quad = \quad P_2$$

$$P_2 S \quad = \quad (4\ 9\ 13\ 8\ 3\ 5\ 5) \quad \hookmapsto \quad (1\ 1\ 1\ 1\ 1\ 1\ 1) \quad = \quad P_3$$

$$P_3 \, S \quad \hookrightarrow \quad (1\,1\,1\,1\,1\,1\,1) \qquad = \quad P_4 \; = \quad P_3.$$

The hidden pattern is a fixed point when they have easy availability of money they have all the states to be on. They visit CSW have bad habits, find more leasure (if no money they will do several types of work to earn) and as they are unaware of the disease they without any fear visit CSWs as they want to show their macho character they smoke and drink and so on. Thus the collective or combined FCM shows the on state of all vectors.

Suppose we consider the on state as the only attribute they have no awareness about the disease say $L = (0\,0\,0\,0\,0\,0\,1)$ and all other attributes are off. We find the effect of L on the dynamical system E.

$$
\begin{aligned}
LE \;&=\; (0\,1\,4\,1\,1\,0\,0) \;\hookrightarrow\; (0\,0\,1\,0\,0\,0\,1) \;=\; L_1 \text{ say} \\
L_1 E \;&=\; (2\,2\,4\,3\,1\,3\,2) \;\hookrightarrow\; (0\,0\,1\,1\,0\,1\,1) \;=\; L_2 \text{ say} \\
L_2 E \;&=\; (2\,6\,9\,4\,1\,4\,2) \;\hookrightarrow\; (0\,1\,1\,1\,0\,1\,1) \;=\; L_3 \text{ say} \\
L_3 E \;&=\; (4\,6\,10\,5\,2\,5\,3) \;\hookrightarrow\; (1\,1\,1\,1\,0\,1\,1) \;=\; L_4 \text{ say} \\
L_4 E \;&\hookrightarrow\; (1\,1\,1\,1\,1\,1\,1),
\end{aligned}
$$

which is a fixed point of the hidden pattern. We can study and interpret the effect of each and every state vector using the C program given in appendix 5.

Also we have only illustrated using 7 attibutes only with 5 experts, one can choose any number of attributes and use any desired number of experts and obtain the results using the C-program.

## 2.4. COMBINED DISJOINT BLOCK FCM AND ITS APPLICATION TO HIV/AIDS PROBLEM

In this section we define for the first time the notion of combined Disjoint Block FCM (CDB FCM) and apply them in the analysis of the socio economic problems of the HIV/AIDS affected migrant labourers.

Let $C_1, \ldots , C_n$ be n nodes / attributes related with some problem. The n may be very large and a non prime. Even though we have C- program to work finding the directed graph and the related connection matrix may be very unwieldy. In such cases we use the notion of combined disjoint block fuzzy cognitive maps. We divide these n attributes into k equal classes and these k equal classes are viewed by k-experts or by the same expert and the corresponding directed graph and the connection matrices are got.

Now these connection matrices is made into a n × n matrix and using the C-program the results are derived. This type of Combined Disjoint Block FCM is know as Combined Disjoint Equal Block FCM. Now some times we may not be in a position to divide the 'n' under study into equal blocks in such cases we use the technique of dividing the n attributes say of some m blocks were each block may not have the same number of attributes, but it is essential that there n attributes are divided into disjoint classes. We use both these techniques in the analysis of the problem. This is the case when n happens to be a prime number. We in the next stage proceed on to define the notion of combined overlap block Fuzzy cognitive maps, for this also we assume two types the overlap is regular and the overlap happens to be irregular. Both the cases are dealt in this book.

**DEFINITION 2.4.1:** *Let $C_1$ … $C_n$ be n distinct attributes of a problem n very large and a non prime. If we divide n into k equal classes i.e., k/n and if n/k = t which are disjoint and if we find the directed graph of each of there k classes of attributes with t attributes each, then their corresponding connection matrices are formed and these connection matrices are joined as blocks to form a n × n matrix.*

*This n × n connection matrix forms the combined disjoint block FCM of equal classes. If the classes are not divided to have equal attributers but if they are disjoint classes we get a n × n connection matrix called the combined disjoint block FCM of unequal classes/ size.*

We illustrate first our models by a simple example. The attributes taken here are related with the HIV/AIDS affected Migrant Labourers. Suppose the 12 attributes $A_1$, $A_2$, $A_3$, $A_4$,…$A_{12}$ given in pages 42-50 is taken i.e.,

$A_1$ - Easy availability of money
$A_2$ - Lack of Education
$A_3$ - Visiting CSWs
$A_4$ - Nature of Profession
$A_5$ - Wrong / Bad company
$A_6$ - Addiction to habit forming substances and visiting CSWs.
$A_7$ - Absence of social responsibility

$A_8$ - Socially free
$A_9$ - Economic Status
$A_{10}$ - More leisure
$A_{11}$ - Machismo / Exaggerated Masculinity
$A_{12}$ - No awareness of the disease.

Using these attributes we give the combined disjoint block fuzzy cognitive map of equal size. We take 3 experts opinion on the 3 disjoint classes so that each class has four attributes / nodes. Let the disjoint classes be $C_1$, $C_2$ and $C_3$ be divided by the following:

$C_1$ = $\{A_1, A_6, A_7, A_{12}\}$, $C_2$ = $\{A_2, A_3, A_4, A_{10}\}$ and $C_3$ = $\{A_5, A_8, A_9, A_{11}\}$.

Now we collect the experts opinion on each of the classes $C_1$, $C_2$ and $C_3$. The directed graph given by the expect on attributes $A_1$, $A_6$, $A_7$ and $A_{12}$ which forms the class $C_1$.

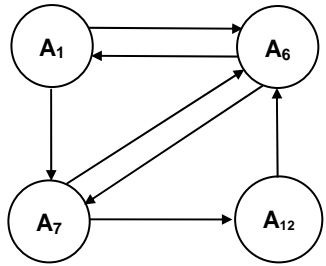

FIGURE 2.4.1

The related connection matrix $B_1$

$$B_1 = \begin{array}{c} \\ A_1 \\ A_6 \\ A_7 \\ A_{12} \end{array} \begin{array}{c} \begin{array}{cccc} A_1 & A_6 & A_7 & A_{12} \end{array} \\ \begin{bmatrix} 0 & 1 & 1 & 0 \\ 1 & 0 & 1 & 0 \\ 0 & 1 & 0 & 1 \\ 0 & 1 & 0 & 0 \end{bmatrix} \end{array}$$

The directed graph given by the expert on $A_2$, $A_3$, $A_4$, $A_{10}$ which forms the class $C_2$.

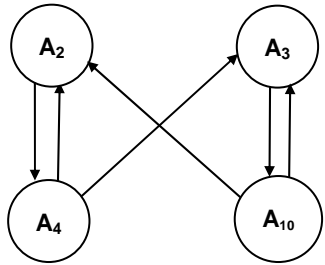

FIGURE 2.4.2

According to this expert the two nodes, lack of education and the nature of profession are very much interrelated. Also he feels most of the educated people have awareness about HIV/AIDS. The related connection matrix $B_2$ is given below:

$$B_2 = \begin{array}{c} \\ A_2 \\ A_3 \\ A_4 \\ A_{10} \end{array} \begin{array}{cccc} A_2 & A_3 & A_4 & A_{10} \\ \begin{bmatrix} 0 & 0 & 1 & 0 \\ 0 & 0 & 0 & 1 \\ 1 & 1 & 0 & 0 \\ 1 & 1 & 0 & 0 \end{bmatrix} \end{array}$$

Now we give the directed graph for the class $C_3$ as given by the expert $C_3 = \{A_5\ A_8\ A_9\ A_{11}\}$.

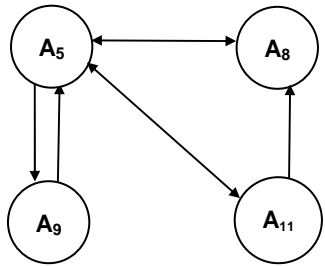

FIGURE 2.4.3

According to this expert bad company induces Macho behaviour. Also when a person has bad habit / bad company he is not bounded socially he without any fear engages himself in doing social evil, so he relates bad company with socially free. If he has no food for the next day, i.e., with this economic status certainly he cannot enjoy or have bad company or bad habits. Usually person with bad company show more of macho behaviour. The related connection matrix

$$B_3 = \begin{array}{c} \\ A_5 \\ A_8 \\ A_9 \\ A_{11} \end{array} \begin{array}{c} \begin{array}{cccc} A_5 & A_8 & A_9 & A_{11} \end{array} \\ \begin{bmatrix} 0 & 1 & 1 & 1 \\ 1 & 0 & 0 & 0 \\ 1 & 0 & 0 & 0 \\ 0 & 1 & 0 & 0 \end{bmatrix} \end{array}$$

Now the combined disjoint block connection matrix of the FCM is given by B

$$\begin{array}{c} \\ A_1 \\ A_6 \\ A_7 \\ A_{12} \\ A_2 \\ A_3 \\ A_4 \\ A_{10} \\ A_5 \\ A_8 \\ A_9 \\ A_{11} \end{array} \begin{array}{c} \begin{array}{cccccccccccc} A_1 & A_6 & A_7 & A_{12} & A_2 & A_3 & A_4 & A_{10} & A_5 & A_8 & A_9 & A_{11} \end{array} \\ \begin{bmatrix} 0 & 1 & 1 & 0 & 0 & 0 & 0 & 0 & 0 & 0 & 0 & 0 \\ 1 & 0 & 1 & 0 & 0 & 0 & 0 & 0 & 0 & 0 & 0 & 0 \\ 0 & 1 & 0 & 1 & 0 & 0 & 0 & 0 & 0 & 0 & 0 & 0 \\ 0 & 1 & 0 & 0 & 0 & 0 & 0 & 0 & 0 & 0 & 0 & 0 \\ 0 & 0 & 0 & 0 & 0 & 0 & 1 & 0 & 0 & 0 & 0 & 0 \\ 0 & 0 & 0 & 0 & 0 & 0 & 0 & 1 & 0 & 0 & 0 & 0 \\ 0 & 0 & 0 & 0 & 1 & 1 & 0 & 0 & 0 & 0 & 0 & 0 \\ 0 & 0 & 0 & 0 & 1 & 1 & 0 & 0 & 0 & 0 & 0 & 0 \\ 0 & 0 & 0 & 0 & 0 & 0 & 0 & 0 & 0 & 1 & 1 & 1 \\ 0 & 0 & 0 & 0 & 0 & 0 & 0 & 0 & 1 & 0 & 0 & 0 \\ 0 & 0 & 0 & 0 & 0 & 0 & 0 & 0 & 1 & 0 & 0 & 0 \\ 0 & 0 & 0 & 0 & 0 & 0 & 0 & 0 & 0 & 1 & 0 & 0 \end{bmatrix} \end{array}$$

Suppose we consider the on state of the attribute easy availability of the money and all other states are off the effect of X = (1 0 0 0 0 0 0 0 0 0 0 0) on the CDB FCM is given by

XB $\hookrightarrow$ (1 1 1 0 0 0 0 0 0 0 0 0) = $X_1$ (Say)

$X_1 B$ $\hookrightarrow$ (1 1 1 1 0 0 0 0 0 0 0 0) = $X_2$ (Say)

$X_2 B$ $\hookrightarrow$ (1 1 1 1 0 0 0 0 0 0 0 0) = $X_3 = X_2$.

($X_2$ is a fixed point of the dynamical system). Thus when one has easy money he has addiction to habit forming substances and visits CSWs he has no social responsibility and he is not aware of how the disease spreads.

Suppose we consider, the on state of the attributes visits CSWs and socially free to be in the on state and all other nodes are in the off state. Now we study the effect on the dynamical system B. Let T = (0 0 0 0 0 1 0 0 0 1 0 0) state vector depicting the on state of CSWs and socially free, passing the state vector T in to the dynamical system B.

$$TB \quad \hookmapsto \quad (0\ 0\ 0\ 0\ 0\ 1\ 0\ 1\ 1\ 1\ 0\ 0) \quad = \quad T_1 \text{ (say)}$$

$$T_1B \quad \hookmapsto \quad (0\ 0\ 0\ 0\ 1\ 1\ 1\ 1\ 1\ 1\ 1\ 1) \quad = \quad T_2 \text{ (say)}$$

$$T_2B \quad \hookmapsto \quad (0\ 0\ 0\ 0\ 1\ 1\ 1\ 1\ 1\ 1\ 1\ 1) \quad = \quad T_3$$

($T_2$ the fixed point of the dynamical system. Thus when he visits CSWs and is socially free one gets the lack of education, money is easily available to him, Absence of social responsibility and no awareness about the disease is in the off state and all other states become on. One can study this dynamical system using C-program given in appendix.

Now we give another application of the Combined disjoint equal block fuzzy cognitive maps to the HIV/AIDS affected migrant labourers. Now we consider the following 15 conceptual nodes $S_1$, …, $S_{15}$ associated with the HIV/AIDS patients and migrancy

$S_1$ - No Awareness of HIV/AIDS in migrant labourers
$S_2$ - Rural living with no education
$S_3$ - Migrancy as truck drivers / daily wagers
$S_4$ - Socially free and irresponsible
$S_5$ - Enjoy life Jolly mood
$S_6$ - Away from family
$S_7$ - More leisure
$S_8$ - No association or union to protect and channelize their money / time
$S_9$ - Visits CSWs
$S_{10}$ - Unreachable by friends and relatives
$S_{11}$ - Easy victims of temptation
$S_{12}$ - Takes to all bad habit / bad company including CSWs
$S_{13}$ - Male Chauvinism

$S_{14}$ - Women as inferior objects

$S_{15}$ - No family binding and respect for wife.

These 15 attributes are divided in 5 classes $C_1, C_2, \ldots, C_5$ with 3 in each class.

Let       $C_1$    =       $\{S_1, S_2, S_3\}$,
          $C_2$    =       $\{S_4, S_5, S_6\}$,
          $C_3$    =       $\{S_7, S_8, S_9\}$,
          $C_4$    =       $\{S_{10}, S_{11}, S_{12}\}$
and       $C_5$    =       $\{S_{13}, S_{14}, S_{15}\}$.

Now we take the experts opinion for each of these classes and take the matrix associated with the combined block disjoint FCMs.

Using the $15 \times 15$ matrix got using the combined connection matrix we derive our conclusions. The experts opinion for the class $C_1 = \{S_1, S_2, S_3\}$ in the form of the directed graph.

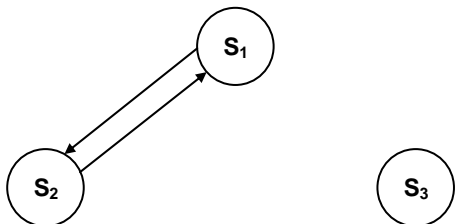

FIGURE 2.4.4

The related connection matrix is given by

$$
\begin{array}{c}
\begin{array}{ccc} S_1 & S_2 & S_3 \end{array} \\
\begin{array}{c} S_1 \\ S_2 \\ S_3 \end{array}
\begin{bmatrix}
0 & 1 & 0 \\
1 & 0 & 0 \\
0 & 0 & 0
\end{bmatrix}
\end{array}
$$

The directed graph for the class $C_2 = \{S_4, S_5, S_6\}$ as given by the expert

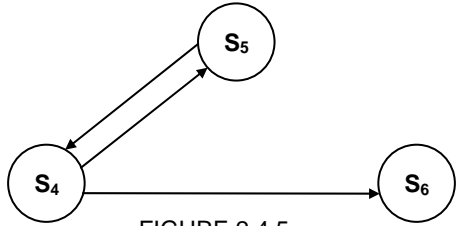

FIGURE 2.4.5

The related connection matrix is given as follows:

$$
\begin{array}{c c c}
 & S_4 & S_5 & S_6 \\
\begin{array}{c} S_4 \\ S_5 \\ S_6 \end{array} &
\left[\begin{array}{c c c}
0 & 1 & 1 \\
1 & 0 & 0 \\
0 & 0 & 0
\end{array}\right]
\end{array}
$$

Now for the class $C_3 = \{S_7, S_8, S_9\}$ the directed graph is given below by the expert

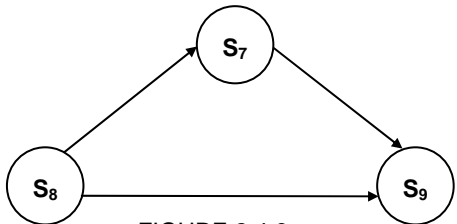

FIGURE 2.4.6

The related connection matrix

$$
\begin{array}{c c c}
 & S_7 & S_8 & S_9 \\
\begin{array}{c} S_7 \\ S_8 \\ S_9 \end{array} &
\left[\begin{array}{c c c}
0 & 0 & 1 \\
1 & 0 & 1 \\
0 & 0 & 0
\end{array}\right]
\end{array}
$$

The directed graph given by the expert for the class $C_4 = \{S_{10}, S_{11}, S_{12}\}$ is given below:

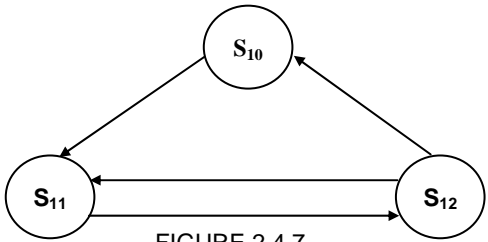

FIGURE 2.4.7

The related connection matrix

$$
\begin{array}{c@{\quad}c}
 & \begin{array}{ccc} S_{10} & S_{11} & S_{12} \end{array} \\
\begin{array}{c} S_{10} \\ S_{11} \\ S_{12} \end{array} &
\left[\begin{array}{ccc}
0 & 1 & 0 \\
0 & 0 & 1 \\
1 & 1 & 0
\end{array}\right]
\end{array}
$$

Now for the class $C_5 = \{S_{13}, S_{14}, S_{15}\}$, the directed graph is as follows:

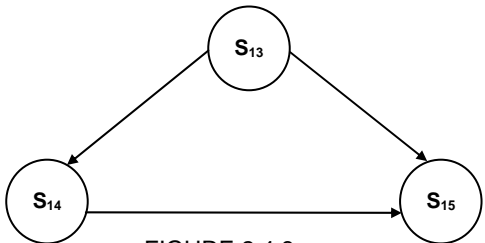

FIGURE 2.4.8

The related connection matrix

$$
\begin{array}{c@{\quad}c}
 & \begin{array}{ccc} S_{13} & S_{14} & S_{15} \end{array} \\
\begin{array}{c} S_{13} \\ S_{14} \\ S_{15} \end{array} &
\left[\begin{array}{ccc}
0 & 1 & 1 \\
0 & 0 & 1 \\
0 & 0 & 0
\end{array}\right]
\end{array}
$$

Now we give the combined block connection matrix related with the 15 attributes, let us denote it by S.

$$\begin{array}{c} \\ S = \end{array} \begin{array}{c} \\ \begin{array}{c}S_1\\S_2\\S_3\\S_4\\S_5\\S_6\\S_7\\S_8\\S_9\\S_{10}\\S_{11}\\S_{12}\\S_{13}\\S_{14}\\S_{15}\end{array} \end{array} \overset{\begin{array}{ccccccccccccccc}S_1 & S_2 & S_3 & S_4 & S_5 & S_6 & S_7 & S_8 & S_9 & S_{10} & S_{11} & S_{12} & S_{13} & S_{14} & S_{15}\end{array}}{\begin{bmatrix} 0&1&0&0&0&0&0&0&0&0&0&0&0&0&0\\ 1&0&0&0&0&0&0&0&0&0&0&0&0&0&0\\ 0&0&0&0&0&0&0&0&0&0&0&0&0&0&0\\ 0&0&0&0&1&1&0&0&0&0&0&0&0&0&0\\ 0&0&0&1&0&0&0&0&0&0&0&0&0&0&0\\ 0&0&0&0&0&0&0&0&0&0&0&0&0&0&0\\ 0&0&0&0&0&0&0&0&1&0&0&0&0&0&0\\ 0&0&0&0&0&0&1&0&1&0&0&0&0&0&0\\ 0&0&0&0&0&0&0&0&0&0&0&0&0&0&0\\ 0&0&0&0&0&0&0&0&0&0&1&0&0&0&0\\ 0&0&0&0&0&0&0&0&0&0&0&1&0&0&0\\ 0&0&0&0&0&0&0&0&0&1&1&0&0&0&0\\ 0&0&0&0&0&0&0&0&0&0&0&0&0&1&1\\ 0&0&0&0&0&0&0&0&0&0&0&0&0&0&1\\ 0&0&0&0&0&0&0&0&0&0&0&0&0&0&0 \end{bmatrix}}$$

Now consider a state vector X = (0 1 0 0 0 1 0 0 0 1 0 0 0 1 0) where the nodes/concepts $S_2$, $S_6$, $S_{10}$ and $S_{14}$ are in the on state. The effect of X on the dynamical system S is given by

$$\begin{array}{lllll} XS & \hookrightarrow & (1\,1\,0\,0\,0\,1\,0\,0\,0\,1\,1\,0\,0\,1\,1) & = & X_1 \text{ say}\\ X_1S & \hookrightarrow & (1\,1\,0\,0\,0\,1\,0\,0\,0\,1\,1\,0\,0\,1\,1) & = & X_2 \text{ say} \end{array}$$

$X_1 = X_2$ is a fixed point of the hidden pattern. We use the C program for CFCM given in the appendix 5 to study the problem and draw conclusions, which is given in the chapter 7 of this book.

## 2.5. COMBINED OVERLAP BLOCK FCM AND ITS USE IN THE ANALYSIS OF HIV/AIDS AFFECTED MIGRANT LABOURERS

Next we give the model of the combined block overlap fuzzy cognitive maps and adapt it mainly on the HIV/AIDS affected migrant labourers relative to their socio economic problems. We adapt it to the model $\{A_1, A_2, \ldots, A_{12}\}$ given in section 2.4 of this book. Let us consider the class

$$
\begin{array}{llll}
C_1 & = & \{A_1\ A_2\ A_3\ A_4\}, & C_2 & = & \{A_3\ A_4\ A_5\ A_6\}, \\
C_3 & = & \{A_5\ A_6\ A_7\ A_8\}, & C_4 & = & \{A_7\ A_8\ A_9\ A_{10}\}, \\
C_5 & = & \{A_9\ A_{10}\ A_{11}\ A_{12}\}\ \text{and} & C_6 & = & \{A_{11}\ A_{12}\ A_1\ A_2\}.
\end{array}
$$

We give the directed graph for each of the classes of attributes $C_1$, $C_2$,…, $C_6$. The directed graph for the four attributes given by $C_1 = \{A_1\ A_2\ A_3\ A_4\}$ as given by the expert is as follows.

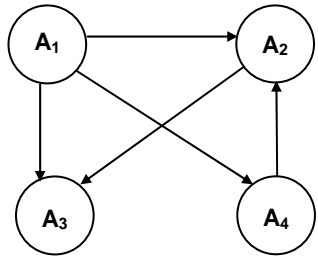

FIGURE 2.5.1

The related connection matrix

$$
M_1 = \begin{array}{c} \\ A_1 \\ A_2 \\ A_3 \\ A_4 \end{array}
\begin{array}{c} \begin{array}{cccc} A_1 & A_2 & A_3 & A_4 \end{array} \\
\begin{bmatrix}
0 & 1 & 1 & 1 \\
0 & 0 & 1 & 0 \\
0 & 0 & 0 & 0 \\
0 & 1 & 0 & 0
\end{bmatrix}
\end{array}
$$

The directed graph for the class of attributes $C_2 = \{A_3\ A_4\ A_5\ A_6\}$ as given by the expert is as follows:

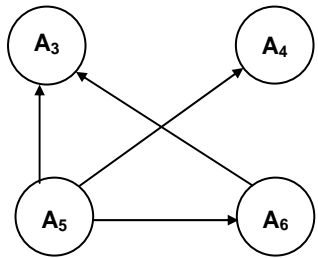

FIGURE 2.5.2

The related connection matrix

$$M_2 = \begin{array}{c} \\ A_3 \\ A_4 \\ A_5 \\ A_6 \end{array} \begin{array}{cccc} A_3 & A_4 & A_5 & A_6 \\ \begin{bmatrix} 0 & 0 & 0 & 0 \\ 0 & 0 & 0 & 0 \\ 1 & 1 & 0 & 1 \\ 1 & 0 & 0 & 0 \end{bmatrix} \end{array}$$

The directed graph associated with the attributes $C_3 = \{A_5\ A_6\ A_7\ A_8\}$ as given by the expert.

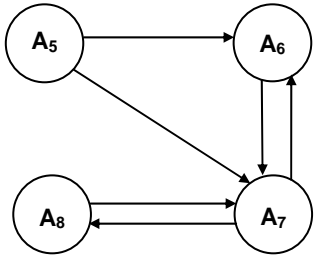

FIGURE 2.5.3

The related connection matrix

$$\begin{array}{c} \\ A_5 \\ A_6 \\ A_7 \\ A_8 \end{array} \begin{array}{cccc} A_5 & A_6 & A_7 & A_8 \\ \begin{bmatrix} 0 & 1 & 1 & 0 \\ 0 & 0 & 1 & 0 \\ 0 & 1 & 0 & 1 \\ 0 & 0 & 1 & 0 \end{bmatrix} \end{array}$$

For the class of attributes $C_4 = \{A_7\ A_8\ A_9\ A_{10}\}$ the directed graph as given by the expert is as follows.

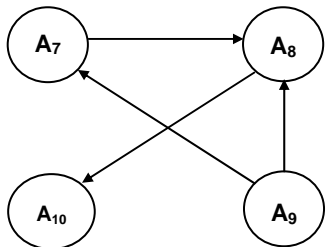

FIGURE 2.5.4

The connection matrix associated with the directed graph.

$$\begin{array}{c} \\ A_7 \\ A_8 \\ A_9 \\ A_{10} \end{array} \begin{array}{c} A_7\ A_8\ A_9\ A_{10} \\ \begin{bmatrix} 0 & 1 & 0 & 0 \\ 0 & 0 & 0 & 1 \\ 1 & 1 & 0 & 0 \\ 0 & 0 & 0 & 0 \end{bmatrix} \end{array}$$

The directed graph associated with the class C = {$A_9$ $A_{10}$ $A_{11}$ $A_{12}$} as given by the expert.

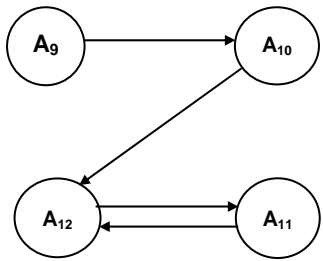

FIGURE 2.5.5

The related connection matrix is

$$\begin{array}{c} \\ A_9 \\ A_{10} \\ A_{11} \\ A_{12} \end{array} \begin{array}{c} A_9\ A_{10}\ A_{11}\ A_{12} \\ \begin{bmatrix} 0 & 1 & 0 & 0 \\ 0 & 0 & 0 & 1 \\ 0 & 0 & 0 & 1 \\ 0 & 0 & 1 & 0 \end{bmatrix} \end{array}$$

Now we give the directed graph of the last class given by the expert

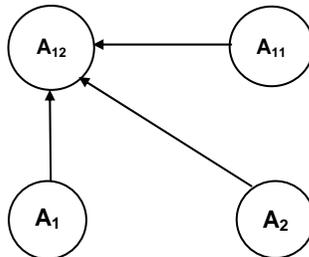

FIGURE 2.5.6

The related connection matrix

$$\begin{array}{c} \\ A_{11} \\ A_{12} \\ A_1 \\ A_2 \end{array} \begin{array}{cccc} A_{11} & A_{12} & A_1 & A_2 \\ \begin{bmatrix} 0 & 1 & 0 & 0 \\ 0 & 0 & 0 & 0 \\ 0 & 1 & 0 & 0 \\ 0 & 1 & 0 & 0 \end{bmatrix} \end{array}$$

The related Combined Block Overlap Fuzzy Cognitive Map (CBOFCM), matrix is as follows. We denote this matrix by W.

$$\begin{array}{c} A_1 \\ A_2 \\ A_3 \\ A_4 \\ A_5 \\ A_6 \\ A_7 \\ A_8 \\ A_9 \\ A_{10} \\ A_{11} \\ A_{12} \end{array} \begin{array}{cccccccccccc} A_1 & A_2 & A_3 & A_4 & A_5 & A_6 & A_7 & A_8 & A_9 & A_{10} & A_{11} & A_{12} \\ \begin{bmatrix} 0 & 1 & 1 & 1 & 0 & 0 & 0 & 0 & 0 & 0 & 0 & 1 \\ 0 & 0 & 1 & 0 & 0 & 0 & 0 & 0 & 0 & 0 & 0 & 1 \\ 0 & 0 & 0 & 0 & 0 & 0 & 0 & 0 & 0 & 0 & 0 & 0 \\ 0 & 1 & 0 & 0 & 0 & 0 & 0 & 0 & 0 & 0 & 0 & 0 \\ 0 & 0 & 1 & 1 & 0 & 2 & 1 & 0 & 0 & 0 & 0 & 0 \\ 0 & 0 & 0 & 1 & 0 & 0 & 1 & 0 & 0 & 0 & 0 & 0 \\ 0 & 0 & 0 & 0 & 0 & 1 & 0 & 2 & 0 & 0 & 0 & 0 \\ 0 & 0 & 0 & 0 & 0 & 0 & 1 & 0 & 0 & 1 & 0 & 0 \\ 0 & 0 & 0 & 0 & 0 & 0 & 0 & 0 & 0 & 1 & 0 & 0 \\ 0 & 0 & 0 & 0 & 0 & 0 & 0 & 1 & 0 & 0 & 0 & 1 \\ 0 & 0 & 0 & 0 & 0 & 0 & 0 & 0 & 0 & 1 & 0 & 2 \\ 0 & 0 & 0 & 0 & 0 & 0 & 0 & 0 & 1 & 0 & 1 & 0 \end{bmatrix} \end{array}$$

Let us consider the state vector X = (1 0 0 0 0 0 0 0 0 0 0 0) where easy availability of money alone is in the on state and all vectors are in the off state. The effect of X on the dynamical system W is give by

$$XW \hookrightarrow (1\ 1\ 1\ 1\ 0\ 0\ 0\ 0\ 0\ 0\ 0\ 1) = X_1 \text{ (say)}$$

$$X_1W \hookrightarrow (1\ 1\ 1\ 1\ 0\ 0\ 0\ 0\ 1\ 0\ 1\ 1) = X_2 \text{ say}$$

$$X_2W \hookrightarrow (1\ 1\ 1\ 1\ 0\ 0\ 0\ 0\ 1\ 1\ 1\ 1) = X_3 \text{ (say)}$$

$$X_3W \hookrightarrow (1\ 1\ 1\ 1\ 0\ 0\ 0\ 1\ 1\ 1\ 1\ 1) = X_4 \text{ (say)}$$

$$X_4W \hookrightarrow (1\ 1\ 1\ 1\ 0\ 0\ 1\ 1\ 1\ 1\ 1\ 1) = X_5 \text{ (say)}$$

$$X_5W \;\mapsto\; (1\ 1\ 1\ 1\ 0\ 1\ 1\ 1\ 1\ 1\ 1\ 1) \;=\; X_6 \text{ (say)}$$

$$X_6W \;\mapsto\; (1\ 1\ 1\ 1\ 0\ 1\ 1\ 1\ 1\ 1\ 1\ 1) \;=\; X_7 \;=\; X_6.$$

The hidden pattern of the CFCM is a fixed point. It says if a person has easy money one acquires all the attributes except $A_5$. Using the C-program one can work for all the solutions and the conclusions given in chapter 7 are based on working of this model.

It is left as an exercise for the reader to compare the dynamical system W which is formed by combined block overlap FCM and the FCM obtained form A given in pages 42 to 52 and 69 to 73.

Now we apply the same model to the attributes already discussed in page 68 to 72 of the book. The attributes associated with it are $P_1, P_2, \ldots, P_{15}$. We format the following classes.

$C_1 = \{P_1\ P_2\ P_3\ P_4\ P_5\},$  $C_2 = \{P_3\ P_4\ P_5\ P_6\ P_7\}$
$C_3 = \{P_6\ P_7\ P_8\ P_9\ P_{10}\},$  $C_4 = \{P_8\ P_9\ P_{10}\ P_{11}\ P_{12}\}$
$C_5 = \{P_{11}\ P_{12}\ P_{13}\ P_{14}\ P_{15}\}$ and  $C_6 = \{P_{13}\ P_{14}\ P_{15}\ P_1\ P_2\}.$

Using these 6 overlapping classes we give the directed graph and their related connection matrices.

The directed graph related with $C_1 = \{P_1\ P_2\ P_3\ P_4\ P_5\}$ given by the expert is as follows:

The related connection matrix is given below:

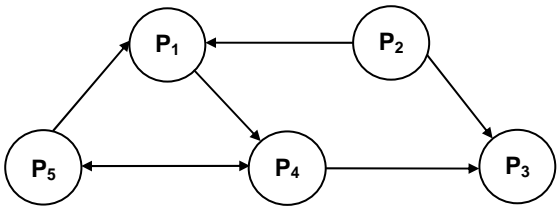

FIGURE 2.5.7

$$\begin{array}{c c c c c c} & P_1 & P_2 & P_3 & P_4 & P_5 \\ \begin{array}{c} P_1 \\ P_2 \\ P_3 \\ P_4 \\ P_5 \end{array} & \left[\begin{array}{c c c c c} 0 & 0 & 0 & 1 & 0 \\ 1 & 0 & 1 & 0 & 0 \\ 0 & 0 & 0 & 0 & 0 \\ 0 & 0 & 1 & 0 & 1 \\ 1 & 0 & 0 & 1 & 0 \end{array}\right] \end{array}$$

Now we give the directed graph of the expert for the class $C_2 = \{P_3 \ P_4 \ P_5 \ P_6 \ P_7\}$.

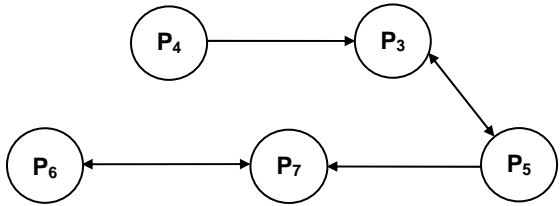

FIGURE 2.5.8

The related connection matrix

$$\begin{array}{c c c c c c} & P_3 & P_4 & P_5 & P_6 & P_7 \\ \begin{array}{c} P_3 \\ P_4 \\ P_5 \\ P_6 \\ P_7 \end{array} & \left[\begin{array}{c c c c c} 0 & 0 & 1 & 0 & 0 \\ 1 & 0 & 0 & 0 & 0 \\ 1 & 0 & 0 & 0 & 1 \\ 0 & 0 & 0 & 0 & 1 \\ 0 & 0 & 0 & 1 & 0 \end{array}\right] \end{array}$$

The directed graph related with the class $C_3 = \{P_6, P_7, P_8, P_9, P_{10}\}$ given by the expert

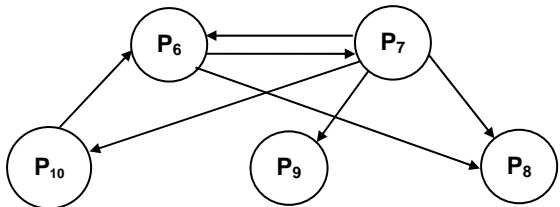

FIGURE 2.5.9

The connection matrix obtained from the directed graph is as follows:

$$
\begin{array}{c}
\phantom{P_{10}}\begin{array}{ccccc} P_6 & P_7 & P_8 & P_9 & P_{10} \end{array} \\
\begin{array}{c} P_6 \\ P_7 \\ P_8 \\ P_9 \\ P_{10} \end{array}
\left[\begin{array}{ccccc}
0 & 1 & 1 & 0 & 0 \\
1 & 0 & 1 & 1 & 1 \\
0 & 0 & 0 & 0 & 0 \\
0 & 0 & 0 & 0 & 0 \\
1 & 0 & 0 & 0 & 0
\end{array}\right]
\end{array}
$$

Now consider the class $C_4 = \{P_8,\ P_9,\ P_{10}\ P_{11}\ P_{12}\}$ the directed graph given by this expert is given below

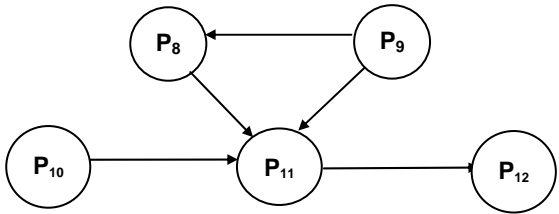

FIGURE 2.5.10

The related connection matrix

$$
\begin{array}{c}
\phantom{P_{10}}\begin{array}{ccccc} P_8 & P_9 & P_{10} & P_{11} & P_{12} \end{array} \\
\begin{array}{c} P_8 \\ P_9 \\ P_{10} \\ P_{11} \\ P_{12} \end{array}
\left[\begin{array}{ccccc}
0 & 0 & 0 & 1 & 0 \\
1 & 0 & 0 & 1 & 0 \\
0 & 0 & 0 & 1 & 0 \\
0 & 0 & 0 & 0 & 1 \\
0 & 0 & 0 & 1 & 0
\end{array}\right]
\end{array}
$$

The directed graph given by the expert related with the class $C_5 = \{P_{11}\ P_{12}\ P_{13}\ P_{14}\ P_{15}\}$ is given below:

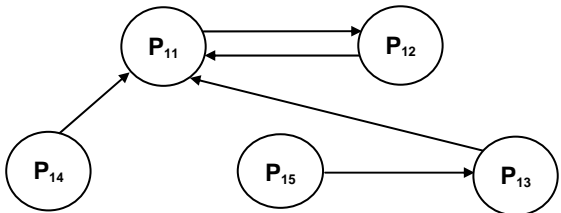

FIGURE 2.5.11

The related connection matrix is as follows:

$$\begin{array}{c} \\ P_{11} \\ P_{12} \\ P_{13} \\ P_{14} \\ P_{15} \end{array} \begin{array}{ccccc} P_{11} & P_{12} & P_{13} & P_{14} & P_{15} \\ \begin{bmatrix} 0 & 1 & 0 & 0 & 0 \\ 1 & 0 & 0 & 0 & 0 \\ 1 & 0 & 0 & 0 & 0 \\ 1 & 0 & 0 & 0 & 0 \\ 0 & 0 & 1 & 0 & 0 \end{bmatrix} \end{array}$$

The directed graph given by the expert for the last class of attributes $C_6 = \{P_{13}\ P_{14}\ P_{15}\ P_1\ P_2\}$

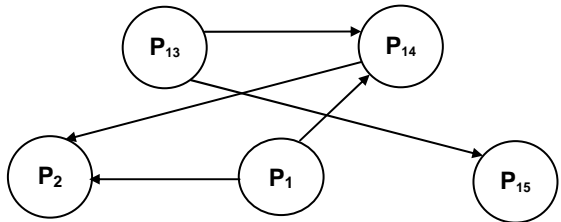

FIGURE 2.5.12

The related collection matrix

$$\begin{array}{c} \\ P_{13} \\ P_{14} \\ P_{15} \\ P_{1} \\ P_{2} \end{array} \begin{array}{ccccc} P_{13} & P_{14} & P_{15} & P_{1} & P_{2} \\ \begin{bmatrix} 0 & 1 & 1 & 0 & 0 \\ 0 & 0 & 0 & 0 & 1 \\ 0 & 0 & 0 & 0 & 0 \\ 0 & 1 & 0 & 0 & 1 \\ 0 & 0 & 0 & 0 & 0 \end{bmatrix} \end{array}$$

Using the six classes of over lapping attributes we formulate the combined block overlap connection matrix V.

$$\begin{array}{c} \quad\quad P_1 \; P_2 \; P_3 \; P_4 \; P_5 \; P_6 \; P_7 \; P_8 \; P_9 \; P_{10} \; P_{11} \; P_{12} \; P_{13} \; P_{14} \, P_{15} \\ \begin{array}{c} P_1 \\ P_2 \\ P_3 \\ P_4 \\ P_5 \\ P_6 \\ P_7 \\ P_8 \\ P_9 \\ P_{10} \\ P_{11} \\ P_{12} \\ P_{13} \\ P_{14} \\ P_{15} \end{array} \left[ \begin{array}{ccccccccccccccc} 0 & 1 & 0 & 1 & 0 & 0 & 0 & 0 & 0 & 0 & 0 & 0 & 0 & 1 & 0 \\ 1 & 0 & 1 & 0 & 0 & 0 & 0 & 0 & 0 & 0 & 0 & 0 & 0 & 0 & 0 \\ 0 & 0 & 0 & 0 & 1 & 0 & 0 & 0 & 0 & 0 & 0 & 0 & 0 & 0 & 0 \\ 0 & 0 & 2 & 0 & 1 & 0 & 0 & 0 & 0 & 0 & 0 & 0 & 0 & 0 & 0 \\ 1 & 0 & 0 & 1 & 0 & 0 & 1 & 0 & 0 & 0 & 0 & 0 & 0 & 0 & 0 \\ 0 & 0 & 1 & 0 & 0 & 0 & 2 & 1 & 0 & 0 & 0 & 0 & 0 & 0 & 0 \\ 0 & 0 & 0 & 0 & 0 & 2 & 0 & 1 & 1 & 1 & 0 & 0 & 0 & 0 & 0 \\ 0 & 0 & 0 & 0 & 0 & 0 & 1 & 0 & 0 & 0 & 1 & 0 & 0 & 0 & 0 \\ 0 & 0 & 0 & 0 & 0 & 0 & 0 & 1 & 0 & 0 & 1 & 0 & 0 & 0 & 0 \\ 0 & 0 & 0 & 0 & 0 & 1 & 0 & 0 & 0 & 1 & 0 & 0 & 0 & 0 & 0 \\ 0 & 0 & 0 & 0 & 0 & 0 & 0 & 0 & 0 & 0 & 0 & 2 & 0 & 0 & 0 \\ 0 & 0 & 0 & 0 & 0 & 0 & 0 & 0 & 0 & 0 & 2 & 0 & 0 & 0 & 0 \\ 0 & 0 & 0 & 0 & 0 & 0 & 0 & 0 & 0 & 0 & 1 & 0 & 0 & 1 & 1 \\ 0 & 1 & 0 & 0 & 0 & 0 & 0 & 0 & 0 & 0 & 1 & 0 & 0 & 0 & 0 \\ 0 & 0 & 0 & 0 & 0 & 0 & 0 & 0 & 0 & 0 & 0 & 0 & 1 & 0 & 0 \end{array} \right] \end{array}$$

Now V gives the $15 \times 15$ matrix of the dynamical system related with the combined block overlap fuzzy cognitive map. Suppose we are interested in studying the effect of state vectors.

Let us consider the state vector Y = (1 0 0 0 1 0 0 0 0 0 0 0 0 0 0) where only the two attributes, No binding with the family and the socially irresponsible co ordinates are in the on state and all other attributes are in the off state. To study the effect of Y on the dynamical system V;

YV $\longmapsto$ (1 1 0 1 1 0 1 0 0 0 0 0 0 1 0) = $Y_1$ (say)

$Y_1$V $\longmapsto$ (1 1 1 1 1 1 1 1 1 1 1 0 1 1 0) = $Y_2$ (say)

$Y_2$V $\longmapsto$ (1 1 1 1 1 1 1 1 1 1 1 1 0 1 0) = $Y_3$ (say)

$Y_3$ V $\longmapsto$ (1 1 1 1 1 1 1 1 1 1 1 1 1 0 1 0) = $Y_4 = Y_3$
(The hidden pattern of the system).

Thus the hidden pattern of the system is a fixed point. When the migrant labourer has no binding with the family and is socially irresponsible we see he has all the other nodes to be working on his nature i.e., they all come to on state only the nodes unreachable by friends or relatives and failure of agriculture is in the off state.

Now using the C-program given in the appendix 5 the reader is given the task of finding the effects of each and every state vector.

Also it is left as an exercise for the reader to compare the matrix P given in page 57-59 and the matrix V and find which approach to this problem is better.

Now we give the final model of this chapter namely blocks of varying size with varying overlap. Consider the 12 attribute model studied relative to the problems of HIV/AIDS affected migrant labourers in page 42-50.

Consider the 12 attributes $A_1$, $A_2$, …, $A_{12}$. Divide them into overlapping blocks

$C_1$ = $\{A_1 A_2 A_3 A_4 A_5\}$, $C_2$ = $\{A_4 A_5 A_6 A_7\}$,
$C_3$ = $\{A_6 A_7 A_8 A_9\}$, $C_4$ = $\{A_8 A_9 A_{10}\}$ and
$C_5$ = $\{A_{10} A_{11} A_{12} A_1 A_2 A_3\}$.

Now we analyze each of these classes of attributes $C_1$ $C_2$ $C_3$ $C_4$ and $C_5$. The directed graph given by the expert for the attributes $\{A_1 A_2 A_3 A_4 A_5\}$.

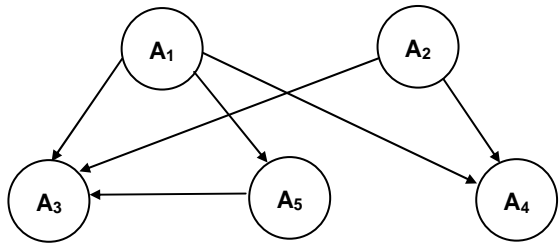

FIGURE 2.5.13

The related connection matrix.

$$
\begin{array}{c}
\phantom{A_1} \begin{array}{ccccc} A_1 & A_2 & A_3 & A_4 & A_5 \end{array} \\
\begin{array}{c} A_1 \\ A_2 \\ A_3 \\ A_4 \\ A_5 \end{array}
\left[ \begin{array}{ccccc}
0 & 0 & 1 & 1 & 1 \\
0 & 0 & 1 & 1 & 0 \\
0 & 0 & 0 & 0 & 0 \\
0 & 0 & 0 & 0 & 0 \\
0 & 0 & 1 & 0 & 0
\end{array} \right]
\end{array}
$$

The directed graph related to the attributes {$A_4$, $A_5$, $A_6$, $A_7$} given by an expert.

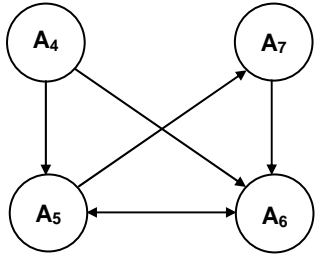

FIGURE 2.5.14

The related connection matrix.

$$
\begin{array}{c}
\quad\ A_4\ A_5\ A_6\ A_7 \\
\begin{array}{c} A_4 \\ A_5 \\ A_6 \\ A_7 \end{array}
\begin{bmatrix}
0 & 1 & 1 & 0 \\
0 & 0 & 1 & 1 \\
0 & 1 & 0 & 0 \\
0 & 0 & 1 & 0
\end{bmatrix}
\end{array}
$$

The directed graph given by the expert related to the attributes $A_6$, $A_7$, $A_8$ and $A_9$.

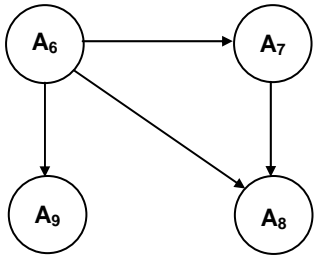

FIGURE 2.5.15

The related connection matrix.

$$
\begin{array}{c}
\quad\ A_6\ A_7\ A_8\ A_9 \\
\begin{array}{c} A_6 \\ A_7 \\ A_8 \\ A_9 \end{array}
\begin{bmatrix}
0 & 1 & 1 & 1 \\
0 & 0 & 1 & 0 \\
0 & 0 & 0 & 0 \\
0 & 0 & 0 & 0
\end{bmatrix}
\end{array}
$$

Now we give the directed graph given by an expert using the attributes $A_8$, $A_9$ and $A_{10}$.

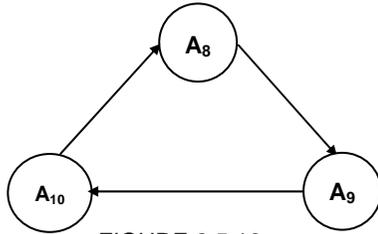

FIGURE 2.5.16

The connection matrix related with this directed graph.

$$
\begin{array}{c}
 & \begin{array}{ccc} A_8 & A_9 & A_{10} \end{array} \\
\begin{array}{c} A_8 \\ A_9 \\ A_{10} \end{array} &
\left[\begin{array}{ccc}
0 & 1 & 0 \\
0 & 0 & 1 \\
1 & 0 & 0
\end{array}\right]
\end{array}
$$

Finally we obtain the connection matrix given by the expert related with the attributes $A_{10}$, $A_{11}$, $A_{12}$, $A_1$, $A_2$ and $A_3$.

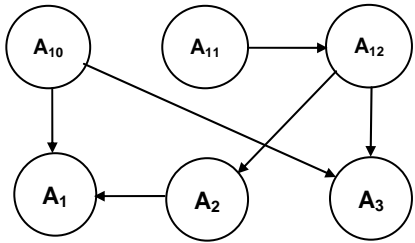

FIGURE 2.5.17

The related connection matrix.

$$
\begin{array}{c}
 & \begin{array}{cccccc} A_{10} & A_{11} & A_{12} & A_1 & A_2 & A_3 \end{array} \\
\begin{array}{c} A_{10} \\ A_{11} \\ A_{12} \\ A_1 \\ A_2 \\ A_3 \end{array} &
\left[\begin{array}{cccccc}
0 & 0 & 0 & 1 & 0 & 1 \\
0 & 0 & 1 & 0 & 0 & 0 \\
0 & 0 & 0 & 0 & 1 & 1 \\
0 & 0 & 0 & 0 & 0 & 0 \\
1 & 0 & 0 & 0 & 0 & 0 \\
0 & 0 & 0 & 0 & 0 & 0
\end{array}\right]
\end{array}
$$

Now we give the related connection combined block over lapping matrix N.

$$\begin{array}{c} \\ A_1 \\ A_2 \\ A_3 \\ A_4 \\ A_5 \\ A_6 \\ A_7 \\ A_8 \\ A_9 \\ A_{10} \\ A_{11} \\ A_{12} \end{array} \begin{array}{cccccccccccc} A_1 & A_2 & A_3 & A_4 & A_5 & A_6 & A_7 & A_8 & A_9 & A_{10} & A_{11} & A_{12} \\ \left[\begin{array}{cccccccccccc} 0 & 0 & 1 & 1 & 1 & 0 & 0 & 0 & 0 & 0 & 0 & 0 \\ 0 & 0 & 1 & 1 & 0 & 0 & 0 & 0 & 0 & 1 & 0 & 0 \\ 0 & 0 & 0 & 0 & 0 & 0 & 0 & 0 & 0 & 0 & 0 & 0 \\ 0 & 0 & 0 & 0 & 1 & 1 & 0 & 0 & 0 & 0 & 0 & 0 \\ 0 & 0 & 1 & 0 & 0 & 1 & 1 & 0 & 0 & 0 & 0 & 0 \\ 0 & 0 & 0 & 0 & 1 & 0 & 1 & 1 & 1 & 0 & 0 & 0 \\ 0 & 0 & 0 & 0 & 0 & 1 & 0 & 1 & 0 & 0 & 0 & 0 \\ 0 & 0 & 0 & 0 & 0 & 0 & 0 & 0 & 1 & 0 & 0 & 0 \\ 0 & 0 & 0 & 0 & 0 & 0 & 0 & 0 & 0 & 1 & 0 & 0 \\ 1 & 0 & 1 & 0 & 0 & 0 & 0 & 1 & 0 & 0 & 0 & 0 \\ 0 & 0 & 0 & 0 & 0 & 0 & 0 & 0 & 0 & 0 & 0 & 1 \\ 0 & 1 & 1 & 0 & 0 & 0 & 0 & 0 & 0 & 0 & 0 & 0 \end{array}\right] \end{array}$$

Now we can analyze the effect of any state vector on the dynamical system N.

Let X = (0 0 0 1 0 0 0 0 0 0 0 0) be the state vector with the only attribute profession in the on state the effect of X on N is given by

XN $\hookrightarrow$ (0 0 0 1 1 1 0 0 0 0 0 0) = $X_1$ (say)

$X_1$N $\hookrightarrow$ (0 0 1 1 1 1 1 1 1 0 0 0) = $X_2$ (say)

$X_2$N $\hookrightarrow$ (0 0 1 1 1 1 1 1 1 0 0 0) = $X_3$ (= $X_2$).

The hidden pattern is a fixed point. Only the vector Easily availability of money lack of education more leisure, machismo / exaggerated masculinity and no awareness of the disease in the off state and all other attributes come to the on state. We see there is difference between resultant vectors when using the dynamical system A and N.

Now we consider the state vector $R_1$ = (1 0 0 0 0 0 1 0 0 1 0 0) in which the attributes easy money, no social responsibility and more leisure i.e., $A_1$, $A_7$ and $A_{10}$ are in the on state and all other attributes are in the off state. The effect of $R_1$ on the dynamical system N is given by

$$R_1N \quad \mapsto \quad (1\ 0\ 1\ 1\ 1\ 1\ 1\ 1\ 0\ 1\ 0\ 0) \quad = \quad R_2 \text{ (say)}$$

$$R_2N \quad \mapsto \quad (1\ 0\ 1\ 1\ 1\ 1\ 1\ 1\ 1\ 1\ 0\ 0) \quad = \quad R_3 \text{ (say)}$$

$$R_3N \quad \mapsto \quad (1\ 0\ 1\ 1\ 1\ 1\ 1\ 1\ 1\ 1\ 0\ 0) \quad = \quad R_4 \quad = \quad R_3$$

$R_3$ is a fixed point. The hidden pattern is a fixed point only the attributes $A_2$, $A_{11}$ and $A_{12}$ remain in the off state all other vector became on.

The reader can compare this with the resultant vector given by A in pages 42 to 50.

Similar studies using the C –program is carried out.

Now we use the attributes $P_1$, $P_2$,…$P_{15}$ given in page 56 and model it using the combined block overlap FCM. First we divide $P_1$, $P_2$,…, $P_5$ into 5 classes

$$
\begin{aligned}
C_1 &= \{P_1\ P_2\ P_3\ P_4\ P_5\ P_6\}, \\
C_2 &= \{P_4\ P_5\ P_6\ P_7\ P_8\ P_9\ P_{10}\}, \\
C_3 &= \{P_7\ P_8\ P_9\ P_{10}\ P_{11}\}, \\
C_4 &= \{P_{11}, P_{12}, P_{13}\}
\end{aligned}
$$

and $\quad C_5 \quad = \quad \{P_{12}\ P_{13}\ P_{14}\ P_{15}\ P_1\ P_2\ P_3\}.$

The expert's directed graph using the attributes $P_1$, $P_2$,…, $P_6$.

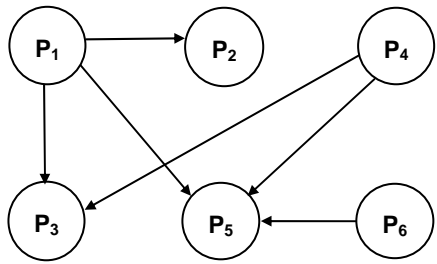

FIGURE 2.5.18

The resulting connection matrix using the directed graph is given below.

$$\begin{array}{c c}
 & \begin{array}{c c c c c c} P_1 & P_2 & P_3 & P_4 & P_5 & P_6 \end{array} \\
\begin{array}{c} P_1 \\ P_2 \\ P_3 \\ P_4 \\ P_5 \\ P_6 \end{array} &
\left[ \begin{array}{c c c c c c}
0 & 1 & 1 & 0 & 1 & 0 \\
0 & 0 & 0 & 0 & 0 & 0 \\
0 & 0 & 0 & 0 & 0 & 0 \\
0 & 0 & 1 & 0 & 1 & 0 \\
0 & 0 & 0 & 0 & 0 & 0 \\
0 & 0 & 0 & 0 & 1 & 0 
\end{array} \right]
\end{array}$$

For the attributes $P_4$, $P_5$,…, $P_9$, $P_{10}$ using the experts opinion we have the following directed graph.

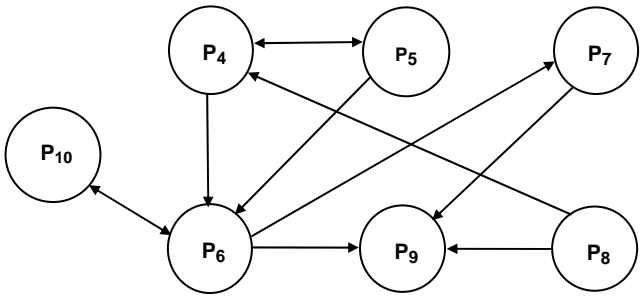

FIGURE 2.5.19

The related connection matrix

$$\begin{array}{c c}
 & \begin{array}{c c c c c c c} P_4 & P_5 & P_6 & P_7 & P_8 & P_9 & P_{10} \end{array} \\
\begin{array}{c} P_4 \\ P_5 \\ P_6 \\ P_7 \\ P_8 \\ P_9 \\ P_{10} \end{array} &
\left[ \begin{array}{c c c c c c c}
0 & 1 & 1 & 0 & 0 & 0 & 0 \\
1 & 0 & 1 & 0 & 0 & 0 & 0 \\
0 & 0 & 0 & 1 & 0 & 1 & 1 \\
0 & 0 & 1 & 0 & 0 & 0 & 0 \\
1 & 0 & 0 & 0 & 0 & 1 & 0 \\
0 & 0 & 0 & 0 & 0 & 0 & 0 \\
0 & 0 & 1 & 0 & 0 & 0 & 0 
\end{array} \right]
\end{array}$$

The directed graph given by the expert for the attributes $P_7$ $P_8…P_{11}$ is given below.

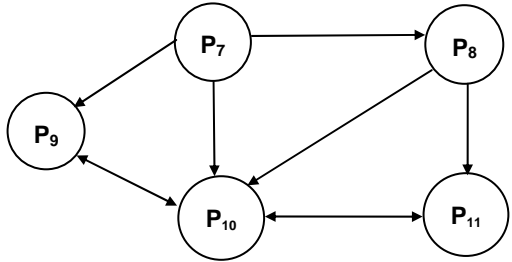

FIGURE 2.5.20

The related connection matrix

$$
\begin{array}{c}
\begin{array}{ccccc} P_7 & P_8 & P_9 & P_{10} & P_{11} \end{array} \\
\begin{array}{c} P_7 \\ P_8 \\ P_9 \\ P_{10} \\ P_{11} \end{array}
\begin{bmatrix}
0 & 1 & 1 & 1 & 0 \\
0 & 0 & 0 & 1 & 1 \\
0 & 0 & 0 & 1 & 0 \\
0 & 0 & 1 & 0 & 1 \\
0 & 0 & 0 & 1 & 0
\end{bmatrix}
\end{array}
$$

Now the directed graph for the three attributes $P_{11}$ $P_{12}$ and $P_{13}$.

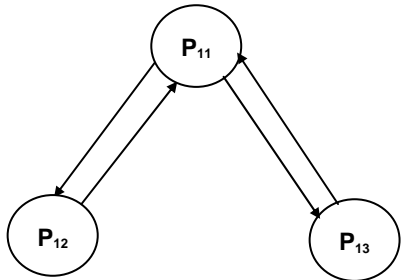

FIGURE 2.5.21

The related connection matrix

$$
\begin{array}{c}
\phantom{P_{11}} \begin{array}{ccc} P_{11} & P_{12} & P_{13} \end{array} \\
\begin{array}{c} P_{11} \\ P_{12} \\ P_{13} \end{array}
\left[
\begin{array}{ccc}
0 & 1 & 1 \\
1 & 0 & 0 \\
1 & 0 & 0
\end{array}
\right]
\end{array}
$$

We give the expert's directed graph using $P_{12}$ $P_{13}$ $P_{14}$ $P_{15}$ $P_1$ $P_2$ $P_3$.

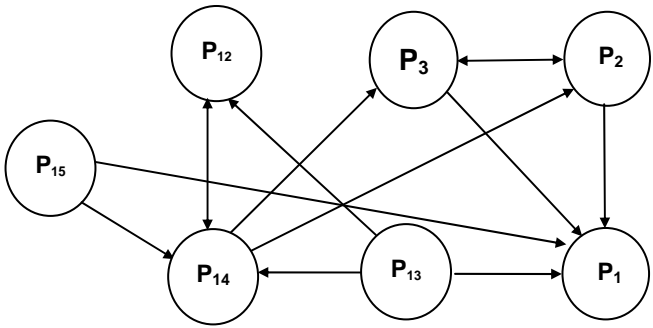

FIGURE 2.5.22

The related connection matrix is

$$
\begin{array}{c}
\phantom{P_{12}} \begin{array}{ccccccc} P_{12} & P_{13} & P_{14} & P_{15} & P_1 & P_2 & P_3 \end{array} \\
\begin{array}{c} P_{12} \\ P_{13} \\ P_{14} \\ P_{15} \\ P_1 \\ P_2 \\ P_3 \end{array}
\left[
\begin{array}{ccccccc}
0 & 0 & 1 & 0 & 0 & 0 & 0 \\
1 & 0 & 1 & 0 & 1 & 0 & 0 \\
1 & 0 & 0 & 0 & 0 & 0 & 1 \\
0 & 0 & 0 & 0 & 1 & 0 & 0 \\
0 & 0 & 0 & 0 & 0 & 0 & 0 \\
0 & 0 & 0 & 0 & 1 & 0 & 1 \\
0 & 0 & 0 & 0 & 1 & 1 & 0
\end{array}
\right]
\end{array}
$$

Let O denote the connection matrix of the combined block overlap FCM.

$$
\begin{array}{c}
\quad\; P_1\; P_2\; P_3\; P_4\; P_5\; P_6\; P_7\; P_8\; P_9\; P_{10}P_{11}\,P_{12}P_{13}P_{14}P_{15} \\
\begin{array}{c}
P_1 \\ P_2 \\ P_3 \\ P_4 \\ P_5 \\ P_6 \\ P_7 \\ P_8 \\ P_9 \\ P_{10} \\ P_{11} \\ P_{12} \\ P_{13} \\ P_{14} \\ P_{15}
\end{array}
\left[
\begin{array}{ccccccccccccccc}
0 & 1 & 1 & 0 & 1 & 0 & 0 & 0 & 0 & 0 & 0 & 0 & 0 & 0 & 0 \\
1 & 0 & 1 & 0 & 0 & 0 & 0 & 0 & 0 & 0 & 0 & 0 & 0 & 0 & 0 \\
1 & 1 & 0 & 0 & 0 & 0 & 0 & 0 & 0 & 0 & 0 & 0 & 0 & 0 & 0 \\
0 & 0 & 1 & 0 & 2 & 1 & 0 & 0 & 0 & 0 & 0 & 0 & 0 & 0 & 0 \\
0 & 0 & 0 & 1 & 0 & 1 & 0 & 0 & 0 & 0 & 0 & 0 & 0 & 0 & 0 \\
0 & 0 & 0 & 0 & 1 & 0 & 1 & 0 & 1 & 1 & 0 & 0 & 0 & 0 & 0 \\
0 & 0 & 0 & 0 & 0 & 1 & 0 & 1 & 1 & 1 & 0 & 0 & 0 & 0 & 0 \\
0 & 0 & 0 & 1 & 0 & 0 & 0 & 0 & 1 & 1 & 1 & 0 & 0 & 0 & 0 \\
0 & 0 & 0 & 0 & 0 & 0 & 0 & 0 & 1 & 1 & 0 & 0 & 0 & 0 & 0 \\
0 & 0 & 0 & 0 & 0 & 1 & 0 & 0 & 0 & 1 & 1 & 0 & 0 & 0 & 0 \\
0 & 0 & 0 & 0 & 0 & 0 & 0 & 0 & 0 & 1 & 0 & 1 & 1 & 0 & 0 \\
0 & 0 & 0 & 0 & 0 & 0 & 0 & 0 & 0 & 0 & 1 & 0 & 0 & 0 & 0 \\
1 & 0 & 0 & 0 & 0 & 0 & 0 & 0 & 0 & 0 & 0 & 0 & 1 & 0 & 1 \\
0 & 0 & 1 & 0 & 0 & 0 & 0 & 0 & 0 & 0 & 0 & 0 & 1 & 0 & 0 \\
1 & 0 & 0 & 0 & 0 & 0 & 0 & 0 & 0 & 0 & 0 & 0 & 0 & 0 & 0
\end{array}
\right]
\end{array}
$$

Now consider a state vector S = (0 1 0 0 0 0 0 0 0 0 0 0 0 0 0 ) the notion male ego alone is in the on state and all other attributes are in the off state.

The effect of S on the dynamical system O is given by

$$SO \;\hookrightarrow\; (1\;1\;1\;0\;0\;0\;0\;0\;0\;0\;0\;0\;0\;0\;0) \qquad = \quad S_1 \text{ (say)}$$

$$S_1O \;\hookrightarrow\; (1\;1\;1\;0\;1\;0\;0\;0\;0\;0\;0\;0\;0\;0\;0) \qquad = \quad S_2 \text{ (say)}$$

$$S_2O \;\hookrightarrow\; (1\;1\;1\;1\;1\;1\;1\;0\;0\;0\;0\;0\;0\;0\;0) \qquad = \quad S_3 \text{ (say)}$$

$$S_3O \;\hookrightarrow\; (1\;1\;1\;1\;1\;1\;1\;1\;0\;1\;1\;0\;0\;0\;0) \qquad = \quad S_4$$

$$S_4\,O \;\hookrightarrow\; (1\;1\;1\;1\;1\;1\;1\;1\;1\;1\;1\;1\;0\;0\;0\;0) \qquad = \quad S_5 \text{ (say)}$$

$$S_5\,O \;\hookrightarrow\; (1\;1\;1\;1\;1\;1\;1\;1\;1\;1\;1\;1\;1\;0\;0) \qquad = \quad S_6 \text{ (say)}$$

$$S_6\,O \;\hookrightarrow\; (1\;1\;1\;1\;1\;1\;1\;1\;1\;1\;1\;1\;1\;1\;0) \qquad = \quad S_7 \,(S_6).$$

Thus the hidden pattern is a fixed point. All vector come to on state except the notion $P_{15}$. Using the C-program the on state of several vectors have been found out. They are used in deriving our conclusions.

An interested reader may compare the different techniques used in this study and derive which model is best suited for this problem.

In the next chapter we proceed on to use the BAM model for the analysis of the HIV/AIDS infected labourers and their psycho, socio economic problems.

Chapter Three

# USE OF BAM TO STUDY THE FACTORS FOR MIGRATION AND VULNERABILITY TO HIV/AIDS

We live in a world of marvelous complexity and variety, a world where events never repeat exactly. Even though events are never exactly the same they are also not completely different. There is a threat of continuity, similarity and predictability that allows us to generalize often correctly from past experience to future events. Neural networks or neuro-computing or brain like computation is based on the wistful hope that we can reproduce at least some of the flexibility and power of the human brain by artificial brains. Neural networks consists of many simple computing elements generally simple non linear summing junctions connected together by connections of varying strength a gross abstraction of the brain which consists of very large number of far more complex neurons connected together with far more complex and far more structured couplings, neural networks architecture cover a wide range.

In one sense every computer is a neural net, because we can view traditional digital logic as constructed from inter connected McCullouch-Pitts neurons. McCullouch-Pitts neurons were proposed in 1943 as models of biological neurons and arranged in networks for a specific purpose of computing logic functions. The problems where artificial neural networks have the most promise are those with a real-world flavor: medical research, signal processing, sociological problems etc.

Neural networks helps to solve these problems with natural mechanisms of generalizations. To over-simplify, suppose we represent an object in a network as a pattern of activation of several units. If a unit or two responds incorrectly the overall pattern stays pretty much of the same, and the network still respond correctly to stimuli when neural networks operate similar inputs naturally produce similar outputs. Most real world perceptual problems have this structure of input-output continuity. The prototype model provides a model for human categorization with a great deal of psychological support. The computational strategy leads to some curious human psychology. For instance in United States people imagine a prototype bird that looks somewhat like a sparrow or a robin. So they learn to judge 'penguins' or "ostriches" as "bad" birds because these birds do not resemble the prototype bird even though they are birds "Badness" shows up in a number of ways.

Neural networks naturally develop this kind of category structure. The problems that neural networks solved well and solved poorly were those where human showed comparable strengths and weaknesses in their cognitive computations". For this reason until quite recently most of the study of neural networks has been carried out by psychologists and cognitive scientists who sought models of human cognitive function. Neural networks deal with uncertainty as humans do, not by deliberate design but as a by product of their parallel distributed structure. Because general statements about both human psychology and the structure of the world embed so deeply in both neural networks and fuzzy systems, it is very appropriate to study the psychological effects of HIV/AIDS patients and their influence on public using this theory. Like social customs these assumptions are obvious only if you grew up with them. Both neural networks and fuzzy systems break with the historical tradition prominent in western thought and we can precisely and unambiguously characterize the world, divide into two categories and then manipulate these descriptions according to precise and formal rules. Huang Po, a Buddhist teacher of the ninth century observed that "To make use of your minds to think conceptually is to leave the substance and attach yourself to form", and "from discrimination between this and that a host of demons blazes forth".

## 3.1 Some Basic Concepts of BAM

Now we go forth to describe the mathematical structure of the Bidirectional Associative Memories (BAM) model. Neural networks recognize ill defined problems without an explicit set of rules. Neurons behave like functions, neurons transduce an unbounded input activation x(t) at time t into a bounded output signal S(x(t)) i.e. Neuronal activations change with time.

Artificial neural networks consists of numerous simple processing units or neurons which can be trained to estimate sampled functions when we do not know the form of the functions. A group of neurons form a field. Neural networks contain many field of neurons. In our text $F_x$ will denote a neuron field, which contains n neurons, and $F_y$ denotes a neuron field, which contains p neurons. The neuronal dynamical system is described by a system of first order differential equations that govern the time-evolution of the neuronal activations or which can be called also as membrane potentials.

$$\dot{x}_i \quad = \quad g_i(X, Y, ...)$$
$$\dot{y}_j \quad = \quad h_j(X, Y, ...)$$

where $\dot{x}_i$ and $\dot{y}_j$ denote respectively the activation time function of the $i^{th}$ neuron in $F_X$ and the $j^{th}$ neuron in $F_Y$. The over dot denotes time differentiation, $g_i$ and $h_j$ are some functions of X, Y, ... where $X(t) = (x_1(t), ... , x_n(t))$ and $Y(t) = (y_1(t), ... , y_p(t))$ define the state of the neuronal dynamical system at time t.

The passive decay model is the simplest activation model, where in the absence of the external stimuli, the activation decays in its resting value

$$\dot{x}_i \quad = \quad x_i$$
and
$$\dot{y}_j \quad = \quad y_j$$

The passive decay rate $A_i > 0$ scales the rate of passive decay to the membranes resting potentials $\dot{x}_i = -A_i x_i$. The default rate is $A_i = 1$, i.e. $\dot{x}_i = -A_i x_i$. The membrane time constant $C_i > 0$ scales the time variables of the activation dynamical system. The default time constant is $C_i = 1$. Thus $C_i \dot{x}_i = -A_i x_i$.

The membrane resting potential $P_i$ is defined as the activation value to which the membrane potential equilibrates in the absence of external inputs. The resting potential is an additive constant and its default value is zero. It need not be positive.

$$P_i \quad = \quad C_i \dot{x}_i + A_i x_i$$
$$I_i \quad = \quad \dot{x}_i + x_i$$

is called the external input of the system. Neurons do not compute alone. Neurons modify their state activations with external input and with feed back from one another. Now, how do we transfer all these actions of neurons activated by inputs their resting potential etc. mathematically. We do this using what are called synaptic connection matrices.

Let us suppose that the field $F_X$ with n neurons is synaptically connected to the field $F_Y$ of p neurons. Let $m_{ij}$ be a synapse where the axon from the $i^{th}$ neuron in $F_X$ terminates. $M_{ij}$ can be positive, negative or zero. The synaptic matrix M is a n by p matrix of real numbers whose entries are the synaptic efficacies $m_{ij}$.

The matrix M describes the forward projections from the neuronal field $F_X$ to the neuronal field $F_Y$. Similarly a p by n synaptic matrix N describes the backward projections from $F_Y$ to $F_X$. Unidirectional networks occur when a neuron field synaptically intra connects to itself. The matrix M be a n by n square matrix. A Bidirectional network occur if $M = N^T$ and $N = M^T$. To describe this synaptic connection matrix more simply, suppose the n neurons in the field $F_X$ synaptically connect to the p-neurons in field $F_Y$. Imagine an axon from the $i^{th}$ neuron in $F_X$ that terminates in a synapse $m_{ij}$, that about the $j^{th}$ neuron in $F_Y$. We assume that the real number $m_{ij}$ summarizes the synapse and that $m_{ij}$ changes so slowly relative to activation fluctuations that is constant.

Thus we assume no learning if $m_{ij} = 0$ for all t. The synaptic value $m_{ij}$ might represent the average rate of release of a neuro-transmitter such as norepinephrine. So, as a rate, $m_{ij}$ can be positive, negative or zero.

When the activation dynamics of the neuronal fields $F_X$ and $F_Y$ lead to the overall stable behaviour the bidirectional networks are called as Bidirectional Associative Memories (BAM). As in the analysis of the HIV/AIDS patients relative to the migrancy we state that the BAM model studied presently and predicting the future after a span of 5 or 10 years may not be the same.

For the system would have reached stability and after the lapse of this time period the activation neurons under investigations and which are going to measure the model would be entirely different.

Thus from now onwards more than the uneducated poor the educated rich and the middle class will be the victims of HIV/AIDS. So for this study presently carried out can only give how migration has affected the life style of poor labourer and had led them to be victims of HIV/AIDS.

Further not only a Bidirectional network leads to BAM also a unidirectional network defines a BAM if M is symmetric i.e. M = $M^T$. We in our analysis mainly use BAM which are bidirectional networks. However we may also use unidirectional BAM chiefly depending on the problems under investigations. We briefly describe the BAM model more technically and mathematically.

An additive activation model is defined by a system of n + p coupled first order differential equations that inter connects the fields $F_X$ and $F_Y$ through the constant synaptic matrices M and N.

$$x_i = -A_i x_i + \sum_{j=1}^{p} S_j(y_j) n_{ji} + I_i \qquad (3.1)$$

$$y_i = -A_j y_j + \sum_{i=1}^{n} S_i(x_i) m_{ij} + J_j \qquad (3.2)$$

$S_i(x_i)$ and $S_j(y_j)$ denote respectively the signal function of the $i^{th}$ neuron in the field $F_X$ and the signal function of the $j^{th}$ neuron in the field $F_Y$.

Discrete additive activation models correspond to neurons with threshold signal functions.

The neurons can assume only two values ON and OFF. ON represents the signal +1, OFF represents 0 or – 1 (– 1 when the representation is bipolar). Additive bivalent models describe asynchronous and stochastic behaviour.

At each moment each neuron can randomly decide whether to change state or whether to emit a new signal given its current activation. The Bidirectional Associative Memory or BAM is a non adaptive additive bivalent neural network. In neural literature the discrete version of the equation (3.1) and (3.2) are often referred to as BAMs.

A discrete additive BAM with threshold signal functions arbitrary thresholds inputs an arbitrary but a constant synaptic connection matrix M and discrete time steps K are defined by the equations

$$x_i^{k+1} = \sum_{j=1}^{p} S_j(y_j^k) m_{ij} + I_i \qquad (3.3)$$

$$y_j^{k+1} = \sum_{i=1}^{n} S_i\left(x_i^{k}\right) m_{ij} + J_i \qquad (3.4)$$

where $m_{ij} \in M$ and $S_i$ and $S_j$ are signal functions. They represent binary or bipolar threshold functions. For arbitrary real valued thresholds $U = (U_1, ..., U_n)$ for $F_X$ neurons and $V = (V_1, ..., V_P)$ for $F_Y$ neurons the threshold binary signal functions corresponds to

$$S_i(x_i^k) = \begin{cases} 1 & \text{if} \quad x_i^k > U_i \\ S_i(x_i^{k-1}) & \text{if} \ x_i^k = U_i \\ 0 & \text{if} \quad x_i^k < U_i \end{cases} \qquad (3.5)$$

and

$$S_j(x_j^k) = \begin{cases} 1 & \text{if} \quad y_j^k > V_j \\ S_j(y_j^{k-1}) & \text{if} \ y_j^k = V_j \\ 0 & \text{if} \quad y_j^k < V_j \end{cases} \qquad (3.6)$$

The bipolar version of these equations yield the signal value -1 when $x_i < U_i$ or when $y_j < V_j$. The bivalent signal functions allow us to model complex asynchronous state change patterns. At any moment different neurons can decide whether to compare their activation to their threshold. At each moment any of the 2n subsets of $F_X$ neurons or 2p subsets of the $F_Y$ neurons can decide to change state. Each neuron may randomly decide whether to check the threshold conditions in the equations (3.5) and (3.6). At each moment each neuron defines a random variable that can assume the value ON(+1) or OFF(0 or -1). The network is often assumed to be deterministic and state changes are synchronous i.e. an entire field of neurons is updated at a time. In case of simple asynchrony only one neuron makes a state change decision at a time. When the subsets represent the entire fields $F_X$ and $F_Y$ synchronous state change results.

In a real life problem the entries of the constant synaptic matrix M depends upon the investigator's feelings. The synaptic matrix is given a weightage according to their feelings. If $x \in F_X$ and $y \in F_Y$ the forward projections from $F_X$ to $F_Y$ is defined by the matrix M. $\{F(x_i, y_j)\} = (m_{ij}) = M, 1 \le i \le n, 1 \le j \le p$.

The backward projections is defined by the matrix $M^T$. $\{F(y_i, x_i)\} = (m_{ji}) = M^T$, $1 \leq i \leq n$, $1 \leq j \leq p$. It is not always true that the backward projections from $F_Y$ to $F_X$ is defined by the matrix $M^T$.

Now we just recollect the notion of bidirectional stability. All BAM state changes lead to fixed point stability. The property holds for synchronous as well as asynchronous state changes.

A BAM system $(F_X, F_Y, M)$ is bidirectionally stable if all inputs converge to fixed point equilibria. Bidirectional stability is a dynamic equilibrium. The same signal information flows back and forth in a bidirectional fixed point. Let us suppose that A denotes a binary n-vector and B denotes a binary p-vector. Let A be the initial input to the BAM system. Then the BAM equilibrates to a bidirectional fixed point $(A_f, B_f)$ as

$$
\begin{array}{lllll}
A & \rightarrow & M & \rightarrow & B \\
A' & \leftarrow & M^T & \leftarrow & B \\
A' & \rightarrow & M & \rightarrow & B' \\
A'' & \leftarrow & M^T & \leftarrow & B' \quad \text{etc.} \\
A_f & \rightarrow & M & \rightarrow & B_f \\
A_f & \leftarrow & M^T & \leftarrow & B_f \quad \text{etc.}
\end{array}
$$

where A', A'', ... and B', B'', ... represents intermediate or transient signal state vectors between respectively A and $A_f$ and B and $B_f$. The fixed point of a Bidirectional system is time dependent.

The fixed point for the initial input vectors can be attained at different times. Based on the synaptic matrix M which is developed by the investigators feelings the time at which bidirectional stability is attained also varies accordingly.

## 3.2 USE OF BAM MODEL TO STUDY THE CAUSE OF VULNERABILITY TO HIV/AIDS AND FACTORS FOR MIGRATION

Now the object is to study the levels of knowledge and awareness relating to STD/HIV/AIDS existing among the migrant labourers in Tamil Nadu; and to understand the attitude, risk behavior and promiscuous sexual practice of migrant labourers. This study was mainly motivated from the data collected by us of the 60 HIV/AIDS infected persons who belonged to the category that comes to be defined as migrant labourers. Almost all of them were natives of (remote) villages and had migrated to the city, typically, "in search of jobs", or because of caste and communal violence. We have noticed how, starting from small villages with hopes and dreams these people had set out to the city, only to succumb to various temptations and finally all their dreams turned into horrid nightmares.

Our research includes probing into areas like: patterns and history of migration work, vulnerabilities and risk exposure in an alien surrounding, 'new' sexual practices or attitudes and discrimination, effect of displacement, coping mechanism etc. We also study the new economic policies of liberalization and globalization and how this has affected people to lose their traditional livelihood and sources of local employment, forcing them into migration.

Our study has been conducted among this informal sector mainly because migrant labourers are more vulnerable to HIV/AIDS infection, when compared to the local population for reasons which include easy money, poverty, powerlessness, inaccessibility to health services, unstable life-style such that insecurity, in jobs, lack of skills, alienation from hometown, lack of community bondage. Moreover, migrant labourers are also not organized into trade unions, as a result of which, they are made victims of horrendous exploitation: they are paid less than minimum wages, they don't receive legal protection, they are unaware of worker's rights issues and essentially lack stability.

Their work periods are rarely permanent; they work as short-term unskilled/ semi-skilled contract labourers or as daily wagers; but they earn well for a day and spend it badly without any social binding or savings or investing on their family members or children. A vast majority around 65% of those interviewed were essentially also part of the 'mobile' population, which was wrapped not only in a single migration from native village to metropolitan city, but also involved in jobs like driving trucks, taxis, etc. which gave them increased mobility. We have also analyzed the patients' feelings about the outreach and intervention programs related to HIV/AIDS and we have sought to comprehend the patterns of marginalization that has increased the predisposition of migrants to HIV and other infectious diseases.

A linguistic questionnaire which was drafted and interviews were conducted for 60 HIV/AIDS patients from the hospitals was the main data used in this analysis. Then the questionnaire was transformed into a Bidirectional Associative Memory (BAM) model (See Appendix 1).

Our sample group consisted of HIV infected migrant labourers whose age group ranged between 20-58 and they were involved in a variety of deregulated labour such as transport or truck drivers, construction labourers, daily wagers or employed in hotels or eateries. We have also investigated the feminization of migration and how women were vulnerable to HIV/AIDS only because of their partners.

Thus we have derived many notable conclusions and suggestion from our study of the socio-economic and psychological aspects of migrant labourers with reference to HIV/AIDS.

### DESCRIPTION OF THE PROBLEM

In view of the linguistic questionnaire we are analyzing the relation among

    a. Causes for migrants' vulnerability to HIV/AIDS
    b. Factors forcing migration
    c. Role of the Government.

We take some subtitles for each of these three main titles.

For the sake of simplicity we are restricted to some major subtitles, which has primarily interested these experts. We use BAM model on the scale [-5,5]. Here we mention that the analysis can be carried out on any other scale according to the whims and fancies of the investigator.

### A: CAUSES FOR MIGRANT LABOURERS VULNERABILITY TO HIV/AIDS

    $A_1$  -  No awareness/ no education
    $A_2$  -  Social Status
    $A_3$  -  No social responsibility and social freedom
    $A_4$  -  Bad company and addictive habits
    $A_5$  -  Types of profession
    $A_6$  -  Cheap availability of CSWs.

### F: Factors forcing people for migration

    $F_1$  -  Lack of labour opportunities in their hometown
    $F_2$  -  Poverty/seeking better status of life
    $F_3$  -  Mobilization of labour contractors

F$_4$ - Infertility of lands due to implementation of wrong research methodologies/failure of monsoon.

## G: Role of the Government

G$_1$ - Alternate job if the harvest fails there by stopping migration

G$_2$ - Awareness clubs in rural areas about HIV/AIDS

G$_3$ - Construction of hospitals in rural areas with HIV/AIDS counseling cell/ compulsory HIV/AIDS test before marriage

G$_4$ - Failed to stop the misled agricultural techniques followed recently by farmers

G$_5$ - No foresight for the government and no precautionary actions taken from the past occurrences.

## Experts opinion on the cause for vulnerability to HIV/AIDS and factors for migration

Taking the neuronal field F$_X$ as the attributes connected with the causes of vulnerability resulting in HIV/AIDS and the neuronal field F$_Y$ is taken as factors forcing migration.

The $6 \times 4$ matrix $M_1$ represents the forward synaptic projections from the neuronal field F$_X$ to the neuronal field F$_Y$.

The $4 \times 6$ matrix $M_1^T$ represents the backward synaptic projections F$_X$ to F$_Y$. Now, taking A$_1$, A$_2$, … , A$_6$ along the rows and F$_1$, …, F$_4$ along the columns we get the synaptic connection matrix $M_1$ which is modeled on the scale [-5, 5]

$$M_1 = \begin{bmatrix} 5 & 2 & 4 & 4 \\ 4 & 3 & 5 & 3 \\ -1 & -2 & 4 & 0 \\ 0 & 4 & 2 & 0 \\ 2 & 4 & 3 & 3 \\ 0 & 1 & 2 & 0 \end{bmatrix}$$

Let X$_K$ be the input vector given as (3 4 -1 -3 -2 1) at the K$^{th}$ time period. The initial vector is given such that illiteracy, lack of awareness, social status and cheap availability of CSWs have stronger impact over migration. We suppose that all neuronal state change decisions are synchronous.

The binary signal vector

$$S(X_K) \qquad = \qquad (1\ 1\ 0\ 0\ 0\ 1).$$

From the activation equation

$$S(X_K)M_1 \qquad = \qquad (9,\ 6,\ 11,\ 7)$$
$$= \qquad Y_{K+1}.$$

From the activation equation

$$S(Y_{K+1}) \qquad = \qquad (1\ 1\ 1\ 1).$$

Now

$$S(Y_{K+1})M_1^T \qquad = \qquad (15,\ 15,\ 1,\ 6,\ 12,\ 3)$$
$$= \qquad X_{K+2}.$$

From the activation equation,

$$S(X_{K+2}) \qquad = \qquad (1\ 1\ 1\ 1\ 1\ 1),$$
$$S(X_{K+2})M_1 \qquad = \qquad (10,\ 12,\ 20,\ 10)$$
$$= \qquad Y_{K+3}.$$
$$S(Y_{K+3}) \qquad = \qquad (1\ 1\ 1\ 1).$$

Thus the binary pair {(1 1 1 1 1 1), (1 1 1 1)} represents a fixed point of the dynamical system. Equilibrium of the system has occurred at the time K + 2, when the starting time was K. Thus this fixed point suggests that illiteracy with unawareness, social status and cheap availability of CSW lead to the patients remaining or becoming socially free with no social responsibility, having all addictive habits and bad company which directly depends on the types of profession they choose.

On the other hand, all the factors of migration also come to on state. Suppose we take only the on state that the availability of CSWs at very cheap rates is in the on state. Say at the $K^{th}$ time we have

$$P_K \qquad = \qquad (0\ 0\ 0\ 0\ 0\ 4),$$
$$S(P_K) \qquad = \qquad (0\ 0\ 0\ 0\ 0\ 1),$$

$$S(P_K)M_1 \quad = \quad (0\ 1\ 2\ 0)$$
$$= \quad Q_{K+1}.$$

$$S(Q_{K+1}) \quad = \quad (0\ 1\ 1\ 0)$$
$$S(Q_{K+1})M_1{}^T \quad = \quad (6\ 8\ 2\ 6\ 7\ 3)$$
$$= \quad P_{K+2}$$

$$S(P_{K+2}) \quad = \quad (1\ 1\ 1\ 1\ 1\ 1)$$
$$S(P_{K+2})M_1 \quad = \quad (10\ 12\ 20\ 10)$$
$$= \quad P_{K+3}$$
$$S(P_{K+3}) \quad = \quad (1\ 1\ 1\ 1).$$

Thus the binary pair {(1 1 1 1 1 1), (1 1 1 1)} represents a fixed point. Thus in the dynamical system given by the expert even if only the cheap availability of the CSW is in the on state, all the other states become on, i.e. they are unaware of the disease, their type of profession, they have bad company and addictive habits, they have no social responsibility and no social fear.

Thus one of the major causes for the spread of HIV/AIDS is the cheap availability of CSWs which is mathematically confirmed from our study. Several other states of vectors have been worked by us for deriving the conclusions.

## Experts Opinion On The Role Of Government And Causes For Vulnerability Of HIV/AIDS

Taking the neuronal field $F_X$ as the role of Government and the neuronal field $F_Y$ as the attributes connected with the causes of vulnerabilities resulting in HIV/ AIDS.

The $6 \times 5$ matrix $M_2$ represents the forward synaptic projections from the neuronal field $F_X$ to the neuronal field $F_Y$.

The $5 \times 6$ matrix $M_2^T$ represents the backward synaptic projection $F_X$ to $F_Y$. Now, taking $G_1, G_2, \ldots, G_5$ along the rows and $A_1, A_2, \ldots, A_6$ along the columns we get the synaptic connection matrix $M_2$ in the scale [-5,5] is as follows:

$$M_2 = \begin{bmatrix} 3 & 4 & -2 & 0 & -1 & 5 \\ 5 & 4 & 3 & -1 & 0 & 4 \\ 1 & 3 & 0 & 1 & 4 & 2 \\ 2 & 3 & -2 & -3 & 0 & 3 \\ 3 & 2 & 0 & 3 & 1 & 4 \end{bmatrix}$$

Let $X_K$ be the input vector (-3 4 -2 -1 3) at the $K^{th}$ instant. The initial vector is given such that Awareness clubs in the rural villages and the Government's inability in foreseeing the conflicts have a stronger impact over the vulnerability of HIV/AIDS. We suppose that all neuronal state change decisions are synchronous.

The binary signal vector

$$S(X_K) \quad = \quad ( 0 1 0 0 1 ).$$

From the activation equation

| | | |
|---|---|---|
| $S(X_K)M_2$ | = | ( 8 6 3 2 1 8 ) |
| | = | $Y_{K+1}$. |

Now,

| | | |
|---|---|---|
| $S(Y_{K+1})$ | = | ( 1 1 1 1 1 ), |
| $S(Y_{K+1}) M_2^T$ | = | ( 9 15 11 3 13 ) |
| | = | $X_{K+2}$. |

Now,

| | | |
|---|---|---|
| $S(X_{K+2})$ | = | (1 1 1 1 1 ), |
| $S(X_{K+2})M_2$ | = | (14 8 -1 0 4 18 ) |
| | = | $Y_{K+3}$. |

Now,

| | | |
|---|---|---|
| $S(Y_{K+3})$ | = | ( 1 1 0 1 1 1 ), |
| $S(Y_{K+3}) M_2^T$ | = | ( 11, 12, 11, 5, 13 ) |
| | = | $X_{K+4}$ |

Thus

| | | |
|---|---|---|
| $S(X_{K+4})$ | = | (1 1 1 1 1 ) |
| | = | $X_{K+2}$ |

and

| | | |
|---|---|---|
| $S(Y_{K+5})$ | = | $Y_{K+3}$. |

The binary pair ((1 1 1 1 1), (1 1 0 1 1 1)) represents a fixed point of the BAM dynamical system. Equilibrium for the system occurs at the time K + 4, when the starting time was K.

This fixed point reveals that the other three conditions cannot be ignored and have its consequences in spreading HIV/ AIDS. Similarly by taking a vector $Y_K$ one can derive conclusions based upon the nature of $Y_K$.

All such conclusions derived from this study are given in the final section.

**Experts Opinion on the Factors of Migration and the Role of Government**

Taking the neuronal field $F_X$ as the attributes connected with the factors of migration and the neuronal field $F_Y$ is the role of Government. The $4 \times 5$ matrix $M_3$ represents the forward synaptic projections from the neuronal field $F_X$ to the neuronal field $F_Y$. The $5 \times 4$ matrix $M_3^T$ represents the backward synaptic projections $F_X$ to $F_Y$. Now, taking $F_1, F_2, F_3, F_4$ along the rows and $G_1, G_2, \ldots, G_5$ along the columns we get the synaptic connection matrix $M_3$ as follows:

$$M_3 = \begin{bmatrix} 4 & 0 & 5 & 3 & 4 \\ 3 & -2 & -4 & 4 & 3 \\ 3 & 0 & 4 & -1 & -2 \\ 2 & 1 & 0 & 5 & 4 \end{bmatrix}$$

Let $X_K$ be the input vector (-2 1 4 -1) at the time K. The initial vector is given such that poverty and mobilization of labour contractors have a greater impact. We suppose that all neuronal state change decisions are synchronous.

The binary signal vector

$$S(X_K) \qquad = \qquad (0\ 1\ 1\ 0).$$

From the activation equation

$$S(X_K)\ M_3 \qquad = \qquad (\ 6\ \text{-}2\ 0\ 3\ 1\ )$$
$$= \qquad Y_{K+1}.$$

Now
$$S(Y_{K+1}\ ) \qquad = \qquad (1\ 0\ 1\ 1\ 1),$$
$$S(Y_{K+1})\ M_3^T \qquad = \qquad (16\ 6\ 4\ 11)$$
$$= \qquad X_{K+2}.$$

Now
$$S(X_{K+2}) \qquad = \qquad (1\ 1\ 1\ 1\ ),$$
$$S\ (X_{K+2})M_3 \qquad = \qquad (12\ \text{-}1\ 5\ 11\ 9)$$
$$= \qquad Y_{K+3}.$$

Now

$$S(Y_{K+3}) \quad = \quad (1\ 0\ 1\ 1\ 1),$$
$$S(Y_{K+3})M_3^T \quad = \quad (16\ 6\ 9\ 11)$$
$$= \quad X_{K+4}.$$

Thus

$$S(X_{K+4}) \quad = \quad (1\ 1\ 1\ 1)$$
$$= \quad X_{K+2},$$
$$S(Y_{K+5}) \quad = \quad Y_{K+3}.$$

The binary pair $\{((1\ 1\ 1\ 1), (1\ 0\ 1\ 1\ 1))\}$ represents a fixed point of the BAM dynamical system. Equilibrium for the system occurs at the time $K+4$, when the starting time was $K$. All the factors point out the failure of the Government in tackling HIV/AIDS. Similarly by taking a vector $Y_K$ one can derive conclusions based upon the nature of $Y_K$.

A complete conclusion about migrant labourers as HIV/AIDS patients based on the BAM model is given in the fourth chapter on conclusions.

Thus these illustrations are given only for the sake of making the reader to understand the workings of the fuzzy models. The calculations of the fixed point using in the BAM-model is given in Appendix 6 where C-program is used to make the computation easy.

Now we use the BAM model in the interval [-4 4] for analyzing the same liquistic questionnaire and we keep the main 3 heads as it is and make changes only the subtitles. These were the subtitles. These were the subtitles proposed by the expert whose opinion is sought.

**A: CAUSES FOR MIGRANT LABOURERS VULNERABILITY TO HIV/AIDS.**

| | | |
|---|---|---|
| $A_1$ | - | No awareness / No education |
| $A_2$ | - | Away from the family for weeks |
| $A_3$ | - | Social status |
| $A_4$ | - | No social responsibility and social freedom |
| $A_5$ | - | Bad company and addictive habits |
| $A_6$ | - | Types of profession where they cannot easily escape from visit of CSWs. |
| $A_7$ | - | Cheap availability of CSWs |

$A_8$ - No union / support group to chanalize their ways of living as saving for future and other monetary benefits from the employer.

$A_9$ - No fear of being watched by friends or relatives.

## F: Factors forcing people for migration

$F_1$ - Lack of labour opportunities in their hometown

$F_2$ - Poverty

$F_3$ - Unemployment

$F_4$ - Mobilization of contract labourers

$F_5$ - Infertility of land failure of agriculture

$F_6$ - Failed Government policies like advent of machinery, no value of small scale industries, weavers problem, match factory problem etc.

## G: Role of Government

$G_1$ - Lack of awareness clubs in rural areas about HIV/AIDS

$G_2$ - No steps to help agricultural cooli only rich farmers are being helped i.e., landowners alone get the benefit from the government helped

$G_3$ - Failed to stop poor yield

$G_4$ - No alternative job provided for agricultural coolie / weavers

$G_5$ - No proper counselling centres for HIV/AIDS in villages

$G_6$ - Compulsory HIV/AIDS test before marriage

$G_7$ - Government does not question about the marriage age of women in rural areas. 90% of marriage is rural areas with no education takes place when women are just below 16 years in 20% of the cases even before they girls attain puberty they are married).

**Expert's opinion on the cause of vulnerability to HIV/AIDS factors for migration**

Taking the neuronal field $F_x$ as the attributes connected with the causes of vulnerability resulting in HIV/AIDS and the neuronal field $F_y$ is taken as factors forcing migration.

The $9 \times 6$ matrix $M_1$ represents the forward synaptic projections from the neuronal field $F_x$ to the neuronal field $F_Y$. The $6 \times 9$ matrix $M^T_1$ represents the backward synaptic projections $F_X$ to $F_Y$.

Now taking $V_2, \ldots , V_9$ along the rows and $F_1, \ldots, F_6$ along the columns we get the synaptic connection matrix $M_1$ which is modelled in the scale [-4, 4].

$$M_1 = \begin{bmatrix} 2 & 0 & 0 & 0 & 0 & 0 \\ 3 & 2 & 2 & 2 & 1 & 3 \\ -2 & 3 & 0 & 2 & 0 & 0 \\ 0 & 0 & 0 & 0 & 0 & 0 \\ 0 & 0 & 0 & 1 & -2 & 0 \\ 4 & 3 & -2 & 3 & 2 & 0 \\ 0 & -1 & 0 & 0 & -2 & 0 \\ 0 & 0 & 0 & 2 & 0 & 0 \\ 0 & -3 & 0 & 0 & -2 & 0 \end{bmatrix}$$

Let $X_K$ be the input vector given as (3 2 1 –1 0 –2 4 -2 1) at the $K^{th}$ time period. The initial vectors as $V_1, V_2, \ldots, V_9$. We suppose that all neuronal state change decisions are synchronous.

The binary signal vector

$\qquad S(X_K) \qquad = \qquad$ (1 1 1 0 0 0 1 0 1)

From the activation equation

$\qquad S (X_K)M_1 \qquad = \qquad$ (3 1 2 4 –3 3)
$\qquad\qquad\qquad\qquad = \qquad Y_{K+1}$
$\qquad S(Y_{K+1}) \qquad = \qquad$ (1 1 1 1 0 1)

Now

$\qquad S (Y_{K+1}) \, M^T_1 \qquad = \qquad$ (2 1 2 3 0 1 8 –1 2 –3)
$\qquad\qquad\qquad\qquad = \qquad X_{K+2}$

$$S\ (X_{K+2}) \quad = \quad (1\ 1\ 1\ 0\ 1\ 1\ 0\ 1\ 0)$$
$$S\ (X_{K+2})\ M_1 \quad = \quad (7,\ 8,\ 0\ ,\ 10,\ 13)$$
$$= \quad Y_{K+3}$$

$$S\ (Y_{K+3}) \quad = \quad (1\ 1\ 0\ 1\ 1\ 1)$$
$$S\ (Y_{K+3})\ M^T_1 \quad = \quad (2,\ 11,\ 3,\ 0,\ -1,\ 12,\ -3\ 2-5)$$
$$= \quad X_{K+4}$$

$$S\ (X_{K+4}) \quad = \quad (1\ 1\ 1\ 0\ 0\ 1\ 0\ 1\ 0)$$
$$S\ (X_{K+4})\ M_1 \quad = \quad (7\ 8\ 0\ 9\ 3\ 3)$$
$$= \quad Y_{K+5}$$

$$S\ (Y_{K+5}) \quad = \quad (1\ 1\ 0\ 1\ 1\ 1)$$
$$S(Y_{K+5})M^T_1 \quad = \quad X_{K+6}$$
$$= \quad X_{K+4.}$$

Thus the binary pair {(1 1 1 0 0 1 0 1 0), (1 1 0 1 1 1)} represents a fixed point of the dynamical system. Equilibrium of the system has occurred at the K + 6 time when the starting point was K.

Thus the fixed point suggests that when cheap availability of CSWs, with no proper awareness and no education, with when they are away from the family with poor social status with no fear of being watched lead the patients to become victims of HIV/AIDS due to their profession and no union/support group to channelize their ways of living. The factors that fuel this are lack of labour opportunity in their hometown, poverty, mobilization of contract labourers, infertility of land failure of agriculture and above all failed government policies like advent of machinery, no value of small scale industries, weavers problem, match factory problem etc.

Using the C program given in the appendix 6 of the book we find the equilibrium of the dynamical system for varying values. The input vector on which the conclusions are made is given in chapter 7.

Suppose we consider the input vector $Y_K = (4\ –1\ 0\ –3\ 0\ 3)$ i.e. the lack of labour opportunities with failed government polices like advent of machinery which replaces human labour like weaver problem and match factor problems and no value for small scale industries to be vector with large positive values and other nodes being either 0 or negative value

| | | |
|---|---|---|
| $S(Y_K)$ | $=$ | $(1\ 0\ 0\ 0\ 0\ 1)$ |
| $S(Y_K)\ M^T_1$ | $=$ | $(2\ 6\ -2\ 0\ 0\ 4\ 0\ 0\ 0)$ |
| | $=$ | $X_{K+1}$ |
| $S(Y_{K+1})$ | $=$ | $(1\ 1\ 0\ 0\ 0\ 1\ 0\ 0\ 0)$ |
| $S(X_{K+1})\ M_1$ | $=$ | $(9\ 5\ 0\ 5\ 3\ 3)$ |
| | $=$ | $Y_{K+2}$ |
| $S(Y_{K+2})$ | $=$ | $(1\ 1\ 0\ 1\ 1\ 1)$ |
| $S(Y_{K+2})\ M^T_1$ | $=$ | $(2,\ 11,\ 3,\ 0,\ -1,\ 12,\ -3,\ 2,\ -5)$ |
| | $=$ | $X_{K+3}$ |
| $S(X_{K+3})$ | $=$ | $(1\ 1\ 1\ 0\ 0\ 1\ 0\ 1\ 0)$ |
| $S(Y_{K+2})M_1$ | $=$ | $(7\ 8\ 0\ 9\ 3\ 3)$ |
| | $=$ | $Y_{K+3}$ |
| $S(Y_{K+3})$ | $=$ | $(1\ 1\ 0\ 1\ 1\ 1)$. |

The stability of the BAM is given by the fixed point, i.e., the binary pair $\{(1\ 1\ 1\ 0\ 0\ 1\ 1\ 0\ 1\ 0),\ (1\ 1\ 0\ 1\ 1\ 1)\}$. The conclusion as before is derived so if the labourers have no employment in their home town and failed government policies result in following problem n the on state of all vectors expect. $F_3$.

Now we proceed on to model using the experts opinion on the factors of migration and the role of government.

Taking the neuronal field $F_X$ as the attributes connected with the factors of migration and the neuronal field $F_Y$ is the role of government. The matrix $M_2$ represents the forward synaptic projections from the neuronal field $F_X$ to the neuronal field $F_Y$. The $7 \times 6$ matrix $M^T_2$ represents the backward synaptic projections $F_X$ to $F_Y$. Now taking the factors $F_1, F_2,\ldots, F_6$ along the rows and $G_1, G_2,\ldots, G_7$ along the columns we get the synaptic connection matrix $M_2$ as follows:

$$M_2 = \begin{array}{c} \\ F_1 \\ F_2 \\ F_3 \\ F_4 \\ F_5 \\ F_6 \end{array} \overset{\begin{array}{ccccccc} G_1 & G_2 & G_3 & G_4 & G_5 & G_6 & G_7 \end{array}}{\begin{bmatrix} 0 & 3 & 4 & 2 & 0 & 0 & 0 \\ 0 & 2 & 3 & 3 & 2 & 0 & 3 \\ 0 & 2 & 3 & 2 & 0 & 0 & 0 \\ 0 & 3 & 2 & 3 & 2 & 0 & 0 \\ 0 & 2 & 2 & 2 & 0 & 0 & 1 \\ 0 & 2 & 1 & 2 & 0 & 0 & 1 \end{bmatrix}}$$

Here the expert feels one of the reasons for marrying off their daughter at a very young age is due to poverty. So they marry them off to old men and at times as second wife and so on.

Suppose we take the fit vector as $X_K$ = (-3 2 0 4 –1 –2) at the $K^{th}$ time period. We suppose that all neuronal state change decisions are synchronous.

The binary signal vector

$$S(X_K) \quad = \quad (0\ 1\ 0\ 1\ 0\ 0)$$

From the activation equation

$$S(X_K)\ M_2 \quad = \quad (0\ 5\ 5\ 6\ 4\ 0\ 3)$$
$$= \quad Y_{K+1}$$

$$S(Y_{K+1}) \quad = \quad (0\ 1\ 1\ 1\ 1\ 0\ 1)$$
$$S(Y_{K+1})\ M^T_2 \quad = \quad (9\ 13\ 7\ 10\ 7\ 6)$$
$$= \quad X_{K+2}$$

$$S(X_{K+2}) \quad = \quad (1\ 1\ 1\ 1\ 1\ 1)$$
$$S(X_{K+2})\ M_2 \quad = \quad (0\ 14\ 15\ 14\ 4\ 0\ 5)$$
$$= \quad Y_{K+3}$$

$$S(Y_{K+3}) \quad = \quad (0\ 1\ 1\ 1\ 1\ 0\ 1).$$

Thus the binary pair {(1 1 1 1 1 1), (0 1 1 1 1 0 1)} represents a fixed point of the dynamical system. Equilibrium of the system has occurred at the time K + 2 when the starting time was K. All attributes come to on state expect the factors Lack of awareness clubs in rural areas about HIV/AIDS and Compulsory HIV/AIDS lest before marriage are just zero, for they have no impact.

While the factor poverty and mobilization of the contract labourer are at the dominance. Let $Y_K$ = (3, -4, -2, -1 0 –1 2) be the input vector at the $K^{th}$ instant. The initial vector is given such that. Next we consider when in the state vector Lack of awareness clubs in rural areas about HIV/AIDS is given the maximum priority followed by marriage age of women is not questioned by the government.

We study the impact of this on migration.

We suppose that all neuronal state change decisions are synchronous.

The binary signal vector.

$$S(Y_K) \qquad = \qquad (1\ 0\ 0\ 0\ 0\ 0\ 1)$$

From the activation equation

$$S(Y_K)\,M_2^T \qquad = \qquad (0\ 3\ 0\ 0\ 1\ 1)$$
$$\qquad\qquad\qquad = \qquad X_{K+1}$$

$$S(X_{K+1}) \qquad = \qquad (0\ 1\ 0\ 0\ 1\ 1)$$
$$S(X_{K+1})\,M_2 \qquad = \qquad (0\ 6\ 6\ 7\ 2\ 0\ 5)$$
$$\qquad\qquad\qquad = \qquad Y_{K+2}$$

$$S(Y_{K+2}) \qquad = \qquad (0\ 1\ 1\ 1\ 1\ 0\ 1)$$
$$S(Y_{K+2})\,M_2^T \qquad = \qquad (9\ 13\ 7\ 10\ 7\ 6)$$
$$\qquad\qquad\qquad = \qquad Y_{K+3}$$

$$S(X_{K+3}) \qquad = \qquad (1\ 1\ 1\ 1\ 1\ 1)$$
$$S(X_{K+3})\,M_2 \qquad = \qquad (0\ 14\ 15\ 14\ 4\ 0\ 5)$$
$$\qquad\qquad\qquad = \qquad Y_{K+4}$$

$$S(Y_{K+4}) \qquad = \qquad (0\ 1\ 1\ 1\ 1\ 0\ 1).$$

Thus the binary pair {(0 1 1 1 1 0 1), (1 1 1 1 1 1)} represents a fixed point of dynamical system. Equilibrium for the system occurs at the time K + 3 when the starting time waste. All nodes except government failure to perform HIV/AIDS test as a compulsory one alone is unaffected and all factors $F_1$, …, $F_6$ become on so these two nodes promote migration in all forms.

Next we study the experts opinion on the role of government and the causes for vulnerability of HIV/AIDs.

Now taking the neuronal field $F_X$ as the role of government. and the neuronal field $F_Y$ as the attributes connected with the causes of vulnerabilities resulting in HIV/AIDS.

The 7 × 9 matrix $M_3$ represents the forward synaptic projections from the neuronal field $F_X$ to the neuronal field $F_Y$ and the 9 × 7 matrix $M_3^T$ represents the backward synaptic projection $F_X$ to $F_Y$.

Now taking $G_1$,…, $G_7$ along the rows and $A_1$,…, $A_9$ along the columns we get the synaptic connection matrix $M_3$ in the scale [-4 4] as follows:

$$M_3 = \begin{bmatrix} 3 & 0 & -2 & 0 & 0 & 4 & 3 & 4 & 2 \\ -1 & 0 & 3 & 0 & 0 & 0 & 2 & 0 & 0 \\ 0 & 3 & 1 & 0 & 1 & 0 & 0 & 0 & 0 \\ -1 & 2 & 3 & 0 & 0 & 0 & 0 & 1 & 0 \\ 3 & 0 & 1 & 0 & 0 & 1 & 1 & 0 & 0 \\ 0 & 0 & 0 & -1 & -1 & 0 & 0 & -1 & -1 \\ 3 & 0 & 1 & 0 & 0 & 0 & 0 & 0 & 0 \end{bmatrix}$$

Let $X_K$ be the input vector (2 –1 +3 –2 –1 0 –1) at the $K^{th}$ instant. The initial vector is given that Failed to stop poor yield and Lack of awareness clubs in rural area for the uneducated about HIV/AIDS, takes up the primary position. We shall assume all neuronal state change decisions are synchronous.

The binary signal vector

$$S(X_K) \qquad = \qquad (1\ 0\ 1\ 0\ 0\ 0\ 0)$$

From the activation equation

$$\begin{aligned}
S(X_K)\,M_3 &= (3, 3\ –1\ 0\ 1\ 4\ 3\ 4\ 2) \\
&= Y_{K+1} \\
S(Y_{K+1}) &= (1\ 1\ 0\ 0\ 1\ 1\ 1\ 1\ 1) \\
S(Y_{K+1})\,M^T_3 &= (16, 1, 4, 2, 5\text{-}3\ 3) \\
&= X_{K+2}
\end{aligned}$$

$$\begin{aligned}
S(X_{K+2}) &= (1\ 1\ 1\ 1\ 1\ 0\ 1) \\
S(X_{K+2})\,M_3 &= (7, 5, 7\ 0\ 1\ 5\ 6\ 5\ 2) \\
&= Y_{K+3}
\end{aligned}$$

$$\begin{aligned}
S(Y_{K+3}) &= (1\ 1\ 1\ 0\ 1\ 1\ 1\ 1\ 1) \\
S(Y_{K+3})\,M_3 &= (5\ 3\ 4\ –1\ 0\ 5\ 6\ 3\ 1) \\
&= X_{K+4}
\end{aligned}$$

$$\begin{aligned}
S(X_{K+4}) &= (1\ 1\ 1\ 0\ 0\ 1\ 1\ 1\ 1) \\
S(X_{K+4})\,M^T_3 &= (14\ 4\ 4\ 5\ 6\ –2\ 4) \\
&= Y_{K+5}
\end{aligned}$$

$$S(Y_{K+5}) = (1\ 1\ 1\ 1\ 1\ 0\ 1).$$

Thus the binary pair {(1 1 1 1 1 0 1), (1 1 10 1 1 1 1 1)} or {(1 1 1 1 1 0 1), (1 1 1 0 0 1 1 1 1)} where (1 1 1 1 1 0 1) is a fixed point and (1 1 1 0 1 1 1 1 1) or (1 1 1 0 0 1 1 1 1) occurs.

No social responsibility on the part of the part of the HIV/AIDS affected migrant labourers is 0 all other co-ordinates becomes on their by proving all other cases enlisted makes these labourers vulnerable to HIV/AIDS or the bad habits/bad company and no social responsibly on the parts of HIV/AIDS affected is 0, all other causes for vulnerability is 'on': all factors relating government becomes on except the node $(G_6)$ – compulsory test for HIV/AIDS before marriage is off. Thus all factors $G_1$, $G_2$, …, $G_5$ and $G_7$ become on.

Now we consider the input vector $Y_K$ = (3 2 –1 –2 0 –3 1 0 –2) at the $K^{th}$ instant. The initial vector is given such that No awareness no education takes a priority followed by away from the family for weeks and cheap availability of CSWs.

The binary signal vector

$$S\,(Y_K) \qquad = \qquad (1\ 1\ 0\ 0\ 0\ 0\ 1\ 0\ 0).$$

From the activation equation

$$S(Y_K)M^T_3 \quad = \quad (6\ 1\ 3\ 1\ 4\ 0\ 3)$$
$$\qquad\quad = \quad X_{K+1}$$
$$S\,(X_{K+1}) \quad = \quad (1\ 1\ 1\ 1\ 1\ 0\ 1)$$
$$S\,(X_{K+1})\,M_3 \quad = \quad (7\ 5\ 7\ 0\ 1\ 5\ 6\ 5\ 2)$$
$$\qquad\quad = \quad Y_{K+2}$$

$$S\,(Y_{K+2}) \quad = \quad (1\ 1\ 1\ 0\ 1\ 1\ 1\ 1\ 1)$$
$$S\,(Y_{K+2})\,M^T_3 \quad = \quad (14,\,4,\,5,\,5\ 6\ –3\ 4)$$
$$\qquad\quad = \quad X_{K+3}$$

$$S\,(X_{K+3}) \quad = \quad (1\ 1\ 1\ 1\ 1\ 0\ 1)$$
$$\qquad\quad = \quad S\,(X_{K+1}).$$

Thus the binary pair {(1 1 1 1 1 0 1), (1 1 1 0 1 1 1 1 1)} is a fixed point of the dynamical system. All factors become on in both the classes of attributes except $G_6$ and $A_4$ i.e., compulsory test for HIV/AIDS remains always off and no social responsibility and social freedom on the part of migrants in no other cause for vulnerability to HIV/AIDS.

Now we obtain the equation. The effect of government influence and the vulnerability of migrant labour using the product of the matrices $M_1$ and $M_2$ and taking the transpose.

$$\left(M_1 \times M_2\right)^T = \begin{bmatrix} 0 & 6 & 8 & 4 & 0 & 0 & 0 \\ 0 & 31 & 33 & 33 & 8 & 0 & 10 \\ 0 & 6 & 8 & 11 & 10 & 0 & 9 \\ 0 & 0 & 0 & 0 & 0 & 0 & 0 \\ 0 & -1 & -2 & -1 & 0 & 0 & -2 \\ 0 & 27 & 29 & 26 & 12 & 0 & 11 \\ 0 & -6 & -7 & -7 & -2 & 0 & -5 \\ 0 & 6 & 4 & 6 & 4 & 0 & 0 \\ 0 & -10 & -13 & -13 & -6 & 0 & -13 \end{bmatrix}$$

We study this in the interval [-14, 14]. This gives the indirect relationship between government policies and the vulnerability of migrant labourers to HIV/AIDS.

The reader is expected to analyse this BAM model using C-program in the appendix 6 and draw conclusions.

Chapter Four

# USE OF NCM TO STUDY THE SOCIO-ECONOMIC PROBLEMS OF HIV/AIDS AFFECTED MIGRANT LABOURERS

It is very unusual in India to see migrant labourers to be the highest HIV/AIDS affected population among the HIV/AIDS affected persons in rural areas. The major reason being we see over 90% of the HIV/AIDS patient are uneducated, economically very poor, and lower middle class from rural areas who take treatment in Tambaram Sanatorium. When badly affected by HIV/AIDS they have no other go either die in the same village or go for treatment in Tambaram Sanatorium hospital. Also it is pertinent to mention that at least 70% of them are ignorant of the mode of spread of disease for they said they felt by visiting CSWs they would be mainly affected by sexually transmitted diseases or venereal diseases. They said they did not know they would become HIV/AIDS victims. So only of the 60 persons 45 of them said after they visited CSWs they took injections and medicines to save themselves from VD/STD. So one of the conjectures which we can derive from their statement is that they would have been doubly careful if they would have been informed about the HIV/AIDS affliction by visiting CSWs.

It is also important to note that most of the HIV/AIDS affected migrant labourers when stay in the place away from their home town are not in a position to channelize their energy in constructive things on the other hand they go for sex which is their only recreation. Also these people do not have association or unions. Thus with no external support they become easy victims of CSWs. Now we just briefly recall the definition of Neutrosophic Cognitive Maps and the notion of neutrosophy.

This chapter introduces the concept of Neutrosophic Cognitive Maps (NCMs). NCMs are a generalization of Fuzzy Cognitive Maps (FCMs). To study or even to define Neutrosophic Cognitive Maps we have to basically define the notions of Neutrosophic Matrix. Further analogous to Fuzzy Relational Maps (FRMs) we define the notion of Neutrosophic Relational Maps, which are clearly a generalization of FRMs. Study of this type is for the first time experimented in this book. When the data under analysis has indeterminate concepts we are not in a position to give it a mathematical expression. Because of the recent introduction of neutrosophic logic by Florentin Smarandache this problem has a solution. Thus we have introduced the additional notion of Neutrosophy in place of Fuzzy theory.

This chapter has four sections. In the first section we just give a very brief introduction to neutrosophy. In section two, we, for the first time, use NCM to analyze the HIV/AIDS affected migrant labourers problem. Section three is the study of the migrants problem using Combined NCM. Section four defines the notion of Combined Disjoint Block NCM and Combined Overlap Block NCM and we apply these two models to this problem.

## 4.1 Basic notions of Neutrosophy and Neutrosophic cognitive maps

In this section we introduce the notion of neutrosophic logic created by Florentin Smarandache, which is an extension/combination of the fuzzy logic in which indeterminacy is included. It has become very essential that the notion of neutrosophic logic play a vital role in several of the real world problems like law, medicine, industry, finance, IT, stocks and share etc. Use of neutrosophic notions will be illustrated/ applied in the later sections of this chapter. Fuzzy theory only measures the grade of membership or the non-existence of a membership in the revolutionary way but fuzzy theory has failed to attribute the concept when the relations between notions or nodes or concepts in problems are indeterminate. In fact one can say the inclusion of the concept of indeterminate situation with fuzzy concepts will form the neutrosophic logic.

As in this book, application of the concept of only fuzzy cognitive maps are dealt which mainly deals with the relation / non-relation between two nodes or concepts but it fails to deal the relation between two conceptual nodes when the relation is an indeterminate one. Neutrosophic logic is the only tool known to us, which deals with the notions of indeterminacy, and here we give a brief description of it. For more about Neutrosophic logic please refer Smarandache.

**DEFINITION 4.1.1:** *In the neutrosophic logic every logical variable x is described by an ordered triple x = (T, I, F) where T is the degree of truth, F is the degree of false and I the level of indeterminacy.*

(A). To maintain consistency with the classical and fuzzy logics and with probability there is the special case where
$$T + I + F = 1.$$
(B). But to refer to intuitionistic logic, which means incomplete information on a variable proposition or event one has
$$T + I + F < 1.$$
(C). Analogically referring to Paraconsistent logic, which means contradictory sources of information about a same logical variable, proposition or event one has
$$T + I + F > 1.$$

Thus the advantage of using Neutrosophic logic is that this logic distinguishes between relative truth that is a truth is one or a few worlds only noted by 1 and absolute truth denoted by $1^+$. Likewise neutrosophic logic distinguishes between relative falsehood, noted by 0 and absolute falsehood noted by $^-0$.

It has several applications. One such given by is as follows:

***Example 4.1.1:*** From a pool of refugees, waiting in a political refugee camp in Turkey to get the American visa, a% have the chance to be accepted – where a varies in the set A, r% to be rejected – where r varies in the set R, and p% to be in pending (not yet decided) – where p varies in P.

Say, for example, that the chance of someone Popescu in the pool to emigrate to USA is (between) 40-60% (considering different criteria of emigration one gets different percentages, we have to take care of all of them), the chance of being rejected is 20-25% or 30-35%, and the chance of being in pending is 10% or 20% or 30%. Then the neutrosophic probability that Popescu emigrates to the Unites States is

*NP (Popescu) = ((40-60) (20-25) $\cup$ (30-35), {10,20,30}), closer to the life.*

This is a better approach than the classical probability, where 40 P(Popescu) 60, because from the pending chance – which will be converted to acceptance or rejection—Popescu might get extra percentage in his will to emigrating and also the superior limit of the subsets sum $60 + 35 + 30 > 100$ and in other cases one may have the inferior sum $< 0$, while in the classical fuzzy set theory the superior sum should be 100 and the inferior sum $\mu$ 0. In a similar way, we could say about the element Popescu that Popescu ((40-60), (20-25) $\cup$ (30-35), {10, 20, 30}) belongs to the set of accepted refugees.

*Example 4.1.2:* The probability that candidate C will win an election is say 25-30% true (percent of people voting for him), 35% false (percent of people voting against him), and 40% or 41% indeterminate (percent of people not coming to the ballot box, or giving a blank vote – not selecting any one or giving a negative vote cutting all candidate on the list). Dialectic and dualism don't work in this case anymore.

*Example 4.1.3:* Another example, the probability that tomorrow it will rain is say 50-54% true according to meteorologists who have investigated the past years weather, 30 or 34-35% false according to today's very sunny and droughty summer, and 10 or 20% undecided (indeterminate).

*Example 4.1.4:* The probability that Yankees will win tomorrow versus Cowboys is 60% true (according to their confrontation's history giving Yankees' satisfaction), 30-32% false (supposing Cowboys are actually up to the mark, while Yankees are declining), and 10 or 11 or 12% indeterminate (left to the hazard: sickness of players, referee's mistakes, atmospheric conditions during the game). These parameters act on players' psychology.

As in this book we use mainly the notion of neutrosophic logic with regard to the indeterminacy of any relation in cognitive maps we are restraining ourselves from dealing with several interesting concepts about neutrosophic logic. As NCMs deals with unsupervised data and the existence or non-existence of cognitive relation, we do not in this book elaborately describe the notion of neutrosophic concepts. See Appendix 9 for more on Neutrosophy. However we just state, suppose in a legal issue the jury or the judge cannot always prove the evidence in a case, in several places we may not be able to derive any conclusions from the existing facts because of which we cannot make a conclusion that no relation exists or otherwise. But existing relation is an indeterminate. So in the case when the concept of indeterminacy exists the judgment ought to be very carefully analyzed be it a civil case or a criminal case. FCMs are deployed only where the existence or non-existence is dealt with but however in our Neutrosophic Cognitive Maps we will deal with the notion of indeterminacy of the evidence also. Thus legal side has lot of Neutrosophic (NCM) applications. Also we will show how NCMs can be used to study factors as varied as stock markets, medical diagnosis, etc. "$I$" denotes the indeterminacy of notion/concept/relation. In this section we for the first time introduce the notion of neutrosophic graphs, illustrate them and give some basic properties. We need the notion of neutrosophic graphs basically to obtain neutrosophic cognitive maps which will be nothing but directed neutrosophic graphs. Similarly neutrosophic relational maps will also be directed neutrosophic graphs. It is no coincidence that graph theory has been independently discovered many times since it may quite properly be regarded as an area of applied mathematics. The subject finds its place in the work of Euler. Subsequent rediscoveries of graph theory were by Kirchhoff and Cayley. Euler (1707-1782) became the father of graph theory as well as topology when in 1936 he settled a famous unsolved problem in his day called the Konigsberg Bridge Problem.

Psychologist Lewin proposed in 1936 that the life space of an individual be represented by a planar map. His view point led the psychologists at the Research center for Group Dynamics to another psychological interpretation of a graph in which people are represented by points and interpersonal relations by lines. Such relations include love, hate, communication and power. In fact it was precisely this approach which led the author to a personal discovery of graph theory, aided and abetted by psychologists L. Festinger and D. Cartwright. Here it is pertinent to mention that the directed graphs of an FCMs or FRMs are nothing but the psychological inter-relations or feelings of different nodes; but it is unfortunate that in all these studies the concept of indeterminacy was never given any place, so in this chapter for the first time we will be having graphs in which the notion of indeterminacy i.e. when two vertex should be connected or not is never dealt with. If graphs are to display human feelings then this point is very vital for in several situations certain relations between concepts may certainly remain an indeterminate. So this section will purely cater to the properties of such graphs the edges of certain vertices may not find its connection i.e., they are indeterminates, which we will be defining as neutrosophic graphs.

The world of theoretical physics discovered graph theory for its own purposes. In the study of statistical mechanics by Uhlenbeck the points stands for molecules and two adjacent points indicate nearest neighbor interaction of some physical kind, for example magnetic attraction or repulsion. But it is forgotten in all these situations we may have molecule structures which need not attract or repel but remain without action or not able to predict the action for such analysis we can certainly adopt the concept of neutrosophic graphs. In a similar interpretation by Lee and Yang the points stand for small cubes in Euclidean space where each cube may or may not be occupied by a molecule. Then two points are adjacent whenever both spaces are occupied.

Feynmann proposed the diagram in which the points represent physical particles and the lines represent paths of the particles after collisions. Just at each stage of applying graph theory we may now feel the neutrosophic graph theory may be more suitable for application.

Now we proceed on to define the neutrosophic graph.

**DEFINITION 4.1.2:** *A neutrosophic graph is a graph in which at least one edge is an indeterminacy denoted by dotted lines.*

**NOTATION***: The indeterminacy of an edge between two vertices will always be denoted by dotted lines.*

***Example 4.1.1:*** The following are neutrosophic graphs:

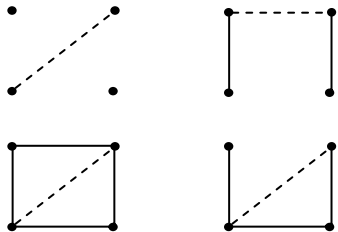

FIGURE: 4.1.1

All graphs in general are not neutrosophic graphs.

***Example 4.1.2:*** The following graphs are not neutrosophic graphs given in Figure 4.1.2:

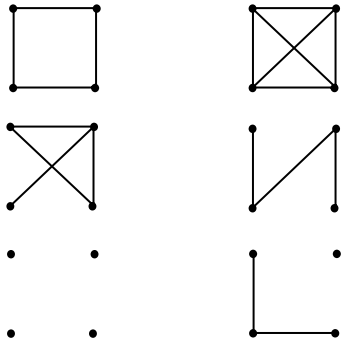

FIGURE: 4.1.2

**DEFINITION 4.1.3:** *A neutrosophic directed graph is a directed graph which has at least one edge to be an indeterminacy.*

**DEFINITION 4.1.4:** *A neutrosophic oriented graph is a neutrosophic directed graph having no symmetric pair of directed indeterminacy lines.*

**DEFINITION 4.1.5:** *A neutrosophic subgraph H of a neutrosophic graph G is a subgraph H which is itself a neutrosophic graph.*

**THEOREM 4.1.1:** *Let G be a neutrosophic graph. All subgraphs of G are not neutrosophic subgraphs of G.*

*Proof:* By an example. Consider the neutrosophic graph given in Figure 4.1.3.

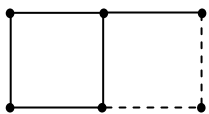

FIGURE: 4.1.3

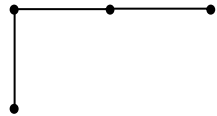

FIGURE: 4.1.4

This has a subgraph given by Figure 4.1.4. which is not a neutrosophic subgraph of G.

**THEOREM 4.1.2:** *Let G be a neutrosophic graph. In general the removal of a point from G need be a neutrosophic subgraph.*

*Proof:* Consider the graph G given in Figure 4.1.5.

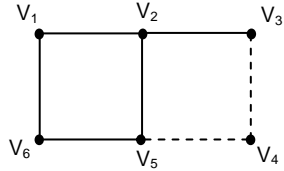

FIGURE: 4.1.5

G \ $v_4$ is only a subgraph of G but is not a neutrosophic subgraph of G. Thus it is interesting to note that this a main feature by which a graph differs from a neutrosophic graph.

**DEFINITION 4.1.6:** *Two graphs G and H are neutrosophically isomorphic if*

   i.    *They are isomorphic*
  ii.    *If there exists a one to one correspondence between their point sets which preserve indeterminacy adjacency.*

**DEFINITION 4.1.7:** *A neutrosophic walk of a neutrosophic graph G is a walk of the graph G in which at least one of the lines is an indeterminacy line. The neutrosophic walk is neutrosophic closed if $v_0 = v_n$ and is neutrosophic open otherwise. It is a neutrosophic trial if all the lines are distinct and at least one of the lines is a indeterminacy line and a path if all points are distinct (i.e. this necessarily means all lines are distinct and at least one line is a line of indeterminacy). If the neutrosophic walk is neutrosophic closed then it is a neutrosophic cycle provided its n points are distinct and $n \geq 3$.*

*A neutrosophic graph is neutrosophic connected if it is connected and at least a pair of points are joined by a path. A neutrosophic maximal connected neutrosophic subgraph of G is called a neutrosophic connected component or simple neutrosophic component of G. Thus a neutrosophic graph has at least two neutrosophic components then it is neutrosophic disconnected. Even if one is a component and another is a neutrosophic component still we do not say the graph is neutrosophic disconnected.*

Neutrosophic disconnected graphs are given in Figure 4.1.6.

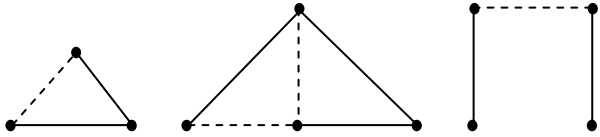

FIGURE: 4.1.6

Graph which is not neutrosophic disconnected is given by Figure 4.1.7.

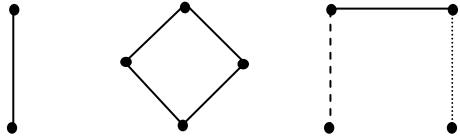

FIGURE: 4.1.7

Several results in this direction can be defined and analysed.

**DEFINITION 4.1.8:** *A neutrosophic bigraph, G is a bigraph, G whose point set V can be partitioned into two subsets $V_1$ and $V_2$ such that at least a line of G which joins $V_1$ with $V_2$ is a line of indeterminacy.*

*This neutrosophic bigraphs will certainly play a role in the study of FRMs and FCMs and in fact give a method of conversion of data from FRMs to FCMs.*

As both the models FRMs and FCMs work on the adjacency or the connection matrix we just define the neutrosophic adjacency matrix related to a neutrosophic graph G given by Figure 4.1.8.

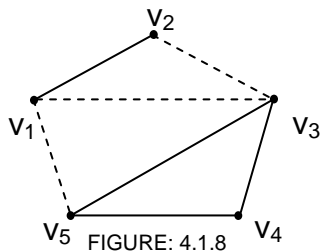

FIGURE: 4.1.8

*The neutrosophic adjacency matrix is N(A)*

$$N(A) = \begin{bmatrix} 0 & 1 & I & 0 & I \\ I & 0 & I & 0 & 0 \\ I & I & 0 & 1 & 1 \\ 0 & 0 & 1 & 0 & 1 \\ I & 0 & 1 & 1 & 0 \end{bmatrix}.$$

*Its entries will not only be 0 and 1 but also the indeterminacy I.*

**DEFINITION 4.1.9:** *Let G be a neutrosophic graph. The adjacency matrix of G with entries from the set (I, 0, 1) is called the neutrosophic adjacency matrix of the graph.*

Now as our main aim is the study of Neutrosophic Cognitive Maps we do not divulge into a very deep study of Neutrosophic Graphs or its properties but have given only the basic and the appropriate notions which are essential for studying of this book.

The notion of Fuzzy Cognitive Maps (FCMs) which are fuzzy signed directed graphs with feedback are discussed and described in Chapter 1. The directed edge $e_{ij}$ from causal concept $C_i$ to concept $C_j$ measures how much $C_i$ causes $C_j$. The time varying concept function $C_i(t)$ measures the non negative occurrence of some fuzzy event, perhaps the strength of a political sentiment, historical trend or opinion about some topics like child labor or school dropouts etc. FCMs model the world as a collection of classes and causal relations between them.

The edge $e_{ij}$ takes values in the fuzzy causal interval [–1, 1] ($e_{ij} = 0$ indicates no causality, $e_{ij} > 0$ indicates causal increase; that $C_j$ increases as $C_i$ increases and $C_j$ decreases as $C_i$ decreases, $e_{ij} < 0$ indicates causal decrease or negative causality $C_j$ decreases as $C_i$ increases or $C_j$, increases as $C_i$ decreases. Simple FCMs have edge value in {-1, 0, 1}. Thus if causality occurs it occurs to maximal positive or negative degree.

It is important to note that $e_{ij}$ measures only absence or presence of influence of the node $C_i$ on $C_j$ but till now any researcher has not contemplated the indeterminacy of any relation between two nodes $C_i$ and $C_j$. When we deal with unsupervised data, there are situations when no relation can be determined between some two nodes. So in this section we try to introduce the indeterminacy in FCMs, and we choose to call this generalized structure as Neutrosophic Cognitive Maps (NCMs). In our view this will certainly give a more appropriate result and also caution the user about the risk of indeterminacy. Now we proceed on to define the concepts about NCMs.

**DEFINITION 4.1.10:** *A Neutrosophic Cognitive Map (NCM) is a neutrosophic directed graph with concepts like policies, events etc. as nodes and causalities or indeterminates as edges. It represents the causal relationship between concepts.*

Let $C_1, C_2, \ldots, C_n$ denote n nodes, further we assume each node is a neutrosophic vector from neutrosophic vector space V. So a node $C_i$ will be represented by $(x_1, \ldots, x_n)$ where $x_k$'s are zero or one or $I$ ($I$ is the indeterminate introduced in Sections 2.2 and 2.3 of the chapter 2) and $x_k = 1$ means that the node $C_k$ is in the on state and $x_k = 0$ means the node is in the off state and $x_k = I$ means the nodes state is an indeterminate at that time or in that situation.

Let $C_i$ and $C_j$ denote the two nodes of the NCM. The directed edge from $C_i$ to $C_j$ denotes the causality of $C_i$ on $C_j$ called connections. Every edge in the NCM is weighted with a number in the set $\{-1, 0, 1, I\}$.

Let $e_{ij}$ be the weight of the directed edge $C_iC_j$, $e_{ij} \in \{-1, 0, 1, I\}$. $e_{ij} = 0$ if $C_i$ does not have any effect on $C_j$, $e_{ij} = 1$ if increase (or decrease) in $C_i$ causes increase (or decreases) in $C_j$, $e_{ij} = -1$ if increase (or decrease) in $C_i$ causes decrease (or increase) in $C_j$. $e_{ij} = I$ if the relation or effect of $C_i$ on $C_j$ is an indeterminate.

**DEFINITION 4.1.11:** *NCMs with edge weight from {-1, 0, 1, I} are called simple NCMs.*

**DEFINITION 4.1.12:** *Let $C_1, C_2, \ldots, C_n$ be nodes of a NCM. Let the neutrosophic matrix N(E) be defined as $N(E) = (e_{ij})$ where $e_{ij}$ is the weight of the directed edge $C_i C_j$, where $e_{ij} \in \{0, 1, -1, I\}$. N(E) is called the neutrosophic adjacency matrix of the NCM.*

**DEFINITION 4.1.13:** *Let $C_1, C_2, \ldots, C_n$ be the nodes of the NCM. Let $A = (a_1, a_2, \ldots, a_n)$ where $a_i \in \{0, 1, I\}$. A is called the instantaneous state neutrosophic vector and it denotes the on-off – indeterminate state position of the node at an instant*

$$a_i \;=\; 0 \text{ if } a_i \text{ is off (no effect)}$$
$$a_i \;=\; 1 \text{ if } a_i \text{ is on (has effect)}$$
$$a_i \;=\; I \text{ if } a_i \text{ is indeterminate(effect cannot be determined)}$$

*for $i = 1, 2, \ldots, n$.*

**DEFINITION 4.1.14:** *Let $C_1$, $C_2$, …, $C_n$ be the nodes of the FCM. Let $\overrightarrow{C_1C_2}$, $\overrightarrow{C_2C_3}$, $\overrightarrow{C_3C_4}$, … , $\overrightarrow{C_iC_j}$ be the edges of the NCM. Then the edges form a directed cycle. An NCM is said to be cyclic if it possesses a directed cyclic. An NCM is said to be acyclic if it does not possess any directed cycle.*

**DEFINITION 4.1.15:** *An NCM with cycles is said to have a feedback. When there is a feedback in the NCM i.e. when the causal relations flow through a cycle in a revolutionary manner the NCM is called a dynamical system.*

**DEFINITION 4.1.16:** *Let $\overrightarrow{C_1C_2}$, $\overrightarrow{C_2C_3}$, $\cdots$, $\overrightarrow{C_{n-1}C_n}$ be cycle, when $C_i$ is switched on and if the causality flow through the edges of a cycle and if it again causes $C_i$, we say that the dynamical system goes round and round. This is true for any node $C_i$, for i = 1, 2,…, n. the equilibrium state for this dynamical system is called the hidden pattern.*

**DEFINITION 4.1.17:** *If the equilibrium state of a dynamical system is a unique state vector, then it is called a fixed point. Consider the NCM with $C_1$, $C_2$,…, $C_n$ as nodes. For example let us start the dynamical system by switching on $C_1$. Let us assume that the NCM settles down with $C_1$ and $C_n$ on, i.e. the state vector remain as (1, 0,…, 1) this neutrosophic state vector (1,0,…, 0, 1) is called the fixed point.*

**DEFINITION 4.1.18:** *If the NCM settles with a neutrosophic state vector repeating in the form*

$$A_1 \rightarrow A_2 \rightarrow … \rightarrow A_i \rightarrow A_1,$$

*then this equilibrium is called a limit cycle of the NCM.*

**METHODS OF DETERMINING THE HIDDEN PATTERN:**

Let $C_1$, $C_2$,…, $C_n$ be the nodes of an NCM, with feedback. Let E be the associated adjacency matrix. Let us find the hidden pattern when $C_1$ is switched on when an input is given as the vector $A_1 = (1, 0, 0,…, 0)$, the data should pass through the neutrosophic matrix N(E), this is done by multiplying $A_1$ by the matrix N(E). Let $A_1N(E) = (a_1, a_2,…, a_n)$ with the threshold operation that is by replacing $a_i$ by 1 if $a_i > k$ and $a_i$ by 0 if $a_i < k$ (k – a suitable positive integer) and $a_i$ by $I$ if $a_i$ is not a integer. We update the resulting concept, the concept $C_1$ is included in the updated vector by making the first coordinate as 1 in the resulting vector. Suppose $A_1N(E) \rightarrow A_2$ then consider $A_2N(E)$ and repeat the same procedure. This procedure is repeated till we get a limit cycle or a fixed point.

**DEFINITION 4.1.19:** *Finite number of NCMs can be combined together to produce the joint effect of all NCMs. If $N(E_1)$, $N(E_2)$,..., $N(E_p)$ be the neutrosophic adjacency matrices of a NCM with nodes $C_1$, $C_2$,..., $C_n$ then the combined NCM is got by adding all the neutrosophic adjacency matrices $N(E_1)$,..., $N(E_p)$. We denote the combined NCMs adjacency neutrosophic matrix by $N(E) = N(E_1) + N(E_2)+...+ N(E_p)$.*

**NOTATION:** Let $(a_1, a_2, \ldots , a_n)$ and $(a'_1, a'_2, \ldots , a'_n)$ be two neutrosophic vectors. We say $(a_1, a_2, \ldots , a_n)$ is equivalent to $(a'_1, a'_2, \ldots , a'_n)$ denoted by $((a_1, a_2, \ldots , a_n) \sim (a'_1, a'_2, \ldots, a'_n)$ if $(a'_1, a'_2, \ldots , a'_n)$ is got after thresholding and updating the vector $(a_1, a_2, \ldots , a_n)$ after passing through the neutrosophic adjacency matrix N(E).

The following are very important:

***Note 1:*** The nodes $C_1$, $C_2$, …, $C_n$ are nodes are not indeterminate nodes because they indicate the concepts which are well known. But the edges connecting $C_i$ and $C_j$ may be indeterminate i.e. an expert may not be in a position to say that $C_i$ has some causality on $C_j$ either will he be in a position to state that $C_i$ has no relation with $C_j$ in such cases the relation between $C_i$ and $C_j$ which is indeterminate is denoted by $I$.

***Note 2:*** The nodes when sent will have only ones and zeros i.e. on and off states, but after the state vector passes through the neutrosophic adjacency matrix the resultant vector will have entries from $\{0, 1, I\}$ i.e. they can be neutrosophic vectors.

The presence of $I$ in any of the coordinate implies the expert cannot say the presence of that node i.e. on state of it after passing through N(E) nor can we say the absence of the node i.e. off state of it the effect on the node after passing through the dynamical system is indeterminate so only it is represented by $I$. Thus only in case of NCMs we can say the effect of any node on other nodes can also be indeterminates. Such possibilities and analysis is totally absent in the case of FCMs.

***Note 3:*** In the neutrosophic matrix N(E), the presence of *I* in the $a_{ij}$ the place shows, that the causality between the two nodes i.e. the effect of $C_i$ on $C_j$ is indeterminate. Such chances of being indeterminate is very possible in case of unsupervised data and that too in the study of FCMs which are derived from the directed graphs.

Thus only NCMs helps in such analysis. Now we shall represent a few examples to show how in this set up NCMs is preferred to FCMs. At the outset before we proceed to give examples we make it clear that all unsupervised data need not have NCMs to be applied to it. Only data which have the relation between two nodes to be an indeterminate need to be modeled with NCMs if the data has no indeterminacy factor between any pair of nodes one need not go for NCMs; FCMs will do the best job.

## 4.2 Use of NCM to analyze the HIV/AIDS affected migrant labourer problems

We just use the basic concept of NCM in our analysis of psychological, socio and economic problems of HIV/AIDS affected migrant labourers. To make the reader compare and comprehend the effect of NCM in the place of FCM, we first analyse using the experts opinion giving them a possibility namely it is not essential to relate or not relate two concepts / nodes and they can also use the symbol *I* to denote that they are not able to determine any form of relation between a pair of nodes i.e., at that time the relation between the nodes is an indeterminate to them.

Now we first study the twelve heads given in page 42-52 and analyzed using FCM.

$A_1$ - Easy Availability of money
$A_2$ - Lack of education
$A_3$ - Visiting CSWs
$A_4$ - Profession
$A_5$ - Wrong / Bad company
$A_6$ - Addiction to habit-forming substances-
visiting CSWs
$A_7$ - Absence of social responsibility

A$_8$ - Socially free
A$_9$ - Economic status
A$_{10}$ - More leisure
A$_{11}$ - Machismo / exaggerated masculinity
A$_{12}$ - No awareness of HIV/AIDS disease.

Now the same expert who spelt out there attributes was now asked to give the opinion suggesting that he has the option to state his inability to associate two concept of nodes such things will be denoted by the letter *I* meaning the relation is an indeterminate.

For the description of each of the 12 attributes please refer page 42-52 of Chapter 2. The neutrosophic graph of the NCM.

Thus we have obtained the experts opinion, which is given the neutrosophic directed graph representation.

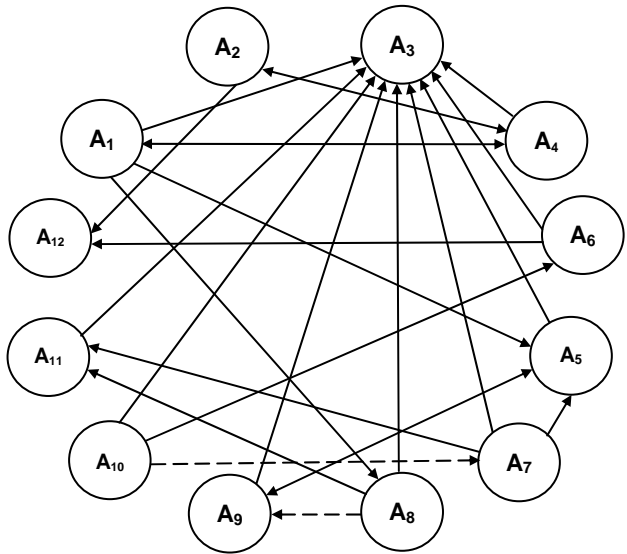

FIGURE: 4.2.1

Let the corresponding neutrosophic connection matrix of the neutrosophic directed graph be given by N(E). Now as the opinion was taken from the same expert we give in this neutrosophic dynamical system the same state vectors for this will help us in the easy comparison of both the systems.

The neutrosophic matrix N(E) is as follows:

$$
\begin{array}{c}
\phantom{A_1}\\
A_1\\
A_2\\
A_3\\
A_4\\
A_5\\
A_6\\
A_7\\
A_8\\
A_9\\
A_{10}\\
A_{11}\\
A_{12}
\end{array}
\begin{array}{cccccccccccc}
A_1 & A_2 & A_3 & A_4 & A_5 & A_6 & A_7 & A_8 & A_9 & A_{10} & A_{11} & A_{12}\\
0 & 0 & 1 & 1 & 1 & 0 & 0 & 1 & 0 & 0 & 0 & 0\\
0 & 0 & 0 & 1 & 0 & 0 & 0 & 0 & 0 & 0 & 0 & 1\\
0 & 0 & 0 & 0 & 0 & 0 & 0 & 0 & 0 & 0 & 0 & 0\\
1 & 1 & 1 & 0 & 0 & 0 & 0 & 0 & 0 & 0 & 0 & 0\\
0 & 0 & 1 & 0 & 0 & 0 & 0 & 0 & 1 & 0 & 0 & 0\\
0 & 0 & 1 & 0 & 0 & 0 & 0 & 0 & 0 & 0 & 0 & 1\\
0 & 0 & 1 & 0 & 1 & 0 & 0 & 0 & 0 & 0 & 1 & 0\\
0 & 0 & 1 & 0 & 0 & 0 & 0 & 0 & I & 0 & 1 & 0\\
0 & 0 & 1 & 0 & 1 & 0 & 0 & 0 & 0 & 0 & 0 & 0\\
0 & 0 & 1 & 0 & 0 & 1 & I & 0 & 0 & 0 & 0 & 0\\
0 & 0 & 1 & 0 & 0 & 0 & 0 & 0 & 0 & 0 & 0 & 0\\
0 & 0 & 0 & 0 & 0 & 1 & 0 & 0 & 0 & 0 & 0 & 0
\end{array}
$$

Let R = (0 0 0 1 0 0 0 0 0 0 0 0 ) be neutrosophic state vector in which only the node $A_4$ that is the types of profession is in on state and all other nodes are in the off state. R is passed into the above neutrosophic connection matrix N(E),

$$
\begin{array}{lllll}
RN\,(E) & \hookrightarrow & (1\ 1\ 1\ 1\ 0\ 0\ 0\ 0\ 0\ 0\ 0\ 0) & = & R_1\\
R_1N\,(E) & \hookrightarrow & (1\ 1\ 1\ 1\ 1\ 0\ 0\ 1\ 0\ 0\ 0\ 1) & = & R_2\\
R_2N\,(E) & \hookrightarrow & (1\ 1\ 1\ 1\ 1\ 1\ 0\ 1\ 1\ 0\ 1\ 1) & = & R_3\\
R_3N\,(E) & \hookrightarrow & (1\ 1\ 1\ 1\ 1\ 1\ 0\ 1\ I\ 0\ 1\ 1\ ) & = & R_4\\
R_4N\,(E) & \hookrightarrow & (1\ 1\ 1\ 1\ 1\ 1\ 0\ 1\ I\ 0\ 1\ 1) & = & R_4
\end{array}
$$

Here the symbol '$\hookrightarrow$' denotes that at every stage the resultant vector has been updated and thresholded.

Thus we see when the attribute profession is in the on state we see the relation between profession and economic status become an indeterminate. For one cannot always say by the type of the profession ones economic status can be determined for it depends on the mouths he has to feed his commitments and other factors which determine the economic status and it has not much to do with the profession.

Now let us consider the state vector S = (1 0 0 0 0 0 0 0 0 0 0 ) i.e., the only neutrosophic vector easy availability of money alone is in the on state and all other vectors are in the off state. To study the effect of S on the neutrosophic dynamical system N(E).

$$SN(E) \hookrightarrow (1\ 0\ 1\ 1\ 1\ 0\ 0\ 1\ 0\ 0\ 0\ 0) = S_1 \text{ say}$$

$$S_1N(E) \hookrightarrow (1\ 1\ 1\ 1\ 1\ 0\ 0\ 1\ I\ 0\ 1\ 0) = S_2 \text{ (say)}$$

$$S_2N(E) \hookrightarrow (1\ 1\ 1\ 1\ 1\ 0\ 0\ 1\ I\ 0\ 1\ 1) = S_3 \text{ (say)}$$

$$S_3N(E) \hookrightarrow (1\ 1\ 1\ 1\ 1\ 1\ 0\ 1\ I\ 0\ 1\ 1) = S_4 \text{ (say)}$$

$$S_4N(E) \hookrightarrow (1\ 1\ 1\ 1\ 1\ 1\ 0\ 1\ I\ 0\ 1\ 1) = S_4$$

$S_4$ is a fixed point. Thus the hidden pattern of the neutrosophic dynamical system is a fixed point easy money makes all other state on except $A_7$ and $A_{10}$ i.e., the absence of social responsibility has however no relation with easy money also the concept more leisure has nothing to do with easy money. However the attribute economic status is an indeterminate for its relation with easy money remains an indeterminate.

The reader is given as an exercise the study the effect of the only attribute easy money in the on state of the fuzzy dynamical system and arrive conclusions based on it and compare it with the FCM. To have an easy comparison between FCMs and NCMs consider the state vector

$$X = (0\ 0\ 0\ 0\ 0\ 0\ 0\ 0\ 0\ 1\ 0).$$

Using X in the neutrosophic connection matrix we get

$$XN(E) \hookrightarrow (0\ 0\ 1\ 0\ 0\ 0\ 0\ 0\ 0\ 1\ 0) = X_1 \text{ (say)}$$

$$X_1N(E) \hookrightarrow (0\ 0\ 1\ 0\ 0\ 0\ 0\ 0\ 0\ 1\ 0) = X_2 = X.$$

Thus we see in both the cases use of FCM and NCM, this case gives the same resultant vector.

Now consider the on state of the nodes profession, wrong/bad company and addiction to habit forming substances and visiting CSWs and all other nodes are in the off state.

Let $Y = (0\ 0\ 0\ 1\ 1\ 1\ 0\ 0\ 0\ 0\ 0).$

The effect of Y on the neutrosophic connection matrix N(E)

$$YN \, (E) \quad \hookrightarrow \quad (1\ 1\ 1\ 1\ 1\ 1\ 0\ 0\ 1\ 0\ 0\ 1) \quad = \quad Y_1 \, \text{(say)}$$

$$Y_1 N(E) \quad \hookrightarrow \quad (1\ 1\ 1\ 1\ 1\ 1\ 0\ 1\ I\ 0\ 1\ 1) \quad = \quad Y_2 \, \text{(say)}$$

$$Y_2 N(E) \quad \hookrightarrow \quad (1\ 1\ 1\ 1\ 1\ 1\ 0\ 1\ I\ 0\ 1\ 1) \quad = \quad Y_3 = Y_2.$$

Thus the neutrosophic hidden pattern of the dynamical system is a fixed point. With the three attributes: profession, wrong/ bad company and addiction of habit forming substances and CSWs are in the on state we see all attributes come to on state except $A_7$ and $A_{10}$ and $A_9$ becomes an indeterminate state the off state of the absence of social responsibility shows that these three nodes have no influence on the dynamical system.

Also we see the economic status is an indeterminate state. However more leisure has no effect on the system. Now we consider the state vector A = (0 0 0 0 0 0 0 0 0 1 1) where only the two attributes Machismo/ Exaggerated masculinity and No awareness are the attributes which are in the on states and all other nodes are in the off state.

Let A $\quad = \quad$ (0 0 0 0 0 0 0 0 0 1 1).

The effect of A on the neutrosophic dynamical system N(E)

$$AN \, (E) \quad \hookrightarrow \quad (0\ 0\ 1\ 0\ 0\ 1\ 0\ 0\ 0\ 0\ 1\ 1) \quad = \quad A_1 \, \text{say}$$

$$A_1 N(E) \quad \hookrightarrow \quad (\ 0\ 0\ 1\ 0\ 0\ 1\ 0\ 0\ 0\ 0\ 1\ 1).$$

When the public lack awareness about the disease and the Machismo/ Exaggerated Masculinity we see they easily visit CSWs and have addiction to habit forming substances and CSWs. Several other on state vectors were sent into the system and the resultant information were recorded for our conclusions given in Chapter Seven.

One vital point obtained using neutrosophy theory is that the economic status is an indeterminate factor while studying the profession; so the truth it imparts is profession and economic status can not be always interlinked.

Only the concept of profession may have some link with the concept of social status for we have taken basically the sample of sixty migrant labourers only from the poor and lower middle class strata so working on the social status in our context will have no significance.

Next we take the view of the same expert now he is given the option that he can denote the notion it is impossible for him interrelate or relation between two concepts is an indeterminacy.

As we have just now pointed out that the persons/ patients whom we had interviewed are largely from the poor or lower middle class. It lies like a myth about the upper middle class or middle class or rich educated people that too from the city who are infected by HIV/AIDS. We are only in a position to conjecture that even if the disease is prevalent in this society it is largely suppressed and no data can ever be got.

For they seek treatment secretly taking confidence of private doctors and paying them a good or sumptuous amount as medical fees. So at every stage it has become our duty to mention that the study and analysis pertains only to poor / very lower middle class with no education and who mainly reside in village or remote village. But all of them had acquired the disease only by migration. Thus our conclusion based on study and analysis regarding the data concerning the remote uneducated poor villagers who have obtained the disease by migration happens to be complete authenticated and implementation of these suggestions will prove to be certainly beneficial to them and avoid the increase in the number of HIV/AIDS patients from this group.

Thus this expert does not include the attribute the economic status as one of the attributes for in his opinion the disease prevalent in all economic status but the publicity or openness is totally lacking in the middle class, upper middle class and rich. They confide it for the reasons best known to them. Thus this expert has already given the justification for working only with seven attributes given in page 54.

The concepts / attributes given by this expert are:

$A'_1$ - Easy availability of money
$A'_2$ - Wrongful company and addictive   habits
$A'_3$ - Visits CSWs
$A'_4$ - Socially irresponsible / free
$A'_5$ - Macho Behaviour
$A'_6$ - More leisure
$A'_7$ - No awareness about disease.

The neutrosophic directed graph given by the expert.

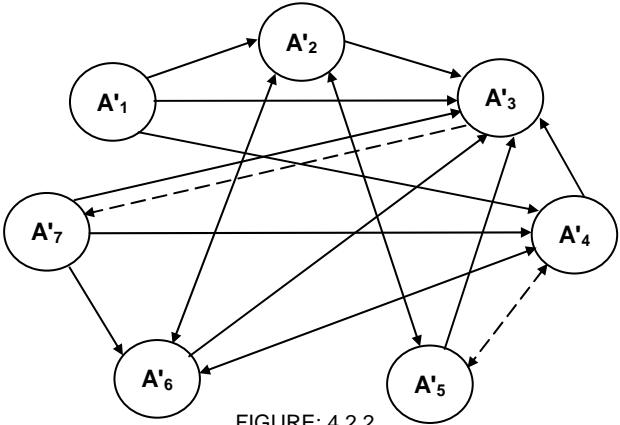

FIGURE: 4.2.2

The associated neutrosophic relational matrix N (A) is given in what follows:

$$N(A) = \begin{array}{c} \\ A'_1 \\ A'_2 \\ A'_3 \\ A'_4 \\ A'_5 \\ A'_6 \\ A'_7 \end{array} \begin{array}{c} A'_1 \; A'_2 \; A'_3 \; A'_4 \; A'_5 \; A'_6 \; A'_7 \\ \begin{bmatrix} 0 & 1 & 1 & 1 & 0 & 0 & 0 \\ 0 & 0 & 1 & 0 & 1 & 1 & 0 \\ 0 & 0 & 0 & 0 & 0 & 0 & I \\ 0 & 0 & 1 & 0 & I & 1 & 0 \\ 0 & 1 & 1 & I & 0 & 0 & 0 \\ 0 & 1 & 1 & 1 & 0 & 0 & 0 \\ 0 & 0 & 1 & 1 & 0 & 1 & 0 \end{bmatrix} \end{array}$$

Also he says that the relation between the attributes visiting CSWs and no awareness about the disease is also an indeterminate. For visiting CSWs need not be only due to the fact he has no awareness about this disease. Their may be cases where he says some other factors which would have forced him to visit CSWs are: extraordinary urge for sex, dejection, absence of wife in case of married men, long stay away from home, family dispute or some inexpressible cause which the client is well aware of. Thus he says he would indicate this relation by *I*. He says the relation between macho behaviour and socially free and the relation socially free or irresponsible and the macho behaviour are indeterminate to him. Thus using his newly given Neutrosophic graph and its related matrix N(A) we proceed to work for the analysis and conclusions i.e., the effect of each of the state vectors. N(A) denotes the neutrosophic connection matrix of the neutrosophic directed graph. Now to find the stability consider the state vector

$$X \quad = \quad (1\ 0\ 0\ 0\ 0\ 0\ 0)$$

that is easy money which is dependent on the profession i.e., $A_1'$ alone is in the on state, and all other attributes are in the off state. Now passing X into N(A) we get

| | | | | |
|---|---|---|---|---|
| XN(A) | ↪ | (1 1 1 1 0 0 0) | = | $X_1$ (say) |
| $X_1$N(A) | ↪ | (1 1 1 1 $I$ 1 $I$) | = | $X_2$ say |
| $X_2$N(A) | ↪ | (1 1 1 1 $I$ 1 $I$) | = | $X_3 = X_2$ (fixed). |

Thus we see the hidden pattern is a fixed point. All the nodes becomes on except the two nodes Macho behaviour and the No awareness about the disease. Both these co-ordinates are indeterminates, thus easy money dependent on profession leaves the nodes, Macho behaviour and the awareness about the disease as indeterminates i.e., the neutrosophic system is not able to state whether there is any influence of the profession on the Macho behaviour or on the awareness about the disease. Further all other attributes come to on state when the person has easy availability of money, which is dependent on his profession.

One can easily compare the neutrosophic result with the fuzzy result. This work is assigned to the interested reader. It is however pertinent to mention that fuzzy method has no effect on these two coordinates; so these two coordinates remain in the off state.

But by using neutrosophy we see that it clearly states we cannot find the relation between macho behaviour and the easy money dependent on the profession. There is always a great difference between saying explicitly no relation exists and the statement that the relation between them is an indeterminate. Thus in our opinion Neutrosophy model in this situation gives a sensitive result than FCMs.

Now let us consider the node 'macho behaviour' to be in the on state i.e., Y = (0 0 0 0 1 0 0) i.e., all other states are in the off state. The effect of Y on the neutrosophic dynamical system N(A) is given by

$$Y \, N(A) \quad \hookrightarrow \quad (0\ 1\ 1\ I\ 1\ 0\ 0) \quad = \quad Y_1 \ (\text{say})$$

$$Y_1 N \, (A) \quad \hookrightarrow \quad (0\ 1\ 1\ I\ 1\ I\ I) \quad = \quad Y_2 \ (\text{say})$$

$$Y_2 N(A) \quad \hookrightarrow \quad (0\ I\ I\ I\ 1\ I\ I) \quad = \quad Y_3 \ (\text{say})$$

$$Y_3 N(A) \quad \hookrightarrow \quad (0\ I\ I\ I\ 1\ I\ I) \quad = \quad Y_4 = Y_3 \ (\text{fixed point})$$

('↪' denotes the resultant vector after passing through the dynamical system has been thresholded and updated.) The hidden pattern of the system is a fixed point. By the use of Neutrosophic model we see that the effect of Macho behaviour is indeterminate on the other coordinates. For a person with Macho behaviour need not visit CSWs or according to the neutrosophic model one does not know i.e., the relation of macho behaviour and the node of visiting CSWs are indeterminate.

Further a person with Macho behaviour need not be a socially irresponsible / free or he may be socially irresponsible/ free thus this adaptation of neutrosophy makes this person/ expert to state that no relation exists between Macho behaviour and the easy money but the existence or the non existence of other relations cannot be said precisely that is these concepts happen to be indeterminates. This is the marked difference between the adaptation fuzzy theory and neutrosophic theory for we see that neutrosophic theory is better or more sensitive in the case of the problem of analyzing the socio economic status of a HIV/AIDS patients who are migrant labourers. We only say for this model the Neutrosophic model is more powerful than the fuzzy cognitive model. The interested reader is expected to work with more state vectors and derive conclusions. However our results are based on the study of the dynamical systems with several state vectors.

We have taken the $3^{rd}$ expert as a refined and a reformed HIV/AIDS patients (For this patient says several of the male in patients in this hospital visit CSWs when they are little better).

$P_1$ – Most of the HIV/AIDS patients do not have binding relationship with family that is why they are often tempted to visit CSW, drink, smoke and waste money on bad company.

$P_2$ – Male chauvinism and macho behaviour.

$P_3$ – Most of the men with male ego and who have no binding in family often think of women as an inferior object and an object of only pleasure.

$P_4$ – Bad company and bad habits.

$P_5$ – Socially irresponsible.

$P_6$ – Several of these persons go to CSW due to uncontrollable emotions that is why some HIV/AIDS patients had acknowledged to him if ever he had self control he would not have erred. That is why most of them do not use safe sex methods even if the CSW advises them.

$P_7$ – Gluttony. This man says that ravenous appetite for food leads to the same kind of urge for sex.

$P_8$ – More leisure.

$P_9$ – According to him that when men/women after a age do not have a goal or an ideal in life.

$P_{10}$ – A person for whom only the body functions, i.e. have no ideal or higher aspirations in life.

$P_{11}$ – CSW.

$P_{12}$ – Most when asked why they visited CSWs the answer will be to 'enjoy life' in a 'jolly mood' 'pass time' and so on.

$P_{13}$ – When mainly they are in a different state they are more free from fear of being noticed so without any sense of shame/ sin they visit CSW.

$P_{14}$ – Most of them reply that they have visited countless CSWs, this is their attitude.

$P_{15}$ – This point he says with his personal experience that most of the agricultural labour have taken up the profession of truck or lorry drivers for livelihood. Thus one of the major reasons is failure of agriculture/poor yield or not up to anticipation. He says in his village alone over 70 families have migrated some to Bombay and some to Madras. Many men have taken the profession of drivers.

With these as nodes we give the directed graph of this expert. Now the same expert is asked the question now with option that he can also say the relation between two nodes to be an indeterminate. Thus we use the same expert for this model to obtain the neutrosophic relation.

FIGURE 4.2.3

Using the neutrosophic graph we obtain the following neutrosophic connection matrix N(P).

$$
\begin{bmatrix}
0 & 1 & 1 & 1 & 1 & 0 & 0 & 0 & 0 & 0 & 1 & I & 0 & 0 & 0 \\
1 & 0 & 1 & 0 & 1 & 0 & 0 & 0 & 0 & 0 & 1 & 1 & 0 & 0 & 0 \\
1 & 1 & 0 & 1 & 0 & 0 & 0 & 0 & 0 & 0 & 1 & 0 & 0 & 0 & 0 \\
0 & 0 & 0 & 0 & 1 & 0 & 0 & 1 & 1 & 0 & 1 & 0 & 0 & 0 & 0 \\
0 & 0 & 0 & 0 & 0 & 0 & 0 & 0 & 0 & 0 & 1 & 1 & 0 & 0 & 0 \\
0 & 0 & 0 & 0 & 0 & 0 & 1 & 1 & 1 & 0 & 1 & 0 & 0 & 1 & 0 \\
0 & 0 & 0 & 0 & 0 & 1 & 0 & 0 & 1 & 1 & 1 & 1 & 0 & 1 & 0 \\
0 & 0 & 0 & 0 & 0 & 0 & 0 & 0 & 1 & 0 & 1 & 1 & 0 & 0 & 0 \\
0 & 0 & 0 & 0 & 0 & 0 & I & 0 & 0 & 0 & 1 & 0 & 0 & 0 & 0 \\
0 & 0 & 0 & 0 & 0 & 0 & 0 & 0 & 1 & 0 & 1 & 0 & 0 & 1 & 0 \\
0 & 0 & 0 & 0 & 0 & 0 & 0 & 0 & 0 & 0 & 0 & 0 & 0 & 0 & 0 \\
0 & 0 & 0 & 0 & 0 & 0 & 0 & 0 & 1 & 0 & 1 & 0 & 1 & 0 & 0 \\
0 & 0 & 0 & I & 0 & 0 & 0 & 0 & 0 & 0 & 1 & 1 & 0 & 0 & I \\
0 & 0 & 0 & 0 & 0 & 0 & 0 & 0 & 0 & 0 & 1 & 0 & 0 & 0 & 0 \\
0 & 0 & 0 & 0 & 0 & 0 & 0 & 0 & 0 & 0 & 0 & 0 & 1 & 0 & 0
\end{bmatrix}
$$

Consider the state vector

$$T \qquad = \qquad (0\ 0\ 0\ 0\ 0\ 0\ 1\ 0\ 0\ 0\ 0\ 0\ 0\ 0\ 0)$$

that is only the attribute $A_7$ is in the on state and all other nodes are in the off state i.e., at the first instance the person is a glutton with regard to food habits, to study the effect of on the neutrosophical system N(P) consider

$$TN(P) \hookrightarrow (0\ 0\ 0\ 0\ 0\ 1\ 1\ 0\ 1\ 1\ 1\ 1\ 0\ 1\ 0) = T_1 \text{ (say)}$$

$$T_1N(P) \hookrightarrow (0\ 0\ 0\ 0\ 0\ 1\ 1\ 1\ 1\ 1\ 1\ 1\ 1\ 1\ 0) = T_2 \text{ (say)}$$

$$T_2\,N(P) \hookrightarrow (0\ 0\ 0\ I\ 0\ 1\ 1\ 1\ 1\ 1\ 1\ 1\ 1\ 1\ I) = T_3 \text{ (say)}$$

$$T_3\,N(P) \hookrightarrow (0\ 0\ 0\ I\ I\ 1\ 1\ 1\ I\ 1\ 1\ 1\ I\ 1\ I) = T_4 \text{ (say)}$$

$$T_4\,N(P) \hookrightarrow (0\ 0\ 0\ I\ I\ 1\ 1\ 1\ I\ 1\ 1\ I\ I\ 1\ I) = T_5 \text{ (say)}$$

$$T_5\,N(P) \hookrightarrow (0\ 0\ 0\ I\ I\ 1\ 1\ 1\ I\ 1\ 1\ I\ I\ 1\ I) = T_5$$

$T_5$ is a fixed point. Thus the hidden pattern of the system is a fixed point with 6 of the coordinates as indeterminates i.e., the on state of $A_7$ cannot precisely say whether these six coordinates become 'on' or remain in the 'off' state; i.e., they are in the indeterminate state. Thus gluttony in food habits has no or nil impact or one relation with "no binding with family" and it has nothing to do with male chauvinism. Gluttony has no effect on women as inferior object all these three coordinates remain in the off state.

Whether there is any relation between gluttony in food habits and bad company and bad habits remains as an indeterminate also the node "socially irresponsible" remains as an indeterminate i.e., it is nether 'on' nor 'off'. But however gluttony in food habits makes the node uncontrollable sexual feelings to be in the on state. One is not able to say whether there exists any relation between gluttony and "more leisure " for the attribute 'more leisure' remains as an indeterminate coordinate. However the on state of the node gluttony has positive effect on the coordinates "no other work for the brain" "Only physically active" and Visits CSWs i.e., all these three coordinates becomes on when a person is a glutton. But weather a glutton "enjoy life or jolly mode" is in the indeterminate state. Also the node unreachable by friends or relatives remains in the indeterminate state. i.e., one is not in a position to find or say weather any relation exists between gluttony and unreachable by friends and relatives. But however to ones surprise a glutton takes pride in visiting countless CSWs for this node becomes on. Finally the attribute, "Failure of agriculture" and its relation with the glutton remain as indeterminate coordinate.

Thus we see a vast difference between the two models neutrosophy cognitive maps and fuzzy cognitive maps for we see the results given by neutrosophy is more sensitive than FCM in the study of this problem.

Next we consider the coordinate male chauvinism is in the on state

$$B = (0\ 1\ 0\ 0\ 0\ 0\ 0\ 0\ 0\ 0\ 0\ 0\ 0\ 0\ 0)$$

i.e., all state vectors are in the off state. We study the effect of S on the neutrosophic dynamical system N(P).

$$S\ N(P) \hookrightarrow (1\ 1\ 1\ 0\ 1\ 0\ 0\ 0\ 0\ 0\ 1\ 1\ 0\ 0\ 0) = S_1\ (say)$$

$$S_1\ N(P) \hookrightarrow (1\ 1\ 1\ 1\ 1\ 0\ 0\ 0\ 0\ 0\ 1\ 1\ 1\ 0\ 0\ ) = S_2\ (say)$$

$$S_2\ N(P) \hookrightarrow (1\ 1\ 1\ 1\ 1\ 0\ 0\ 1\ 1\ 0\ 1\ 1\ 1\ 1\ I) = S_3\ (say)$$

$$S_3\ N(P) \hookrightarrow (1\ 1\ 1\ 1\ 1\ 0\ I\ 1\ 1\ 0\ 1\ 1\ 1\ 1\ I) = S_4\ (say)$$

$$S_4\ N(P) \hookrightarrow (1\ 1\ 1\ 1\ 1\ I\ I\ 1\ 1\ I\ 1\ 1\ 1\ 1\ I) = S_5\ (say)$$

$$S_5\ N(P) \hookrightarrow (1\ 1\ 1\ 1\ 1\ I\ I\ I\ 1\ 1\ I\ 1\ 1\ I\ I) = S_6\ (say)$$

$$S_6\ N(P) \hookrightarrow (1\ 1\ 1\ 1\ 1\ I\ I\ I\ 1\ I\ 1\ I\ I\ I\ I) = S_2\ (say)$$

$$S_7\ N(P) \hookrightarrow (1\ 1\ 1\ 1\ 1\ I\ I\ I\ 1\ I\ 1\ I\ I\ I\ I) = S_2$$

a fixed point. If male chauvinist alone is in the on state we see in the resultant vector the nodes $P_1$, $P_3$, $P_4$, $P_5$ and $P_{11}$ become on and the nodes $P_6$, $P_7$, $P_8$, $P_{10}$, $P_{12}$, $P_{13}$, $P_{14}$ and $P_{15}$ are indeterminate. We see that the male egoist has 'no binding with his family', 'thinks women as 'inferior objects',' 'bad company and formation of bad habits,' 'socially irresponsible,' 'have no other work for the brain' and 'visits CSWs' so all these nodes become on. But we see that when a male chauvinist has uncontrollable sex feelings or he is 'glutton' or he has 'more leisure', 'only physically active' and 'enjoys life – jolly mood', 'unreachable by friends and relatives' 'pride in visiting CSWs' or 'failure of agriculture', these nodes become indeterminate by which we cannot say there no relevance but we can not or we do not know whether there is relevance.

## 4.3 Combined NCM to study the migrant problems

Just in section 4.2 we have defined the combined NCM. In order to obtain the combined effect we use combined NCMs, here we use 5 experts opinion we assign to them the 10 attributes which we have taken as a compulsory one in order to use it in the simple combined NCMs.

The 10 attributes given to them after consultation with several experts are as follows: since the model is mainly based only to study the relationship between the catching of HIV/AIDS and migrancy. We give the nodes $A_1, \ldots A_{10}$ and briefly describe what each of the nodes $A_1, A_2, \ldots, A_{10}$ signifies

$A_1$ – Easy money – no channelizing of mind / body after work
$A_2$ – Living away from family for weeks
$A_3$ – Free from fear of being watched by friends and relatives
$A_4$ – Bad habits / company
$A_5$ – Lack of proper family binding
$A_6$ – No proper counseling for saving money/ proper spending
$A_7$ – Lack of any proper knowledge about HIV/AIDS.
$A_8$ – CSWs only recreation
$A_9$ – No proper motivation / higher aims in life
$A_{10}$ – Age a big disadvantage.

We briefly describe each of these attributes relative to the context of migrant labourers.

**$A_1$ – easy money - no  channelizing of body / mind after work**

As far as the migrant labourers are concerned that too specially affected by the HIV/AIDS we see that they are mainly daily wagers or given wage after completing a task i.e., going on trips to deliver goods or people or weekly wagers i.e., paid weekly. Whatever be their job as far as the 60 migrant labourers are concerned 58 were not the people who worked for a month to earn their pay monthly. So for them money was easily earned, further after hours of labour their only two main recreations were liquor and visiting CSWs. Apart from this they had no further thinking.

This is very clear from the interviews. When we say easy money they after eating and stay they have enough money left over to spend on CSWs and liquor. None of them ever had the intension or the knowledge of saving money for when they were with the disease most of them said they had no money even to buy a cup of tea to drink or to buy and eat what they liked. If they ever had the habit of saving they would not have been in this plight. Some of them really felt their present poor status and their past living how lavishly they spent money.

## A$_2$ - Living away from family for weeks

Majority of these migrant labourers lived away from their family for at least a week or more and some of them had their profession to be truck / lorry drivers so they had a superstition that their body would be relived of heat if they indulged in sex.

So most of the truck /lorry drivers said they visited CSWs to lessen the heat of their body. Also as they basically lacked self control they easily visited CSWs for sex. As most of these men were in their twenties and early thirties and also married since 3 to 5 years they easily became victims of temptation . Thus at this point it is pertinent to mention that even if the patient was just in his forties must have been infected at his early thirties and late twenties only. Age played a major havoc in the life of the migrant labourers to catch HIV/AIDS.

## A$_3$ - Free from fear of being watched by friends and relatives

These migrant labourers when they are away from their families i.e., when they are away from their home town they are free in a sense that they are fully free from being watched reported by their friends and relatives.

For visiting a CSW according to the Indian society that too in the Tamil society is a crime, an heinous act, an act very much condemned by the family people. So these labourers when they want to or visit the CSWs they are free from the fear of being watched by friends and relatives as they live in a new place and with an unknown identify so it favours them to act with courage.

Once they become used to, a stage comes when they have become really addicted to for it is no exaggeration on the part of us for we saw out of 60 HIV/AIDS infected migrant labourers only 1 has not visited CSWs. 58 of them said with pride that they had visited in numerable times the CSWs to enjoy life or in jolly mood. Thus we see an unknown places provides a comfortable atmosphere to visit CSWs.

### A₄ - Bad habits / company

Only 10% of the HIV/AIDS affected migrant labourers accepted they did all acts of their own of the 60, 54 of them said they smoke, drink or visit CSWs because of their friends. Almost all of them said that it was their friends who spoilt them. Thus we see that bad habits / bad company is also one of the reasons for these migrant labourers to catch HIV/AIDS.

### A₅ - Lack of proper family binding

This node is mainly put for with a view if the migrant labourers afflicted by HIV/AIDS have some binding on the family then certainly

1. They would disclose the disease they suffer from
2. They would test their wife and children for the same
3. They would not infect their wife
4. Above all they would not visit numerous CSWs.

Hence we are justified in taking this node as a related one with the problem

### A₆ - No proper counseling for saving money / proper spending.

Most of the migrant labourers are only daily wagers they do not have an association or any form of union which can promote their welfare by some welfare schemes like building houses, educating children, buying land etc. or do not know the art of saving the money they earn by a hard labour.

### A₇ - Lack of any knowledge about HIV/AIDS.

Most of the migrant labourers whom we have interviewed that too hailing from villages are uneducated so they do not know the basic knowledge about it, like the way it is caused, the way it spreads from one to another, the enormous physical suffering due to HIV/AIDS and above all it has no cure but is a livable disease.

For certainly if one is fully aware of the problem in catching HIV/AIDS the persons concerned would by very careful. For 58 of the 60 HIV/AIDS patients we interviewed only after knowing that they are infected with HIV/AIDS they have come to know about the disease.

For this reason we see out of 60 at least 52 of them wanted to commit suicide when they came to know that they had HIV/AIDS. Thus lack of knowledge about HIV/AIDS had made them HIV/AIDS patients.

## A$_8$ - CSWs are only recreation.

Most of the HIV/AIDS affected migrant labourers openly acknowledged that after the days work and after taking liquor only there was urge for sex.

So they were forced to visit CSWs. Also as they stayed away from the family and as they were mostly in places which had alien language or atmosphere, they were forced to visit CSWs for sex. They had no other aim or goal.

That is why 98% of them said they had visited numerous CSWs to enjoy life in jolly mood and so on and so forth.

## A$_9$ - No proper motivation / higher aims in life.

Most of the migrant labourers were from remote villages with no education so they had no higher aspirations in life they did not even feel they have to earn and it was equally important to save or they did not know the value of education or the moral values or material values so very easily they became victims of all bad habits so they had no higher aims in life.

Everything started with smoke, liquor and CSWs and ended with the same. For even while we were interviewing most of them did not aspire for educating their children or improve the state of their living conditions. They were so much absorbed with their living style and satisfied by it.

Only because they were infected by HIV/AIDS which not only took away their health, even doing routine but it made them absolutely immobile and very sick for any thing even could not eat or walk.

If there were strong supportive groups or the remote uneducated villagers were empowered to improve themselves certainly this plight can be changed without any difficulty.

### A₁₀ - Age a big disadvantage

Most of the migrant labourers come out from their native village for seeking employment only in the age group from 17 to 35 where 17 yrs is the approximate least and 35 years is the approximate greatest there may be age group little lesser than the least or little greater than the greatest. But we are mostly interested in taking the maximum clustering age group.

That is why most of the patients affected by HIV/AIDS among the migrant labourers were in the age group 26 to 35; for it may take for the infection to show its symptoms in a period of 3 years to 13 years depending on the health conditions of these persons their food habits and also other habits. Thus we see age in which they are exposed plays a major criteria in getting the infection as they are vulnerable to all types of temptations including CSWs. Thus they catch the disease mainly due to ignorance catalyzed by age.

Thus if there were effective awareness program given from the age of 13 years in villages this could be prevented.

Thus we have taken the age as a major criteria for being infected by HIV/AIDS.

Now having explained the 10 attributes we give the expects opinion in the neutrosophic model to finally obtain the combined effect.

The neutrosophic directed graph given by the first expert.

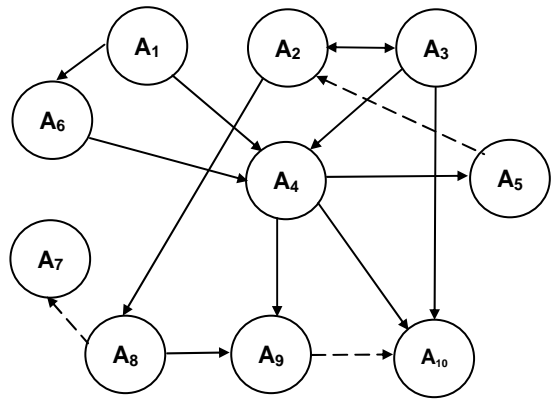

FIGURE: 4.3.1

The related neutrosophic connection matrix N(E₁)

$$N(E_1) = \begin{array}{c} \\ A_1 \\ A_2 \\ A_3 \\ A_4 \\ A_5 \\ A_6 \\ A_7 \\ A_8 \\ A_9 \\ A_{10} \end{array} \begin{array}{c} A_1 \ A_2 \ A_3 \ A_4 \ A_5 \ A_6 \ A_7 \ A_8 \ A_9 \ A_{10} \\ \begin{bmatrix} 0 & 0 & 0 & 1 & 0 & 1 & 0 & 0 & 0 & 0 \\ 0 & 0 & 1 & 0 & 0 & 0 & 0 & 1 & 0 & 0 \\ 0 & 0 & 0 & 1 & 0 & 0 & 0 & 0 & 0 & 1 \\ 0 & 0 & 0 & 0 & 1 & 0 & 0 & 0 & 1 & 1 \\ 0 & I & 0 & 0 & 0 & 0 & 0 & 0 & 0 & 0 \\ 0 & 0 & 0 & 1 & 0 & 0 & 0 & 0 & 0 & 0 \\ 0 & 0 & 0 & 0 & 0 & 0 & 0 & 0 & 0 & 0 \\ 0 & 0 & 0 & 0 & 0 & 0 & I & 0 & 1 & 0 \\ 0 & 0 & 0 & 0 & 0 & 0 & 0 & 0 & 0 & I \\ 0 & 0 & 0 & 0 & 0 & 0 & 0 & 0 & 0 & 0 \end{bmatrix} \end{array}$$

Take any state vector say X = (0 0 1 0 0 0 0 0 0 0 ) i.e., the only node $A_3$ i.e., "free from being watched is in the on state and all other nodes are in the off state. We see the effect of X on the neutrosophical dynamical system $N(E_1)$

$$
\begin{array}{llll}
XN(E_1) & \hookrightarrow & (0\ 0\ 1\ 1\ 0\ 0\ 0\ 0\ 0\ 1) & = & X_1 \text{ (say)} \\
X_1N(E_1) & \hookrightarrow & (0\ 0\ 1\ 1\ 1\ 1\ 0\ 0\ 1\ 1) & = & X_2 \text{ (say)} \\
X_2N(E_1) & \hookrightarrow & (0\ I\ 1\ 1\ 1\ 0\ 0\ 0\ 1\ 1) & = & X_3 \text{ (say)} \\
X_3N(E_1) & \hookrightarrow & (0\ I\ 1\ 1\ 1\ 0\ 0\ I\ 1\ 1) & = & X_4 \text{ (say)} \\
X_4N(E_1) & \hookrightarrow & (0\ I\ 1\ 1\ 1\ 0\ I\ I\ I\ 1) & = & X_5 \text{ (say)} \\
X_5N(E_1) & \hookrightarrow & (0\ I\ 1\ 1\ 1\ 0\ I\ I\ I\ 1) & = & X_6 = X_5
\end{array}
$$

$X_5$ is the fixed point. Thus we get the hidden pattern of the neutrosophic dynamical system to be a fixed point which asserts if a person is free from being watched by friends and relatives and is alone, is in the on state then the states $A_4$, $A_5$ and $A_{10}$ become on and $A_2$, $A_7$, $A_8$ and $A_9$ become as indeterminable coordinates. i.e. if a migrant labour in his age group from 20-35 is free from the fear of being catched by relatives and friends certainly he has or develops bad habits and bad friends, he certainly lacks family binding age is the cause of it i.e., age is the big disadvantage for being so.

However the coordinates living away from the family is an indeterminate for it is neither 'off' nor 'on'. Lack of any proper knowledge about HIV/AIDS is also an indeterminate coordinate for one cannot say the on or off state of it. Likewise one cannot say he visits CSWs as only recreation or does not visit CSWs for it also is an indeterminate node. For this expert cannot say for certain he visits or he does not visit CSWs, $A_8$ also remains as an indeterminate coordinate. $A_9$ that is 'No proper motivation / higher aims in life' is neither off nor on its is also an indeterminate coordinate for one is not able to state whether free from being watched has any relation or not.

Next we consider the on state of the coordinate $A_5$. $A_{10}$ that is 'age is a big disadvantage' and $A_5$ – 'Lack of proper family binding', all other coordinates remain in the 'off' state. Let $Y = (0\ 0\ 0\ 0\ 1\ 0\ 0\ 0\ 0\ 1)$ be the instantaneous state vector we study the effect of Y on the dynamical system $N(E_1)$

$$YN(E_1) \quad \hookrightarrow \quad (0\,I\,0\,0\,1\,0\,0\,0\,0\,1) \quad = \quad Y \text{ (say)}$$
$$Y_1N(E_1) \quad \hookrightarrow \quad (0\,I\,I\,0\,1\,0\,0\,I\,0\,1) \quad = \quad Y_2 \text{ (say)}$$
$$Y_2N(E_1) \quad \hookrightarrow \quad (0\,I\,I\,I\,1\,0\,I\,I\,I\,1) \quad = \quad Y_3 \text{ (say)}$$
$$Y_3N(E_1) \quad \hookrightarrow \quad (0\,I\,I\,I\,1\,0\,I\,I\,I\,1) \quad = \quad Y_4 = Y_3$$
$$\text{(a fixed point).}$$

Thus the hidden pattern is a fixed point and one sees however this has no effect on easy money but according to this expert $A_2$, $A_3$, $A_4$, $A_7$, $A_8$, and $A_9$ are indeterminates for one cannot say the effect or no effect of these coordinates but they remain as indeterminate one there by making the analyzer to think about it. Now let us consider the state vector in which only the nodes $A_6$ and $A_8$ are in the on state and all other nodes are in the off state. The effect of Z $= (0\ 0\ 0\ 0\ 0\ 1\ 0\ 1\ 0\ 0\ )$ on the neutrosophical dynamical system $N(E_1)$ is as follows:

$$Z\,N(E_1) \quad \hookrightarrow \quad (0\,0\,0\,1\,0\,1\,I\,1\,1\,0) \quad = \quad Z_1 \text{ (say)}$$
$$Z_1\,N(E_1) \quad \hookrightarrow \quad (0\,0\,0\,1\,1\,1\,I\,1\,I\,1) \quad = \quad Z_2 \text{ (say)}$$
$$Z_2\,N(E_1) \quad \hookrightarrow \quad (0\,0\,0\,1\,1\,1\,I\,1\,I\,I) \quad = \quad Z_3 \text{ (say)}$$
$$Z_3\,N(E_1) \quad \hookrightarrow \quad (0\,0\,0\,1\,1\,1\,I\,1\,I\,I)) \quad = \quad Z_4 = Z_3$$
$$\text{(a fixed point).}$$

Thus the hidden pattern of the dynamical systems is a fixed point in which the on state of the vector $A_6$ and $A_8$ makes the on state of the nodes $A_4$, $A_5$ i.e., when one has no proper counseling and visits CSWs for recreation certainly he has bad habits and enjoys bad company and has no proper binding with the family members and the attributes $A_7$, $A_9$ and $A_{10}$ are indeterminates one cannot say for certain about the on or off state for these nodes viz Lack of any proper knowledge about HIV/AIDS, No proper motivation / higher aims in life and age a big disadvantage remain as an indeterminate. Thus we see this method gives the on state, off state and the indeterminate state of attributes. We see that living away from family for weeks; easy money and no channelizing of mind and body after work and free from fear of being watched by friends and relatives remain in the off state their by asserting that there is no relation between the nodes no proper counseling for saving money and visiting CSWs has no relation with the easy money, no fear of being monitored by friends and relatives and living away from families for weeks. But there exists nodes, which are indeterminate when $A_6$ and $A_8$ are in the on state.

Now for the same set of nodes/attributes we obtain the second experts opinion. The neutrosophic directed graph given by him.

FIGURE: 4.3.2

and the related connection matrix $N(E_2)$ is given below

$$
\begin{array}{c@{\ }c}
 & \begin{array}{cccccccccc} A_1 & A_2 & A_3 & A_4 & A_5 & A_6 & A_7 & A_8 & A_9 & A_{10} \end{array} \\
\begin{array}{c} A_1 \\ A_2 \\ A_3 \\ A_4 \\ A_5 \\ A_6 \\ A_7 \\ A_8 \\ A_9 \\ A_{10} \end{array} &
\left[\begin{array}{cccccccccc}
0 & 0 & 0 & 1 & 0 & 0 & 0 & 0 & 0 & 0 \\
I & 0 & 1 & 0 & 1 & 0 & I & 0 & 0 & 0 \\
0 & 0 & 0 & 1 & 0 & 0 & 0 & 0 & 0 & 0 \\
I & 0 & 0 & 0 & 0 & 1 & 0 & 0 & 0 & 0 \\
0 & 0 & 0 & 0 & 0 & 0 & 0 & 0 & 0 & 0 \\
0 & 0 & 0 & 0 & 0 & 0 & 1 & 0 & 0 & 0 \\
0 & 0 & 0 & 0 & 0 & 0 & 0 & 1 & 0 & 0 \\
0 & 0 & 0 & 0 & 0 & 0 & 0 & 0 & 0 & 0 \\
0 & 0 & 0 & 0 & 0 & 0 & 0 & 0 & 0 & 1 \\
0 & 0 & 0 & 0 & 0 & 0 & 0 & 1 & 0 & 0
\end{array}\right]
\end{array}
$$

We see the effect of the same vectors as that the first expert for comparison and take the same state vector X = (0 0 1 0 0 0 0 0 0 0), that is the node free from being watched by friends and relatives is in the on state. Now we study the effect of the state vector X on the dynamical system N ($E_2$).

$$XN(E_2) \hookrightarrow (0\ 0\ 1\ 1\ 0\ 0\ 0\ 0\ 0\ 0) = X_1 \text{ (say)}$$
$$X_1N(E_2) \hookrightarrow (I\ 0\ 1\ 1\ 0\ 1\ 0\ 0\ 0\ 0) = X_2 \text{ (say)}$$
$$X_2N(E_2) \hookrightarrow (I\ 0\ 1\ 1\ 0\ 1\ 1\ 0\ 0\ 0) = X_3 \text{ (say)}$$
$$X_3N(E_2) \hookrightarrow (I\ 0\ 1\ 1\ 0\ 1\ 1\ 1\ 0\ 0) = X_4 \text{ (say)}$$
$$X_4N(E_2) \hookrightarrow (I\ 0\ 1\ 1\ 1\ 1\ 1\ 1\ 0\ 0) = X_5 \text{ (say)}$$
$$X_5N(E) \hookrightarrow (I\ 0\ 1\ 1\ 1\ 1\ 1\ 1\ 0\ 0) = X_6 \text{ (say)} = X_5.$$

$X_5$ is a fixed point. The hidden pattern of the dynamical system is a fixed point. When the attribute free from being watched by friends and relatives is in the on state we see the attribute $A_1$ "Easy money and no channelizing of mind / body after work" is in an indeterminate state $I$. However living away from the family for weeks with $A_2$ and $A_9$ are in the off state and all other attributes $A_4$, $A_5$, …, $A_8$ become on there by saying when they are not watched they are free so they carry out all bad activity and their age and no motivation in life acts as a catalyst for all these acts. Now we study the effect of state vector Z = (0 0 0 0 0 1 0 1 0 0), the nodes $A_6$ and $A_8$ i.e., no proper counseling for saving money or spending in the proper channel and CSWs only recreation are in the on state and all other attributes are off. We study the effect of Z on the neutrosophical dynamical system N($E_2$).

$$ZN(E_2) \hookrightarrow (0\ 0\ 0\ 0\ 0\ 1\ 1\ 1\ 1\ 0) = Z_1 \text{ (say)}$$
$$Z_1N(E_2) \hookrightarrow (0\ 0\ 0\ 0\ 0\ 1\ 1\ 1\ 0\ 0) = Z_2 \text{ (say)} = Z_1.$$

$Z_2$ is a fixed point. Thus we see when there is no proper counseling of spending the money they earn as duty wages and CSWs are their only recreating the attributes $A_1$, $A_2$, $A_3$, $A_4$. $A_9$ and $A_{10}$ are in the off state that is easy money, living away from family, free from being observed, bad habits and company lack of proper family binding etc. remain in the off state.

Several other state vectors has been substituted and analyzed in the case of the dynamical system N ($E_2$) given by the second expert. Next we consider the effect of the neutrosophic dynamical system as given by the third expert. The directed graph given by the $3^{rd}$ expert related with the 10 attributes $A_1,\ldots,A_{10}$ are

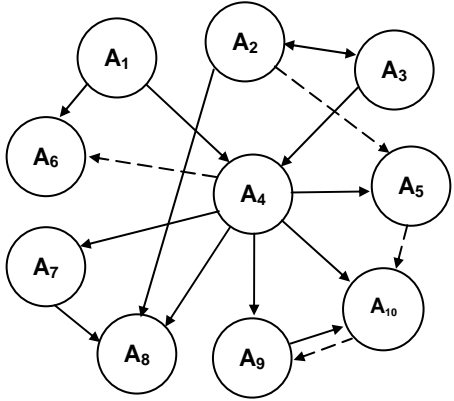

FIGURE: 4.3.3

We give the related $10 \times 10$ connection matrix N ($E_3$)

$$
\begin{array}{c c}
 & \begin{array}{cccccccccc} A_1 & A_2 & A_3 & A_4 & A_5 & A_6 & A_7 & A_8 & A_9 & A_{10} \end{array} \\
\begin{array}{c} A_1 \\ A_2 \\ A_3 \\ A_4 \\ A_5 \\ A_6 \\ A_7 \\ A_8 \\ A_9 \\ A_{10} \end{array} &
\left[
\begin{array}{cccccccccc}
0 & 0 & 0 & 1 & 0 & 1 & 0 & 0 & 0 & 0 \\
0 & 0 & 1 & 0 & I & 0 & 0 & 1 & 0 & 0 \\
0 & 1 & 0 & 1 & 0 & 0 & 0 & 0 & 0 & 0 \\
0 & 0 & 0 & 0 & 1 & I & 1 & 1 & 1 & 1 \\
0 & 0 & 0 & 0 & 0 & 0 & 0 & 0 & 0 & I \\
0 & 0 & 0 & 0 & 0 & 0 & 0 & 0 & 0 & 0 \\
0 & 0 & 0 & 0 & 0 & 0 & 0 & 1 & 0 & 0 \\
0 & 0 & 0 & 0 & 0 & 0 & 0 & 0 & 0 & 0 \\
0 & 0 & 0 & 0 & 0 & 0 & 0 & 0 & 0 & 1 \\
0 & 0 & 0 & 0 & 0 & 0 & 0 & 0 & I & 0 \\
\end{array}
\right]
\end{array}
$$

Now we consider the state vector U = (0 0 1 0 0 0 0 0 0 0) where only $A_3$ is in the on state and all other nodes are in the off state. The effect of U on the neuronal dynamic system N($E_3$):

$$
\begin{array}{llll}
\text{UN} (E_3) & \hookrightarrow & (0\ 1\ 1\ 1\ 0\ 0\ 0\ 0\ 0\ 0) & = \quad U_1 \text{ (say)} \\
U_1 N(E_3) & \hookrightarrow & (0\ 1\ 1\ 1\ 1\ I\ 1\ 1\ 1\ 1) & = \quad U_2 \text{ (say)} \\
U_2 N(E_3) & \hookrightarrow & (0\ 1\ 1\ 1\ I\ I\ 1\ 1\ I\ 1) & = \quad U_3 \text{ (say)} \\
\end{array}
$$

$U_3N(E_3)$    ↪    $(0\ 1\ 1\ 1\ I\ I\ 1\ 1\ I\ I)$    =    $U_4$ (say)

$U_4N(E_3)$    ↪    $(0\ 1\ 1\ 1\ I\ I\ 1\ 1\ I\ I)$    =    $U_5$ =  $U_4$.

$U_4$ is a fixed point. The hidden pattern of the dynamical system is a fixed point when the HIV/AIDS migrant labourer was free from being observed by friends and relatives the attributes $A_2$, $A_4$, $A_7$ and $A_8$ become as on states and the nodes $A_5$, $A_6$, $A_9$ and $A_{10}$ were in an indeterminate state with $A_1$ alone in the off state. Their by meaning easy money had ill effect on them that is when they were free from being watched by friends/relatives. However bad they were for certain living away from the family, they had all bad habits and company, however lack of family binding and  no proper counseling to spend money remained in the indeterminate state that is nothing can be said about their family ties or the knowledge of spending earned money property but certainly they are unaware of HIV/AIDS and CSWs was their only recreation came to on state. Nothing could be said about their age criteria or higher motivation in life.

Now consider the effect of the state vector V = (0 0 0 0 0 1 0 1 0 0) on the neuronal dynamical system $N(E_3)$ i.e., only the nodes $A_6$ and  $A_8$ are in the on state and all other nodes are in the off state that is bad habits/bad company and no proper counseling of how to spend properly the money they earn are in the on state. The effect of V on $N(E_3)$ is

$VN (E_3)$    ↪    $(0\ 0\ 0\ 0\ 0\ 1\ 0\ 1\ 0\ 0)$    =    $V_1 = V$.

Thus we see both states happen to be the fixed point. Thus the on state of these two vectors have nil effect on the system. Next we consider the on state of the node $A_2$  and all other nodes are in the off state. The effect of the state vector T = (0 1 0 0 0 0 0 0 0 0) on the neuronal dynamical system N $(E_3)$

$TN (E_3)$    ↪    $(0\ 1\ 1\ 0\ I\ 0\ 0\ 1\ 0\ 0)$    =    $T_1$ (say)

$T_1N(E_3)$    ↪    $(0\ 1\ 1\ 1\ I\ 0\ 0\ 1\ 0\ 0)$    =    $T_2$ say

$T_2N(E_3)$    ↪    $(0\ 1\ 1\ 1\ I\ I\ 1\ 1\ 1\ 1)$    =    $T_3$ (say)

$T_3N(E_3)$    ↪    $(0\ 1\ 1\ 1\ I\ I\ 1\ 1\ I\ I)$    =    $T_4$ (say)

$T_4N(E_3)$    ↪    $(0\ 1\ 1\ 1\ I\ I\ 1\ 1\ I\ I)$    =    $T_5 (= T_4)$

a fixed point of the neutrosophic dynamical system). Suppose for in the study of the HIV/AIDS affected migrant labourers only the node living away from the family is in the on state and all other nodes are in the off state. We see in the resultant state vector the nodes $A_3$, $A_4$, $A_7$ and $A_8$ become on state $A_1$ remaining in the off state. The vectors $A_5$, $A_6$, $A_9$ and $A_{10}$ are indeterminates, conclusions are given for this resultant vector earlier. Thus we see whether they are free from being watched or live away from the family the net result in their activities is the same.

Next we give the fourth experts opinion on the set of 10 attributes $A_1$, $A_2$,…, $A_9$, $A_{10}$. The neutrosophic directed graph is as follows:

FIGURE: 4.3.4

The related neutrosophic relational matrix N ($E_4$) is given below.

$$
\begin{array}{c|cccccccccc}
 & A_1 & A_2 & A_3 & A_4 & A_5 & A_6 & A_7 & A_8 & A_9 & A_{10} \\
\hline
A_1 & 0 & 0 & 0 & 1 & 0 & 0 & 0 & 1 & 0 & 0 \\
A_2 & 0 & 0 & 0 & I & 0 & 0 & 0 & 0 & 0 & 0 \\
A_3 & 0 & 0 & 0 & 0 & 0 & 0 & 0 & 0 & 1 & 0 \\
A_4 & 0 & 0 & 0 & 0 & 1 & 0 & 0 & 0 & 0 & 0 \\
A_5 & 0 & 0 & 0 & 1 & 0 & 0 & 0 & 1 & 0 & 0 \\
A_6 & 0 & 0 & 0 & 0 & 0 & 0 & 0 & 0 & 0 & I \\
A_7 & 0 & 0 & 0 & 0 & 0 & 0 & 0 & 1 & 0 & 0 \\
A_8 & 0 & 0 & 0 & 0 & 0 & 0 & 0 & 0 & 0 & 0 \\
A_9 & 0 & 0 & 0 & 0 & 0 & 0 & 0 & 1 & 0 & 1 \\
A_{10} & 0 & 0 & 0 & 0 & 0 & 0 & 0 & 0 & 0 & 0 \\
\end{array}
$$

Consider the state vector P = (0 0 1 0 0 0 0 0 0 0) where only the concept $A_3$ is in the on state and all other nodes are in the off state that is free from being observed by friends and relatives is in the on state. The effect of the state vector P on the neutrosophic dynamical system $N(E_4)$

$$P \, N(E_4) \quad \hookrightarrow \quad (0\ 0\ 1\ 0\ 0\ 0\ 0\ 0\ 1\ 0) \quad = \quad P_1 \text{ (say)}$$

$$P_1 N(E_4) \quad \hookrightarrow \quad (0\ 0\ 1\ 0\ 0\ 0\ 0\ 1\ 1\ 1) \quad = \quad P_2 \text{ (say)}$$

$$P_2 \, N(E_4) \quad \hookrightarrow \quad (0\ 0\ 1\ 0\ 0\ 0\ 0\ 1\ 1\ 1) \quad = \quad P_3 = (P_2)$$

Thus we get the fixed point of the neutrosophic dynamical system. We see the nodes $A_1$, $A_2$, $A_4$, $A_5$, $A_6$ and $A_7$ are in the on state. According to this expert $A_8$, $A_9$ and $A_{10}$ become on there by indicating that CSWs are only recreation, no higher motivation in life and age stands an catalyst for it.

Next we consider the effect of the state vector Q = (0 0 0 0 0 1 0 1 0 0). The effect of Q on the dynamical system $N(E_4)$

$$QN(E_4) \quad \hookrightarrow \quad (0\ 0\ 0\ 0\ 0\ 1\ 0\ 1\ 0\ I) \quad = \quad Q_1 \text{ (say)}$$

$$Q_1 N(E_4) \quad \hookrightarrow \quad (0\ 0\ 0\ 0\ 0\ 1\ 0\ 1\ 0\ I) \quad = \quad Q_2 = Q$$
$$\text{(fixed point)}$$

Thus we see when $A_6$ and $A_8$ are in the on state the resultant vector the question of whether age is an disadvantage becomes an indeterminacy.

Now consider the on state of the attribute $A_1$ alone and all other nodes are off. The state vector R = (1 0 0 0 0 0 0 0 0 0) that is only easy money node is alone in the on state and all other nodes are off. The effect of R on the dynamical system $N(E_4)$ is

$$RN(E_4) \quad \hookrightarrow \quad (1\ 0\ 0\ 1\ 0\ 0\ 0\ 1\ 0\ 0) \quad = \quad R_1$$

$$R_1 N(E_4) \quad \hookrightarrow \quad (1\ 0\ 0\ 1\ 1\ 0\ 0\ 1\ 0\ 0) \quad = \quad R_2$$

$$R_2 N(E_4) \quad \hookrightarrow \quad (1\ 0\ 0\ 1\ 1\ 0\ 0\ 1\ 0\ 0) \quad = \quad R_3 (R_2 \text{ say})$$

Fixed point of the dynamical system. We see the nodes $A_4$, $A_5$ and $A_8$ become on that is when they have easy money they have lack of family binding, no proper counseling to spend/save money and visiting CSWs is a recreation according to this expert.

Next consider the fifth experts opinion on the 10 attributes $A_1$, $A_2$, …, $A_{10}$.

The neutrosophic directed graph given by him.

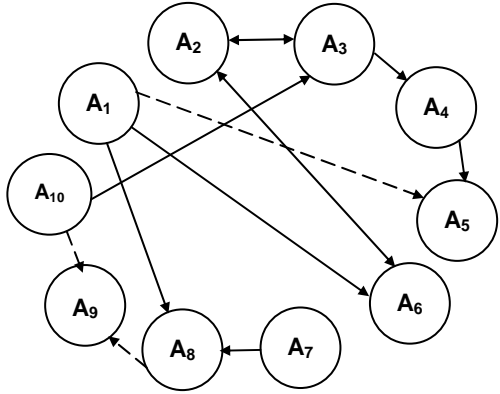

FIGURE: 4.3.5

The related neutrosophic connection matrix N (E$_5$)

$$
\begin{array}{c}
 \\
A_1 \\
A_2 \\
A_3 \\
A_4 \\
A_5 \\
A_6 \\
A_7 \\
A_8 \\
A_9 \\
A_{10}
\end{array}
\begin{array}{c}
A_1\ A_2\ A_3\ A_4\ A_5\ A_6\ A_7\ A_8\ A_9\ A_{10} \\
\left[
\begin{array}{cccccccccc}
0 & 0 & 0 & 0 & I & 1 & 0 & 1 & 0 & 0 \\
0 & 0 & 1 & 0 & 0 & 1 & 0 & 0 & 0 & 0 \\
0 & 1 & 0 & 0 & 0 & 0 & 0 & 0 & 0 & 0 \\
0 & 0 & 0 & 0 & 1 & 0 & 0 & 0 & 0 & 0 \\
0 & 0 & 0 & 0 & 0 & 0 & 0 & 0 & 0 & 0 \\
0 & 1 & 0 & 0 & 0 & 0 & 0 & 0 & 0 & 0 \\
0 & 0 & 0 & 0 & 0 & 0 & 0 & 1 & 0 & 0 \\
0 & 0 & 0 & 0 & 0 & 0 & 0 & 0 & I & 0 \\
0 & 0 & 0 & 0 & 0 & 0 & 0 & 0 & 0 & 0 \\
0 & 0 & 1 & 0 & 0 & 0 & 0 & 0 & I & 0
\end{array}
\right]
\end{array}
$$

Consider the state vector X in which A$_3$ that is free from being watched alone is in the on state and all other nodes are in the off state. Let X = (0 0 1 0 0 0 0 0 0 0). The effect of the state vector X on the neuronal dynamical system N (E$_5$).

$$XN(E_5) \quad \hookrightarrow \quad (0\ 1\ 1\ 0\ 0\ 0\ 0\ 0\ 0\ 0) \quad = \quad X_1 \text{ (say)}$$

$$X_1N(E_5) \quad \hookrightarrow \quad (0\ 1\ 1\ \ 0\ 0\ 1\ 0\ 0\ 0\ 0) \quad = \quad X_2 \text{ (say)}$$

$$X_2N\,(E_5) \quad \hookrightarrow \quad (0\ 1\ 1\ 0\ 0\ 1\ 0\ 0\ 0\ 0) \quad = \quad X_3 \ = \ (X_2)$$

X$_2$ is a fixed point. In the resultant vector only the nodes A$_2$ and A$_6$ becomes on state and all other nodes remain off, there by indicating that they live away from family and have no proper counseling to spend the money they earn. Now we analyze the effect of the state vector Y in which the attributes $A_6$ and $A_7$ are in the on state and all other nodes are in the off state i.e. Y = 90 0 0 0 0 1 1 0 0 0). The effect of state vector Y on the dynamical system $N(E_5)$

$$YN(E_5) \quad \hookrightarrow \quad (0\ 1\ 0\ 0\ 0\ 1\ 0\ 1\ 0\ 0) \quad = \quad Y_1 \text{ (say)}$$

$$Y_1N(E_5) \quad \hookrightarrow \quad (0\ 1\ 1\ 0\ 0\ 1\ 0\ 1\ I\ 0) \quad = \quad Y_2 \text{ (say)}$$

$$Y_2N(E_5) \quad \hookrightarrow \quad (0\ 1\ 1\ 0\ 0\ 1\ 0\ 1\ I\ 0) \quad = \quad Y_3 = Y_2$$

a fixed point of the dynamical system. Thus we see the resultant vector had made the nodes $A_2$, $A_3$ to be in the on state and the $A_9$ node is an indeterminate and all the other nodes $A_1$, $A_4$, $A_5$, $A_7$ and $A_{10}$ are in the off state; their by indicating that with no proper counseling of spending the money they earn coupled with CSWs their only recreation tells that they are away from the family with no fear of being observed by friends and relatives however the node no proper motivation/no higher aims in life remain as an indeterminate state.

We have obtained the effect of the individual neutrosophic cognitive maps now we find the effect of this 5 NCMs given by the five experts where $N(E) = N(E_1) + N(E_2) + N(E_3) + N(E_4) + N(E_5)$ and we study the effect of each of the state vectors and draw conclusions based on our analysis

|        | $A_1$ | $A_2$ | $A_3$ | $A_4$ | $A_5$ | $A_6$ | $A_7$ | $A_8$ | $A_9$ | $A_{10}$ |
|--------|-------|-------|-------|-------|-------|-------|-------|-------|-------|----------|
| $A_1$  | 0 | 0 | 0 | 3 | I | 3 | 0 | 2 | 0 | 0 |
| $A_2$  | I | 0 | 4 | I | I | 1 | I | 2 | 0 | 0 |
| $A_3$  | 0 | 2 | 0 | 3 | 0 | 0 | 0 | 0 | 2 | 1 |
| $A_4$  | I | 0 | 0 | 0 | 3 | I | 1 | 1 | 2 | 2 |
| $A_5$  | 0 | I | 0 | 1 | 0 | 0 | 0 | 1 | 0 | I |
| $A_6$  | 0 | 1 | 0 | 0 | 0 | 0 | 1 | 0 | 0 | I |
| $A_7$  | 0 | 0 | 0 | 0 | 0 | 0 | 0 | 3 | 0 | 0 |
| $A_8$  | 0 | 0 | 0 | 0 | 0 | 0 | 1 | 0 | 1 | 0 |
| $A_9$  | 0 | 0 | 0 | 0 | 0 | 0 | 0 | 1 | 0 | 1 |
| $A_{10}$ | 0 | 0 | 0 | 0 | 0 | 0 | 0 | 1 | I | 0 |

$N(E)$ is the combined neutrosophic connection matrix associated with the combined neutrosophic cognitive map. Now we proceed

on to study the effect of each and every state vector, which the expert is interested.

We make the modification in the thresholding of the state vector $(a_1, .., a_n)$ if $a_i \geq 2$ put 1, if $a_i < 2$ put 0. $a_i \geq 2I$ put I, $a_i < 2I$ put 0, $I + 1 = 0$, $2I + 1 = I$ and $2 + I = 1$. (The thresholding of the resultant vectors is sometimes based on the norms given by the expert). Suppose $C = (0\ 0\ 1\ 0\ 0\ 0\ 0\ 0\ 0\ 0)$ be the state vector in which only the attribute $A_3$ is in the on state and all other vectors are in the off state to study the effect on the dynamical system N(E).

$$CN(E) \hookrightarrow (0\ 2\ 1\ 3\ 0\ 0\ 0\ 0\ 2\ 1)$$
$$= (0\ 1\ 1\ 1\ 0\ 0\ 0\ 0\ 1\ 0) = C_1 \text{ (say)}$$
$$C_1N\ (E) \hookrightarrow (I\ 1\ 1\ 1\ 1\ 0\ 0\ 1\ 1\ 1) = C_2 \text{ (say)}$$
$$C_2N\ (E) \hookrightarrow (I\ 1\ 1\ 0\ 0\ I\ 1\ 1\ 1\ 1) = C_3 \text{ (say)}$$
$$C_3N\ (E) \hookrightarrow (I\ 1\ 1\ 0\ 0\ I\ 1\ 1\ 1\ 1) \text{ a fixed point}$$

which say $A_1$ and $A_6$ are indeterminates and $A_2$, $A_7$, $A_8$, $A_9$ and $A_{10}$ come to on state with $A_4$ and $A_5$ are in the off state. Thus we see easy money and no proper counseling of how to spend and save the earned money remains indeterminate state. Easy money and no channelizing of mind, lack of knowledge about HIV/AIDS, CSWs only recreation, no proper motivation for higher aims in life and age a big disadvantage come to on state.

We analyze different state vectors using the Combined Neutrosophic Cognitive Map (CNCM) and obtain conclusions based on our study.

## 4.4 Combined block NCMs and their applications

Now for the first time we define four new types of combined neutrosophic cognitive maps. These four techniques are mainly used when the number of attributes under consideration are large. Further it is difficult to obtain the neutrosophic directed graphs when the number of associated attributes is a big number. The directed graphs are difficult to be constructed using large number of attributes. So we in this book introduce the four new types of combined neutrosophic cognitive map techniques. First we define the new notion of combined disjoint block neutrosophic cognitive maps.

**DEFINITION 4.4.1:** *Let $A_1,..., A_n$ be n attributes for which the neutrosophic cognitive map is to be constructed (n large and a non prime) associated with some problem which is under investigation. For the combined disjoint block model we divide the n attributes into equal number of disjoint classes say $C_1,...,C_t$. i.e., t/ n. Now each $C_i$ contains n/ t number of chosen $A_t$; $1 \le t \le n$. (i = 1, 2,..., t). To each of the classes $C_i$ we using the experts opinion obtain the neutrosophic directed graph. Using these neutrosophic directed graphs we obtain their related neutrosophic connection matrices. This matrices are arranged in such a way we obtain the $n \times n$ connection matrix. This connection matrix represents the matrix associated with the Combined Disjoint Block Neutrosophic Cognitive Maps (CDBNCM).*

*It is important and interesting to note that when the division of n attributes $A_1,..., A_n$ are arbitrary but disjoint into say t classes say $C_1, C_2, ..., C_t$ where $C_i$ may have different number of elements in them, where even n can be a prime number then we call such structure also as combined neutrosophic block disjoint cognitive maps here only the blocks need not contain equal number of elements.*

We illustrate as well as exhibit these models in the study of the real world problem related with the HIV/AIDS migrant labourers. Now we adopt this mathematical model in the analysis of the conceptual nodes $P_1, P_2,..., P_5$ given in chapter 2 page 57.

Suppose we study with equal classes of division. Say $C_1$ $C_2$ and $C_3$ where $C_1 = \{P_1, P_2, P_3, P_4, P_5\}$, $C_2 = \{P_6, P_7, P, P_9, P_{10}\}$ and $C_3 = \{P_{11}, P_{12}, P_{13} P_{14} P_5\}$.

This expert has already taken 15 attributes now we obtain the neutrosophic directed graph using experts opinions for the classes $C_1$, $C_2$ and $C_3$.

The neutrosophic directed graph given by the expert related with the attributes $P_1, P_2, P_3, P_4$ and $P_5$ is as follows:

The related neutrosophic connection matrix is

$$\begin{array}{c c c c c c} & P_1 & P_2 & P_3 & P_4 & P_5 \\ \begin{matrix} P_1 \\ P_2 \\ P_3 \\ P_4 \\ P_5 \end{matrix} & \begin{bmatrix} 0 & 1 & 0 & I & 0 \\ 0 & 0 & 1 & I & 0 \\ 0 & 0 & 0 & 0 & 0 \\ 0 & 0 & 0 & 0 & 0 \\ 0 & 0 & 0 & 1 & 0 \end{bmatrix} \end{array}$$

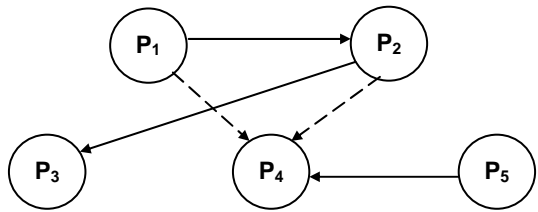

FIGURE: 4.4.1

Now we obtain the experts opinion on the attributes $P_6$, $P_7$, $P_8$, $P_9$ and $P_{10}$. The related neutrosophic directed graph is

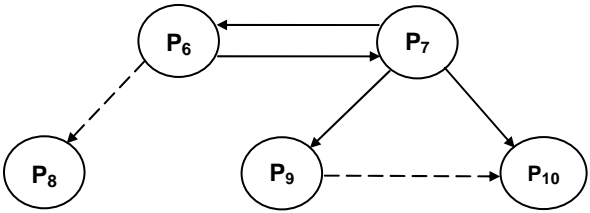

FIGURE: 4.4.2

The related neutrosophic connective matrix

$$\begin{array}{c c c c c c} & P_6 & P_7 & P_8 & P_9 & P_{10} \\ \begin{matrix} P_6 \\ P_7 \\ P_8 \\ P_9 \\ P_{10} \end{matrix} & \begin{bmatrix} 0 & 1 & I & 0 & 0 \\ 1 & 0 & 0 & 1 & 1 \\ 0 & 0 & 0 & 0 & 0 \\ 0 & 0 & 0 & 0 & I \\ 0 & 0 & 0 & 0 & 0 \end{bmatrix} \end{array}$$

The directed graph associated with the attributes $P_{11}$ $P_{12}$ $P_{13}$ $P_{14}$ and $P_{15}$ as given by the expert is

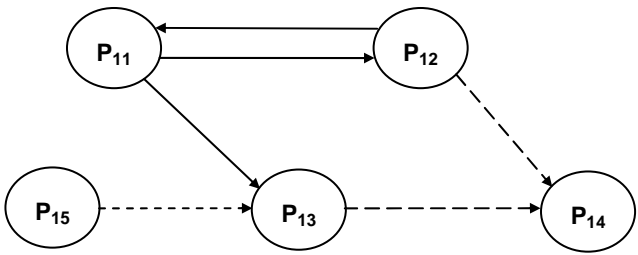

FIGURE: 4.4.3

The related neutrosophic connection matrix

$$\begin{array}{c c} & \begin{array}{c c c c c} P_{11} & P_{12} & P_{13} & P_{14} & P_{15} \end{array} \\ \begin{array}{c} P_{11} \\ P_{12} \\ P_{13} \\ P_{14} \\ P_{15} \end{array} & \left[ \begin{array}{c c c c c} 0 & 1 & 1 & 0 & 0 \\ 1 & 0 & 0 & I & 0 \\ 0 & 0 & 0 & I & 0 \\ 0 & 0 & 0 & 0 & 0 \\ 0 & 0 & I & 0 & 0 \end{array} \right] \end{array}$$

Now we give the related combined block disjoint neutrosophic connection matrix N (A) of the CDBNCM.

$$\begin{array}{c c} & \begin{array}{c c c c c c c c c c c c c c c} P_1 & P_2 & P_3 & P_4 & P_5 & P_6 & P_7 & P_8 & P_9 & P_{10} & P_{11} & P_{12} & P_{13} & P_{14} & P_{15} \end{array} \\ \begin{array}{c} P_1 \\ P_2 \\ P_3 \\ P_4 \\ P_5 \\ P_6 \\ P_7 \\ P_8 \\ P_9 \\ P_{10} \\ P_{11} \\ P_{12} \\ P_{13} \\ P_{14} \\ P_{15} \end{array} & \left[ \begin{array}{c c c c c c c c c c c c c c c} 0 & 1 & 0 & I & 0 & 0 & 0 & 0 & 0 & 0 & 0 & 0 & 0 & 0 & 0 \\ 0 & 0 & 1 & I & 0 & 0 & 0 & 0 & 0 & 0 & 0 & 0 & 0 & 0 & 0 \\ 0 & 0 & 0 & 0 & 0 & 0 & 0 & 0 & 0 & 0 & 0 & 0 & 0 & 0 & 0 \\ 0 & 0 & 0 & 0 & 0 & 0 & 0 & 0 & 0 & 0 & 0 & 0 & 0 & 0 & 0 \\ 0 & 0 & 0 & 1 & 0 & 0 & 0 & 0 & 0 & 0 & 0 & 0 & 0 & 0 & 0 \\ 0 & 0 & 0 & 0 & 0 & 0 & 1 & I & 0 & 0 & 0 & 0 & 0 & 0 & 0 \\ 0 & 0 & 0 & 0 & 0 & 1 & 0 & 0 & 1 & 1 & 0 & 0 & 0 & 0 & 0 \\ 0 & 0 & 0 & 0 & 0 & 0 & 0 & 0 & 0 & 0 & 0 & 0 & 0 & 0 & 0 \\ 0 & 0 & 0 & 0 & 0 & 0 & 0 & 0 & 0 & I & 0 & 0 & 0 & 0 & 0 \\ 0 & 0 & 0 & 0 & 0 & 0 & 0 & 0 & 0 & 0 & 0 & 0 & 0 & 0 & 0 \\ 0 & 0 & 0 & 0 & 0 & 0 & 0 & 0 & 0 & 0 & 0 & 1 & 1 & 0 & 0 \\ 0 & 0 & 0 & 0 & 0 & 0 & 0 & 0 & 0 & 0 & 1 & 0 & 0 & I & 0 \\ 0 & 0 & 0 & 0 & 0 & 0 & 0 & 0 & 0 & 0 & 0 & 0 & 0 & I & 0 \\ 0 & 0 & 0 & 0 & 0 & 0 & 0 & 0 & 0 & 0 & 0 & 0 & 0 & 0 & 0 \\ 0 & 0 & 0 & 0 & 0 & 0 & 0 & 0 & 0 & 0 & 0 & 0 & I & 0 & 0 \end{array} \right] \end{array}$$

Using N(A) the dynamical system of the given model we study the effect of state vectors.

Suppose S = (0 1 0 0 0 0 0 0 0 0 0 0 0 0 0) i.e., only the attribute male ego is in the on state and all other nodes are in the off state. To find the effect of S on the dynamical system N(A)

$$SN(A) \hookrightarrow (0\ 1\ 1\ I\ 0\ 0\ 0\ 0\ 0\ 0\ 0\ 0\ 0) = S_1 \text{ (say)}$$

$$S_1N(A) \hookrightarrow (0\ 1\ 1\ I\ 0\ 0\ 0\ 0\ 0\ 0\ 0\ 0\ 0) = S_2 = S_1 \text{ (say)}$$

Thus a person with male ego presumably treats women as inferior object and one cannot say then they have bad company or bad habits. The node $A_4$ remains as an indeterminate and all other nodes remains in the off state. Suppose we consider the state vector $Y = (1\ 0\ 0\ 0\ 0\ 0\ 1\ 0\ 0\ 0\ 1\ 0\ 0\ 0)$ that is the nodes $A_1$, $A_7$ and $A_{11}$ are in the on state and all other nodes are in the off state. The effect of $Y$ on the dynamical system $N(A)$

$$YN(A) \hookrightarrow (1\ 1\ 0\ I\ 0\ 1\ 1\ 0\ 1\ 1\ \ 1\ 0\ 0\ 0\ 0) = Y_1 \text{(say)}$$

$$Y_1N(A) \hookrightarrow (1\ 1\ 1\ I\ 0\ 1\ 1\ I\ 1\ 0\ 1\ 1\ 0\ 0\ ) = Y_2 \text{(say)}$$

$$Y_2N(A) \hookrightarrow (1\ 1\ 1\ I\ 0\ 1\ 1\ I\ 1\ 0\ 1\ 1\ 1\ I\ 0) = Y_3 \text{(say)}$$

$$Y_3N(A) \hookrightarrow (1\ 1\ 1\ I\ 0\ 1\ 1\ I\ 1\ 0\ 1\ 1\ 1\ I\ 0) = Y_4 = Y_3 \text{(say).}$$

The hidden pattern of the dynamical system is a fixed point in which the nodes $A_2$, $A_3$, $A_6$, $A_9$, $A_{12}$ and $A_{13}$ come to on state. The nodes $A_5$, $A_{10}$ and $A_{15}$ are 0 so they are in the off state. The nodes $A_4$ and $A_8$ are in the indeterminate state. Thus when no binding with family, gluttony food habits and visits CSWs the patient is a male egoist, who treat women as inferior objects, they are men with unusual desire for sex, with only physical activity, wants to enjoy life and be in jolly mood, he has no close friends/ relatives around him. One cannot say whether they enjoy bad company or they have more leisure or they are sexual perverts.

Thus we see that the combined Neutrosophic disjoint block cognitive maps works well when we consider different combination of attributes from different class as the state vector. This would be still more effective if the classes we reshuffled and the maps and directed graphs and obtain their related matrices are added. Now we apply it so another model given in page 42 to 52 of chapter 2.

Let $A_1$, $A_2$,…, $A_{12}$ be the 12 attributes related with a HIV/AIDS migrant labourer and the social and economic problem he faces. Let us divide these 12 attributes in classes so that each class has 3 attributes $C_1 = \{A_1\ A_2\ A_3\}$, $C_2 = \{A_4,\ A_5\ A_6\}$, $C_3 = \{A_7\ A_8\ A_9\}$ and $C_4 = \{A_{10}\ A_{11}\ A_{12}\}$.

Now we obtain an expert's opinion on each of these classes in the form of directed graph which is transformed on to a connection matrix.

The neutrosophic directed graph as given by the expert for the attributes $A_1$, $A_2$, $A_3$ from the class $C_1$.

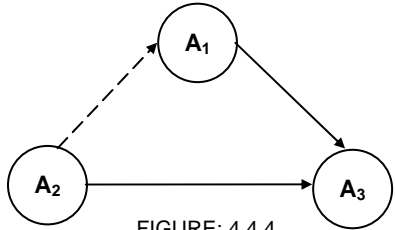

FIGURE: 4.4.4

The related connection matrix

$$
\begin{array}{c}
\phantom{A_1} \begin{array}{ccc} A_1 & A_2 & A_3 \end{array} \\
\begin{array}{c} A_1 \\ A_2 \\ A_3 \end{array}
\begin{bmatrix}
0 & 0 & 1 \\
I & 0 & 1 \\
0 & 0 & 0
\end{bmatrix}
\end{array}
$$

The neutrosophic directed graph as given by the expert for the class $C_2 = \{A_4 \ A_5 \ A_6\}$.

The related connection neutrosophic matrix is

$$
\begin{array}{c}
\phantom{A_4} \begin{array}{ccc} A_4 & A_5 & A_6 \end{array} \\
\begin{array}{c} A_4 \\ A_5 \\ A_6 \end{array}
\begin{bmatrix}
0 & I & I \\
0 & 0 & 1 \\
0 & 0 & 0
\end{bmatrix}
\end{array}
$$

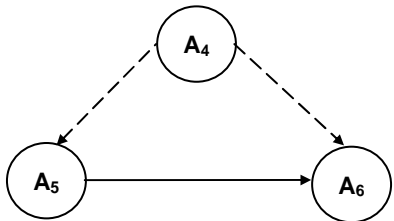

FIGURE: 4.4.5

The neutrosophic directed graph related to the class $C_3 = \{A_7, A_8, A_9\}$ given by the expert.

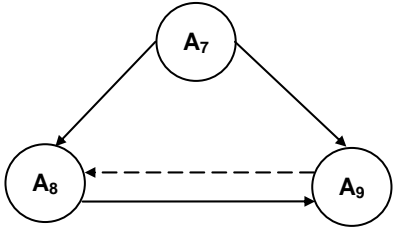

FIGURE: 4.4.6

The related neutrosophic connection matrix

$$
\begin{array}{c}
\quad\quad A_7\, A_8\, A_9 \\
\begin{array}{c} A_7 \\ A_8 \\ A_9 \end{array}
\left[\begin{array}{ccc}
0 & 1 & 1 \\
0 & 0 & 1 \\
0 & I & 0
\end{array}\right]
\end{array}
$$

The directed neutrosophic graph associated with the class $C_4 = \{A_{10}\ A_{11}\ A_{12}\}$ given by the expert is as follows:

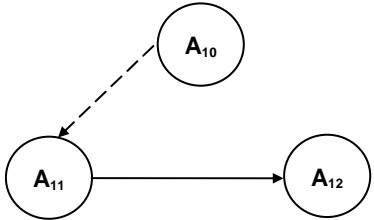

FIGURE: 4.4.7

The related neutrosophic connection matrix is given below

$$
\begin{array}{c}
\quad\quad A_{10}\, A_{11}\, A_{12} \\
\begin{array}{c} A_{10} \\ A_{11} \\ A_{12} \end{array}
\left[\begin{array}{ccc}
0 & I & 0 \\
0 & 0 & 1 \\
0 & 0 & 0
\end{array}\right]
\end{array}
$$

Now we find the combined disjoint block matrix associated with the classes $C_1$, $C_2$, $C_3$ and $C_4$. Let N (B) denote the associated neutrosophic $12 \times 12$ matrix related with the attributes $A_1$, $A_2$,…, $A_{12}$.

$$
\begin{array}{c}
\quad\ A_1\ A_2\ A_3\ A_4\ \ A_5\ A_6\ A_7\ A_8\ A_9\ A_{10}\ A_{11}\ A_{12} \\
\begin{array}{c}
A_1 \\ A_2 \\ A_3 \\ A_4 \\ A_5 \\ A_6 \\ A_7 \\ A_8 \\ A_9 \\ A_{10} \\ A_{11} \\ A_{12}
\end{array}
\left[
\begin{array}{cccccccccccc}
0 & 0 & 1 & 0 & 0 & 0 & 0 & 0 & 0 & 0 & 0 & 0 \\
I & 0 & 1 & 0 & 0 & 0 & 0 & 0 & 0 & 0 & 0 & 0 \\
0 & 0 & 0 & 0 & 0 & 0 & 0 & 0 & 0 & 0 & 0 & 0 \\
0 & 0 & 0 & 0 & I & I & 0 & 0 & 0 & 0 & 0 & 0 \\
0 & 0 & 0 & 0 & 0 & 1 & 0 & 0 & 0 & 0 & 0 & 0 \\
0 & 0 & 0 & 0 & 0 & 0 & 0 & 0 & 0 & 0 & 0 & 0 \\
0 & 0 & 0 & 0 & 0 & 0 & 0 & 1 & 1 & 0 & 0 & 0 \\
0 & 0 & 0 & 0 & 0 & 0 & 0 & 0 & 1 & 0 & 0 & 0 \\
0 & 0 & 0 & 0 & 0 & 0 & 0 & I & 0 & 0 & 0 & 0 \\
0 & 0 & 0 & 0 & 0 & 0 & 0 & 0 & 0 & 0 & I & 0 \\
0 & 0 & 0 & 0 & 0 & 0 & 0 & 0 & 0 & 0 & 0 & 1 \\
0 & 0 & 0 & 0 & 0 & 0 & 0 & 0 & 0 & 0 & 0 & 0
\end{array}
\right]
\end{array}
$$

Now let us find the effect of the state vector (0 0 1 0 1 0 0 0 1 1 0 0) = X on the dynamical system N (B) that is when the nodes $A_3$, $A_5$, $A_9$ and $A_{10}$ are in the on states and all other nodes are in the off state. Consider

$$XN (B) \hookrightarrow (0\,0\,1\,0\,1\,1\,0\,0\,1\,1\,I\,0) \quad = \quad X_1 \text{ (say)}$$

$$X_1N(B) \hookrightarrow (0\,0\,1\,0\,1\,0\,0\,I\,1\,1\,I\,0) \quad = \quad X_2 \text{ (say)}$$

$$X_2N(B) \hookrightarrow (0\,0\,1\,0\,1\,0\,0\,I\,1\,1\,I\,I) \quad = \quad X_3 \text{ (say)}$$

$$X_3N(B) \hookrightarrow (0\,0\,1\,0\,1\,0\,0\,I\,1\,1\,I\,I) \quad = \quad X_3 \text{ (fixed point)}$$

$X_3$ is a fixed point. Thus we see the effect this state vector on the dynamical system that is the resultant state vector is such that the nodes $A_1$, $A_2$, $A_4$, $A_6$ and $A_7$ are in the off state. The nodes $A_8$, $A_{11}$ and $A_{12}$ are in the indeterminate state. Apart from this no node comes to on state a very rare situation.

Now we illustrate for the same two models by the combined disjoint block neutrosophic cognitive models where now the block length are not all equal. Each of the block length need not be the same in general. Now we adopt this model in the study of the migrant labourers affected by HIV/AIDS and their related socio economic problems. The 15 attributes $P_1$, $P_2$,…, $P_{15}$ given in page 57 of chapter two is taken.

Now let us consider the classes $C_1 = \{P_1\ P_2\ P_3\ P_4\} = C_2 = \{P_5\ P_6,\ P_7\}$, $C_3 = \{P_8\ P_9\ P_{10}\ P_{11}\ P_{12}\}$ and $C_4 = \{P_{13}\ P_{14},\ P_{15}\}$.

Now using the experts opinion we give the related directed graph. The neutrosophic directed graph related with the attributes $C_1 = \{P_1, P_2, P_3, P_4\}$.

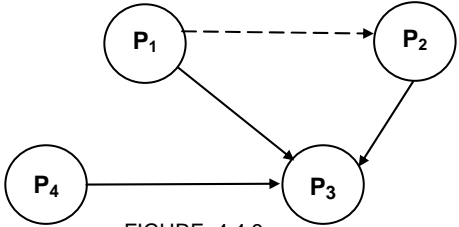

FIGURE: 4.4.8

The related connection matrix.

$$
\begin{array}{c c}
 & \begin{array}{c c c c} P_1 & P_2 & P_3 & P_4 \end{array} \\
\begin{array}{c} P_1 \\ P_2 \\ P_3 \\ P_4 \end{array} &
\left[\begin{array}{c c c c}
0 & I & 1 & 0 \\
0 & 0 & 1 & 0 \\
0 & 0 & 0 & 0 \\
0 & 0 & 1 & 0
\end{array}\right]
\end{array}
$$

Now neutrosophic directed graph by the expert related with the class $C_2 = \{P_5, P_6, P_7\}$.

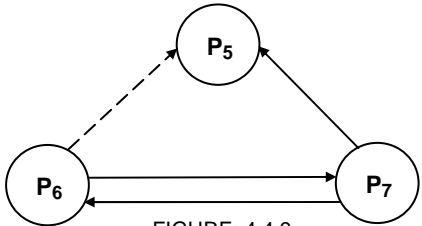

FIGURE: 4.4.9

The related neutrosophic connection matrix

$$
\begin{array}{c c}
 & \begin{array}{c c c} P_5 & P_6 & P_7 \end{array} \\
\begin{array}{c} P_5 \\ P_6 \\ P_7 \end{array} &
\left[\begin{array}{c c c}
0 & 0 & 0 \\
I & 0 & 1 \\
1 & 1 & 0
\end{array}\right]
\end{array}
$$

For the class $C_3 = \{P_8\ P_9\ P_{10}\ P_{11}\ P_{12}\}$ we have the experts opinion given by the following neutrosophic directed graph.

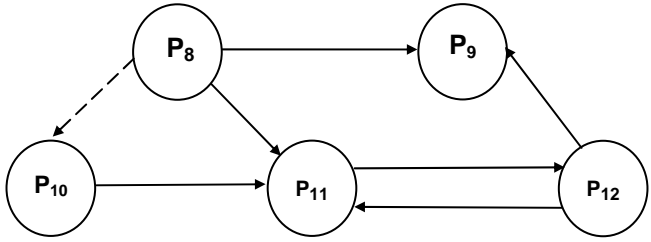

FIGURE: 4.4.10

The related neutrosophic connection matrix.

$$\begin{array}{c} \\ P_8 \\ P_9 \\ P_{10} \\ P_{11} \\ P_{12} \end{array} \begin{array}{ccccc} P_8 & P_9 & P_{10} & P_{11} & P_{12} \\ \left[\begin{array}{ccccc} 0 & 1 & I & 1 & 0 \\ 0 & 0 & 0 & 0 & 0 \\ 0 & 0 & 0 & 1 & 0 \\ 0 & 0 & 0 & 0 & 1 \\ 0 & 1 & 0 & 1 & 0 \end{array}\right] \end{array}$$

Now we give the neutrosophic directed graph of the class $C_5$ as given by the expert.

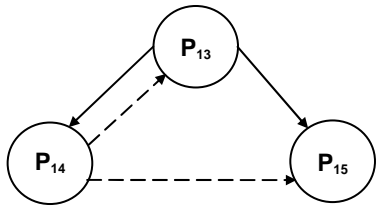

FIGURE: 4.4.11

The corresponding neutrosophic connection matrix

$$\begin{array}{c} \\ P_{13} \\ P_{14} \\ P_{15} \end{array} \begin{array}{ccc} P_{13} & P_{14} & P_{15} \\ \left[\begin{array}{ccc} 0 & 1 & 1 \\ I & 0 & I \\ 0 & 0 & 0 \end{array}\right] \end{array}$$

Now we using the 5 neutrosophic relation matrices obtain the related combined block disjoint neutrosophic cognitive maps matrix of different length. Let us denote it by N(D)

|        | $P_1$ | $P_2$ | $P_3$ | $P_4$ | $P_5$ | $P_6$ | $P_7$ | $P_8$ | $P_9$ | $P_{10}$ | $P_{11}$ | $P_{12}$ | $P_{13}$ | $P_{14}$ | $P_{15}$ |
|--------|----|----|----|----|----|----|----|----|----|----|----|----|----|----|----|
| $P_1$    | 0 | $I$ | 1 | 0 | 0 | 0 | 0 | 0 | 0 | 0 | 0 | 0 | 0 | 0 | 0 |
| $P_2$    | 0 | 0 | 1 | 0 | 0 | 0 | 0 | 0 | 0 | 0 | 0 | 0 | 0 | 0 | 0 |
| $P_3$    | 0 | 0 | 0 | 0 | 0 | 0 | 0 | 0 | 0 | 0 | 0 | 0 | 0 | 0 | 0 |
| $P_4$    | 0 | 0 | 1 | 0 | 0 | 0 | 0 | 0 | 0 | 0 | 0 | 0 | 0 | 0 | 0 |
| $P_5$    | 0 | 0 | 0 | 0 | 0 | 0 | 0 | 0 | 0 | 0 | 0 | 0 | 0 | 0 | 0 |
| $P_6$    | 0 | 0 | 0 | 0 | $I$ | 0 | 1 | 0 | 0 | 0 | 0 | 0 | 0 | 0 | 0 |
| $P_7$    | 0 | 0 | 0 | 0 | 1 | 1 | 0 | 0 | 0 | 0 | 0 | 0 | 0 | 0 | 0 |
| $P_8$    | 0 | 0 | 0 | 0 | 0 | 0 | 0 | 0 | 1 | $I$ | 1 | 0 | 0 | 0 | 0 |
| $P_9$    | 0 | 0 | 0 | 0 | 0 | 0 | 0 | 0 | 0 | 0 | 0 | 0 | 0 | 0 | 0 |
| $P_{10}$   | 0 | 0 | 0 | 0 | 0 | 0 | 0 | 0 | 0 | 0 | 1 | 0 | 0 | 0 | 0 |
| $P_{11}$   | 0 | 0 | 0 | 0 | 0 | 0 | 0 | 0 | 0 | 0 | 0 | 1 | 0 | 0 | 0 |
| $P_{12}$   | 0 | 0 | 0 | 0 | 0 | 0 | 0 | 0 | 1 | 0 | 1 | 0 | 0 | 0 | 0 |
| $P_{13}$   | 0 | 0 | 0 | 0 | 0 | 0 | 0 | 0 | 0 | 0 | 0 | 0 | 0 | 1 | 1 |
| $P_{14}$   | 0 | 0 | 0 | 0 | 0 | 0 | 0 | 0 | 0 | 0 | 0 | 0 | $I$ | 0 | $I$ |
| $P_{15}$   | 0 | 0 | 0 | 0 | 0 | 0 | 0 | 0 | 0 | 0 | 0 | 0 | 0 | 0 | 0 |

Consider the state vector $X = (0\ 0\ 1\ 0\ 0\ 1\ 0\ 0\ 1\ 0\ 0\ 0\ 1)$ i.e., the nodes $P_3$, $P_7$, $P_{11}$ and $P_{15}$ are in the on state.

The effect of X on the neutrosophic dynamical system N (D)

$$XN \quad \hookrightarrow \quad (0\ 0\ 1\ 0\ 1\ 1\ 1\ 0\ 0\ 0\ 1\ 1\ 0\ 0\ I) \quad = \quad X_1 \text{ (say)}$$

$$X_1 N(D) \hookrightarrow \quad (0\ 0\ 1\ 0\ 0\ 1\ 1\ 0\ 1\ 0\ 1\ 1\ 0\ 0\ 1) \quad = \quad X_2 \text{ (say)}$$

$$X_2 N(D) \hookrightarrow \quad (0\ 0\ 1\ 0\ 0\ 1\ 1\ 0\ 1\ 0\ 1\ 1\ 0\ 0\ 1) \quad = \quad X_3\ (=X_2)$$

$X_2$ is a fixed point. Thus the hidden pattern for this state vector is a fixed point.

We see when the nodes women as inferior objects, food habits, visits CSWs and failure of agriculture to be in the on state and all other nodes are in the off states we see that the nodes $P_1$, $P_2$, $P_4$, $P_5$, $P_8$, $P_{10}$, $P_{13}$ and $P_{14}$ remain in the off state there by indicating that no binding with family, male chauvinism, bad company and bad habits, socially irresponsible, more leisure, only physically active, unreachable by friends or relatives and sexual perverts remain unaffected only the nodes $P_6$ and $P_9$ come to on state; that is uncontrollable sexual feelings and no work for brain are the nodes which influence $P_3$, $P_7$, $P_{11}$ and $P_{15}$. One can use any state vector and derive the conclusion.

Now we model the 12 attributes given in page 42 to 53 the node / concepts related with an HIV/AIDS affected migrant labourer. Let $A_1$, $A_2$,…, $A_{12}$ be divided into classes of different cardinalities i.e., each block does not contain same number of attributes.

Let
$$C_1 = \{A_1\ A_2\ A_3\ A_4\ A_5\},$$
$$C_2 = \{A_6\ A_7\ A_8\}$$
and $$C_3 = \{A_9\ A_{10}\ A_{11}\ A_{12}\}.$$

Now using the opinion of the expert we give the neutrosophic directed graph and their related connection matrix.

The neutrosophic directed graph related with the class $C_1 = \{A_1, A_2\ A_3\ A_4\ A_5\}$ is an follows.

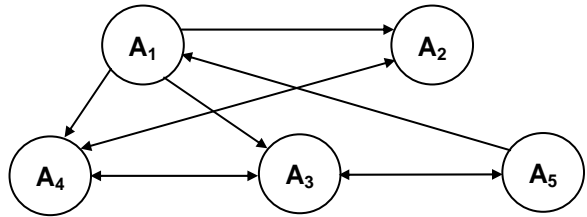

FIGURE: 4.4.12

The related connection matrix

$$
\begin{array}{c}
\quad\quad A_1\ A_2\ A_3\ A_4\ A_5 \\
\begin{array}{c}
A_1 \\ A_2 \\ A_3 \\ A_4 \\ A_5
\end{array}
\begin{bmatrix}
0 & I & 1 & I & 0 \\
0 & 0 & 0 & 1 & 0 \\
0 & 0 & 0 & 1 & 1 \\
0 & 1 & 1 & 0 & 0 \\
1 & 0 & 1 & 0 & 0
\end{bmatrix}
\end{array}
$$

Now we give the neutrosophic directed graph given by the expert for the class $C_2 = \{A_6\ A_7\ A_8\}$.

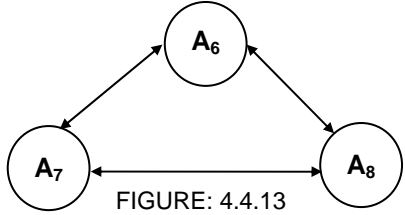

FIGURE: 4.4.13

The related connection matrix

$$\begin{array}{c}
\phantom{A_6}\begin{array}{ccc} A_6 & A_7 & A_8 \end{array} \\
\begin{array}{c} A_6 \\ A_7 \\ A_8 \end{array}
\begin{bmatrix} 0 & 1 & 1 \\ 1 & 0 & 1 \\ 1 & 1 & 0 \end{bmatrix}
\end{array}$$

Now the directed graph and the connection matrix for the class $C_3$ = ($A_9$ $A_{10}$ $A_{11}$ $A_{12}$) as given by the expert

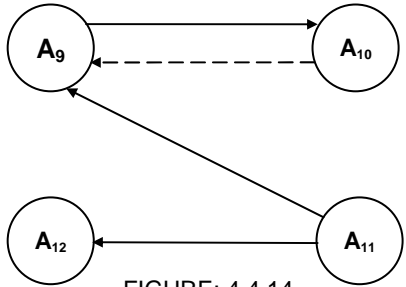

FIGURE: 4.4.14

The related connection matrix is as follows.

$$\begin{array}{c}
\phantom{A_9}\begin{array}{cccc} A_9 & A_{10} & A_{11} & A_{12} \end{array} \\
\begin{array}{c} A_9 \\ A_{10} \\ A_{11} \\ A_{12} \end{array}
\begin{bmatrix} 0 & 1 & 0 & 0 \\ I & 0 & 0 & 0 \\ 1 & 0 & 0 & 1 \\ 0 & 0 & 0 & 0 \end{bmatrix}
\end{array}$$

Now using these three connection matrices we give the related connection matrix of the combined block disjoint neutrosophic cognitive maps of varying length

$$\begin{array}{c@{\ }c}
& \begin{matrix} A_1 & A_2 & A_3 & A_4 & A_5 & A_6 & A_7 & A_8 & A_9 & A_{10} & A_{11} & A_{12} \end{matrix} \\
\begin{matrix} A_1 \\ A_2 \\ A_3 \\ A_4 \\ A_5 \\ A_6 \\ A_7 \\ A_8 \\ A_9 \\ A_{10} \\ A_{11} \\ A_{12} \end{matrix} &
\begin{bmatrix}
0 & I & 1 & I & 0 & 0 & 0 & 0 & 0 & 0 & 0 & 0 \\
0 & 0 & 0 & 1 & 0 & 0 & 0 & 0 & 0 & 0 & 0 & 0 \\
0 & 0 & 0 & 1 & 1 & 0 & 0 & 0 & 0 & 0 & 0 & 0 \\
0 & 1 & 1 & 0 & 0 & 0 & 0 & 0 & 0 & 0 & 0 & 0 \\
1 & 0 & 1 & 0 & 0 & 0 & 0 & 0 & 0 & 0 & 0 & 0 \\
0 & 0 & 0 & 0 & 0 & 0 & 1 & 1 & 0 & 0 & 0 & 0 \\
0 & 0 & 0 & 0 & 0 & 1 & 0 & 1 & 0 & 0 & 0 & 0 \\
0 & 0 & 0 & 0 & 0 & 1 & 1 & 0 & 0 & 0 & 0 & 0 \\
0 & 0 & 0 & 0 & 0 & 0 & 0 & 0 & 0 & 1 & 0 & 0 \\
0 & 0 & 0 & 0 & 0 & 0 & 0 & 0 & I & 0 & 0 & 0 \\
0 & 0 & 0 & 0 & 0 & 0 & 0 & 0 & 1 & 0 & 0 & 1 \\
0 & 0 & 0 & 0 & 0 & 0 & 0 & 0 & 0 & 0 & 0 & 0
\end{bmatrix}
\end{array}$$

Let $N(E)$ denotes the relational connection matrix. Let us consider the state vector $X = (0\ 1\ 0\ 0\ 0\ 0\ 1\ 0\ 0\ 0\ 1\ 0)$ that is $A_2$, $A_7$ and $A_{11}$ are in the on state and all other attributes are in the off state; the effect of $X$ on the neutrosophic dynamical system $N(E)$

$$\begin{aligned}
XN(E) &\hookrightarrow (0\ 1\ 0\ 1\ 0\ 1\ 1\ 1\ 1\ 0\ 11) \\
X_1N(E) &\hookrightarrow (0\ 1\ 1\ 1\ 0\ 1\ 1\ 1\ 1\ 0\ 1\ 1) = X_2 \text{ (say)} \\
X_2N(E) &\hookrightarrow (0\ 1\ 1\ 1\ 1\ 1\ 1\ 1\ 1\ 1\ 1\ 1) = X_3 \text{ say} \\
X_3N(E) &\hookrightarrow (1\ 1\ 1\ 1\ 1\ 1\ 1\ 1\ 1\ I\ 1\ 1) = X_4 \\
X_4N(E) &\hookrightarrow (1\ 1\ 1\ 1\ 1\ 1\ 1\ 1\ 1\ I\ I\ 1\ 1) = X_5 \\
X_5N(E) &\hookrightarrow (1\ 1\ 1\ 1\ 1\ 1\ 1\ 1\ 1\ I\ I\ 1\ 1) = X_6 = (X_5)
\end{aligned}$$

Thus the hidden pattern of the dynamical system is a fixed point. All the state vector become on when the attributes $A_2$, $A_7$ and $A_{11}$ are in the on state so these three nodes are sufficient to make all the nodes on.

Now we proceed onto define the two new notions about combined block over lap neutrosophic cognitive maps of same length and different length:

**DEFINITION 4.4.2:** *Consider a set of nodes / attributes $A_1$, $A_2$,..., $A_n$ to be analyzed using the combined overlap block neutrosophic maps using blocks of same length, let us divide the attributes $A_1$,*

*$A_2$..., $A_n$ in to m classes with a overlap of s elements ie $A_1$, $A_2$,..., $A_n$ is divided into classes $C_1$,..., $C_m$ where we have number of elements in each $C_i$ is s where there are common elements between the $C_i$'s and $C_{i+1}$'s. Now we using the experts opinion on each of these classes obtain the neutrosophic directed graph and their related neutrosophic connection matrix relative to COBNCM.*

*We arrange these neutrosophic connection matrices into a n × n combined neutrosophic connection matrix.*

*Now instead of making the number of elements in each $C_i$ to be the same we can have different number of elements in each $C_i$ and also for any distinct $C_i$ and $C_{i+1}$ we may have different number of elements in $C_i \cap C_{i+1}$. The neutrosophic directed graph is drawn for each of these $C_i$ and their corresponding neutrosophic connection matrix is obtained. These neutrosophic connection matrices are combined and a n × n neutrosophic matrix is obtained which is termed as the combined overlap block neutrosophic matrix relative to the combined overlap block neutrosophic cognitive maps (COBNCM) of varying block sizes.*

Now we illustrate these two models by the neutrosophic relation between the HIV/AIDS migrant labourers and their related socio economic problem. Let $P_1$, $P_2$,..., $P_{15}$ be the attributes given in chapter two in page 57. Now we divide them into classes $C_1$, ..., $C_5$ where

$C_1$ = {$P_1$ $P_2$ $P_3$ $P_4$ $P_5$},  $C_2$ = {$P_4$ $P_5$ $P_6$ $P_7$ $P_8$},
$C_3$ = {$P_7$ $P_8$ $P_9$ $P_{10}$ $P_{11}$},  $C_4$ = {$P_{10}$ $P_{11}$ $P_{12}$ $P_{13}$ $P_{14}$}
and $C_5$ = {$P_{13}$ $P_{14}$ $P_{15}$ $P_1$ $P_2$}.

Now using the experts opinion we obtain the neutrosophic directed graph for each of these classes. The neutrosophic directed graph for $C_1$ as given by the expert is

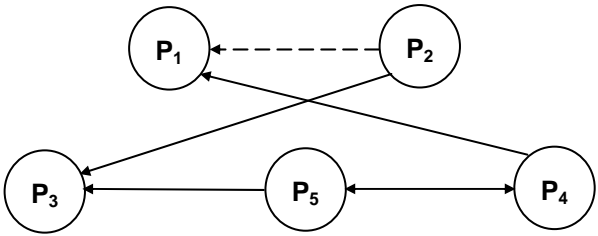

FIGURE: 4.4.15

The related connection matrix is

$$
\begin{array}{c}
\quad\quad P_1\ P_2\ P_3\ P_4\ P_5 \\
\begin{array}{c}
P_1 \\ P_2 \\ P_3 \\ P_4 \\ P_5
\end{array}
\begin{bmatrix}
0 & 0 & 0 & 0 & 0 \\
I & 0 & 1 & 0 & 0 \\
0 & 0 & 0 & 0 & 0 \\
1 & 0 & 0 & 0 & 1 \\
0 & 0 & 1 & 1 & 0
\end{bmatrix}
\end{array}
$$

The experts opinion on the class $C_2$ in the form of the neutrosophic directed graph

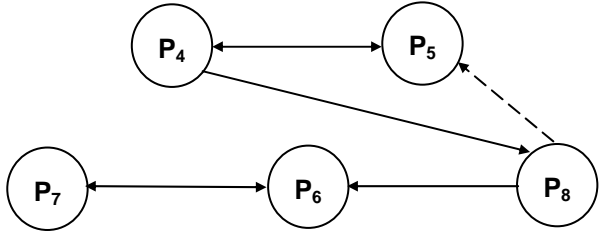

FIGURE: 4.4.16

The related connection matrix is

$$
\begin{array}{c}
\quad\quad P_4\ P_5\ P_6\ P_7\ P_8 \\
\begin{array}{c}
P_4 \\ P_5 \\ P_6 \\ P_7 \\ P_8
\end{array}
\begin{bmatrix}
0 & 1 & 0 & 0 & 1 \\
1 & 0 & 0 & 0 & 0 \\
0 & 0 & 0 & 1 & 0 \\
0 & 0 & 1 & 0 & 0 \\
0 & I & 1 & 0 & 0
\end{bmatrix}
\end{array}
$$

The experts opinion on the class $C_3 = \{P_7\ P_8\ P_9\ P_{10}\ P_{11}\}$ in the form of the neutrosophic directed graph

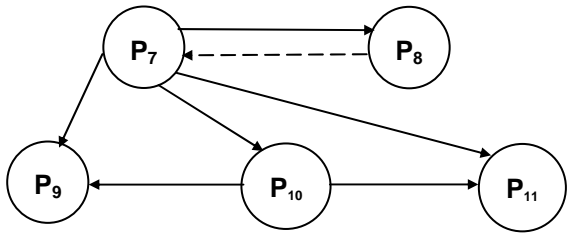

FIGURE: 4.4.17

The related neutrosophic connection matrix

$$
\begin{array}{c}
\phantom{P_{7}} \\
P_{7} \\
P_{8} \\
P_{9} \\
P_{10} \\
P_{11}
\end{array}
\begin{array}{c}
P_7 \ P_8 \ P_9 \ P_{10} \ P_{11} \\
\begin{bmatrix}
0 & 1 & 1 & 1 & 1 \\
I & 0 & 0 & 0 & 0 \\
0 & 0 & 0 & 0 & 0 \\
0 & 0 & 1 & 0 & 1 \\
0 & 0 & 0 & 0 & 0
\end{bmatrix}
\end{array}
$$

The directed neutrosophic graph related to the class $C_4 = \{P_{10} \ P_{11} \ P_{12} \ P_{13} \ P_{14}\}$.

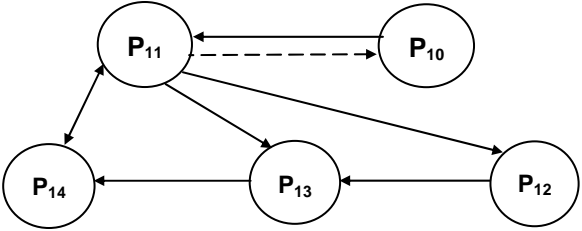

FIGURE: 4.4.18

The related neutrosophic connection matrix

$$
\begin{array}{c}
\phantom{P_{10}} \\
P_{10} \\
P_{11} \\
P_{12} \\
P_{13} \\
P_{14}
\end{array}
\begin{array}{c}
P_{10} \ P_{11} \ P_{12} \ P_{13} \ P_{14} \\
\begin{bmatrix}
0 & 1 & 0 & 0 & 0 \\
I & 0 & 1 & 1 & 1 \\
0 & 0 & 0 & 1 & 0 \\
0 & 0 & 0 & 0 & 1 \\
0 & 1 & 0 & 0 & 0
\end{bmatrix}
\end{array}
$$

Now we give the directed graph for $C_5 = \{P_{13} \ P_{14} \ P_{15} \ P_1 \ P_2\}$.

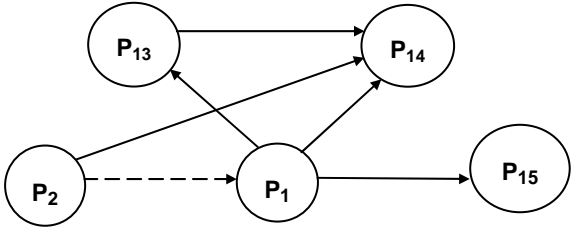

FIGURE: 4.4.19

The related connection matrix

$$
\begin{array}{c}
\phantom{P_{13}} \\
P_{13} \\
P_{14} \\
P_{15} \\
P_{1} \\
P_{2}
\end{array}
\begin{array}{c}
P_{13}\; P_{14}\; P_{15}\; P_{1}\; P_{2} \\
\begin{bmatrix}
0 & 1 & 0 & 0 & 0 \\
0 & 0 & 0 & 0 & 0 \\
0 & 0 & 0 & 0 & 0 \\
1 & 1 & 1 & 0 & 0 \\
0 & 1 & 0 & I & 0
\end{bmatrix}
\end{array}
$$

Now we give the combined connection matrix related with the combined overlap block NCM.

| | $P_1$ | $P_2$ | $P_3$ | $P_4$ | $P_5$ | $P_6$ | $P_7$ | $P_8$ | $P_9$ | $P_{10}$ | $P_{11}$ | $P_{12}$ | $P_{13}$ | $P_{14}$ | $P_{15}$ |
|---|---|---|---|---|---|---|---|---|---|---|---|---|---|---|---|
| $P_1$ | 0 | 0 | 0 | 0 | 0 | 0 | 0 | 0 | 0 | 0 | 0 | 0 | 1 | 1 | 1 |
| $P_2$ | $I$ | 0 | 1 | 0 | 0 | 0 | 0 | 0 | 0 | 0 | 0 | 0 | 0 | 1 | 0 |
| $P_3$ | 0 | 0 | 0 | 0 | 0 | 0 | 0 | 0 | 0 | 0 | 0 | 0 | 0 | 0 | 0 |
| $P_4$ | 1 | 0 | 0 | 0 | 2 | 0 | 0 | 1 | 0 | 0 | 0 | 0 | 0 | 0 | 0 |
| $P_5$ | 0 | 0 | 1 | 2 | 0 | 0 | 0 | 0 | 0 | 0 | 0 | 0 | 0 | 0 | 0 |
| $P_6$ | 0 | 0 | 0 | 0 | 0 | 0 | 1 | 0 | 0 | 0 | 0 | 0 | 0 | 0 | 0 |
| $P_7$ | 0 | 0 | 0 | 0 | 0 | 1 | 0 | 0 | 0 | 0 | 0 | 0 | 0 | 0 | 0 |
| $P_8$ | 0 | 0 | 0 | 0 | $I$ | 1 | $I$ | 0 | 0 | 0 | 0 | 0 | 0 | 0 | 0 |
| $P_9$ | 0 | 0 | 0 | 0 | 0 | 0 | 0 | 0 | 0 | 1 | 1 | 0 | 0 | 0 | 0 |
| $P_{10}$ | 0 | 0 | 0 | 0 | 0 | 0 | 0 | 0 | 1 | 0 | 2 | 0 | 0 | 0 | 0 |
| $P_{11}$ | 0 | 0 | 0 | 0 | 0 | 0 | 0 | 0 | 0 | $I$ | 0 | 1 | 1 | 1 | 0 |
| $P_{12}$ | 0 | 0 | 0 | 0 | 0 | 0 | 0 | 0 | 0 | $I$ | 0 | 0 | 1 | 0 | 0 |
| $P_{13}$ | 0 | 0 | 0 | 0 | 0 | 0 | 0 | 0 | 0 | 0 | 0 | 0 | 0 | 2 | 0 |
| $P_{14}$ | 0 | 0 | 0 | 0 | 0 | 0 | 0 | 0 | 0 | 0 | 0 | 0 | 0 | 0 | 0 |
| $P_{15}$ | 0 | 0 | 0 | 0 | 0 | 0 | 0 | 0 | 0 | 0 | 1 | 0 | 0 | 0 | 0 |

Let us denote this neutrosophic connection matrix for the combined overlap block NCM by N(O). Now consider the state vector X = (0 1 0 0 0 0 0 1 0 1 0 0 1 0 1) where the attributes $A_2$, $A_8$, $A_{10}$ $A_{13}$ and $A_{15}$ to be in the on state and all other nodes are in the off state. The effect of X on the dynamical system N (O) is given by

X N(O) ↪ ($I$ 1 1 0 $I$ 1 $I$ 1 1 1 1 0 1 1 1) = X₁ (say)

$$X_1 \, N(O) \hookrightarrow \quad (I \, 1 \, 0 \, I \, I \, 0 \, 0 \, 1 \, 1 \, 1 \, 1 \, 1 \, 1 \, 1 \, 1) \quad = \quad X_2 \text{ (say)}$$

$$X_2 \, N(O) \hookrightarrow \quad (I \, 1 \, 0 \, I \, I \, 0 \, I \, 1 \, 1 \, 1 \, 1 \, 1 \, 1 \, 1 \, 1) \quad = \quad X_3 = (X_2)$$

Thus the hidden pattern is a fixed point of the dynamical system. When male chauvinism, more leisure, only physically active, unreachable by friends and relatives and failure of agriculture is in the on state that is the HIV/AIDS patient all the five attributes only and all other attributes are in the off state we see, the attribute whether he has no binding with the family is an indeterminate, also whether he has bad habits and bad company is also an indeterminate. The concept of socially irresponsible and he is glutton remains as indeterminates. The nodes women as inferior objects, uncontrollable sex feelings are in the off state. All other states become on that is he has no work for the brain, visits CSWs, enjoys life and they are sexual perverts.

Several other state vectors can be substituted and conclusions arrived at based on them. Now we depict the same model with the attributes $A_1$, $A_2$,…, $A_{12}$ given in pages 42 to 52. Let us take the attributes $A_1$, $A_2$,…, $A_{12}$ and divide them into blocks of equal length but with over lap.

$$
\begin{aligned}
C_1 &= \{A_1, A_2, A_3 \, A_4\}, & C_2 &= \{A_3, A_4, A_5 \, A_6\}, \\
C_3 &= \{A_5, A_6, A_7 \, A_8\}, & C_4 &= \{A_7, A_8, A_9 \, A_{10}\} \\
C_5 &= \{A_9, A_{10}, A_{11} \, A_{12}\} \text{ and} & C_6 &= \{A_{11}, A_{12}, A_1 \, A_2\}
\end{aligned}
$$

and using the classes $C_1$, $C_2$,…, $C_6$ find using the experts opinion the neutrosophic graph and their related connection matrix using which we would find the combined overlap block NCM of blocks of equal sizes.

The neutrosophic graph related with class $C_1 = \{A_1, A_2, A_3, A_4\}$ given by the expert is as follows.

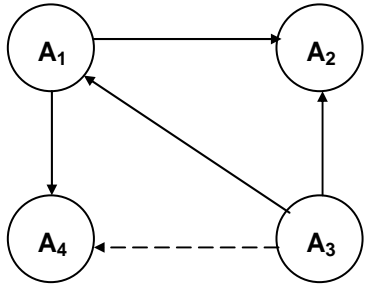

FIGURE: 4.4.20

The related connection matrix is

$$
\begin{array}{c c c c c}
 & A_1 & A_2 & A_3 & A_4 \\
\begin{array}{c} A_1 \\ A_2 \\ A_3 \\ A_4 \end{array} &
\left[ \begin{array}{c c c c}
0 & 1 & 0 & 1 \\
0 & 0 & 0 & 0 \\
1 & 1 & 0 & I \\
0 & 0 & 0 & 0
\end{array} \right]
\end{array}
$$

Next we give the neutrosophic directed graph for $C_2 = \{A_3 \ A_4 \ A_5 \ A_6\}$ given by the expert.

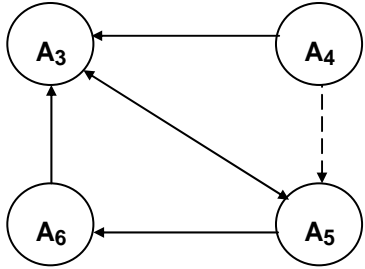

FIGURE: 4.4.21

The relational matrix associated with the above graph.

$$
\begin{array}{c c c c c}
 & A_3 & A_4 & A_5 & A_6 \\
\begin{array}{c} A_3 \\ A_4 \\ A_5 \\ A_6 \end{array} &
\left[ \begin{array}{c c c c}
0 & 0 & 1 & 0 \\
1 & 0 & I & 0 \\
1 & 0 & 0 & 1 \\
1 & 0 & 0 & 0
\end{array} \right]
\end{array}
$$

The directed graph given by the expert related to the class $C_3 = \{A_5 \ A_6 \ A_7 \ A_8\}$.

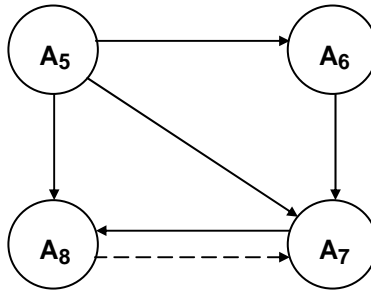

FIGURE: 4.4.22

The related matrix

$$
\begin{array}{c@{\quad}c}
 & \begin{array}{cccc} A_5 & A_6 & A_7 & A_8 \end{array} \\
\begin{array}{c} A_5 \\ A_6 \\ A_7 \\ A_8 \end{array} &
\left[\begin{array}{cccc}
0 & 1 & 1 & 1 \\
0 & 0 & 1 & 0 \\
0 & 0 & 0 & 1 \\
0 & 0 & I & 0
\end{array}\right]
\end{array}
$$

The directed graph associated with $\{A_7\ A_8\ A_9\ A_{10}\}$ is as follows

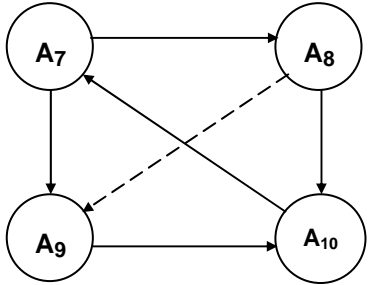

FIGURE: 4.4.23

The related connection matrix

$$
\begin{array}{c@{\quad}c}
 & \begin{array}{cccc} A_7 & A_8 & A_9 & A_{10} \end{array} \\
\begin{array}{c} A_7 \\ A_8 \\ A_9 \\ A_{10} \end{array} &
\left[\begin{array}{cccc}
0 & 1 & 1 & 0 \\
0 & 0 & I & 1 \\
0 & 0 & 0 & 1 \\
1 & 0 & 0 & 0
\end{array}\right]
\end{array}
$$

The directed graph given by the expert for the attributes $\{A_9\ A_{10}\ A_{11}\ A_{12}\}$.

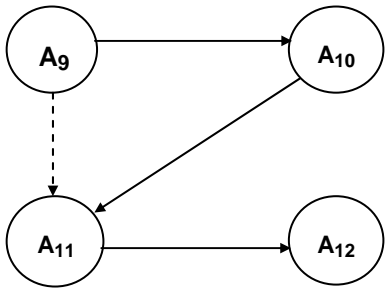

FIGURE: 4.4.24

The related connection matrix

$$
\begin{array}{c}
\quad\quad A_9\; A_{10}\; A_{11}\; A_{12} \\
\begin{array}{c} A_9 \\ A_{10} \\ A_{11} \\ A_{12} \end{array}
\begin{bmatrix}
0 & 1 & I & 0 \\
0 & 0 & 1 & 0 \\
0 & 0 & 0 & 1 \\
0 & 0 & 0 & 0
\end{bmatrix}
\end{array}
$$

For the attributes $C_6 = \{A_{11}, A_{12}, A_1\; A_2\}$ and the related graph is as follows:

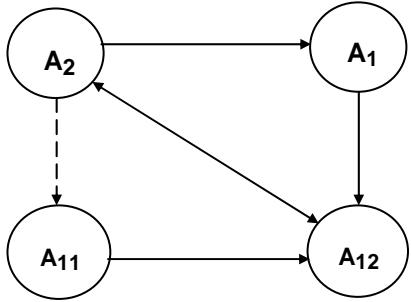

FIGURE: 4.4.25

The related connection matrix

$$
\begin{array}{c}
\quad\quad A_{11}\; A_{12}\; A_1\; A_2 \\
\begin{array}{c} A_{11} \\ A_{12} \\ A_1 \\ A_2 \end{array}
\begin{bmatrix}
0 & 1 & 0 & 0 \\
0 & 0 & 0 & 1 \\
0 & 1 & 0 & 0 \\
I & 1 & 1 & 0
\end{bmatrix}
\end{array}
$$

Using all these connection matrices we obtain the combined block overlap neutrosophic cognitive map.

Let us denote this $12 \times 12$ matrix O(A).

$$\begin{array}{c} \\ A_1 \\ A_2 \\ A_3 \\ A_4 \\ A_5 \\ A_6 \\ A_7 \\ A_8 \\ A_9 \\ A_{10} \\ A_{11} \\ A_{12} \end{array} \begin{array}{cccccccccccc} A_1 & A_2 & A_3 & A_4 & A_5 & A_6 & A_7 & A_8 & A_9 & A_{10} & A_{11} & A_{12} \\ 0 & 1 & 0 & 1 & 0 & 0 & 0 & 0 & 0 & 0 & 0 & 1 \\ 1 & 0 & 0 & 0 & 0 & 0 & 0 & 0 & 0 & 0 & I & 1 \\ 1 & 1 & 0 & I & 1 & 0 & 0 & 0 & 0 & 0 & 0 & 0 \\ 0 & 0 & 1 & 0 & I & 0 & 0 & 0 & 0 & 0 & 0 & 0 \\ 0 & 0 & 1 & 0 & 0 & 2 & 1 & 1 & 0 & 0 & 0 & 0 \\ 0 & 0 & 1 & 0 & 0 & 0 & 1 & 0 & 0 & 0 & 0 & 0 \\ 0 & 0 & 0 & 0 & 0 & 0 & 0 & 2 & 1 & 0 & 0 & 0 \\ 0 & 0 & 0 & 0 & 0 & 0 & 1 & 0 & I & 0 & 0 & 0 \\ 0 & 0 & 0 & 0 & 0 & 0 & 0 & 0 & 0 & 2 & I & 1 \\ 0 & 0 & 0 & 0 & 0 & 0 & 1 & 0 & 0 & 0 & 1 & 0 \\ 0 & 0 & 0 & 0 & 0 & 0 & 0 & 0 & 0 & 0 & 0 & 2 \\ 0 & 1 & 0 & 0 & 0 & 0 & 0 & 0 & 0 & 0 & 0 & 0 \end{array}$$

Using this matrix we obtain the effect of any state vector.

Suppose X = {1 0 0 0 0 1 0 0 0 0 1 0}, that is the nodes $A_1$, $A_6$ and $A_{11}$ are in the on state and all other nodes are in the off state. To find the effect of X on the dynamical system O(A)

| | | | | |
|---|---|---|---|---|
| XO (A) | ↪ | (1 1 1 1 0 1 1 0 0 0 1 1) | = | $X_1$ (say) |
| $X_1$ O(A) | ↪ | (1 1 1 $II$ 1 1 1 1 0 $I$ 1} | = | $X_2$ (say) |
| $X_2$ O(A) | ↪ | (1 1 1 $II$ 1 1 1 $I$ 1 1 1) | = | $X_3$ (say) |
| $X_3$ O(A) | ↪ | (1 1 1 $II$ 1 1 1 $II$ 1 1) | = | $X_4$ (say) |
| $X_4$ O(A) | ↪ | (1 1 1 $II$ 1 1 1 $II$ 1 1) | = | $X_5$ (= $X_4$) |

Thus the resultant vector is a fixed point which has made all nodes on or indeterminate. Now as the final illustrations we describe the newly defined combined block overlap NCM with blocks of varying length. Now we consider the same case where the attributes associated with the HIV/AIDS migrant labourer is taken as $P_1 P_2, \ldots, P_{15}$. (Page 57 of chapter II).

Now we divide the 15 attributes into classes

$C_1$ = $(P_1 P_2 P_3 P_4 P_5)$, $C_2$ = $(P_4 P_5 P_6 P_7 P_8 P_9 P_{10})$,
$C_3$ = $(P_9 P_{10} P_{11} P_{12})$, $C_4$ = $(P_{10} P_{11} P_{12} P_{13} P_{14} P_{15})$
and $C_5$ = $(P_{15} P_1 P_2 P_3)$.

Now using these attributes in these 5 classes we give the neutrosophic directed graph given by the experts. The directed graph related with the attributes $\{P_1\ P_2\ P_3\ P_4\ P_5\}$.

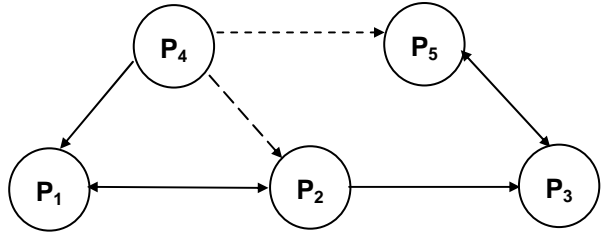

FIGURE: 4.4.26

The related matrix

$$
\begin{array}{c c}
 & \begin{array}{c c c c c} P_1 & P_2 & P_3 & P_4 & P_5 \end{array} \\
\begin{array}{c} P_1 \\ P_2 \\ P_3 \\ P_4 \\ P_5 \end{array} &
\begin{bmatrix}
0 & 1 & 0 & 0 & 0 \\
1 & 0 & 1 & 0 & 0 \\
0 & 0 & 0 & 0 & 1 \\
1 & I & 0 & 0 & I \\
0 & 0 & 1 & 0 & 0
\end{bmatrix}
\end{array}
$$

The directed graph given by the expert related with the attributes $C_2 = \{P_4\ P_5\ P_6\ P_7\ P_8\ P_9\ P_{10}\}$.

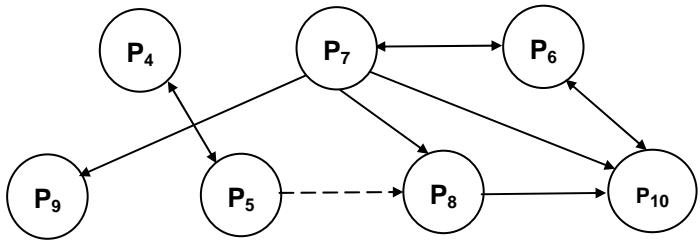

FIGURE: 4.4.27

The related connection matrix is

$$\begin{array}{c} \quad\quad P_4\ P_5\ P_6\ P_7\ P_8\ P_9\ P_{10} \\ \begin{array}{c} P_4 \\ P_5 \\ P_6 \\ P_7 \\ P_8 \\ P_9 \\ P_{10} \end{array} \begin{bmatrix} 0 & 1 & 0 & 0 & 0 & 0 & 0 \\ 1 & 0 & 0 & 0 & I & 0 & 0 \\ 0 & 0 & 0 & 1 & 0 & 0 & 1 \\ 0 & 0 & 1 & 0 & 1 & 1 & 1 \\ 0 & 0 & 0 & 0 & 0 & 0 & 1 \\ 0 & 0 & 0 & 0 & 0 & 0 & 0 \\ 0 & 0 & 1 & 0 & 0 & 0 & 0 \end{bmatrix} \end{array}$$

using the experts opinion we have the following directed graph for the attributes $C_3 = \{P_9\ P_{10}\ P_{11}\ P_{12}\}$.

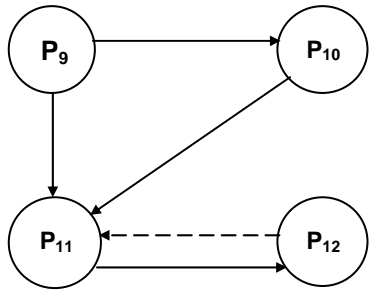

FIGURE: 4.4.28

The related connection matrix

$$\begin{array}{c} \quad\quad P_9\ P_{10}\ P_{11}\ P_{12} \\ \begin{array}{c} P_9 \\ P_{10} \\ P_{11} \\ P_{12} \end{array} \begin{bmatrix} 0 & 1 & 1 & 0 \\ 0 & 0 & 1 & 0 \\ 0 & 0 & 0 & 1 \\ 0 & 0 & I & 0 \end{bmatrix} \end{array}$$

The directed graph related to the attributes in $C_4$; $C_4 = \{P_{10}\ P_{11}\ P_{12}\ P_{13}\ P_{14}\ P_{15}\}$.

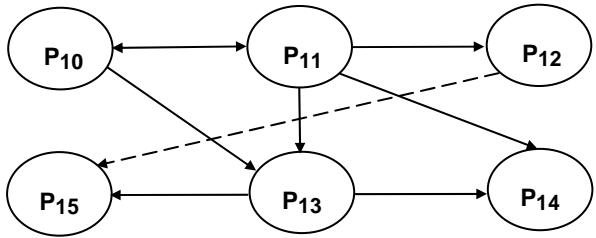

FIGURE: 4.4.29

The related connection matrix

$$
\begin{array}{c}
\phantom{P_{10}}\begin{array}{cccccc} P_{10} & P_{11} & P_{12} & P_{13} & P_{14} & P_{15} \end{array}\\
\begin{array}{c} P_{10} \\ P_{11} \\ P_{12} \\ P_{13} \\ P_{14} \\ P_{15} \end{array}
\left[\begin{array}{cccccc}
0 & 1 & 0 & 1 & 0 & 0 \\
1 & 0 & 1 & 1 & 1 & 0 \\
0 & 0 & 0 & 0 & 0 & I \\
0 & 0 & 0 & 0 & 1 & 1 \\
0 & 0 & 0 & 0 & 0 & 0 \\
0 & 0 & 0 & 0 & 0 & 0
\end{array}\right]
\end{array}
$$

The directed graph related with the final class of attributes $C_5 = \{P_{15}\, P_1 P_2\, P_3\}$.

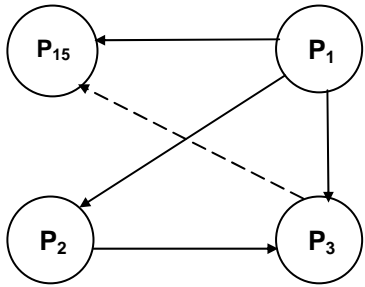

FIGURE: 4.4.30

The related connection matrix

$$
\begin{array}{c}
\phantom{P_{15}}\begin{array}{cccc} P_{15} & P_1 & P_2 & P_3 \end{array}\\
\begin{array}{c} P_{15} \\ P_1 \\ P_2 \\ P_3 \end{array}
\left[\begin{array}{cccc}
0 & 0 & 0 & 0 \\
1 & 0 & 1 & 1 \\
0 & 0 & 0 & 1 \\
I & 0 & 0 & 0
\end{array}\right]
\end{array}
$$

Now we give the connection matrix of the combined block overlap NCM with blocks of varying length

$$
\begin{array}{c}
\begin{array}{cccccccccccccccc}
 & P_1 & P_2 & P_3 & P_4 & P_5 & P_6 & P_7 & P_8 & P_9 & P_{10} & P_{11} & P_{12} & P_{13} & P_{14} & P_{15}
\end{array} \\
\begin{array}{c}
P_1 \\ P_2 \\ P_3 \\ P_4 \\ P_5 \\ P_6 \\ P_7 \\ P_8 \\ P_9 \\ P_{10} \\ P_{11} \\ P_{12} \\ P_{13} \\ P_{14} \\ P_{15}
\end{array}
\begin{bmatrix}
0 & 2 & 1 & 0 & 0 & 0 & 0 & 0 & 0 & 0 & 0 & 0 & 0 & 0 & 1 \\
I & 0 & 2 & 0 & 0 & 0 & 0 & 0 & 0 & 0 & 0 & 0 & 0 & 0 & 0 \\
0 & 0 & 0 & 0 & 1 & 0 & 0 & 0 & 0 & 0 & 0 & 0 & 0 & 0 & I \\
1 & I & 0 & 0 & I & 0 & 0 & 0 & 0 & 0 & 0 & 0 & 0 & 0 & 0 \\
0 & 0 & 1 & 1 & 0 & 0 & 0 & I & 0 & 0 & 0 & 0 & 0 & 0 & 0 \\
0 & 0 & 0 & 0 & 0 & 0 & 1 & 0 & 0 & 1 & 0 & 0 & 0 & 0 & 0 \\
0 & 0 & 0 & 0 & 0 & 1 & 0 & 1 & 1 & 1 & 0 & 0 & 0 & 0 & 0 \\
0 & 0 & 0 & 0 & 0 & 0 & 0 & 0 & 0 & 0 & 0 & 0 & 0 & 0 & 0 \\
0 & 0 & 0 & 0 & 0 & 0 & 0 & 0 & 0 & 1 & 0 & 0 & 0 & 0 & 0 \\
0 & 0 & 0 & 0 & 0 & 1 & 0 & 0 & 0 & 0 & 1 & 0 & 1 & 0 & 0 \\
0 & 0 & 0 & 0 & 0 & 0 & 0 & 0 & 0 & 1 & 0 & 1 & 1 & 1 & 0 \\
0 & 0 & 0 & 0 & 0 & 0 & 0 & 0 & 0 & 0 & 0 & 0 & 0 & 0 & I \\
0 & 0 & 0 & 0 & 0 & 0 & 0 & 0 & 0 & 0 & 0 & 0 & 0 & 1 & 1 \\
0 & 0 & 0 & 0 & 0 & 0 & 0 & 0 & 0 & 0 & 0 & 0 & 0 & 0 & 0 \\
0 & 0 & 0 & 0 & 0 & 0 & 0 & 0 & 0 & 0 & 1 & 0 & 0 & 0 & 0
\end{bmatrix}
\end{array}
$$

Now we study the effect of any state vector on the dynamical system $N(O_1)$ where $N(O_1)$ represents a $15 \times 15$ matrix.

For consider the sate vector $A_1$, $A_5$, $A_9$ $A_{11}$ and $A_{15}$ in the on state and all vectors with off state

$$X \quad = \quad (1\ 0\ 0\ 0\ 1\ 0\ 0\ 0\ 1\ 0\ 1\ 0\ 0\ 0\ 1)$$

The effect of X on the dynamical system $N(O_1)$ is given by

$$XN\,(O_1) \quad \hookrightarrow \quad (1\ 1\ 1\ 1\ 1\ 0\ 0\ I\ 1\ 1\ 1\ 1\ 1\ 1\ 1) \quad = \quad X_1 \text{ (say)}$$

$$X_1\,N(O_1) \quad \hookrightarrow \quad (1\ 1\ 1\ 1\ 1\ 1\ 1\ 0\ I\ 1\ 1\ 1\ 1\ 1\ 1\ 1) \quad = \quad X_2 \text{ (say)}$$

$$X_2\,N(O_2) \quad \hookrightarrow \quad (1\ 1\ 1\ 1\ 1\ 1\ 1\ 1\ I\ 1\ 1\ 1\ 1\ 1\ 1\ 1) \quad = \quad X_3 \text{ (say)}$$

$$X_3\,N\,(O_1) \quad \hookrightarrow \quad (1\ 1\ 1\ 1\ 1\ 1\ 1\ 1\ I\ 1\ 1\ 1\ 1\ 1\ 1\ 1)$$

The hidden pattern is a fixed point of the system. Using the resultant vector any reader can interpret the resultant vector.

Next we model the problem given in page 42 to 53 of chapter two using the combined block overlap NCF.

Let $A_1$, $A_2$, ..., $A_{12}$ be taken as in page 42 to 53. Now divide $A_1$, $A_2$,..., $A_{12}$ into overlapping classes of different lengths.

Let $C_1 = \{A_1\ A_2\ A_3\ A_4\ A_5\ A_6\}$, $C_2 = \{A_5\ A_6\ A_7\ A_8\}$ and $C_3 = \{A_8\ A_9\ A_{10}\ A_{11}\ A_{12}\ A_1, A_2\}$.

Now using the experts opinion we find the related connection matrix from the neutrosophic graph given by them for the 3 classes of attributes $C_1\ C_2\ C_3$.

The directed graph related to $C_1 = \{A_1\ A_2\ A_3\ A_4\ A_5\ A_6\}$ is

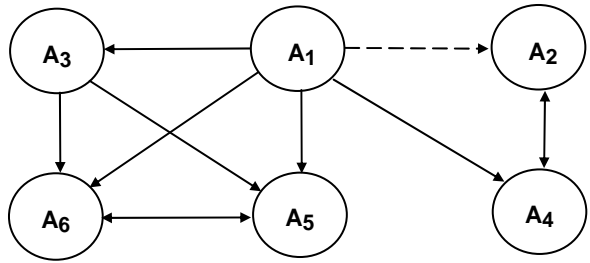

FIGURE: 4.4.31

The related neutrosophic connection matrix

$$\begin{array}{c} \\ A_1 \\ A_2 \\ A_3 \\ A_4 \\ A_5 \\ A_6 \end{array} \begin{array}{c} \begin{array}{cccccc} A_1 & A_2 & A_3 & A_4 & A_5 & A_6 \end{array} \\ \begin{bmatrix} 0 & I & 1 & 1 & 1 & 1 \\ 0 & 0 & 0 & 1 & 0 & 0 \\ 0 & 0 & 0 & 0 & 1 & 1 \\ 0 & 1 & 0 & 0 & 0 & 0 \\ 0 & 0 & 0 & 0 & 0 & 1 \\ 0 & 0 & 0 & 0 & 1 & 0 \end{bmatrix} \end{array}$$

The neutrosophic directed graph related with $C_2 = \{A_5\ A_6\ A_7\ A_8\}$ as given by the expert.

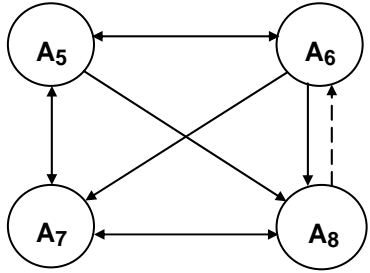

FIGURE: 4.4.32

The related connection matrix.

$$
\begin{array}{c}
\quad\ A_5\ A_6\ A_7\ A_8 \\
\begin{array}{c}
A_5 \\ A_6 \\ A_7 \\ A_8
\end{array}
\left[
\begin{array}{cccc}
0 & 1 & 1 & 1 \\
1 & 0 & 1 & 1 \\
1 & 0 & 0 & 1 \\
0 & I & 1 & 0
\end{array}
\right]
\end{array}
$$

Now using the class $C_3 = \{A_8\ A_9\ A_{10}\ A_{11}\ A_{12}\ A_1\ A_2\}$ we form the neutrosophic directed graph using the experts opinion.

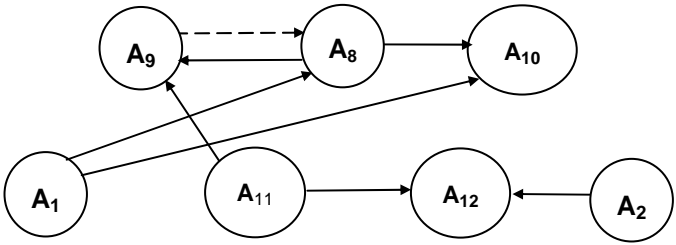

FIGURE: 4.4.33

Now we give the relational connection matrix associated with the directed graph.

$$\begin{array}{c@{\quad}ccccccc}
& A_8 & A_9 & A_{10} & A_{11} & A_{12} & A_1 & A_2 \\
A_8 & 0 & 1 & 1 & 0 & 0 & 0 & 0 \\
A_9 & I & 0 & 0 & 0 & 0 & 0 & 0 \\
A_{10} & 0 & 0 & 0 & 0 & 0 & 0 & 0 \\
A_{11} & 0 & 1 & 0 & 0 & 1 & 0 & 0 \\
A_{12} & 0 & 0 & 0 & 0 & 0 & 0 & 0 \\
A_1 & 1 & 0 & 1 & 0 & 0 & 0 & 0 \\
A_2 & 0 & 0 & 0 & 0 & 1 & 0 & 0
\end{array}$$

Now we find the matrix related with combined block over lap NCM.

$$\begin{array}{c@{\quad}cccccccccccc}
& A_1 & A_2 & A_3 & A_4 & A_5 & A_6 & A_7 & A_8 & A_9 & A_{10} & A_{11} & A_{12} \\
A_1 & 0 & I & 1 & 1 & 1 & 1 & 0 & 1 & 0 & 1 & 0 & 0 \\
A_2 & 0 & 0 & 0 & 1 & 0 & 0 & 0 & 0 & 0 & 0 & 0 & 1 \\
A_3 & 0 & 0 & 0 & 0 & 1 & 1 & 0 & 0 & 0 & 0 & 0 & 0 \\
A_4 & 0 & 1 & 0 & 0 & 0 & 0 & 0 & 0 & 0 & 0 & 0 & 0 \\
A_5 & 0 & 0 & 0 & 0 & 0 & 2 & 1 & 1 & 0 & 0 & 0 & 0 \\
A_6 & 0 & 0 & 0 & 0 & 2 & 0 & 1 & 1 & 0 & 0 & 0 & 0 \\
A_7 & 0 & 0 & 0 & 0 & 1 & 0 & 0 & 0 & 0 & 0 & 0 & 0 \\
A_8 & 0 & 0 & 0 & 0 & 0 & I & 1 & 0 & 1 & 1 & 0 & 0 \\
A_9 & 0 & 0 & 0 & 0 & 0 & 0 & 0 & I & 0 & 0 & 0 & 0 \\
A_{10} & 0 & 0 & 0 & 0 & 0 & 0 & 0 & 0 & 0 & 0 & 0 & 0 \\
A_{11} & 0 & 0 & 0 & 0 & 0 & 0 & 0 & 0 & 1 & 0 & 0 & 1 \\
A_{12} & 0 & 0 & 0 & 0 & 0 & 0 & 0 & 0 & 0 & 0 & 0 & 0
\end{array}$$

Now we can study the effect of each and every state vector. Let us denote this $12 \times 12$ by matrix by N(R).

Let us consider the state vector X = (1 1 0 0 0 1 0 1 0 0 0 1). The effect of X on the system N (R) is given by

$$XN\,(R) \;\hookrightarrow\; (1\ 1\ 1\ 1\ 1\ 1\ 1\ 1\ 1\ 1\ 0\ 1) \;=\; X_1 \text{ (say)}$$

$$X_1\,N\,(R) \;\hookrightarrow\; (1\ 1\ 1\ 1\ 1\ 1\ 1\ 1\ 1\ 1\ 0\ 1) \;=\; X_2\,(= X_1 \text{ say})$$

This is a fixed point. Thus when the nodes $A_1$, $A_2$, $A_6$, $A_8$ and $A_{12}$ are in the on state in the input vector we obtain a resultant in which all nodes come to on state except $A_{11}$. Hence these five nodes have effect on all nodes except $A_{11}$.

The interested reader can construct a C program for NRM and CNRM as in case of FRMs and CFRMs given in appendix 7 and 8 of this book.

Chapter Five

# USE OF FRM AND NRM TO STUDY THE HIV/AIDS AFFECTED MIGRANT LABOURERS AND THEIR SOCIO ECONOMIC CONDITIONS

In this chapter we study the socio economic conditions and its relation to the HIV/AIDS affected migrant labourers. The main observation from our study is that the migration of the labourers has resulted due to unemployment problems mainly the failure of agriculture. Thus these people with least knowledge about the urban life and the evils of urban society come for livelihood and with a mission for earning for the family but unfortunately they learn all the evils of the urban society and specially became victims of HIV/AIDS as their only and sole recreation and past time after work hours is sex by visiting CSWs.

Urban people are cleaver and well aware of HIV/AIDS so they go for protected sex on the contrary these rural male go for unprotected sex by not only infecting themselves but their wife's and children in due course of time. This is one of the major reasons that only rural poor uneducated flood the hospital and we saw most of them in a very pathetic state. Majority of them (90% of them) came to know about HIV/AIDS and its incurability only after they had become infected by it. With a great fear to let the family know they kept it in dark even from their wives and unfortunately infected them.

Thus over 95% of the women who have been infected by HIV/AIDS from rural areas were only infected by their husbands. As the men to woman transmission is very fast i.e., 24 times greater than woman to man transmission these infected women became full blown HIV/AIDS patients in a span of 6months to 5 years, and they came for treatment only when the disease had become chromic in them.

Thus here we study the problem using FRMs and NRMs. This chapter has six sections. In section one we define we define FRM and study this model using FRM. In this section several experts opinion are used and the combined FRM model is also used to find the conclusions. In section two we for the first time define the notion of linked FRM and its applications. We define for the first time the new notions of combined block disjoint FRM and combined block overlap FRM and use them to study this model.

This is done in section three of this chapter. In section four we introduce NRM and study the problem using them. In section five we introduce the new notion of linked neutrosophic relational maps and apply it in the case of the social problems faced by the HIV/AIDS patients. For the first time we define the notions of combined block disjoint neutrosophic relational maps and combined block overlap neutrosophic relational maps and study them and adopt them in this model. This analysis is carried out in section 6 of this chapter.

## 5.1 Use of Fuzzy Relational Maps in the Study of Relation between HIV/AIDS Migrants and their Socio Economic Conditions

The notion of FRM is new it was defined in [110, 111]. We just recall the basic definition of FRM. We introduce the notion of Fuzzy relational maps (FRMs); they are constructed analogous to FCMs described and discussed in the earlier sections. In FCMs we promote the correlations between causal associations among concurrently active units.

But in FRMs we divide the very causal associations into two disjoint units, for example, the relation between a teacher and a student or relation between an employee or employer or a relation between doctor and patient and so on. Thus for us to define a FRM we need a domain space and a range space which are disjoint in the sense of concepts. We further assume no intermediate relation exists within the domain elements or node and the range spaces elements. The number of elements in the range space need not in general be equal to the number of elements in the domain space.

Thus throughout this section we assume the elements of the domain space are taken from the real vector space of dimension n and that of the range space are real vectors from the vector space of dimension m (m in general need not be equal to n). We denote by R the set of nodes $R_1,\ldots, R_m$ of the range space, where R = $\{(x_1,\ldots, x_m) \mid x_j = 0 \text{ or } 1\}$ for j = 1, 2,…, m. If $x_i = 1$ it means that the node $R_i$ is in the on state and if $x_i = 0$ it means that the node $R_i$ is in the off state. Similarly D denotes the nodes $D_1, D_2,\ldots, D_n$ of the domain space where D = $\{(x_1,\ldots, x_n) \mid x_j = 0 \text{ or } 1\}$ for i = 1, 2,…, n. If $x_i = 1$ it means that the node $D_i$ is in the on state and if $x_i = 0$ it means that the node $D_i$ is in the off state.

Now we proceed on to define a FRM.

**DEFINITION 5.1.1:** *A FRM is a directed graph or a map from D to R with concepts like policies or events etc, as nodes and causalities as edges. It represents causal relations between spaces D and R .*

*Let $D_i$ and $R_j$ denote that the two nodes of an FRM. The directed edge from $D_i$ to $R_j$ denotes the causality of $D_i$ on $R_j$ called relations. Every edge in the FRM is weighted with a number in the set {0, ±1}.*

*Let $e_{ij}$ be the weight of the edge $D_iR_j$, $e_{ij} \in \{0, ±1\}$. The weight of the edge $D_i R_j$ is positive if increase in $D_i$ implies increase in $R_j$ or decrease in $D_i$ implies decrease in $R_j$ ie causality of $D_i$ on $R_j$ is 1. If $e_{ij} = 0$, then $D_i$ does not have any effect on $R_j$ . We do not discuss the cases when increase in $D_i$ implies decrease in $R_j$ or decrease in $D_i$ implies increase in $R_j$ .*

**DEFINITION 5.1.2:** *When the nodes of the FRM are fuzzy sets then they are called fuzzy nodes. FRMs with edge weights {0, ±1} are called simple FRMs.*

**DEFINITION 5.1.3:** *Let $D_1, \ldots, D_n$ be the nodes of the domain space D of an FRM and $R_1, \ldots, R_m$ be the nodes of the range space R of an FRM. Let the matrix E be defined as E = $(e_{ij})$ where $e_{ij}$ is the weight of the directed edge $D_iR_j$ (or $R_jD_i$), E is called the relational matrix of the FRM.*

<u>*Note*</u>: It is pertinent to mention here that unlike the FCMs the FRMs can be a rectangular matrix with rows corresponding to the domain space and columns corresponding to the range space. This is one of the marked difference between FRMs and FCMs.

**DEFINITION 5.1.4:** *Let $D_1$, …, $D_n$ and $R_1$,…, $R_m$ denote the nodes of the FRM. Let $A = (a_1,…,a_n)$, $a_i \in \{0, \pm1\}$. A is called the instantaneous state vector of the domain space and it denotes the on-off position of the nodes at any instant.*

*Similarly let $B = (b_1,…, b_m)$ $b_i \in \{0, \pm1\}$. B is called instantaneous state vector of the range space and it denotes the on-off position of the nodes at any instant $a_i = 0$ if $a_i$ is off and $a_i = 1$ if $a_i$ is on for i= 1, 2,…, n Similarly, $b_i = 0$ if $b_i$ is off and $b_i = 1$ if $b_i$ is on, for i= 1, 2,…, m.*

**DEFINITION 5.1.5:** *Let $D_1$, …, $D_n$ and $R_1$,…, $R_m$ be the nodes of an FRM. Let $D_iR_j$ (or $R_j\ D_i$) be the edges of an FRM, $j = 1, 2,…, m$ and i= 1, 2,…, n. Let the edges form a directed cycle. An FRM is said to be a cycle if it posses a directed cycle. An FRM is said to be acyclic if it does not posses any directed cycle.*

**DEFINITION 5.1.6:** *An FRM with cycles is said to be an FRM with feedback.*

**DEFINITION 5.1.7:** *When there is a feedback in the FRM, i.e. when the causal relations flow through a cycle in a revolutionary manner, the FRM is called a dynamical system.*

**DEFINITION 5.1.8:** *Let $D_i\ R_j$ (or $R_j\ D_i$), $1 \leq j \leq m$, $1 \leq i \leq n$. When $R_i$ (or $D_j$) is switched on and if causality flows through edges of the cycle and if it again causes $R_i$ (or$D_j$), we say that the dynamical system goes round and round.*

*This is true for any node $R_j$ (or $D_i$) for $1 \leq i \leq n$, (or $1 \leq j \leq m$). The equilibrium state of this dynamical system is called the hidden pattern.*

**DEFINITION 5.1.9:** *If the equilibrium state of a dynamical system is a unique state vector, then it is called a fixed point. Consider an FRM with $R_1$, $R_2$,…, $R_m$ and $D_1$, $D_2$,…, $D_n$ as nodes.*

*For example, let us start the dynamical system by switching on $R_1$ (or $D_1$). Let us assume that the FRM settles down with $R_1$ and $R_m$ (or $D_1$ and $D_n$) on, i.e. the state vector remains as (1, 0, …, 0, 1) in R (or 1, 0, 0, … , 0, 1) in D), This state vector is called the fixed point.*

**DEFINITION 5.1.10:** *If the FRM settles down with a state vector repeating in the form*

$$A_1 \rightarrow A_2 \rightarrow A_3 \rightarrow \ldots \rightarrow A_i \rightarrow A_1 \ (or \ B_1 \rightarrow B_2 \rightarrow \ldots \rightarrow B_i \rightarrow B_1)$$

*then this equilibrium is called a limit cycle.*

**Methods of Determining the Hidden Pattern**

Let $R_1, R_2, \ldots, R_m$ and $D_1, D_2, \ldots, D_n$ be the nodes of a FRM with feedback. Let E be the relational matrix. Let us find a hidden pattern when $D_1$ is switched on i.e. when an input is given as vector $A_1 = (1, 0, \ldots, 0)$ in $D_1$, the data should pass through the relational matrix E. This is done by multiplying $A_1$ with the relational matrix E. Let $A_1E = (r_1, r_2, \ldots, r_m)$, after thresholding and updating the resultant vector we get $A_1 E \in R$. Now let $B = A_1E$ we pass on B into $E^T$ and obtain $BE^T$. We update and threshold the vector $BE^T$ so that $BE^T \in D$. This procedure is repeated till we get a limit cycle or a fixed point.

**DEFINITION 5.1.11:** *Finite number of FRMs can be combined together to produce the joint effect of all the FRMs. Let $E_1, \ldots, E_p$ be the relational matrices of the FRMs with nodes $R_1, R_2, \ldots, R_m$ and $D_1, D_2, \ldots, D_n$, then the combined FRM is represented by the relational matrix $E = E_1 + \ldots + E_p$.*

Of the 60 HIV/AIDS infected migrants labourers interviewed only one was infected by his wife who was a cook and 3 persons were infected on their jolly trips to enjoy life. Thus all the 56 of them were infected only in their work place i.e., a place away from their hometown. It is also more important to mention that all the 59 of them were infected only by CSWs and several of them acknowledged that they visited CSWs innumerable number of times.

Thus their social conditions in urban life made them visit CSWs and become HIV/AIDS victims. The urban society so detached from the rural ones for the busy state of city never bothers about the happenings in the neighbourhood. This alfooness, uncared lonely life had lead the rural people to go for CSWs.

Also being from poor or lower middle class they stay in very ill furnished places and when the chance of them having a T.V even in their house is impossible it is still an impossibility to have a T.V. set so that they can spend their leisure in watching it. Thus their only spare time job is visiting CSWs.

For some when they live in an urban place which is not Tamil Nadu, the problem of language also plays a role forcing them to CSWs, as they cannot have the satisfaction of watching other language movies. As they are basically uneducated and poor do not have the facility of visiting tourist spots or place of importance to while away their spare time in good and constructive activities.

Thus after the work hours for majority of them CSWs was the only pastime. Being unaware of HIV/AIDS had unprotected sex coupled with all other bad habits became a victim of HIV/AIDS. The spread of HIV/AIDS was faster as most of them suffered from STD/VD, so the infection to HIV/AIDS to them was direct through blood.

Now using FRMs and NRMs we study the socio economic problem of HIV/AIDS affected persons.

We first take the opinion of an expert and study the problem using FRM. We take the following attributes as the nodes of the domain space) of Fuzzy Relational Maps which gives the problems faced by the migrant labourer from rural areas. $D = \{D_1, D_2, D_3, D_4, D_5, \ldots, D_8\}$ where $D_1, D_2, \ldots, D_8$ are described.

$D_1$ - Poverty
$D_2$ - Failure of Agriculture
$D_3$ - No awareness about HIV/AIDS
$D_4$ - More leisure after work hours
$D_5$ - Away from home for weeks
$D_6$ - No association for migrant labourers
$D_7$ - No alternatives provided by government.
$D_8$ - No fear of being observed by friends or relatives in urban life style.

The attributes taken for the range space are $\{R_1 \ R_2 \ R_3 \ R_4 \ R_5\}$

$R_1$ - Migration to City (Cause of)
$R_2$ - Visiting as labourers CSWs
$R_3$ - Bad Company
$R_4$ - Easy victims of bad habits
$R_5$ - Easily affected by HIV/AIDS.

We take an experts view who is a migrant labourer affected / infected by HIV/AIDS.

The directed graph related with this expert's opinion is given in the next page.

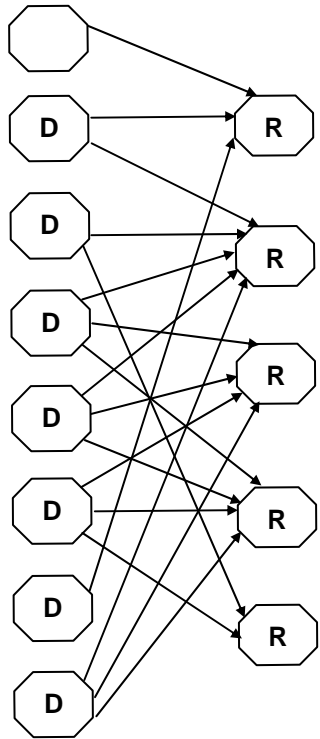

FIGURE: 5.1.1

The expert himself agrees that mainly no one observes them they are totally free so only they choose easily take up to all bad habits. The expert was in late twenties.

Also he says most of them are unaware of methods to save money they earn and spend in the same way. Some say they get some advance or an amount deducted by their employer, which is given to them when they go home for a week or 3 to four days. So they freely spend the money they earn daily. Except the truck drivers majority of these labourers are only daily wagers. They say the place they stay after work hours is very inconducible so only they go for CSWs.

Almost 95% (58 out of 60) of these persons drink alcohol daily. They also acknowledge that when they visit CSWs they are fully drunk. Thus this opinion of his was sought after discussion. He felt this was the problem faced by migrant labourers of his age group. Some were just married or married for 5 to 6 years so they invariably went to CSWs for sex.

We give the related matrix E of the directed graph given in figure 5.1.1.

$$
\begin{array}{c}
\quad\quad R_1\ R_2\ R_3\ R_4\ R_5 \\
\begin{array}{c}
D_1 \\ D_2 \\ D_3 \\ D_4 \\ D_5 \\ D_6 \\ D_7 \\ D_8
\end{array}
\left[
\begin{array}{ccccc}
1 & 0 & 0 & 0 & 0 \\
1 & 1 & 0 & 0 & 0 \\
0 & 1 & 0 & 0 & 1 \\
0 & 1 & 1 & 1 & 0 \\
0 & 1 & 1 & 1 & 0 \\
0 & 0 & 1 & 1 & 1 \\
1 & 0 & 0 & 0 & 0 \\
0 & 1 & 1 & 1 & 0
\end{array}
\right]
\end{array}
$$

Let E denote the connection matrix of the FRM. Now we study the hidden pattern of the dynamical system E. Suppose we consider the on state of the node $D_2$ i.e., failure of agriculture and all other nodes are in the off state.

Let

$$X \quad = \quad (0\ 1\ 0\ 0\ 0\ 0\ 0\ 0).$$

The effect of X on the dynamical system E is given by

$$XE \quad \hookrightarrow \quad (1\ 1\ 0\ 0\ 0) \quad\quad = \quad Y_1 \text{ say}$$

$$Y_1E^T \quad \hookrightarrow \quad (1\ 1\ 1\ 1\ 1\ 1\ 0\ 11) \quad = \quad X_1 \text{ say}$$

$$X_1E \quad \hookrightarrow \quad (1\ 1\ 1\ 1\ 1) \quad\quad = \quad Y_2 \text{ say}$$

$$Y_2E^T \quad \hookrightarrow \quad (11\ 1\ 1\ 1\ 1\ 1\ 1) \quad = \quad X_2 \text{ say}$$

$$X_2E \quad \hookrightarrow \quad (1\ 1\ 1\ 1\ 1) \quad\quad = \quad Y_3 \text{ say}$$

$$Y_3E^T \quad \hookrightarrow \quad (11\ 1\ 1\ 1\ 1\ 1\ 1) \quad = \quad X_3 = X_2 \text{ a fixed point}$$

of the hidden pattern. Thus

$$X_3E \quad \hookrightarrow \quad (1\ 1\ 1\ 1\ 1).$$

Thus we get the pair of resultant vectors to be $\{(1\ 1\ 1\ 1\ 1\ 1\ 1\ 1),\ (11\ 1\ 1\ 1)\}$. Thus the only node failure of agriculture has lead to the on state of all others nodes except poverty and no alternatives provided by the government. For when they go for job in urban area due their earning poverty may not exists. Once they have taken up jobs in cities it is immaterial whether the government is providing them with some alternatives or not for they have taken their own course of action. Once they face failure due to agriculture they see all the range space coordinates come to on state.

Now consider the on state of the state vector "easy victims of HIV/AIDS" i.e., $R_4$ is in the on state. The effect of $R_4$ on the dynamical system E is given by the following procedure

Let

$$T \quad = \quad (0\ 0\ 0\ 1\ 0)$$

$$TE^T \quad \hookrightarrow \quad (0\ 0\ 0\ 1\ 1\ 1\ 0\ 1) = \quad S \text{ (say)}$$

$$SE \quad \hookrightarrow \quad (0\ 1\ 1\ 1\ 1) \quad = \quad T_1 \text{ say}$$

$$T_1E^T \quad \hookrightarrow \quad (0\ 1\ 1\ 1\ 1\ 1\ 0\ 1) = \quad S_1 \text{ (say)}$$

$$S_1E \quad \hookrightarrow \quad (1\ 1\ 1\ 1\ 1) \quad = \quad T_2 \text{ (say)}$$

$$T_2E^T \quad \hookrightarrow \quad (1\ 1\ 1\ 1\ 1\ 1\ 1\ 1) = \quad S_2 \text{ (say)}$$

$$S_2E \quad \hookrightarrow \quad (1\ 1\ 1\ 1\ 1).$$

Thus the resultant vector is a fixed point given by {(1 1 1 1 1 1 1 1), (1 1 1 1 1)}. Thus when the coordinate easy victims of the disease HIV/AIDS is in the on state all other co-ordinates come to the on state i.e., they are easy victims as they are away from the family, poverty, no awareness, more leisure, no fear of being observed by the family failure of the agriculture. Several such conclusions can be obtained using the hidden pattern of the dynamical system.

Next we use a second expert's opinion to study the problem: The causes of migrant labourers vulnerability to HIV/AIDS and the role of the government. We model this using Fuzzy Relational Maps. The role of government is taken as the domain space and that of the migrant labourers vulnerability to HIV/AIDS as range space.

**ROLE OF GOVERNMENT**

$D_1$ -  Awareness program about HIV/AIDS in rural areas where most people are uneducated is insufficient.

$D_2$ - No alternative job provided to rural agricultural labourers when agriculture fails.

$D_3$ - Failure to stop mislead agriculture techniques.

$D_4$ - No for sight by the government to take precautionary measure from the past experience.

$D_5$ - Failure to stop functioning of CSWs. When CSWs flourish no steps taken to test them for HIV/AIDS and give them free medical aid and help them to use condoms.

$D_6$ - Government does not help rural uneducated agricultural labourers even after migration.

The attributes given by the expert that causes the migrant labourer vulnerability to HIV/AIDS.

$R_1$ - No education / No awareness.

$R_2$ - Cheap availability of CSWs.

$R_3$ - Away from family for weeks.

$R_4$ - Superstition about sex / profession.

$R_5$ - No union to channelize them and advice about HIV/AIDS.

$R_6$ - Not watched by people in city and unaware of the disease.

$R_7$ - No job in native place.

$R_8$ - Free availability of alcohol.

$R_9$ - Infertility of land, so labour contractors take advantage of poverty.

$R_{10}$ - No proper medical aid by government / counselling.

The directed graph related with this experts opinion.

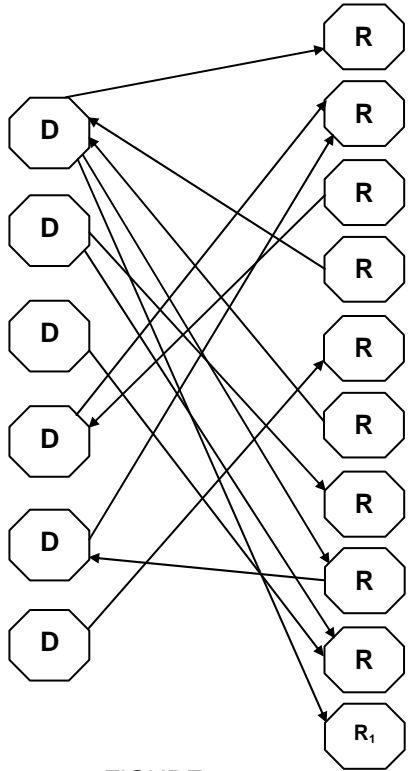

FIGURE: 5.1.2

The related connection matrix

$$
\begin{array}{c}
\phantom{D_1} \\
D_1 \\
D_2 \\
D_3 \\
D_4 \\
D_5 \\
D_6
\end{array}
\begin{array}{c}
R_1\ R_2\ R_3\ \ R_4\ R_5\ \ R_6\ R_7\ R_8\ R_9\ R_{10} \\
\begin{bmatrix}
1 & 0 & 0 & 1 & 0 & 1 & 0 & 1 & 0 & 1 \\
0 & 0 & 0 & 0 & 0 & 0 & 1 & 0 & 1 & 0 \\
0 & 0 & 0 & 0 & 0 & 0 & 0 & 0 & 1 & 0 \\
0 & 1 & 1 & 0 & 0 & 0 & 0 & 0 & 0 & 0 \\
0 & 1 & 0 & 0 & 0 & 0 & 0 & 1 & 0 & 0 \\
0 & 0 & 0 & 0 & 1 & 0 & 0 & 0 & 0 & 0
\end{bmatrix}
\end{array}
$$

Let F denote the $6 \times 10$ connection matrix. The effect of the state vector

$$X \quad = \quad (0\ 0\ 1\ 0\ 0\ 0),$$

i.e., the attributes $A_3$ is in the on state and all other nodes are in the off state

| | | | | |
|---|---|---|---|---|
| XF | ↪ | $(0\ 0\ 0\ 0\ 0\ 0\ 0\ 0\ 1\ 0)$ | = | Y (say) |
| $YF^T$ | ↪ | $(0\ 1\ 1\ 0\ 0\ 0)$ | = | $X_1$ say |
| $X_1\,F$ | ↪ | $(0\ 0\ 0\ 0\ 0\ 0\ 1\ 0\ 1\ 0)$ | = | $Y_1$ say |
| $Y_1F^T$ | ↪ | $(1\ 1\ 1\ 0\ 0\ 0)$ | = | $X_2$ say |
| $X_2F$ | ↪ | $(1\ 0\ 0\ 1\ 0\ 1\ 1\ 1\ 1\ 1)$ | = | $Y_2$ say |
| $Y_2F^T$ | ↪ | $(1\ 1\ 1\ 0\ 1\ 0)$ | = | $X_3$ say |
| $X_3\,F$ | ↪ | $(1\ 1\ 1\ 1\ 0\ 1\ 1\ 1\ 1\ 1)$ | = | $Y_3$ |
| $Y_3F^T$ | ↪ | $(1\ 1\ 1\ 1\ 1\ 0)$ | = | $X_4$ |
| $X_4F$ | ↪ | $(1\ 1\ 1\ 1\ 0\ 1\ 1\ 1\ 1\ 1)$ | = | $Y_4$ say |
| $Y_4F^T$ | ↪ | $(1\ 1\ 1\ 1\ 1\ 0)$ | = | $X_5 = X_4.$ |

Thus we get a binary pair, which is a fixed point of the dynamical system i.e., $\{(1\ 1\ 1\ 1\ 1\ 0), (1\ 1\ 1\ 1\ 0\ 1\ 1\ 1\ 1\ 1)\}$. Thus because of the government failure to stop mislead agricultural technique all nodes except $D_6$ in the domain space and in the range space $R_5$ and all other states are onto the on. For more please use the program in C given in the appendix 7 of the book and work for the resultant vectors. However we have worked with several state vectors to obtain the conclusions given in chapter VII.

To study the HIV/AIDS and its relation with migrant labour we first enlist the Risk/ problems faced by the migrant labourers in the new place where he has temporarily migrated for job are:

| | | |
|---|---|---|
| $M_1$ | – | Bad company and Bad Habits. |
| $M_2$ | – | CSWs. |
| $M_3$ | – | No proper union for the migrant labourers. |
| $M_4$ | – | Living in totally different set up (different from the village life. |
| $M_5$ | – | lack of education. |
| $M_6$ | – | Away from family for work. |
| $M_7$ | – | No fear of being observed. |
| $M_8$ | – | No higher motivation in life. |

$M_9$    –    No way to engage themselves after work hours.

$M_{10}$   –    Their place of stay after work hours in very inconducive.

Let us consider these ten attributes $M_1 \ldots, M_{10}$ as the nodes related with the domain space. We describe each of the attributes now.

## $M_1$ – Bad company and Bad Habits

The migrants from rural areas become easy victims to bad company and bad habits. They after work hours take up to alcohol as it is a way to help them to go to sleep as the family members are away also the interview showed of the 60 interviewed 58 of them were drunkards i.e., they drink everyday and they took proud in saying so. Thus they being already addicted to bad habits in the absence of family members become more prone to them

## $M_2$ – CSWs

The sadness about Tamil Nadu and especially cities in it like Chennai we do not have a red light area like Bombay the CSWs are present all over the city and the seekers of them know their places and rates. Thus these uneducated rural men seek after them after working hours, as they are available for cheap rates. So it is very easy for them to visit CSWs for sexual satisfaction mainly when they are away from their family for weeks.

## $M_3$ – No proper union for the migrant labourers

The migrant labourers may be sales people, working as cleaners or servers in hotels, truck drivers, construction labourers, bore pump labourers, labourer who lay road, etc. These labourers do not have any form of job security or job related protection they do all types of even risky jobs for money as they are in dire need of some job to keep up their living. They do not have union which can raise their voice at time of problems or any form of counselling or help as the types of migrant labourers is heterogeneous not of a particular type. Further government also does not show any interest about these people for they are very poor, from depressed classes and uneducated so they show no concern over them. These people come for their very living.

**M₄ – Living in totally different set up**

96% of the migrant labourers are from very remote or rural areas with no basic education. The life in the villages are entirely different from the city for, in city life no one ever knows what is happening to the neighbour or who is their neighbour, Thus it has become for the rural people to feel very lonely and left out in city after work hours so they seek for alcohol and CSWs. Thus the daring distinction in the life style makes several of the migrant labourers to be more prone to all sorts of bad habits ultimately leaving them as HIV/AIDS affected ones.

**M₅ – Lack of education**

Of the 60 migrant labourers we have interviewed at least 3% had never entered school premises and 80% only studied less than $8^{th}$ std and they lack education. Even those who said they had studied up to $5^{th}$ std or $6^{th}$ std or $7^{th}$ std. humbly acknowledged they don't know to read and have forgotten so they do not have any form of education. They were totally ignorant of the disease HIV/AIDS until they became the victims of it. But it is very important to note that the city people are very well aware of the disease and are very careful to take care of themselves that they don't become affected by HIV/AIDS in contrast to the rural uneducated migrant labourers.

**M₆ – Away from family**

These migrant labourer are away from the family for weeks; they do not have any enterprise after the work hours. In most cases the work hours end by six in the evening or even earlier. They after this time are fully free and do not and cannot watch T.V. or engage themselves in any form of productive occupation the majority of them drink alcohol as it is easily available now as government sells it and go for CWS; for they have least fear of being watched for no one ever notices them in the city.

**M₇ – No fear of being observed**

The migrant labourers go to CSWs when they are totally drunk. They think they can do anything for no one can question them. Thus we see the migrant labourers who live in city for some days or weeks are prone to all bad habits and at onc stage or other become socially irresponsible. When we say they become socially irresponsible elements we wish to state that they not only infect themselves but their family members with HIV/AIDS by their careless and reckless ways of living; like visiting CSWs in a drunken condition and not following any protected sex. We say they are socially irresponsible for from the interviews we come to know that of the 60 migrant labourers whom we interviewed all the 58 of them took pride in saying they had visited countless number of CSWs and they took further pride in saying that they had all habits like consuming alcohol, continuous smokers and some adids of ganga or drugs. Thus when we say they have no fear of being observed by other we mean all these, which has come from live interviews from them.

### $M_8$ – No higher motivation in Life

When we say these HIV/AIDS affected migrant labourers have no higher motivation in life, for they do not have aspiration of educating their children or providing a better status to their family members but on the contrary what they earn they spend on themselves that too not on good cloths are comfort for them but they spend on smoke, alcohol and on the CSWs. Thus these HIV/AIDS patients who are migrant labourers live for the day enjoy the day they don't take any effort in improving the family or their children. They speak of their poverty, if poverty is acute how do they have mind to spend it in wrong means. Thus we have put this tile after lot of contemplation and from the interviews.

### $M_9$ – No way to engage themselves after work hours

Most of them complain of their inability to do anything after work hours what they mainly do after work hours is they self sympathize themselves by saying they have worked for the day so to forget the body pain at the first stage they consume liquor. Once they are drunk invariably they are exposed to CSWs either by their friends or even in the place where they take the drink. Thus they being unaware of how the HIV/AIDS spreads they become easy victims of the same. So their only recreation after work hours is CSWs. This was clearly acknowledged by 60% of them for the other migrant labourers whom we had interviewed were Truck drivers who have sex on the high ways from the CSWs posted for the same purpose of trapping these truck drivers.

**M$_{10}$ – Their place of stay after work hours is very inconducive**

They stay after work hours in places where no comfort or any form of living by any human is possible so some say they drink to get sleep. The sanitary conditions are still poorer. They pay a rent of Rs.50 to Rs.150 a month for their stay. Since the place of stay after work hours is very uncomfortable they spend out. They are not able to think that they can pay a better rent and stay in a better place by using the money they spend on CSWs, smoke and drink (alcohol). They feel that these 3 are impossible to be sacrificed in comparison with their place of stay. So they prefer in stay in such dingy corners in the midst of mosquitoes and with no proper sanitation.

We take a few attributes related with the government.

G$_1$ –   Lack of fore sight in the part of government to provide job opportunities for agricultural labourers when they have no job or failure of monsoon.

G$_2$ –   No plans of government ever reach the very poor and the uneducated rural people they survive only by physical labour.

G$_3$ –   Till date government has not given any form protection /incentive / adult education / training in other handicrafts of the migrant labourers. In most cases the migrant labourers are victims of harassment, ill treatment by the rich farmers/employer.

G$_4$ –   Easy availability of liquor.

G$_5$ –   Availability of CSWs at very a cheap rate every where in the city.

G$_6$ –   Government's failure to take any steps to form a support group or mobilization of them to empower them leading to better living conditions and protection from HIV/AIDS.

G$_7$ –   Awareness program about HIV/AIDS has never reached these people so they visit CSWs without fear and is one of the causes of infecting even the CSWs.

We take the seven attributes $G_1$, $G_2$,…, $G_7$ as those nodes of the range space and $M_1$, $M_2$, …, $M_{10}$ as the nodes of the domain space and obtain the directed graph of the FRM using an experts opinion.

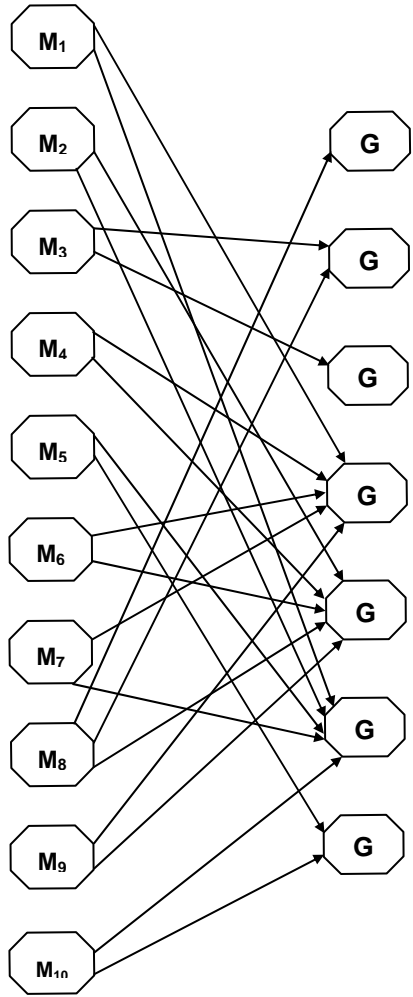

FIGURE: 5.1.3

Using the directed graph given by this expert we obtain the related connections matrix. Let us denote the $10 \times 7$ connection matrix of the FRM by $F_1$

$$
F_1 = \begin{array}{c} \\ M_1 \\ M_2 \\ M_3 \\ M_4 \\ M_5 \\ M_6 \\ M_7 \\ M_8 \\ M_9 \\ M_{10} \end{array}
\begin{array}{c} G_1\ G_2\ G_3\ \ G_4\ \ G_5\ G_6\ G_7 \\
\begin{bmatrix}
0 & 0 & 0 & 1 & 0 & 1 & 0 \\
0 & 0 & 0 & 0 & 1 & 1 & 0 \\
0 & 1 & 1 & 0 & 0 & 0 & 0 \\
0 & 0 & 0 & 1 & 1 & 0 & 0 \\
0 & 0 & 0 & 0 & 0 & 1 & 1 \\
0 & 0 & 0 & 1 & 1 & 0 & 0 \\
0 & 0 & 0 & 1 & 0 & 1 & 0 \\
1 & 1 & 0 & 0 & 1 & 0 & 0 \\
0 & 0 & 0 & 1 & 1 & 0 & 0 \\
0 & 0 & 0 & 0 & 0 & 1 & 1
\end{bmatrix}
\end{array}
$$

Suppose we assume the stability of the dynamical system. Now we study the effect of the state vector on the system $F_1$ taken from the domain space or from the range space.

Suppose X = (0 0 0 0 0 0 0 1 0 0) is take as the state vector only the vector a migrant uneducated labouer has no higher motivation about life and all other states are off. We now study effect of X on the dynamical system F and also the effect of X on the range space.

$$XF_1 \quad \hookrightarrow \quad (1\ 1\ 0\ 0\ 1\ 0\ 0) \qquad = \quad Y$$

$$YF^T \quad \hookrightarrow \quad (0\ 1\ 1\ 1\ 0\ 1\ 0\ 1\ 1\ 0\,) = \quad X_1$$

say '$\hookrightarrow$' denote the vector is updated and thresholded at each stage.

$$X_1F_1 \quad \hookrightarrow \quad (1\ 1\ 1\ 1\ 1\ 1\ 0) \qquad = \quad Y_1$$

$$Y_1F_1{}^T \quad \hookrightarrow \quad (1\ 1\ 1\ 1\ 1\ 1\ 1\ 1\ 1\ 1) \quad = \quad X_2 = (X_1)$$

$$X_2F \quad \hookrightarrow \quad (\,1\ 1\ 1\ 1\ 1\ 1\ 1\,) \qquad = \quad Y_2$$

$$Y_2F^T \quad \hookrightarrow \quad (1\ 1\ 1\ 1\ 1\ 1\ 1\ 1\ 1\ 1) \quad = \quad X_3 = \quad X_2.$$

Thus the hidden pattern of the dynamical system is a fixed point, when the migrant labour has no higher values life all the nodes in both the domain and range space become on given by the pair {(1 1 1 1 1 1 1 1 1 1), (1 1 1 1 1 1 1)}. Which show the importance of the impact of this node on other nodes.

Now we study the state vector P = (1 0 0 0 0 0 0 0 0 1) i.e., the states vector denotes the on state of the attributes $M_1$ and $M_{10}$ and all other nodes remain off. The effect of P on the dynamical system F is given by

$$PF_1 \quad \hookrightarrow \quad (0\ 0\ 0\ 1\ 0\ 1\ 1) \quad = \quad Q\ (say)$$

$$Q_1F_1^T \quad \hookrightarrow \quad (1\ 1\ 0\ 1\ 1\ 1\ 1\ 1\ 1\ 1) \quad = \quad P_1\ (say)$$

$$P_1F_1 \quad \hookrightarrow \quad (1\ 1\ 0\ 1\ 1\ 1\ 1) \quad = \quad Q_1\ (say)$$

$$Q_1\ F_1^T \quad \hookrightarrow \quad (1\ 1\ 1\ 1\ 1\ 1\ 1\ 1\ 1\ 1) \quad = \quad P_2\ (say)$$

$$P_2\ F_1 \quad \hookrightarrow \quad (1\ 1\ 1\ 1\ 1\ 1\ 1) \quad = \quad Q_2\ (say)$$

$$Q_2F_1^T \quad \hookrightarrow \quad (1\ 1\ 1\ 1\ 1\ 1\ 1\ 1\ 1\ 1) \quad = \quad P_3\ (P_2).$$

Thus the resultant binary pair {(1 1 1 1 1 1 1 1 1 1), (1 1 1 1 1 1 1)} is a fixed point of the system. We see all nodes both in the domain space come to on state when the migrant labourer has bad habits, bad company and the place of stay of work hours is inconducive.

Now we consider the effect of the state vector of the range space on the dynamical system $F_1$. Let Z = (0 1 0 0 0 0 0) i.e., only the attribute $G_2$ is in the on state and all other nodes are in the off state. The effect of Z on the system $F_1$ is given by

$$ZF_1^T \quad \hookrightarrow \quad (0\ 0\ 1\ 0\ 0\ 0\ 0\ 1\ 0\ 0) \quad = \quad A\ (say)$$

$$A\ F_1 \quad \hookrightarrow \quad (1\ 1\ 1\ 0\ 1\ 0\ 0) \quad = \quad Z_1\ (say)$$

$$Z_1\ F_1^T \quad \hookrightarrow \quad (0\ 1\ 1\ 1\ 0\ 1\ 0\ 1\ 1\ 0) \quad = \quad A_1\ say$$

$$A_1F_1 \quad \hookrightarrow \quad (1\ 1\ 1\ 1\ 1\ 1\ 0) \quad = \quad Z_2(say)$$

$$Z_2\ F_1^T \quad \hookrightarrow \quad (1\ 1\ 1\ 1\ 1\ 1\ 1\ 1\ 1\ 1) \quad = \quad A_2\ say$$

$$A_2F_1 \quad \hookrightarrow \quad (1\ 1\ 1\ 1\ 1\ 1\ 1) \quad = \quad Z_3(say)$$

$$Z_3\ F_1^T \quad \hookrightarrow \quad (1\ 1\ 1\ 1\ 1\ 1\ 1\ 1\ 1\ 1) \quad = \quad Z_3\ say.$$

Thus we see the hidden pattern of the system is a fixed point. Thus the fixed binary pair is {(1 1 1 1 1 1 1 1 1 1), (1 1 1 1 1 1 1)}. Several such illustrations can be obtained using different sets of state vectors and conclusions derived on the resultant vector. Also we can make use of the C program given in the appendix of the book so that the calculations can be made further simple.

Next we seek the opinion of another expert for the same data so that after obtaining say 4 experts opinion we would use the combined FRM and derive conclusions based on these experts opinion. Now we consider some four experts opinion over the same set of attributes given by the domain space as $M_1$, $M_2$, … , $M_{10}$ and that of the range space R given by $G_1$, $G_2$, … , $G_7$. We give the directed graph given by the expert using the set of attributes $\{M_1, M_2, …, M_{10}\}$ and $\{G_1, G_2,…, G_7\}$ in the following:

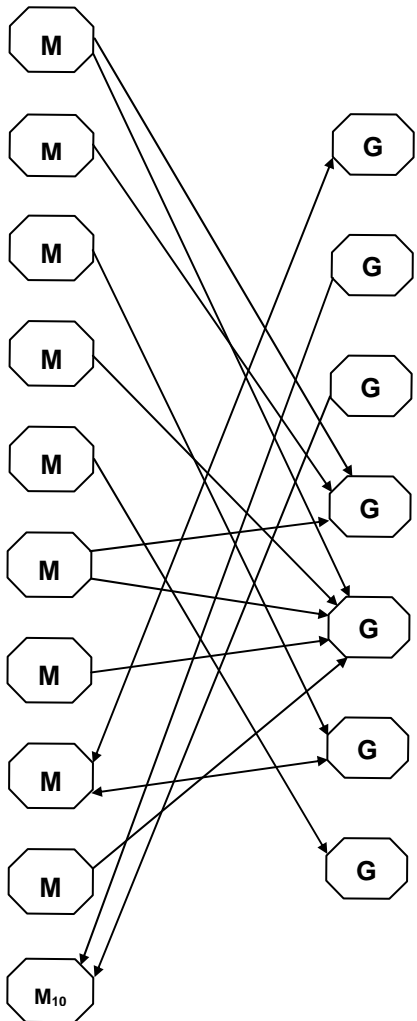

FIGURE: 5.1.4

The related $10 \times 7$ matrix $F_2$ of the graph is us following

$$F_2 = \begin{array}{c} \\ M_1 \\ M_2 \\ M_3 \\ M_4 \\ M_5 \\ M_6 \\ M_7 \\ M_8 \\ M_9 \\ M_{10} \end{array} \begin{array}{ccccccc} G_1 & G_2 & G_3 & G_4 & G_5 & G_6 & G_7 \\ \left[\begin{array}{ccccccc} 0 & 0 & 0 & 1 & 1 & 0 & 0 \\ 0 & 0 & 0 & 1 & 0 & 0 & 0 \\ 0 & 0 & 0 & 0 & 0 & 1 & 0 \\ 0 & 0 & 0 & 0 & 1 & 0 & 0 \\ 0 & 0 & 0 & 0 & 0 & 0 & 1 \\ 0 & 0 & 0 & 1 & 1 & 0 & 0 \\ 0 & 0 & 0 & 0 & 1 & 0 & 0 \\ 1 & 0 & 0 & 0 & 0 & 1 & 0 \\ 0 & 0 & 0 & 0 & 1 & 0 & 0 \\ 0 & 1 & 1 & 0 & 0 & 0 & 0 \end{array}\right] \end{array}$$

Now we study the effect of the dynamical system $F_2$ on the state vectors.

Suppose we have to the find the effect of the state vector X = (0 0 0 0 0 0 0 1 0 0) i.e., node attribute $M_8$ alone is in the on state and all other nodes are in the off state. The effect of X on the dynamical system $F_2$ is given by

$$XF_2 \quad \hookrightarrow \quad (1\ 0\ 0\ 0\ 0\ 1\ 0) \quad = \quad Y \text{ (say)}$$
$$YF_2^T \quad \hookrightarrow \quad (0\ 0\ 1\ 0\ 0\ 0\ 0\ 1\ 0\ 0) = X_1 \text{ say}$$
$$X_1F_2 \quad \hookrightarrow \quad (1\ 0\ 0\ 0\ 0\ 1\ 0) \quad = \quad Y_1 \text{ (say)}$$
$$Y_1 F_2^T \quad \hookrightarrow \quad (0\ 0\ 1\ 0\ 0\ 0\ 0\ 1\ 0\ 0) = X_2 = X_1.$$

Thus the hidden pattern of the dynamical system is a fixed point given by the binary pair {(0 0 1 0 0 0 0 1 0 0), (1 0 0 0 0 1 0)} we see the on state of $M_8$ makes the on state of $G_1$, $G_6$ and $M_3$.

Now we study the effect of the state vector P = (1 0 0 0 0 0 0 0 0 0) i.e., only the nodes $M_1$ is in the on state and all other nodes are in the off state. The effect of P on the dynamical system $F_2$ is given by

$$PF_2 \quad \hookrightarrow \quad (0\ 0\ 0\ 1\ 1\ 0\ 0) \quad = \quad Q$$

$$QF_2^T \quad \hookrightarrow \quad (1\ 1\ 0\ 1\ 0\ 1\ 1\ 0\ 1\ 0) \quad = \quad P_1$$
$$P_1\ F_2 \quad \hookrightarrow \quad (0\ 0\ 0\ 1\ 1\ 0\ 0) \quad\quad = \quad Q_1 = Q_1$$
$$Q_1\ F_2^T \quad \hookrightarrow \quad (1\ 1\ 0\ 1\ 0\ 1\ 1\ 0\ 1\ 0) \quad = \quad P_1.$$

Thus the fixed point is a binary pair given by {(1 1 0 1 0 1 1 0 1 0), (0 0 0 1 1 0 0)} that is the effect of bad company and bad habits makes the nodes $M_2$, $M_4$, $M_6$, $M_7$ and $M_9$ on state in the domain space and this node has no impact on no proper union, lack of education, no higher motivation and incondusive place of stay.

We have worked with several state vector to arrive at the conclusions. Now we consider the state vector B = (1 0 0 0 0 0 0) where only the node $G_1$ in the on state and all other nodes are in the off state. The effect of B on the dynamical system $F_2$ is given by

$$BF_2^T \quad \hookrightarrow \quad (0\ 0\ 0\ 0\ 0\ 0\ 0\ 1\ 0\ 0) \quad = \quad A.$$
$$AF_2 \quad \hookrightarrow \quad (1\ 0\ 0\ 0\ 0\ 1\ 0) \quad = \quad B_1$$
$$B_1\ F_2^T \quad \hookrightarrow \quad (1\ 0\ 1\ 0\ 0\ 0\ 0\ 1\ 0\ 0) \quad = \quad A_1$$
$$A_1\ F_2 \quad \hookrightarrow \quad (1\ 0\ 0\ 1\ 1\ 1\ 0) \quad = \quad A_2$$
$$B_2\ F_2^T \quad \hookrightarrow \quad (1\ 1\ 1\ 1\ 0\ 1\ 1\ 1\ 1\ 0) \quad = \quad A_2$$
$$A_2\ F_2 \quad \hookrightarrow \quad (1\ 0\ 0\ 1\ 1\ 1\ 0) \quad\quad = \quad B_3 = B_2$$

$$B_3\ F_2^T \quad \hookrightarrow \quad\quad A_3 \quad (= A_2).$$

Thus the hidden pattern is a fixed point of the dynamical system. It is given by the binary pair {(1 1 1 1 0 1 1 1 1 0), (1 0 0 1 1 1 0)}.

The node lack of foresight is so powerful it has made all nodes in the domain space on expert lack of education.

Now using the same set of attributes we once again calculate the effect of state vectors using a third experts opinion. We take the domain space D = {$M_1$, $M_2$,…, $M_{10}$} and the range space R = {$G_1$, $G_3$,…, $G_7$}.

The directed graph as given by the expert is as follows:

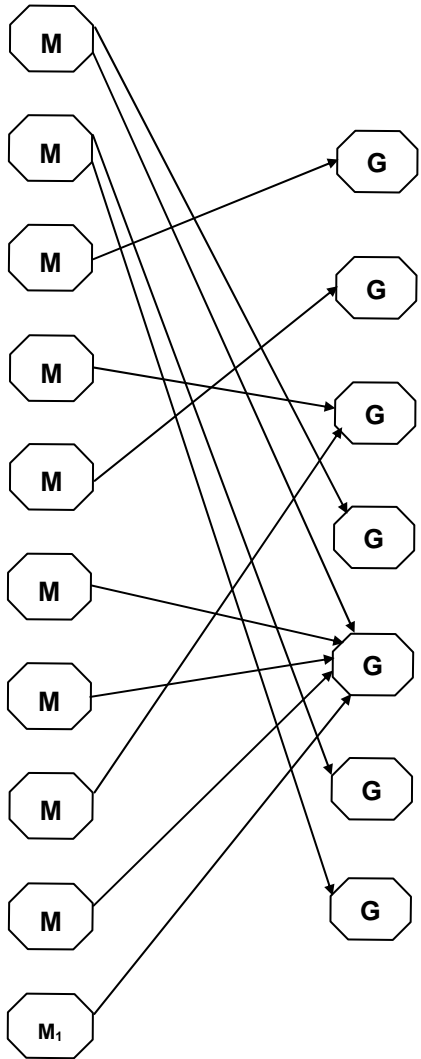

FIGURE: 5.1.5

The connection matrix $F_3$ associated with the directed graph given above is as follows

$$
\begin{array}{c}
\phantom{F_2 = M_1}\begin{array}{ccccccc} G_1 & G_2 & G_3 & G_4 & G_5 & G_6 & G_7 \end{array} \\
F_2 = \begin{array}{c} M_1 \\ M_2 \\ M_3 \\ M_4 \\ M_5 \\ M_6 \\ M_7 \\ M_8 \\ M_9 \\ M_{10} \end{array}
\begin{bmatrix}
0 & 0 & 0 & 1 & 1 & 0 & 0 \\
0 & 0 & 0 & 0 & 0 & 1 & 1 \\
1 & 0 & 0 & 0 & 0 & 0 & 0 \\
0 & 0 & 1 & 0 & 0 & 0 & 0 \\
0 & 1 & 0 & 0 & 0 & 0 & 1 \\
0 & 0 & 0 & 0 & 1 & 0 & 0 \\
0 & 0 & 0 & 0 & 1 & 0 & 0 \\
0 & 0 & 1 & 0 & 0 & 0 & 0 \\
0 & 0 & 0 & 0 & 1 & 0 & 0 \\
0 & 0 & 0 & 0 & 1 & 0 & 0
\end{bmatrix}
\end{array}
$$

$F_3$ is a $10 \times 7$ matrix consider any state vector $X = (0\ 0\ 0\ 0\ 0\ 1\ 0\ 0\ 0\ 0)$ where only the attribute $M_6$ is in the on state and all other vectors are in the off state the effect of X on the dynamical system $F_3$

$$
\begin{aligned}
XF_3 &\hookrightarrow (0\ 0\ 0\ 0\ 1\ 0\ 0) & &= & Y \\
YF_3^{T} &\hookrightarrow (1\ 0\ 0\ 0\ 0\ 1\ 1\ 0\ 1\ 1) & &= & X_1 \text{ say} \\
X_1 F_3 &\hookrightarrow (0\ 0\ 0\ 1\ 1\ 0\ 0) & &= & Y_1 \text{ say} \\
Y_1 F_3^{T} &\hookrightarrow (1\ 0\ 0\ 0\ 0\ 1\ 1\ 0\ 1\ 1) & &= & X_2 = X.
\end{aligned}
$$

Thus we see the hidden pattern of the dynamical system is a fixed point given by the binary pair $\{(1\ 0\ 0\ 0\ 0\ 1\ 1\ 0\ 1\ 1), (0\ 0\ 0\ 1\ 1\ 0\ 0)\}$. When the migrant labourer is away from the family is in the on state we see the nodes $M_1$, $M_7$, $M_9$ and $M_{10}$ come to on state and all other nodes in the domains are in the off state in the range space for easy availability of liquor and CSWs at a cheap rate.

Next we consider the state vector $Y = (0\ 0\ 0\ 0\ 1\ 0\ 0)$ from the range space. Effect of Y on the dynamical system $F_3$ is given by

$$
\begin{aligned}
YF_3^{T} &\hookrightarrow (1\ 0\ 0\ 0\ 0\ 1\ 1\ 0\ 1\ 1) & &= & X \text{ (say)} \\
X F_3 &\hookrightarrow (0\ 0\ 0\ 1\ 1\ 0\ 0) & &= & Y_1 \text{ say} \\
Y_1 F_3^{T} &\hookrightarrow (1\ 0\ 0\ 0\ 0\ 1\ 1\ 0\ 1\ 1) & &= & X_1 \text{ (Say)} \\
X_1 F_3 &\hookrightarrow (0\ 0\ 0\ 1\ 1\ 0\ 0) & &= & Y_2 (= Y_1 \text{ say}).
\end{aligned}
$$

Thus the hidden pattern of the dynamical system is a fixed point the binary $\{(1\ 0\ 0\ 0\ 0\ 1\ 1\ 0\ 1\ 1), (0\ 0\ 0\ 1\ 1\ 0\ 0)\}$. We see when $G_5$ that is the availability of CSWs at a very cheap rate is in the on state we see the node $G_4$ comes to the on state that the cheap availability of liquor comes to force. Thus CSWs and drinking alchol goes hand in hand. Further the nodes $M_1$, $M_6$, $M_7$, $M_9$ and $M_{10}$ come to on state, that is the migrants seek bad company, they are away from their family, no fear of being obsevered, no way of engaging after work hours and their place of stay is incondusive.

Now we seek the fourth expert opinion on the same set of attributes D = {$M_1$, $M_2$, $M_3$, …, $M_{10}$} and R = {$R_1$, $R_2$, $R_3$,…, $R_7$}.

The directed graph of the given system is as follows:

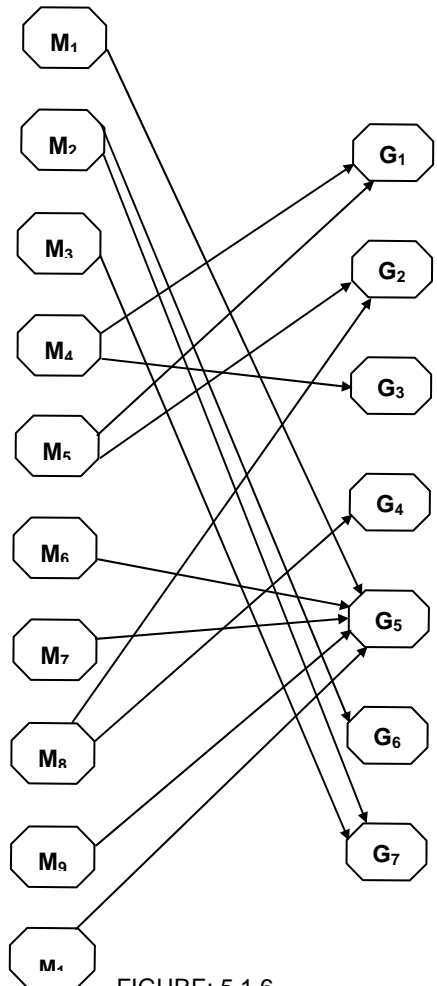

FIGURE: 5.1.6

The related connection matrix is given by $F_4$

$$
\begin{array}{c@{\;}c}
 & \begin{array}{ccccccc} G_1 & G_2 & G_3 & G_4 & G_5 & G_6 & G_7 \end{array} \\
\begin{array}{c} M_1 \\ M_2 \\ M_3 \\ M_4 \\ M_5 \\ M_6 \\ M_7 \\ M_8 \\ M_9 \\ M_{10} \end{array} &
\left[\begin{array}{ccccccc}
0 & 0 & 0 & 0 & 1 & 0 & 0 \\
0 & 0 & 0 & 0 & 0 & 1 & 1 \\
0 & 0 & 0 & 0 & 0 & 0 & 1 \\
1 & 0 & 1 & 0 & 0 & 0 & 0 \\
1 & 1 & 0 & 0 & 0 & 0 & 0 \\
0 & 0 & 0 & 0 & 1 & 0 & 0 \\
0 & 0 & 0 & 0 & 1 & 0 & 0 \\
0 & 1 & 0 & 1 & 0 & 0 & 0 \\
0 & 0 & 0 & 0 & 1 & 0 & 0 \\
0 & 0 & 0 & 0 & 1 & 0 & 0
\end{array}\right]
\end{array}
$$

Let us consider the state vector $X = (0\ 0\ 0\ 0\ 0\ 0\ 0\ 1\ 0\ 0)$ where only the attribute $M_8$ is in the on state and all other attributes are in the off state. The effect of $X$ on the dynamical system $F_4$ is given by

$$
\begin{array}{llll}
XF_4 & \hookrightarrow & (0\ 1\ 0\ 1\ 0\ 0\ 0) & = & Y \text{ say} \\
YF_4{}^T & \hookrightarrow & (0\ 0\ 0\ 0\ 1\ 0\ 0\ 1\ 0\ 0) & = & X_1 \text{ say} \\
X_1 F_4 & \hookrightarrow & (1\ 1\ 0\ 1\ 0\ 0\ 0) & = & Y_1 \text{ say} \\
Y_1 F_4{}^T & \hookrightarrow & (0\ 1\ 0\ 1\ 1\ 0\ 0\ 1\ 0\ 0) & = & X_2 \text{ say} \\
X_2 F_4{}^T & \hookrightarrow & (1\ 1\ 1\ 1\ 0\ 1\ 1) & = & Y_2 \text{ say} \\
Y_2 F_4{}^T & \hookrightarrow & (0\ 1\ 0\ 1\ 1\ 0\ 0\ 1\ 0\ 0) & = & X_3\ (X_3 = X_2) \\
X_3 F_4 & \hookrightarrow & (1\ 1\ 1\ 1\ 0\ 1\ 1) & = & Y_3\ (= Y_2).
\end{array}
$$

Thus the hidden pattern of the dynamical system is a fixed point given by the binary pair $\{(0\ 0\ 0\ 0\ 1\ 0\ 0\ 1\ 0\ 0), (1\ 1\ 1\ 1\ 0\ 1\ 1)\}$.

We see when no higher motivation in life is in the on state we see the nodes $M_2$, $M_4$ and $M_5$ come to on state there by indicating they visits CSWs, they lack education and live on a totally different set up.

Also $G_1$, $G_2$, $G_3$, $G_4$, $G_5$ and $G_6$ come to on state that is government lacks foresight, it has no plans to help the poor,

availability of CSWS and liquor and its failure to spread awareness.

Now using the three experts opinion relative to the matrices $F_2$, $F_3$ and $F_4$ we obtain the combined FRM. The Combined FRM (CFRM), F is given by $F = F_2 + F_3 + F_4$

$$F = \begin{array}{c}  \\ M_1 \\ M_2 \\ M_3 \\ M_4 \\ M_5 \\ M_6 \\ M_7 \\ M_8 \\ M_9 \\ M_{10} \end{array} \begin{array}{c} G_1 \ G_2 \ G_3 \ G_4 \ G_5 \ G_6 \ G_7 \\ \begin{bmatrix} 0 & 0 & 0 & 2 & 3 & 0 & 0 \\ 0 & 0 & 0 & 1 & 0 & 2 & 2 \\ 1 & 0 & 0 & 0 & 0 & 1 & 1 \\ 1 & 0 & 2 & 0 & 1 & 0 & 0 \\ 1 & 2 & 0 & 0 & 0 & 0 & 1 \\ 0 & 0 & 0 & 1 & 2 & 0 & 0 \\ 0 & 0 & 0 & 1 & 2 & 0 & 0 \\ 1 & 1 & 1 & 1 & 0 & 1 & 0 \\ 0 & 0 & 0 & 0 & 3 & 0 & 0 \\ 0 & 1 & 1 & 0 & 2 & 0 & 0 \end{bmatrix} \end{array}$$

Using the dynamical system F of the CFRM we find the effect of each of the state vector. Also we have given in the appendix a C-program to work for the stability of the system i.e., finidng the hidden pattern be it be a fixed point or a limit cycle.

Now let us consider the state vector

$$X = (0\ 0\ 0\ 0\ 0\ 0\ 0\ 1\ 0\ 0)$$

i.e., only the node $M_8$ is in the on state and all others nodes are in the off state. The effect of X on the dynamical system F is given by

| | | | | |
|---|---|---|---|---|
| XF | $\hookrightarrow$ | (1 1 1 1 0 1 0) | = | Y say |
| $YF^T$ | $\hookrightarrow$ | (1 1 1 1 1 1 1 1 0 1) | = | $X_1$ say |
| $X_1 F$ | $\hookrightarrow$ | (1 1 1 1 1 1 1) | = | $Y_1$ say |
| $Y_1F^T$ | $\hookrightarrow$ | (1 1 1 1 1 1 1 1 1 1). | | |

Thus using the CFRM the hidden pattern of the dynamical system s a fixed point given by the binary pair {(1 1 1 1 1 1 1 1 1 1 1), (1 1 1 1 1 1 1)}. Thus all the nodes become on both from the range and the domain space.

This shows in the combined opinion all states become on when only the attribut $M_8$ is in the on state. Thus the only node no higher motivation in life leads to the on state of all attributes and the government to motivate them.

Let us consider the state vector A = (1 0 0 0 0 0 0) where only the attribute $G_1$ is in the on state and all other vectors are in the off state. Now we study the effect of A on the CFRM i.e.,

$$AF^T \quad \hookrightarrow \quad (0\ 0\ 1\ 1\ 1\ 0\ 0\ 1\ 0\ 0) \quad = \quad B \text{ say}$$

$$BF \quad \hookrightarrow \quad (1\ 1\ 1\ 1\ 1\ 1\ 1) \quad\quad = \quad A_1 \text{ say}$$

$$A_1\,F^T \quad \hookrightarrow \quad (1\ 1\ 1\ 1\ 1\ 1\ 1\ 1\ 1\ 1).$$

Thus when the only attribute $G_1$ is in the on state it gives the hidden pattern of the system to be a fixed point which is given by the binary pair {(1 1 1 1 1 1 1 1 1 1 1), (1 1 1 1 1 1 1)}.

Thus we see when lack of foresight on the part of the government to provide job opportunities for agricultural labourer leads to the on state of all attributes in both the domain and range space.

Now we study the social problems faced by a HIV/AIDS affected migrant labourer. It is unfortunate to state at least 80% of these HIV/AIDS affected migrant labourers were first the worst victims of V.D. and other sexually transmitted diseases. Thus even before they were tested to have HIV/AIDS they had taken treatment in different places by different doctors.

By the time they were confirmed to have HIV/AIDS most of them at least 75% have totally or mostly drained all their resources for taking treatment and medical aid. Several of them addressed their doctors in a very abusive language for they felt that the doctors have cheated them.

Further the prolonged delay in diagonizing the disease has made the disease in them chronic and very serious or to be more precise the patients had become very seriously ill and even immobile. Thus certainly a few died just even before they took any proper treatment for HIV/AIDS.

Now we study the basic reason for this are two fold; one the doctors take this as an opportune movement to exploit them monetarily though they are well aware of the fact that these patients suffer from HIV/AIDS. The second part even if some of the doctors would have informed them they fearing their close relatives and their social status in the village would not have disclosed the fact but suffer silently and lonely to death until they became very serious.

Also HIV/AIDS is a disease, which does not make the sufferer die soon but it takes sufficient time in making the patient invalid and very serious for even over six months. Thus fearing society they don't even disclose about the disease to their wife and they infect their wife also. Thus the chain of infection continues.

It is important to know in several of the villages if they come to know some one in their village is affected by HIV/AIDS the whole family is denied water, provision ration even say simple cup of lea from a tea shop the stigma is so powerful they hide the disease from the public and from the family due to the fear of segregation.

This mainly shows that if the awareness program among the rural uneducated people is powerful certainly they will not make the mistake of getting the disease through CSWs, even granted they visited the CSWs they would certainly adopt safe sex methods due to fear of being infected by the disease. Thus we feel that awareness program to reach the illiterate would alone will stop the spread of HIV/AIDS in rural uneducated areas. For of the 60 people interviewed only one said he knew about HIV/AIDS before he was infected.

All the migrant labourers were infected by HIV/AIDS only through sexual transmissions. Thus a prevention program is not a very difficult one. It is not a heterogeneous problem in other countries where spread in via several ways. It is surprisingly noted that these patients themselves acknowledged that they were advised by CSWs to use safe sex methods only they did not follow it.

So it has become true that several of the CSWs were spreading awareness about HIV/AIDS but only these uneducated rural poor migrant labourers became victims due to ignorance and arrogance.

We say arrogance for they themselves acknowledged that they did not listen to her (CSW). Some of them said they will start to respect women; some said for all their acts they are now undergoing the punishment.

Most of them had become so invalid and poor; they were fully submissive and in fact depressed and dejected. Several of them recalled how they wasted money and now how they do not have money even to buy and drink a cup of tea.

Now let us consider the attributes leading a HIV/AIDS migrant patients to socially behave irresponsible.

The attributes, which promote the migration and migrant labourers becoming victim of HIV/AIDS, is taken as a domain space.

$M_1$ - Government help never reaches the rural poor illiterate

$M_2$ - Availability of cheap liquor

$M_3$ - Cheap availability of CSWs

$M_4$ - No job opportunities

$M_5$ - No proper health center

$M_6$ - No school even for primary classes

$M_7$ - Poverty

$M_8$ - No proper road or bus facilities

$M_9$ - Living conditions questionably poor

$M_{10}$ - Government unconcern over their development in any plans so only introduction of machine for harvest etc has crippled agricultural labourer.

Now we proceed onto give the attributes related with the migrant labourers who suffer and are ultimate victims of HIV/AIDS.

$P_1$ - No education / No help by government

$P_2$ - Awareness program never reaches them

$P_3$ - No responsibility of parents to educate children

$P_4$ - Marriage age of girls very low say even at 11 years they are married

$P_5$ - addiction to cheap liquor

$P_6$ - Smoke and visit CSWs

$P_7$ - Very questionable living condition so migrate.

Now we take three experts opinion and obtain the directed graph of the FRM so that making use of these three experts we can find the combined FRM and its effect on state vectors to draw conclusions.

The directed graph given by the first expert for the FRM

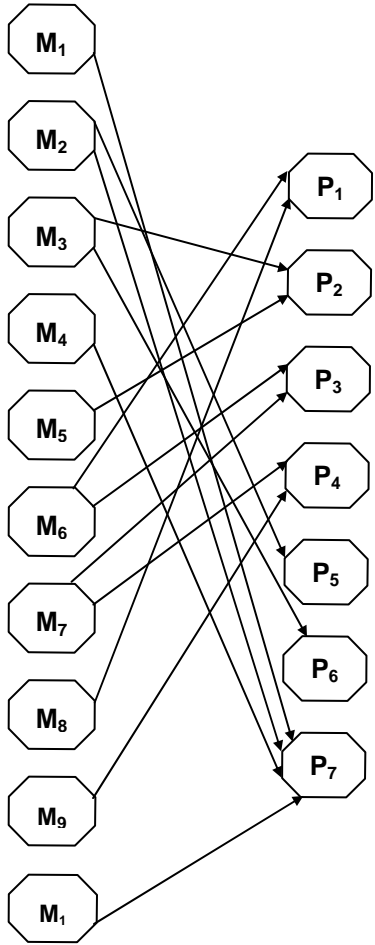

FIGURE: 5.1.7

The related $10 \times 7$ connection matrix of the above given directed graph is denoted by $G_1$ is as follows

$$G_1 = \begin{array}{c} \\ M_1 \\ M_2 \\ M_3 \\ M_4 \\ M_5 \\ M_6 \\ M_7 \\ M_8 \\ M_9 \\ M_{10} \end{array} \begin{array}{ccccccc} P_1 & P_2 & P_3 & P_4 & P_5 & P_6 & P_7 \\ \left[\begin{array}{ccccccc} 0 & 0 & 0 & 0 & 0 & 0 & 1 \\ 0 & 0 & 0 & 0 & 1 & 1 & 0 \\ 0 & 1 & 0 & 0 & 0 & 1 & 0 \\ 0 & 0 & 0 & 0 & 0 & 0 & 1 \\ 0 & 1 & 0 & 0 & 0 & 0 & 0 \\ 1 & 0 & 1 & 0 & 0 & 0 & 0 \\ 0 & 0 & 1 & 1 & 0 & 0 & 0 \\ 1 & 0 & 0 & 0 & 0 & 0 & 0 \\ 0 & 0 & 0 & 1 & 0 & 0 & 0 \\ 0 & 0 & 0 & 0 & 0 & 0 & 1 \end{array}\right] \end{array}$$

Let us consider the effect of the state vector $X = (0\ 0\ 0\ 0\ 1\ 0\ 0\ 0\ 0\ 0)$ on the dynamical system $G_1$ where $M_5$ alone is in the on state all other vectors are in the off state

$$\begin{array}{llll} XG_1 & \hookrightarrow & (0\ 1\ 0\ 0\ 0\ 0\ 0) & = & Y \\ YG_1^T & \hookrightarrow & (0\ 0\ 1\ 0\ 1\ 0\ 0\ 0\ 0\ 0) & = & X_1 \\ X_1G_1 & \hookrightarrow & (0\ 1\ 0\ 0\ 0\ 1\ 0) & = & Y_1 \\ Y_1\,G_1^T & \hookrightarrow & (0\ 1\ 1\ 0\ 1\ 0\ 0\ 0\ 0\ 0) & = & X_2 \text{ say} \\ X_2\,G_1 & \hookrightarrow & (0\ 1\ 0\ 0\ 1\ 1\ 0) & = & Y_2 \text{ say} \\ Y_2\,G_1 & \hookrightarrow & (0\ 1\ 1\ 0\ 1\ 0\ 0\ 0\ 0\ 0) & = & X_3 = (X_2 \text{ say}). \end{array}$$

The hidden pattern of the dynamical system is a fixed point given by the binary pair $\{(0\ 1\ 0\ 0\ 0\ 1\ 0), (0\ 1\ 1\ 0\ 1\ 0\ 0\ 0\ 0\ 0)\}$. We see when only the node no proper health center is in the on state we see the nodes $M_2$ and $M_3$ come to on state which implies they drink cheap liquor and visit CSWs for they do not know its bad effect on health. Also $P_2$ and $P_6$ come to on state that awareness program never reaches them and they visit CSWs and also smoke. Now suppose we consider the state vector $Y = (0\ 0\ 1\ 0\ 0\ 0\ 0)$ where $P_3$ alone is in the on state and all other nodes are in the off state the effect of $Y$ on the dynamical system $G_1$ is given by

$$\begin{array}{llll} YG_1^T & \hookrightarrow & (0\ 0\ 0\ 0\ 0\ 1\ 1\ 0\ 0\ 0) & = & X_1 \text{ say} \\ X_1\,G_1 & \hookrightarrow & (1\ 0\ 1\ 1\ 0\ 0\ 0) & = & Y_1 \text{ say} \\ Y_1\,G_1^T & \hookrightarrow & (0\ 0\ 0\ 0\ 0\ 1\ 1\ 1\ 1\ 0) & = & X_2 \text{ say} \\ X_2\,G_1 & \hookrightarrow & (1\ 0\ 1\ 1\ 0\ 0\ 0) & = & Y_2\ (Y_1 \text{ say}). \end{array}$$

Thus the hidden pattern of the dynamical system is a fixed point given by the binary pair {(0 0 0 0 0 1 1 1 1 0), (1 0 1 1 0 0 0)}. When only the node no responsibilty on the part of the parents to educate the children is in the on state and all nodes are in the off state we see $P_1$ and $P_4$ becomes on their by indicating the parents themselves have no education/ government does not help them automatically the marriage age of the girl children is very low. The nodes, which becomes on in the domain space are $M_6$, $M_7$, $M_8$ and $M_9$ that is no school even for primary classes, poverty, no proper road facilities and living condition of them are questionably poor. Now we consider the second experts opinion over the same set of attributes. The directed graph given by the expert relating the attributes {$M_1$, $M_2$,…, $M_{10}$} and {$P_1$, $P_2$, $P_3$,…,$P_7$} is as follows:

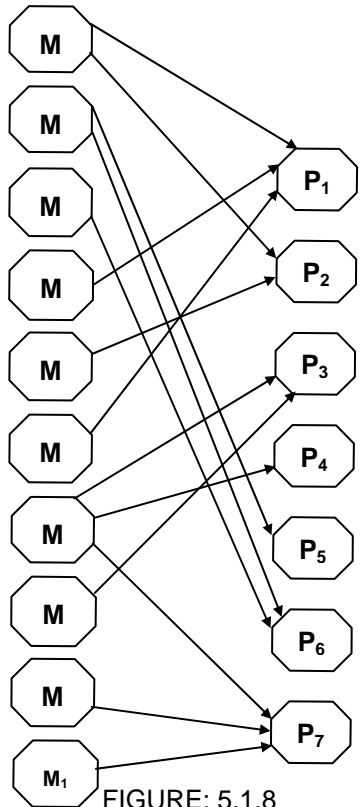

FIGURE: 5.1.8

The related $10 \times 7$ connection matrix $G_2$ given below

$$
\begin{array}{c}
\phantom{M_{10}} \quad P_1 \; P_2 \; P_3 \; P_4 \; P_5 \; P_6 \; P_7 \\
\begin{array}{c}
M_1 \\ M_2 \\ M_3 \\ M_4 \\ M_5 \\ M_6 \\ M_7 \\ M_8 \\ M_9 \\ M_{10}
\end{array}
\begin{bmatrix}
1 & 1 & 0 & 0 & 0 & 0 & 0 \\
0 & 0 & 0 & 0 & 1 & 1 & 0 \\
0 & 0 & 0 & 0 & 0 & 1 & 0 \\
1 & 0 & 0 & 0 & 0 & 0 & 0 \\
0 & 1 & 0 & 0 & 0 & 0 & 0 \\
1 & 0 & 0 & 0 & 0 & 0 & 0 \\
0 & 0 & 1 & 1 & 0 & 0 & 1 \\
0 & 0 & 1 & 0 & 0 & 0 & 0 \\
0 & 0 & 0 & 0 & 0 & 0 & 1 \\
0 & 0 & 0 & 0 & 0 & 0 & 1
\end{bmatrix}
\end{array}
$$

Now using the connection matrix $G_2$ one can study the effect of the state vector on the dynamical system.

Let
$$X \qquad = \qquad (0\,0\,0\,0\,1\,0\,0\,0\,0\,0)$$

be the state vector which has only the attribute $M_5$ to be in the on state and all other nodes are in the off state.

The effect of X on the dynamical system $G_2$ is given by

$$XG_2 \quad \hookrightarrow \quad (0\,1\,0\,0\,0\,0\,0) \qquad = \qquad \text{Y (say)}$$

$$YG_2^T \quad \hookrightarrow \quad (1\,0\,0\,0\,1\,0\,0\,0\,0\,0) \quad = \quad X_1 \text{ say}$$

$$X_1\,G_2 \quad \hookrightarrow \quad (1\,1\,0\,0\,0\,0\,0) \qquad = \qquad Y_1 \text{ say}$$

$$Y_1\,G_2^T \quad \hookrightarrow \quad (1\,0\,0\,1\,1\,1\,0\,0\,0\,0) \quad = \quad X_2 \text{ say}$$

$$X_2\,G_2 \quad \hookrightarrow \quad (1\,1\,0\,0\,0\,0\,0) \qquad = \quad Y_2\,(= Y_1).$$

Thus the hidden pattern of the dynamical system is a fixed point given by the binary pair {(1 0 0 1 1 1 0 0 0 0), (1 1 0 0 0 0 0)}.

Thus when no proper health center exists in a village we see the nodes $M_{11}$, $M_4$ and $M_6$ come to on state in the domain space and $P_1$ and $P_2$ comes to on state in the range space there by indicating they have no job opportunities, no school for even primary education and government help never reaches them and no education and awareness programs never reach them. So government has not even extended the help of building a good health center for these rural uneducated poor.

Now we consider the state vector $Y = (0\ 0\ 1\ 0\ 0\ 0\ 0)$ only $P_3$ attribute is in the on state and all other state vectors are in the off state.

The effect of $Y$ on the dynamical system $G_2$ is given by

$$YG^T_2 \quad \hookrightarrow \quad (0\ 0\ 0\ 0\ 0\ 0\ 1\ 1\ 0\ 0) \quad = \quad X \text{ say}$$

$$XG_2 \quad \hookrightarrow \quad (0\ 0\ 1\ 1\ 0\ 0\ 1) \quad\quad\quad = \quad Y_1 \text{ say}$$

$$Y_2G^T_2 \quad \hookrightarrow \quad (0\ 0\ 0\ 0\ 0\ 0\ 1\ 1\ 1\ 1) \quad = \quad X_1 \text{ say}$$

$$X_1G_2 \quad \hookrightarrow \quad (0\ 0\ 1\ 1\ 0\ 0\ 1) \quad\quad\quad = \quad Y_2\ (=Y_1).$$

Thus the hidden pattern of the dynamical system is a fixed point. The fixed binary pair is given by $\{(0\ 0\ 0\ 0\ 0\ 0\ 1\ 1\ 1\ 1), (0\ 0\ 1\ 1\ 0\ 0\ 1)\}$.

We see only when the node no responsibility of parents to educate the children in the range space we see the nodes $P_4$ and $P_7$ become on and all other nodes remain off that is thus the on state of this node makes the parents get their girl children get married at a very young age and their living conditions are questionably poor so they migrate hence can not find means to educate their children as they constantly move. Also the nodes $M_7$, $M_8$, $M_9$ and $M_{10}$ come to on state in the domain space there by poverty, no road facility (that is even to send their children for studies to neigbouring school) living conditions is poor and government unconcern to help there agriculture labourer with alternatives come to on state.

Now we proceed on the get the third experts opinion using same attributes $\{M_1, M_2,\ldots, M_{10}\}$ for the domain space and $\{P_1, P_2,\ldots, P_7\}$ for the range space. The directed graph given by the expert is

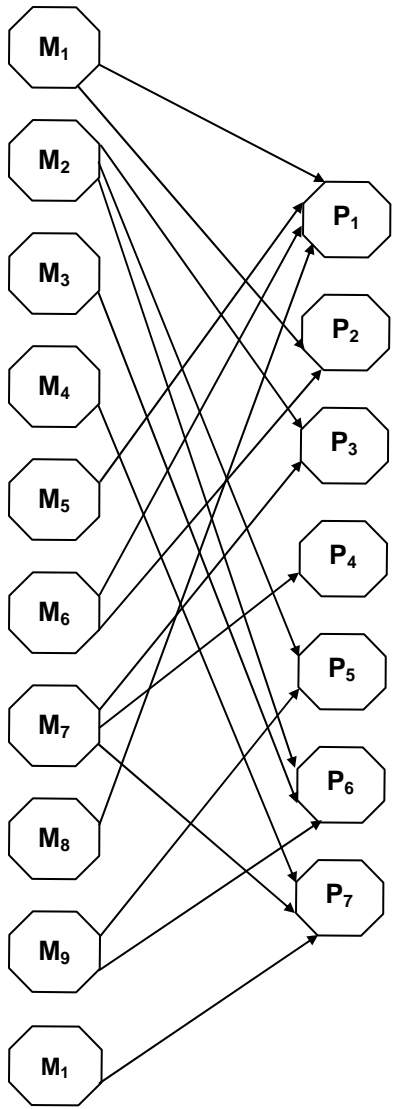

FIGURE: 5.1.9

The related connection $10 \times 7$ matrix $G_3$ is as follows:

$$
\begin{array}{c}
\phantom{M_{10}} \quad P_1 \ P_2 \ P_3 \ P_4 \ P_5 \ P_6 \ P_7 \\
\begin{array}{c}
M_1 \\ M_2 \\ M_3 \\ M_4 \\ M_5 \\ M_6 \\ M_7 \\ M_8 \\ M_9 \\ M_{10}
\end{array}
\begin{bmatrix}
1 & 1 & 0 & 0 & 0 & 0 & 0 \\
0 & 0 & 1 & 0 & 1 & 1 & 0 \\
0 & 0 & 0 & 0 & 0 & 1 & 0 \\
0 & 0 & 0 & 0 & 0 & 0 & 1 \\
1 & 0 & 0 & 0 & 0 & 0 & 0 \\
1 & 1 & 0 & 0 & 0 & 0 & 0 \\
0 & 0 & 1 & 1 & 0 & 0 & 1 \\
1 & 0 & 0 & 0 & 0 & 0 & 0 \\
0 & 0 & 0 & 0 & 1 & 1 & 0 \\
0 & 0 & 0 & 0 & 0 & 0 & 1
\end{bmatrix}
\end{array}
$$

The effect of the state vector X = (0 0 0 0 0 1 0 0 0 0) where only the attribute $M_6$ is in the on state and all other nodes are in the off state. Effect of X on $G_3$ is given by

$$
\begin{aligned}
XG_3 &\hookrightarrow (1\ 1\ 0\ 0\ 0\ 0\ 0) & &= Y \\
YG^T_3 &\hookrightarrow (1\ 0\ 0\ 0\ 1\ 1\ 0\ 1\ 0\ 0) & &= X_1 \\
X_1\,G_3 &\hookrightarrow (1\ 1\ 0\ 0\ 0\ 1) & &= Y_1 \text{ say} \\
Y_1\,G^T_3 &\hookrightarrow (1\ 0\ 0\ 1\ 1\ 1\ 1\ 1\ 0\ 1) & &= X_2 \text{ say} \\
X_2\,G_3 &\hookrightarrow (1\ 1\ 1\ 1\ 0\ 0\ 1) & &= Y_2 \text{ say} \\
Y_2\,G^T_3 &\hookrightarrow (1\ 1\ 0\ 1\ 1\ 1\ 1\ 1\ 0\ 1) & &= X_3 \text{ say} \\
X_3\,G_3 &\hookrightarrow (1\ 1\ 1\ 1\ 1\ 1\ 1) & &= X_4 \text{ (say)} \\
X_4\,G_3 &\hookrightarrow (1\ 1\ 1\ 1\ 1\ 1\ 1) & &= Y_4 \text{ (= } Y_3).
\end{aligned}
$$

Thus the hidden pattern of the dynamical system is a fixed point given by the binary pair {(1 1 1 1 1 1 1 1 1 1), (1 1 1 1 1 1 1)}. When there is no school even for primary classes alone is int the on state, we see all the nodes both in the domain space and the range space come to on state there by indicating the strong impact on the absence of even a primary school.

Now suppose we consider the state vector $Y_1$ = (0 0 1 0 0 0 0) i.e., only the attribute $P_3$ is in the on state and all other attributes in the off state the effect of $Y_1$ on the dynamical system is given by

$Y_1 G^T_3 \hookrightarrow (0\ 1\ 0\ 0\ 0\ 0\ 1\ 0\ 0\ 0) = X_1$ say

$Y_1\ G_3 \hookrightarrow (0\ 0\ 1\ 1\ 1\ 1\ 1) = Y_2$ (say)

$Y_2\ G^T_3 \hookrightarrow (0\ 1\ 1\ 1\ 0\ 0\ 1\ 0\ 1\ 1) = X_2$ say

$X_2\ G_3 \hookrightarrow (0\ 0\ 1\ 1\ 1\ 1\ 1) = Y_3\ (= Y_2$ say)

$Y_3\ G^T_3 \hookrightarrow (0\ 1\ 1\ 1\ 0\ 0\ 1\ 0\ 1\ 1) = X_3\ (= X_2$ say).

Thus the hidden pattern of the dynamical system is a fixed point given by the binary pair {(0 0 1 1 1 1 1), (0 1 1 1 0 0 1 0 1 1)}.

Only when the node parents not responsible enough to educate their children alone is in the on state in the range space and all other nodes are in the off state we see that $P_4$, $P_5$, $P_6$ and $P_7$ come to on state in the range space and $M_2$, $M_3$, $M_4$, $M_7$, $M_9$ and $M_{10}$ come to on state their by indicating the significance of education.

Now we obtain the combined FRM using the matrices $G_1$, $G_2$ and $G_3$; $G = G_1 + G_2 + G_3$ where G is the connection matrix related with the combined FRM G is also a $10 \times 6$ matrix having {$M_1$, $M_2$,…, $M_{10}$} and {$P_1$, $P_2$,…, $P_7$} as the attributes

|        | $P_1$ | $P_2$ | $P_3$ | $P_4$ | $P_5$ | $P_6$ | $P_7$ |
|--------|-------|-------|-------|-------|-------|-------|-------|
| $M_1$  | 2 | 2 | 0 | 0 | 0 | 0 | 1 |
| $M_2$  | 0 | 0 | 1 | 0 | 2 | 2 | 0 |
| $M_3$  | 0 | 1 | 0 | 0 | 0 | 3 | 0 |
| $M_4$  | 1 | 0 | 0 | 0 | 0 | 0 | 2 |
| $M_5$  | 1 | 2 | 0 | 0 | 0 | 0 | 0 |
| $M_6$  | 3 | 1 | 1 | 0 | 0 | 0 | 0 |
| $M_7$  | 0 | 0 | 1 | 1 | 0 | 0 | 2 |
| $M_8$  | 2 | 0 | 2 | 0 | 0 | 0 | 0 |
| $M_9$  | 0 | 0 | 0 | 1 | 1 | 1 | 1 |
| $M_{10}$ | 0 | 0 | 0 | 0 | 0 | 0 | 3 |

Consider the state vector X = (0 0 0 0 0 0 0 1 0 0) only the attributes $M_8$ is in the on state and all vectors are in the off state the effect of X on the dynamical system G is given by

$XG \hookrightarrow (1\ 0\ 1\ 0\ 0\ 0\ 0) = Y$

$$YG^T \quad \hookrightarrow \quad (1\ 1\ 1\ 1\ 1\ 1\ 0\ 1\ 0\ 0) \quad = \quad X_1\ (\text{say})$$

$$X_1\ G \quad \hookrightarrow \quad (1\ 1\ 1\ 0\ 1\ 1\ 1) \quad = \quad Y_1\ \text{say}$$

$$Y_1\ G^T \quad \hookrightarrow \quad (1\ 1\ 1\ 1\ 1\ 1\ 1\ 1\ 1\ 1) \quad = \quad X_2\ \text{say}$$

$$X_2\ G \quad \hookrightarrow \quad (1\ 1\ 1\ 1\ 1\ 1\ 1) \quad = \quad Y_2$$

$$Y_2\ G^T \quad \hookrightarrow \quad (1\ 1\ 1\ 1\ 1\ 1\ 1\ 1\ 1\ 1) \quad = \quad X_3\ (= X_2).$$

Thus we see the hidden pattern of the dynamical system is a fixed point given by the binary pair {(1 1 1 1 1 1 1 1 1 1), (1 1 1 1 1 1 1)}. We see when the villiage has np proper road or bus route we see all the nodes both from the domain space come to on state there by making the essence of the road or bus facility from one place too another. Now let us consider the on state the node $P_7$ and all other nodes are in the off state, the effect of given by

$$YG^T \quad \hookrightarrow \quad (1\ 0\ 0\ 1\ 0\ 0\ 1\ 0\ 11) \quad = \quad X\ \text{say}$$

$$XG^T \quad \hookrightarrow \quad (1\ 1\ 1\ 1\ 1\ 1\ 1) \quad = \quad Y_1\ \text{say}$$

$$Y_1\ G^T \quad \hookrightarrow \quad (1\ 1\ 1\ 1\ 1\ 1\ 1\ 1\ 1\ 1) \quad = \quad X_1\ \text{say}.$$

Thus we see the binary pair {(1 1 1 1 1 1 1), (1 1 1 1 1 1 1 1 1 1)} is a fixed point. When $P_7$ node alone is in the on state that is very questionable living conditions so they migrate in the range space and all other nodes are in the off state we see all the nodes both in the range and domain space become on in the CFCM.

Using the C-program given in the appendix 8 for the combined FRM one can easily solve the problem.

## 5.2 Defintion of Linked FRMs and its applications

Here we introduce the new notion of linked FRMs. We may have a problems to which we may be able to associate attributes which may be divided into spaces say D, R and S we would be in a position to associate say D and R and R and S and find it impossible to draw the directed graph relating D and S but we know they are indirectly related in such cases, the method of linked FRM provides us a means to relate these two using the connection matrices associated with D and R and R and S.

**DEFINITION 5.2.1:** *Let P be some problem under study. Suppose the attributes associated with P be divided into 3 classes say D, R and S where D has n atiributes ($D_1$, …, $D_n$) R has m attributes say ($R_1$, $R_2$,…, $R_m$) and S has t attributes say ($S_1$, $S_2$,…, $S_t$). We are in a position to get the directed graph relating the spaces D and R and suppose the connection matrix associated with it is a n $\times$ m matrix say A. Suppose we are in a position to associate the attributes in the classes R and S through a directed graph and the related connection matrix be B which is a m $\times$ t matrix. But we are aware of the fact that the attributes in the classes D and S are related but we are not in a position to exibit them in the form of the directed graph so it may not possible for us to get the related matrix.*

*Now we can relate the attributes in D and S by calculating the product of the matrices A and B. A.B = C gives the n $\times$ t matrix which links D and S, C will serve as the relational matrix. We call such FRMs as linked FRMs we can link any finite number of FRMs. For suppose $A_1$, $A_2$,…, $A_t$ denote t classes of attributes. We can at each stage link them by the above process. The relation when the linked FRM is given we call the directed graph representation of the directed graph as the linked directed graph; for we have linked to get the directed graph.*

We now illustrate these linked FRM by the problem of migrant labourers affected with HIV/AIDS.

The linked FRM inter relate the indirect relation between two sets of attribute which cannot be related by directed graphs. Now we consider the 3 sets of attributes related with the migrant labourers given in chapter III under the three heads

| | | |
|---|---|---|
| A | – | Causes for migrant labourers vulnerability to HIV/AIDS. |
| F | – | Factors forcing for migration. |
| G | – | Role of government. |

Now using the two FRMs relating causes for migrant labourers vulnerability to HIV/AIDS and its relation to factors forcing people for migration cannot be found directly so using the FRM relating the two sets of attributes Factors forcing people for migration and the Role of government. Now using these two FRMs and the FRM relating the causes for migrant labourers vulnerability to HIV/AIDS and its relation with the role of government one can find the relation between A and F. The FRMs directed graph relating attributes A = {$A_1$, $A_2$, …, $A_6$} and G = {$G_1$, $G_2$, $G_3$ , $G_4$ , $G_5$}

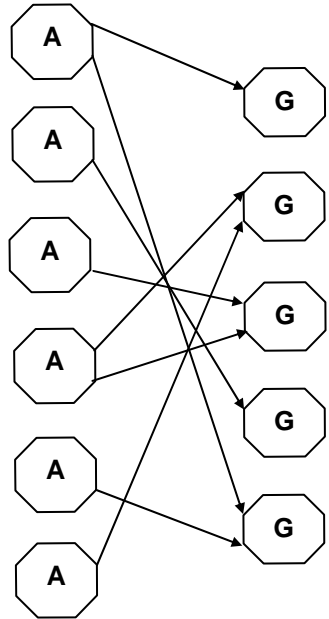

FIGURE: 5.2.1

The related $6 \times 5$ connection matrix of the FRM denoted $L_1$

$$
\begin{array}{c c}
 & \begin{array}{c c c c c} G_1 & G_2 & G_3 & G_4 & G_5 \end{array} \\
\begin{array}{c} A_1 \\ A_2 \\ A_3 \\ A_4 \\ A_5 \\ A_6 \end{array} &
\left[ \begin{array}{c c c c c}
1 & 0 & 0 & 0 & 1 \\
0 & 0 & 0 & 1 & 0 \\
0 & 0 & 1 & 0 & 0 \\
0 & 1 & 1 & 0 & 0 \\
0 & 0 & 0 & 0 & 1 \\
0 & 1 & 0 & 0 & 0
\end{array} \right]
\end{array}
$$

Suppose we have the state vector $X = (0\ 0\ 0\ 1\ 0\ 0)$ ie the attribute $A_4$ alone is in the on state and all other nodes are in the off state. The effect of $X$ on the dynamical system $L_1$ is given by

$$XL_1 \quad \hookrightarrow \quad (0\ 0\ 1\ 0\ 0) \qquad = \quad Y \text{ say}$$

$$YL_1^T \quad \hookrightarrow \quad (0\ 0\ 1\ 1\ 0\ 0) \qquad = \quad X_1 \text{ say}$$

$$X_1 L_1 \quad \hookrightarrow \quad (0\ 1\ 1\ 0\ 0) \qquad = \quad Y_1 \text{ say}$$

$Y_1 L_1 \hookrightarrow (0\ 0\ 1\ 1\ 0\ 1) = X_2$ say

$X_2 L_1 \hookrightarrow (0\ 1\ 1\ 0\ 0) = Y_3 = (Y_2)$.

Thus the hidden pattern of the dynamical system is a fixed point given by the binary pair {(0 1 1 0 0), (0 0 1 1 0 1)} which implies bad company and bad addictive habits, when alone is in the on state we see the nodes $A_3$ and $A_6$ come to on state and in the range space the nodes $G_2$ and $G_3$ come to on state. Thus we see the nodes no social responsibility and cheap availability of CSWs come to on state. So that even persons with bad company and bad habits can save themselves once they are aware about it.

Next we consider the effect of the state vector in the domain space where only the node $G_1$ is in the on state and all other vectors are in the off state. The effect of $Y = (1\ 0\ 0\ 0\ 0)$ on the dynamical system $L_1$ is given by

$YL^T_1 \hookrightarrow (1\ 0\ 0\ 0\ 0\ 0) = X$

$XL_1 \hookrightarrow (1\ 0\ 0\ 0\ 1) = Y_1$ say

$Y_1 L^T_1 \hookrightarrow (1\ 0\ 0\ 0\ 1\ 0) = X_1$ say

$X_1 L_1 \hookrightarrow (1\ 1\ 0\ 0\ 1) = Y_2$ say

$Y_2 L^T_1 \hookrightarrow (1\ 0\ 0\ 1\ 1\ 1) = X_2$ say

$X_2 L_1 \hookrightarrow (1\ 1\ 1\ 0\ 1) = Y_3$ (say)

$Y_3 L^T_1 \hookrightarrow (1\ 0\ 1\ 1\ 1\ 1) = X_3$ say

$X_3 L_1 \hookrightarrow (1\ 1\ 1\ 0\ 1) = Y_4 = (Y_3)$.

Thus the hidden pattern of the dynamical system is a fixed point given by the binary pair {(1 1 1 0 1), (1 0 1 1 1 1)}. We see when the node alternative job if the harvest fails there by stopping migration is in the on state we see in the resultant statevectors $G_2$, $G_3$ and $G_5$ come to on state there by stating awareness clubs in rural areas about HIV/AIDS in villages come on together government had failes to take precautionary steps to help them from past experience. Further 11 nodes $A_1$, $A_3$, $A_4$, $A_5$, $A_6$ come to on states forcing them as they have no education so the profession they choose make them seek bad companay and this helped by the cheap availability of CSWs.

The directed graph of the FRM relating the factors forcing migration and the role of government is given by the following directed graph.

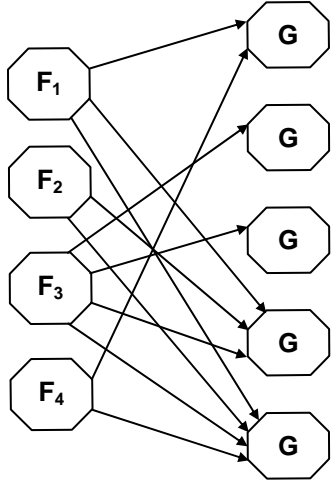

FIGURE: 5.2.2

According to expert, when the mobilization labour contractors with no awareness about HIV/AIDS as the awareness program in rural areas especially when they are uneducated and poor never reaches them, also if there were some health centers or hospitals in the rural areas the people in those places would have become aware of HIV/AIDS but this also does not what takes place is the mobilization of labour contractors; who in the city become easy victims of HIV/AIDS with the catalyst of liquor and cheap availability of CSWs.

So this expert has justified his directed graph in this way. As we have no right to modify these opinion we give it as it is even if it is meaningless to some the related connection matrix be denoted by $L_2$

$$L_2 = \begin{array}{c} \\ F_1 \\ F_2 \\ F_3 \\ F_4 \end{array} \begin{array}{ccccc} G_1 & G_2 & G_3 & G_4 & G_5 \\ \begin{bmatrix} 1 & 0 & 0 & 1 & 1 \\ 0 & 0 & 0 & 1 & 1 \\ 0 & 1 & 1 & 1 & 1 \\ 1 & 0 & 0 & 0 & 1 \end{bmatrix} \end{array}$$

Let us consider the vector X = (1 0 0 0) ie only the attribute $F_1$ is in the on state and all other attributes are in the off state. The effect of X on the dynamical system $L_2$

$$XL_2 \hookrightarrow (1\ 0\ 0\ 1\ 1) = Y \text{ say}$$

$$YL^T_2 \hookrightarrow (1\ 1\ 1\ 1) = X_1 \text{ say}$$

$$X_1 L_2 \hookrightarrow (1\ 1\ 1\ 1\ 1) = Y_1 \text{ say}$$

$$Y_1 L^T_2 \hookrightarrow (1\ 1\ 1\ 1) = X_2\ (=X_1).$$

Thus the hidden pattern is a fixed point given by the binary pair {(1 1 1 1), (1 1 1 1 1)} we see if lack of labour opportunity alone is in the home town all other nodes both from the domain and range space become on their by forcing the government to take better steps to find means to stop unemployment by alternative means which show vital role played by that node on all other nodes both from the domain and range space.

Now suppose we consider the state vector in which the attribute $G_5$ alone is in the on state and all other attributes are in the off state the effect of the state vector Y = (0 0 0 0 1) on the dynamical system $L_2$

$$YL^T_2 \hookrightarrow (1\ 1\ 1\ 1) = X \text{ (say)}$$

$$XL_2 \hookrightarrow (1\ 1\ 1\ 1\ 1) = Y_1 \text{ (say)}$$

$$Y_1 L^T_2 \hookrightarrow (1\ 1\ 1\ 1) = X_1 = X.$$

Thus the hidden pattern is a fixed point of the dynamical system given by the binary {(1 1 1 1), (1 1 1 1 1)}. Thus when the attributes no foresight for the government and no precautionary actions taken from the past occurrences in on state all attributes in both the domain space and range space come to on state there by showing the importance of that attribute.

Using the C program in the appendix 7 of the book one can work out for the conclusions.

Now using the FRM matrices $L_1$ and $L_2$ we link the relation between the class of attributes {$A_1$, $A_2$, …, $A_6$} and {$F_1$, $F_2$, $F_3$, $F_4$}. Let L denote the link relating these two attributes given by $L_1$ × $L^T_2$ ~ L this links the attribute {$A_1$, …, $A_6$} with {$F_1$, $F_2$, $F_3$, $F_4$}.

$$L = \begin{bmatrix} 1 & 1 & 1 & 1 \\ 1 & 1 & 1 & 0 \\ 0 & 0 & 1 & 0 \\ 0 & 0 & 1 & 0 \\ 1 & 1 & 1 & 1 \\ 0 & 0 & 1 & 0 \end{bmatrix}$$

Here '~' represents the entries $a_{ij}$ in $L_1 L^{-2}_2$ are replaced by 0 if $a_{ij} \leq 0$; $a_{ij} = 1$ if $a_{ij}$ is greater than or equal to 1.

Now we can using this linked matrix draw the relational map connecting $\{A_1, A_2, \ldots, A_6\}$ and $\{F_1\ F_2\ F_3\ F_4\}$ which we will be calling as the linked relational graph, unlike in a linked relational graph we will have the relations between only two nodes to exist only implicitly and not explicitly.

The linked relational graph relating $\{A_1, \ldots, A_6\}$ and $\{F_1\ F_2\ F_3\ F_4\}$ obtained from the matrix L is as follows.

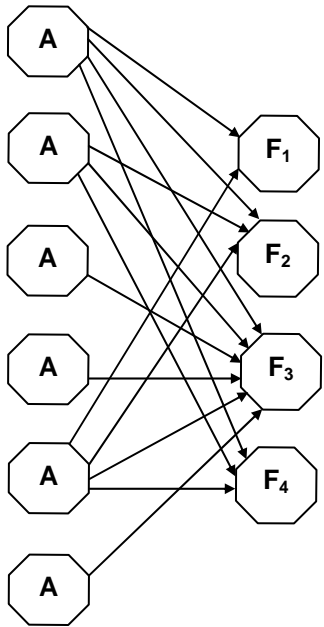

FIGURE: 5.2.3

Suppose we consider the state vector X = {0 1 0 0 0 0} ie the attribute $A_2$ alone is in the on state and all other nodes are in the off state, effect of X on the linked dynamical system L is given by

| XL | $\hookrightarrow$ | (1 1 1 0) | = | Y (say) |
|----|----|----|----|----|
| $YL^T$ | $\hookrightarrow$ | (1 1 1 1 1 1 1) | = | $X_1$ (say) |
| $X_1 L$ | $\hookrightarrow$ | (1 1 1 1) | = | $Y_1$ |
| $Y_1 L^T$ | $\hookrightarrow$ | (1 1 1 1 1 1 1) | = | $X_2 (=X_1)$ say. |

Thus the hidden pattern is a fixed point of the dynamical system given by the binary pair {(1 1 1 1 1 1 1), (1 1 1 1)}. This resultant vector shows how powerful is the influence of the node $A_2$ on the whole dynamical systems.

Now consider the state vector X = (0 0 1 0 0 0), ie only the attribute $A_3$ is in the on state and all other attributes are in the off state. The effect of X on the dynamical system L is given by

| XL | $\hookrightarrow$ | (0 0 1 0) | = | Y say |
|----|----|----|----|----|
| $YL^T$ | $\hookrightarrow$ | (1 1 1 1 1 1) | = | $X_1$ say |
| XL | $\hookrightarrow$ | (1 1 1 1) | = | $Y_1$say |
| $Y_1 L^T$ | $\hookrightarrow$ | (1 1 1 1 1 1) | = | $X_2 (=X_1)$. |

Thus the fixed point is a binary pair {(1 1 1 1 1 1), (1 1 1 1)}. Here also the node $A_3$ no social responsibility has so much of influence on nodes of both the domain and rage space, there by insisiting if a migrant labourer is socially irresponsible all evil take over ultimately leaving him as an HIV/AIDS patient.

Next we proceed on to study the model with some different types of attributes associated the HIV/AIDS affected migrant labourers.

Thus we can link seemingly unrelated issues but which are implicitly related into a model. Hence this study is very helpful when one wishes to study the effect of one set of attributes over another set of attributes, which are not explicitly related. Now we consider another three sets of attributes related with the HIV/AIDS migrant labourers under three separate heads.

**Causes of migration and their vulnerability to HIV/AIDS**

$C_1$ - Easy victims of temptation so fall a pray to CSWs

$C_2$ - Already have all bad habits

$C_3$ - No education so have no principle and goals in life

$C_4$ - They are in a new surrounding and all alone

$C_5$ - No work after days work

**Causes of migrant labourers becoming HIV/AIDS infected**

$A_1$ - Visiting CSWs

$A_2$ - Slaves of bad habits like smoke, alchol and CSWs

$A_3$ - Unaware of HIV/.AIDS

$A_4$ - No education so have not learnt about how HIV/AIDS spreads as the awareness program is rural area is very poor

$A_5$ - Away from family

$A_6$ - No friends / relatives to observe their activities

$A_7$ - No higher motivation very free after work so, only recreation is CSWs.

**Factors acting as Catalyst in the spread of HIV/AIDS among migrant labourers in he city**

$T_1$ - Cheap availability of CSWs in all places

$T_2$ - Free sale of liquor

$T_3$ - Free from being observed

$T_4$ - Bad company / bad friends

$T_5$ - Easy money no knowledge of saving / beter investment

$T_6$ - Men in rural areas are rarely aware of the mode of spread of HIV/AIDS.

Now using the experts opinion we give the directed graph of the FRM relating the attributes $\{A_1, \ldots, A_7\}$ and $\{T_1, \ldots, T_6\}$ in the following page.

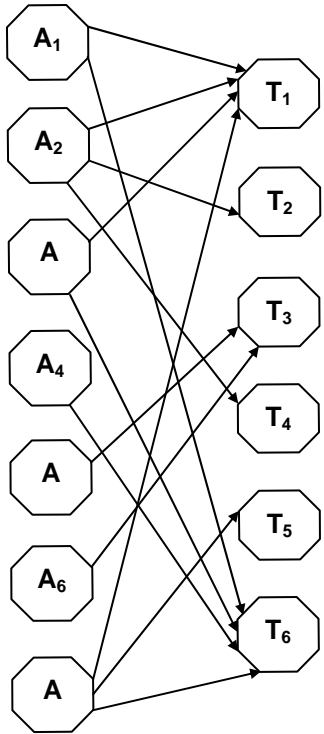

FIGURE: 5.2.4

The connection matrix related with the FRM be given by $L_1$

$$
L_1 = \begin{array}{c} \\ A_1 \\ A_2 \\ A_3 \\ A_4 \\ A_5 \\ A_6 \\ A_7 \end{array}
\begin{array}{c} T_1\ T_2\ T_3\ T_4\ T_5\ T_6 \\ \begin{bmatrix} 1 & 0 & 0 & 0 & 0 & 1 \\ 1 & 1 & 0 & 1 & 0 & 0 \\ 1 & 0 & 0 & 0 & 0 & 1 \\ 0 & 0 & 0 & 0 & 0 & 1 \\ 0 & 0 & 1 & 0 & 0 & 0 \\ 0 & 0 & 1 & 0 & 0 & 0 \\ 1 & 0 & 0 & 0 & 1 & 1 \end{bmatrix} \end{array}
$$

We now study the effect of the state vector X = (0 0 0 0 1 0 0) where only the node $A_5$ is in the on state all nodes are in the off state.

$$XL_1 \quad \hookrightarrow \quad (0\ 0\ 1\ 0\ 0\ 0) \quad = \quad Y \text{ say}$$

$$YL^T_1 \quad \hookrightarrow \quad (0\ 0\ 0\ 0\ 1\ 1\ 0) \quad = \quad X_1 \text{ say}$$

$$X_1 L_1 \quad \hookrightarrow \quad (0\ 0\ 1\ 0\ 0\ 0) \quad = \quad Y_1 \ (=Y).$$

Thus the hidden pattern of the dynamical system is a fixed point given by the binary pair {(0 0 1 0 0 0), (0 0 0 0 11 0)}. Thus when the node away from the family is in the on state we see that the resultant pair $A_6$ comes to on state that is they are free from the fear of being observed by friends and they are free to act in any way. All other nodes remain in the off state.

Now we consider the effect of the state vector Y = (0 0 0 1 0 0) ie only the node $T_4$ is in the on state and all other vectors are in the off state; the effect of $T_4$ on $L_1$ is given by

$$YL^T_1 \quad \hookrightarrow \quad (0\ 1\ 0\ 0\ 0\ 0\ 0) \quad = \quad X \text{ say}$$

$$XL_1 \quad \hookrightarrow \quad (1\ 1\ 0\ 1\ 0\ 0) \quad = \quad Y_1 \text{ say}$$

$$Y_1L^T_1 \quad \hookrightarrow \quad (1\ 1\ 1\ 0\ 0\ 0\ 1) \quad = \quad X_1 \text{ say}$$

$$Y_1 L_1 \quad \hookrightarrow \quad (1\ 1\ 0\ 1\ 1\ 1) \quad = \quad Y_2 \text{ say}$$

$$Y_2L^T_1 \quad \hookrightarrow \quad (1\ 1\ 1\ 1\ 0\ 0\ 1) \quad = \quad X_2 \text{ say}$$

$$X_2L_1 \quad \hookrightarrow \quad (1\ 1\ 0\ 1\ 1\ 1) \quad = \quad Y_3 = (Y_2).$$

Thus the hidden pattern of the dynamical system is a fixed point given by the binary pair {(1 1 0 1 1 1), (1 1 1 1 0 0 1)}. When the node $T_4$ is in the on state that is bad company and bad friends we see in the resultant pair all nodes become on both in the range and the domain space expert $T_3$ and $A_5$ and $A_6$ there by indicating free from being observed has no effect when they are under the spell of bad company and friends. Also they freely visit CSWs, they are slaves of all bad and addictive habits, unaware of HIV/AIDS, no education so do not know about HIV/AIDS and as they have no higher aspiration only CSWs are the main source of past time after work hours.

Now we obtain the FRM relating the sets of attributes $\{C_1, C_2, \ldots, C_5\}$ and $\{A_1, A_2, \ldots, A_7\}$. The directed graph of the FRM is given below.

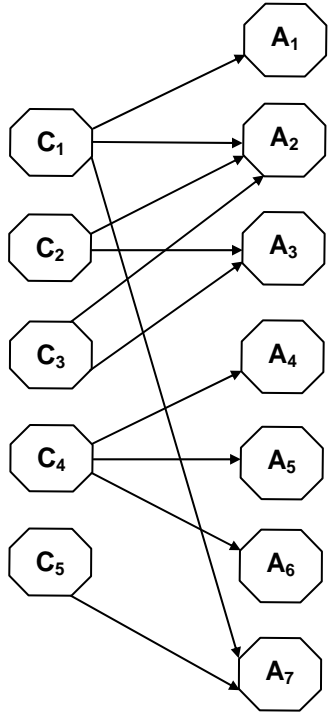

FIGURE: 5.2.5

The connection matrix related with the directed graph be denoted $L_2$

$$L_2 = \begin{array}{c} \\ C_1 \\ C_2 \\ C_3 \\ C_4 \\ C_5 \end{array} \begin{array}{c} \begin{array}{ccccccc} A_1 & A_2 & A_3 & A_4 & A_5 & A_6 & A_7 \end{array} \\ \left[ \begin{array}{ccccccc} 1 & 1 & 0 & 0 & 0 & 0 & 1 \\ 0 & 1 & 1 & 0 & 0 & 0 & 0 \\ 0 & 1 & 1 & 0 & 0 & 0 & 0 \\ 0 & 0 & 0 & 1 & 1 & 1 & 0 \\ 0 & 0 & 0 & 0 & 0 & 0 & 1 \end{array} \right] \end{array}$$

Let us consider the state vector X = (1 0 0 0 0) only the attribute $C_1$ is in the on state and all other vectors are in the off state the effect of X on the dynamical system is given by

$$XL_2 \quad \hookrightarrow \quad (1\ 1\ 0\ 0\ 0\ 01) \quad = \quad Y$$

$$YL^T_2 \quad \hookrightarrow \quad (1\ 1\ 1\ 0\ 1) \quad = \quad X_1 \text{ say}$$

$$X_1 L_2 \quad \hookrightarrow \quad (1\ 1\ 1\ 0\ 0\ 0\ 1) \quad = \quad Y_1 \text{ say}$$

$$Y_1 L^T_2 \quad \hookrightarrow \quad (1\ 1\ 1\ 0\ 1) \quad = \quad X_2\ (= X_1).$$

Thus the hidden pattern of the FRM is a fixed point given by the binary pair {(1 1 1 0 1), (1 1 1 0 0 0 1)}. We see when the node that they are easy victims of temptation and fall a pray is in the on state we see in the resultant vectors the nodes we see the nodes $C_2$, $C_3$ and $C_5$ come to on state that is they are already slaves of bad habits, they have no education and priciple and no work after their days duty. Also in the range space A1, A2, A3 and A7 come to on state there by indicating they are slaves of bad habits, unaware of HIV/AIDS and visiting CSWs is their only recreation.

Suppose we consider the state vector Y = (0 0 1 0 0 0 0) the attribute $A_3$ to be in the on state the effect of Y on the dynamical system $L_2$ is given by

$$YL^T_2 \quad \hookrightarrow \quad (0\ 1\ 1\ 0\ 0) \quad = \quad X$$

$$XL_2 \quad \hookrightarrow \quad (0\ 1\ 1\ 0\ 0\ 0\ 0) \quad = \quad Y_1 \text{ say}$$

$$Y_1 L^T_2 \quad \hookrightarrow \quad (1\ 1\ 1\ 0\ 0) \quad = \quad X_1 \text{ say}$$

$$X_1 L_2 \quad \hookrightarrow \quad (1\ 1\ 1\ 0\ 0\ 01) \quad = \quad Y_2 \text{ say}$$

$$Y_2 L_2 \quad \hookrightarrow \quad (1\ 1\ 1\ 0\ 1) \quad = \quad X_2 \text{ say}$$

$$Y_2 L_2 \quad \hookrightarrow \quad (1\ 1\ 1\ 0\ 0\ 0\ 1) \quad = \quad Y_3 = Y_2.$$

Thus the hidden pattern of the dynamical system is a fixed point given by the binary pair {(1 1 1 0 0 0 1), (1 1 1 0 1)}. We see when the node unaware of HIV/AIDS is on the resultant is the same as when the node easy victims of temptation so they fall a pray to CSWs is on.

Now we are not in a position to relate the sets of attributes {$C_1$ $C_2$ $C_3$ $C_4$ $C_5$} and {$T_1$ $T_2$ $T_3$ $T_4$ $T_5$ $T_6$} and we inter link them by taking the product of the matrix $L_2 L_1$. Let $L \sim L_2 L_1$ we see L is a $5 \times 6$ matrix given as below.

$$\begin{array}{c} \\ C_1 \\ C_2 \\ C_3 \\ C_4 \\ C_5 \end{array} \begin{array}{cccccc} T_1 & T_2 & T_3 & T_4 & T_5 & T_6 \\ \begin{bmatrix} 1 & 1 & 0 & 1 & 1 & 1 \\ 1 & 1 & 0 & 1 & 0 & 0 \\ 1 & 1 & 0 & 1 & 0 & 0 \\ 0 & 0 & 1 & 0 & 0 & 0 \\ 1 & 0 & 0 & 0 & 1 & 1 \end{bmatrix} \end{array}$$

Thus we get the matrix linking the attributes $\{C_1, C_2, \ldots, C_5\}$ and $\{T_1, T_2, \ldots, T_6\}$. We can now get the relational graph which will be called as the link directed graph.

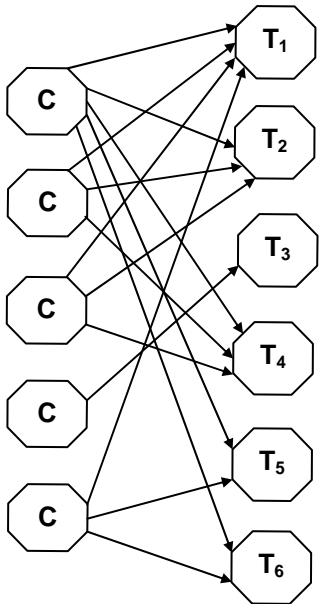

FIGURE: 5.2.6

Now we can study of any of the state vectors in the dynamical system L. Let us consider the state vector $X = (0\ 0\ 1\ 0\ 0\ 0)$ ie only the attribute $T_3$ is in the on state and all other nodes are in the off state. The effect of X on the dynamical system L is given by

$$XL^T \quad \hookrightarrow \quad (0\ 0\ 0\ 1\ 0) \quad = \quad Y$$

$$YL \quad \hookrightarrow \quad (0\ 0\ 1\ 0\ 0\ 0) \quad = \quad X_1 = X$$

Thus the fixed point is the binary pair given by $\{(1\ 0\ 0\ 0\ 1\ 0), (0\ 0\ 1\ 0\ 0\ 0)\}$. When free from being observed is in the on state we see no attributes from that space become on but $A_1$ and $A_5$ become on that is they visit CSWs and no one ever sees them.

Now let us consider the state vector $Y = (1\ 0\ 0\ 0\ 0)$ that is only the attribute $C_1$ is in the on state and all other vectors are in the off state the effect of Y on the dynamical system L

$$YL \quad \hookrightarrow \quad (1\ 1\ 0\ 1\ 1\ 1) \quad = \quad X$$

$$XL^T \quad \hookrightarrow \quad (1\ 1\ 1\ 0\ 1) \quad = \quad Y \text{ say}$$

$$Y_1 L \quad \hookrightarrow \quad (1\ 1\ 0\ 1\ 1\ 1) \quad = \quad X_1 = X.$$

Thus the hidden pattern is a fixed point given by the binary pair $\{(1\ 1\ 0\ 1\ 1\ 1), (1\ 1\ 1\ 0\ 1)\}$. When easy victims of temptation is in the on state we see in the resultant vector $C_2$, $C_3$ and $C_5$ come to on state there by confirming they have all bad habits no education and no work after bad habits no education and no work after days duty and $T_1$, $T_2$, $T_4$, $T_5$ and $T_6$ come to on state they actions are catalyzed by cheap availability of CSWs, free sale of liquor, bad company and bad habits, easy money and no habits of saving and they are unaware of how HIV/AIDS spreads. Thus we can use any state vector and obtain conclusion based on it.

It is left as an exercise for the reader to form a program in C language to find linked FRMs.

## 5.3 Use of CDBFRM to analyse the problem of HIV/AIDS affected migrants

We now define of combined disjoint block FRM of equal and varied lengths or sizes. We illustrate these only in the context of the HIV/AIDS affected migrant labourers.

**DEFINITION 5.3.1:** *Consider a FRM model with n elements {$D_1$, $D_2$,…, $D_n$} in the domain space and {$R_1$,…, $R_m$} in the range space. Now if we divide m into say t blocks and n into r blocks such that m = st and n = pr then we see that the t blocks are disjoint and also r blocks are disjoint.*

*We get many $p \times r$ matries, which is formed into a $n \times m$ matrix which we call as the combined disjoint block FRM of equal length/sizes.*

*If instead of dividing the n attributes of the domain space and matributes of the range space into equal classes each we can also divide them into unequal classes but continue to be disjoint. This case we call as the combined disjoint block FRM of unequal length/sizes.*

**DEFINITION 5.3.2:** *Consider FRM with $D_1$, $D_2$,…, $D_n$ as the attributes related with the domain space and $R_1$, $R_2$,…, $R_m$ be the attributes related with the range space.*

*Now divide the attributes related with the domain space into equal number of classes with overlap and the range space is also divided into classes with equal number of elements but with overlap. We form the $n \times m$ matrix using these blocks of matrices. This matrix corresponds to the combined overlap block FRM of blocks of equal length with over lap.*

*We could also divide the data into blocks of different length and unequal overlap and sizes.*

Depending on the data here we would illustrate both these models only with the problem of HIV/AIDS affected migrant labourers affected with HIV/AIDS.

Now using the FRM model given in chapter V page 201 to 202 using the attributes of the domain space {$D_1$, …, $D_6$} and that of the range space {$R_1$, $R_2$, …, $R_{10}$} we give the disjoint block decomposition of them and give the related directed graphs and their associated connection matrices.

Let the data be divided into two disjoint classes

$$C_1 = \{(D_1 \ D_2 \ D_3), (R_1 \ R_2 \ R_3 \ R_4 \ R_5)\} \text{ and}$$

$$C_2 = \{(D_4 \ D_5 \ D_6), (R_6, R_7, R, R_9, R_{10})\}.$$

The directed graph for the class $C_1$ as given by the expert is as follows:

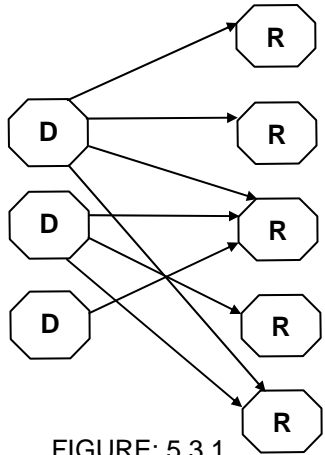

FIGURE: 5.3.1

The related connection matrix is given by

$$
\begin{array}{c}
\begin{array}{ccccc} R_1 & R_2 & R_3 & R_4 & R_5 \end{array} \\
\begin{array}{c} D_1 \\ D_2 \\ D_3 \end{array}
\begin{bmatrix}
1 & 1 & 1 & 0 & 1 \\
0 & 0 & 1 & 1 & 1 \\
0 & 0 & 1 & 0 & 0
\end{bmatrix}
\end{array}
$$

The directed graph relative to $C_2 = \{(D_4\ D_5\ D_6),\ (R_6,\ R_7,\ R,\ R_9,\ R_{10})\}$ as given by the expert is as follows:

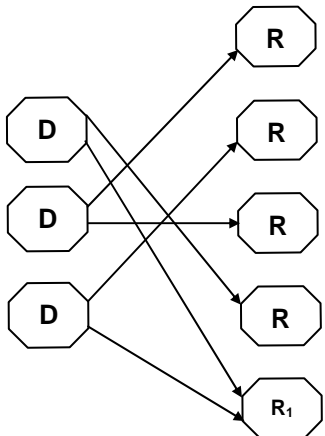

FIGURE: 5.3.2

The related connection matrix is given by

$$
\begin{array}{c}
\phantom{D_4} \begin{array}{ccccc} R_6 & R_7 & R_8 & R_9 & R_{10} \end{array} \\
\begin{array}{c} D_4 \\ D_5 \\ D_6 \end{array}
\begin{bmatrix}
0 & 0 & 0 & 1 & 1 \\
1 & 0 & 1 & 0 & 0 \\
0 & 1 & 0 & 0 & 1
\end{bmatrix}
\end{array}
$$

Using these two connection matrices we now get the matrix associated with the combined disjoint block FRM

$$
\begin{array}{c}
\phantom{D_1} \begin{array}{cccccccccc} R_1 & R_2 & R_3 & R_4 & R_5 & R_6 & R_7 & R_8 & R_9 & R_{10} \end{array} \\
\begin{array}{c} D_1 \\ D_2 \\ D_3 \\ D_4 \\ D_5 \\ D_6 \end{array}
\begin{bmatrix}
1 & 1 & 1 & 0 & 1 & 0 & 0 & 0 & 0 & 0 \\
0 & 0 & 1 & 1 & 1 & 0 & 0 & 0 & 0 & 0 \\
0 & 0 & 1 & 0 & 0 & 0 & 0 & 0 & 0 & 0 \\
0 & 0 & 0 & 0 & 0 & 0 & 0 & 0 & 1 & 1 \\
0 & 0 & 0 & 0 & 0 & 1 & 0 & 1 & 0 & 0 \\
0 & 0 & 0 & 0 & 0 & 0 & 1 & 0 & 0 & 1
\end{bmatrix}
\end{array}
$$

Let C (F) denote the $6 \times 10$ matrix. The effect of the state vector X = (0 0 1 0 0 0) on the dynamical system C(F) only when the $D_3$ node is in the an state and all other nodes are in the off state

$$
\begin{array}{llll}
XC(F) & \hookrightarrow & (0\ 0\ 1\ 0\ 0\ 0\ 0\ 0\ 0\ 0) & = & Y \\
YC(F)^T{}_T & \hookrightarrow & (1\ 1\ 1\ 0\ 0\ 0) & = & X_1 \\
X_1\, C\ (F) & \hookrightarrow & (1\ 1\ 1\ 1\ 1\ 0\ 0\ 0\ 0\ 0) & = & Y_1 \\
Y_1\ (C(F))^T & \hookrightarrow & (1\ 1\ 1\ 0\ 0\ 0) & = & X_2 = X_1
\end{array}
$$

So

$$
\begin{array}{llll}
X_2\ (C(F)) & \hookrightarrow & (1\ 1\ 1\ 1\ 1\ 0\ 0\ 0\ 0\ 0) & = & Y_2 = Y_1.
\end{array}
$$

Thus the hidden pattern is a fixed point given by the binary pair {(1 1 1 0 0 0), (1 1 1 1 1 0 0 0 0 0)}. We see when failure to stop mislead agriculture techniques alone are in the on state, in the resultant vector the nodes $D_1$, $D_2$ come to on state and $R_1$, $R_2$, $R_3$, $R_4$ and $R_5$ come to on state there by indicating awareness program by government and government has not provided alternative for agriculture and in the range space no education, cheap availablitity of CSWs, away from family for weeks, supersistion about sex and no union to channelize and makle them aware of HIV/AIDS comes to on state.

Let us consider the state vector $Y = (0\ 0\ 0\ 0\ 0\ 0\ 1\ 0\ 0\ 0)$ ie only the node $R_7$ is in the on state and all other vectors are in the off state. The effect of Y on C(F) is given by

$$YC(F)^T \hookrightarrow (0\ 0\ 0\ 0\ 0\ 1) = X$$
$$XC(F) \hookrightarrow (0\ 0\ 0\ 0\ 0\ 1\ 0\ 0\ 1) = Y_1 \text{ say}$$
$$Y_1(C(F))^T \hookrightarrow (0\ 0\ 0\ 1\ 0\ 1) = X_1 \text{ say}$$
$$X_1\ (C(F)) \hookrightarrow (0\ 0\ 0\ 0\ 0\ 1\ 0\ 1\ 1) = Y_2 \text{ say}$$
$$Y_2\ (C(F)) \hookrightarrow (0\ 0\ 0\ 1\ 0\ 1) = X_2\ (= X_1).$$

Thus the fixed point is a binary pair given by {(0 0 0 1 0 1), (0 0 0 0 0 1 0 1 1)}. When only the node $R_7$, no job in the native place is in the on state we see in the resultant vector $R_9$ and $R_{10}$ on the range space become on, infertility of land so land contractors take advantages they also do not get proper medical counseling by government. Also in the domain space $D_4$ and $D_6$ come to on state there by indicating the short sightedness of government regarding the problems and government do not help these agricultural labourers.

Now let us consider the model given in page 204-205 with the domain space {$M_1\ M_2,\ \ldots,\ M_9$} which the expert wishes to combine the two attributes $M_9$ and $M_{10}$ into a single attribute $M_9$. The attributes of the range space are taken as {$G_1\ G_2,\ \ldots,\ G_6$} where $G_1$ and $G_2$ are combined together to give the attribute $G_1$ so $G_3$ is renamed as $G_2$, $G_4$ as $G_3$, $G_5$ as $G_4$, $G_6$ as $G_5$ and $G_7$ as $G_6$.

Now using the disjoint block decomposition we divide these two sets of attributes into 3 classes

$$C_1 = \{(M_1,\ M_2\ M_3),\ (G_1\ G_2)\},$$
$$C_2 = \{(M_4\ M_5\ M_6),\ (G_3\ G_4)\} \text{ and }$$
$$C_3 = \{(M_7\ M_8\ M_9),\ (G_5\ G_6)\}$$

Now using the experts opinion we give the related directed connection matrices.

The directed graph for the class $C_1$ given by the expert is as follows:

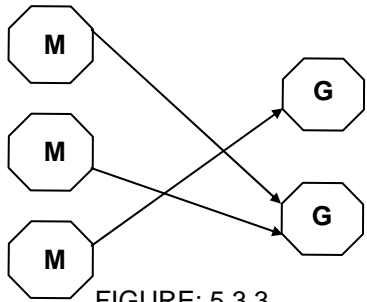

FIGURE: 5.3.3

The associated relational matrix is

$$\begin{array}{c} \\ M_1 \\ M_2 \\ M_3 \end{array} \begin{array}{cc} G_1 & G_2 \\ \begin{bmatrix} 0 & 1 \\ 0 & 1 \\ 1 & 0 \end{bmatrix} \end{array}$$

The directed graph associated given by the expert for the class $C_2 = \{(M_4\ M_5\ M_6)\ (G_3\ G_4)\}$.

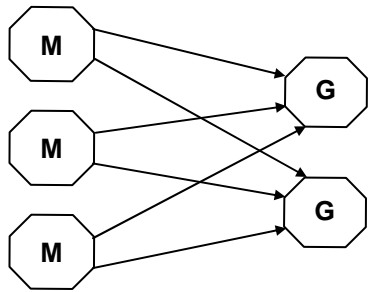

FIGURE: 5.3.4

The associated relational matrix.

$$\begin{array}{c} \\ M_4 \\ M_5 \\ M_6 \end{array} \begin{array}{cc} G_3 & G_4 \\ \begin{bmatrix} 1 & 1 \\ 1 & 1 \\ 1 & 1 \end{bmatrix} \end{array}$$

The directed graph related with $C_3 = \{(M_7 \, M \, M_9), (G_5 \, G_6)\}$.

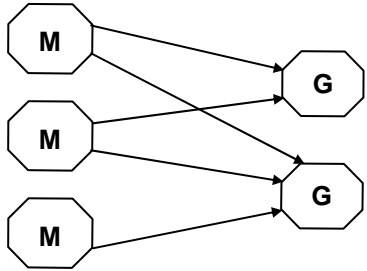

FIGURE: 5.3.5

The associated connection matrix

$$
\begin{array}{c}
\phantom{M_7} \quad G_5 \; G_6 \\
\begin{array}{c} M_7 \\ M_8 \\ M_9 \end{array}
\begin{bmatrix} 1 & 1 \\ 1 & 1 \\ 0 & 1 \end{bmatrix}
\end{array}
$$

The connection matrix associated with the combined disjoint block FRM is a $9 \times 6$ matrix and is denoted by C(G).

$$
\begin{array}{c}
\phantom{M_1} \quad G_1 \; G_2 \; G_3 \; G_4 \; G_5 \; G_6 \\
\begin{array}{c} M_1 \\ M_2 \\ M_3 \\ M_4 \\ M_5 \\ M_6 \\ M_7 \\ M \\ M_9 \end{array}
\begin{bmatrix}
0 & 1 & 0 & 0 & 0 & 0 \\
0 & 1 & 0 & 0 & 0 & 0 \\
1 & 0 & 0 & 0 & 0 & 0 \\
0 & 0 & 1 & 1 & 0 & 0 \\
0 & 0 & 1 & 1 & 0 & 0 \\
0 & 0 & 1 & 1 & 0 & 0 \\
0 & 0 & 0 & 0 & 1 & 1 \\
0 & 0 & 0 & 0 & 1 & 1 \\
0 & 0 & 0 & 0 & 1 & 1
\end{bmatrix}
\end{array}
$$

Let us consider the effect of the state vector $X = (0\ 0\ 0\ 1\ 0\ 0\ 0\ 0\ 1)$ of the dynamical system C(G) where the attributes $M_4$ and $M_9$ are in the on state and all other vectors are in the off state.

| XC(G) | ↪ | (0 0 110 1) | = | Y say |
|---|---|---|---|---|
| Y(C(G))$^T$ | ↪ | (0 0 0 1 1 1 1 1 1) | = | $X_1$ say |
| $X_1$ C(G) | ↪ | (0 0 0 1 1 1 1) | = | $Y_1$ say |
| $Y_1$ C(G) | ↪ | (0 0 0 1 1 1 1 1 1) | = | $X_2 (= X_1)$. |

Thus the fixed point is a binary pair {(0 0 0 1 1 1 1 1 1), (0 0 1 1 1 1)}. When $M_4$ and $M_9$ are in the on state we see in the resultant vector $M_1$, $M_2$ and $M_3$ are in off state but $M_5$, $M_6$, $M_7$ and $M_8$ come to on state in the domain space and $G_3$, $G_4$, $G_5$ and $G_6$ become on and $G_1$ and $G_2$ remain in the off state.

The combined disjoint block FRM would be more powerful and useful only when the attributes can be still split into disjoint disassociated classes other wise it would not be of much use as over lap relation would not find their role.

Now we give the model when the blocks are of unequal size. Now we consider the model with the sets of attributes given by $\{A_1, A_2,…, A_6\}$ and $\{G_1, G_2, G_3, G_4, G_5\}$. Let us divide it into classes

$C_1 = \{A_1 A_2), (G_1, G_4)\}$ and $C_2 = \{(A_3 A_4 A_5 A_6) (G_2 G_3 G_5)\}$.

The related directed graph is given by the expert is as follows:

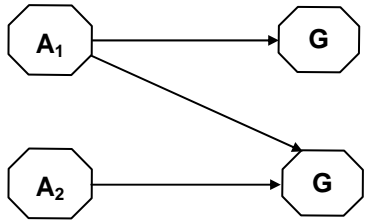

FIGURE: 5.3.6

The connection matrix associated with it

$$\begin{array}{c} \\ A_1 \\ A_2 \end{array} \begin{array}{cc} G_1\ G_4 \\ \begin{bmatrix} 1 & 1 \\ 0 & 1 \end{bmatrix} \end{array}$$

The related directed graph of $C_2$ is given by the expert is as follows:

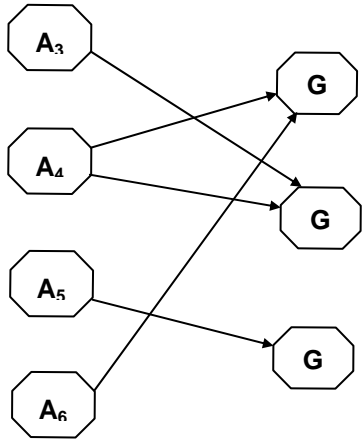

FIGURE: 5.3.7

The related connection matrix

$$\begin{array}{c} & \begin{array}{ccc} G_2 & G_3 & G_5 \end{array} \\ \begin{array}{c} A_3 \\ A_4 \\ A_5 \\ A_6 \end{array} & \left[\begin{array}{ccc} 0 & 1 & 0 \\ 1 & 1 & 0 \\ 0 & 0 & 1 \\ 1 & 0 & 0 \end{array}\right] \end{array}$$

Now we give the associated connection $6 \times 5$ matrix by S

$$\begin{array}{c} & \begin{array}{ccccc} G_1 & G_2 & G_3 & G_4 & G_5 \end{array} \\ \begin{array}{c} A_1 \\ A_2 \\ A_3 \\ A_4 \\ A_5 \\ A_6 \end{array} & \left[\begin{array}{ccccc} 1 & 0 & 0 & 1 & 0 \\ 0 & 0 & 0 & 1 & 0 \\ 0 & 0 & 1 & 0 & 1 \\ 0 & 1 & 1 & 0 & 0 \\ 0 & 0 & 0 & 0 & 1 \\ 0 & 1 & 0 & 0 & 0 \end{array}\right] \end{array}$$

S is the matrix of the combined disjoint block FRM is unequal sizes. Consider the state vector X = (1 0 0 1 0 0) in which the attributes $A_1$ and $A_4$ are in the on state and all other vectors are in the off state the effect of X on the system S is given by

$$\begin{aligned}
XS &\hookrightarrow (1\ 1\ 0\ 1\ 0) &&= Y \\
YS^T &\hookrightarrow (1\ 1\ 0\ 1\ 0\ 1) &&= X_1 \\
X_1\,S &\hookrightarrow (1\ 1\ 1\ 1\ 0) &&= Y_1 \\
Y_1\,S^T &\hookrightarrow (1\ 1\ 1\ 1\ 0\ 1) &&= X_2 \\
X_2\,S &\hookrightarrow (1\ 1\ 1\ 1\ 1) &&= Y_2 \\
Y_2\,S^T &\hookrightarrow (1\ 1\ 1\ 1\ 1\ 1) &&= X_3 \\
X_3\,S &\hookrightarrow (1\ 1\ 1\ 1\ 1) &&= Y_3\ (=Y_2).
\end{aligned}$$

Thus the hidden pattern is a fixed point given by the binary pair $\{(1\ 1\ 1\ 1\ 1\ 1)\,,(1\ 1\ 1\ 1\ 1)\}$.

Now consider the state vector $Y = (0\ 0\ 0\ 1\ 0)$ with only $G_4$ in the on state and all other vectors are in off state. The effect

$$\begin{aligned}
YS^T &\hookrightarrow (1\ 1\ 0\ 0\ 0\ 0) &&= X\ \text{say} \\
XS &\hookrightarrow (1\ 0\ 0\ 1\ 0) &&= Y_1 \\
Y_1 S^T &\hookrightarrow (1\ 1\ 0\ 0\ 0\ 0) &&= X_1\ (= X).
\end{aligned}$$

Thus the fixed binary pair is $\{(1\ 1\ 0\ 0\ 0\ 0), (1\ 0\ 0\ 1\ 0)\}$ which makes the on state of $G_1$ only in the range space and the on state of $A_1$ and $A_2$ in the domain space.

Now we consider the model given in page 198 of this chapter Let the domain space be taken as $\{(D_1, D_2,.., D_8), (R_1\ R_2\ R_3\ R_4\ R_5)\}$ consider the classes $C_1$ and $C_2$ where $C_1 = \{(D_1\ D_2\ D_3\ (R_1\ R_2)\}$ and $C_2\ \{(D_4\ D_5\ D_6\ D_7\ D_8), (R_3\ R_4\ R_5)\}$ Now we obtain the experts opinion and give the directed graph of them.

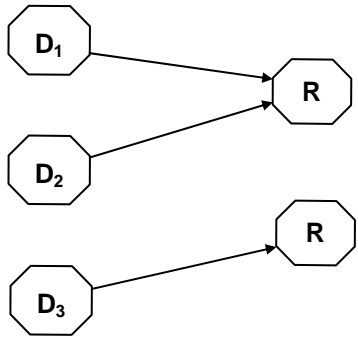

FIGURE: 5.3.8

The related connection matrix is given by

$$
\begin{array}{cc}
 & R_1 \ \ R_2 \\
\begin{array}{c} D_1 \\ D_2 \\ D_3 \end{array} &
\left[\begin{array}{cc}
1 & 0 \\
1 & 0 \\
0 & 1
\end{array}\right]
\end{array}
$$

The directed graph related to $C_2$ as given by the expert is

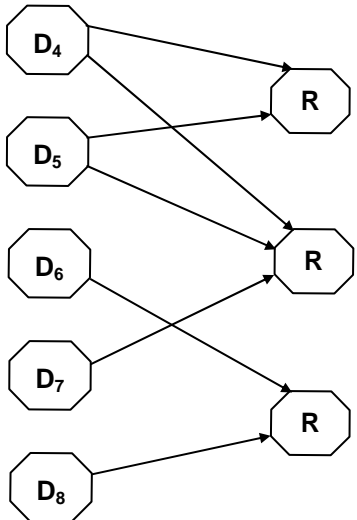

FIGURE: 5.3.9

The related connection matrix is given by

$$
\begin{array}{cccc}
 & R_3 \ \ R_4 \ \ R_5 \\
\begin{array}{c} D_4 \\ D_5 \\ D_6 \\ D_7 \\ D_8 \end{array} &
\left[\begin{array}{ccc}
1 & 1 & 0 \\
1 & 1 & 0 \\
0 & 0 & 1 \\
0 & 1 & 0 \\
1 & 0 & 1
\end{array}\right]
\end{array}
$$

The connection matrix associated with the combined disjoint block FRM of different sizes is given by the $8 \times 5$ matrix and is denoted by G

$$
\begin{array}{c}
\begin{array}{ccccc} R_1 & R_2 & R_3 & R_4 & R_5 \end{array} \\
\begin{array}{c} D_1 \\ D_2 \\ D_3 \\ D_4 \\ D_5 \\ D_6 \\ D_7 \\ D_8 \end{array}
\begin{bmatrix}
1 & 0 & 0 & 0 & 0 \\
1 & 0 & 0 & 0 & 0 \\
0 & 1 & 0 & 0 & 0 \\
0 & 0 & 1 & 1 & 0 \\
0 & 0 & 1 & 1 & 0 \\
0 & 0 & 0 & 0 & 1 \\
0 & 0 & 0 & 1 & 0 \\
0 & 0 & 1 & 0 & 1
\end{bmatrix}
\end{array}
$$

Let us consider the effect of the state vector X = (0 0 0 0 1 0 0 0) where $D_5$ is in the on state and all other attributes are in the off state. The effect of X on G is given by

$$
\begin{array}{llll}
XG & \hookrightarrow & (0\ 0\ 1\ 1\ 0) & = & Y \\
YG^T & \hookrightarrow & (0\ 0\ 0\ 1\ 1\ 0\ 1\ 1) & = & X_1 \text{ (say)} \\
X_1\,G & \hookrightarrow & (0\ 0\ 1\ 1\ 1) & = & Y_1 \\
Y_1\,G^T & \hookrightarrow & (0\ 0\ 0\ 1\ 1\ 1\ 1\ 1) & = & X_1 \text{ say} \\
X_1\,G & \hookrightarrow & (0\ 0\ 1\ 1\ 1) & = & Y_2 = Y_1.
\end{array}
$$

Thus the hidden pattern is a fixed point given by the binary pair {(0 0 1 1 1), (0 0 0 1 1 1 1 1)}, which shows its strong effect on the dynamical system.

Consider the state vector
$$Y = (1\ 0\ 0\ 0\ 0).$$

$$
\begin{array}{llll}
YG^T & \hookrightarrow & (1\ 1\ 0\ 0\ 0\ 0\ 0\ 0) & = & X \\
XG & \hookrightarrow & (1\ 0\ 0\ 0\ 0) & = & Y_1 (= Y).
\end{array}
$$

Thus the fixed binary pair is given by {(1 0 0 0 0), (1 1 0 0 0 0 0 0)}.

Consider the FRM given in pages 204 - 209 of this chapter. The attributes are taken as $\{M_1, M_2, \ldots, M_{10}\}$ and $\{G_1, G_2, \ldots, G_7\}$. The combined block disjoint FRM of varying sizes is calculated using the classes

$$C_1 = \{(M_1, M_2\ M_3\ M_4)\ (G_1\ G_2\ G_3)\}$$
$$C_2 = \{(M_5\ M_6\ M_7\ M_8)\ (G_4\ G_5)\}\ \text{and}$$
$$C_3 = \{(M_9\ M_{10}),\ (G_6\ G_7)\}.$$

The directed graph of class $C_1$ as given by the expert is as follows:

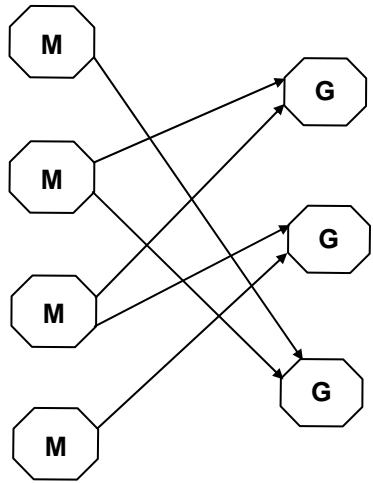

FIGURE: 5.3.10

The related connection matrix

$$
\begin{array}{c}
\\
M_1 \\
M_2 \\
M_3 \\
M_4
\end{array}
\begin{array}{ccc}
G_1 & G_2 & G_3 \\
\left[\begin{array}{ccc}
0 & 0 & 1 \\
1 & 0 & 1 \\
1 & 1 & 0 \\
0 & 1 & 0
\end{array}\right]
\end{array}
$$

The directed graph relative to $C_2$ as given by the expert is given in the following:

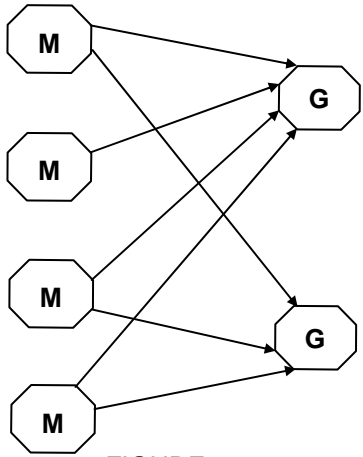

FIGURE: 5.3.11

The related connection matrix

$$
\begin{array}{c}
\phantom{M_5}\begin{array}{cc} G_4 & G_5 \end{array} \\
\begin{array}{c} M_5 \\ M_6 \\ M_7 \\ M_8 \end{array}
\begin{bmatrix} 1 & 1 \\ 0 & 1 \\ 1 & 1 \\ 1 & 1 \end{bmatrix}
\end{array}
$$

The directed graph related to $C_3$ given by the expert is

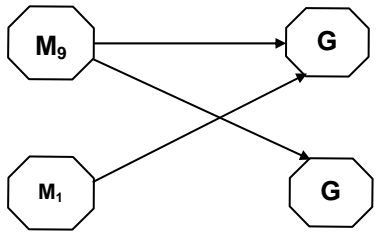

The related connection matrix

$$
\begin{array}{c}
\phantom{M_9}\begin{array}{cc} G_6 & G_7 \end{array} \\
\begin{array}{c} M_9 \\ M_{10} \end{array}
\begin{bmatrix} 1 & 1 \\ 0 & 1 \end{bmatrix}
\end{array}
$$

Now we give the relational matrix associated with the combined Block disjoint FRM is a $10 \times 7$ matrix denoted by (CP)

$$
\begin{array}{c}
\phantom{M_{10}} \begin{array}{ccccccc} G_1 & G_2 & G_3 & G_4 & G_5 & G_6 & G_7 \end{array} \\
\begin{array}{c} M_1 \\ M_2 \\ M_3 \\ M_4 \\ M_5 \\ M_6 \\ M_7 \\ M_8 \\ M_9 \\ M_{10} \end{array}
\begin{bmatrix}
0 & 0 & 1 & 0 & 0 & 0 & 0 \\
1 & 0 & 1 & 0 & 0 & 0 & 0 \\
1 & 1 & 0 & 0 & 0 & 0 & 0 \\
0 & 1 & 0 & 0 & 0 & 0 & 0 \\
0 & 0 & 0 & 1 & 1 & 0 & 0 \\
0 & 0 & 0 & 0 & 1 & 0 & 0 \\
0 & 0 & 0 & 1 & 1 & 0 & 0 \\
0 & 0 & 0 & 1 & 1 & 0 & 0 \\
0 & 0 & 0 & 0 & 0 & 1 & 1 \\
0 & 0 & 0 & 0 & 0 & 0 & 1
\end{bmatrix}
\end{array}
$$

Let us consider the state vector $X = (0\ 0\ 0\ 0\ 0\ 0\ 0\ 1\ 0\ 0)$ where the node $M_8$ alone is in the on state and all other nodes are in the off state. The effect of X on the dynamical system C(P)

$$X\,C(P) \quad \hookrightarrow \quad (0\ 0\ 0\ 1\ 1\ 0\ 0) \quad = \quad Y$$

$$Y C(P)^T \quad \hookrightarrow \quad (0\ 0\ 0\ 0\ 1\ 1\ 1\ 1\ 0\ 0) \quad = \quad X_1$$

$$X_1\,C(P) \quad \hookrightarrow \quad (0\ 0\ 0\ 1\ 1\ 0\ 0) \quad = \quad Y_1 = Y.$$

Thus the fixed point is a binary pair given by $\{(0\ 0\ 0\ 0\ 1\ 1\ 1\ 1\ 0\ 0), (0\ 0\ 0\ 1\ 1\ 0\ 0)\}$ it makes on the attributes $M_5$ to $M_{10}$, which shows the strong impact of $M_8$ on the other nodes.

We can calculate the result and state vector of any state vector using the C-program given in the appendex of this book. Next we consider the combined overlap block FRM used in modeling the HIV/AIDS migrant laboures problem. Let us consider the attributes given page 198 of this chapter. The classes of attributes used are $\{D_1, D_2, \ldots, D\}$ and $\{R_1, R_2, \ldots, R_5\}$. We divide these into overlap blocks $C_1, C_2, C_3$ and $C_4$,

$$C_1 = \{(D_1\ D_2\ D_3\ D_4), (R_1\ R_2\ R_3)\},$$
$$C_2 = \{(D_3\ D_4\ D_5\ D_6), (R_2\ R_3\ R_4)\}$$
$$C_3 = \{(D_5\ D_6\ D_7\ D_8), (R_3\ R_4\ R_5)\} \text{ and}$$
$$C_4 = \{(D_7\ D_8\ D_1\ D_2), (R_4\ R_5\ R_1)\}.$$

Now we obtain the directed graph of these four FRMs. The directed graph related to the class $C_1$.

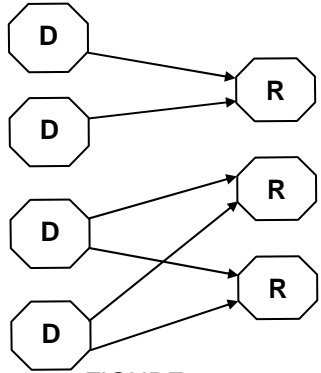

FIGURE: 5.3.13

The related relation matrix

$$
\begin{array}{c}
\quad\ \ R_1\ R_2\ R_3 \\
\begin{array}{c}
D_1 \\ D_2 \\ D_3 \\ D_4
\end{array}
\left[
\begin{array}{ccc}
1 & 0 & 0 \\
1 & 0 & 0 \\
0 & 1 & 1 \\
0 & 1 & 1
\end{array}
\right]
\end{array}
$$

The directed graph given by the expert of class $C_2$ is as follows:

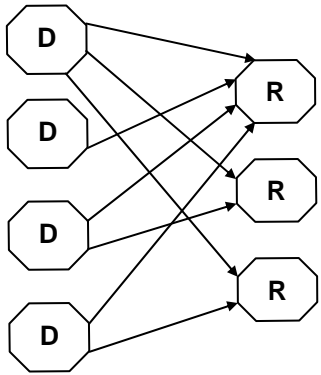

FIGURE: 5.3.14

The related relational matrix

$$\begin{array}{c@{\quad}ccc} & R_2 & R_3 & R_4 \\ D_3 & 1 & 1 & 1 \\ D_4 & 1 & 0 & 0 \\ D_5 & 1 & 1 & 0 \\ D_6 & 1 & 0 & 1 \end{array}$$

Next we consider the directed graph given by the expert for $C_3$.

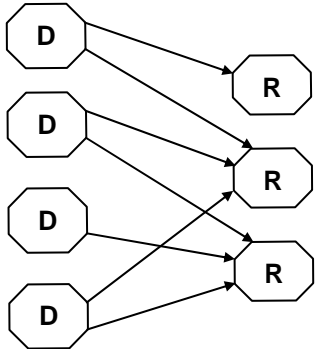

FIGURE: 5.3.15

The related connection matrix is

$$\begin{array}{c@{\quad}ccc} & R_3 & R_4 & R_5 \\ D_5 & 1 & 1 & 0 \\ D_6 & 1 & 1 & 0 \\ D_7 & 0 & 0 & 1 \\ D_8 & 1 & 1 & 0 \end{array}$$

The directed graph for the final class $\{(D_7\ D_8\ D_1\ D_2),\ (R_4,\ R_5\ R_1)$.

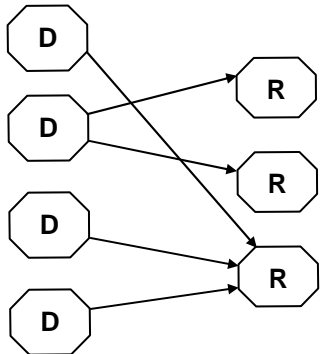

FIGURE: 5.3.16

The related connection matrix

$$
\begin{array}{c}
 & \begin{array}{ccc} R_4 & R_5 & R_6 \end{array} \\
\begin{array}{c} D_7 \\ D_8 \\ D_9 \\ D_{10} \end{array} &
\left[\begin{array}{ccc}
0 & 0 & 1 \\
1 & 1 & 0 \\
0 & 0 & 1 \\
0 & 0 & 1
\end{array}\right]
\end{array}
$$

Now using these four connection matrices we obtain the connection matrix relative to the combined overlap block FRM which is $8 \times 5$ matrix and is denoted by C(T).

$$
\begin{array}{c}
 & \begin{array}{ccccc} R_1 & R_2 & R_3 & R_4 & R_5 \end{array} \\
\begin{array}{c} D_1 \\ D_2 \\ D_3 \\ D_4 \\ D_5 \\ D_6 \\ D_7 \\ D_8 \end{array} &
\left[\begin{array}{ccccc}
2 & 0 & 0 & 0 & 0 \\
2 & 0 & 0 & 0 & 0 \\
0 & 2 & 2 & 1 & 0 \\
0 & 2 & 1 & 0 & 0 \\
0 & 1 & 2 & 1 & 0 \\
0 & 1 & 1 & 1 & 0 \\
1 & 0 & 0 & 0 & 2 \\
0 & 0 & 1 & 2 & 0
\end{array}\right]
\end{array}
$$

Let X = (0 0 0 0 1 0 0 0) be the state vector with the node $D_5$ to be in the on state all other nodes are in the off state. Effect of X on the dynamical system (CT) is given by

$$
\begin{array}{lll}
\text{XC(T)} & \hookrightarrow \ (0\ 1\ 1\ 1\ 0) & = \quad \text{Y} \\
\text{Y (C(T))}^T & \hookrightarrow \ (0\ 0\ 1\ 1\ 1\ 1\ 0\ 1) & = \quad X_1 \\
X_1 \ \text{(CT))} & \hookrightarrow \ (0\ 1\ 1\ 1\ 0) & = \quad Y_1 \ (Y_1 \ \text{say).}
\end{array}
$$

Thus the hidden pattern is a fixed point given by the binary pair {(0 1 1 1 0), (0 0 1 1 1 1 10 1)}, the reader is expected to analyse the binary pair.

Consider the state vector Y = (1 0 0 0 0) ie only the attribute $R_1$ is in the on state all other nodes are in the off state. Effect of Y on the dynamical system C(T).

$$Y(C(T))^T \hookrightarrow (1\ 1\ 0\ 0\ 0\ 0\ 1\ 0) = X$$
$$XC(T) \hookrightarrow (1\ 0\ 0\ 0\ 1) = Y_1$$
$$Y_1(C(T))^T \hookrightarrow (1\ 1\ 0\ 0\ 0\ 0\ 10) = X_1\ (= X_1).$$

Thus the fixed point is a binary pair {(1 1 0 0 0 0 1 0), (10 0 0 1)}. This node does not have a very strong impact on the other dynamical system.

Now we consider the combined overlap Block FRM to study the model given in page 204 - 205 of this chapter. Take the sets of attributes. {$M_1\ M_2,\ldots, M_{10}$}, ($G_1\ G_2,\ldots, G_7$)}.

We divide them into overlapping blocks given by $C_1\ C_2$ and $C_3$ where

$$C1 = \{(M_1, M_2\ M_3\ M_4)\ (G_1\ G_2\ G_3)\},$$
$$C2 = \{(M_4, M_5\ M_6\ M_7)\ (G_3\ G_4\ G_5)\}\ \text{and}$$
$$C3 = \{(M_7, M_8\ M_9\ M_{10})\ (G_5\ G_6\ G_7)\}.$$

The directed graph given by the expert for the class $C_1$ is

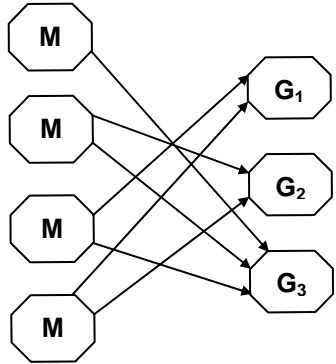

FIGURE: 5.3.17

The related relational connection matrix

$$\begin{array}{c} \\ M_1 \\ M_2 \\ M_3 \\ M_4 \end{array} \begin{array}{ccc} G_1 & G_2 & G_3 \\ \begin{bmatrix} 0 & 0 & 1 \\ 0 & 1 & 1 \\ 1 & 0 & 1 \\ 1 & 1 & 0 \end{bmatrix} \end{array}$$

Now the directed graph given by the expert for class $C_2$.

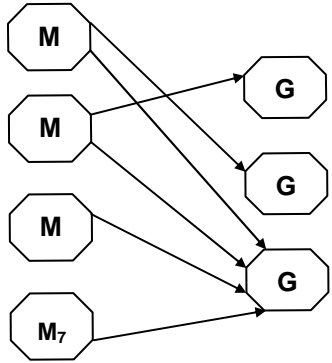

FIGURE: 5.3.18

The related connection matrix is

$$
\begin{array}{c c c}
 & G_3 & G_4 & G_5
\end{array}
$$
$$
\begin{array}{c}
M_4 \\
M_5 \\
M_6 \\
M_7
\end{array}
\begin{bmatrix}
0 & 1 & 1 \\
1 & 0 & 1 \\
0 & 0 & 1 \\
0 & 0 & 1
\end{bmatrix}
$$

Now the directed graph given by the expert for the last class.

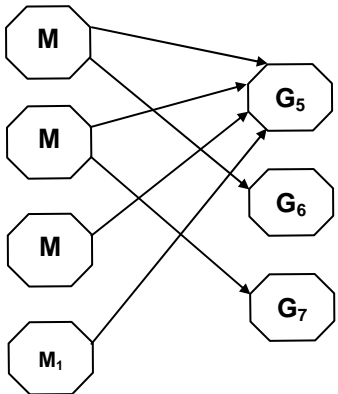

FIGURE: 5.3.19

The relational connection matrix is given by

$$
\begin{array}{c c c c}
 & G_5 & G_6 & G_7 \\
M_7 & \begin{bmatrix} 1 & 1 & 0 \\ 1 & 0 & 1 \\ 1 & 0 & 0 \\ 1 & 0 & 0 \end{bmatrix} \\
M_8 & & \\
M_9 & & \\
M_{10} & &
\end{array}
$$

Now using these three connection matrices give the $10 \times 7$ connection matrix associated with the combined block over lap FRM where the blocks are of equal size, we denote this matrix by C(Q)

$$
\begin{array}{c c c c c c c c}
 & G_1 & G_2 & G_3 & G_4 & G_5 & G_6 & G_7 \\
M_1 & \begin{bmatrix} 0 & 0 & 1 & 1 & 1 & 0 & 0 \\ 0 & 1 & 1 & 0 & 1 & 0 & 0 \\ 1 & 0 & 1 & 0 & 1 & 0 & 0 \\ 1 & 1 & 0 & 1 & 1 & 0 & 0 \\ 0 & 0 & 1 & 0 & 1 & 1 & 1 \\ 0 & 0 & 0 & 0 & 1 & 0 & 0 \\ 0 & 0 & 0 & 0 & 1 & 1 & 0 \\ 0 & 0 & 0 & 0 & 1 & 0 & 1 \\ 0 & 0 & 0 & 0 & 1 & 0 & 0 \\ 0 & 0 & 0 & 0 & 1 & 0 & 0 \end{bmatrix} \\
M_2 & \\
M_3 & \\
M_4 & \\
M_5 & \\
M_6 & \\
M_7 & \\
M_8 & \\
M_9 & \\
M_{10} &
\end{array}
$$

Let X = (0 0 0 0 0 0 0 1 0 0) ie only $M_8$ node is in the on state all other nodes are in the off state. The effect of X on the system C(Q) is given by

$$
\begin{aligned}
XC(Q) &\hookrightarrow (0\ 0\ 0\ 0\ 1\ 0\ 1) &=& \quad Y \\
Y(C(Q))^T &\hookrightarrow (1\ 1\ 1\ 1\ 1\ 1\ 1\ 1\ 1\ 1) &=& \quad X_1 \text{ (say)} \\
X_1\, C(Q) &\hookrightarrow (1\ 1\ 1\ 1\ 1\ 1\ 1) &=& \quad Y_1 \\
Y_1\, C(Q)^T &\hookrightarrow X_2 = X_1.
\end{aligned}
$$

Thus the fixed point is a binary pair given by {(1 1 1 1 1 1 1), (1 1 1 1 1 1 1 1 1 1)}, which shows the strong impact of node M8 on the total dynamical system as all nodes come to on state in both the domain and range space. Let us consider the state vector Y = (0 0 1 0 0 0 0 0) ie only the state $G_3$ is in the on state and all other vectors are in the off state effect of Y on C(Q)

$$
YC(Q)^T \quad \hookrightarrow \quad (0\ 1\ 1\ 10\ 1\ 0\ 0\ 0\ 0) \quad = \quad X
$$

$$ZC(Q) \quad \hookrightarrow \quad (1\ 1\ 1\ 1\ 1\ 0\ 0) \quad\quad = \quad Y_1 \text{ say}$$
$$Y_1 C(Q)^T \quad \hookrightarrow \quad (1\ 1\ 1\ 1\ 1\ 1\ 1\ 1\ 1\ 1) \ = \quad X_1$$
$$X_1 C(Q) \quad \hookrightarrow \quad (1\ 1\ 1\ 1\ 1\ 1\ 1) \quad\quad = \quad Y_2 = Y_1.$$

Thus the hidden pattern is a fixed point given by the binary pair $\{(1\ 1\ 1\ 1\ 1\ 1\ 1\ 1\ 1\ 1), (1\ 1\ 1\ 1\ 1\ 1\ 1)\}$. That is the combined overlap block FRM is very sensitive and it makes all the nodes to be in the on state. Next we consider the combined overlap block FRM of different sizes. Let us consider the attributes given in page of this chapter, with the sets of attributes $\{(D_1\ D_2\ D_3\ D_4\ D_5\ D_6), (R_1\ R_2, \ldots, R_{10})\}$. Let us form the classes

$$C_1 = \{(D_1\ D_2\ D_3\ D_4), (R_1\ R_2\ R_3\ R_4\ R_5)\}$$
$$C_2 \ \{(D_3\ D_4\ D_5), (R_3\ R_4\ R_5\ R_6)\} \text{ and}$$
$$C_3 = \{(D_4\ D_5\ D_6\ D_1\ D_2), (R_6\ R_7\ R\ R_9\ R_{10})\}.$$

Now we study the model implications. The directed graph of the FRM related to class $C_1$ as given by the expert is as follows

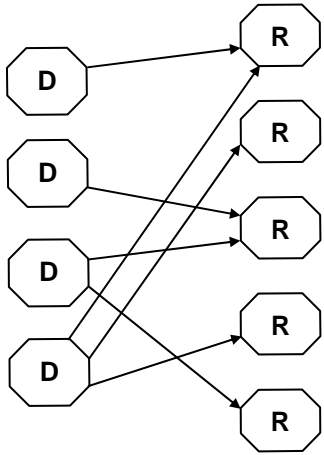

FIGURE: 5.3.20

The related relational matrix
$$R_1\ R_2\ R_3\ R_4\ R_5$$

$$\begin{array}{c} D_1 \\ D_2 \\ D_3 \\ D_4 \end{array} \begin{bmatrix} 1 & 0 & 0 & 0 & 0 \\ 0 & 0 & 1 & 0 & 0 \\ 0 & 0 & 1 & 0 & 1 \\ 1 & 1 & 0 & 1 & 0 \end{bmatrix}$$

The directed graph given by the expert for the class $C_2$ is

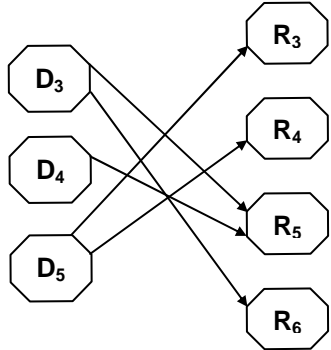

FIGURE: 5.3.21

The related connection matrix is

$$\begin{array}{cccc} & R_3 & R_4 & R_5 & R_6 \end{array}$$
$$\begin{array}{c} D_3 \\ D_4 \\ D_5 \end{array} \begin{bmatrix} 0 & 0 & 1 & 1 \\ 0 & 0 & 1 & 0 \\ 1 & 1 & 0 & 0 \end{bmatrix}$$

The directed graph given by the expert for of the last class $C_3$.

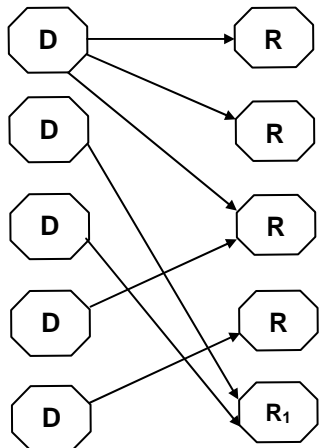

FIGURE: 5.3.22

The related connection matrix of he class $C_3$ is

$$
\begin{array}{c@{\;}c}
& \begin{array}{ccccc} R_6 & R_7 & R_8 & R_9 & R_{10} \end{array} \\
\begin{array}{c} D_4 \\ D_5 \\ D_6 \\ D_1 \\ D_2 \end{array} &
\left[\begin{array}{ccccc}
1 & 1 & 1 & 0 & 0 \\
0 & 0 & 0 & 0 & 1 \\
0 & 0 & 0 & 0 & 1 \\
0 & 0 & 1 & 0 & 0 \\
0 & 0 & 0 & 1 & 0
\end{array}\right]
\end{array}
$$

Now we give the connection matrix C(T) of the combined overlap block FRM. This is a $6 \times 10$ matrix

$$
\begin{array}{c@{\;}c}
& \begin{array}{cccccccccc} R_1 & R_2 & R_3 & R_4 & R_5 & R_6 & R_7 & R_8 & R_9 & R_{10} \end{array} \\
\begin{array}{c} D_1 \\ D_2 \\ D_3 \\ D_4 \\ D_5 \\ D_6 \end{array} &
\left[\begin{array}{cccccccccc}
1 & 0 & 0 & 0 & 0 & 0 & 0 & 1 & 0 & 0 \\
0 & 0 & 1 & 0 & 0 & 0 & 0 & 0 & 1 & 0 \\
0 & 0 & 1 & 1 & 1 & 1 & 0 & 0 & 0 & 0 \\
1 & 1 & 0 & 2 & 0 & 1 & 1 & 1 & 0 & 0 \\
0 & 0 & 1 & 1 & 0 & 0 & 0 & 0 & 0 & 1 \\
0 & 0 & 0 & 0 & 0 & 0 & 0 & 0 & 0 & 1
\end{array}\right]
\end{array}
$$

Now we study the effect of the dynamical system on any state vector X = (0 0 0 0 1 0) ie only the node $D_5$ is in the on state all nodes are in the off state.

$$
\begin{array}{llll}
\text{X(CT))} & \hookrightarrow & (0\,0\,1\,1\,0\,0\,0\,0\,0\,1) & = & Y \\
\text{Y(C(T))}^T & \hookrightarrow & (0\,1\,1\,1\,1\,1) & = & X_1 \\
\text{X}_1\text{ (CT))} & \hookrightarrow & (1\,1\,1\,1\,1\,1\,1\,1\,1\,1) & = & Y_1 \\
\text{Y}_1\text{ C(T)}^T & \hookrightarrow & (1\,1\,1\,1\,1\,1\,1).
\end{array}
$$

Thus the fixed point is the binary pair {(1 1 1 1 1 1 1 1 1 1), (1 1 1 1 1 1)}. Thus the node D5 has a very strong influence on all nodes making all of them to on state.

Let us consider the state vector Y = (0 0 1 0 0 0) where the attribute $D_3$ alone is in the on state and all nodes are in the off state. The effect of Y on C(T) is given by

$$YC(T) \quad \hookrightarrow \quad (0\ 0\ 1\ 1\ 1\ 1\ 0\ 0\ 0\ 0) \quad = \quad X$$

$$X\ (C(T))^T \quad \hookrightarrow \quad (0\ 1\ 1\ 1\ 1\ 0) \quad\quad = \quad Y_1$$

$$Y_1\ C(T) \quad \hookrightarrow \quad (1\ 1\ 1\ 1\ 1\ 1\ 1\ 1\ 1\ 1) \quad = \quad X_1$$

$$X_1\ C(T)^T \quad \hookrightarrow \quad (1\ 1\ 1\ 1\ 1\ 1).$$

Thus the binary pair is a fixed point given by $\{(1\ 1\ 1\ 1\ 1\ 1)$, $(1\ 1\ 1\ 1\ 1\ 1\ 1\ 1\ 1\ 1)\}$. Both the nodes $D_3$ and $D_5$ have same impact on the dynamical system for the fixed point in both the cases are the same.

We study another set of attributes the HIV/AIDS affected migrant labourers and their related problem. Let us consider the sets of attributes $\{(A_1\ A_2\ A_3,\ \ldots,\ A_6),\ (G_1\ G_2, \ldots, G_5)\}$ is divided into over lapping classes by

$$C_1 = \{(A_1\ A_2\ A_3\ A_4),\ (G_1\ G_2, G_3)\} \text{ and}$$
$$C_2\ = \{(A_4\ A_5\ A_6),\ (G_3\ G_4\ G_5)\}$$

The directed graph relating to the class $C_1$ as given by the expert is as follows:

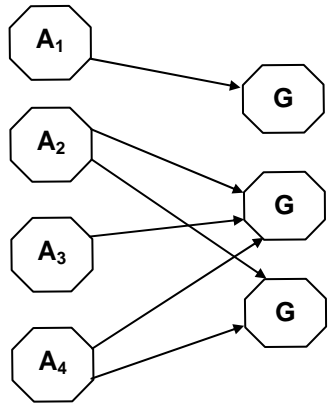

FIGURE: 5.3.23

The related connection matrix

$$G_1\ \ G_2\ \ G_3$$

$$\begin{array}{c}A_1\\A_2\\A_3\\A_4\end{array}\begin{bmatrix}1&0&0\\0&1&1\\0&1&0\\0&1&1\end{bmatrix}$$

The directed graph of the FRM related to class $C_2$.

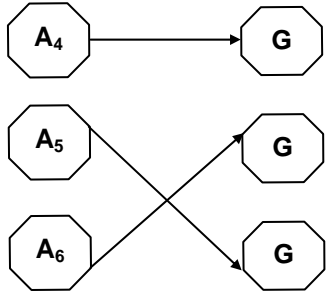

FIGURE: 5.3.24

The related connection matrix

$$\begin{array}{cccc}&G_3&G_4&G_6\\A_4&\begin{bmatrix}1&0&0\\A_5&0&0&1\\A_6&0&1&0\end{bmatrix}\end{array}$$

The matrix related to the combined overlap block FRM denoted by C(V)

$$\begin{array}{cccccc}&G_1&G_2&G_3&G_4&G_5\\A_1&\begin{bmatrix}1&0&0&0&0\\A_2&0&1&1&0&0\\A_3&0&1&0&0&0\\A_4&0&1&2&0&0\\A_5&0&0&0&0&1\\A_6&0&0&0&1&0\end{bmatrix}\end{array}$$

Consider the state vector $X = (0\ 0\ 0\ 1\ 0\ 0)$ in which only the node $A_4$ is in the on state and rest of the nodes in the off state. The effect of X on C(V)

$$X\,C(V) \quad \hookrightarrow \quad (0\;1\;1\;0\;0) \quad = \quad Y$$

$$Y(C(V))^T \quad \hookrightarrow \quad (0\;1\;1\;1\;0\;0) \quad = \quad X_1$$

$$X_1\,C(V) \quad \hookrightarrow \quad (0\;1\;1\;0\;0) \quad = \quad Y_1\,(=Y).$$

The hidden pattern of the dynamical system in a fixed point given by {(0 1 1 0 0), (0 1 1 1 0 0)}, the effect of $G_1$ in the dynamical system is nothing in the view of this expert.

Next consider the state vector $Y = (1\;0\;0\;0\;0)$ ie only the node $G_1$ is in the on state all vectors are in off state effect of the vector Y on the dynamical system C(V) is given by

$$YC(V)^T \quad \hookrightarrow \quad (1\;0\;0\;0\;0\;0) \quad = \quad X$$

$$X\,C(V) \quad \hookrightarrow \quad (1\;0\;0\;0\;0) \quad = \quad Y_1 \;=Y$$

Thus we see the hidden pattern is a fixed point given by the binary pair {(1 0 0 0 0), (1 0 0 0 0 0)}. The impact of G1 in the dynamical system is nothing in the view of this expert.

## 5.4. NRMs to Analyze the Problem of HIV/AIDS Migrant Labourers

Neutrosophic Cognitive Maps (NCMs) promote the causal relationships between concurrently active units or decides the absence of any relation between two units or the indeterminance of any relation between any two units. But in Neutrosophic Relational Maps (NRMs) we divide the very causal nodes into two disjoint units. Thus for the modeling of a NRM we need a domain space and a range space which are disjoint in the sense of concepts. We further assume no intermediate relations exist within the domain and the range spaces. The number of elements or nodes in the range space need not be equal to the number of elements or nodes in the domain space.

Throughout this section we assume the elements of a domain space are taken from the neutrosophic vector space of dimension n and that of the range space are neutrosophic vector space of dimension m. (m in general need not be equal to n). We denote by R the set of nodes $R_1,\ldots, R_m$ of the range space, where $R = \{(x_1,\ldots, x_m) \mid x_j = 0 \text{ or } 1 \text{ for } j = 1, 2, \ldots, m\}$.

If $x_i = 1$ it means that node $R_i$ is in the on state and if $x_i = 0$ it means that the node $R_i$ is in the off state and if $x_i = I$ in the resultant vector it means the effect of the node $x_i$ is indeterminate or whether it will be off or on cannot be predicted by the neutrosophic dynamical system.

It is very important to note that when we send the state vectors they are always taken as the real state vectors for we know the node or the concept is in the on state or in the off state but when the state vector passes through the Neutrosophic dynamical system some other node may become indeterminate i.e. due to the presence of a node we may not be able to predict the presence or the absence of the other node i.e., it is indeterminate, denoted by the symbol I, thus the resultant vector can be a neutrosophic vector.

**DEFINITION 5.4.1:** *A Neutrosophic Relational Map (NRM) is a Neutrosophic directed graph or a map from D to R with concepts like policies or events etc. as nodes and causalities as edges. (Here by causalities we mean or include the indeterminate causalities also). It represents Neutrosophic Relations and Causal Relations between spaces D and R .*

*Let $D_i$ and $R_j$ denote the nodes of an NRM. The directed edge from $D_i$ to $R_j$ denotes the causality of $D_i$ on $R_j$ called relations. Every edge in the NRM is weighted with a number in the set {0, +1, –1, I}. Let $e_{ij}$ be the weight of the edge $D_i R_j$, $e_{ij} \in \{0, 1, -1, I\}$. The weight of the edge $D_i R_j$ is positive if increase in $D_i$ implies increase in $R_j$ or decrease in $D_i$ implies decrease in $R_j$ i.e. causality of $D_i$ on $R_j$ is 1. If $e_{ij} = -1$ then increase (or decrease) in $D_i$ implies decrease (or increase) in $R_j$. If $e_{ij} = 0$ then $D_i$ does not have any effect on $R_j$. If $e_{ij} = I$ it implies we are not in a position to determine the effect of $D_i$ on $R_j$ i.e. the effect of $D_i$ on $R_j$ is an indeterminate so we denote it by I.*

**DEFINITION 5.4.2:** *When the nodes of the NRM take edge values from {0, 1, –1, I} we say the NRMs are simple NRMs.*

**DEFINITION 5.4.3:** *Let $D_1, …, D_n$ be the nodes of the domain space D of an NRM and let $R_1, R_2,…, R_m$ be the nodes of the range space R of the same NRM. Let the matrix N(E) be defined as N(E) = $(e_{ij})$ where $e_{ij}$ is the weight of the directed edge $D_i R_j$ (or $R_j D_i$ ) and $e_{ij} \in \{0, 1, -1, I\}$. N(E) is called the Neutrosophic Relational Matrix of the NRM.*

The following remark is important and interesting to find its mention in this book.

**Remark**: Unlike NCMs, NRMs can also be rectangular matrices with rows corresponding to the domain space and columns corresponding to the range space. This is one of the marked difference between NRMs and NCMs. Further the number of entries for a particular model which can be treated as disjoint sets when dealt as a NRM has very much less entries than when the same model is treated as a NCM.

Thus in many cases when the unsupervised data under study or consideration can be spilt as disjoint sets of nodes or concepts; certainly NRMs are a better tool than the NCMs.

**DEFINITION 5.4.4:** *Let $D_1, \ldots, D_n$ and $R_1, \ldots, R_m$ denote the nodes of a NRM. Let $A = (a_1, \ldots, a_n)$, $a_i \in \{0, 1, -1\}$ is called the Neutrosophic instantaneous state vector of the domain space and it denotes the on-off position of the nodes at any instant. Similarly let $B = (b_1, \ldots, b_n)$ $b_i \in \{0, 1, -1\}$, B is called instantaneous state vector of the range space and it denotes the on-off position of the nodes at any instant, $a_i = 0$ if $a_i$ is off and $a_i = 1$ if $a_i$ is on for $i = 1, 2, \ldots, n$. Similarly, $b_i = 0$ if $b_i$ is off and $b_i = 1$ if $b_i$ is on for $i = 1, 2, \ldots, m$.*

**DEFINITION 5.4.5:** *Let $D_1, \ldots, D_n$ and $R_1, R_2, \ldots, R_m$ be the nodes of a NRM. Let $D_i R_j$ (or $R_j D_i$) be the edges of an NRM, $j = 1, 2, \ldots, m$ and $i = 1, 2, \ldots, n$. The edges form a directed cycle. An NRM is said to be a cycle if it possess a directed cycle. An NRM is said to be acyclic if it does not possess any directed cycle.*

**DEFINITION 5.4.6:** *A NRM with cycles is said to be a NRM with feedback.*

**DEFINITION 5.4.7:** *When there is a feedback in the NRM i.e. when the causal relations flow through a cycle in a revolutionary manner the NRM is called a Neutrosophic dynamical system.*

**DEFINITION 5.4.8:** *Let $D_i R_j$ (or $R_j D_i$) $1 \leq j \leq m$, $1 \leq i \leq n$, when $R_j$ (or $D_i$) is switched on and if causality flows through edges of a cycle and if it again causes $R_j$ (or $D_i$) we say that the Neutrosophical dynamical system goes round and round. This is true for any node $R_j$ (or $D_i$) for $1 \leq j \leq m$ (or $1 \leq i \leq n$). The* *equilibrium state of this Neutrosophical dynamical system is called the Neutrosophic hidden pattern.*

**DEFINITION 5.4.9:** *If the equilibrium state of a Neutrosophical dynamical system is a unique Neutrosophic state vector, then it is called the fixed point. Consider an NRM with $R_1$, $R_2$, ..., $R_m$ and $D_1$, $D_2$,..., $D_n$ as nodes. For example let us start the dynamical system by switching on $R_1$ (or $D_1$). Let us assume that the NRM settles down with $R_1$ and $R_m$ (or $D_1$ and $D_n$) on, or indeterminate on, i.e. the Neutrosophic state vector remains as (1, 0, 0,..., 1) or (1, 0, 0,...I) (or (1, 0, 0,...1) or (1, 0, 0,...I) in D), this state vector is called the fixed point.*

**DEFINITION 5.4.10:** *If the NRM settles down with a state vector repeating in the form $A_1 \rightarrow A_2 \rightarrow A_3 \rightarrow ... \rightarrow A_i \rightarrow A_1$ (or $B_1 \rightarrow B_2 \rightarrow ... \rightarrow B_i \rightarrow B_1$) then this equilibrium is called a limit cycle.*

**METHODS OF DETERMINING THE HIDDEN PATTERN IN A NRM**

Let $R_1$, $R_2$,..., $R_m$ and $D_1$, $D_2$,..., $D_n$ be the nodes of a NRM with feedback. Let N(E) be the Neutrosophic Relational Matrix. Let us find the hidden pattern when $D_1$ is switched on i.e. when an input is given as a vector; $A_1 = (1, 0, ..., 0)$ in D; the data should pass through the relational matrix N(E).

This is done by multiplying $A_1$ with the Neutrosophic relational matrix N(E). Let $A_1N(E) = (r_1, r_2..., r_m)$ after thresholding and updating the resultant vector we get $A_1E \in R$, Now let $B = A_1E$ we pass on B into the system $(N(E))^T$ and obtain $B(N(E))^T$. We update and threshold the vector $B(N(E))^T$ so that $B(N(E))^T \in D$.

This procedure is repeated till we get a limit cycle or a fixed point.

**DEFINITION 5.4.11:** *Finite number of NRMs can be combined together to produce the joint effect of all NRMs. Let $N(E_1)$, $N(E_2)$,..., $N(E_r)$ be the Neutrosophic relational matrices of the NRMs with nodes $R_1$,..., $R_m$ and $D_1$,...,$D_n$, then the combined NRM is represented by the neutrosophic relational matrix $N(E) = N(E_1) + N(E_2) +...+ N(E_r)$.*

In this section we for the first time study the socio economic problems of the HIV/AIDS affected migrant labourers. This study is pertinent for when we were seeking experts opinion some made it very clear that they could not say the relationship i.e., they were unaware of whether a relation existed between two nodes under the investigation. As fuzzy relational maps of fuzzy cognitive maps does not help us in marking indeterminancy we are justified in adopting NRM (i.e., Neutrosophic relational maps for investigation). Now we use the NRM model for the attributes given in page 198 of this chapter. Let us take for the domain space D of the NRM $\{D_1, D_2,\ldots, D_8\}$ and the range space R = $\{R_1, R_2,\ldots R_5\}$ already explained in this chapter.

Now using the experts opinion we give the neutrosophic directed graph.

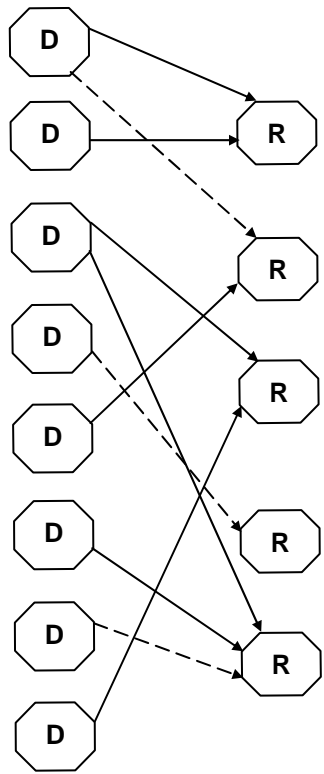

FIGURE: 5.4.1

We using the neutrosophic directed graph give the related neutrosophic connection matrix N(A)

$$
\begin{array}{c c c c c c}
 & R_1 & R_2 & R_3 & R_4 & R_5 \\
D_1 & \begin{bmatrix} 1 & I & 0 & 0 & 0 \\ D_2 & 1 & 0 & 0 & 0 & 0 \\ D_3 & 0 & 1 & 0 & 0 & 1 \\ D_4 & 0 & 0 & 0 & I & 0 \\ D_5 & 0 & 1 & 0 & 0 & 0 \\ D_6 & 0 & 0 & 0 & 0 & 1 \\ D_7 & 0 & 0 & 0 & 0 & I \\ D_8 & 0 & 0 & 1 & 0 & 0 \end{bmatrix}
\end{array}
$$

Let us consider the state vector X = (0 0 0 1 0 0 0 0) ie only the attribute $D_4$ is in the on state and all other attributes are in the off state XN(A)

$$
\begin{array}{llll}
\text{XN(A)} & \hookrightarrow & (0\ 0\ 0\ I\ 0) & = & \text{Y} \\
\text{YN(A)}^{\text{T}} & \hookrightarrow & (0\ 0\ 0\ 1\ 0\ 0\ 0\ 0) & = & X_1\ (=X)
\end{array}
$$

Thus we see that this node has only impact to make the node $D_4$ an indeterminate. Now consider the state vector Y = (0 0 0 0 1) where the node $R_5$ alone is in the on state all other state vectors are in the off state.

$$
\begin{array}{llll}
\text{YN(A)}^{\text{T}} & \hookrightarrow & (0\ 0\ 1\ 0\ 0\ 1\ I\ 0) & = & \text{X} \\
\text{XN(A)} & \hookrightarrow & (0\ 1\ 0\ 0\ 1) & = & Y_1\ (= Y) \\
Y_1\text{N(A)}^{\text{T}} & \hookrightarrow & (I\ 0\ 1\ 0\ 1\ 1\ I\ 0) & = & X_1\ \text{(say)} \\
X_1\ \text{N(A)} & \hookrightarrow & (I\ 1\ 0\ 0\ 1) & = & Y_2\ \text{say} \\
Y_2\ \text{N(A)}^{\text{T}} & \hookrightarrow & (I\ I\ 1\ 0\ 1\ 1\ I\ 0) & = & X_2\ \text{say} \\
\text{XN(A)} & \hookrightarrow & (I\ 1\ 0\ 0\ 1) & = & Y_3\ (=Y_2)
\end{array}
$$

Thus the hidden pattern of the dynamic system is a fixed point given by the binary pair {$(I\ I\ 1\ 0\ 1\ 1\ I\ 0)$, $(I\ 1\ 0\ 0\ 1)$}. The node easily affected by HIV/AIDS leads to the on state of $D_3$, $D_5$ and $D_6$ and indeterminate state of $D_1$, $D_2$ and $D_7$. Also $R_1$ is indeterminate and $R_2$ becomes on. Now consider the state vector X = (1 0 0 0 0 0 0 0 ) where only the attribute $D_1$ is in the on state and all other attributes are in the off state. Now the effect of X on the dynamical system N(A) is given by

$$
\begin{array}{llll}
\text{XN(A)} & \hookrightarrow & (1\ I\ 0\ 0\ 0) & = & \text{Y} \\
\text{Y(N(A))}^{\text{T}} & \hookrightarrow & (I\ 1\ I\ 0\ I\ 0\ 0\ 0) & = & X_1 \\
X_1\text{N (A)} & \hookrightarrow & (1\ I\ 0\ 0\ I) & = & Y_2 \\
Y_2\text{N(A)}^{\text{T}} & \hookrightarrow & (1\ 1\ I\ 0\ I\ I\ I\ 0) & = & X_2 \\
X_2\text{N(A)} & \hookrightarrow & (1\ I\ 0\ 0\ I) & = & Y_3\ (= Y_2).
\end{array}
$$

Thus we get a fixed point given by the binary pair {(1 *I* 0 0 I), (*I* 1 *I* 0 *I I I* 0)} which is the hidden pattern of the dynamical system. Several nodes become indeterminate when D1 alone is in the on state in the state vector. Now we further illustrate the problem of HIV/AIDS affected migrant labourers using NRMs. Take the attributes {$M_1$, $M_2$, …, $M_{10}$} and {$P_1$, $P_2$,…, $P_7$} given in page 222 of this chapter. Now we calculate the combined NRM using experts opinion. We would use 5 experts to find the related NRM for these pair of attributes.

The neutrosophic directed graph given by the first expert for the NRM.

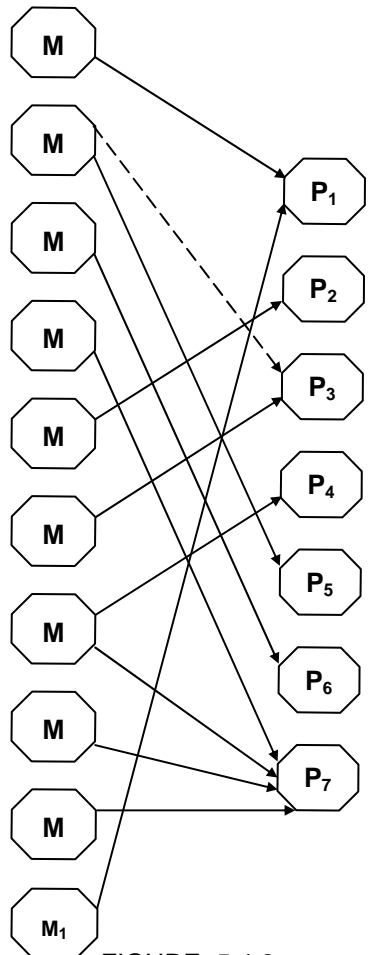

FIGURE: 5.4.2

Now we give the related connection $10 \times 7$ matrix $N(S_1)$

$$
\begin{array}{c}
 \\
M_1 \\
M_2 \\
M_3 \\
M_4 \\
M_5 \\
M_6 \\
M_7 \\
M_8 \\
M_9 \\
M_{10}
\end{array}
\begin{array}{c}
P_1\ P_2\ P_3\ P_4\ P_5\ P_6\ P_7 \\
\left[\begin{array}{ccccccc}
1 & 0 & 0 & 0 & 0 & 0 & 0 \\
0 & 0 & I & 0 & 1 & 0 & 0 \\
0 & 0 & 0 & 0 & 0 & 1 & 0 \\
0 & 0 & 0 & 0 & 0 & 0 & 1 \\
0 & 1 & 0 & 0 & 0 & 0 & 0 \\
0 & 0 & 1 & 0 & 0 & 0 & 0 \\
0 & 0 & 0 & 1 & 0 & 0 & 1 \\
0 & 0 & 0 & 0 & 0 & 0 & 1 \\
0 & 0 & 0 & 0 & 0 & 0 & 1 \\
0 & 0 & 0 & 0 & 0 & 0 & 1
\end{array}\right]
\end{array}
$$

Let us consider the effect of the state vector X = (0 0 0 0 1 0 0 0 0 0) on the dynamical system, in this state vector X only the node $M_5$ is in the on state and all other nodes are in the off state.

$$
\begin{aligned}
&\text{XN } (S_1) &&\hookrightarrow &&(0\ 1\ 0\ 0\ 0\ 0\ 0) &&= &&\text{Y} \\
&N(S_1)^T &&\hookrightarrow &&(0\ 0\ 0\ 0\ 1\ 0\ 0\ 0\ 0\ 0).
\end{aligned}
$$

Thus we get the fixed point. No awareness program by the government ever reaches them as their villages do not have any health centers.

Next let us consider the state vector Y = (0 0 1 0 0 0 0) where only the node $P_3$ is in the on state and all other nodes are in the off state effect of Y on the dynamical system $N(S_1)$ is

$$
\begin{aligned}
&\text{YN}(S_1)^T &&\hookrightarrow &&(0\ I\ 0\ 0\ 0\ 1\ 0\ 0\ 0\ 0) &&= &&\text{X say} \\
&\text{XN}(S_1) &&\hookrightarrow &&(0\ 0\ 1\ 0\ 1\ 0\ 0) &&= &&Y_1 \\
&Y_1 N(S_1)^T &&\hookrightarrow &&(0\ 0\ 0\ 0\ 0\ 1\ 0\ 0\ 0\ 0) &&= &&X_1 \\
&X_1\ N(S_1)^T &&\hookrightarrow &&(0\ 0\ 1\ 0\ 0\ 0) &&= &&Y_2 \\
&Y_2\ N(S_1)^T &&\hookrightarrow &&(0\ I\ 0\ 0\ 0\ 1\ 0\ 0\ 0\ 0) &&= &&X_2 \\
&X_2\ N(S_1) &&\hookrightarrow &&(0\ 0\ 1\ 0\ 1\ 0\ 0) &&= &&Y_3 \\
&Y_3\ N(S_1) &&\hookrightarrow &&(0\ 0\ 0\ 0\ 0\ 1\ 0\ 0\ 0\ 0).
\end{aligned}
$$

Thus we see the hidden pattern is a limit cycle fluctuating between the binary pairs {{(0 0 1 0 0 0 0), (0 *I* 0 0 0 1 0 0 0 0)} and {(0 0 1 0 1 0 0), (0 0 0 0 0 1 0 0 0 0)}.

Now we consider the second experts opinion over the same set of attributes and obtain the directed graph given by him.

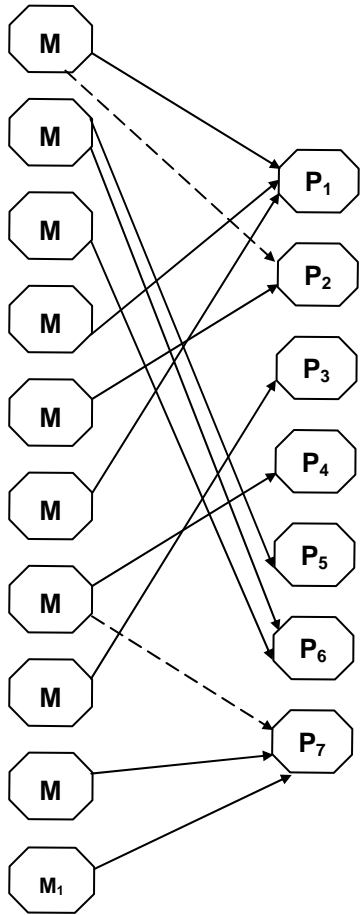

FIGURE: 5.4.3

We obtain the relational neutrosophic $10 \times 7$ matrix $N(S_2)$

$$\begin{array}{c} \\ M_1 \\ M_2 \\ M_3 \\ M_4 \\ M_5 \\ M_6 \\ M_7 \\ M_8 \\ M_9 \\ M_{10} \end{array} \begin{array}{ccccccc} P_1 & P_2 & P_3 & P_4 & P_5 & P_6 & P_7 \\ \begin{bmatrix} 1 & I & 0 & 0 & 0 & 0 & 0 \\ 0 & 0 & 0 & 0 & 1 & 1 & 0 \\ 0 & 0 & 0 & 0 & 0 & 1 & 0 \\ 1 & 0 & 0 & 0 & 0 & 0 & 0 \\ 0 & 1 & 0 & 0 & 0 & 0 & 0 \\ 1 & 0 & 0 & 0 & 0 & 0 & 0 \\ 0 & 0 & 0 & 1 & 0 & 0 & I \\ 0 & 0 & 1 & 0 & 0 & 0 & 0 \\ 0 & 0 & 0 & 0 & 0 & 0 & 1 \\ 0 & 0 & 0 & 0 & 0 & 0 & 1 \end{bmatrix} \end{array}$$

Now using thus connection matrix $N(S_2)$ we can find the effect of any state vector on the dynamical system Let X = (0 0 0 0 1 0 0 0 0 0) be the state vector in which the attribute $M_5$ alone is in the on state and all other nodes are in the off state. Let us study the effect of the state vector X on $N(S_2)$.

$$XN(S_2) \quad \hookrightarrow \quad (0\ 1\ 0\ 0\ 0\ 0\ 0) \quad = \quad Y$$

$$YN(S_2)^T \quad \hookrightarrow \quad (I\ 0\ 0\ 0\ 1\ 0\ 0\ 0\ 0\ 0) \quad = \quad X_1$$

$$X_1\ N(S_2) \quad \hookrightarrow \quad (\pm\ 0\ 0\ I\ 1\ I\ 0\ 0\ 0\ 0) \quad = \quad X_2$$

$$Y_1\ N(S_2)^T \quad \hookrightarrow \quad (I\ 0\ 0\ I\ 1\ I\ 0\ 0\ 0\ 0) \quad = \quad X_2$$

$$X_2\ N(S_2) \quad \hookrightarrow \quad (I\ 0\ 0\ 0\ 0\ 0\ 0) \quad = \quad Y_2 = Y_1.$$

Thus the hidden pattern of the dynamical system is given by the fixed point which is the binary pair {$(I\ 0\ 0\ 0\ 0\ 0\ 0\ 0)$, $(I\ 0\ 0\ I\ 1\ I\ 0\ 0\ 0\ 0)$} which has only indeterministic influence on the dynamical system.

Now let us consider the state vector Y = (0 0 1 0 0 0 0) ie only the attribute $P_3$ is in the on state and all other vectors are in the off state the effect of Y on the dynamical system $N(S_2)$ is given by

$$YN(S_2)^T \quad \hookrightarrow \quad (0\ 0\ 0\ 0\ 0\ 0\ 0\ 1\ 0\ 0) \quad = \quad X$$

$$YN(S_2) \quad \hookrightarrow \quad (0\ 0\ 1\ 0\ 0\ 0\ 0).$$

Thus the fixed point is the binary pair {(0 0 1 0 0 0 0), (0 0 0 0 0 0 0 1 0 0)} has only influence on $M_8$ and according to this expert no connection with other nodes.

Now we proceed on to get the third experts opinion using the same set of attributes $\{M_1, M_2, \ldots, M_{10}\}$ and $\{P_1, P_2, \ldots, P_7\}$.

The directed graph related to the third experts opinion is

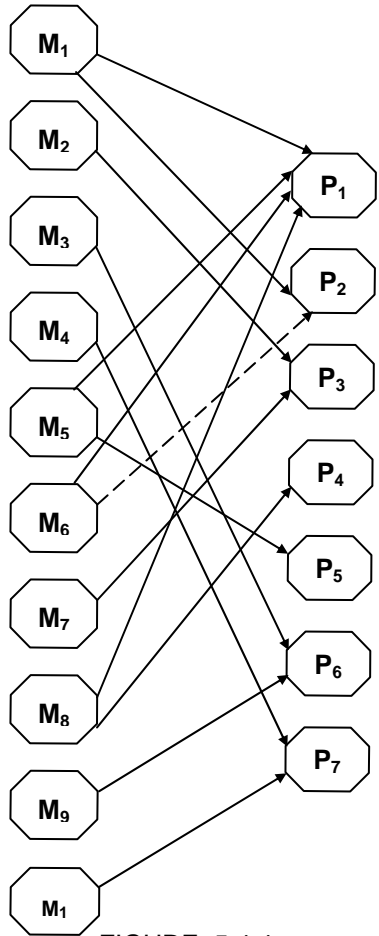

FIGURE: 5.4.4

We give the relational matrix associated with the neutrosophic directed graph and denote it by $N(S_3)$

$$
\begin{array}{c}
\phantom{M_{10}} \begin{array}{ccccccc} P_1 & P_2 & P_3 & P_4 & P_5 & P_6 & P_7 \end{array} \\
\begin{array}{c}
M_1 \\ M_2 \\ M_3 \\ M_4 \\ M_5 \\ M_6 \\ M_7 \\ M_8 \\ M_9 \\ M_{10}
\end{array}
\begin{bmatrix}
1 & 1 & 0 & 0 & 0 & 0 & 0 \\
0 & 0 & 1 & 0 & 0 & 0 & 0 \\
0 & 0 & 0 & 0 & 0 & 1 & 0 \\
0 & 0 & 0 & 0 & 0 & 0 & 1 \\
1 & 0 & 0 & 0 & 1 & 0 & 0 \\
1 & I & 0 & 0 & 0 & 0 & 0 \\
0 & 0 & 1 & 0 & 0 & 0 & 0 \\
1 & 0 & 0 & 1 & 0 & 0 & 0 \\
0 & 0 & 0 & 0 & 0 & 1 & 0 \\
0 & 0 & 0 & 0 & 0 & 0 & 1
\end{bmatrix}
\end{array}
$$

Let is consider the state vector X = (0 0 0 0 1 0 0 0 0 0) only the node $M_5$ is in the on state and all other nodes are in the off state. The effect of X on $N(S_3)$

$$XN(S_3) \quad \hookrightarrow \quad (1\ 0\ 0\ 0\ 1\ 0\ 0) \qquad = \quad Y$$

$$YN(S_3)^T \quad \hookrightarrow \quad (1\ 0\ 0\ 0\ 1\ 1\ 0\ 1\ 0\ 0) = \quad X_1$$

$$X_1\, N(S_3) \quad \hookrightarrow \quad (1\ 0\ 0\ 1\ 1\ \ 0\ 0) \qquad = \quad Y_1$$

$$Y_1\, N(S_3) \quad \hookrightarrow \quad (1\ 0\ 0\ 0\ 1\ 1\ 1\ 0\ 1\ 0\ 0) = \quad X_2 = X_1.$$

Thus the hidden pattern of the dynamical system in a fixed point given by the binary pair {(1 0 0 1 1 0 0), (1 0 0 0 1 1 0 1 0 0)} which makes on the nodes $P_1$, $P_4$ and $P_5$ in the range space.

Now we proceed on to consider the state vector Y = (0 0 1 0 0 0 0) i.e., only the node $P_3$ is in the on state and all other nodes are in the off state. The effect of Y on he dynamical system $N(S_3)$ is given by

$$YN(S_3)^T \quad \hookrightarrow \quad (0\ 1\ 0\ 0\ 0\ 0\ 1\ 0\ 0\ 0) = \quad X$$

$$XN(S_3) \quad \hookrightarrow \quad (0\ 0\ 1\ 0\ 0\ 0\ 0).$$

Thus the fixed point of the state vector Y is given by the binary pair {(0 0 1 0 0 0 0), (0 1 0 0 0 0 1 0 0 0)}, the node P3 has least effect on other nodes.

Next we consider the opinion of the forth expert for the same set of attributes $\{(M_1, M_2,\ldots, M_{10}\}$ $(P_1, P_2,\ldots, P_7\}$. The directed graph related to the forth experts opinion is as follows:

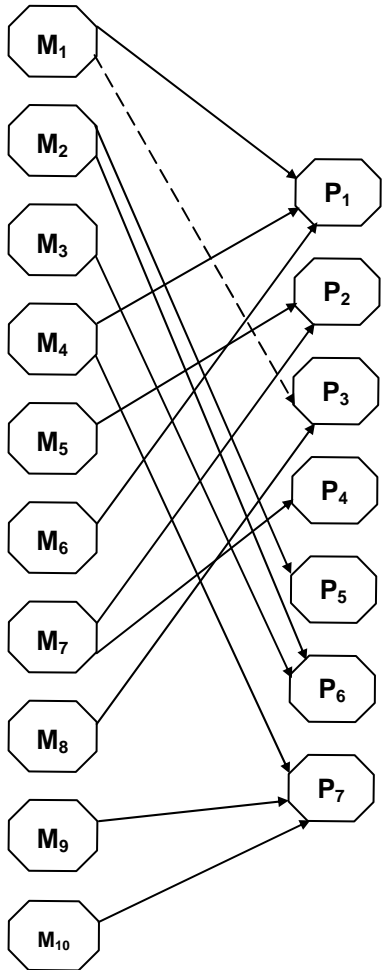

FIGURE: 5.4.5

The related neutrosophic $10 \times 7$ matrix $N(S_4)$ is given in the following page

$$\begin{array}{c}\begin{array}{ccccccc} P_1 & P_2 & P_3 & P_4 & P_5 & P_6 & P_7 \end{array}\\ \begin{array}{c} M_1 \\ M_2 \\ M_3 \\ M_4 \\ M_5 \\ M_6 \\ M_7 \\ M_8 \\ M_9 \\ M_{10} \end{array} \begin{bmatrix} 1 & 0 & I & 0 & 0 & 0 & 0 \\ 0 & 0 & 0 & 0 & 1 & 1 & 0 \\ 0 & 0 & 1 & 0 & 0 & 0 & 0 \\ 1 & 0 & 0 & 0 & 0 & 0 & 0 \\ 0 & 1 & 0 & 0 & 0 & 0 & 0 \\ 1 & 0 & 0 & 0 & 0 & 0 & 0 \\ 0 & 1 & 0 & 1 & 0 & 0 & 0 \\ 0 & 0 & 1 & 0 & 0 & 0 & 0 \\ 0 & 0 & 0 & 0 & 0 & 0 & 1 \\ 0 & 0 & 0 & 0 & 0 & 0 & 1 \end{bmatrix}\end{array}$$

Let us consider the state vector X = (0 0 0 0 1 0 0 0 0) in which only the node $M_5$ is in the on state and all other nodes are in the off state. The effect of X on the dynamical system $N(S_4)$.

$$\begin{array}{llll} XN(S_4) & \hookrightarrow & (0\ 1\ 0\ 0\ 0\ 0\ 0) & = & Y \\ YN(S_4)^T & \hookrightarrow & (0\ 0\ 0\ 0\ 1\ 0\ 1\ 0\ 0\ 0) & = & X_1 \\ X_1 N(S_4) & \hookrightarrow & (0\ 1\ 0\ 1\ 0\ 0\ 0) & = & Y_1 \\ Y_1 N(S_4)^T & \hookrightarrow & (0\ 0\ 0\ 0\ 1\ 0\ 1\ 0\ 0\ 0) & = & X_2 = X_1. \end{array}$$

Thus the fixed point of the NRM is a binary pair {(0 1 0 1 0 0 0), (0 0 0 0 1 0 1 0 0 0)}, which shows $M_5$ has main effect only $M_7$ and other nodes remain unaffected.

Next we consider the on state of the node $P_3$ alone and all the othr nodes are in the off state. Let Y = (0 0 1 0 0 0 0) be the state vector Effect of Y on the dynamical system $N(S_4)$ is given by

$$\begin{array}{llll} YN(s_4)^T & \hookrightarrow & (I\ 0\ 1\ 0\ 0\ 0\ 0\ 1\ 0\ 0) & = & X \\ XN(S_4) & \hookrightarrow & (I\ 0\ 1\ 0\ 0\ 0\ 0) & = & Y_1 \\ Y_1 N(S_4)^T & \hookrightarrow & (I\ 0\ 1\ I\ 0\ I\ 0\ 1\ 0\ 0) & = & X_1 \\ X_1 N(S_4)^T & \hookrightarrow & (I\ 0\ 1\ 0\ 0\ 0\ 0) & = & Y_2 = Y_1. \end{array}$$

Thus the binary pair associated with the fixed point in {($I$ 0 1 0 0 0 0), ($I$ 0 1 $I$ 0 $I$ 0 1 0 0)}. $P_3$ makes $P_1$ alone as an indeterminate in the range space. Thus according to this expert $P_3$ has least effect. Also $M_1$, $M_4$ and $M_6$ become indeterminate in the domain space and nodes $M_3$ and $M_8$ come to on state.

Now we consider the opinion of the fifth expert on the same set of attributes $\{(M_1, M_2, ..., M_{10}) (P_1 P_2, ..., P_7)\}$.

We give the associated neutrosophic directed graph given by the expert.

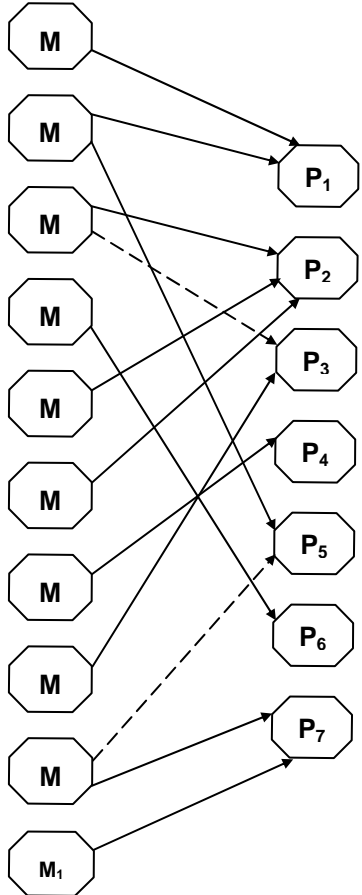

FIGURE: 5.4.6

The related connection matrix of the NRM. Let us denote this matrix by $N(S_5)$.

$$
\begin{array}{c}
\quad\quad P_1\; P_2\; P_3\; P_4\;\; P_5\; P_6\; P_7 \\
\begin{array}{c}
M_1 \\ M_2 \\ M_3 \\ M_4 \\ M_5 \\ M_6 \\ M_7 \\ M_8 \\ M_9 \\ M_{10}
\end{array}
\left[
\begin{array}{ccccccc}
1 & 0 & 0 & 0 & 0 & 0 & 0 \\
1 & 0 & 0 & 0 & 1 & 0 & 0 \\
0 & 1 & I & 0 & 0 & 0 & 0 \\
0 & 0 & 0 & 0 & 0 & 1 & 0 \\
0 & 1 & 0 & 0 & 0 & 0 & 0 \\
0 & 1 & 0 & 0 & 0 & 0 & 0 \\
0 & 0 & 0 & 1 & 0 & 0 & 0 \\
0 & 0 & 1 & 0 & 0 & 0 & 0 \\
0 & 0 & 0 & 0 & I & 0 & 1 \\
0 & 0 & 0 & 0 & 0 & 0 & 1
\end{array}
\right]
\end{array}
$$

Let us consider the state vector X = (0 0 0 0 1 0 0 0 0 0) ie the node $M_5$ alone in the on state and all other nodes are in the off state. The effect of X on the dynamical system $N(S_5)$.

$$
\begin{array}{llll}
X\,N\,(S_5) & \hookrightarrow & (0\;1\;0\;0\;0\;0\;0) & = \quad Y \\
Y\,N\,(S_5)^T & \hookrightarrow & (0\;0\;1\;0\;1\;1\;0\;0\;0\;0) & = \quad X_1 \\
X_1\,N(S_5) & \hookrightarrow & (0\;1\;I\;1\;0\;0\;0) & = \quad Y_1 \\
Y\,N\,(S_5)^T & \hookrightarrow & (0\;0\;0\;0\;1\;1\;1\;\;I\;0\;0) & = \quad X_2 \\
X_2\,N(S_5) & \hookrightarrow & (0\;1\;I\;1\;0\;0\;0) & = \quad Y_2 = Y_1.
\end{array}
$$

Thus the hidden pattern of the dynamical system is a fixed point given by the binary pair {(0 1 $I$ 1 0 0 0), (0 0 0 0 1 1 1 $I$ 0 0)}. Interpret the resultant vector!

Now let us consider the state vector Y = (0 0 1 0 0 0 0) where only the node $M_3$ is in the on state and all other nodes are in the off state. Effect of Y on the dynamical system $N(S_5)$

$$
\begin{array}{llll}
Y\,N\,(S_5)^T & \hookrightarrow & (0\;0\;I\;0\;0\;0\;0\;1\;0\;0) & = \quad X \\
X\,N\,(S_5) & \hookrightarrow & (0\;I\;1\;0\;0\;0\;0) & = \quad Y_1 \\
Y_1\,N(S_5) & \hookrightarrow & (0\;0\;0\;0\;I\;I\;0\;1\;0\;0) & = \quad X_1 \\
X_1\,N(S_5) & \hookrightarrow & (0\;I\;1\;0\;0\;0\;0) & = \quad Y_2 = Y_1.
\end{array}
$$

Thus the hidden pattern is a fixed point given by the binary pair {(0 $I$ 1 0 0 0 0), (0 0 0 0 $I$ $I$ 0 1 0 0 )}. This node has mainly indeterministic effect on the dynamical system except on the node $M_8$ that becomes on.

Now we form the combined NRM using the fire experts opinion. The related connection matrix of the combined NRM is given by $N(S) = N(S_1) + N(S_2) + N(S_3) + N(S_4) + N(S_5)$ ie sum of the 5 matrices.

$$
\begin{array}{c}
\\ M_1 \\ M_2 \\ M_3 \\ M_4 \\ M_5 \\ M_6 \\ M_7 \\ M_8 \\ M_9 \\ M_{10}
\end{array}
\begin{array}{c}
\begin{array}{ccccccc} P_1 & P_2 & P_3 & P_4 & P_5 & P_6 & P_7 \end{array} \\
\left[\begin{array}{ccccccc}
5 & 0 & I & 0 & 0 & 0 & 0 \\
1 & 0 & 0 & 0 & 4 & 2 & 0 \\
0 & 1 & 0 & 0 & 0 & 3 & 0 \\
2 & 0 & 0 & 0 & 0 & 1 & 2 \\
1 & 4 & 0 & 0 & 1 & 0 & 0 \\
3 & 0 & 1 & 0 & 0 & 0 & 0 \\
0 & 1 & 0 & 4 & 0 & 0 & 0 \\
1 & 0 & 3 & 1 & 0 & 0 & 1 \\
0 & 0 & 0 & 0 & I & 1 & 4 \\
0 & 0 & 0 & 0 & 0 & 0 & 5
\end{array}\right]
\end{array}
$$

Now we consider the effect of any state vector on the combined NRM $N(S)$. Suppose $X = (0\,0\,0\,0\ 1\,0\,0\,0\,0)$ ie only the node $M_5$ is in the on state and all other nodes are in the off state. The effect of X on the CNRM is given by

We threshold the CNRM as $a_I \geq 2$ we put 1 if $a_I < 2$ we put o if $a_I$ is $I + 1$, then $a_I = 0$

$XN(S) \quad \hookrightarrow \quad (0\,1\,0\,0\,0\,0\,0) \qquad = \quad Y$

$YN(S)^T \quad \hookrightarrow \quad (0\,0\,0\,0\,1\,0\,0\,0\,0)$.

Thus the fixed pt is $\{(0\,0\,0\,0\,1\,0\,0\,0\,0)\}$…

Let $Y = (0\,0\,1\,0\,0\,0)$ ie only the node $P_3$ is in the on state and all other nodes are in the off state effect of Y on $N(S)$

$YN(S)^T \quad \hookrightarrow \quad (I\,0\,0\,0\,0\,0\,0\,1\,0\,0) \quad = \quad X$

$XN(S)^T \quad \hookrightarrow \quad (I\,0\,1\,0\,0\,0\,0) \qquad = \quad Y_1$

$Y_1\,N(S)^T \quad \hookrightarrow \quad (I\,I\,0\,1\,I\,1\,0\,1\,0\,0) \quad = \quad X_1$

$X_1\,N(S) \quad \hookrightarrow \quad (I\,I\,1\,0\,I\,I\,1) \qquad = \quad Y_2$

$$Y_2 \, N(S)^T \quad \hookrightarrow \quad (I\,I\,I\,I\,I\,I\,0\,1\,1\,1) \quad = \quad X_2$$

$$X_2 \, N(S) \quad \hookrightarrow \quad (I\,I\,1\,1\,I\,1) \quad\quad = \quad Y_3$$

$$Y_3 \, N(S)^T \quad \hookrightarrow \quad (I\,I\,I\,0\,I\,I\,I\,0\,1\,1) \quad = \quad X_3$$

$$X_3 \, N(S)^T \quad \hookrightarrow \quad (I\,I\,1\,I\,I\,I\,1) \quad\quad = \quad Y_4.$$

Thus the resultant vector is a fixed point given by the binary pair $\{(I\ I\ I\ I\ I\ I\ 0\ 1\ 1\ 1), (I\ I\ 1\ I\ I\ I\ 1)\}$. Most of the nodes of the dynamical system in the CNRM are indeterminate.

## 5.5 Linked Neutrosophic Relational Maps and its application to Migrant Problems

By the term migrant problems we mean the socio economic problems faced by the HIV/AIDS affected migrant labourer. To study this problem we were forced to construct a new neutrosophic model called the linked neutrosophic relational maps. Neutrosophic relational maps cannot always be obtained given two classes of attributes, which are formed for a particular problem. The inter relation between three or more classes at times becomes so implicit and one finds it difficult to obtain the neutrosophic bigraph by using them. In such cases we can make use of the linked neutrosophic relational maps.

**DEFINITION 5.5.1:** *Let P be some problem, which is under investigation. Let us have 3 sets of disjoint attributes related with P, i.e., say $(C_1, C_2, \ldots, C_n)$, $(B_1, B_2, \ldots, B_5)$ and $(A_1, A_2, \ldots, A_m)$ be the 3 classes of n, t and m attributes respectively. Suppose the expert gives the two neutrosophic digraphs relating the classes of attributes $\{(C_1, \ldots, C_n)$ and $(B_1, \ldots, B_t)\}$ and $\{(B_1, \ldots, B_t)$ and $(A_1, \ldots, A_m)\}$ and the classes of attributes $\{(C_1, \ldots, C_n)$ and $(A_1, \ldots, A_m)\}$ are related implicitely and not able to give the relation, then form the neutrosophic relational matrices R and S related to the two bigraph (say) R is a n × t matrix and S a t × m matrix. Then the product of the matrix R × S gives a n × m matrix which gives the relational matrix between $\{(C_1, C_2, \ldots, C_n)$ and $(A_1, \ldots, A_m)\}$. We call this matrix the linked relational matrix and the neutrosophic graph drawn using this linked relational matrix, called the linked neutrosophic directed graph.*

Now suppose we have for any problem P some r set of disjoint classes of attributes are available and each class has say

$n_1$, $n_2$,…, $n_p$ elements in them; we can always obtain the relation between desired ones by linking them. Thus we can get linked matrices and these linked matrices will give the linked relational neutrosophic graphs.

We illustrate this by the problem of HIV/AIDS affected migrant labourers and their socio economic conditions. Let us consider the three classes of attributes relating the socio-economic and psychological aspects of migrant labourers with reference to HIV/AIDS given under the three broad heads

A – Causes for migrants vulnerability to HIV/AIDS
F – Factors forcing migration
G – Role of government

given in chapter 3 page 105-106 of this book. The subtitles under these three heads are are also given in the same pages.

Suppose the classes (A and G) and (F and G) related by the neutrosophic relational maps and the classes (A and F) could not be related explicitly by the expert. We use the method of linked neutrosophic relational maps and obtain the linked neutrosophic bigraph using the linked relational matrices.

The neutrosophic directed graph using the attributes A and G.

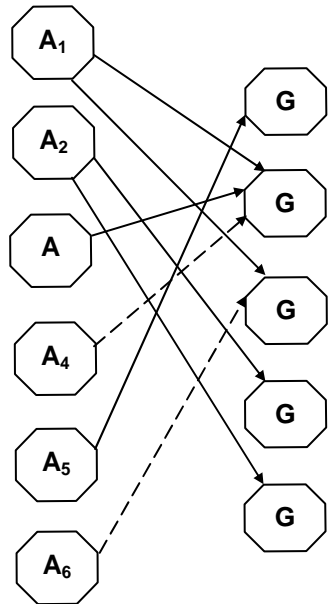

FIGURE: 5.5 .1

According to this expert when they have no education / no awareness it implies the awareness class in rural area about HIV/AIDS is not effective.

It also implies the health center if they are in the village does not propagate the awareness about HIV/AIDS or any other related program about HIV/AIDS. However the expert says relation between migrant labourers who are victims of bad company and bad habit and the awareness clubs in rural areas about HIV/AIDS is an indeterminate. He says only a socially irresponsible man will be unaware of HIV/AIDS. Also according to this expert the relation between cheap availability of CSWs and construction of hospitals to spread HIV/AIDS awareness and related program is an indeterminate. However he inter-relates. Types of profession with alternative job if harvest fails, there by stopping migration.

We give the related neutrosophic connection matrix N(L)

$$
\begin{array}{c}
\quad\; G_1\; G_2\;\; G_3\;\; G_4\; G_5 \\
\begin{array}{c}
A_1 \\ A_2 \\ A_3 \\ A_4 \\ A_5 \\ A_6
\end{array}
\left[
\begin{array}{ccccc}
0 & 1 & 1 & 0 & 0 \\
0 & 0 & 0 & 1 & 1 \\
0 & 1 & 0 & 0 & 0 \\
0 & I & 0 & 0 & 0 \\
1 & 0 & 0 & 0 & 0 \\
0 & 0 & I & 0 & 0
\end{array}
\right]
\end{array}
$$

Suppose the node $G_2$ alone is in the on state and all other nodes in the off state let x = (0 1 0 0 0). The effect of X on N(L) on the dynamical system is given by

$$
\begin{aligned}
XN(L)^T &\hookrightarrow (1\ 0\ 1\ I\ 0\ 0) &=& \quad Y \\
YN(L) &\hookrightarrow (0\ 1\ 1\ 0\ 0) &=& \quad X_1 \\
Y_1\,N(L)^T &\hookrightarrow (1\ 0\ 1\ I\ 0\ I) &=& \quad Y_1 \\
Y_1\,N(L) &\hookrightarrow (0\ 1\ I\ 0\ 0) &=& \quad X_2 \\
X_2\,N(L)^T &\hookrightarrow (I\ 0\ 1\ I\ 0\ 1) &=& \quad Y_2 = Y_1.
\end{aligned}
$$

Thus the hidden pattern of the neutrosophic dynamical system is a fixed point given by the binary pair {(1 0 1 $I$ 0 I) (0 1 $I$ 0 0)}. The reader is expected to anlyze the fixed point.

Suppose we consider the node $A_3$ to be in the in the on state and all other nodes are in the off state. Denote the state vector by X = (0 0 1 0 0 0)

$$XN (L) \quad \hookrightarrow \quad (0\ 1\ 0\ 0\ 0) \quad = \quad Y$$

$$YN(L)^T \quad \hookrightarrow \quad (1\ 0\ 1\ I\ 0\ 0) \quad = \quad X_1$$

$$X_1\ N(L) \quad \hookrightarrow \quad (0\ 1\ 1\ 0\ 0) \quad = \quad Y_1$$

$$Y_1\ N(L)^T \quad \hookrightarrow \quad (1\ 0\ 1\ I\ 0\ I) \quad = \quad X_2$$

$$X_2\ N(L) \quad \hookrightarrow \quad (0\ 1\ 1\ 0\ 0) \quad = \quad Y_2 = Y_1.$$

So the hidden pattern of the dynamical system is a fixed point given by the pair {(0 1 1 0 0), (1 0 1 $I$ 0 I)}, which makes A1 and A6 on and A4 is an indeterminate. G2 and G3 become on in the range space.

It is left as an exercise for the reader to obtain a C-program to determine the resultant vector given the state vector and the relational neutrosophic matrix.

Now using the same expert we obtain the neutrosophic bigraph relating the attributes {$F_1$ $F_2$ $F_3$ $F_4$} and {$G_1$ $G_2$ $G_3$ $G_4$ $G_5$}.

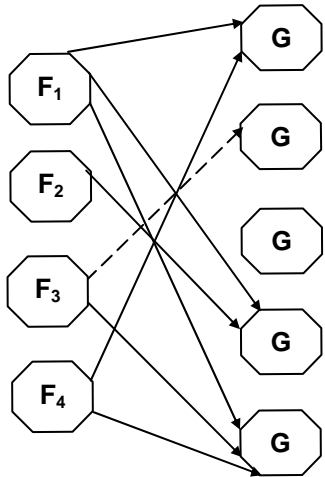

FIGURE: 5.5.2

The relational connection neutrosophic matrix.

$$
\begin{array}{c}
\quad\ G_1\ G_2\ G_3\ G_4\ G_5 \\
\begin{array}{c} F_1 \\ F_2 \\ F_3 \\ F_4 \end{array}
\begin{bmatrix}
1 & 0 & 0 & 1 & 1 \\
0 & 0 & 0 & 1 & 0 \\
0 & I & 0 & 0 & 1 \\
1 & 0 & 0 & 0 & 1
\end{bmatrix}
\end{array}
$$

Let us denote this neutrosophic relational matrix by $N(L_1)$. Suppose consider the state vector $X = (0\ 0\ 1\ 0\ 0)$. Thus we see the government construction of hospitals in rural areas with HIV/AIDS counseling cell / compulsory HIV/AIDS test before marriage or any other program has no effect on any of the factors forcing people for migration which is very true for when people suffer with unemployment.

Hence live under poverty and a third party man takes advantage of the situation and mobilizes them as contract labour and when the yield is so poor due to inferlity of land because of modern techniques.

Next we consider the state vector $Y = (0\ 1\ 0\ 0)$ ie only the attribute $F_2$ is in the on state and all other nodes are in the off state. The effect of Y on the dynamical system $N(L_1)$

$$
\begin{array}{llllll}
Y N(L_1) & \hookrightarrow & (0\ 0\ 0\ 1\ 0) & = & X \\
X N(L_1)^T & \hookrightarrow & (1\ 1\ 0\ 0) & = & Y_1 \\
Y_1 N(L_1) & \hookrightarrow & (1\ 0\ 0\ 1\ 1) & = & X_1 \\
X_1 N(L_1)^T & \hookrightarrow & (1\ 1\ 1\ 1) & = & Y_2 \\
Y_2 N(L_1) & \hookrightarrow & (1\ I\ 0\ 11) & = & X_2 \\
X_2 N(L_1)^T & \hookrightarrow & (!\ 1\ I\ 1) & = & Y_3 \\
Y_3 N(L_1) & \hookrightarrow & (1\ I\ 0\ 11) & = & X_3 & = & X_2.
\end{array}
$$

Thus the hidden pattern of the neutrosophic system is a fixed point given by the binary pair $((1\ I\ 0\ 1\ 1), (1\ 1\ I\ 1)\}$ all state become on and indeterminate thus the node $F_2$ has a large influence on the dynamical system.

Now using the neutrosophic relational matrices $N(L)$ and $N(L_1)$ we can obtain the linked connection matrix relating the class of attributes $\{ A_1, A_2, A_3, A_4, A_5, A_6\}$ and $\{F_1\ F_2\ F_3\ F_4\}$.

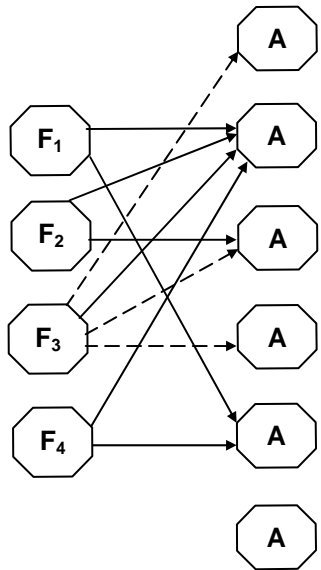

FIGURE: 5.5.3

The product of the matrix $N(L)$ $N(L_1)^T \sim N(P)$ $N(P)$ is obtained after thresholding the entries by the rule $a_{ij} \geq 1$ is replaced by 0.

$$N(P) = \begin{array}{c} \\ A_1 \\ A_2 \\ A_3 \\ A_4 \\ A_5 \\ A_6 \end{array} \begin{array}{cccc} F_1 & F_2 & F_3 & F_4 \\ \left[ \begin{array}{cccc} 0 & 0 & I & 0 \\ 1 & 1 & 1 & 1 \\ 0 & 0 & I & 0 \\ 0 & 0 & I & 0 \\ 1 & 0 & 0 & 1 \\ 0 & 0 & 0 & 0 \end{array} \right] \end{array}$$

From matrix $N(P)$ we can obtain the linked relational map relating $\{A_1, \ldots, A_6\}$ and $\{F_1\ F_2\ F_3\ F_4\}$. However the attribute $A_6$ remain unrelated seen from figure 5.5.3. Now we can study the effect of state vector on this linked dynamical system. Now we study the effect of the state vector say $X = (0\ 0\ 0\ 1\ 0\ 0)$ ie only the attribute $A_4$ is in the on state and all other attributes are in the off state.

$$XN\,(P) \quad \hookrightarrow \quad (0\ 0\ I\ 0\ ) \quad = \quad Y$$

$$YN(P)^T \quad \hookrightarrow \quad (I\ I\ I\ I\ 0\ 0\ ) \quad = \quad X_1$$

$$X_1\,N(P) \quad \hookrightarrow \quad (I\ I\ I\ I) \quad = \quad Y_1$$

$$Y_1\,N\,(P) \quad \hookrightarrow \quad (I\ I\ I\ 1\ I\ 0) \quad = \quad X_2$$

$$X_2\,N(P) \quad \hookrightarrow \quad (I\ I\ I\ I).$$

that is when the attribute $A_a$ is in the on state one is not in a position to state whether the states of the range space that is the attributes $A_1$, $A_2$ $A_3$ and $A_5$ is being influenced or not so they reach a state which is indeterminate i.e., the indeterminate state only $A_6$ is in the off state. Thus we see only using the neutrosophic model we are in a position to say the models influence on certain state vectors leave even indeterminate state that is one cannot say the state is off or the state is on. Thus this provides a very new insight in the study and analysis of the problem.

Next we consider a state vector say $Y = (0\ 1\ 0\ 0)$ i.e., only the node $F_2$ is in the on state and all other vectors are in the off state. The effect of $Y$ on the dynamical system $N(P)$ is given by

$$YN\,(P)^T \quad \hookrightarrow \quad (0\ 1\ 0\ 0\ 0\ 0) \quad = \quad X$$

$$XN\,(P) \quad \hookrightarrow \quad (1\ 1\ 1\ 1) \quad = \quad Y_1$$

$$Y_1\,N(P) \quad \hookrightarrow \quad (I\ 1\ I\ I\ I\ 1\ 0) \quad = \quad X_1$$

$$X_1\,N(P) \quad \hookrightarrow \quad (1\ 1\ I\ 1) \quad = \quad Y_2$$

$$Y_2N(P) \quad \hookrightarrow \quad (I\ 1\ I\ I\ 1\ 0) \quad = \quad X_2\,(= X_1).$$

Thus the hidden pattern is a fixed point given by the binary pair $\{(1\ 1\ I\ 1),\ (I\ 1\ I\ I\ 1\ 0)\}$. $F_2$ has a strong influence on the dynamical system.

Now we proceed on to study some other model using the linked neutrosophic NRM. We analyze for the model given in page 239-240 of this chapter Let $C = \{C_1\ C_2\ C_3\ C_4\ C_5\}$, $A = \{A_1, A_2,\ldots, A_7\}$ and $T = \{T_1, T_2,\ldots, T_6\}$.

Now using the experts opinion we obtain the directed graph of the NRM relating the attributes $\{A_1, A_2, \ldots, A_7\}$ and $\{T_1, T_2, \ldots, T_6\}$ which is given by the following diagram:

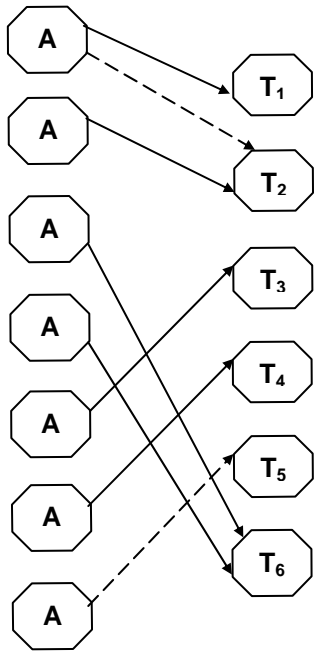

FIGURE: 5.5.4

Now using the neutrosophic directed graph we give the related connection neutrosophic matrix N(L).

$$
\begin{array}{c c}
& \begin{array}{c c c c c c} T_1 & T_2 & T_3 & T_4 & T_5 & T_6 \end{array} \\
\begin{array}{c} A_1 \\ A_2 \\ A_3 \\ A_4 \\ A_5 \\ A_6 \\ A_7 \end{array} &
\left[ \begin{array}{c c c c c c}
1 & I & 0 & 0 & 0 & 0 \\
0 & 1 & 0 & 0 & 0 & 0 \\
0 & 0 & 0 & 0 & 0 & 1 \\
0 & 0 & 0 & 0 & 0 & 1 \\
0 & 0 & 1 & 0 & 0 & 0 \\
0 & 0 & 0 & 1 & 0 & 0 \\
0 & 0 & 0 & 0 & I & 0
\end{array} \right]
\end{array}
$$

Using N(L) we obtain the effect of any given state vector. Let $X = (0\ 0\ 0\ 1\ 0\ 0\ 0)$ be the state vector in which only the node $A_4$ is in the on state and all other nodes are in the off state. Effect of X on the dynamical system N(L) is given by

$$
\begin{array}{llll}
XN(L) & \hookrightarrow & (0\ 0\ 0\ 0\ 0\ 1) & = \quad Y \\
YN(L)^T & \hookrightarrow & (0\ 0\ 1\ 1\ 0\ 0\ 0) & = \quad X_1
\end{array}
$$

$$X_1 \, N(L) \quad \hookrightarrow \quad (0\ 0\ 0\ 0\ 0\ 1) \quad = \quad Y_1 = Y_1.$$

Thus the hidden pattern is a fixed point given by the binary pair $\{(0\ 0\ 0\ 0\ 0\ 1)\ (0\ 0\ 1\ 1\ 0\ 0\ 0)\}$, which has a very low influence on the system.

Next consider the state vector $X = (0\ 0\ 0\ 0\ 1\ 0\ 0)$ where only the node which is on is $A_5$. The effect of X on the dynamical system N(L) is given by

$$XN(L) \quad \hookrightarrow \quad (0\ 0\ 1\ 0\ 0\ 0) \quad = \quad Y$$
$$YN(L)^T \quad \hookrightarrow \quad (0\ 0\ 0\ 0\ 1\ 0\ 0) \quad = \quad X_1 = X.$$

Thus the fixed point of the dynamical system is a binary pair given by $\{(0\ 0\ 1\ 0\ 0\ 0),\ (0\ 0\ 0\ 0\ 1\ 0\ 0)\}$, according to the expert the effect on one node on the other node is very poor. Let us now consider the effect of the state vector $Y = (0\ 0\ 0\ 1\ 0\ 0)$ ie only the node $T_4$ is in the on state and all other nodes are in the off state. Effect of Y on the dynamical system N(L) is given by

$$YN(L)^T \quad \hookrightarrow \quad (0\ 0\ 0\ 0\ 0\ 1\ 0\ 0) \quad = \quad X$$
$$XN\ (L) \quad \hookrightarrow \quad (0\ 0\ 0\ 1\ 0\ 0) \quad = \quad Y_1\ (= Y).$$

Thus the binary pair is $\{(0\ 0\ 0\ 1\ 0\ 0),\ (0\ 0\ 0\ 0\ 0\ 1\ 0\ 0)\}$. Same conclusion as above. Now we give the directed graph of the NRM relating the attributes $C = \{C_1\ C_2\ C_3\ C_4\ C_5\}$ and $A = \{A_1\ A_2,\ldots, A_6\}$ given by the expert.

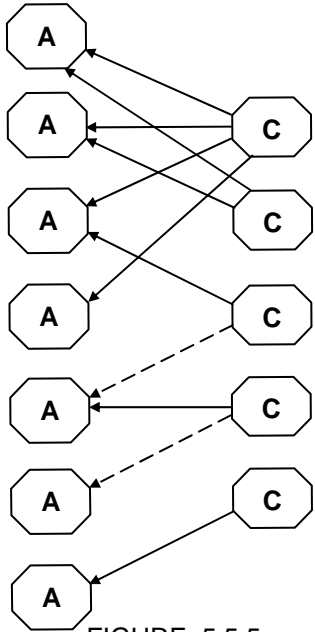

FIGURE: 5.5.5

Using the directed neutrosophic graph we now obtain the related neutrosophic connection matrix N(L).

$$
\begin{array}{c}
\quad\quad A_1 \ A_2 \ A_3 \ A_4 \ \ A_5 \ A_6 \ A_7 \\
\begin{array}{c} C_1 \\ C_2 \\ C_3 \\ C_4 \\ C_5 \end{array}
\begin{bmatrix}
1 & 1 & 1 & 1 & 0 & 0 & 0 \\
1 & 1 & 0 & 0 & 0 & 0 & 0 \\
0 & 0 & 1 & 0 & I & 0 & 0 \\
0 & 0 & 0 & 0 & 1 & I & 0 \\
0 & 0 & 0 & 0 & 0 & 0 & 1
\end{bmatrix}
\end{array}
$$

Let us consider the state vector $X = (0\ 0\ 1\ 0\ 0)$ ie only the attribute $C_3$ is in the on state and all other vectors are in the off state. The effect of X on the dynamical sytem N(4) is given by

$$
\begin{array}{llll}
X N (L_1) & \hookrightarrow & (0\ 0\ 1\ 0\ I\ 0\ 0) & = \quad Y \\
Y N (L_1) & \hookrightarrow & (1\ 0\ 1\ I\ 0) & = \quad X_1 \\
X_1\ N(L_1) & \hookrightarrow & (1\ 1\ 1\ 1\ I\ I\ 0) & = \quad Y_1 \\
Y_1\ N(L_1) & \hookrightarrow & (1\ 1\ 1\ I\ 0) & = \quad X_2 \\
X_2\ N(L_1) & \hookrightarrow & (1\ 1\ 1\ 1\ I\ I\ 0) & = \quad Y_2 = (Y_1).
\end{array}
$$

Thus the hidden pattern of the dynamical system is a fixed point given by the binary pair $\{(1\ 1\ 1\ I\ 0),\ (1\ 1\ 1\ 1\ I\ I\ 0)\}$. The effect of the node $C_3$ is strong for it makes other nodes on or indeterminate.

Now we proceed on to study the effect of state vector $X = (1\ 0\ 0\ 0\ 0)$ where only the attribute $C_1$ is in the on state and all other vectors are in the off state. The effect of X on the dynamical system is given by

$$
\begin{array}{llll}
X N (L_1) & \hookrightarrow & (1\ 1\ 1\ 1\ 0\ 0\ 0) & = \quad Y \\
Y N (L_1)^T & \hookrightarrow & (1\ 1\ 1\ 0\ 0) & = \quad X_1 \\
X_1\ N(L_1) & \hookrightarrow & (1\ 1\ 1\ 1\ I\ 0\ 0) & = \quad Y_1 \\
Y_1 N(L_1)^T & \hookrightarrow & (1\ 1\ I\ I\ 0) & = \quad X_2 \\
X_2\ N(L_1) & \hookrightarrow & (1\ 1\ I\ 1\ I\ I\ 0) & = \quad Y_2 \\
Y_2\ N(L_2)^T & \hookrightarrow & (1\ 1\ I\ I\ 0) & = \quad X_3 = \quad X_1.
\end{array}
$$

Thus the hidden pattern of the dynamical system is the fixed point given by the binary pair $\{(1\ 1\ I\ I\ 0),\ (1\ 1\ I\ 1\ I\ I\ 0)\}$. The influence of the node $C_1$ also on the dynamical system is very strong.

Suppose consider the state vector $Y = (0\ 0\ 1\ 0\ 0\ 0\ 0)$ ie only the attribute $A_3$ alone is in the on state and all other attributes are in the off state. The effect of Y on the dynamical system $N(L_1)$ is given by

$$
\begin{array}{llll}
Y N (L_1)^T & \hookrightarrow & (1\ 0\ 1\ 0\ 0) & = & X \\
X N (L_1) & \hookrightarrow & (1\ 1\ 1\ 1\ I\ 0\ 0) & = & Y_1 \\
Y_1 N(L_1)^T & \hookrightarrow & (1\ 1\ I\ I\ 0) & = & X_1 \\
X_1 N (L_1) & \hookrightarrow & (1\ 1\ 1\ 1\ I\ I\ 0) & = & Y_2 \\
Y_2 N(L_1)^T & \hookrightarrow & (1\ 1\ I\ I\ 0). &
\end{array}
$$

Thus the hidden pattern of the dynamical system is a fixed point given by the binary pair $\{(1\ 1\ I\ I\ 0),\ (1\ 1\ 1\ 1\ I\ I\ 0)\}$. The effect of $A_3$ on the system is more or less same as the effect of the node $C_1$ on the system.

The expert feels that he is not in a position to interrelate the attributes (Cand T so now we use the method of linked NRM to find the link between $(C_1,\ldots,\ C_5)$ and $(T_1,\ldots,\ T_6)$ by laking the product of the matrices $N(L_1) \times N(L)$ which gives a $5 \times 6$ matrix relating C and T. Let $N(L_1) \times N(L) \sim N(L_2)$ whre '~' denote the entries in the product of the matrices $N(L_1) \times N(L)$ is replaced by 1 if the entries $a_{ij} \geq 1$ and by 0 if $a_{ij} \leq 0$.

|       | $T_1$ | $T_2$ | $T_3$ | $T_4$ | $T_5$ | $T_6$ |
|-------|-------|-------|-------|-------|-------|-------|
| $C_1$ | 1 | $I$ | 0 | 0 | 0 | 1 |
| $C_2$ | 0 | $I$ | 0 | 0 | 0 | 0 |
| $C_3$ | 0 | 0 | $I$ | 0 | 0 | 1 |
| $C_4$ | 0 | 0 | 1 | $I$ | 0 | 0 |
| $C_5$ | 0 | 0 | 0 | 0 | $I$ | 0 |

Now we can obtain the relational neutrosophic directed graph which will link the attributes $\{C_1,\ldots,\ C_5\}$ and $\{T_1,\ldots,\ T_6\}$. We called this neutrosophic graph as the linked directed neutrosophic graph for it is not the directed graph given by the expert but only implicitly, the relations are calculated using the aid of the linked matrix.

Let as now study the effect of the state vector $X = (0\ 0\ 1\ 0\ 0\ 0)$ i.e., only the attribute $T_3$ is in the on state and all other nodes are in the off state. The effect of X on the linked dynamical system $N(L_2)$ is given by

$$\begin{array}{lllll}
XN\,(L_2)^T & \hookrightarrow & (0\ 0\ I\ 0\ 1) & = & Y \\
YN\,(L_2) & \hookrightarrow & (0\ 0\ 1\ I\ 0\ I) & = & X_1 \\
X_1\,N(L_2)^T & \hookrightarrow & (I\ 0\ I\ I\ 0) & = & Y_1 \\
Y_1\,N(L_2) & \hookrightarrow & (I\ I\ 1\ I\ 0\ I) & = & X_2\ \text{say} \\
X_2\,N(L_2)^T & \hookrightarrow & (I\ I\ I\ I\ 0) & = & Y_2 \\
Y_2\,N\,(L_2) & \hookrightarrow & (I\ I\ 1\ I\ 0\ I) & = & X_3 = \quad X_2.
\end{array}$$

Thus the hidden pattern of the dynamical system is given by the fixed point, which is given by the binary pair, $\{(I\ I\ I\ I\ 0),\ (I\ I\ 1\ I\ 0\ I)\}$ i.e., several nodes in the domain and range space become indeterminate with the influence of $T_3$.

Now let us consider the state vector $Y = (1\ 0\ 0\ 0\ 0)$ ie only the attribute $C_1$ is in the on state and all other vectors are in the off state the effect of $Y$ on the dynamical system $N(L_2)$ is given by

$$\begin{array}{lllll}
YN(L_2) & \hookrightarrow & (1\ I\ 0\ 0\ 0\ 01) & = & X \\
XN\,(L_2) & \hookrightarrow & (1\ I\ 1\ 0\ 0) & = & Y_1 \\
Y_1\,N\,(L_2) & \hookrightarrow & (1\ I\ I\ 0\ 0\ 1) & = & X_1 \\
X_1\,N\,(L_2) & \hookrightarrow & (1\ I\ I\ I\ 0) & = & Y_2 \\
Y_2\,N\,(L_2) & \hookrightarrow & (1\ I\ I\ I\ 0\ I) & = & X_2 \\
X_2\,N(L_2)^T & \hookrightarrow & (1\ I\ I\ I\ 0) & = & Y_3 = Y_2.
\end{array}$$

Thus we get the hidden pattern of the dynamical system to be the fixed point given by the binary pair $\{(1\ I\ I\ I\ 0),\ (1\ I\ I\ I\ 0\ I)\}$. The node C1 has strong impact on the dynamical system $N(L_2)$.

## 5.6 Combined Disjoint Block NRM and Combined Overlap Block NRM

When the problem under investigation has several attributes it becomes very difficult to handle them as a whole and all the more it becomes difficult to get the expert opinion and the related neutrosophic directed graph so we face the problem of the management of all the attributes together so in this chapter we give new methods by which the problem can be managed section by section or part by part without sacrificing any thing. We define here four methods.

**DEFINITION 5.6.1:** *Consider a NRM model with n elements say {D₁, D₂,…, Dₙ} in the domain space and {R₁, R₂,…, Rₘ} elements in the range space. Suppose we imagine that both m and n are fairly large numbers. Now we divide the n elements (ie attributes) of the domain space into t blocks each (assuming n is not a prime) such that no element finds its place in two blocks and all the elements are necessarily present in one and only one block such that each block has the same number of elements.*

*Similarly we divide the m attributes of the range space in r blocks so that no element is found in more than one block ie every element is in one and only one block. Thus each block has equal number of elements in them (As m is also a assurance to be a large number which is not a prime. Now if we take one block from the t block and one block from the r-block and form the neutrosophic directed graph using an expert's opinion. It is important to see no block repeats itself either from the r blocks or from the t blocks i.e., each block from the rage space as well as the domain space occurs only once.*
*Now we take the directed graph of each pair of blocks (one block taken from the domain space and one block taken from the range space) and their related neutrosophic connection matrix. We use all these connection matrices and form the n × m neutrosophic matrix, which is called the connection matrix of the Combined Disjoint Block of the NRM of equal sizes.*

*Now instead of dividing them into equal sizes blocks we can also divide them into different size blocks so that still they continue to be disjoint.*
*Thus if we use equal size block we mention so other wise we say the Combined Disjoint Block of the NRM of varying sizes.*

Now we give the definition of combined over lap blocks NRM.

**DEFINITION 5.6.2:** *Let us consider a NRM where the attributes related with the domain space say D = {D₁,…, Dₙ} is such that n is very large and the attributes connected with the range space say R = R₁,…, Rₘ} is such that m is very large. Now we want to get the neutrosophic directed graph, it is very difficult to manage the directed graph got using the expert opinion for the number of nodes is large equally, large is the number of edges also it may happen n and m may be primes and further it may be the case where some of the attributes may be common.*

*In such a case we adopt a new technique called the combined over lap block of equal length. Take the n attributes of the domain space divide them into blocks having same number of element also see that the number elements overlap between any two blocks is either empty the same number of element; carryout the similar procedure in case of the range space elements also.*

*Now pair the blocks by taking a block from the range elements and a block from the domain elements draw using the experts opinion the neutrosophic directed graph. See that each blocks in the pairs occur only once. Now using the neutrosophic directed graphs we obtain every pairs neutrosophic relational matrices using these matrices we form the $n \times m$ matrices which is the matrix related with the combined overlap block of the NRM of equal size if we very the number of elements in each of these blocks we call such a model as the combined overlap block of the NRM with varying sizes of blocks.*

It has become pertinent to mention here that we can by all means divide a NCM, if the attributes so formed into disjoint blocks of all them as combined disjoint block of the NCM. On similar lines we can divide NCM into overlap blocks of equal size and using this form the combined overlap block NCM. Already this concept for FCM have been used by Bart Kosko in his book [62]. Here we only define and study it in case of NCMs. As our main motivation in this book is to analyze the problem of HIV/AIDS affected migrant labour from rural areas of Tamil nadu we illustrate all these models only by HIV/AIDS affected migrant laboures and the problems faced by them. In several case we draw the conclusion based on our analysis which is given in chapter 7 of this book.

Now using the NRM model given in chapter V using the attributes of the domain space $\{D_1 \ D_2, \ldots, D_6\}$ and that of the range space $\{R_1, \ldots, R_{10}\}$.

We give he disjoint block decomposition of them and give the related directed neutrosophic graphs and their associated neutrosophic connection matrices.

$$C_1 = \{D_1 \ D_2 \ D_3), (R_1 \ R_2 \ R_3 \ldots R_5)\} \text{ and}$$
$$C_2 = \{(D_4 \ D_5 \ D_6), (R_6 \ R_7 \ R_8 \ R_9 \ R_{10})\}$$

The neutrosophic directed graph for the class $C_1$ given by the expert is as follows.

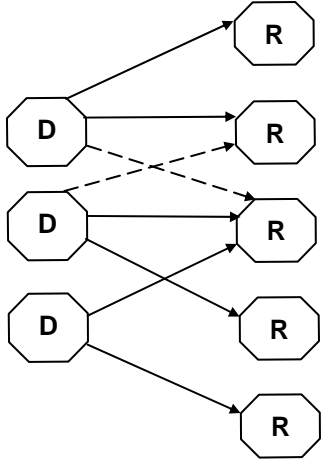

FIGURE: 5.6.1

The related neutrosophic connection matrix is

$$\begin{array}{c} \\ D_1 \\ D_2 \\ D_3 \end{array} \begin{array}{ccccc} R_1 & R_2 & R_3 & R_4 & R_5 \\ \begin{bmatrix} 1 & 1 & I & 0 & 0 \\ 0 & I & 1 & 1 & 0 \\ 0 & 0 & 1 & 0 & 1 \end{bmatrix} \end{array}$$

The directed graph for the class $C_2 = \{(D_4\ D_5\ D_6),\ (R_6\ R_7\ R_8\ R_9\ R_{10})\}$ given by the expert is as follows.

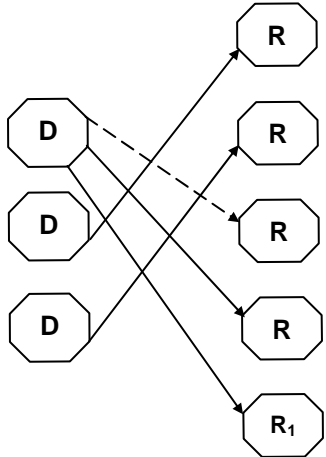

FIGURE: 5.6.2

The related neutrosophic connection matrix is

$$
\begin{array}{c}
 \\
D_4 \\
D_5 \\
D_6
\end{array}
\begin{array}{c}
R_6\ R_7\ R_8\ R_9\ R_{10} \\
\begin{bmatrix}
0 & 0 & I & 1 & 1 \\
1 & 0 & 0 & 0 & 0 \\
0 & 1 & 0 & 0 & 0
\end{bmatrix}
\end{array}
$$

Using these two connection matrices we obtain the $6 \times 10$ matrix associated with the combined disjoint block NRM which we denote by C(N).

$$
\begin{array}{c}
 \\
D_1 \\
D_2 \\
D_3 \\
D_4 \\
D_5 \\
D_6
\end{array}
\begin{array}{c}
R_1\ R_2\ R_3\ R_4\ R_5\ R_6\ R_7\ R_8\ R_9\ R_{10} \\
\begin{bmatrix}
1 & 1 & I & 0 & 0 & 0 & 0 & 0 & 0 & 0 \\
0 & I & 1 & 1 & 0 & 0 & 0 & 0 & 0 & 0 \\
0 & 0 & 1 & 0 & 1 & 0 & 0 & 0 & 0 & 0 \\
0 & 0 & 0 & 0 & 0 & 0 & 0 & I & 1 & 1 \\
0 & 0 & 0 & 0 & 0 & 1 & 0 & 0 & 0 & 0 \\
0 & 0 & 0 & 0 & 0 & 0 & 1 & 0 & 0 & 0
\end{bmatrix}
\end{array}
$$

Now we analyze the effect of the state vector X = (0 0 1 0 0 0) i.e., only the attribute $D_3$ alone is in the on state and all other nodes are in the off state.

$$
\begin{array}{llll}
XC\,(N) & \hookrightarrow & (0\ 0\ 1\ 0\ 1\ 0\ 0\ 0\ 0\ 0) & = & Y \\
YC\,(N)^T & \hookrightarrow & (I\ 1\ 1\ 0\ 0\ 0) & = & X_1 \\
X_1\,C(N) & \hookrightarrow & (I\,I\ 1\ 1\ 1\ 0\ 0\ 0\ 0\ 0) & = & Y_1 \\
Y_1\,C(N)^T & \hookrightarrow & (I\ 1\ 1\ 0\ 0\ 0) & = & X_2\,(= X_1)
\end{array}
$$

The hidden pattern is a fixed point given by the binary pair {(I 1 1 0 0 0), (I I 1 1 1 00 0 0 0)}, which has influence in the system which is evident by the resultant vector.

Let us now consider the state vector Y = (0 0 0 0 0 0 1 0 0 0) i.e., only the node $R_7$ is in the on state and all other nodes are in the off state. Effect on Y on the neutrosophic dynamical system C(N) is given by

$$
\begin{array}{llll}
YC\,(N)^T & \hookrightarrow & (0\ 0\ 0\ 0\ 0\ 1) & = & X \\
XC\,(N) & \hookrightarrow & (0\ 0\ 0\ 0\ 0\ 0\ 1\ 0\ 0\ 0)
\end{array}
$$

Thus we see the binary pair {(0 0 0 0 0 0 1 0 0 0). Thus pair {(0 0 0 0 0 0 1 0 0 0), (0 0 0 0 0 1) is a fixed point of the dynamical system. This pair has no influence on the dynamical system. Now let us consider the model given in page with the domain space {$M_1$ $M_2$,…, $M_9$} where this expert wishes to combine the attributes $M_9$ and $M_{10}$.

The attributes of the range space are taken as {$G_1$, …, $G_6$} where $G_1$ and $G_2$ are combined together to form the attribute $G_2$ so $G_3$ is labeled as $G_2$, $G_4$ as $G_3$ and so on $G_6$ as $G_5$. Now we use the method of disjoint block decomposite and divide these sets into 3 classes

$$C_1 = \{(M_1\ M_2\ M_3),\ (G_1\ G_2)\},$$
$$C_2 = \{(M_4\ M_5\ M_6),\ (G_3\ (G_4)\}\ and$$
$$C_3 = \{(M_7\ M_8\ M_9),\ (G_5\ (G_6)\}.$$

Now using the expert's opinion we give the related directed neutrosophic graphs for the classes $C_1$, $C_2$ and $C_3$ respectively.

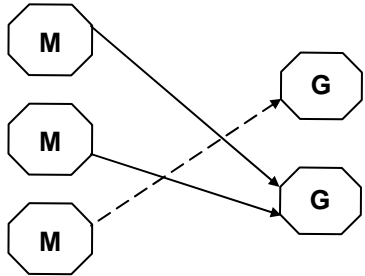

FIGURE: 5.6.3

The associated neutrosophic connection matrix is

$$\begin{array}{c} \\ M_1 \\ M_2 \\ M_3 \end{array} \begin{array}{cc} G_1 & G_2 \\ \begin{bmatrix} 0 & 1 \\ 0 & 1 \\ I & 0 \end{bmatrix} \end{array}$$

The directed neutrosophic graph given by the expert related with the class $C_2$.

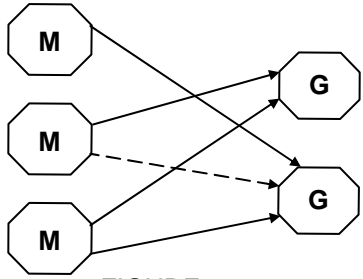

FIGURE: 5.6.4

The related connection matrix.

$$
\begin{array}{c c}
 & G_3 \quad G_4 \\
\begin{matrix} M_4 \\ M_5 \\ M_6 \end{matrix} &
\begin{bmatrix} 0 & 1 \\ 1 & I \\ 1 & 1 \end{bmatrix}
\end{array}
$$

Now using the experts opinion give the related neutrosophic directed graph for the class $C_3$.

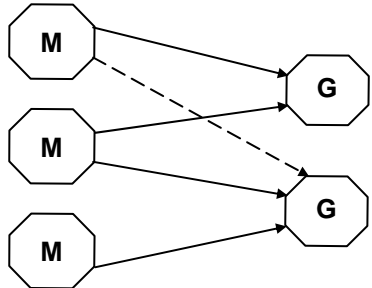

FIGURE: 5.6.5

The related connection matrix

$$
\begin{array}{c c}
 & G_5 \quad G_6 \\
\begin{matrix} M_7 \\ M_8 \\ M_9 \end{matrix} &
\begin{bmatrix} 1 & I \\ 1 & 1 \\ 0 & 1 \end{bmatrix}
\end{array}
$$

Now using these three connection matrices we get the related relational matrix of the combined disjoint block NRM which is denoted by C(M) and C(M) is a $9 \times 6$ matrix.

$$\begin{array}{c} \\ M_1 \\ M_2 \\ M_3 \\ M_4 \\ M_5 \\ M_6 \\ M_7 \\ M_8 \\ M_9 \end{array} \begin{array}{cccccc} G_1 & G_2 & G_3 & G_4 & G_5 & G_6 \\ \left[\begin{array}{cccccc} 0 & 1 & 0 & 0 & 0 & 0 \\ 0 & 1 & 0 & 0 & 0 & 0 \\ I & 0 & 0 & 0 & 0 & 0 \\ 0 & 0 & 0 & 1 & 0 & 0 \\ 0 & 0 & 1 & I & 0 & 0 \\ 0 & 0 & 1 & 1 & 0 & 0 \\ 0 & 0 & 0 & 0 & 1 & I \\ 0 & 0 & 0 & 0 & 1 & 1 \\ 0 & 0 & 0 & 0 & 0 & 1 \end{array}\right] \end{array}$$

Let us consider the effect of the state vector X = (0 0 0 1 0 0 0 0 1) of the dynamical system where only the attributes $M_4$ and $M_9$ are in the on state and all other nodes are in the off state. The effect of X on C(M) is given by

$$\begin{array}{lllll} \text{XC (M)} & \hookrightarrow & (0\,0\,0\,1\,0\,1) & = & \text{Y} \\ \text{Y C(M)}^\text{T} & \hookrightarrow & (0\,0\,0\,1\,I\,1\,I\,1\,1) & = & \text{X}_1 \\ \text{X}_1\,\text{C(M)} & \hookrightarrow & (0\,0\,I\,1\,I\,1) & = & \text{Y}_1 \\ \text{Y}_1\,\text{C(M)} & \hookrightarrow & (0\,0\,0\,1\,I\,I\,I\,I\,1) & = & \text{X}_2 \\ \text{X}_2\,\text{C(M)} & \hookrightarrow & (0\,0\,I\,I\,I\,1) & = & \text{Y}_2 \\ \text{Y}_2\,\text{C (M)} & \hookrightarrow & (0\,0\,0\,1\,I\,I\,I\,I\,I\,I) & = & \text{X}_3 \\ \text{X}_3\,\text{C (M)} & \hookrightarrow & (0\,0\,I\,I\,I\,I) & = & \text{Y}_3 = \text{Y}_2 \\ \text{Y}_2\,\text{C(M)} & \hookrightarrow & \text{X}_3 & & \end{array}$$

Thus the binary pair is a fixed point given by {(0 0 $I$ $I$ $I$ I), (0 0 0 1 $I$ $I$ $I$ $I$ 1)} which shows a very strong influence on the system.

Now we illustrate the model in which the blocks are disjoint but the number of blocks are not same. We consider the model with the sets of attributes given by {$A_1$, $A_2$, …,$A_6$} and {$G_1$ $G_2$… $G_5$}. Let us divide into equivalence classes.

$$C_1 = \{(A_1\ A_2)\ (G_1\ G_4)\} \text{ and}$$
$$C_2 = \{(A_3\ A_4\ A_5\ A_6),\ (G_2\ G_3\ G_5)\}$$

The directed neutrosophic graph related with the class $C_1$.

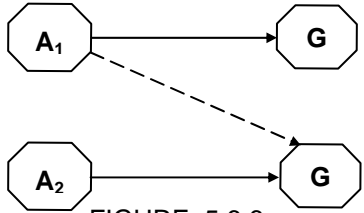

FIGURE: 5.6.6

The connection neutrosophic matrix associated with the above graph is

$$\begin{array}{c} \\ A_1 \\ A_2 \end{array} \begin{array}{cc} G_1 & G_4 \\ \begin{bmatrix} 1 & I \\ 0 & 1 \end{bmatrix} \end{array}$$

The related neutrosophic direct graph given the expert associated with the class $C_2$ is given by

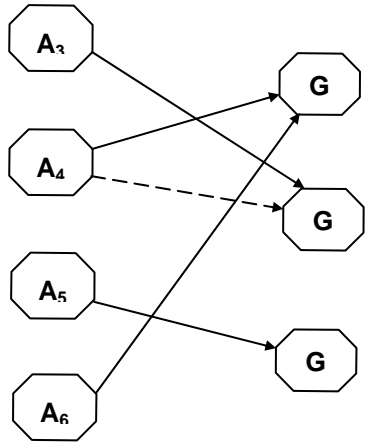

FIGURE: 5.6.7

The related connection matrix.

$$\begin{array}{c} \\ A_3 \\ A_4 \\ A_5 \\ A_6 \end{array} \begin{array}{ccc} G_2 & G_3 & G_5 \\ \begin{bmatrix} 0 & 0 & 0 \\ 1 & I & 1 \\ 0 & 0 & 1 \\ 1 & 0 & 0 \end{bmatrix} \end{array}$$

Now we give the associated connection $6 \times 5$ matrix C(S) given by

$$
\begin{array}{c}
\\
A_1 \\
A_2 \\
A_3 \\
A_4 \\
A_5 \\
A_6
\end{array}
\begin{array}{c}
\begin{array}{ccccc}
G_1 & G_2 & G_3 & G_4 & G_5
\end{array} \\
\left[
\begin{array}{ccccc}
1 & 0 & 0 & I & 0 \\
0 & 0 & 0 & 1 & 0 \\
0 & 0 & 0 & 0 & 0 \\
0 & 1 & I & 0 & 0 \\
0 & 0 & 0 & 0 & 1 \\
0 & 1 & 0 & 0 & 0
\end{array}
\right]
\end{array}
$$

C(S) is the neutrosophic matrix of the combined disjoint block NRM of unequal sizes, consider the state vector X = (1 0 0 1 0 0) in which the attributes $A_1$ and $A_4$ are in the on state and all other nodes are in the off state. Effect of X on the dynamical system C(S) is given by

$$X \ C(S) \quad \hookrightarrow \quad (1\ 1\ I\ I\ 0) \quad = \quad Y$$

$$Y \ C(S)^T \quad \hookrightarrow \quad (1\ I\ 0\ 1\ 0\ 1) \quad = \quad X_1$$

$$X_1 \ C(S) \quad \hookrightarrow \quad (1\ 1\ I\ I\ 0) \quad = \quad Y_1 \ (= Y)$$

Thus the hidden pattern of the dynamical system is the fixed point given by the binary pair {(1 1 $I$ $I$ 0) (1 $I$ 0 1 0 1)}. The nodes $A_1$ and $A_4$ together have a very strong impact on C(S). Now consider the state vector Y = (0 0 0 1 0) ie only the attribute $G_4$ is in the on state and all other vectors are in off state. The effect of Y on the dynamical system C(S) is given by

$$YC(S)^T \quad \hookrightarrow \quad (I\ 1\ 0\ 0\ 0\ 0) \quad = \quad X$$

$$X \ (C\ (S)) \quad \hookrightarrow \quad (I\ 0\ 0\ I\ 0) \quad = \quad Y_1$$

$$Y_1 \ C(S)^T \quad \hookrightarrow \quad (I\ 1\ 0\ 0\ 0\ 0) \quad = \quad X_1 = X$$

Thus the hidden point is the fixed point given by the binary pair {($I$ 1 0 0 0 0), ($I$ 0 0 1 0)}. The effect of the node G4 on the dynamical system C(S) is not very strong. Now let us consider yet a new model given in page 198 of this chapter let the domain space of the NRM be taken as {($D_1$ $D_2$,…, $D_8$) } and that of the range space be {$R_1$ $R_2$…, $R_5$)} consider the classes $C_1$ and $C_2$ where $C_1$ = {($D_1$ $D_2$ $D_3$), ($R_1$ $R_2$)} and $C_2$ = {($D_4$ $D_5$ $D_6$ $D_7$ $D_8$), ($R_3$ $R_4$ $R_5$)}. Now we obtain the experts opinion and give the directed graph related to them.

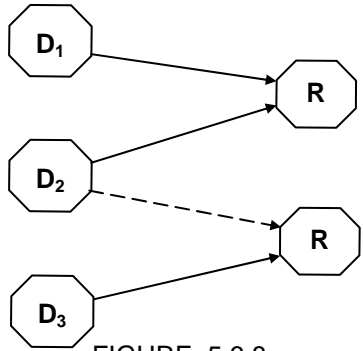

FIGURE: 5.6.8

The related neutrosophic connection matrix is

$$
\begin{array}{c}
\quad\quad R_1 \; R_2 \\
\begin{array}{c} D_1 \\ D_2 \\ D_3 \end{array}
\left[\begin{array}{cc}
1 & 0 \\
1 & I \\
0 & 1
\end{array}\right]
\end{array}
$$

The directed graph related to the class $C_2$ is as follows

The related connection matrix is given by

$$
\begin{array}{c}
\begin{array}{ccc} R_3 & R_4 & R_5 \end{array} \\
\begin{array}{c} D_4 \\ D_5 \\ D_6 \\ D_7 \\ D_8 \end{array}
\begin{bmatrix}
1 & 1 & 0 \\
1 & I & 0 \\
0 & 0 & 1 \\
0 & 1 & 0 \\
0 & 0 & I
\end{bmatrix}
\end{array}
$$

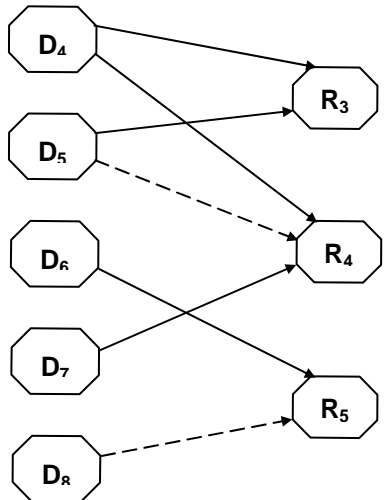

FIGURE: 5.6.9

Now we obtain the related neutrosophic connection matrix of the combined disjoint block NRM of different sizes given by the 8 × 5 matrix and is denoted by C(G)

$$
\begin{array}{c}
\begin{array}{ccccc} R_1 & R_2 & R_3 & R_4 & R_5 \end{array} \\
\begin{array}{c} D_1 \\ D_2 \\ D_3 \\ D_4 \\ D_5 \\ D_6 \\ D_7 \\ D_8 \end{array}
\begin{bmatrix}
1 & 0 & 0 & 0 & 0 \\
1 & I & 0 & 0 & 0 \\
0 & 1 & 0 & 0 & 0 \\
0 & 0 & 1 & I & 0 \\
0 & 0 & 1 & I & 0 \\
0 & 0 & 0 & 0 & 1 \\
0 & 0 & 0 & 1 & 0 \\
0 & 0 & 0 & 0 & I
\end{bmatrix}
\end{array}
$$

Now we obtain the effect of a state vector X = (0 0 0 0 1 0 0 0) where only the node $D_5$ is in the on state and all other nodes are in the off state. The resultant vector of X using the dynamical system C(G)

$$
\begin{array}{llll}
\text{XC (G)} & \hookrightarrow & (0\ 0\ 1\ I\ 0) & = & \text{Y} \\
\text{Y C (G)}^{\text{T}} & \hookrightarrow & (0\ 0\ 0\ I\ I\ 0\ 0\ 0) & = & \text{X}_1 \\
\text{X}_1\ \text{C(G)} & \hookrightarrow & (0\ 0\ 1\ I\ 0) & = & \text{Y}_1\ (=\text{Y})
\end{array}
$$

Thus the fixed point is given by the binary pair {(0 0 1 $I$ 0), (0 0 0 1 1 0 0 0)}, which has a least effect on the system.

Let us consider the state vector Y = (1 0 0 0 0) that is $R_1$ is in the on state and all nodes in the rage space are in the off state

$$
\begin{array}{llll}
\text{YC (G)}^{\text{T}} & \hookrightarrow & (1\ 1\ 0\ 0\ 0\ 0\ 0\ 0) & = & \text{X} \\
\text{XC (G)} & \hookrightarrow & (1\ I\ 0\ 0\ 0) & = & \text{Y}_1
\end{array}
$$

$$YC(G)^T \quad \hookrightarrow \quad (1\ I\ I\ 0\ 0\ 0\ 0\ 0\ 0) \quad = \quad X_1$$
$$X_1\ C(G) \quad \hookrightarrow \quad (1\ I\ 0\ 0\ 0) \quad = \quad Y_2\ (=Y_1)$$

Thus the hidden pattern of the dynamical system is a fixed point given by the binary pair $\{(1\ I\ 0\ 0\ 0), (1\ I\ I\ 0\ 0\ 0\ 0\ 0)\}$ according to this expert the effect of the node $R_1$ is least on the system. Consider the NRM given in the attributes are taken as $\{M_1\ M_2,\ldots M_{10}\}$ and $\{G_1\ G_2,\ldots G_7\}$ the combined block disjoint NRM of varying sizes is calculated using the classes $C_1 = \{M_1\ M_2\ M_3\ M_4\}\ (G_1\ G_2\ )\}$, $C_2 = \{\ M_5\ M_6\ M_7\ M_8\}\ (G_4\ G_5\ )\}$ and $C_3 = \{(M_9\ M_{10}),\ (G_6,\ G_7)\}$.

The neutrosophic directed graph of class $C_1$ given by the expert is as follows.

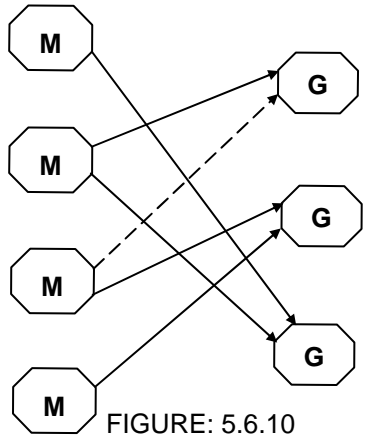

FIGURE: 5.6.10

The related neutrosophic connection matrix

$$\begin{array}{c c} & \begin{matrix} G_1 & G_2 & G_3 \end{matrix} \\ \begin{matrix} M_1 \\ M_2 \\ M_3 \\ M_4 \end{matrix} & \begin{bmatrix} 0 & 0 & 1 \\ 1 & 0 & 1 \\ I & 1 & 0 \\ 0 & 1 & 0 \end{bmatrix} \end{array}$$

Now using the expert opinion we give directed graph for class $C_2$.

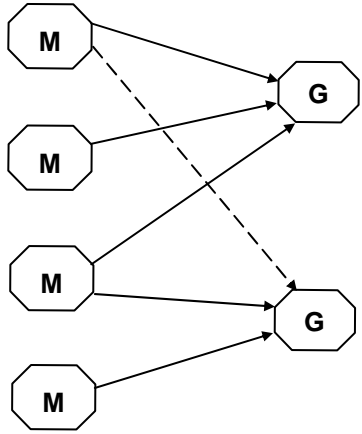

FIGURE: 5.6.11

The related connection matrix

$$\begin{array}{c} \quad\quad G_4 \;\; G_5 \\ \begin{array}{c} M_5 \\ M_6 \\ M_7 \\ M_8 \end{array} \begin{bmatrix} 1 & I \\ 1 & 0 \\ 1 & 1 \\ 0 & 1 \end{bmatrix} \end{array}$$

Now the directed graph given by the expert is

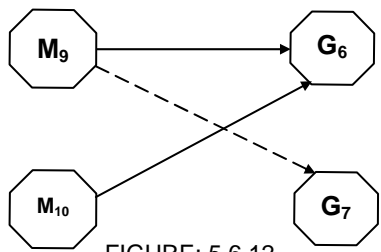

FIGURE: 5.6.12

The related neutrosophic matrix

$$\begin{array}{c} \quad\quad G_6 \;\; G_7 \\ \begin{array}{c} M_9 \\ M_{10} \end{array} \begin{bmatrix} 1 & I \\ 1 & 0 \end{bmatrix} \end{array}$$

Now we give the relational connection matrix of the combined disjoint block NRM of varying sizes. We denote this $10 \times 7$ matrix by N(W)

$$
\begin{array}{c}
\\
M_1 \\
M_2 \\
M_3 \\
M_4 \\
M_5 \\
M_6 \\
M_7 \\
M_8 \\
M_9 \\
M_{10}
\end{array}
\begin{array}{c}
G_1\ G_2\ G_3\ G_4\ G_5\ G_6\ G_7 \\
\left[\begin{array}{ccccccc}
0 & 0 & 1 & 0 & 0 & 0 & 0 \\
1 & 0 & 1 & 0 & 0 & 0 & 0 \\
I & 1 & 0 & 0 & 0 & 0 & 0 \\
0 & 1 & 0 & 0 & 0 & 0 & 0 \\
0 & 0 & 0 & 1 & I & 0 & 0 \\
0 & 0 & 0 & 1 & 0 & 0 & 0 \\
0 & 0 & 0 & 1 & 1 & 0 & 0 \\
0 & 0 & 0 & 0 & 1 & 0 & 0 \\
0 & 0 & 0 & 0 & 0 & 1 & I \\
0 & 0 & 0 & 0 & 0 & 1 & 0
\end{array}\right]
\end{array}
$$

Let us consider the state vector $X = (0\ 0\ 0\ 0\ 0\ 0\ 0\ 1\ 0\ 0)$ where only node $M_8$ alone is in the on state and all other nodes are in the off state. The effect of X on the dynamical system C(W) is given by

$$
\begin{array}{lll}
XC(W) \hookrightarrow (0\ 0\ 0\ 0\ 1\ 0\ 0) & = & Y \\
YC(W)^T \hookrightarrow (0\ 0\ \ 0\ 0\ I\ 0\ 1\ 0\ 0) & = & X_1 \\
X_1\ X(W) \hookrightarrow (0\ 0\ 0\ I\ 1\ 0\ 0) & = & Y_1 \\
Y_1\ (C(W))^T \hookrightarrow (0\ 0\ 0\ 0\ I\ I\ I\ I\ 1\ 0\ 0) & & \\
X_2\ C(W) \hookrightarrow (0\ 0\ 0\ I\ I\ 0\ 0) & = & Y_2 \\
Y_2\ C(W)^T \hookrightarrow (0\ 0\ 0\ 0\ I\ I\ I\ 1\ 0\ 0) & = & X_3\ (=X_2).
\end{array}
$$

Thus the fixed point is given by the binary pair $\{(0\ 0\ 0\ I\ I\ 0\ 0), (0\ 0\ 0\ 0\ I\ I\ I\ 1\ 0\ 0)\}$. The reader is expected to give the effect of the node $M_8$ on the system.

Next we consider the combined overlap block NRM used in modelling the HIV/AIDS migrant labourers problem. Let us consider the attributes given 198. The classes of attributes used are $\{D_1, D_2,\ldots, D_8\}$ and $\{R_1, R_2,\ldots, R_5\}$ We divide these into over lap blocks $C_1, C_2\ C_3$ and $C_4$ where $C_1 = \{(D_1\ D_2\ D_3\ D_4), (R_1\ R_2$

$R_3$)}, $C_2$ = {($D_3$ $D_4$ $D_5$ $D_6$), ($R_2$ $R_3$ $R_4$)}, $C_3$ = {($D_5$ $D_6$, $D_7$, $D_8$), ($R_3$ $R_4$ $R_5$)} and $C_4$ = {($D_7$ $D_8$ $D_1$, $D_2$), ($R_4$ $R_5$ $R_1$)}.

Now we obtain using the experts opinion the directed graph of these NRMs. The directed graph related to the class $C_1$.

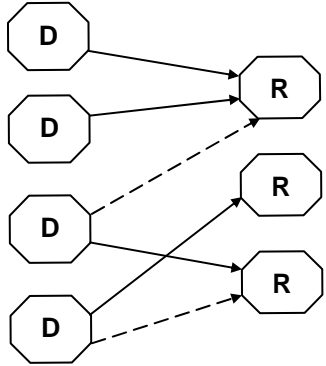

FIGURE: 5.6.13

The related connection matrix

$$\begin{array}{c} \\ D_1 \\ D_2 \\ D_3 \\ D_4 \end{array} \begin{array}{ccc} R_1 & R_2 & R_3 \\ \begin{bmatrix} 1 & 0 & 0 \\ 1 & 0 & 0 \\ I & 0 & 1 \\ 0 & 1 & I \end{bmatrix} \end{array}$$

The directed graph given by the expert related to the class $C_2$.

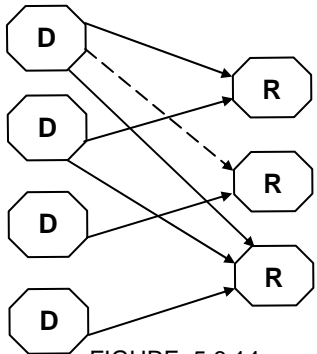

FIGURE: 5.6.14

The related connection matrix is given as follows

$$\begin{array}{c} \quad\quad R_2 \;\; R_3 \;\; R_4 \\ \begin{array}{c} D_3 \\ D_4 \\ D_5 \\ D_6 \end{array} \left[ \begin{array}{ccc} 1 & I & 1 \\ 1 & 0 & 1 \\ 0 & 1 & 0 \\ 0 & 0 & 1 \end{array} \right] \end{array}$$

Directed graph given by the expert for the class $C_3$ is as follows.

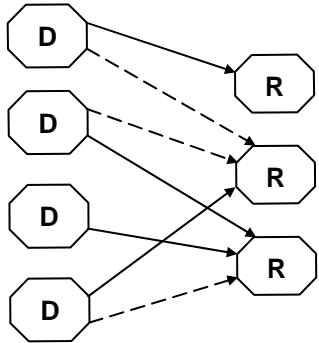

FIGURE: 5.6.15

The related connection matrix of the above directed graph

$$\begin{array}{c} \quad\quad R_3 \;\; R_4 \; R_5 \\ \begin{array}{c} D_5 \\ D_6 \\ D_7 \\ D_8 \end{array} \left[ \begin{array}{ccc} 1 & I & 0 \\ 0 & I & 1 \\ 0 & 0 & 1 \\ 0 & 1 & I \end{array} \right] \end{array}$$

The directed graph given by the expert for the last class $C_4$.

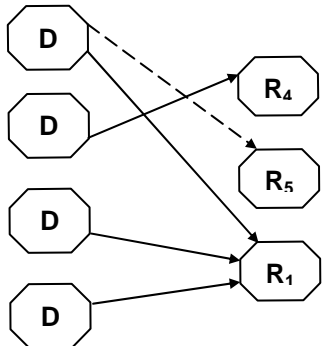

FIGURE: 5.6.16

The related neutrosophic matrix

$$
\begin{array}{c}
 \\
D_7 \\
D_8 \\
D_1 \\
D_2
\end{array}
\begin{array}{ccc}
R_4 & R_5 & R_6 \\
\left[\begin{array}{ccc}
0 & I & 1 \\
1 & 0 & 0 \\
0 & 0 & 1 \\
0 & 0 & 1
\end{array}\right]
\end{array}
$$

Now using these four connection matrices we obtain the relational neutrosophic matrix of the combined block overlap NRM denoted by N(V)

$$
\begin{array}{c}
 \\
D_1 \\
D_2 \\
D_3 \\
D_4 \\
D_5 \\
D_6 \\
D_7 \\
D_8
\end{array}
\begin{array}{ccccc}
R_1 & R_2 & R_3 & R_4 & R_5 \\
\left[\begin{array}{ccccc}
2 & 0 & 0 & 0 & 0 \\
2 & 0 & 0 & 0 & 0 \\
I & 1 & I & 1 & 0 \\
0 & 2 & I & 1 & 0 \\
0 & 0 & 2 & 1 & 0 \\
0 & 0 & 0 & I & 1 \\
1 & 0 & 0 & 0 & I \\
0 & 0 & 0 & 2 & I
\end{array}\right]
\end{array}
$$

Let X = (0 0 0 0 1 0 0 0) be the state vector with the node $D_5$ to be in the on state and all other nodes are in the off state the effect of X on the dynamical system N(V)

$$
\begin{aligned}
XN(V) &\hookrightarrow (0\,0\,1\,1\,0) &&= Y \\
YN(V)^T &\hookrightarrow (0\,0\,I\,I\,1\,I\,0\,1) &&= X_1 \\
X_1\,N(V) &\hookrightarrow (I\,I\,I\,1\,I) &&= Y_1 \\
Y_1\,N(V)^T &\hookrightarrow (I\,I\,I\,I\,1\,I\,I\,1) &&= X_2 \\
X_2\,(N(V)) &\hookrightarrow (I\,I\,I\,1\,I) &&= Y_2\,(=Y_1).
\end{aligned}
$$

Thus the fixed point is a binary pair given by {($I\,I\,I\,1\,I$), ($I\,I\,I\,I\,1\,I\,I\,1$)}which has a very strong influence on the dynamical system as D5 makes all nodes to on state or in an indeterminate state.

Next we consider the state vector Y = (1 0 0 0 ) that is only the node $R_1$ is in the on state and all other nodes are in the off state. The effect of Y on the dynamical system N(V) is given by

$$Y(N(V))^T \quad \hookrightarrow \quad (1\ 1\ I\ 0\ 0\ 0\ 1\ 0) \qquad = \quad X$$

$$XN(V) \quad \hookrightarrow \quad (1\ I\ I\ I\ I) \qquad = \quad Y_1$$

$$Y_1(N(V)) \quad \hookrightarrow \quad (1\ 1\ I\ I\ I\ I\ I\ I) \qquad = \quad X_1$$

$$X_1\ N(V) \quad \hookrightarrow \quad (1\ I\ I\ I\ I) \qquad = \quad Y_2 = Y_1.$$

Thus the hidden pattern of the dynamical system is a fixed point given by the binary pair $\{(1\ I\ I\ I\ I)\ (1\ 1\ I\ I\ I\ I\ I\ I)\}$ as before this node has very strong impact on the system.

Next we consider the combined overlap block NRM to study the given model in page 205 - 208 of this chapter. Take the sets of attributes given by $\{M_1\ M_2\ldots,\ M_{10}\}$ and $(G_1\ G_2\ldots\ G_7)$. We divide them into over lapping blocks given by the classes $C_1$, $C_2$ and $C_3$. where $C_1$ $\{(M_1\ M_2\ M_3\ M_4),\ (G_1\ G_2\ G_3)\}$ $C_2 = \{(M_4\ M_5\ M_6\ M_7), (G_3\ G_4\ G_5)\}$ and $C_3 = \{(M_7\ M_8\ M_9\ M_{10}),\ (G_5\ G_6\ G_7)\}$.

The directed graph for the class $C_1$ as given by the expert as follows

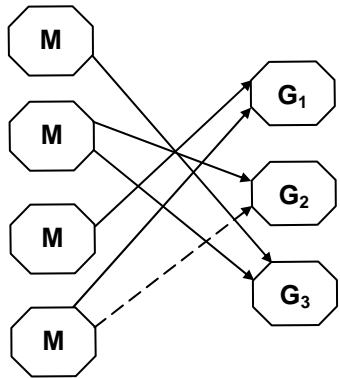

FIGURE: 5.6.17

The related connection matrix of the NRM is

$$
\begin{array}{c}
\phantom{M_1}\ \ G_1\ \ G_2\ \ G_3 \\
\begin{array}{c}
M_1 \\ M_2 \\ M_3 \\ M_4
\end{array}
\begin{bmatrix}
0 & 0 & 1 \\
0 & 1 & 1 \\
1 & 0 & 0 \\
1 & I & 0
\end{bmatrix}
\end{array}
$$

The neutrosophic directed graph given by the expert for the class $C_2$.

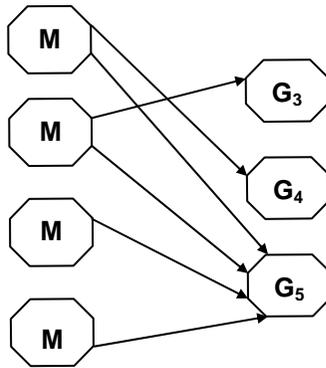

FIGURE: 5.6.18

The related connection matrix of the above directed graph is as follows

$$
\begin{array}{c}
\begin{array}{ccc} G_3 & G_4 & G_5 \end{array} \\
\begin{array}{c} M_4 \\ M_5 \\ M_6 \\ M_7 \end{array}
\begin{bmatrix}
0 & 1 & 1 \\
1 & 1 & 0 \\
0 & 0 & 1 \\
0 & 0 & 1
\end{bmatrix}
\end{array}
$$

Now we using the experts opinion give the directed graph of the last class.

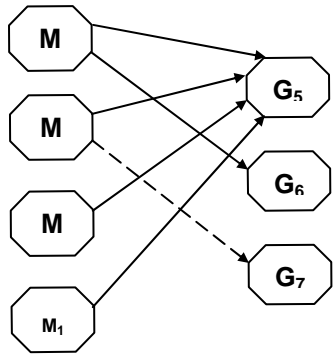

FIGURE: 5.6.19

The related neutrosophic connection matrix is given below

$$\begin{array}{c} \begin{array}{ccc} G_5 & G_6 & G_7 \end{array} \\ \begin{array}{c} M_7 \\ M_8 \\ M_9 \\ M_{10} \end{array} \begin{bmatrix} 1 & 1 & 0 \\ I & 0 & 1 \\ 1 & 0 & 0 \\ 1 & 0 & 0 \end{bmatrix} \end{array}$$

Now using there three connection matrices we obtain the 10 x 7 connection neutrosophic matrix associated with the combined block overlap NRM where the blocks are of equal size. We denote this matrix by N(Q).

$$\begin{array}{c} \begin{array}{ccccccc} G_1 & G_2 & G_3 & G_4 & G_5 & G_6 & G_7 \end{array} \\ \begin{array}{c} M_1 \\ M_2 \\ M_3 \\ M_4 \\ M_5 \\ M_6 \\ M_7 \\ M_8 \\ M_9 \\ M_{10} \end{array} \begin{bmatrix} 0 & 0 & 1 & 0 & 0 & 0 & 0 \\ 0 & 1 & I & 0 & 0 & 0 & 0 \\ 1 & 0 & 0 & 0 & 0 & 0 & 0 \\ 1 & I & 0 & 1 & 1 & 0 & 0 \\ 0 & 0 & 1 & 1 & 0 & 0 & 0 \\ 0 & 0 & 0 & 0 & 1 & 0 & 0 \\ 0 & 0 & 0 & 0 & 1 & 1 & 1 \\ 0 & 0 & 0 & 0 & I & 0 & 1 \\ 0 & 0 & 0 & 0 & 1 & 0 & 0 \\ 0 & 0 & 0 & 0 & 1 & 0 & 0 \end{bmatrix} \end{array}$$

Let $X = (0\ 0\ 0\ 0\ 0\ 0\ 0\ 1\ 0\ 0)$ be the state vector in which only the node $M_8$ is in the on state and all other nodes are in the off state. The effect of X on the dynamical system N(Q).

$$\begin{array}{llll}
XN(Q) & \hookrightarrow & (0\ 0\ 0\ 0\ 1\ 0\ 1) & = & Y \\
YN(Q)^T & \hookrightarrow & (0\ 0\ 0\ 1\ 0\ 1\ 1\ 1\ 1\ 1) & = & X_1 \\
X_1 N(Q) & \hookrightarrow & (1\ I\ 0\ 1\ 1\ 1\ 1) & = & Y_1 \\
Y_1 N(Q)^T & \hookrightarrow & (0\ I\ 1\ I\ 1\ 1\ 1\ 1\ 1\ 1) & = & X_2 \\
X_2 N(Q) & \hookrightarrow & (I\ I\ I\ I\ I\ 1\ 1) & = & Y_2.
\end{array}$$

It is left as exercise to reader to find the fixed point.

Let us consider the state vector $Y = (0\ 0\ 1\ 0\ 0\ 0\ 0)$ ie only the node $G_3$ is on the on state and all other nodes are in the off state. The effect of Y on N(Q)

$$YN(Q)^T \quad \hookrightarrow \quad (1\ 1\ 0\ 1\ 0\ 0\ 0\ 0\ 0) = \quad X$$

$$XN(Q) \quad \hookrightarrow \quad (1\ 0\ 1\ 1\ 1\ \ 0\ 0) \quad = \quad Y_1$$

$$Y_1\,(NQ)^T \quad \hookrightarrow \quad (1\ 1\ 1\ 1\ 1\ 1\ 1\ 1\ 1\ 1) = \quad X_1$$

$$Y_1\,N(Q) \quad \hookrightarrow \quad (1\ 1\ 1\ 1\ 1\ 1\ 1).$$

Thus the hidden pattern of the dynamical system is a fixed point given by the binary pair {(1 1 1 1 1 1 1 1 1), (1 1 1 1 1 1 1 1 1)}. The node $G_3$ has a very strong influence on the system so all nodes come ti on state. We see the combined overlap block NRM is every sensitive to the situation.

Next we consider the combined overlap block FRM of different sizes. Let us consider the attributes given in page 201-202 with set of attributes {($D_1$ $D_2$ $D_3$ $D_4$), ($R_1$ $R_2$ $R_3$ $R_4$ $R_5$)}, $C_2$ = {($D_3$ $D_4$ $D_5$ ), ($R_3$ $R_4$ $R_5$ $R_6$)}, $C_3$ = {($D_4$ $D_5$ $D_6$ $D_1$, $D_2$), ($R_6$ $R_7$ $R_8$ $R_{10}$)}. Now we analyze the model.
The directed given by the expert for the class $C_1$ of the NRM.

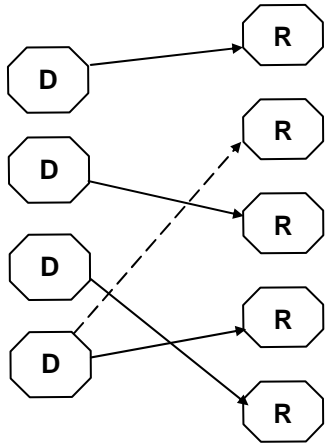

FIGURE: 5.6.20

The related relation matrix is as follows

$$\begin{array}{c@{\quad}ccccc}
 & R_1 & R_2 & R_3 & R_4 & R_5 \\
D_1 & \begin{bmatrix} 1 \\ 0 \\ 0 \\ 0 \end{bmatrix} & \begin{matrix} 0 \\ 0 \\ 0 \\ I \end{matrix} & \begin{matrix} 0 \\ 1 \\ 0 \\ 0 \end{matrix} & \begin{matrix} 0 \\ 0 \\ 0 \\ 1 \end{matrix} & \begin{matrix} 0 \\ 0 \\ 1 \\ 0 \end{bmatrix}
\end{array}$$

Directed graph given by the expert relative to the class $C_2$ is

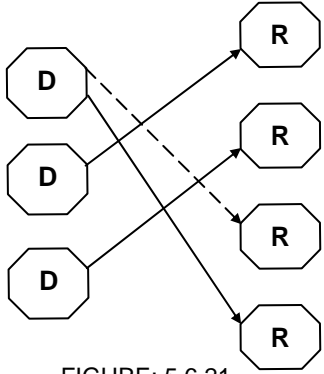

FIGURE: 5.6.21

The related matrix of the directed graph is as follows.

$$\begin{array}{c} \\ D_3 \\ D_4 \\ D_5 \end{array} \begin{array}{cccc} R_3 & R_4 & R_5 & R_6 \\ \begin{bmatrix} 0 & 0 & I & 1 \\ 1 & 0 & 0 & 0 \\ 0 & 1 & 0 & 0 \end{bmatrix} \end{array}$$

Now the expert opinion in he form of direct graph for the class $C_3$.

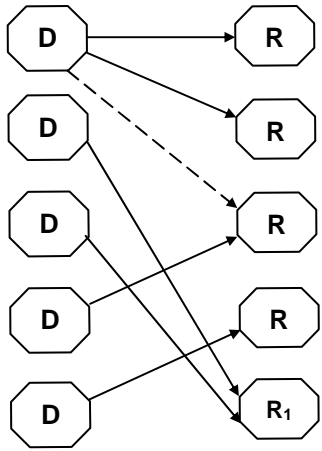

FIGURE: 5.6.22

The related connection matrix for this neutrosophic directed bigraph is

$$\begin{array}{c} \\ D_4 \\ D_5 \\ D_6 \\ D_1 \\ D_2 \end{array} \begin{array}{c} R_6 \; R_7 \; R_8 \; R_9 \; R_{10} \\ \left[ \begin{array}{ccccc} 1 & 1 & I & 0 & 0 \\ 0 & 0 & 0 & 0 & 1 \\ 0 & 0 & 0 & 0 & 1 \\ 0 & 0 & 1 & 0 & 0 \\ 0 & 0 & 0 & 1 & 0 \end{array} \right] \end{array}$$

Now using this matrices we give the connection matrix N((T) which gives the related combined overlap block NRM.

This is a $6 \times 10$ matrix as given below:

$$\begin{array}{c} \\ D_1 \\ D_2 \\ D_3 \\ D_4 \\ D_5 \\ D_6 \end{array} \begin{array}{c} R_1 \; R_2 \; R_3 \; R_4 \; R_5 \; R_6 \; R_7 \; R_8 \; R_9 \; R_{10} \\ \left[ \begin{array}{cccccccccc} 1 & 0 & 0 & 0 & 0 & 0 & 0 & 1 & 0 & 0 \\ 0 & 0 & 1 & 0 & 0 & 0 & 0 & 0 & 1 & 0 \\ 0 & 0 & 0 & 0 & 0 & 1 & 0 & 0 & 0 & 0 \\ 0 & I & 1 & 1 & 0 & 1 & 1 & I & 0 & 0 \\ 0 & 0 & 0 & 1 & 0 & 0 & 0 & 0 & 0 & 1 \\ 0 & 0 & 0 & 0 & 0 & 0 & 0 & 0 & 0 & 1 \end{array} \right] \end{array}$$

Now we consider the effect of the state vector X = (0 0 0 0 1 0) on the dynamical system where only the node $D_5$ is in the on state and all nodes are in the off state

$$\begin{array}{lllll}
\text{XN (T)} & \hookrightarrow & (0\,0\,0\,10\,0\,0\,0\,0\,1) & = & \text{Y} \\
\text{YN(T)}^\text{T} & \hookrightarrow & (0\,0\,0\,1\,1\,1) & = & \text{X}_1 \\
\text{X}_1\,\text{N(T)} & \hookrightarrow & (0\,I\,1\,1\,0\,1\,1\,I\,0\,1) & = & \text{Y}_1 \\
\text{Y}_1\,\text{N(T)}^\text{T} & \hookrightarrow & (I\,1\,1\,1\,1\,1) & = & \text{X}_2 \\
\text{X}_2\,\text{N(T)} & \hookrightarrow & (I\,I\,1\,10\,1\,1\,0\,1\,1) & = & \text{Y}_2 \\
\text{Y}_2\,(\text{N(T)})^\text{T} & \hookrightarrow & (I\,1\,1\,1\,1\,1) & = & \text{X}_3 = (=\text{X}_2).
\end{array}$$

Thus the hidden pattern of the dynamical system is a fixed point given by the binary pair $\{(I\ 1\ 1\ 1\ 1\ 1),\ (I\ I\ 1\ 1\ 0\ 1\ 1\ 0\ 1\ 1)\}$.

The influence of $D_5$ on the dynamical system is very strong so it affects most of the nodes. Let us consider the state vector Y = (0 0 1 0 0 0) where the attribute $D_3$ alone is in the on state and all other nodes are the off state.

The effect of Y on the dynamical system (N(T)) is given by

$$
\begin{array}{llll}
YN\,(T) & \hookrightarrow & (0\ 0\ 0\ 0\ 0\ 1\ 0\ 0\ 0\ 0) & = \quad X \\
XN\,(T)^T & \hookrightarrow & (0\ 0\ 1\ 1\ 0\ 0) & = \quad Y_1 \\
Y_1\,N(T) & \hookrightarrow & (0\ I\ 1\ 1\ 0\ 1\ 1\ I\ 0\ 0) & = \quad X_1 \\
X_1\,N(T)^T & \hookrightarrow & (I\ 1\ 1\ 1\ 1\ 0) & = \quad Y_2 \\
Y_2\,N(T) & \hookrightarrow & (I\ I\ 1\ 1\ 0\ 11\ I\ 1\ 1) & = \quad X_2 \\
X_2\,N(T)^T & \hookrightarrow & (I\ 1\ 1\ 1\ 1\ 1) & = \quad Y_3 \\
Y_3\,N(T) & \hookrightarrow & (I\ I\ 1\ 1\ 0\ 1\ 1\ I\ 1\ 1) & = \quad X_3\ (=X_2).
\end{array}
$$

Thus the hidden pattern of the dynamical system is the fixed point given by the binary pair. $\{(I\ 1\ 1\ 1\ 1\ 1),\ (I\ I\ 1\ 1\ 0\ 11I\ 1\ 1)\}$. The $D_3$ also has a strong effect on the dynamical system for it influences all nodes expect $R_5$.

We now study another set of attributes related with the HIV/AIDS affected migrant labours. Let us consider the sets attributes $\{(A_1\ A_2\ ,.... \ A_6),\ (G_1\ G_2,... \ G_5)\}$ of the NRM. This is divided into two overlapping classes $C_1$ and $C_2$

$$C_1 = \{(A_1\ A_2\ A_3\ A_4),\ (G_1\ G_2\ G_3)\}$$
and
$$C_2 = \{(A_4\ A_5\ A_6)\ (G_3\ G_4\ G_5)\}.$$

The directed graph as given by the expert for the attributes in class $C_1$.

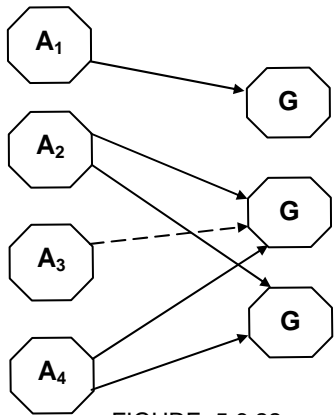

FIGURE: 5.6.22

The related relational neutrosophic matrix is

$$
\begin{array}{c@{}c}
 & \begin{array}{ccc} G_1 & G_2 & G_3 \end{array} \\
\begin{array}{c} A_1 \\ A_2 \\ A_3 \\ A_4 \end{array} &
\left[ \begin{array}{ccc}
1 & 0 & 0 \\
0 & 1 & 1 \\
0 & I & 0 \\
0 & 1 & 1
\end{array} \right]
\end{array}
$$

The directed graph given by the expert for the class $C_2 = \{(A_4\ A_5\ A_6)\ (G_3\ G_4\ G_5)\}$ is as follows:

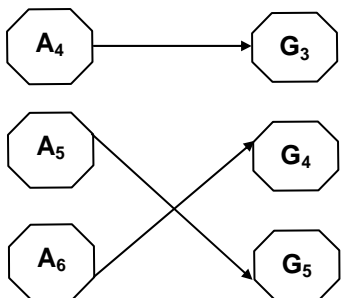

FIGURE: 5.6.23

The related connection matrix

$$\begin{array}{c} \begin{array}{ccc} G_3 & G_4 & G_5 \end{array} \\ \begin{array}{c} A_4 \\ A_5 \\ A_6 \end{array} \begin{bmatrix} 1 & 0 & 0 \\ 0 & 0 & 1 \\ 0 & 1 & 0 \end{bmatrix} \end{array}$$

The matrix related to the combined block overlap NRM denoted by N(V).

$$\begin{array}{c} \begin{array}{ccccc} G_1 & G_2 & G_3 & G_4 & G_5 \end{array} \\ \begin{array}{c} A_1 \\ A_2 \\ A_3 \\ A_4 \\ A_5 \\ A_6 \end{array} \begin{bmatrix} 1 & 0 & 0 & 0 & 0 \\ 0 & 1 & 0 & 0 & 0 \\ 0 & I & 0 & 0 & 0 \\ 0 & 1 & 2 & 0 & 0 \\ 0 & 0 & 0 & 0 & 1 \\ 0 & 0 & 0 & 1 & 0 \end{bmatrix} \end{array}$$

Consider the state vector $X = (0\ 0\ 0\ 10\ 0)$ in which only the node $A_4$ is in the on state and all other nodes are in the off state. The effect of X on N(V)

$$\begin{array}{llll} XN(U) & \hookrightarrow & (0\ 1\ 1\ 0\ 0) & = & Y \\ YN\ (U)^T & \hookrightarrow & (0\ 1\ I\ 1\ 0\ 0) & = & X_1 \\ X_1\ N(U) & \hookrightarrow & (0\ 1\ 1\ 0\ 0) & = & Y_1 = Y. \end{array}$$

Thus the hidden pattern of the dynamical system is a fixed point given by the binary pair $\{(0\ 1\ 1\ 0\ 0\ ),\ (0\ 1\ I\ 1\ 0\ 0)\}$. It has meager effect on the system

Next we consider the state vector $Y = (1\ 0\ 0\ 0\ 0)$ i.e., only the node $G_1$ is in the on state all other nodes are in the off state. The effect of Y on the dynamical system N(U) is given by

$$\begin{array}{llll} G_1\ N(U)^T & \hookrightarrow & (1\ 0\ 0\ 0\ 0\ 0) & = & X \\ XN\ (U) & \hookrightarrow & (1\ 0\ 0\ 0\ 0) & = & Y_1\ (=Y). \end{array}$$

Thus the hidden pattern of the dynamical system is a binary pair given by $\{(1\ 0\ 0\ 0\ 0)\ (1\ 0\ 0\ 0\ 0\ 0)\}$. The effect is least on the system.

It is important to mention here that once the connection matrix of the combined block disjoint NRM, or combined block overlap NRM is formed it is a matter of routine to determine the fixed point or limit cycle. Once again we leave it for the reader to modify the FRM-program to the case of NRM only in it is matter of an addition $I$ also in the coordinate where the indeterminate $I$ behaves as $I.I = I$ but $1 + I = 0$ r + $I$ if r > 1 is 1 and $1 + rI$ if r > 1 is $I$. Thus the interested reader to requested to construct a C program.

Chapter Six

# INTERVIEWS OF 60 MIGRANT LABOURERS LIVING WITH HIV/AIDS

This chapter is entirely devoted to giving the condensed data from the interviews taken by us of 60 people living with HIV/AIDS who represent our sample. All the 60 male members of our sample had become victims of HIV while migrating or in migration or acquired the disease from a place that was not their hometown. Also we state that some of the interviewed patients did not talk with us freely; they displayed reluctance in speaking to us about their disease and a lot of indeterminacy prevailed while talking to them.

The very analysis of the migrants and AIDS are mainly based on the interviews  it has become essential to give the interviews of them. The reader is advised to be non-critical over the interview for two reasons. One as much as possible we have given the interview verbatim not even altering the sequence of their expressions and talks. Secondly it is seen in some places the interviews of some migrant labourers are little contradictory. As we do not want to even polish or make any changes in the interviews given by them the observations made by us are put only as additional observations. These interviews may throw light on the living conditions and the socio-economic problems faced by the migrant labourers.

At some places in the interviews we find contradicting answers, yet we have presented it as it is, for when we questioned them twice or tried to make some confirmation about facts, they refused to be frank and they just bluntly said, "Do not trouble us". One has to understand that because of this it has become impossible for us to make our data doubly confirmed. One of the major observations we place at the outset is that several of the patients refused to give the name of their native place. We have also respected their desire for anonymity and as a result in this book, we have not mentioned any of the names of the patients or their relatives. Most of them do not know the difference between HIV and AIDS. For at least 80% of them, on coming to know that they were HIV/AIDS infected, their first reaction was to attempt suicide. We also state that this was due to two factors: one, the stigma that people associate with AIDS; two, because of the idea that AIDS has no cure.

Also all of them had the assumption that the only cause for the disease was commercial sex workers. All the 60 persons interviewed had been infected not from their hometown but only from outside their hometowns.

The major places of transmission of HIV/AIDS infection are from Andhra Pradesh and Mumbai. Over 85 of the interviewed members (or their husbands) had the habit of alcoholism also. 60 of them are willing to come out and talk about AIDS in front of the public and counsel others. We observed that in the case of the truck-drivers, their occupational mobility and amount of money easily earned in their hand lead to their having frequent sexual encounters, to the extent that a lot of the male patients say that the number of CSWs they visited was "numerous", "countless" etc. Most of them feared that they were ill treated by their relatives, neighbours and people in their villages. Those infected with the disease do not belong, as often misconstrued, to the poor strata of society alone. 60% of the interviewed people living with HIV/AIDS had become poor only after getting the HIV infection and spending on the treatment. Most of them belonged to the middle class or upper middle class and from the rich society.

For the purpose of this project we used an elaborate questionnaire that contained 191 questions. The questionnaire, in expert jargon, can be described as a fuzzy neutrosophic linguistic questionnaire (Appendix 1).

All the interviews were taken at the Tambaram Sanatorium from people living with HIV/AIDS. Most of those interviewed are in the advanced stages of the disease; and several of them are inpatients at the Tambaram Sanatorium.

We have given the interviews verbatim. Also accompanying the interviews are the additional observations that we, as interviewers, could gather. We have not tried even to change the sequences in which they spoke.

This is mainly done for the sake of simplicity to a reader. Here are the summaries of the interviews.

## 1

On 20-4-02 we interviewed a 25-year-old man affected by HIV/AIDS. He is uneducated. He had migrated from Andhra Pradesh and was employed as a sanitary worker in a hospital at Chennai. He is married and his wife is also affected with HIV/AIDS. They have no children. As his disease was in advanced stages only his mother who was taking care of him answered most of the questions. He was an alcoholic and he used to visit commercial sex workers (CSWs). He confesses that it was his friends who spoiled him. He visited commercial sex workers in Andhra Pradesh and has thereby acquired the disease. He belongs to the lower middle class and lives in a rented house. He claims that 75% of people are aware of HIV/ AIDS. The major symptoms which forced him to take the test for HIV were

    i.   Partial paralysis or immobility of legs and hands
   ii.   No control over urination
  iii.   Itching and scabies

His in-laws were unaware of his infection with HIV. Only his wife knows about it. He feels self-control and leading a good life alone can stop the spread of the disease. The friends who were with him, though they were addicted to drinking liquor, they did not visit commercial sex workers. He wept for not having any children and his inability to have children. He feels people of his village are reserved with him. His only wish is that the people living with HIV/AIDS should be saved.

*Additional Observation: This patient broke down during the course of the interview. He was very ill. With the exception of one or two questions which he answered in tears, his mother answered most of them. Both he and his mother are least bothered about the fact that his wife has contracted the disease from him, neither have they informed about this to their in-laws. His wife had once sadly remarked to his mother that they had wilfully ruined her life.*

## 2

A 44-year-old lorry driver who is an inpatient in the Tambaram Sanatorium was interviewed on 16-4-02. He has answered 112 questions out of the total 191. He is married and has 2 children: a boy and a girl. He is illiterate. He has not informed any of his relatives about his HIV+ status, and none of them visit him. He is depressed for he has no money even to buy and drink tea. He was quite well off and rich, but due to his addictive habits and this disease he has been reduced to poverty. His first symptom was fever and loss of weight. His wife has not acquired the disease and she was tested negative. He has visited several commercial sex workers. He maintains that the number is "countless." He has paid them from Rs.50 to Rs.100 per visit. Since he had to stay away from home for a week or a month sometimes he visited these roadside commercial sex workers. His life-story is very pathetic and heart-rending. He began his life as a child labourer. He started to earn when he was eight. His father had two wives and his mother had died when he was 3 years old. He admits he took up jobs to save himself from the cruelty of his stepmother. He adds that his stepmother had illegal relationship with another man and she swindled his father's property worth about 10 lakhs. She had beaten him with a broom when he was just five years old.

He also suffers from other Sexually Transmitted Diseases (STDs). He has informed his wife about his disease and has taken treatment for the last 10 years. He feels that his relatives dislike him. Both his children are well aware of the disease he suffers from. He advises that before the marriage both the man and woman should be tested for HIV/AIDS. He also says that people do not talk about HIV/AIDS in public as they feel it is a social stigma. He feels the advertisements in the poster form with written material are of no use to patients like him who do not know to read. He feels that to make people aware they should go from village to village and perform plays simultaneously spreading awareness of HIV/AIDS.

He says that the villagers would not allow people affected with AIDS to be in the village and they would be chased out. He has not seen any TV advertisement about HIV/AIDS. He also said that all lorry drivers were shown a video show on HIV/AIDS at Lucas TVS near 100 ft. Road in Chennai. He learnt several things from the video show, but he confesses that before the show itself, he was very infected with the disease. So when he was seeing the video show his heart was beating rapidly with fear. He has taken treatment from a doctor in Cochin by paying Rs.10, 000 but that treatment did not cure him. He had all habits like alcoholism and smoking. His employer knows about his disease and has paid him Rs.5, 000 for treatment. His only wish is to marry off his daughter without any social stigma of saying that her father is an AIDS patient.

He confesses that as soon as he came to know that he had HIV/AIDS, he tried to commit suicide. Only the patients around him in the hospital saved him. His other symptoms apart from STDs are stomach pain, headache, fever etc. He tells us that he was in the special ward for 10 days. His assessment is that no one from the hospital has ever returned home. All of them die here. He wishes to die with least suffering. He adds that he feels happy to talk with us.

***Additional Observation:*** *This man longs for affection. He is determined not to think of women badly. He is very frank, he never wanted to leave us, and spoke with us for over two hours. He wants everyone to take the HIV/AIDS test before marriage.*

### 3

An inpatient of the Tambaram hospital aged around 40 years who was a lorry driver by profession was interviewed on 20-4-02. He is from the poor strata of the society. He was suffering from fever on the day of the interview. He is educated up to the 6th standard. His wife (who is educated up to the 8th standard) is also affected by HIV/AIDS. He has widely traveled to all the states and when he was away from home he had regularly visited commercial sex workers. He has answered only very few questions. Shockingly, his main wish is that the commercial sex workers must be banned. Several of his answers were indeterminates, some were conveyed only by signs if asked for confirmation, he would answer "don't trouble me."

He says that the government should take this step and the police should cooperate for this. His parents are not aware of his infection with the disease. He has two children who have not been tested for the disease. His wife was tested and it has been confirmed that she is HIV+. He tried to commit suicide as soon as he came to know that he has HIV/AIDS.

He has not confessed or talked about his disease or his visits to CSWs to his wife. He was the only son of his parents. His marriage was arranged and he is married for the past 10 years. The commercial sex workers had advised him to use safe sex methods but only he had not followed it. The couple however seemed to be composed and the wife freely talked about the problems.

*Additional Observation: He seemed to be very reserved and tensed but later he became composed. He acknowledges about his promiscuous sexual activities in the presence of his wife.*

**4**

A 31-year-old man an inpatient at the Tambaram Sanatorium was interviewed on 20-04-2002. He is by profession a driver and is educated up to the 10th standard. His wife is also an HIV/AIDS patient undergoing treatment in the same hospital. They have two sons and who were tested for the disease; the children are not infected with HIV/AIDS.

He had taken up the job at a very young age since he had no money to continue his education. He suffers from common cold, cough and fever. He was married in 1995. He lives in poverty in a small hut. His basic symptom was breathlessness.

All his relatives are aware of his disease. He had visited commercial sex workers before and after marriage. He says that for the past 5 years he knows about AIDS. He feels sex education must be imparted to children from eighth standard onwards. He has changed his religion from Hinduism to Christianity for solace and comfort. He prays and believes in Jesus Christ.

His wife who was with him was listening to the interview but did not talk anything to us because this man was dominating her and he did not like her talking with us.

*Additional Observation: He seems to be very reserved, but accepts that he visited sex-workers before and after marriage. He feels depressed as he had got the disease due to ignorance. His wife's parents are highly hysterical and extremely affected because their daughter has become infected. They have become enraged to the extent of manhandling the patient.*

*His wife was relatively passive. She says that her husband could not study after 10th standard as he had the family responsibility of protecting his parents and he was her relative. Might be her expression showed to us that she felt that if he was educated he would have been saved from this disease and would have taken up a good Government job.*

# 5

This 30-year-old illiterate inpatient has answered around 140 of the questions. We interviewed him at the Tambaram Hospital on 19-04-02. He is by profession a weaver. He is from a remote village in Tamil Nadu. He has acquired the disease at Andhra Pradesh or Mumbai. He had visited north-Indian commercial sex workers in Mumbai and also visited Andhra Pradesh sex workers. He had done these acts as a part of enjoyment with friends. He is married and has 2 children. He was married at 22 years. For the past two years he is in the hospital. The main symptoms that he was suffering from are diarrhoea, fever, and blisters all over the body, scabies and pain in the limbs. He had regularly visited CSWs right from the age of 20. He left home after a quarrel and went to Andhra Pradesh, where he visited commercial sex workers. He says that as a result, he acquired STD and he had visited doctors and got himself cured. Then he had gone to Mumbai for four months just on a pleasure trip with his friends. There he had visited north-Indian commercial sex workers and paid them Rs.20 to Rs.50.

He had met a commercial sex worker in Mumbai whom he wanted to marry but she had refused saying that she had been cheated and she never wants to cheat anybody. She was beautiful so he was attracted towards her. At Nagari in Andhra Pradesh he had visited CSWs and paid them from Rs.10 to Rs.50. He confesses that he used to see many pornographic movies.

He never had any sexual relations within his native village. But his group of friends indulged in these activities. He said some of the commercial sex workers had visible symptoms of STDs i.e. blisters on the legs and thighs. When he had questioned one of them about it, she said that it was simply due to heat. He says he had doubt about it even then. From Mumbai he once again visited only Andhra Pradesh. He was in Andhra Pradesh for 5 years. Ultimately he again went to his hometown and got married. Actually he wanted to postpone his marriage with a view to enjoy life more. But his parents and relatives forced him into a marriage.

He is a Vanniyar by caste and a Hindu by religion. His parents and in-laws are unaware of the fact that he is suffering from HIV/AIDS. He has not even disclosed these things to his wife. He has occasionally used safe sex methods but not always. He used to feel bad while buying condoms but he used to secretly buy it. He has not used it.

His first stage of the disease was blisters and scabies all over the body. At 26 years he became aware of AIDS. He had great urge for sex and he used to take sex-enhancing tablets that cost around Rs.40 each. Even after he was aware of his disease he had sex with his wife. He has not tested whether his wife has HIV/AIDS. He belongs to the upper middle class and has the bad habits of smoking and drinking liquor. He has not informed it to his parents with a fear that they will not spend money on him. He says only fear prevents people from talking about HIV/AIDS openly. He also suggested a immediate ban on commercial sex workers and suggested that both the society and the police should work together to ban it. In his opinion, the only way to stop HIV/AIDS is to ban commercial sex work.

***Additional Observation****: This patient conceals how he got money to live alone or how he maintained himself after he ran from his house in anger. He is very careful in divulging facts to us. He mentions that the cost of the CSWs varies according to their beauty. We talked with him for over ninety minutes but unlike other patients he could hide certain facts and never contradict himself till the end. He was unrepentant and does not regret or feel guilty. He has the fear that the society will disown him. He does not bother about his twin girls (under 2 years) or his wife who is just in her early twenties.*

**6**

This 26-year-old man was interviewed on 11-03-2002, at the Tambaram Sanatorium. He says he is an agricultural coolie but after talking with him we learnt that he was carrying on a very successful paddy business once. He is married and has two children. He maintains his children and his wife are not affected by HIV/AIDS. It was surprising when he said that his father-in-law who was a lorry driver had HIV/AIDS. He further says that only his in-laws were aware of his disease and his parents are not aware of his disease.

He did not want to say his native place however by slip of tongue he said his village was near Tirupathur. He had visited commercial sex workers even before marriage with a view to enjoy life. Only here he says that he sold his house and he was doing paddy business. He does not wish to talk out or answer several of the questions. There were several indeterminate facets and nothing could be made of it. He suspects that he could have got the disease from the commercial sex workers of Andhra Pradesh. Several times in this interview he said, "don't trouble me by questioning." He was rich and reduced to being a coolie because of his habits. He acknowledges that he has spent a lot of money on his friends and on commercial sex workers. He is unaware of methods to prevent the disease. He had all sorts of addictive habits except drug addiction.

He is of the view the government's advertisements are enough for one to correct oneself. He was aware about AIDS but he had indulged in these "mistakes" carelessly.

***Additional Observation****: All the time this patient kept saying: "Don't trouble me." His only contention is that the Government has not taken proper steps to provide medicine for HIV/AIDS. He blames the government. He does not accept that HIV/AIDS is incurable. He acknowledges that his father-in-law, who was a lorry driver, died of AIDS. He even sold his house and together with his friends he spent a lot of money on commercial sex workers. He does not show any care about his wife and children. Though we wanted to find out more about his family, he was pretty reluctant and indifferent and did not disclose any thing about them.*

### 7

This 25-year man was interviewed on 1-3-2002, he is an inpatient of the Tambaram Sanatorium. He is a runaway from home, he ran away at eight years. He belongs to the Chettiar community and works as a cook. He has studied up to 4th standard. He is unmarried. Now he stays in Chennai–19. Both his parents are alive. No one in the family knows about his whereabouts or his disease.

He said he was forced to have sex with an elderly woman while he was a vendor on the train at the Nellore station. Thus, as a platform dweller he was exposed to several bad habits and hardships which was not his choice. He says he has had sex with countless number of women.

He admits that the commercial sex workers provided him with safe sex methods however he declined to use them. He admits he was unaware of HIV/AIDS. He had all habits like drinking, smoking, visiting CSWs and drugs from the age of 20. First he had come into contact with CSWs from Hyderabad. He is contradictory in this interview.

For he also states that at 12 years he started to use drugs and smoke and had acquired these habits only due to his friends. His initial symptom was scabies. He used to see pornographic movies. Only after getting drunk he approached the CSWs. He further explains that at 21 years he had visited CSWs in the Vijayawada area. Surprisingly, he expresses his intention to come out in public to spread awareness of AIDS.

***Additional Observation****: We come to know that this man was a victim of the worst forces of social exploitation. His story is terribly pathetic, he has made a living as a 'platform case' with no one to take care of him, and with all nefarious and anti-social elements of the society.*

## 8

This man is an inpatient of the Tambaram Hospital. He hails from an agriculturist family in Salem. He was quite well off but due to the failure of agriculture and the low yield he took up the profession of driver. He has studied up to the 10th standard. He is unmarried. He can drive truck, bus or lorry. He was earning around Rs.3000 per month.

He feels to stop HIV/AIDS the government should first eradicate poverty. He says this because whichever commercial sex worker he has met, they have said that they have taken up that profession mainly due to poverty to support their children and other dependents.

He is of the view that this is one of the main ways of eradicating HIV/AIDS because the CSWs are the major source of spread of this infection. He was the one to point out that all notices and posters were of no use for the majority of the HIV/AIDS patients are uneducated so they would not be in a position to even read it and gain awareness. He joined this job at 20 years. He started to take alcohol and started smoking when he took up the profession of bus driver.

First he and the conductor had sex with women who traveled from one village to another. These ladies were working as labourers in agriculture fields and traveled from one village to other. He claims that his friend the conductor is no more. He suggests that his friend could have died of AIDS. He says that he knew he was affected with HIV/AIDS in August 2001.

When asked, "How his state of mind was?" he is the only person who said that he was not sad because he was reaping the punishment for his bad deeds. He is not interested in getting married. He wants to live only up to the stage till he can be self-sufficient; once he is totally invalid he wants to be killed: that is he professes mercy killing of terminally ill patients.

Several women with whom he had sex were women whose husbands had gone outstation for work for 10 months and so these ladies wanted sex and they used people like him. They used to give liquor/arrack to these women before they had sex. He is of the opinion that the women were good and only his friend the conductor who used to be with him was not so and was already infected by HIV/AIDS. He would have also infected these women. The conductor also had the habit of smoking and was an alcoholic.

This man has also talked about their husbands to these women and some said they had no sex satisfaction from their husbands. He did not want to disclose the names of his parents. He insists that sex education must be given to school-children. He is also aware of the wiring job and his job as a driver was only temporary. As he felt that nobody should get this disease he does not want to marry. When asked about the symptoms of the disease he said not only loss of weight, but also blisters, scabies and itching.

He has never disclosed all these to anyone. Since we have asked him several times for this interview and because of our many questions he has talked his heart. He feels that when the patients become invalid i.e. when they are not able to walk and eat they should be killed (mercy killing).

When we told him that killing was not proper and was a sin, he said that Gandhi himself had said that it is better to kill than see anybody suffer and therefore it was okay for we are certain that these AIDS patients cannot be cured. He is a unique case for he has never felt sad about his disease or suffering. He has not taken any other treatment.

***Additional Observation****: Some of the messages which he conveyed are shocking. He says that apart from visiting commercial sex workers, he also had 'casual' sex with women whose husbands were away from home for long periods of time and also with housewives. He is one of the patients who feels that mercy killing is better.*

*This man even quotes Mahatma Gandhi to justify mercy killing at a stage where the affected person cannot eat or do any thing of his/her own.*

*We were surprised by the attitude of this man: he feels that he is undergoing the punishments for his sins; as a result he doesn't feel sad or bad. He is very composed and relaxed. He feels happy that his friends still visit him and they share the same cigar.*

*He states that he often laughs and he is determined not to feel sad for that happiness and laughter would cure him. He is of the opinion that he is keeping good health.*

## 9

On 19-2-02 we interviewed a 36-year-old man with HIV/AIDS who is an inpatient at the Tambaram Sanatorium. He is an agriculturist earning around Rs.80, 000 per year. He is married and is educated up to the 6th standard. He is a native of Mellur. He is a Kallar by caste.

He came to know about the Tambaram Sanatorium because his two close friends are taking treatment and two other friends had died of HIV/AIDS. Both his parents are alive and his in-laws i.e. father-in-law and mother-in-law are also alive. He belongs to the upper middle class and his brother is working abroad. He has two children aged 6½ and 4 years. He has spent lots of money on his friends. His family members do not know; only his wife has come to know of his being HIV+ from the hospital.

His wife looks at him as though he is a guilty person after she came to know that he has disease. She is not friendly with him. He started visiting CSWs at the age of 21. He has taken the test for AIDS at the age of 26, 27 and 28 because he wanted to go abroad.

When he was tested for the second time he was found to be HIV+. He could not go abroad but they have deleted his card as he was infected by HIV/AIDS. He said it is the "defect of his age" that he had such habits. He had spent around Rs.100/- for sex with commercial sex workers in Mumbai where he stayed for a year to get visa.

Even after knowing that he had the disease he had unprotected sex with his wife. His major and only symptom was itching and he at times brought out blood during coughing. He has tested his wife without her knowledge. His wife is also infected with HIV. He has seen pornographic movies especially English movies. He feels that today the government will give one conclusion and another conclusion tomorrow so he does not want to say anything. He feels specially drivers are prone to AIDS.

He has spent 8 lakhs rupees for treatment in Karaikudi and that one Dr. Marutha Pandi has cheated him. He came to know that he was affected by HIV/AIDS disease only after his marriage. He had not tested his children for AIDS as he had the disease only in the year 1999 but the children were born much before that. He came to know about safe sex methods only in the 34th year of his life. He has answered around 123 questions.

*Additional Observation: His wife does not like him now and she is very aloof and this has aggravated his mental make up and hence his physical conditions. He is taking treatment without the knowledge of anybody from his family. He has informed his family that he comes to Tambaram to visit his friends.*

## 10

On 19-3-02, we interviewed a 32-year-old man, who was in the lorry business, he had been educated up to the 10th standard. He is infected with HIV/AIDS and is taking treatment at the Tambaram Sanatorium. He is married and has a five-year-old child. Both his wife and the child were tested and are free from HIV/AIDS. He is a Mudaliar by caste.

After the disease he converted from his native religion to Christianity. He has a house of his own and he is from the upper-middle class strata. He was found to be HIV seropositive in December 2001; he further says that he did the 'mistake' around 18-19 years of age, before his marriage. His native place is Vandavasi. He contracted the disease at Villupuram; because he wanted to have a 'jolly' time. After he came to know that he had the disease he consumed 35 sleeping tablets and still he did not die. He adds that he didn't die because the sleeping tablets were past the expiration date. His main symptom was diarrhoea. He did not know about HIV/AIDS when he had sex with commercial sex workers. He expresses his idea that sex education for school-children is a must.

*Additional Observation: He feels that he liked the interview and he says the interview may be useful to the public.*

## 11

On 1-4-02 we met a 58-year-old HIV+ man staying as an inpatient in the Tambaram Hospital. His primary profession was agriculture, then he became a postman since he was an ex-service man. He is married and has 2 sons and 2 daughters. He belongs to the upper middle class society. His educational qualification is E.S.L.C. He is also suffering from tuberculosis.

He has changed his religion to Christianity for he has a complete faith that Jesus will cure and forgive him. He feels the hospital is not clean. He asks "When the government spends lakhs for cancer why not for AIDS?" He says that Hindu religious leaders (like the Christians) should give counseling and solace to AIDS patients. When Hindu leaders can give voice and support for cancer patients why not for AIDS; was the question he raised.

He feels that Andhra Pradesh CSWs have more prevalence of AIDS (though we do not know how far it is correct!) Perhaps he wants to convey that there are more commercial sex workers hailing from Andhra Pradesh as compared to Tamil Nadu. He also says that the Andhra Pradesh CSWs are very young girls and they have serious STD disease and symptoms. He had come into contact with a commercial sex worker in Pondicherry when he went there for a conference. Once again he went to a conference in Kallakurchi where he had sex with a commercial sex worker. He did it with a view to enjoy life. Also he had an illegal affair when he was in the military. He has answered over 46 questions. He has not answered the question whether his wife is affected with AIDS. The answer to this question is an indeterminate for we could not get anything out of him about his wife!

*Additional Observation: He has converted from his native religion. He drove home the point that in Hinduism when religious leaders talk for cancer, but why not for HIV/AIDS? He reports that he had unprotected sex under the influence of alcohol. Though he answered a lot of questions, he did not wish to answer any of the questions on his wife or children or parents or in-laws.*

## 12

We interviewed a 24-year-old unmarried driver on 5-3-02 who earned around Rs.2000 p.m. He is from Pasumapatti. No one else in his family is affected with HIV/AIDS. He has answered around 59 questions. His brother who is also a driver has several times advised him not to fall into bad habits. He has been in this job for the past 8 years.

When he was a cleaner itself he had a friend who was a driver who went to CSWs and induced him to do so. So without a thought he involved himself in a similar fashion. That driver friend of his is no more. So far this man has spent around two lakhs rupees for the treatment.

He had a love affair with a Kallar caste-Hindu girl but as he belonged to a lower caste they drove him away. She makes a living by working as a servant in several households. He admits that he has visited countless women because of love failure.

He has visited several sex workers. His marriage was fixed but stopped due to the disease because he felt that he should not spoil others. He does not know how much his brother has spent on him to cure this disease. His friend who has made mistakes with him is in the same hospital. He gave the ward number and the name of his friend also.

He said he would come in public and talk to people for increasing awareness of HIV/AIDS. He often felt his family has developed itself after a lot of hard work but they are spending lots of money on his treatment. He does not know to read so he did not know about the disease. When he was told about the disease he felt very sad. Also when he visited the sex workers they advised him to use safe sex methods to which he did not agree.

***Additional Observation****: He tells his story of jilted love. He continues to speak about the woman he loved. It is also a case of peer-group negative influence which has led him astray: this is best evidenced when he tells that his friend who committed the mistakes along with him is also in the same hospital.*

### 13

On 5-3-02 we interviewed an illiterate 38-year-old inpatient at the Tambaram Sanatorium. He is a native of Kallai village, Perambalur district. He is married and has switched his job from working in a country medicine shop to a shop in Burma Bazaar. He earlier worked in the herbal-leaves packing industry. Only when he was away from the family he developed several bad habits and one such was visiting commercial sex workers also. His parents are not kind with him. The in-laws look at him as if he is a criminal and therefore nobody stays with him. He adds that if his health recovers he would go back to work in Burma Bazaar and earn around Rs.200 to Rs.300 a day in business.

His family is in difficulty even to eat properly. He admits that he had sexual liaisons with women who came to work in the shop. His wife was very angry, she does not talk to him and uses disrespectful words and ill-treats him. If the government helps his children it is sufficient for him. He expresses his feeling that when he had sex with other women and with CSWs he was not aware of the consequences. He has two children.

His wife is also affected with HIV. His children are not affected by HIV/AIDS. He justifies that since he was away from his wife he had to seek sex from CSWs. He has no support. His only sorrow is about who would take care of his children. He says he thought of committing suicide. But he felt he could give the children to the government and then die. He has his own house. He knew he was being unfaithful to his wife yet as he was alone in Red hills, Chennai while working in a tea shop, he regularly sought after commercial sex workers.

***Additional Observation***: *This patient cares very much about his children. There is a very clear self-appraisal of his behaviour.*

## 14

The patient was a farmer who later became a lorry driver. He is around 30 years and has been educated up to the 10th standard. His wife is 26 years old and she has passed the 11th standard. They have two children. They belong to the Gounder caste and they are a middle-class family. The marriage was an arranged one. She was already related to her husband. His wife is also affected by the disease. His children were tested and they are not HIV positive. He did not give the name of his exact native place but simply said that he hails from a village near Coimbatore.

In the year 2001, he was tested and was found to be HIV positive. He first told his parents about his disease. He had visited commercial sex workers in Mumbai, and had spent around Rs.500 on them. He says he feels very sad when he sees his children and wife and repents for his actions.

When he was 25 years he knew about HIV/AIDS. The sex workers offered him condoms. He said that he has used safe sex methods yet he is surprised as to how he got the disease. Most of the CSWs had said to him that they did sex-work for the sake of livelihood. He feels free and compulsory screening for AIDS must be done and sex education can be started at school in the 8th standard. Only his friends support him. His friend is also affected with AIDS and is in the same hospital with stomach pain and diarrhoea. Both mother and father-in-law know about this and they treat him kindly. He is an alcoholic and smoker. He is willing to spread awareness of HIV/AIDS in the public.

*Additional Observation: This man says that everything is immaterial now because he has become a victim. He feels no one should be affected the way he has been affected. He does not feel sad about his disease.*

## 15

On 10-2-02 we interviewed a 35-year-old driver living with HIV/AIDS. He is married and is educated up to the 7th standard. His wife is 28 years old and is educated up to the 10th standard. He is an alcoholic and smoker. He is also affected by tuberculosis. He first had sex with commercial sex workers when he was about 20 years. He was unaware of AIDS when he visited commercial sex workers. His marriage was an arranged one.

He was very reluctant to speak and did not even answer the question whether his wife was tested for HIV/AIDS. He was silent and did not respond to several questions and the answer for them remained indeterminate.

*Additional Observation: He has four children: 3 females and a male. He says he had sex not only with CSWs but also with family women. He did not disclose how much he spent on them. When he sought after the CSWs, he did not know about HIV/AIDS but he was only aware of VD. After he acquired this disease only he knew "that HIV/AIDS was caused by a poisonous germ" and he says "there is no way to kill it or cure it."*

*He did not know the difference between HIV and AIDS but he says it comes because of 'ladies'. He tells us some of his friends are also ill with HIV/AIDS and they take medicines regularly. He advises that youngsters should not seek after CSWs and that is the only way to be free from HIV/AIDS.*

## 16

On 10-2-02 we met a 28-year-old driver in the Tambaram Hospital who was affected by AIDS. He said he is a smoker but he is not an alcoholic because he drinks only beer.

He is from the Pillai (Hindu) caste. He is a tuberculosis patient and he says he weighs only 30 kg. He has completed primary schooling. He was unable to walk when we met him. He had been married just 6 months ago. His marriage was an arranged one. He knew about HIV/AIDS.

He did not answer many of our questions. He says that he has visited his friends in this hospital who were affected by AIDS. He said he does not know how the disease is spread. He acknowledged that he got the disease because of going as a lorry driver to different states and distant places.

*Additional Observation: He is aware and knows about HIV/AIDS. He has come to visit his friends who were admitted here and thereby he had learnt about HIV/AIDS. He also told us that his friends got the infection because they had visited CSWs but he does not know how he got the disease: he is unique in a way for he does not even accept his visit to the CSW. He says "one gets the disease through ladies (CSWs)."*

*He was saying discordantly that the mistake-doers are the politicians! We could not get the clarification for this statement.*

### 17

We met an inpatient in the Tambaram Sanatorium aged around 41 years who was a driver. He has not passed SSLC. His native place is Bangalore. He was a driver for the past 25 years and is a Hindu by religion. He said he got the disease because of his insolence.

He admitted to being addicted to alcohol and smoking. He said he had visited commercial sex workers and felt very ashamed to speak about them. Only his wife knows about it. He says he was like a tiger in those days, but now he has become like a dog.

He is very sad. His wife wept very bitterly. He did not want to talk more about this.

*Additional Observation: This man has taken up the job of driving at the very early age of 16 itself. He is highly depressed and sensitive. He says he got the disease because of his arrogance.*

*From the age of 25 he was a regular visitor of CSWs. His wife says it is her destiny to have a HIV/AIDS affected husband. He has not disclosed about the disease to others for the main reason that his daughter's life would be affected. He says even the neighbouring homes do not know about this.*

## 18

This man did not even disclose his name. He is 27 years and has completed SSLC. A driver by profession, he is suffering from tuberculosis also.

He acquired the disease in Mumbai and had spent around Rs.60 per CSW. He felt like committing suicide when he came to know he was affected by HIV/AIDS. He is married and has not tested his wife for HIV/AIDS. He said he couldn't say about the disease in his village for people in the village will not give even water. His village people are not aware of any safe sex methods or the cause for HIV/AIDS. He feels he cannot even say about the disease to his friends.

*Additional Observation: He is of the opinion that he got the disease because he committed mistakes. He does not disclose what were the 'mistakes'.*

*He has a child. He does not believe as much in God as he believes in Medicine. He is also of the opinion that his disease cannot be disclosed to friends. Further he maintains that he is unaware of safe-sex methods.*

*He holds the opinion that to live with some self-respect it is better he does not disclose about the disease to anyone: a clear indication of his stigmatized identity within the community.*

## 19

This 24-year-old inpatient was interviewed on 3-4-02. He is a lorry driver and also owns 3 lorries. He had sex with commercial sex workers right from 20 years of age. He is unmarried. He is not able to distinguish between HIV and AIDS.

He belongs to Karoor. He acknowledges that he got the disease from CSWs. He was unaware about HIV/AIDS.

He repents because he is not even able to walk now. He has spent around Rs.200/- for each visit to these commercial sex workers. He says his family people know about the disease. He feels he did all this only as a "jolly enjoyment." In the end he admits to having visited the commercial sex worker even at the age of 18. He is contradictory.

He claims the disease is only due to CSW and he repents desperately for ruining his life by deluding himself that he was enjoying it.

*Additional Observation: He thinks he can cure himself completely and go home. He is not even in a position to walk. It is surprising that he has used safe sex methods and he does not know how he has got this disease. He says "men are affected faster than women." He has lost weight drastically. He shouts in anguish "CSWs must be asked to work as coolies and their profession must be banned!"*

## 20

We met a 28-year-old man, a native of Thiruvannamalai who is an inpatient at the Tambaram Sanatorium. He is a lorry driver by profession and is unmarried. He had visited CSWs right from the age of 23. He is unaware of HIV/AIDS. He did not use safe sex methods though he was aware of condoms. He has spent Rs.30 to Rs.40 on CSWs. He has visited around 50 to 100 commercial sex workers. He had the symptoms of STDs and met a doctor and was cured. He says as the disease spreads because of the commercial sex workers they must be banned. He feels everybody thinks it is very indecent to talk about the disease and socially it is impossible to accept it.

*Additional Observation: It took six months to know that he had HIV/AIDS. His family members wept when they came to know he had HIV/AIDS. He had STDs even before he was confirmed to be HIV+. He says the CSWs should not exist. He feels the government must find the medicine to cure the disease. He feels that his relatives are unchanged even after the disease. He also says that in the Tambaram Sanatorium, Siddha medicine is also given. He also wanted free medicines to be distributed to all AIDS patients.*

## 21

A 47-year-old man living with HIV/AIDS was interviewed on 5-03-02. He is a bus-driver by profession. He took up this job at the age of 18. He has 3 children. From the age of 18 he has had sex. He says he was not aware about HIV/AIDS. At the age of 20 he married, it was not a love marriage. He has sought after addictive habits like drugs and drink. He hails from the Manaparai Town.

The major symptoms were ulcers in the mouth, diarrhoea etc. After a long time his family doctor advised him to go for the HIV/AIDS test. He feels now that his living is of no use and wonders how he is going to earn and live. He does not wish to live. He said he has no money even to buy and eat anything or even to have a cup of tea. He says he has not wasted any of his wife's property but he has wasted the money he earned as tips. He says he was not aware of any safe sex methods.

*Additional Observation: He said that today's youth are well versed with the knowledge of computers and that they should be given proper sex education. He says none of his relatives know about the fact that he is suffering from AIDS. He is in a way unique for he says a man can have any number of women provided the women do not have any disease and the man is also free of any disease. We are not able to find any other facts about his real state of mind except that he feels insecure and dejected.*

## 22

On 2-4-02 we met a 34-year-old agricultural labourer from Tannampadi. He is educated up to 10th standard. He had acquired the disease when he was working in a company. His symptoms were fever and breathlessness; he also had tuberculosis. He is married and has a child aged 2 years. The child and his wife are also HIV/AIDS infected. They have been tested to be HIV+.

He first said he had affair only with family women. Then after some time he confesses that he had visited commercial sex workers and spent Rs.150 for every visit. He then acknowledges that after a year he had sex with a CSW he found changes in his body. He admits that even before marriage he had sex with CSWs. He has read pornographic books. He is not aware of the difference between HIV and AIDS. For the past two months he is immobile i.e. not even able to walk or sit.

*Additional Observation: He hails from a joint family. He says his father had two wives: one wife is dead and another is living. He has not disclosed about the disease to anyone for the fear of ill-treatment. None of his relatives know about the disease. He has not only spent on CSW but he says he had affairs with other women also. He feels dejected for getting the disease because of his own bad activities. The minute he came to know that he is infected by HIV/AIDS, his first question was whether the disease is curable or not. We came to know during the course of the interview, that he has been bedridden for the past two months,*

*even now he is immobile. He never spoke about the treatment given to his children and wife.*

## 23

We met a HIV/AIDS affected man aged about 26 years in the Tambaram Sanatorium. He has studied up to S.S.L.C. and he is working as a salesman in Chennai in Old Washermanpet. He is not married. His father died when he was just 6 years old. He is suffering with HIV/AIDS for the past four years. He is now affected by tuberculosis also. He is a native of Pondicherry but was born and brought up in Chennai. He tells that from the age of 23 years he started to visit commercial sex workers. He maintains that he was well aware of the fact that the disease spreads through CSWs. He was cautious and guarded in divulging information to us. He admits that he was a regular customer at a lodge in Thiruthani that provided commercial sex workers. He has seen all types of movies and only recently he has stopped. Earlier he used to smoke. He feels he is better and improved now. He has never visited quacks or other doctors for treatment.

He had a set of friends and enjoyed life in this manner. He does not like relatives, even his own brothers and sisters and says they are of no use. He tells us that he will now advice his friends to go for the HIV/AIDS test as they have also visited the same set of commercial sex workers. He feels that he used to select commercials sex workers with beautiful face and body. He later admits that this was the mistake he has made.

He also told us that if he doubts these commercial sex workers about having any disease they would say, "you need not visit us." They have said to him that they did not have money to go and test for the disease at this juncture. He feels that the government when it permits the CSWs to flourish it should take the responsibility to test them and provide the information to them.

We were able to get the true feelings of this patient and we came to know of how he had contracted the disease. Each time, he used to spend from Rs.200- Rs.300 for the sex-workers.

## 24

We interviewed a lorry driver aged 46 years (an inpatient of Tambaram Hospital). He is from Ottanchathiram in the Dindigul district. He used to drive from Orissa border to Kerala and return back. He openly confessed that he had visited commercial sex workers only in Andhra Pradesh on his way to Kerala. They had to spend Rs.50 on these women who are aged between 18 to 20 years.

He is married and has a son. First when he came to know he had HIV/AIDS he was shocked but he says he consoled himself saying that after all one day we are going to die. He feels that he is infected for the past one year. His wife is also infected and she has HIV/AIDS. But his son is not infected.

He has taken the Siddha treatment from Kerala and feels it is of no use. His friend who is a driver is also affected with HIV/AIDS. He has recently changed his religion to Christianity, for the past two months he is a Christian.

He proudly shows that he is wearing a cross and Jesus would cure him. He is educated up to the 5th standard. He was born as a twin and his twin-brother owns lorries. None of his brothers visit him. He says he was well aware that because he was a lorry driver he would get this disease. He also had tuberculosis. He has all habits like alcoholism, drugs, smoking and visiting CSWs.

## 25

We interviewed a 42-year-old man who is an inpatient of the Tambaram Sanatorium. He is married and he belongs to the middle class. He says his children are not having this disease. His primary symptoms were cold and fever. All his relatives look at him with contempt and hatred and he feels bad about this.

He first came into contact with sex workers when he stayed away from his village. Later he had worked in a petty tea-shop in Chennai and stayed there. He used to eat and sleep in the same shop. He was an agriculturist labourer, because of the failure of agriculture, he had come to Chennai in search of a livelihood.

*Additional Observation: He is from a joint family: all of them used to go and work as hired-labourers in agricultural fields. Now they don't get any job, so he came alone to Chennai in search of livelihood. He says he contracted the disease only in Chennai.*

*He further adds that only a friend from a neighbouring village had brought him to a CSW in Chennai. He was unaware of safe-sex methods. Because of more leisure and easy money he had become a victim of HIV/AIDS. He has no education, his*

*profession and the absence of social responsibility are the reasons that led him to acquire HIV/AIDS.*

## 26

On 2-3-02 we interviewed a 34-year-old rich Muslim man from Chengelput who is an illiterate agriculturist. He is married and has a son. All his brothers and sisters are agriculturists. He was unaware of HIV/AIDS. He bluntly acknowledges that he had visited commercial sex workers before and after marriage. He has told his wife about his activities only after he acquired HIV/AIDS. She is also infected and has lost weight. He underwent a test for the disease in Mumbai. Only then he came to this hospital. He had flu and tuberculosis. He confesses that only because of him his wife had contracted the disease. They both have skin lesions and scabies on and off. He has seen pornographic magazines. The commercial sex workers whom he met in Mumbai were from Tamil Nadu, he also says his friends have the disease and most of his friends with AIDS have died. He wishes for a free test for HIV/AIDS by the government.

*Additional Observation: Nine years after his marriage, he had contracted the disease from the CSW. He says that he went to them, for, at times his wife has refused to have sex with him. He has been to Mumbai and visited CSWs even there. He has used safe sex methods yet he does not know how he has contracted the disease. According to him, in Tamil Nadu the Tambaram Sanatorium is the only good hospital for HIV/AIDS patients. He is of the opinion that several people will not give even water to drink to a person if they come to know that he/she is an AIDS patient. He has not confided about the disease to any of his friends, only very close relatives know about it. He also adds that he is not alcoholic.*

## 27

On 7-3-02 a 34-year-old man living with HIV/AIDS and an inpatient of the Tambaram Sanatorium was interviewed. He has completed his SSLC and he belongs to Pattukottai, and he is a Hindu by religion. He is married and his wife is also educated up to the school final. They have 3 children. He knows two languages: Hindi and Tamil.

He is from a middle class family and has an own house. He knows that he is infected for the past five years. He has concealed this to his wife. He acquired the infection during a visit to Mumbai where he had gone with his friends. His wife is tested and was found to be HIV negative. The children were all tested to be negative.

He says that HIV/AIDS spreads only through women. He had spent money and had sex with over a dozen women in Mumbai. His only symptom was fever.

*Additional Observation: He had got this disease in Mumbai five years ago and he has not yet disclosed this to his wife. He was found to have HIV in Pattukottai Hospital. In his own words, he had "enjoyed" life in Mumbai.*

*His friends are well and unaffected by this disease. For several questions he simply sighed, that was the only answer. However he says HIV/AIDS affected children can be accommodated in a separate hostel! This was surprising. He had gone to Mumbai in order to go to a foreign country. He says that he was basically an agriculturist coolie. He was later working as a cycle mechanic. He too wants to ban commercial sex work.*

## 28

On 7-3-02 a 40-year-old man living with HIV/AIDS, (an inpatient from the Tambaram hospital); who was SSLC qualified, driver by profession was interviewed. He is a Hindu Gounder from Pennagaram town in Dharmapuri district. He has three children. He has sold all his property. He can talk Hindi, Tamil, Kannada and Telugu. He is an alcoholic and is all alone now. He had these habits when he was a driver and even before his marriage. He was tested to be HIV+ only recently. He has spent a lot of money as doctors said he was suffering from brain fever. He hails from an agriculturist family.

*Additional Observation: He basically hails from an agriculturist family. He was treated for 10 days in a private hospital where they said that he was suffering from brain fever. Since he was not cured he went to Salem. In Salem also he was given glucose and they also said that he was suffering from brain-fever. Then he went to Bangalore. A doctor at Bangalore, treated him for ten days and asked him to get treatment from Tambaram. He confesses that he does not know anything about HIV/AIDS. In his*

*opinion, no lorry driver visits CSW without condoms. He thinks in a drunken state he might have forgotten to use it. He has not yet tested his wife and children for HIV/AIDS. He says earlier he was careless only now he repents for it. He complements the government especially the Tambaram Sanatorium for the medical advancement.*

<div align="center">

## 29

</div>

A 31-year-old man living with HIV/AIDS from Tiruchi and now admitted as an inpatient of the Tambaram Sanatorium was interviewed. He has studied up to the ninth standard and is by profession a lorry driver. His hereditary occupation is agriculture. He is married and has two children. He openly acknowledges that he has visited commercial sex workers when he was driving lorries and tractors. He feels he had made the mistake of not utilizing safe sex methods. He belongs to the Chettiar community.

He says he was aware of safe sex methods even around the age of 10 to 15 years. His wife or children have not taken the test for HIV/AIDS. He lives in his own house. He has lost weight from 69 kgs to 50 kgs. He says he had visited commercial sex workers in Andhra Pradesh, Mumbai, Delhi, Assam, Tiruchi and Kerala. He did not answer most of the questions.

*Additional Observation: This man has two brothers who are agriculturists. After knowing that he has acquired the disease they hardly visit him. He also adds that his mother-in-law is very affectionate with him. It is unfortunate that he has so far not taken any steps to test his wife or children for the disease. In Tiruchi he was tested and the doctor immediately advised him to go to Tambaram and take treatment.*

*In his opinion, several of his driver friends who have HIV/AIDS say that they are suffering from tuberculosis etc. but they do not disclose the truth. He has contracted the disease from Andhra Pradesh (we learn further by his statements that the CSWs are very cheap in this region and that they offer their services for Rs. 20 to Rs. 50. They are also aged between 16 to 19 years.*

*He also adds that most of these women suffer from STD). He adds that he has visited CSWs in Delhi, Mumbai, Assam and Kerala. He first disclosed that he had HIV/AIDS to his mother. He has taken Siddha treatment also. He has not seen the advertisements on TV. He requests the government to help the*

*AIDS patients and their families and give free medicine to them and free education to their children.*

## 30

A driver aged about 46 years, now an inpatient of the Tambaram hospital was interviewed. According to him, ignorance was the cause for HIV/AIDS. He refused to speak about his children, parents or wife. He had visited commercial sex workers while he went on outstation trips and he spent from Rs.50 to Rs.100 on these CSWs per visit. He did not divulge the name of his native place. For the past one year he is in the hospital. He has seen pornographic magazines also.

Now he has changed his religion to Christianity as none of the native Gods like Mari Aatha, Kali Amman cured him, he mainly believes that Jesus will at least forgive him though not cure him!

***Additional Observation**: Even his brothers do not visit him now. His friends are free from the disease. He did not disclose the states to which he traveled and where he might have acquired the disease but he acknowledges that he has visited CSWs. He is also an alcoholic. He answered only 17 questions. He has even refused to answer the questions like "Are you married? Father's name? Village name? Hereditary occupation?" etc.*

## 31

An illiterate fabrication worker aged 35 years was interviewed on 18-4-02. He is from Ambattur, Chennai, and he is a Nadar by caste. He has 3 daughters. His marriage was an arranged one. Before 5 years he had some blisters in private parts, he visited doctors and was cured. After this he visited several commercial sex workers from Andhra Pradesh, Mumbai, Goa etc. Also his mouth was ulcerated so could not eat and he also had a TB attack. His children are healthy and his wife is well but does not visit him. He states that when the blisters and scabies did not heal he doubted that he might have HIV/AIDS. So at that time itself he made up his mind that he would have HIV/AIDS.

He feels that because he drank and he did all 'bad' acts without rhyme or reason he is sick. He also thinks that it is a social stigma to suffer from this disease. He told that the free test by the government to detect HIV/AIDS is welcome. He is of the opinion that pornographic cinema and lewd pictures are also one of the reasons for the spread of HIV/AIDS. He says the government and the police are supporting commercial sex workers; hence the spread of the disease cannot be stopped. He admits to being well off.

*Additional Observation: He feels that some doctors are kind and considerate and a few of them are irritated and rude. Whenever his wife went to her mother's house, he used to go with his friends to drink and enjoy life. He says those friends are healthy and they do not visit him. He took Siddha medicine and his condition became very serious. He is of the opinion that most of the patients when they are little better in health they drink and have sex and that is why they become more serious. He has got this information from the other patients. He reports that if the patients obey the doctor's advice they will not become more serious. He is of the view that all CSWs must be tested for the disease and they should be given free medicine. He says that in Parry's Corner in Chennai, women from Kerala, Karnataka and Andhra Pradesh roam in the nights to get customers. These women are not suspected even by the police. They are hand in glove with them. In his opinion, these women are the main cause for the disease in Tamil Nadu, they fall on the customer and seduce him. He reasons that as men are drunk most of the nights and so became a easy prey by not only losing their belongings but by acquiring the disease. According to him, AIDS spreads because of pornographic movies and easily available CSWs (who just charge from Rs.20 to Rs. 50).*

## 32

On 13-2-02 a 38-year-old man, an inpatient of the Tambaram Sanatorium living with HIV/AIDS was interviewed. He is from Tirunelveli, and is a Nadar by caste. He is unmarried and has worked in Mumbai for the past 20 years. He was already taking treatment in Bangalore. He felt little better and left for Mumbai.

He took treatment in Mumbai also, which he says was not effective. He admits to having visited commercial sex workers. He is educated up to the 8th standard. He started such activities at the age of 20. He says "if you are in Mumbai you will know everything." But he says he knew about HIV/AIDS and also tells us that the disease is incurable.

He confesses that he knew he would be affected by HIV/AIDS if he visited commercial sex workers; still he went. He is of the opinion that he was well when he did exercise and had good food, but once again he took to the same habit and he felt his health condition became from bad to worse. He has scabies and blisters in the hands and legs. When he had a tuberculosis attack four years ago, it took six months for him to be cured.

***Additional Observation****: He is the only one who says he had visited CSWs in spite of knowing fully well he would get HIV/AIDS.*

## 33

On 25-4-02 we met an illiterate, 27-year-old who was an inpatient from the Tambaram Hospital. His job is embroidery. He is married for the past 6 years and has 2 girl children. He is from the Tiruvallur district. When his wife came to know about his HIV infection she left him saying that after becoming well he could join her. He is taking care of himself and also adds that they have given him a separate plate, tumblers and place. He bluntly says he does not know how HIV/AIDS spreads. He is highly dejected and says it is better to die. He visited Mumbai on a business trip and he acquired the disease only in Mumbai from the commercial sex workers. He has had sex with commercial sex workers from the age of 15.

***Additional Observation****: This is the only case were the wife has left the husband and gone to her mother's house when she came to know that he has HIV/AIDS.*

*He is unaware of how HIV/AIDS spreads. He has also contracted the disease in Mumbai. He reports to us that all the CSWs whom he visited gave him condom but he only did not use it for he did not think he would get the disease.*

## 34

This 47 year-old inpatient of the Tambaram Sanatorium is a lorry driver by profession: he hails from Coimbatore. He is married but he has no children and is very poor. He is not aware of this disease but he says that he had blisters and scabies in the private parts right from the age of 15. None of his relatives except his wife and brother-in-law know about it. Both of them are very kind and caring with him. When he went on long trips he visited commercial sex workers.

*Additional Observation: He confesses that he had sex with several CSWs. He says that he always had money to visit CSWs though he was poor. He has taken country medicine for 60 days but he is of the opinion that it was of no use for it did not cure him.*

## 35

This inpatient of the Tambaram Sanatorium aged 32 years is a farmer living with HIV/AIDS. He has only studied up to the 7th standard. He has visited Mumbai, Delhi and Goa. He was taking treatment in Kerala and has spent lots of money. He is not married. He is a smoker and an alcoholic. He is from Chettiar Thoppu near Chidambaram. He says he did not get married because he did not want to infect any woman. He has been treated before also for the same disease. Earlier he had faith in God but as it did not make him feel better now he has no faith in God. He is the only one in his family and village who is affected by HIV/AIDS and every relative knows about this in his village. People know about it: some consoled him and others look at him badly. He had visited commercial sex workers right from the age of 20. He accepts that his age was the reason for his behaviour.

He had gone to the sex workers only from the neighboring town or from other place i.e., not from his own place. He knows about safe sex methods yet he was not bothered about it. He says in Goa he spent Rs.100 on each CSW he had visited and that in Mumbai the rate is according to her beauty. Further he adds that all types of CSWs were available in Goa. He is an alcoholic and smoker. He feels he has to die. None of his friends have this disease.

*Additional Observation: First when he came to know about the disease he thought he could commit suicide. He cries as he utters this. He says the number of commercial sex workers he has visited is numerous. He further adds if he had been aware of the disease he would have been very cautious.*

## 36

On 3-4-02 we interviewed a 34-year-old inpatient from the Tambaram hospital. He is married and his wife is no more. He says that his wife also had sexually transmitted diseases (STDs). He feels she might have died of some venereal disease. He has two sons. They belong to an agriculturist family. His father had deserted his wife shortly after this man's birth. This man's mother visits him.

He accepts he has had sex with many CSWs, for which he has no account. He has gone to several places in different states selling timber and on his way i.e., on roadsides he has had sex with commercial sex workers. He acknowledges that he knows about safe-sex methods but he did not use them for he was not comfortable with condoms. He is uneducated. He feels people do not talk about the disease for the fear of social stigma. He accepts that he did not feel sad when he came to know he had the disease. He tells us that both his children are healthy and that his wife would have contracted the disease from him and that could be the reason why she died.

***Additional Observation****: This is an instance of a broken family for his father is a Vanniyar and his mother a Dalit. After the birth of him they separated and his mother alone brought him up. Only his mother sees him, father does not know about his disease and he does not wish to disclose this to his father. He feels on Sundays street dramas can be performed about HIV/AIDS for awareness. He accepts that soon after the infection, the affected do not take treatment but instead they go on seeking astrological advice, star's position, other superstitions, which makes the disease advance and leads to chronic stages.*

## 37

On 27-4-02 we met a 50 year-old illiterate man affected by HIV/AIDS. The wife is aged about 45 years and her husband is 50 years. They are illiterate. They belong to Pachaiyur near Krishnagiri. He has 3 children. His wife is also HIV+.

The woman has no close relatives. First she had fever and was very serious so she was tested. She was found to have HIV/AIDS so the doctors advised her to test for her husband and he too was found to be HIV seropostive. He has tuberculosis and throat infection. She feels sad about their plight. He has left his job for the last 4 years.

He used to stay away from home to sell clothes and at that time he had visited commercial sex workers. She is unaware of the difference between HIV/AIDS and does not know how to prevent the disease. She has spent over Rs.10,000/- to cure him. He did not talk much or even show any interest in this interview.

*Additional Observation: He did not go for any work for the past 4 years; as he was unable to work and never used to eat well and often was vomiting. His wife had spent lot of money on treatment. She had an attack of fever which could not be cured. So she was advised by doctors for blood test and was found to be HIV+. Then as per the request of the doctors her husband was also tested and was found to have HIV/AIDS. She was dejected and has submitted to fate. Some of the villagers in her village are also affected by this disease but none of them disclose it for fear of isolation and social stigma.*

## 38

On 18-4-02 we interviewed a 26-year-old man and his 22-year-old wife. Both are affected by HIV/AIDS. She is by profession a farmer and has two children. Her husband is HIV/AIDS infected, and is affected with tuberculosis and is admitted in the same hospital. She acknowledges that her husband had stayed away from her. We had interviewed the couple simultaneously.

Being an agriculturist labourer she may get once a while a job which was never certain. Both of them hail from Karur. As soon as he came to know he was affected with HIV/AIDS he wanted to consume rat poison and die. Their marriage was an arranged one, and she is related to him. They are economically backward. He adds that he had gone to Andhra Pradesh and only there he had contracted the disease from CSWs.

Also he had an affair with a married woman from his place who belonged to the Chettiar caste. He had sex with Telugu women and one of his friends from Erode had introduced him to these women.

He admits that he sees movies and that these movies kindle base emotions. They go for movies with CSWs first, then only take separate room in nearby lodges according to their social status. He emphasizes that he has never visited his village CSWs. When asked about his symptoms he explains that he had often experienced burning sensation in the private parts.

*Additional Observation: They are from an agriculturist family. She expects that the government must cure her and her husband*

*by taking care of their medical expenses. The only cause for HIV/AIDS according to her is that family people have illegal sex with others. She says her only support (physical and monetary) is her mother. Her father-in-law does not know about their disease.*

*Her husband says he wanted to hang himself or consume rat poison. He has taken treatment for STD/VD from several doctors. He was accompanied with 3 friends whenever he visited CSWs. He says they used to take the CSW with them to cinema-theater after which they booked a room in the lodge.*

*He sadly adds that he has not seen his three friends after he had been admitted in this hospital. Even before his marriage he had sex with CSWs because of friends.*

## 39

We interviewed a 35-year-old man from Chidambaram who was a cook by profession and an outpatient of the Tambaram Sanatorium. When he had migrated to Mayawaram on a business trip he had married a widowed woman from the Naicker (Vanniyar) community, she was already having a child.

Often the child was very sick, when they tested the child it had HIV/AIDS and presently the child's health condition is very serious. His wife also died in 1999 in the General Hospital of HIV/AIDS. He got married in 1996.

He has two children: one child died just after birth and the other is the child of his wife's first husband. His wife was an assistant to cooks, she used to wash vessels etc. The man says that she had relationship with men not for money but mainly to have some 'manly' support, but they had cheated her.

He accepts that she had enormous faith in him and he supported her till the end. He adds that the priests in the church of the village also requested him to marry her saying that she was good. She has not disclosed much of this to him. Only her parents have told him these things.

*Additional Observation: This is the only case in our sample study where the man has got the infection from his wife. He had never visited CSWs and does not have any other habits. He is special in a way, for he says in marriage more than seeing the horoscope one must first perform HIV/AIDS test.*

## 40

On 18-3-02 we interviewed a 41-year-old inpatient educated up to the fourth standard. He is married and they have 3 children. His occupation is agriculture. He belongs to Vellore. According to him, when he went to Nellore he had sex with commercial sex workers. He has paid about Rs.500 for that visit. He knows both Tamil and Telugu. He feels it is a social stigma to talk about HIV/AIDS to others. He feels more advertisements must be given to spread awareness of HIV/AIDS. He admits that his wife has not taken the HIV/AIDS test so far.

*Additional Observation: He was aware of condoms he only has made the mistake by not using it. Further the CSWs whom he visited did not offer him any safe-sex methods. It is observed that Andhra Pradesh based CSWs are not as aware as those in Mumbai (who first and foremost insist on condom usage.) He had an upper middle class background.*

*He says earlier he used to get lot of agricultural yield and good profit as an agriculturist, but nowadays the major amount is spent as pay for coolies. His blood test for AIDS was first performed at Tirupathi. He admits that this was mainly done to keep this disease a secret from his family members and friends. The doctors who treated him said they would cure him but they did not cure him. He adds that he has seen HIV/AIDS advertisements in TV. But still he is not aware of how this disease spreads.*

*None of his relatives visit him. Another point is, after he got the disease he did not have sex with his wife.*

## 41

We interviewed one of two brothers living with HIV/AIDS on 18-2-02, at the Tambaram Sanatorium. They are educated up to 10th standard. They are Hindu Naidus living in Ulundhukottai. They divulged that they were infected because they had unprotected sex with commercial sex workers. First, when they were tested they had only VD, they were careless and now both of them are affected by HIV/AIDS. One of the brothers acknowledges that he has tested his wife and she is suffering from VD. He says first that his wife was very ill with some symptoms like white discharge, then they tested themselves. He blurts out, "She knew the loose character of me yet she married me." He has visited only the commercial sex workers from Nellore when he was hardly 16 years old. He is unaware of how HIV spreads. He acknowledges that he has spent from Rs.5 to Rs.15 then, but recently he had to spend Rs.50 on the commercial sex workers. He adds that to prevent HIV/AIDS the CSWs must be fully banned. He accepts that his brother who has always accompanied him in all acts was losing weight badly so only he was also tested based on this suspicion. He confesses that he and his brother did all acts together that is why both of them in the family are affected by AIDS.

***Additional Observation***: *According to him, CSWs must be rehabilitated with work according to their ability. Three members in this family are affected: this man, his brother and his wife. They have only one child, they have not tested for HIV/AIDS for their child. Their parents are not informed, as they do not like to inform them. He adds that his brother has contracted the disease because after this man had sex with a CSW, his brother would have sex with the same CSW. They had sex even when they were in 10th standard.*

## 42

On 19-8-2002 we interviewed an uneducated farmer aged 35 years. He is a Hindu Gounder from Villupuram. He earned from Rs.20,000 to 30,000 earlier, these days the yield is low and he earns much lesser. He knows Tamil, Telugu and Hindi. He first left his village and came to Chennai and was doing some daily labour for 6 days then he went to Mumbai and spent some days there. Then he visited Bangalore for 15 days. He first said that he had never visited commercial sex workers or had no habits except alcoholism. Also he said because he had tuberculosis he would have got HIV/AIDS. After talking with him he says his wife is not affected by HIV/AIDS and acknowledges that he has visited CSWs first in Mumbai then in Bangalore and later confessed he got AIDS only through them.

***Additional Observation:*** *He is an agriculturist coolie but later he had become a farmer by purchasing 2 acres of land. His hereditary occupation is agriculture. In his opinion because he had a love failure in his life and his wife was not up to his satisfaction. He adds that some of the doctors are very irritated with him. He is of the view that except this hospital there is no* *place where HIV/AIDS treatment is properly given so he feels this is the only best hospital. He has not disclosed this disease to his village people or to his relatives. He says that the medicine which they give are of value from Rs.1000 to Rs.2000, he feels very sad because some of the patients take these medicines which is given free and throw it away. Some of them do not take the medicine regularly.*

*He says some CSWs are really in a pitiable state whereas some of them do their job mainly to earn more money easily. He is absolutely unaware of any safe sex methods. He emphasizes, "Children should be openly taught that by illegal sex they would get HIV/AIDS. Only then they will be careful." He said that both his legs and hands were paralyzed that is why he was first hospitalized.*

### 43

On 8-2-2002 we met a 37-year-old man who is now an inpatient of the Tambaram hospital. He is unmarried. He has studied up to 6th standard. Only after taking medicine from this hospital he is in a position to walk a bit. He is unaware of HIV/AIDS. His mother took him to a village doctor but it was no use.

He had contracted the disease when he went to Namakkal to buy a cow; he had visited commercial sex workers in Namakkal. Only once he had visited them. He is a milk vendor by profession living in North Chennai, a native of Pangalam-North.

*Additional Observation: It is felt from his talk that if this man had known that Namakkal is a place where HIV/AIDS is highly prevalent he would not have visited CSWs at Namakkal. He said they supply CSWs in the lodge in Namakkal and he has to only pay the rent for the room. He confesses he was drunk that is why this problem. He is unaware of anything about AIDS. He is a bachelor and when he had the disease he used country medicine for cure. He has hope that if he takes the medicine regularly from the Tambaram Sanatorium he would be cured.*

### 44

A 26-year-old illiterate man who is a 'Thavil' artist by profession from Mayawaram was interviewed by us in February 2002. He was married and on the first night itself because of some problems with the bride he left home to perform a Kutcheri and after finishing the performance he and his friends visited commercial sex workers and he feels only there he would have contracted the disease. He is a Pillai by caste. He was very sad when he came to know that he had HIV/AIDS. From then onwards whenever he went out to give performances he visited commercial sex workers by paying them from Rs.20 to Rs. 30 per visit. He started this habit at the age of 22. He is unaware about AIDS.

*Additional Observation: He says because of his stammering and baldness, his marriage was a failure even on the first day. As he does not know to cook he will soon get married again.*

## 45

On 6-2-02 we met a 35-year-old lorry driver at the Tambaram Sanatorium. He is married and is a Christian. He hails from Thirunelveli district in Tamil Nadu. He was tested to have HIV/AIDS in the Palayamkottai General Hospital. He has no children. As a driver he has traveled widely and had visited CSWs from all states. He feels that he has got the disease because of his extreme bad habits and uncontrollable activities. He answered only 16 questions. So no special facts can be derived for he did not even answer the question whether he has tested his wife.

## 46

A Naidu man aged 30 and an inpatient of the Tambaram hospital was interviewed on 8-3-02. He is unmarried. He has visited CSWs in Padappai and Navaloor and spent Rs.100 and Rs.200 for each visit. No one in his family knows about this, as he is the head of the family. He feels he does not need any government support. He does not see movies. He is not aware of any safe-sex methods. He says he had met numerous/ uncountable CSWs, who were Tamils and not Naidus. He will not come in public to speak about HIV/AIDS because he thinks that he will lose respect.

*Additional Observation: He has studied up to 5th standard. He is contented that the doctors and nurses look after him well and only after admitting in this hospital he feels better. He has taken country medicine. No one in his family knows about his seropositivity. He has 9 siblings, 7 sisters and 2 brothers. Contradicting this, he later says his brothers and sisters know*

*about his infection. He is the eldest in the family. "Do not touch" is the motto of his village people when it comes to AIDS patients.*

*He has visited numerous CSWs, he says that some of them were very thin, and most of the CSWs are married, some of them say that their husbands are a tuberculosis patient, drunkard etc. He explains that the CSWs have a namesake husband, so that they are not troubled by other men and that man extends security to them and takes money from the CSWs.*

## 47

A 40-year-old lorry driver was interviewed on 1-4-02. He has studied up to 7th standard and hails from a remote village in Tamil Nadu. He was reluctant to give the name of his village. He is married for the past 10 years. He has visited CSWs from the age of 20 for jollyness and happiness. He has two children. The children and his wife are free from HIV/AIDS. He is ignorant of the disease as well as the safe sex methods. When he was away from home for a week or a fortnight he necessarily visited CSWs from other states.

***Additional Observation****: He is taking treatment for the past 1½ years. His wife is a close blood relative of him. Still he was all alone in the hospital. He wants the government to find medicine to cure the patients. He did not want to disclose his community. He also acknowledges that he visited CSW to be happy and to enjoy life.*

## 48

An ITI trained 31-year-old man living with HIV/AIDS was interviewed by us at Tambaram Hospital. He is a motor mechanic hailing from Tirupur. He is married. He thinks that when he was in Mumbai he should have got the disease. His marriage was arranged and was done without asking him whether he was satisfied. His wife is unaware of the disease, but his mother knows everything. All his friends are alcoholics. When he visited CSWs he was ignorant of HIV/AIDS.

***Additional Observation****: In our sample of 60 we see that several persons from Tirupur are affected by HIV/AIDS this indirectly hints at the fact that such industrial cities have higher concentrations of HIV/AIDS due to large presence of migrants.*

*He has got himself tested at Madurai. Many HIV/AIDS patients do not test for the disease at their hometown but only a little away from their hometown in order to keep the information away from the family.*

## 49

A hotel manager aged 28 years and who has studied up to 8th standard was interviewed on 28-1-02. He is presently an inpatient of the Tambaram Sanatorium. He started to have sex from the age of 23. He also drives auto. He is not married. He had sex with commercial sex workers belonging to Andhra Pradesh, Kerala and Pondicherry. He did not know how he contracted HIV/AIDS. He says he had fever often and also had blisters in the private parts. He used to take treatment for the fever and STD and repeatedly he had the same type of attack. He has spent from Rs.50 to 100 each on the commercial sex workers for each visit.

***Additional Observation***: *He had used safe sex methods at times and not always. He also admits that he has STD. He adds that doctors in Tambaram Sanatorium take care of him properly. He has also taken treatment from Kerala, but it was not effective that is why some of them advised him to get treatment from the Tambaram hospital. He acknowledges he has committed the mistake willfully. None of his friends know about this. He accepts that the CSWs always advised him for safe sex methods.*

## 50

A 28-year-old youth hailing from Tenkasi was interviewed on 6-3-02. He worked as a construction labourer. He is married and has 2 children. He has tested the children for HIV/AIDS and they do not have it. He has tested his wife and she has HIV/AIDS. His wife says she has got the disease from him. Only in the Kuttralam platform he had sex with CSWs by paying Rs.50 to Rs.100 for each visit from the age of 20. He had blisters and scabies in the private parts. He has beaten his wife whenever she questioned him on where he stayed during the nights. His mother-in-law is unaware of the disease.

***Additional Observation***: *He is an inpatient for the past 2 years. His wife interferes to say that her husband lies and also her in-laws cannot understand the situation, they only will blame her* *saying that their son will never make mistakes and that he has got HIV/AIDS from her. She alleges that her in-laws are ready to get their son married, even for the second time, in spite of the fact that he is an HIV/AIDS patient. Both of them have also taken country medicine for cure. This man's parents are so superstitious they think that by going from temple to temple and seeking the advice of fortune tellers, would cure him of the disease. She had felt that she can die but she lives for the sake of her two daughters. They say good-willed and service-minded people should adopt the children orphaned due to HIV/AIDS.*

## 51

A 44-year-old, driver by profession was interviewed on 15-02-02. He is from Neyveli. He is married and has a son and a daughter. He says that from '80-'84 he had regularly visited commercial sex workers. He confesses that it is the only mistake he has made. He has visited CSWs from Maharashtra, Rajasthan etc. When he disclosed about his activities his wife has shouted badly at him. He had taken treatment from Kerala and spent Rs.10,000 on the same.

*Additional Observation: He was an outpatient of the hospital for one year. Blood test was first performed at his village and he was asked to immediately to go to the Tambaram Sanatorium. He admits that he used to spend Rs. 20 to Rs. 30 on CSWs.*

*He explains when the men went on trips, the CSWs would be under the trees to offer their services. He also adds that the CSWs are of three types: those who take up the trade due to poverty; those who do it for earning; for easy money, and because they like to do it and those who take it for they have failed in love. He tells us that in those years 80-84, he did not know about this disease. Now, he is so depressed that he wishes to commit suicide.*

*His wife does not come and see him. He says he is like an orphan. None of his relatives come to see him. What people feel about the disease he says is "good" (i.e., the social stigma it holds with it). He says information on the HIV/AIDS disease and the way it spreads, can be given as a course for study in schools.*

## 52

On 13-02-2002 we interviewed a 35-year-old driver from the Kadamadai town in the Dharmapuri district. He is married, his wife was tested and she was found to be free from HIV/AIDS. He acknowledges that because he is a driver he has widely visited CSWs. He has all habits like drinking and smoking. He has spent from Rs.50 to Rs.100 for each visit to the CSWs. But he refused to name the states in which he had acquired the disease.

*Additional Observation: He is a driver in a private concern. His family occupation as well as his wife's family occupation is only agriculture. Because of very poor yield and no proper payment they have sought different jobs. He says he has asked the CSW whether they have AIDS, and they have said "we do not have AIDS".*

## 53

On 19-2-02 we interviewed a man aged 32 years who has studied up to 9th standard. He is married and has three children. He resides in Thiruvarur. He is also affected by tuberculosis.

He is a lorry driver by profession. He tells us that he has visited commercial sex workers once a week. He took to these habits because he was advised by his doctor that after the birth of their 3rd child, his wife was very weak and sick, so he should not have frequent sex with her.

*Additional Observation: He has stopped sending his children to school: his 10 year-old girl and the 9 year-old boy. He has not tested his wife for HIV/AIDS. He shares his view that when the doctor advised him to go to the 28th ward he says the disease increased. He fainted twice and fell from the cot and once again fainted near the bathroom.*

*According to him, this disease worsens with fear and tension. He thinks that because his wife is weak he would have got the disease from her. Yet, he has not tested her for HIV/AIDS and she is, according to him free from the disease. He has an enormous complexity of thoughts! This man feels that for HIV/AIDS treatment only the Tambaram Sanatorium gives the best treatment and no other hospital in South India gives such a good treatment. He spoke about the TNSACS in Egmore. He feels that to observe HIV/AIDS patients who are terminally ill is like living in a torture chamber.*

*He had stomach pain, then they said that he had jaundice, then he began to lose weight, stamina and strength. Then they said he had tuberculosis. Then they said it was not TB but only*

*HIV/AIDS. But still he has doubts whether the doctors have checked him properly. Doctors have told him, "If you have made this mistake, even 20 years ago, you have the chance of getting HIV/AIDS." He wants to be cured and then continue to work as driver.*

## 54

On 13-02-02 we interviewed a man living with HIV/AIDS at the Tambaram Sanatorium. He was reserved and did not disclose his name or age or his native place. He was working as a labourer in a small hotel.

He is unmarried and does not know any trade. Might be it looks like he should have been a runaway from home. He did not mention about his parents or his native place. He advises the youth to be self controlled. Later he said after much conversation, that from his childhood he has lived all alone in the platform. There was no one to question or advice him. He has had sex with hundreds of CSWs and has spent Rs.100/- per visit. He has gone to Bangalore, Vijyawada and Mumbai on work.

*Additional Observation: He did not disclose name, age, caste, religion or his native place. He also said that doctors look after him well. He is very affected for he says he had no one to reform or bring him up properly. As he was a platform dweller he was misused by several people. He was very much abused physically also. He looked very ill with unattended scabies and scars on the body with flies buzzing around him.*

## 55

This patient did not answer several of the questions. He mentioned that he has children from which we have to conclude that he is married. Also he never gave the name of any of his close relatives be it his father or mother. He has paid from Rs.10 to Rs.50 for every visit to the CSWs. He has gone out of his native place and to other states where he had visited commercial sex workers. For the past five years he is in the hospital.

First when he was 24 years he spent about Rs.150/- on a commercial sex worker. He says the number of women he has visited is numerous. Two of his close friends have died. He does not believe in God. In his opinion, there is no difference between the educated and the uneducated in getting HIV/AIDS or in visiting commercial sex workers.

*Additional Observation: He did not answer the question whether his children are studying in schools or not. He says that government should punish and control men seeking CSW and patronizing illegal sex. He is of the opinion that people think even if they talk with the HIV/AIDS patients, they will get HIV/AIDS. In his opinion more than their fear of HIV/AIDS they fear the social stigma attached with it.*

## 56

On 2-3-02 we met a 27-year-old man in the Tambaram Sanatorium. He has studied till the 7th standard and his hereditary occupation is agriculture. He hails from Tirupathur. He is a Hindu and is married. A child was born and died at two years. He has a second child who is two years old. His wife is pregnant with the third child. His marriage was an arranged one. He is a driver. He has not tested his wife. He acknowledges openly that even after marriage whenever he went to different places he visited commercial sex workers. Only after marriage he has visited several commercial sex workers and had paid them about Rs.50 per visit. He confesses that he does not know that he would get HIV/AIDS by this. He adds that no one in his family knows about his disease. He further states that he had visited an Andhra Pradesh woman who was a commercial sex worker and she was good looking so whenever he went on business trips to Andhra Pradesh he used to visit only her and now the last time when he went to see her they said she had died. He thinks that this was the mistake he made and he feels that he got the disease only from her.

*Additional Observation: Even now he has not tested his wife who is pregnant though the child might also be infected. He does not show any interest in this. For headache first he used tablets from medical shops. He used to get fever often and only therefore he visited an MBBS doctor, who gave him medicine and he felt a little better. Then once again he became sick so he visited an MD doctor. That doctor said 'only after testing your blood I will start the treatment.' This doctor said that he has HIV/AIDS, so advised him to go to the Tambaram Sanatorium. The doctor had said "you may live 2 to 3 years, only if you go to Tambaram, otherwise it is*

*very difficult." So he came to Tambaram. He acknowledges that several of the CSWs have asked him to use condoms.*

*Only he has refused to use it! He is of the opinion that because of ignorance, he has contracted the disease. His wife has shouted at him that he has contracted some sexual diseases from others. He had hidden this from her. He explains that if he admits that he suffers from HIV/AIDS nobody in the village would even give him water, they will treat him badly. That is why he has said that he is suffering from tuberculosis. He next points out that in this large world with so many people, is there not a single scientist to find medicine to cure the disease which kills thousands of people! He advises those who get married to have their blood tested. He has no hesitancy to say what he feels, for we have come to work for their sake only.*

<div style="text-align:center">**57**</div>

A 37-year-old lorry driver was interviewed on 10-4-02. He is from Gujarat. He is married for the past 4 years. He has visited commercial sex workers and paid them from Rs.50 to Rs.100. He has gone to Mumbai, Kerala and Andhra Pradesh.

***Additional Observation****: He started to visit CSWs at the age of 24. He is affected by STD/VD. He has taken 9 injections for cure. He had contracted the disease on his trips from the roadside CSWs. After this attack of disease he never used to get down from the vehicle which he was driving. He has no children. The drivers say that when the lorry travels on the road, the CSWs flash torchlights in the dark, so that the driver gets down and has sex. He says they pay from Rs.50 to Rs. 100. He accepts that he had so far visited 20 to 25 CSWs. He adds that he is separated from his wife for the past three years.*

*He lost weight that is why his blood was tested and he was found to have tuberculosis and HIV. He says the HIV/AIDS patient become very tired even if they walk 20 feet. He feels, "HIV/AIDS is the punishment for their sins." He comments that these days even at 12 years, a boy seeks sex. He adds very sadly that he is an orphan.*

<div style="text-align:center">**58**</div>

A 20-year-old with a pass in SSLC was interviewed on 10-4-2002. He is from Mumbai. He has spent Rs. 100 per visit on the commercial sex workers. His native place is Villupuram. His profession is stitching bags. He has got this disease only from Mumbai.

*Additional Observation: He is not very social for he has not talked with other people around him at the hospital. He could not talk with us as his father was curiously and sincerely observing him.*

## 59

On 8-3-02 we met a 32-year-old man living with HIV/AIDS. He belongs to the Hindu Vanniyar community and is from Ulundoorpettai. He is a construction labourer. He had gone to Mumbai when he was 18 years old on a business trip with friends. He has spent lots of money on the CSWs. His wife is also infected with HIV.

*Additional Observation: He feels that he did mistakes because of his age. He also tells that he has used safe-sex methods. He wonders how he has got the disease in spite of adopting safe-sex methods. But later he says sometimes he used condoms and at other times he did not use it. According to him, only the CSWs supply men with condoms. He has never purchased it. He has spent lots of money on CSWs. Even after marriage, he had frequented CSWs. Only his wife was benevolent enough to admit him in this hospital. His parents do not know that he has the disease.*

*The child and wife are also affected with the disease. He says that after seeing the HIV/AIDS patients in the hospital he has decided never to make mistakes. He has talked to CSWs, some of them had said that they were forcefully brought into this trade, and they have even asked him to kindly free them. But he had replied that since the police was also supporting the trade he could not do anything.*

## 60

We interviewed a 32-year-old SSLC educated man living with HIV/AIDS. He did not want to disclose his name. He is a driver by profession. He has spent around Rs.100/- on CSWs per visit. He knows Tamil, Kannada, Telugu and English. He is married and has a child.

*Additional Observation: For eight months he had not disclosed the fact that he was HIV seropositive to anyone. He suffered silently. He has spent over Rs.5,000 on country medicines.*

*He has hidden the news of his disease from his family because no one in his family must suffer the social stigma due to him.*

We had just given the verbatim translations and additional observations of the interviews of the 60 HIV/AIDS patients who had acquired the disease by migration.

These interviews are as per their statements and the following drawbacks may be found. When we talked to the patients in the first place they were frightened, they were not open, they are extremely reluctant to talk. Some of them did not like to answer all the questions. Very often they told us to skip the chronology of the questionnaire. They preferred to talk whatever was on their minds.

Secondly, several times we had to take the interview of the slightly healthy person with AIDS by bringing him/her out of his bed, for some of the patients become aggressive even if a neighbouring patient talks about HIV/AIDS and the cause of it. They were afraid that we will publicize their name and so on. So it took several days to become friendly. Next we had to give them complete confidence that none of the news will be given with identity or any publicity. We also explained to them that it was purely for research purpose.

Thirdly they were so much depressed that it was very difficult to make them talk. The guilt was predominant in them, however we wish to state that most of them, after a few minutes of conversation, accepted and gave the true cause of how they got HIV/AIDS. Fourthly, they were so poor that even if they wished to have a cup of tea they could not get it for they did not even possess that paltry amount so patients who could walk were taken to tea shops and we spent several minutes with them over a cup of tea. It was a very bitter experience when some of them wept to us stating how lavishly they have spent money and how poor their plight was now.

Some of them were reluctant to talk so not much of information could be gathered. Many of them did not give their name and some did not give their parents or in laws name or the name of their village or even the district. Now we give a brief description of our feeling, our ways of approach, our discussions etc.

At the outset we wish to state that our interviewers too were a little hesitant. Some of my students got fever the first day when they met the patients. One of the students just ran away never to return to work. Later, they confessed that it was because they had been very frightened. Though we have given them enough courage and guidance and instructions, in the beginning stages, fear was dominant in the students also. Many weeks into the project, we heard our students tell us that they had become the best of friends with the AIDS patients.

Only after interviewing the people with AIDS, we too came to know several of the symptoms, external changes and internal changes of the disease.

We had a ready-made questionnaire which had over 191 questions. It is important to note that several questions were evolved by discussions thereby making the answers and questions very spontaneous, not only from the questionnaire but also from our discussions, which gave us great insight into the understanding of AIDS. We had gradations in the answers for some questions. Gradations cannot be given to all questions like "Are you married"? "How many children?" Profession? Educational qualifications? But we could grade the questions like what made them visit CSWs? How many? Reasons for it? How they felt when they come to know about the disease? etc. The answers were graded or it was marked in some cases like indeterminates.

From their facial reactions, change of tones we analyze and grade their answers.

The mode of communication was very difficult for even their Tamil dialect was very different from our Tamil dialect mainly because they hailed from remote areas. Some of them could not follow our questions. At other points of time, we could observe that at times the patients became self-conscious, and they couldn't speak out their mind.

Also almost all of the male HIV/AIDS patients had all habits like alcoholism, smoking and visiting CSWs and most of them have become victims of this disease because of ignorance. They denied any knowledge about HIV/AIDS. Some of them had taken Siddha treatment from Kerala.

Majority of the HIV/AIDS patients selected for our study on migrant labourers are drivers. A striking point is that all of them

(except one) wanted the government to find medicine for curing AIDS.

Only one patient was of the opinion that discovery of any medicine to cure HIV/AIDS would result in the further decline of morality. So he did not want any cure to be found for this disease.

The faith in God of the people living with HIV/AIDS is distinctly different. One category of people have changed faith from their native religion to Christianity. Some have lost faith in God: they have become atheists, some of them say prayer is merely a ritual, others say that they have more faith in modern medicine than in God. Some of the patients are highly philosophical when they speak of their death and of God, some others have resigned it to fate.

Some more awareness programs must be given to those in the rural areas. A few of the persons interviewed by us said that as they were uneducated they could not read the advertisement on the bus or auto. So some of them suggested "figurative methods" of advertisement. Some of them could not even sit or walk or stand. Everyday after the interviews we met and recorded our opinions in a diary and this is also used in our research.

A few of HIV/AIDS infected patients preferred mercy killing at the terminal and most advanced stages of the disease. They are unaware of the present days steps taken by government and also of the many awareness drives.

All of them wanted to ban commercial sex workers. They felt that the government, the public and the police should cooperate for it. One man said that the CSWs must be rehabilitated to do some other work, another said that they should become coolies. Another man said that poverty must be abolished for AIDS to be annihilated. A few of the affected women with AIDS said that the CSWs exist because of the demand; if there was no demand for commercial sex workers, what would they do?

So for this the only solution is that men should have control over emotions and at no time seek commercial sex workers; then automatically the CSWs would have to seek different ways for their livelihood/ existence. On the contrary the demand is so high that is why they thrive well. Only one woman with HIV/AIDS, said that she never knew that there was something like commercial sex workers.

Another observation we made was the symptoms the HIV/AIDS patients suffer. Some of the patients were of the opinion that the government should find some medicines to control their emotions before they find medicine for HIV/AIDS.

Most of HIV/AIDS infected patients in our sample for this study were poor, lower middle class or middle class. To purchase this disease itself they had to spent from Rs.20 to Rs.100 for each visit to the CSWs!

Only after acquiring the disease they had become poorer for they had spent at the least Rs.20,000 for medical check up, transportation, medicine, etc.

We observed that the HIV+ men in our sample, had, at no point of time in their lives exercised control over their sexual behaviours. This is made clear from the way they were answering the question "How many CSWs you have visited" - innumerable. Uniformly they were all suffering from depression. Yet another striking behavior of the HIV/AIDS affected women who were with their husbands in the hospital was that they were composed and calmer than the women were alone. This alone shows how much the women who have been affected by their husbands react. Some other women were so depressed they wept all the time.

The following were the external observations made by us on the majority of the patients:

(1) More than 51 of them had skin ailments like scabies on their hands and legs, change of color of the skin in the cheeks and forehead, boils on the eyelids, big boils on the body, some bandaged, some with medicine applied, some not properly taken care of with flies hovering on.

(2) About 52 of them had Tuberculosis, cough, common cold and breathlessness.

(3) Around 31 of the interviewed patients were very thin and not even in a position a walk 20 ft. without difficultly.

(4) In their talk and act the helplessness of their state and their depression was conveyed.

(5) Several were all alone with nobody to look after them: for in some cases the relatives did not know, as they have not informed and for some their relatives have deserted them because of the disease.

(6) The youngest women patient was a 16 year old married at 11 and widowed at 11½.

(7) The youngest male patient in our study was a 20 year old.

The emotion of these 101 patients will be studied based on the data collected.

Chapter Seven

# CONCLUSIONS

The conclusions of our study are mainly made from the analysis of the data using fuzzy theory and also from the feedback we received from the 60 HIV/AIDS patients whom we have interviewed as a part of our research study on HIV/AIDS. In fact for several of the interviews, we held free-wheeling discussions, allowing the patients to speak at length to us, we made them pour out to us all their worries and trauma, we shared all their suggestions.

In fact, the format of our interview was designed so as to make it a very interesting session, which brought out the personality and nature of the patient whom we interviewed. We hope, that it is not wrong to point out at this juncture, that all those interviewed were not only HIV seropositive, but suffered from symptoms of full-blown AIDS. Next, all the patients whose data is presented here were interviewed at the Tambaram Sanatorium. As a result, our data is restricted to the type of people who visit the Tambaram Sanatorium in maximum numbers: the poor and middle classes. This is an unfortunate drawback of our study on HIV/AIDS. The disease by itself knows no economic strata. There is migrancy among the rich and the elite too, but our study does not cater to those angles.

Also, it is pertinent to mention here that the results which we have presented here are purely technical and in no place, our feelings or our opinion is given. We therefore wish to make it clear that this work does not carry any of our bias. Thus, this chapter on conclusions is based on the analysis of the data. This chapter has four sections.

7.1. Observations from our study
7.2. Suggestions given by the HIV/AIDS patients
7.3. Conclusions derived from the interviews
7.4. Conclusions and suggestions based on our

research study

The only angle of the study of HIV/AIDS is in relation with migration; specially migration from the rural areas. So we do not say the same can be applied to patients who have not acquired the disease due to migration—once again by the term 'migration' we mean that the disease was from travel/mobility be it any jolly trip or a business-related travel or travel as profession.

## 7.1. Observations from our study

The following observations were made from our study. We first enlist our observations and then we make the derivations from these observations by consulting with the experts in those fields. We have consulted people of religious groups, agricultural researchers, trade union leaders, doctors, sociologists, psychologists and NGOs. The raw observations are:

1. Over 74% of the HIV/AIDS patients interviewed were agricultural coolies (i.e. labourers working in fields of others, who were hired for daily wage).

   They attributed their migration to the following reasons:

   a. Around 8 of them said that due to caste clashes and other inexpressible reasons they were forced to come away from their villages and migrate to the cities where they had to seek their livelihood. They said that this was because after the caste clashes, they were not given job by their oppressor caste landlords, and so, to eke out their living as well as for fear of further assault on themselves they have migrated to the city and changed several jobs varying from construction labourers to cleaners and workers in hotels, petty shops or tea shops.

   b. Nearly 70% of the migrant labourers with HIV/AIDS infection attributed their migrancy to the failure of agriculture. When we asked about it, they attributed this to the failure of monsoon, also it was the failure of agriculture which meant yield was much less than the expected. So, the landlords were not in a position to pay them well, as their pay was in proportion to yield, that too in terms of harvested grains.

c. Some of them were purely employed as 'Saanars' where their work was to clean and manure the coconut trees and pluck coconuts and they complained of no work for in several places the coconut trees had dried or did not give any yield.

d. The other migrant labourers, around 30% were truck drivers or car drivers or construction labourers who had to widely travel for their livelihood, and for many of them agriculture was the only hereditary occupation.

2. 93% of them were living below the poverty line presently (at the time of the interview), though over 80% of the persons had belonged to lower/upper middle class.

3. Over 96% of them had all the addictive habits like drinking, smoking and visiting CSWs.

4. Most of them, that is, around 92% were treated by quacks and over 68% have taken treatment from Kerala were a few of them were tested as HIV negative but really had HIV seropositivity. Thus they lost most of their money in wrong treatments. Also they spent money in visiting temples, fortune-tellers, astrologers etc. which only made the disease chronic and more advanced.

5. Most of them confessed that they feel better after admitting themselves in the Tambaram Sanatorium. Most of them wanted free medicine and treatment for AIDS.

6. Most of the men were not interested in making their wives to undergo the HIV/AIDS testing. Neither were they willing to test their children. Several said that by doing so the disease would be made public and relatives and friends would come to know of it and hence this would make it impossible for them to live in the villages.

7. Around 40% of them had changed their religious belief; that is they had changed from their native religion/Islam to Christianity. Some had given up on their own and started to believe in Jesus Christ though their family members had not yet converted. These conversions could be attributed to the Christian nuns and priests who often came and consoled these AIDS patients. 7 of them had lost faith in God for they said God had no power to cure so they felt it is a waste to pray.

8.   Almost all the men interviewed (59 out of the 60) accepted that they had visited CSWs; only one man claimed that he had contracted the disease from his wife.

9.   In most of the cases, the CSWs offered safe sex methods to these men, only they did not use it because of their macho behaviour: they were intoxicated, they did not want to use it etc. Shockingly, over one-third of our sample study did not know the modes of communication of HIV/AIDS. This was another probable reason why they had not used condoms.

10.  Over 85% were not interested in anything and were depressed and frustrated not about their sufferings or the symptoms of the disease but the social stigma and their increased isolation in a society which had alienated them. They had no one to look after them, no one to whom they could turn to for solace or comfort.

11.   Their high-risk behavior is attributed to the following:

> a.   Long stay away from family.
> b.   Lewd/vulgar and porn movies inducing sexual aberration.
> c.   No social fear as they lived away from home and relatives.
> d.   They earned well, so to have pleasure, sought after commercial sex workers.
> e.   Drinking habits and other addictions.
> f.   No motivation in life, unambitious, leading risky, unplanned lives.
> g.   Illiteracy.
> h.   Absence of any goal in life.
> i.   Negative peer group influence, bad company: Friends who lead them astray and consequent behavioural change.
> j.   Absence of moral, ethical/cultural values of life.
> k.   Availability of CSWs for paltry rates.
> l.   After work hours the only mode of relaxation was sex/CSWs.
> m.   The vulgar and obscene cinema songs which kindle only the lower and base emotions.

n.  Complete freedom with no one to question and had no conscience and could not think of the consequences.

o.  The sudden change of atmosphere and the urban settings provided all support for them to fall a prey to commercial sex workers.

p.  They earn daily and spend daily and are not even interested in wife or children or their own dress or looks. Such a 'vagabond' lifestyle resulted in their condemnable sexual behaviour.

12.  The migrant labourers had no union or a proper peer-group leader who could advise or channelize their behaviour. They are apart from the mainstream, they are not organized workers, and so they do not enjoy either the benefits or the legal safeguards that other workers enjoy or do not know the modes of saving.

13.  Around 60% of them felt relieved after talking with us and explicitly said this to us. They felt consoled and some of them cried, displayed emotions of anger, despair and dejection. Some 20 were reluctant and they never opened their hearts fully and did not give the truth at all: even their village's name, father's name, caste, age etc. The remaining were frustrated and left their life to fate ("vidhi") so they were not openly coming out to express their state of mind.

14.  Some wanted the government to help them with free medicine and basic needs like hot water, good food, good bedding etc. and recreation like TV., radio, divine music, religious consolation and psychological counseling.

15.  85 of them felt that the Tambaram Sanatorium is the only hospital in south India which gives costly, effective medicine that makes them feel better.

16.  All of them except one man wanted the government to immediately find a cure. That man said that a cure for AIDS would cause a further degeneration of morality. Another man who was interviewed was ashamed that there are thousands of scientists in India and so far no one has found the medicine to cure AIDS.

17. Some 26 of them expressed their preference for mercy killing. They felt that when they are not in the state to get up or walk they should be killed. (i.e., death wish of terminally ill patients).

18. 93% of them said that the doctors at the Tambaram Sanatorium were taking good care of them.

19. None of them in our sample study claimed to have any spastic children.

20. It is a surprise that most of the infected men were in a position to locate from which place/CSW they had contracted the disease. Also, most of them had harboured the suspicion that they might be suffering from HIV/AIDS.

21. The first disclosure was never done to their spouse. In most of the cases where the disclosure was done, it was to their brother or to their mother.

22. Only one patient out of the 60 knew about the free counseling for HIV/AIDS, and this man too did not know the phone number. It therefore points to us that we must build up the awareness campaign. [Note: This research was carried out in the pre-Pulli Raja days.]

23. Some of them had married for a second time very recently after the death of their first wife.

24. According to our study the maximum age group in which migrant labourers are HIV/AIDS patients are from the age group 31-34 followed by 24-30. This is the stage in which the persons are highly affected with HIV/AIDS. They are in the advanced stages of the disease. From this it is observed that it is in keeping with the fact that these persons must have contracted the disease around the age group just from 20 and above. The same is confirmed from three different analysis of RTD and CETD matrices carried out by varying the number of age groups to four and five in section 2.1 in chapter 2.

25. We also make another observation that in the age group 26-30 the migrant workers who are HIV/AIDS affected patients not only had CSWs for sex but also other married women as heterosexual partners. Further the study indicates the migrant workers who are affected with HIV/AIDS in the age group 31-35 are the ones who visit quacks for the treatment followed by the age group 20-25. This could be attributed to the fact that those in these older age groups are mostly married, and hence they face problems of anonymity and disclosure; apart from stigma and discrimination within their family.

26. It is important to note that the migrant HIV/AIDS affected persons in the age group 31-35 are prone to bad company and alcohol. However the data shows that they don't have other women as heterosexual partners but visit only CSWs regularly.

We have given here only the major and main effects from our study which cannot be got from the collected data.

27. The HIV/AIDS patients among migrant labourers in the age group 41-48 are now least involved in CSWs, other women, quacks, alcohol, smoke etc. However this may not be true in case of educated city dwellers or rich people living with HIV/AIDS.

This is because after that age they are not migrants. They would have settled into a life of secure and steady incomes. After this age, they would not indulge either in cross-border migration or rural-urban migration. That is why the mathematical study shows a sharp negative trend at such an age.

28. Taking the type of profession described in section 2.2 of chapter 2 we see that based on the profession, is with profession in the 'on state' and all other states in the off state then we analyze using FCM model the resultant gives all the states except no social responsibility and more leisure is in the 'off' state.

Thus depending on profession this dynamical system predicts the on state of all attributes associated with the HIV/AIDS patients like easy money, no education, bad company, bad habits, socially free and economic status. Also if easy money, no social responsibility and more leisure are in the 'on state' then all the 12 attributes given in section 2.2 come to the "on state".

Thus when migrant labourer earns money easily with no sense of social responsibility and has lot of leisure then mostly the migrant patient is endowed with all bad habits.

We are able to see that the exaggerated notions of masculinity, high levels of carelessness and risk-taking, and an attitude which doesn't takes stock of the disadvantages of such a lifestyle usually lead men into contracting the AIDS infection.

Also the study indicates most of them are unaware about the disease. All these are established mathematically using the FCM model. It is important to mention here that in certain cases i.e. when certain nodes are on all other nodes become on their by making the importance that attribute or that particular concept in the HIV/AIDS patients strong influence on all other nodes.

## 7.2 Suggestions given by the HIV/AIDS patients

These were the suggestion given by these HIV/AIDS patients whom we have interviewed.

Here we do not discuss or analyze about the problems of implementing the suggestions but only give it as stated by them; it is up to the concerned to take a note of them.

1.   They felt intensive awareness program about HIV/AIDS must be given to them. Several of them felt that this education was more important and urgently needed than adult education as most of them made mistakes due to lack of knowledge about the spread of HIV/AIDS. Only after knowing that they were affected by HIV/AIDS they came to know about safe and protected sex. Some of the patients are yet to know how HIV spreads.

2.   As advertisements are given about products and events and such other things, advertisements regarding the awareness for HIV/AIDS can be given both in Radio and TV. Also in cinema theaters advertisements on safe sex can be given. Moreover, concrete steps should be taken to ensure that the radio as a sector is given a lot of weightage. This is because, almost all the truck-drivers, lorry drivers and the migrants spend a lot of their time on the road, and elsewhere listening to radio. The new breed of FM culture, has increased this a hundredfold. So, it is in the radio sector that we have to bring in the maximum of our change. Further, unlike the television which is planned for family viewing, the radio is often listened to, only by one person at a time (in the family), so privacy while listening is guaranteed.

3.   More health-care centers must be given in villages with a section of doctors to counsel about STD and HIV/AIDS.

4.    When agriculture fails mainly land owners are helped by the government but till date no steps concerning the daily coolies (or hired farmhands) who work in these fields are helped. Also steps have to be taken to secure their futures and also ideas and suggestions for alternative solutions could be implemented.

5.    Free medical check up for HIV/AIDS can be conducted as womenfolk especially those living in these village, whenever their husbands are infected and they fail to take any treatment until the diseases is in the full blown state.

6.    When questioned about embracing Christianity the converted male patients said that most of the Hindu religious leaders work for cancer patients, disabled etc. but no Hindu religious leader has so far come publicly to give support or solace to the HIV/AIDS patients. "Also Hindu faith does not pardon you it only says you bear the karma, but Christianity pardons you so to get pardoned, we have embraced Christianity." To an extent, their change to Christianity basically results from the fact that the patients had been approached by Christian missionaries, and these religious people had given them a patient hearing. Also the practice of confession, made the mental burden of the AIDS patients, seem light and easy.

7.    The Government should provide job opportunities for the HIV/AIDS patients by giving them a certain percentage of reservation in employment; this percentage can be fixed with the percentage of AIDS patients in the general population in India. This shall open up a great amount of change, remove the stigma that is associated with the disease, and also remove existing patterns of workplace discrimination against HIV/AIDS.

8.    In the villages, the HIV seropositive patients are ill-treated so badly that even basic needs like water, provisions etc. are denied not only to them, but even to the family members of HIV/AIDS patients. That is why they fear to disclose the disease even to their wives or children as a result of which the infected men decide not to carry out the HIV/AIDS test on their wives or children, because they fear that the information shall spread everywhere very soon. Such offenders who disgrace and insult the AIDS patients must be made to undergo stringent punishments.

9. Counseling must be given to the general public and society at large to treat AIDS patients with respect. Such sensitization programs are the need of the hour. We feel the grave importance of this situation because, even in their death beds the patients think of nothing else, but lament at the stigma which they are made to undergo.

10. Some of the HIV/AIDS patients who had committed mistakes in the age group of 20 to 25 years said that cinema was one of the major reasons for their fault and that it had certain scenes that had made them lose self control and seek for commercial sex workers. They were against screening of such movies. In fact one of the patient requested the government to take steps to ensure that some medicine should be found to control their sexual desire and the urge for release and gratification, instead of finding medicine for AIDS.

11. All of them shared the view that free HIV/AIDS test and medicine by the government would be of use and would enable the concerned to be tested without any difficulty or even monetary bearing.

12. Almost all of them felt that sex education in schools was a necessity for that alone would make them aware of this dangerous, incurable, and socially a so-called stigmatized disease.

13. The uneducated HIV/AIDS patients lamented that they could not read and several of them had no TV so by some other means the uneducated should be given awareness and programs to educate them about HIV/AIDS and the way it is spread.

Among the 60 persons we have taken for our sample study, over 70% were educated only up to the primary school i.e. 2nd, 3rd, 4th, or 5th std. Categorically all of them emphasized that HIV/AIDS can be completely eradicated if the commercial sex workers did not exist. All of them agreed upon the fact that they were infected mainly by them. This lack of awareness is very frustrating.

Some of the infected patients (we have reliably learnt) still continue to visit CSWs in nearby areas and they practice unprotected sex, this is highly frustrated for they would be infected CSWs/ or other women by the number of infected increases steadily.

14. It is of extreme importance that the AIDS patient who suffer from depression given medical treatment to alleviate the depression. This depression, which some of the doctors say is a side effect of the anti-retroviral therapy, is a cause for the increasing number of suicides at the Tambaram Sanatorium.

An important but a very delicate observation is that most of the persons contracted the disease from Andhra Pradesh, New Delhi, Mumbai, Goa, Kerala. Majority of the infections were contracted at Mumbai, next at Andhra Pradesh where the CSWs are available at the cheapest rates.

We could peculiarly observe that even those who had never visited CSWs in their own hometowns visited CSWs on their visits to other cities. In other words, the migrancy/mobility further aggravated their high-risk behaviour and increased their vulnerability to AIDS.

## 7.3 Conclusion derived from the Interviews

At the first instance we clarify that the interview was carried out on a small sample from the Tambaram Sanatorium, where the inpatients were only from villages, uneducated and have acquired the disease by migration only. So by no means this study will reflect on the general behaviour or the socio-psychology of HIV/AIDS patients. Also it cannot be related with educated city dwellers or the urban rich.

For, in our opinion to pick up a case study of educated and posh city dwellers for HIV/AIDS that too from rich, upper and middle class society may be impossible for the disease is well concealed, the treatment is taken secretly and they are unapproachable for they fear that the news might be spread and they are righteously bothered about their anonymity. In the cities, the stigma associated with HIV/AIDS carries a heavy price.

So our study pertains only to the ignorant people from rural and semi-urban areas who lack the awareness and thereby infect themselves, their wife and children.

The conclusions given here are from the interview only (exclusive to those derived mathematically).

1. At least 70% of the CSW had advised them to use safe sex methods, only these men did not use condoms. On one hand, it shows that the NGOs have been successful in educating the CSWs. On the other hand, this practice is very dangerous, because these men risk not only themselves, but also put the health of the CSWs at stake, by having unprotected sex.

2.    Several of the patients get admitted for treatment for HIV/AIDS and after they are little better they once again begin to drink liquor, smoke and go to the CSWs and at the same time also have sex with their wives.

Some of the patients also throw away the costly medicines given to them. It is important to take this into account, and to build awareness and counselling programmes at such hospitals and care centers in order to improve the patient's health.

3.    At least 10% of them said they used safe sex methods yet they had contracted the disease.

This make us to think in the following ways:

a.    Their state is not normal and alert, for almost all of them are in an intoxicated and drunken state when they visit the CSWs. So, they admit of having forgotten to wear the condoms, etc.

b.    Some of the condoms might have been of low quality.

c.    Certainly the NGOs outreach and other awareness programs have educated the CSWs about the importance of safe sex; unfortunately, the CSWs cannot do anything when the customers are so reluctant to use condoms.

Finally from the interviews, observations and analysis of these patients the following has to be mentioned:

Men take pride in saying that they have addictions to alcoholism, smoking and visiting CSWs. It has become part and parcel of their male ego to proudly display their aggressive macho behaviour and to say that they have all the addictive habits instead of feeling bad for being a slave to it. Of the 60 men interviewed only one said that he did not have any bad habits. Unless men shed their false ego and lead a life of reality most of the awareness program about HIV/AIDS will only be a failure.

This study and analysis pertains only to illiterate villagers who are migrants; so when we are speaking about the increased percentage of unprotected sex, or about the aggressive and notable masculine/macho behaviour we are only speaking about this section of people.

We have found out that there are high levels of ignorance about the disease. While one section doesn't know how AIDS spreads; there is another vast majority that believes in various baseless myths about controlling the spread of HIV/AIDS; the most common examples are: eating high amounts of non-vegetarian food and having a strong body would prevent spread of AIDS, injecting penicillin and having unprotected sex would stop infection from one individual to another, washing the private parts with soda after having sex would kill the HIV virus etc. It is time awareness programs are made to ensure that such myths are shattered. Yet another factor, which is relevant and needs our attention is that most of the patients have been under the influence of alcohol when they had visited the CSWs. When there are five lakh HIV/AIDS patients in Tamil Nadu, and given the fact that there is an established link between alcoholism and high-risk behaviour, we, as a state are not taking the right steps. The future looks bleak: the government itself sells liquor and condoms; indirectly encouraging licentious behaviour.

There is another line of thought which criticizes the safe-sex promotion campaign. It is actually very unreasonable for the government to go about and say 'practice all vices' but please be careful about it. It is unethical and shows a complacent tendency. Several of these patients said that the CSW are supported by some elements among the police and by the public demand. Unless the government takes stringent measures to prevent trafficking in humans, and also ensures that the existing CSWs are rehabilitated, the way to an absence of spread of AIDS would not be possible. Current trends of the government plans give a stress to 'abstinence' and 'prevention', but what is needed is a more total revamping of our views.

It is a pity that we as a society fail to inculcate values, or decorum or social behavior or propriety. Our attitudes and the existing hypocrisy surrounding the discussion related to sex puts all our future at stake.

We had interviewed 60 HIV/AIDS patients. The first-hand information even before the mathematical analysis of the data is:

a.   The couples were more composed, relaxed and less depressed than the individual patients, except only one couple where the wife slippered the husband in our presence so hysterically because she was very dejected.

b.   Also it was seen from the data that 41 cases the affected persons were drivers (25) or their wives who had contracted the disease had their husband to be drivers. Of the 60, 45 had their hereditary profession to be agriculture or agriculture labourers.

c.   It is beyond comprehension to state that 59 of the 60 men were infected by the CSWs, that too not of their own locality or even district. Over 35% of them had CSWs mainly from Mumbai. A few of them openly acknowledge that Mumbai gives the best CSWs and so even in their long trips they halt the journey at Mumbai for the purpose of visiting the CSWs.

Nearly 22% had also met CSWs in Andhra Pradesh and the main reason attributed for it is that they get CSWs there for a very cheap rate varying from Rs.20 to Rs.50. So in their opinion Andhra Pradesh provides cheaper CSWs. One person was infected from Saudi Arabia, one from Malaysia and one from Singapore. Within Tamil Nadu, there were migrants who had been infected at Pondicherry, Kallakuruchi, Kutrallam, Padappai, Navallur, Madurai, Tiruthanni, Namakkal and Chennai.

d.   It is pertinent to note the commendable services of the NGOs and over 70% of the CSWs had advised these men to use condoms, however it is only the men who have failed to use it and thereby became victims of HIV/AIDS. Thus almost all the sixty men interviewed were fully aware of the condom only they did not use it for the best reasons known to them. We find that the men do not mind spending any amount of money to visit the CSW (and there by buy the disease) but they have an inbuilt and vehement attitude against condom usage. Even the condoms are offered the men do not use it, instead they offer a variety of excuses for not using it.

At this juncture, we also feel the need to stress for the importance of the female condom. Though it is more costlier, and more difficult in its practicality, it has become very essential for its sale and availability among the CSWs to be increased. It would in the long run, certainly reduce the alarming rates of the spread of the disease.

e.   Apart from these, 59 of the 60 men proudly acknowledged that they have all bad habits like alcohol, CSWs, smoking, etc. It is a pity to state out of the 60 men 54 of them said without any hesitation that they visited "innumerable", "countless" CSWs, which displays their atrocious macho behaviour.

## 7.4. Conclusions And Suggestions Based On Our Research Study

1.   It is of foremost importance that the future generations are sensitized about the AIDS disease. The most recent statistics about the AIDS epidemic shows that half of all the new HIV infections are in the age group of 15 to 24 years.

Moreover, over 50% of the adolescents in this age group have serious misconceptions about HIV/AIDS. Elaborate statistics also indicate that over 80 percentage of the teenagers who have sex do not use condoms. News reports have pointed out that children came in school-uniforms to take the AIDS test. Thus it has become a need of the hour for the Government to take AIDS awareness programs to the classrooms.

Steps in this direction have been taken by the TNSACS and the Government of Tamil Nadu: the most recent step in this direction was the implementation of the Directorate of Teacher Education, Research and Training to UNICEF through TNSACS for AIDS education in 3,800 odd government and aided schools in the state at a cost of Rs 50 lakh.

What is also needed very much is the dropping of the taboo connected with sex. On one hand we criticize the western lifestyle. On the other hand, we find that most of the present day youngsters are aping it. So, they must be taught our cultural values. Yet, this must not be at the cost of masking them into degrees of ignorance.

For instance, even today, it is a practice in rural areas for teachers never to teach the required lessons of sex education, and human biology which is a part of the syllabus.

Worse, in many of the schools, the very pages which contain references to the sexual reproduction in humans are stapled, and then only the textbooks are handed over to the students. What this does is doubly detrimental to the general mindset of the students. Either they grow up unawares of this scientific knowledge and are instead full of misconceptions, or their curiosity is ignited and they begin to experiment putting themselves into risk.

In order to popularize AIDS awareness among this age group, it is necessary to not only train their teachers but also to work through the various and preexisting organizations in these schools and colleges: like the NCC (National Cadet Crops) and NSS (National Service Scheme). While such grassroot students organization have taken many an issue such as Rainwater

Harvesting, they must dedicate themselves to the work of AIDS awareness.

2. Further, when we inquired the students themselves into the topic of sex education, over 87% of them expressed their wish that the lessons must be imparted not by their school teachers, but by trained experts. Majority of them responded enthusiastically to the idea of using peer-group educators who are just a few years older. They attributed this to the lack of inhibition, the ease of discussion and above all to the spontaneity. The students whom we interviewed expressed a preference to having these lessons in non-school settings.

3. Next, we must work together to see that the stigma and discrimination relating to AIDS must be totally eradicated. This has been the theme for 2003-4. When we conducted these interviews, we were very moved by the pathetic plight of these AIDS patients. In the last stages of their life, they have no one to care or cater to them. We could see that they yearned for human comfort and consolation, and that they cherished even the few hours we spent with them. Our children should be taught the importance of being humane, and we must all give back to society what we have taken from it. A day or two must be compulsorily allotted in the school calendar, and children must be taken to the AIDS care centers where they get to talk to the patients and show their solidarity. While children are made to visit orphanages, old age homes, homes for the disabled and such other places, it has not occurred to us to take these children and make them meet the HIV/AIDS patients. This shall make them empathic with the AIDS patients and at the same time make each one of them to secure for themselves a future of commitment. In the long run this shall remove the social stigma that is associated with AIDS.

4. Only pictorial representations have the maximum impact on the minds of people: this is a sheer psychology. That is why, even in Montessori schools, in the first stage the children are taught with toys and pictures. So, if we want to get something into the mind of adults, specially those who are uneducated and impulsive, the advertisement must be in the form of pictures (moving/still) and high-impact visuals and dialogues (when it is a video) only that will have lots of positive impact. This was pointed out to the interviewer by several of the HIV/AIDS patients, who said that being illiterate they had no idea of what HIV/AIDS posters or written advertisements conveyed. In this context we need to adopt a multi-pronged approach to implement the awareness. Methods of advocacy could include individual and group discussions, campaigning, exhibitions, video shows, rally, stickers, booklets, handbills, posters, flipcharts, peer and volunteer training, street plays, songs, etc.

5.   It is very unfortunate to note down, that while the media bears the onus for taking this awareness programs to great levels, at the same time, the same media must also take the blame for the vulgarity and obscene perversions that are represented. The standard dailies, the magazines, the cinema songs and movies, serials all of them present a very skewed version of life.

It is also sad to note that in order to boost their sales, a few magazines take to publishing soft porn photographs, they start columns on sex advice etc. While it is true that such columns do dispel many myths and spread awareness, to a large extent they also provide lewd reading matter apart from increasing the sales and subscriptions of such magazines.

We feel a separate study can be carried out by media-trend experts to understand the effect of such suggestive, soft porn, lewd and vulgar print media and the obscene portrayal of women in semi-naked states in the television.

While habits like alcoholism, smoking or even visiting commercial sex workers has the capability to only pollute one person at a time, the media has the capacity to lead astray thousands and lakhs of people, especially adults, at the same time.

6.   Because the cinema plays such an important role in shaping the popular culture, and because actresses and actors are seen as icons and role-models by the youth of today, we have to take steps to include them as ambassadors of AIDS awareness.

On the Indian level, we have popular celebrities like Shabana Azmi and Nafisa Ali who lend their voice for the cause of AIDS, on the international level the UNAIDS has roped in many well-known and highly renowned icons to spread the message, even a former US President like Bill Clinton, or the richest man in the world Bill Gates are involved in AIDS related charities. Likewise, at the Tamil Nadu level we can rope in celebrities to spread the message and spread the awareness. It can have a tremendous and immediate impact among the youth.

7. Sometimes the films end up giving a very wrong view of life. Not only are some of the song and dance sequences the worst in vulgarity, but some of the messages that the movies convey are downright condemnable.

For instance, a common topic that is often handled is the "rags-to-riches" story, of how runaway boys from villages, or those migrating from the villages to the cities, come under a stroke of luck and become rich and powerful. It is true that there are one-in-ten-million such stories, but the glamour with which it is portrayed has caused many a person to run away from home, and many a youth to come to the city with nothing but hope. Likewise, many people have ingrained the notion that it is okay to visit CSWs, if they do not have complete gratification within their marriage: because it is what the cinemas and serials portray.

8. Since this study of ours dealt with migration and mobility, we have found out that there are a considerable number of patients who have been infected by HIV during their visit to tourist spots. Thus, greater awareness and outreach programs much be carried out in such places.

One of our observation and suggestion is that steps must be taken to have AIDS awareness and safe-sex promotion advertisements to be placed in major railway stations; and also advertised on the giant electronic display boards along with the other routine advertisements. Over the course of a day, there are lakhs of people who traverse through the trains of India, and so, such an outreach program would be very beneficial and would target a great human mass at a time.

Our surveyors have however found that even at the Chennai Central Station, such AIDS awareness ads are absent.

9. We feel that psychologically the AIDS patients at the care centers need enormous amounts of mental support and counseling. They are very depressed, driven to suicide, and hopelessly weepy and unhappy most of the time. Many of them even expressed their wishes for mercy killing. Few of them complained that they were left like orphans. So, such patients must be rehabilitated. We also suggest that AIDS patients staying in the hospital must be formed and linked into support groups.

10. We come to know that several of the AIDS patients (over 60%) had attempted suicide after hearing about their infection. Also, most of the men had not divulged to their wives the fact that they were HIV positive. We suggest that it would be a good social investment if, on the lines of European and American countries, we design counselling programs to be given to these patients: On what they could do after being confirmed HIV+, how they can break this information to their families, where to seek support, how to take medical help, etc.

11. At this final juncture, we also wish to point out the importance of the Panchayats and local bodies in spreading the message of AIDS awareness at the grassroots level.

For centuries, the Panchayats have been the backbone of the Indian nation. It is these village institutions that have controlled and developed the structures of administration, it is this system that ensures the day-to-day affairs in villages. Mahatma Gandhi envisaged a *gram rajya,* in our country. We feel that such an oldest institution can be successfully used to have one of the most effective intervention programs.

Such an intervention strategy, where the Panchayats spearhead the message of AIDS awareness, will ensure not only a stop to the spread of the disease, but it would even remove the stigma associated with the disease that results in brazen incidents of discrimination and abuse in rural areas.

12. There has to be essential social changes coming. We feel that HIV/AIDS is not just a disease, it is a social problem. So, it needs to be tackled not only medically, but must be tackled socially. In order to solve this problem people from all strata of life must come together.

And we must intervene to prevent all the various and deep-rooted causes of the diseases. Within a patriarchal society as our own, we have to give great importance to women empowerment, for women continue to remain the most innocent, hapless and vulnerable victims of the disease.

Because of the extent to which patriarchy has been institutionalized, even women meekly accept their positions of subordination. They are unable to raise any of their voices against their husband's misbehaviours: be it visiting CSWs, alcoholism, or tobacco addiction.

In the long run, even they wrongly perceive masculinity and aggressiveness to be embodied in having such habits. Within the rural context, their levels of knowledge about the disease are abysmally low, and even after their infection most of them are unaware of how the disease spreads.

14. We could observe that the HIV/AIDS patients are more affected by the social stigma associated with the disease, than the physical symptoms of the disease itself. This is because the discrimination they suffer from, is not because of the stigma associated with a particular disease but because of the stigma associated with sex (which is one of the most known reasons for transmission of disease).

So, until this concept of dealing with sexuality is no longer treated with crime and guilt, we cannot expect the discriminatory attitude to be erased. It is the responsibility of all of us, as citizens, as members of a civil society to give respectability to these affected patients, and to treat them with love and care. People must not traumatize the victims, instead they must work towards creating a better society. At this juncture we wish to point out that on the lines of the various criminal laws, we can also enact stringent laws and penalty measures to those who discriminate against AIDS patients.

Next, we must also have the proactive approach of ensuring that the people in general mingle with the AIDS patients without hesitation or discrimination. We had observed from our study that because of this stigmatization and isolation, there is only a consequent deterioration and degeneration of the society.

The neglect and isolation of the male AIDS patients only increases their high risk behaviour, making some of them seek CSWs, or drown themselves in alcohol and drugs.

In this too, the suffering of the women AIDS patients is very regrettable. While their husbands are infected the woman sacrifice their entire lives, they live with them, take care of them, and remain with them till their death. But we always observe that when the women die, they are left to die as orphans. Even their own families, the families where they took birth, neglect them. This is very heart-rending; and we must evolve and take initiative to take care of such women patients.

15. An alarming trend that we could notice was the number of people who had sought treatment from quacks, and the high number of people who had taken to astrology, rituals and special prayers to cure themselves of the disease. Because of taking what the patients refer to as the Kerala medicine (spending very high sums of money: as much as 20,000), their health had in fact deteriorated and their condition had become chronic. Because of the marginalization, they do not have access to proper health care facilities in their new surroundings.

That is one main reason why many of them come to the Government Hospital of Thoracic Medicine only in the very advanced stages of AIDS. This can be countered by mounting campaigns that target the AIDS patients and tell them where to get the standard and reliable treatment.

In the context of Tamil Nadu, where more than half of the country's AIDS patients reside, AIDS has become a generalized disease and so, it has become of high importance to develop such campaigns. At least we should work towards debunking the quacks, who squander all the wealth of the AIDS patients.

16. At the village level, programs have to be imparted by the government and voluntary organizations on the alternatives to agriculture and about sustainable agriculture. Many of the rural people who have migrated attributed their migration to the failure of agriculture. The newly developed hybrid varieties of crops have failed, and very often the yield is just sufficient to pay the hired labourers, there is an absence of profit. For most of the poor families in the rural areas, migration has become a survival strategy. While the green revolution in the north ushered in economic change, in the Tamil Nadu situation such a long-lasting change has not been brought. It is essential to educate the people about switching over to other means of sustenance, and at the same time, caution them about the dangers of a migrant life-style.

17. One of our final points is about religion. We could observe that a considerable section of the affected patients had a change of religious belief: one, they stopped believing in God, two, they converted to another faith (in this case, Christianity).

While this loss of faith is attributed to the fact that they lost all faith in God when they learnt that this disease had no cure. In fact many of the patients, sound very rationalistic when they say, "We think Medicine, and not God will cure us". Those who have converted to Christianity after their illness, claim that while the Hindu religion offers no pardon and blames everything on Karma, the Christian religion has a confession and a pardon for all of them.

Some other patients pointed out that in the Christian faith, topics like AIDS was discussed even in the churches and sermons were given. One patient asked as: "When the Hindu religious leaders can speak about Cancer, why don't they lend a voice for the cause of AIDS?"

## 7.5 Adaptation of neutrosophic theory in the analysis of the migrant labourers of rural Tamil Nadu affected by HIV/AIDS

We got the following suggestions and conclusions. We have used the tool of NCM, NRM linked NRM, combined NCM combined NRM, combined disjoint block NRM and combined overlap block NRM.

1.      3% (2 of the 60 persons whom we interviewed) of the migration is an indeterminancy for they do not migrate for job or better living conditions not even attracted by the city poshness but for trying employment in oil countries and countries like Singapore and Malaysia; for already they are very rich and so especially they go to Bombay for visa passport and have caught HIV/AIDS from Bombay.

2.      The migrants who came out on caste clashes still lived away from home town as orphans for the fear of their kith and kin being tortured. Thus several of attributes associated with them remained as indeterminate for they cannot follow the majority pattern but only an exceptional one. Also some ran away with a ransom amount and lead a life of vagabond and there by caught HIV/AIDS but never disclosed how they had such large amount of money!

3.      Also some ran away with a ransom amount and lead a life of vagabond and there by caught HIV/AIDS but never disclosed how they had such large amount of money! The direct role of governments help reaching the poor again labourers remain as an indeterminacy for they fear the head of the village who after getting the cash or kind does not distribute to the real needy by keeps to himself. The poor labourers due to fear cannot question him. There are several stray incidents if questioned the consequence they face are very dire. These factors do not come out at any stage. Government does not take any special steps to see it reaches the grass root level people for they are unconcerned about the welfare of them. The government fears to question or punish these village heads for they are the vote bankers. This vicious circle  is an indeterminacy and needs more attention. So the data supplied by government and the real help received by them can never be correlated in any way.

4. At least 55 out of 60 interviewed had taken treatment several times for STD/VD from several doctors. They by no means were aware of the fact that when they with STD had sex with CSWs who had HIV/AIDS they would directly be infected in the blood stream itself by HIV/AIDS. But they after visited CSWs once again took treatment. This factor remains an indeterminate for when we interviewed the problems and the sufferings they faced due to STD/ VD was very dominating. This remains as an unsurveyed area for we feel 46 out 60 had been infected and became very serious with HIV/AIDS because they are suffering with STD/VD. Thus we suggest it is very much important on the part of all doctors who treat patients for STD/VD to advice them to be proper i.e. not have unprotected sex till they are fully cured of at least open of wounds and ulcers in their private parts. That is why several became full blown HIV/AIDS within a very short span of 3 to 4 years. Certainly this would make those who visit CSWs to be not only careful but save themselves from HIV/AIDS. Unless this is implemented by some forceful means it is impossible even to lessen the number of HIV/AIDS affected people especially in rural areas.

5. Thus it has become a primary importance that whether a counseling center for HIV/AIDS is present in rural areas it is an essentially that a Health Center to treat them for STD/VD and all sexually transmitted disease must be built where the danger of people being easily affected by HIV/AIDS mainly men due to STD/VD should be counseled and must be given all problems they face due to STD/VD when they have unprotected sex with CSWs. For this alone can save the poor ignorant uneducated migrant labour in rural areas from HIV/AIDS. Free treatment for STD/VD must be given with a elaborate advice/ awareness program about HIV/AIDS.

6. The fact whether HIV/AIDS is due to poverty is an indeterminacy but in case of rural uneducated labourer migrancy is a root cause of HIV/AIDS. Thus it is suggested awareness programs about HIV/AIDS be given to rural uneducated in a very different way so that it is not only understandable by them but it makes a strong impact on them so that they are careful to save themselves and others from HIV/AIDS.

7.     Further we saw the compulsory HIV/AIDS test before marriage was an indeterminacy in the investigation. One cannot say for certain it would serve a great purpose in case of male migrants but certainly save their life partners.

8.     The marriage age of women in rural uneducated villages is from 11 years on wards. The main reasons are

1.   Poverty.

2.   No proper school in the nearby area.

3.   No health centers to take care of women psychological problems.

4.   Women treated only as a burden till married to someone.

5.   The large amount of money spent as dowry is also a cause of early marriage.

6.   Thus it is suggested that severe punishment be given to both who give dowry as well as those who receive dowry in marriages.

7.   The law should be implemented powerfully on the marriage age.

8.   These girls should be educated at least up to $10^{th}$ std and the education must be made absolutely free.

9.   The giving in marriage of young girls to very old men as well as second wife should be stopped New laws must be formulated for the same.

Unless these suggestions are strictly followed the status of rural uneducated women in general cannot be changed and in particular the number of HIV/AIDS victims in this group will be exponentially increasing, for they are the passive victims of HIV/AIDS for no fault of them.

9.     Rural uneducated men must change and must be ashamed of their habits which has made them catch HIV/AIDS.

10.    It is an indeterminacy whether HIV/AIDS disease had made them loose faith in god or change faiths. But data shows that 8%  had lost faith in god 10% changed their faith to Christianity.

11.    It is unfortunate to note till date native religious leaders had failed to speak about HIV/AIDS. It is further noted counseling to HIV/AIDS patients are given in churches.

12.    It remains as an indeterminate whether the poor quality of grains and other food products accelerates the HIV/AIDS  victims ie a person who should become a full blown AIDS in 10 years becomes a victims by 2 to 3 years.

13.    The very poor living conditions in the rural uneducated society is one of the cause for spread of HIV/AIDS. ( more or less similar situation as in South Africa ).

14.    This study shows women can be easily educated about awareness of HIV/AIDS but it is very difficult to educate or change rural uneducated men about awareness of HIV/AIDS for CSWs have become the machinery of spreading about HIV/AIDS but these rural migrant labourers failed to take up their advice and become HIV/AIDS victims. This is clear from our analysis. The only reason for this is the society is a male dominated one where women are only objects catering to the needs of their pleasure and labour.

15.    It is important to campaign in all HIV/AIDS awareness programs in rural areas that as most of the disease like cancer, diabetes are not so fully curable so is also HIV/AIDS. This is a livable disease provided people affected by HIV/AIDS lead a hygienic and healthy life and not continue to be slaves of bad habits. This alone will not only remove the social stigma associated with it but help people from neglecting themselves and their family. So that this will voluntarily make them to come and take up blood test and treatment if found to be affected by HIV/AIDS. So the very approach of awareness program in rural areas must not be scaring them there by making the victims as untouchables but the program should make this disease as one among other diseases. Thus is the only way in rural areas to accept an HIV/AIDS patients.

16.     Counseling units in rural areas about HIV/AIDS also about STD / VD are totally absent. For this the unemployed youth can be trained and motivated to do the job of counseling this will play a double role one employment for the unemployed youth as well as the successful awareness program in the rural villages.

17.     Women empowerment groups must be established in all rural villages with government support so that women are not that much explored and made passive victims of HIV/AIDS. It is sad to note most of the women associations union fail to speak about HIV/AIDS infected women on the contrary talk about several other factors which do not in any way help them. For they discover they are HIV/AIDS victims only after a chronic stage and these patients are terminally ill discarded and disowned by the family.

96% of the HIV/AIDS infected migrants interviewed were married men the cause of the infection being only CSWs. Just like South Africa men these married men have multiple sex partners. It remains as an indeterminate whether these CSWs whom they visited regularly were regular customers or whether they had "affairs" with any other women.

So if women are empowered it will certainly decrease the number of men being infected by HIV/AIDS; for they would fearlessly spread about awareness of HIA/AIDS.

Appendix One

# QUESTIONNAIRE

## MANUAL FOR FIELD RESEARCHERS

The points to be had in mind while interviewing the AIDS patients.

❑ Names of the AIDS patients, will not be disclosed and their identities shall be protected.

❑ The interviews of the AIDS patients shall be used by us purely for research purposes.

❑ The patients shall be let to speak out freely and fully. The interviewer shall not intercept when the AIDS patients are speaking (to make them tell what the interviewers have in their mind.)

❑ They shall always move with the AIDS patients with kindness and respect, care and compassion.

❑ When interviewing the patients for the first time, no questions disliked by the patients, shall be asked to them, unless they voluntarily come forward to give the answers.

❑ The interviewers shall select their questions according to the circumstances/ persons.

❑ While gathering information from the AIDS patients, who reside in the hospitals permanently, and from the out-patients of the hospitals an uniform technique shall not be used.

❑ Answers to questions with '*' will be graded and if the answers are indeterminate a note of it will be made.

❑ When these questions are found not adequate, more questions can be framed as per the requirements or depending upon the situation.

❑ The interview shall be opened and concluded depending upon how openly the AIDS patients speak out.

❑ The interview with AIDS patients shall begin with the casual remark that "the AIDS is only a disease which causes death just as all other disease which cause death".

❑ When situations arise in which the AIDS patients conceal certain information / misrepresent or answer that is not determinate, the field researchers shall adopt some psychological techniques and note them down separately.

Interviewee            :

Date and place of interview      :

Circumstance of the interview   :
❑ Bedridden          ❑ Could walk
❑ Could not speak             ❑ Looked well
❑ Looked ill          ❑ Was very feeble
❑ Could only sit      ❑ Was depressed
❑ Normal             ❑ Other observation:
❑ Somebody spoke        ❑ Both of them spoke
❑ Cannot say anything from external appearance
   (indeterminate)

## GENERAL PARTICULARS

1. Name :
2. Age :
3. Sex :
4. Educational qualification :
5. Native place and Address :
6. Present Address :
7. Caste :
8. Religion :

## JOB-RELATED INFORMATION

9. Profession :
10. The office in which you worked/are working :
11. How many person are employed in the office?:
12. Name of the profession :
13. At what age did you join the profession? :
14. Whether the job was permanent or temporary ? :
15. What skills do you know? :
16. Do you do this job hereditarily? :
17. Have you any other sources of income? :
18. How much do you earn as monthly income? :

19. How did you spend
    your income?        :
    ❑ On self                    ❑ On family
    ❑ Towards drink/CSW/smoke    ❑ For accommodation
    ❑ All                        ❑ None

20. Have you gone abroad?
    If yes, the name of country?
    The purpose of your trip?
    The income you earned?    :

21. Whichever places you went
    from your native place
    to earn money?        :

22. Have you informed or not
    about your AIDS disease to :
    ❑ Friends                    ❑ None
    ❑ Siblings (brothers/sisters) ❑ In-laws
    ❑ Children                   ❑ Spouse
    ❑ Other relatives            ❑ Strangers
    ❑ Boss                       ❑ All

23. Do you think that you have
    the right to conceal your
    AIDS disease from
    your employer?        :
    ❑ Yes                        ❑ No
    ❑ Partially                  ❑ Silent
    ❑ No reaction                ❑ Just looks sad
    ❑ Confused                   ❑ Suppresses feeling

24. Do you want to talk about
    your AIDS disease to your
    fellow employees?
    (or, it is not necessary?)    :

MARRIAGE

25. Are you married?        :
26. Your age at marriage    :

27. Love or arranged marriage :

28. How was the marriage held?:

☐ Hindu                    ☐ Muslim
☐ Christian                ☐ None
☐ Unreactive               ☐ Other

29. Are you satisfied with
    your marriage?        :

☐ Fully                    ☐ Partially
☐ Silent                   ☐ Other reaction
☐ Not                      ☐ No reaction

30. If not fully satisfied, why
    did you agree to marriage?  :

31. Were you aware of AIDS
    before marriage?        :

32. Some people have more
    than one wife.
    Do you have so…?         :

**WIFE**

33. Name             :
34. Educational qualification   :
35. Age              :
36. Occupation            :
37. Govt.. /private sector job?  :
38. Is wife alive?       :
    If no, did you remarry?   :

**CHILDREN**

39. How many children?      :
40. Up to what standard the
    children have studied ?   :

41. Age gap between children?   :

42. What jobs the children
    are doing? (If over 18)    :

43. Do your children know
    that you have this disease?   :

44. How you and your children
    were living   before you
    contracted HIV?       :

45. How are your children
    behaving towards you after
    they became aware of your
    AIDS disease?        :

    ❑ Cold and distant      ❑ Never visit
    ❑ Kind              ❑ No concern
    ❑ As usual           ❑ Other reaction

46. Did you contract AIDS before/
    after your child's birth?  :

47. After which child, did you
    contract AIDS?        :

48. Even after knowing that
    you were infected with AIDS
    did you have the child?   :

49. Does that child have AIDS?   :

50. How you had been/have
    been maintaining the child?  :

51. What you expect the Govt.
    to do for your children?  :

52. Whoever, and in whichever
    manner, can AIDS education
    be imparted to children?:

    ❑ Advertisements in TV   ❑ Advt. in Paper
    ❑ Songs/ dramas        ❑ Posters
    ❑ Does not know        ❑ Radio-broadcast
    ❑ Teaching in School    ❑ Others

53. What do you think about
    giving education to the
    AIDS–affected children
    together with others?      :
    - ❑ No                     ❑ Yes
    - ❑ No reaction            ❑ Never
    - ❑ Sign on the face       ❑ Other reaction
    - ❑ Indeterminate
54. Are your children spastic?      :

**PARTICULARS OF PARENTS /IN-LAWS**

55. Name                   :
56. Age                    :
57. Educational Qualification      :
58. Occupation             :
    Govt. or Pvt. employment?   :
59. Whether alive?            :
60. How many wives for father? :
61. Can you tell about the
    properties your parents
    and ancestors left for you?   :
62. Financial position: Before/After:
    - ❑ Rich                   ❑ Poor
    - ❑ In debt                ❑ Upper middle class
    - ❑ Lower middle class     ❑ Others
63. When they gave their
    daughter in marriage to you,
    what did they provide you?   :
64. When your parents/in-laws
    came to know that you have
    AIDS what was their
    psychological reaction?   :
    - ❑ Cold                   ❑ No reaction
    - ❑ Scornful               ❑ Blaming
    - ❑ As usual               ❑ Depressed

- ❑ Dejected
- ❑ Sad

- ❑ Aversion
- ❑ Cannot say

## PARTICULARS OF RESIDENCE

65. Is it own house or rented?        :
66. What is the nature of the
    house you live?        :
67. If own house, are there
    any tenants living in it?  :
68. If rented house what rent
    are you paying?        :
69. If apartment whether
    government housing
    or private buildings? :
70. Whether nuclear family
    or a joint family ?        :
71. Are you living in your
    ancestral home?        :
72. If joint family, with
    whom are you living?    :
73. What type of locality
    is your house situated?    :
74. Whether any of your
    friends have this disease?    :
75. How many of your friends
    have contracted AIDS?    :
76. What is their condition now? :
77. Have you spoken with your
    friends about AIDS?        :
    - ❑ Never
    - ❑ Fully
    - ❑ Aggressive
    - ❑ Despair
    - ❑ Dejection

    - ❑ Sometimes
    - ❑ Silent
    - ❑ Pathetic
    - ❑ Anger
    - ❑ Cannot say

❑ Shows changes in facial expression like fear, sadness

78. What is the nature of your friends? What are their characteristics?

| | |
|---|---|
| ❑ Careless | ❑ Responsible |
| ❑ Irresponsible | ❑ Commercial |
| ❑ Escapist | ❑ Cannot say |
| ❑ Not honest | ❑ Sadist |

79. Can you tell about the friendship between you and your friends in the present context? (if no reply is given we tick under these heads)

| | |
|---|---|
| ❑ Smoke | ❑ Drink |
| ❑ Visits CSWs | ❑ Sensitive |
| ❑ Very Commercial | ❑ Commercial |
| ❑ Very sensitive | ❑ Silent |
| ❑ Happy | ❑ Unhappy |
| ❑ Voice chokes | ❑ Looked in |
| ❑ Cannot say | contemplation |

80. Can you tell about the occupation of your friends? :

**COMMERCIAL SEX WORKERS**

81. Can you tell the reasons for seeking improper sex partner?:

| | |
|---|---|
| ❑ Urge for sex | ❑ Not satisfied |
| ❑ Male Ego | ❑ Silent |
| ❑ Doesn't like to answer | ❑ No reaction |
| ❑ Looks sad | ❑ Other reaction |
| ❑ Cannot say | |

82. How did your improper sexual companion became acquainted with you? :

83. What kind of expenses did you meet when you sought

improper sex partner?    :

84. Can you give the particulars
    regarding your expenses in
    relation to the CSWs
    who you visited?      :

85. Have you taken any of your
    friends when you went
    again to the CSW?

86. What do you think about the
    living standard of the CSW?
    - ❏ Forced into the trade        ❏ Destitute
    - ❏ Poverty                      ❏ Family burden
    - ❏  As per liking               ❏ For money
    - ❏ Others                       ❏ No answer
    - ❏ Cannot say

87. Did you fear getting AIDS
    when you went to the CSW? :
    - ❏ Depressed                    ❏ Sad
    - ❏ Contempt                     ❏ No reaction
    - ❏ Anger                        ❏ Fear
    - ❏ Despair                      ❏ Cannot say

88. When you have sexual
    relation with her, will you
    ask her if she has AIDS?  :
    - ❏ Yes                          ❏ No
    - ❏ Fear to ask                  ❏ Guilty look
    - ❏ Cannot ask                   ❏ Look of despair
    - ❏ Look of emptiness            ❏ Look of repentance
    - ❏ Look of contempt             ❏ Others

89. Have you consulted  your
    friends about the changes
    occurred in your body after
    you had sex with CSW?   :
    - ❏ Yes/ No                      ❏ Cannot say
    - ❏ Fear                         ❏ Sad look

❑ Tries to hide expression    ❑ Guilty expression
❑ Unhappy expression    ❑ No reaction
❑ Depression seem    ❑ Change of voice
❑ Change in posture    ❑ Change in topic
❑ Hides all emotions    ❑ Silent
❑ Closes eye    ❑ Looks vacant

90. How many times you had sexual relation with a CSW before you came to know about AIDS?   :

91. Can you tell me the words spoken to you by the CSW during your visit?   :

92. What was your mental condition after you returned from the CSW?   :

❑ Sad    ❑ No reaction
❑ Happy    ❑ Guilty
❑ Not guilty    ❑ Unhappy
❑ Others    ❑ Dejected
❑ Not able to determine

93. What were your thoughts when you had sex with your spouse after returning from the CSW?   :

❑ No feeling    ❑ Sad
❑ Guilty    ❑ Unreactive
❑ Silent    ❑ Looks vacant
❑ Desperate    ❑ Others
❑ Cannot say

**STATE OF MIND OF THE PERSON**

94. Can you tell how you contracted AIDS disease?:

95. In which year, month and

day, the symptoms of AIDS
began to show up?    :

96. After changes in the body
became visible have you
realized that it was AIDS?   :

97. Can you tell the gradual
changes that happened in
the body as symptom of AIDS?:

98. How was your physical
condition before you
contracted HIV/AIDS?   :

Changes in mind/emotion

   ❏ Solitary            ❏ Blank
   ❏ Guilty              ❏ Sad
   ❏ Unhappy         ❏ Irritated
   ❏ Desperate       ❏ Depressed
   ❏ Suicidal          ❏ Others
   ❏ Cannot say

99.   To whom did you first reveal
the news of your infection?:

100.  How the people around
you moved with before
you contracted AIDS? :

101.  How your relations
behaved with you after
they became aware of
your AIDS infection?   :

102.  To what type of physicians
did you go for consultation
(Allopathy/Siddha/Ayurveda) :

103.  What is the reasons for
going to him?      :

104.  Please tell the difference
in your mental conditions

between the occasion you went to the doctor for consultation for other ailments and when you went to consult on AIDS? :

- ❏ No feeling
- ❏ High heartbeat
- ❏ Sad
- ❏ Desperate
- ❏ Confused
- ❏ Cannot say
- ❏ Ashamed
- ❏ Hides Feeling
- ❏ Disgraced
- ❏ Others
- ❏ Ignorance

105. Can you point out the difference between the behaviour of your relatives towards you at the early stages of your AIDS disease and their behaviour after your disease became chronic? :

106. Which hospital you opted among the government hospitals? Why? :

107. Can you tell us about the manner in which your doctor is taking care? :

- ❏ Well
- ❏ Unconcerned
- ❏ Very caring
- ❏ Discriminatory
- ❏ Cannot say
- ❏ Not so nicely
- ❏ No reaction
- ❏ Committed
- ❏ Uptight

108. Whom you doubt as the person who caused this damage to you? :

109. Do you want to disclose to others about your AIDS disease?(or, do you want

to conceal it?)          :

110. Please tell the reasons?    :

111. After knowing about your
     infection, did you have
     sex with the same person?:

112. Does anyone else in
     your family have AIDS?
     Give the particulars.  :

113. When you came to know
     that you had AIDS, what
     was your mental state?    :
     ❑ Sad                     ❑ Ashamed
     ❑ Highly depressed        ❑ Suicidal
     ❑ No feeling              ❑ Shocked
     ❑ High heartbeat          ❑ Wept uncontrollably
     ❑ Unhappy                 ❑ Disgraced
     ❑ Lost Happiness          ❑ Courageous
     ❑ Cannot say              ❑ Others

114. What is your present
     mental condition after
     your working in the
     early stages?       :

115. At the time when AIDS
     was said to be a horrible
     disease some said not to
     worry about it. What do
     you think about this?  :

116. The patients, who suffer
     from diseases like heart
     and kidney ailment,
     invite quick death
     by bemoaning in a fear
     psychology thinking always
     about the disease. What
     do you think about them? :

- ❑ Mercy killing justification ❑ Worst in case of AIDS
- ❑ Sad ❑ No reaction
- ❑ Cannot say ❑ Others

117. Did your pulse rate shoot up when you talked or thought about AIDS at the early stages of the disease? :

118. How long the impact of the AIDS, talks about it, the reading about it and the thoughts about it, were working upon you? :
- ❑ Never ❑ Always
- ❑ No reaction ❑ Passive
- ❑ Silent ❑ Others
- ❑ Cannot say

119. What was your mental condition when you read or heard news on AIDS in newspapers, radio and television after you contracted AIDS? :
- ❑ Sad ❑ Ashamed
- ❑ Highly depressed ❑ Suicidal
- ❑ No feeling ❑ Shocked
- ❑ High heartbeat ❑ Wept uncontrollably
- ❑ Unhappy ❑ Disgraced
- ❑ Lost Happiness ❑ Courageous
- ❑ Cannot say ❑ Others

120. What do you expect the mental state of other AIDS patients when they see AIDS advertisements : (pick from above)

121. How and in which manner the daily life pattern was affected after you

contracted AIDS?     :

- ❑ Socially       ❑ Economically
- ❑ Mentally       ❑ Materially
- ❑ Spiritually     ❑ Physically
- ❑ None           ❑ All
- ❑ Cannot say

122. If you had been educated already about AIDS, could you have escaped from this disease?     :

123. If not, what reason can you tell?     :

## GENERAL QUESTIONS

124. Have you consulted your family members regarding AIDS?     :

125. Are there AIDS patients in your family/relations?  :

126. Do you know anybody, known to you who has HIV?:

127. What you think about how your family members/ relations/friends should behave in taking care of the AIDS patients?  :

- ❑ Kind          ❑ Harsh
- ❑ Open          ❑ Considerate
- ❑ Cannot say     ❑ Others

128.* What is your view on the opinion that if medicine is discovered to fully cure AIDS, social deterioration will further be worsened? :

129.* What is your opinion on the concept that the government should call out all the AIDS patients for free medical examination/checkup?     :

130.* What do you think about those who say: AIDS will not touch me. Don't talk about it to me?        :

131.* What is your view on the statement that "it is the social duty of everyone to help AIDS patients?"     :

132.  What do you think about teaching sex education?           :

133.  What do you think about yellow journalism and pornographic/blue films? :

134.  Different views prevail in different countries on the use of contraceptives. What is your opinion? :

135.  Whether one has or has not the right to know whether his colleagues in his workplace has AIDS?:

136.* It is natural to have a fear while speaking, reading and commenting about AIDS. Is there a way to change this altitude?   :

137. Can you differentiate between HIV/AIDS/STD from one another?     :

138. Do you know how many years will it take for HIV to convert into AIDS?  :

139.* Have you had a discussion with your colleagues and friends on AIDS?     :

140. About what you would talk during such discussion?     :

141.* What is your mental state on seeing an AIDS patient?     :

142. Will you explain to the public about AIDS if you are given an opportunity to do so?             :

143. Many people speak differently about AIDS. You might have heard/ spoken to some of them. Can you recall any one unforgettable incident?     :

144. Is the AIDS patient alone responsible for contracting the disease?   :

145. "When anyone sees/hears an AIDS patient, s/he is emotionally surcharged. S/he is subjected to mental affliction/grief/ anger/insult/frustration and behaves so. Can you tell the reason for such behaviour?    :

146.  "The foremost duty of the government is to ascertain how many citizens are suffering from AIDS. But no such precise assessment has been made in India." What is your opinion on this statement?    :

147.  Why no free medical test is not being conducted on all citizens to find out whether anyone has AIDS like for smallpox/ diabetes/leprosy?    :

148.* We talk about diabetes, leprosy and cancer in all places. But, we do not talk about AIDS. Why?:

149.  What do you think about those who say "Increasing of AIDS disease is not my problem. It is the problem of concerned persons".    :

150.* "When you know your friend has AIDS how do you move with him" why? :

151.  What expenses you want the government to meet for AIDS control?    :

152.  What further measures do you expect the government

to take to control AIDS? :

153. Which of the measures
the government is under-
taking to prevent AIDS is
most attractive to you? :

154. In which medical system
you have faith: Allopathy,
Siddha/Ayurveda/Unani?
The reason for it? :

155. Do you think that AIDS will
spread through mosquito
biting/kissing on hand/lips
and eating in same plate? :

156. Which advertisements for
AIDS awareness is mostly
liked by you? :

157. What do you think when
you see in private the
advertisements on AIDS? :
❑ Fear                    ❑ Guilt
❑ Unconcerned             ❑ Sadness
❑ Cannot say              ❑ Others

158. What do you think when
you see advertisements
in the company of your
family members? :
❑ Fear                    ❑ Guilt
❑ Unconcerned             ❑ Sadness
❑ Cannot say              ❑ Others

159. What do you think of AIDS
awareness advertisements? :

160. Can you suggest good
advertisement techniques
for AIDS awareness/

publicity?            :

161. What further steps the central/state governments must take to eradicate AIDS?    :

162. Should all people agree to undergo AIDS test if the government comes forward to do the test freely?

163. What kind of expense you want the government to bear for the AIDS control? :

164. Do you think that the AIDS patients should/shouldn't be employed in govt./pvt. services?            :

165. "Any disease can be contracted in the world. But the AIDS should never come" what do you think about this?        :

166.* "AIDS patients should be given top priority in medical treatment above all other patients". What is your view on this opinion? :

167. Do you think that the AIDS patient should also be present when the doctors are discussing about the nature of AIDS?        :

168.* What do you think about interviewing you by people like us?            :

169. What other questions you
     want us to ask you, other
     than those we asked you? :

170. Please tell what part of our
     questions, you like/don't
     like the most?           :

## ADVERTISEMENTS

171.*Through which source you
     came to know about AIDS
     before you contracted it? :

172. If you came to know it
     through advertisements,
     what type was the media? :

173. Our govt. has provided
     telephone facilities in 77
     places throughout the nation
     to tell confidentially doubts
     questions therein. Have you
     used this facility? Do you
     know about this facility?   :

174. Do you know about AIDS?
     By whatever further steps
     the government can take
     the advts. to the masses? :

175.* Now you regret that the
     immoral behaviour which
     you started one day has
     now brought dreaded  AIDS
     disease to you. What do you
     think about it now?      :

176.*What is your view about

fully eradicating prostitution
from our country?     :

177. How can the NGO/social
     organizations make the
     advertisements on AIDS
     reach the public?    :

178. In what ways can your
     religion publicize AIDS?   :

**FOOD HABITS**

179. What kind of food habits
     your family was used to
     vegetarian/meat-eating?  :

180. Have you changed your
     food habits after this
     disease contracted you?  :

181. If you have changed, how?:

182. How is your physical
     conditions after you changed
     your food habits?

**RELIGION**

183. Do you profess the same
     religion hereditarily? If not,
     in which generation you
     converted to this religion? :

184. Have you converted into
     any other religion after
     you contracted AIDS?  :

185. Did you convert? From
     which religion?      :

186. What services you rendered
     to people of that religion
     after your conversion? Or
     was it for personal solace?:
187. What reasons you can
     give for your conversion?  :
188. Do you think the religion is
     giving you consolation after
     your conversion? Or is it a
     testing tool?         :
189. Did both of your parents
     belong to a same religion? :
190. Are you continuing to visit
     the religious institutions
     after you contracted AIDS?:
191. Do you continue to believe
     in God after getting AIDS? :

Appendix Two

# TABLE OF STATISTICS

| SL. No. | AGE/-SEX | EDU. | NATIVE PLACE | PRESENT ADDRESS | PLACE OF INFECTION | HEREDITARY OCCUPATION | PRESENT OCCUPATI-ON |
|---|---|---|---|---|---|---|---|
| 1. | 25-M | Nil | Andhra Pradesh | Ashok Nagar, Chennai | Andhra Pradesh | Sanitary Worker | Sanitary Worker |
| 2. | 44-M | Nil | Neyveli | Medavakka-m | Assam-Ker-ala route | Nil | Lorry Driver |
| 3. | 40-M | 6th | Velavalampatti | Hospital | Many states | Agriculture | Lorry Driver |
| 4. | 31-M | 10th | Namakkal | Hospital | Many states | Agriculture | Driver |
| 5. | 30-M | - | Kumarapalaya-m | Hospital | Mumbai, Nagari | Weavers | Nil |
| 6. | 26-M | 7th | Tiruppathur | Hospital | Same locality | Agriculture | Agriculture |
| 7. | 25-M | 4th | Andhra Pradesh | Chennai - 19 | Hyderabad | Nil | Cook |
| 8. | 28-M | 10th | Salem | Cuddalore | Many states | Agriculture | Driver: bus, lorry, truck |
| 9. | 36-M | 6th | Melore | Hospital | Mumbai | Agriculture | Agriculture |
| 10. | 32-M | 10th | Vandavasi | Hospital | Villupuram | Agriculture | Lorry Driver |
| 11. | 58-M | 8th | Did not divulge | Hospital | Pondicherry & Kallakurchi | Agriculture | Agriculture & Postman |
| 12. | 24-M | Nil | Pasumapatti | Hospital | Many states | Agriculture | Lorry driver |
| 13. | 38-M | 3rd | Perambalur | Hospital | Chennai | Agriculture | Several jobs |
| 14. | 30-M | 10th | Near Coimbatore | Hospital | Mumbai | Agriculture | Lorry driver |
| 15. | 35-M | 7Th | Did not divulge | Hospital | Many states | Agriculture | Driver |
| 16. | 28-M | 5th | Did not divulge | Hospital | Did not divulge | Nil | Driver |
| 17. | 41-M | 10th | Bangalore | Hospital | During trips | Nil | Driver |
| 18. | 27-M | 10th | Did not divulge | Hospital | Mumbai | Agriculture | Lorry Driver |
| 19. | 24-M | 6th | Karoor | Hospital | Outside TN | Agriculture | Driver |
| 20. | 28-M | 10th | Tiruvannamal-ai | Hospital | Many states | Nil | Driver |
| 21. | 47-M | Nil | Tiruchi | Hospital | Many states | Agriculture | Driver |
| 22. | 34-M | 10th | Tannampadi | Hospital | Nearby village | Nil | Works in a shop |
| 23. | 26-M | 10th | Old Washermanpet | Old Washer-m-anpet | Tiruthanni | Nil | Salesman |
| 24. | 46-M | 5th | Ottancathiram | Hospital | Other states | Nil | Driver |
| 25. | 42-M | Nil | Did not divulge | Hospital | Chennai | Agriculture | Teashop |
| 26. | 34-M | Nil | Chengelput | Hospital | Mumbai | Agriculture | Agriculture |
| 27. | 34-M | 10th | Pattukottai | Hospital | Mumbai | Agriculture | Mechanic |
| 28. | 40-M | 10th | Dharmapuri | Hospital | Other states | Nil | Driver |
| 29. | 31-M | 9th | Tiruchi | Hospital | Many states | Agriculture | Driver |
| 30. | 46-M | 5th | Did not divulge | Hospital | Outstation | Nil | Driver |
| 31. | 35-M | 5th | Ambattur | Ambattur | AP, Kerala | Nil | Fabric. work |
| 32. | 38-M | 8th | Tirunelveli | Mumbai | Mumbai | Nil | Hotel Server |
| 33. | 27-M | Nil | Tiruvallur | Hospital | Mumbai | Nil | Embroidery |
| 34. | 47-M | 10th | Coimbatore | Hospital | Many states | Nil | Lorry Driver |
| 35. | 32-M | 7th | Chettiarthoppu | Hospital | New Delhi Mumbai, Goa | Agriculture | Agriculture |
| 36. | 34-M | Nil | Did not divulge | Mylapore | Other states | Agriculture | Timber seller |
| 37. | 55-M | Nil | Krishnagiri | Hospital | Other states | Nil | Cloth seller |
| 38. | 26-M | 5th | Karoor | Hospital | Outstation | Agriculture | Agriculture |
| 39. | 35-M | 8th | Chidambaram | Hospital | Thro' spouse | Nil | Cook |
| 40. | 41-M | 5th | Vellore | Hospital | Nellore | Agriculture | Agriculture |
| 41. | 29-M | 10th | Ulundukottai, T.N | Hospital | Nellore | Agriculture | Agriculture |
| 42. | 35-M | Nil | Villupuram | Hospital | Mumbai | Agriculture | Agriculture |
| 43. | 37-M | 6th | Pangkalam-no-rth | North Madras | Namakkal | Nil | Milk vendor |
| 44. | 26-M | Nil | Mayavaram | Virudachal-am | Local areas | Nil | Violinist |
| 45. | 35-M | 8th | Tirunelveli | Hospital | Other states | Agriculture | Driver |
| 46. | 30-M | 5th | Did not divulge | Hospital | Padappai Navalloor | Agriculture | Not known |
| 47. | 40-M | 7th | Did not divulge | Hospital | Other states | Agriculture | Lorry driver |

| 48. | 31-M | ITI | Tirupur | Hospital | Mumbai | Agriculture | Mechanic |
|---|---|---|---|---|---|---|---|
| 49. | 28-M | 8th | Did not divulge | Hospital | Kerala, AP Pondicherry | Nil | Hotel Manager & Auto Driver |
| 50. | 28-M | Nil | Tenkasi | Hospital | Kuttralam Rail. Platform | Agriculture | Agriculture |
| 51. | 44-M | Nil | Neyveli | Hospital | Maharashtr-a Rajasthan | Nil | Lorry driver |
| 52. | 35-M | 8th | Dharmapuri | Hospital | Other states | Agriculture | Driver |
| 53. | 32-M | 9th | Tiruvarur | Hospital | Other states | Agriculture | Lorry driver |
| 54. | 23-M | Nil | Did not divulge | Platform life | Mumbai, Vijayawada Bangalore | Nil | Hotel worker |
| 55. | 25-M | Nil | Did not divulge | Hospital | Other states | Nil | Nil |
| 56. | 27-M | 7th | Tiruppathur | Hospital | Andhra Pradesh | Agriculture | Nil |
| 57. | 37-M | 9th | Gujarat | Chennai | Other states | Nil | Driver |
| 58. | 20-M | 10th | Mumbai | Mumbai | Mumbai | Nil | Nil |
| 59. | 32-M | 5th | Ulundoorpettai | Hospital | Mumbai | Agriculture | Contract Labourer |
| 60. | 32-M | 10th | Did not divulge | Hospital | Many states | Agriculture | Driver |

Appendix Three

# C-PROGRAM FOR CETD AND RTD MATRIX

**PROGRAM IN C LANGUAGE TO FIND THE ROW SUMS OF THE CETD MATRIX AND REFINED TIME DEPENDENT MATRIX**

```
#include<stdio.h>
#include<math.h>
void main()
{
float a[8][8],x[8][8],me[8],std[8];
float t[8],c1=0.0,c2=0.0,c=0.0,max[8];
float ce[8][8],row[8],c3=0.0,min,r,s;
  int i,j,m,n,e[8][8],p,q,u,v;
  clrscr();
  printf("enter the no. of rows and columns\n");
  scanf("%d %d",&m,&n);
  printf("enter the initial matrix\n");
  for(i=0;i<m;i++)
  {
    for(j=0;j<n;j++)
      scanf("%f",&a[i][j]);
    printf("\n");
  }
  printf("enter the time interval in each
row\n");
  for(i=0;i<m;i++)
    scanf("%f",&t[i]);
  printf("enter the epsilon\n");
  scanf("%f",&s);
  for(i=0;i<m;i++)
```

```c
{for(j=0;j<n;j++)
    ce[i][j]=0.0;}
for(i=0;i<m;i++)
{for(j=0;j<n;j++)
    a[i][j]=a[i][j]/t[i];}
for(j=0;j<n;j++)
{for(i=0;i<m;i++)
    c1=c1+a[i][j];
  me[j]=c1/m;
  c1=0.0;}
for(j=0;j<n;j++)
{for(i=0;i<m;i++)
    c2=c2+pow((a[i][j]-me[j]),2);
  std[j]=sqrt(c2/m);
  c2=0.0;}
for(i=0;i<m;i++)
{
  for(j=0;j<n;j++)
    x[i][j]=(a[i][j]-me[j])/std[j];
}
for(i=0;i<m;i++)
{
  c=fabs(x[i][0]);
  for(j=0;j<n;j++)
  {
    if(fabs(x[i][j])>c)
    c=fabs(x[i][j]);
  }
  max[i]=c;
}
min=max[0];
for(i=1;i<m;i++)
{
  if(max[i]<min)
  min=max[i];
}
printf("ATD matrix:\n");
for(i=0;i<m;i++)
```

```c
  {
    for(j=0;j<n;j++)
      printf("%f\t",a[i][j]);
    printf("\n");
  }
  getch();
  printf("\n%f",min);
  printf("\n");
  for(p=0;p<m;p++)
  {
    for(q=0;q<n;q++)
     {
       printf("\n\n");
       if(fabs((fabs(x[p][q])-min))<s)
       {
    r=fabs(x[p][q]);
    printf("RTD matrix for alpha=
%f\n\n",r);
    for(i=0;i<m;i++)
    {for(j=0;j<n;j++)
      { if(x[i][j]>=r)
          e[i][j]=1;
        if(fabs(x[i][j])<r)
          e[i][j]=0;
        if(x[i][j]<=(-r))
          e[i][j]=-1;
        printf("%d\t",e[i][j]);
      }
     printf("\n");
    }
    for(u=0;u<m;u++)
    {for(v=0;v<n;v++)
        ce[u][v]=ce[u][v]+e[u][v];}
      }getch();
    }
  }
  printf("\nCETD matrix is:\n\n");
  for(i=0;i<m;i++)
```

```
  {
    for(j=0;j<n;j++)
      printf("%f\t",ce[i][j]);
    printf("\n");
  }
  getch();
  for(i=0;i<m;i++)
  {
    c3=0.0;
    for(j=0;j<n;j++)
      c3=c3+ce[i][j];
    row[i]=c3;
  }
  printf("\nRow sums of the CETD matrix;\n");
  for(i=0;i<m;i++)
    printf("row %d= %f\n",i+1,row[i]);
  getch();
}
```

Appendix Four

# C-PROGRAM FOR FCM

**TO FIND THE HIDDEN PATTERN OF THE GIVEN VECTOR**
**GIVEN THE CONNECTION MATRIX OF THE FCM**

```c
#include<stdio.h>
#include<math.h>
void main()
{
  int a[8][8],b[8][8],y[8];
  int x[8],x1[8],i,j,k,n,t=1,c,s=0,u,v;
  clrscr();
  printf("enter the order of the matrix: ");
  scanf("%d",&n);
  printf("enter the initial matrix:\n\n");
  for(i=0;i<n;i++)
  {
    for(j=0;j<n;j++)
      scanf("%d",&a[i][j]);
    printf("\n");
  }
  printf("\nenter the input vector:\n");
  for(i=0;i<n;i++)
    scanf("%d",&x[i]);
  for(i=0;i<n;i++)
  {
    x1[i]=x[i];
  }
  for(i=0;i<n;i++)
    y[i]=0;
  for(k=1;t==1;k++)
  {
    for(i=0;i<n;i++)
      b[s][i]=x[i];
    s++;
    for(j=0;j<n;j++)
    {
      c=0;
      for(i=0;i<n;i++)
    c=c+(x[i]*a[i][j]);
      y[j]=c;
    }
    for(i=0;i<n;i++)
    {
      if(y[i]>0)
```

```
   y[i]=1;
      if(y[i]<0)
   y[i]=-1;
      if(y[i]==0)
   y[i]=0;
   }
   printf("\n\n");
   for(i=0;i<n;i++)
   {
      if(abs(y[i]+x1[i])<=1)
   y[i]=y[i]+x1[i];
   }
   for(j=0;j<s;j++)
   {
      t=0;
      for(i=0;i<n;i++)
      {
   if(b[j][i]!=y[i])
     t=1;
   x[i]=y[i];
      }
      if(t==0)
      break;
   }
   }
   printf("\nthe input vector is:\n");
   for(i=0;i<n;i++)
     printf("%d\t",x1[i]);
   for(i=0;i<n;i++)
     b[s][i]=y[i];
   for(j=0;j<=s;j++)
   {
      v=0;
      for(i=0;i<n;i++)
      {
   if(b[j][i]!=y[i])
     v=1;
      }
      if(v==0)
      break;
   }
   u=j;
   printf("\n");
   for(j=u;j<=s;j++)
   {
     printf("\nthe  input  vector  for  the  %d
iteration:",j+1);
     for(i=0;i<n;i++)
        printf("%d\t",b[j][i]);
     printf("\n");
   }
}
```

Appendix Five

# C-PROGRAM FOR CFCM

**PROGRAM IN C LANGUAGE TO FIND THE FIXED POINT FOR
COMBINED FUZZY COGINTIVE MAPS (CFCM)**

```c
#include<stdio.h>
#include<math.h>
void main()
{
  int a[8][8],b[8][8],y[8],m[8][8];
  int x[8],x1[8],i,j,k,n,t=1,c,s=0,u,v,l;
  clrscr();
  printf("enter the number of experts: ");
  scanf("%d",&l);
  printf("enter the order of the matrix: ");
  scanf("%d",&n);
  for(i=0;i<n;i++)
  {
    for(j=0;j<n;j++)
      a[i][j]=0;
  }
  for(k=0;k<l;k++)
  {
    printf("enter the %d expert's matrix\n",k+1);
    for(i=0;i<n;i++)
    {
      for(j=0;j<n;j++)
    {
      scanf("%d",&m[i][j]);
      a[i][j]=a[i][j]+m[i][j];
    }
    }
  }
  printf("\nenter the input vector:\n");
  for(i=0;i<n;i++)
    scanf("%d",&x[i]);
  printf("\nThe given matrix is:\n");
  for(i=0;i<n;i++)
  {
    printf("\n");
    for(j=0;j<n;j++)
      printf("%d\t",a[i][j]);
  }
```

```
for(i=0;i<n;i++)
{
  x1[i]=x[i];
}
for(i=0;i<n;i++)
  y[i]=0;
for(k=1;t==1;k++)
{
  for(i=0;i<n;i++)
    b[s][i]=x[i];
  s++;
  for(j=0;j<n;j++)
  {
    c=0;
    for(i=0;i<n;i++)
  c=c+(x[i]*a[i][j]);
    y[j]=c;
  }
  for(i=0;i<n;i++)
  {
    if(y[i]>0)
  y[i]=1;
    if(y[i]<0)
  y[i]=-1;
    if(y[i]==0)
  y[i]=0;
  }
  printf("\n\n");
  for(i=0;i<n;i++)
  {
    if(abs(y[i]+x1[i])<=1)
  y[i]=y[i]+x1[i];
  }
  for(j=0;j<s;j++)
  {
    t=0;
    for(i=0;i<n;i++)
    {
  if(b[j][i]!=y[i])
    t=1;
  x[i]=y[i];
    }
    if(t==0)
    break;
  }
}
printf("\nthe input vector is:\n");
for(i=0;i<n;i++)
  printf("%d\t",x1[i]);
for(i=0;i<n;i++)
  b[s][i]=y[i];
for(j=0;j<=s;j++)
  {
    v=0;
    for(i=0;i<n;i++)
    {
```

```
   if(b[j][i]!=y[i])
      v=1;
      }
       if(v==0)
       break;
    }
  u=j;
  printf("\n");
  for(j=u;j<=s;j++)
  {
    printf("\nthe   input   vector   for   the   %d
iteration:",j+1);
    for(i=0;i<n;i++)
      printf("%d\t",b[j][i]);
    printf("\n");
  }
}
```

Appendix Six

# C-PROGRAM FOR BAM

**PROGRAM IN C LANGUAGE TO FIND THE FIXED POINT OF THE GIVEN INPUT VECTOR IN BIDIRECTIONAL ASSOCIATIVE MEMORIES MATRIX**

```c
#include<stdio.h>
void main()
{
  int m[8][8],x[8],sx[8],y[8],sy[8];
int x1[8],y1[8],sx1[8],sy1[8],i,j,a,b,t;
  int p,q,u=0,v=0;
  clrscr();
  printf("enter the no. of rows\n");
  scanf("%d",&p);
  printf("enter the no. of columns\n");
  scanf("%d",&q);
  printf("enter the matrix\n");
  for(i=0;i<p;i++)
  {
    for(j=0;j<q;j++)
      scanf("%d",&m[i][j]);
    printf("\n");
  }
  printf("enter the input vector\n");
  for(i=0;i<p;i++)
    scanf("%d",&x[i]);
  printf("\n");
  printf("\nenter the initial signal fn\n");
  for(i=0;i<q;i++)
    scanf("%d",&sy[i]);
  for(i=0;i<p;i++)
  {
    if(x[i]<=0)
      sx[i]=0;
    if(x[i]>0)
      sx[i]=1;
  }
  for(t=0;(u==0)||(v==0);t++)
  {
    u=1;
    v=1;
    printf("\n\n");
```

```c
    for(j=0;j<q;j++)
    {
      a=0;
      for(i=0;i<p;i++)
        a=a+sx[i]*m[i][j];
      y1[j]=a;
    }
    for(i=0;i<q;i++)
    {
      if(y1[i]<0)
        sy1[i]=0;
      if(y1[i]==0)
        sy1[i]=sy[i];
      if(y1[i]>0)
        sy1[i]=1;
    }
    for(i=0;i<p;i++)
    {
      b=0;
      for(j=0;j<q;j++)
        b=b+sy1[j]*m[i][j];
      x1[i]=b;
    }
    for(i=0;i<p;i++)
    {
      if(x1[i]<0)
        sx1[i]=0;
      if(x1[i]==0)
        sx1[i]=sx[i];
      if(x1[i]>0)
        sx1[i]=1;
    }
    for(i=0;i<p;i++)
    {
      if(sx[i]!=sx1[i])
        u=0;
      sx[i]=sx1[i];
    }
    for(i=0;i<q;i++)
    {
      if(sy[i]!=sy1[i])
   v=0;
      sy[i]=sy1[i];
    }
  }
  printf("\ntotal no of iterations\n");
  printf("%d",2*t);
  printf("\n\nthe fixed points are:\n");
  for(i=0;i<p;i++)
    printf("%d\t",sx1[i]);
  printf("\n");
  for(i=0;i<q;i++)
    printf("%d\t",sy1[i]);
}
```

Appendix Seven

# C-PROGRAM FOR FRM

**PROGRAM IN C LANGUAGE TO FIND THE FIXED POINT FOR A GIVEN INPUT VECTOR FOR FRM**

```
#include<stdio.h>
#include<math.h>
void main()
{
  int a[8][8],b[8][8],c[8][8],y[8],y1[8],e[8][8];
  int
x[8],x1[8],i,j,k,t=1,t1=1,d,c1,c2,s=0,u,v,u1,v1,s
1=0,p,q,p1;
  clrscr();
  printf("enter the number of rows and columns of
the matrix: ");
  scanf("%d %d",&p,&q);
  printf("enter the initial matrix:\n\n");
  for(i=0;i<p;i++)
  {
    for(j=0;j<q;j++)
      scanf("%d",&a[i][j]);
    printf("\n");
  }
  printf("\nThe given matrix is:\n");
  for(i=0;i<p;i++)
  {
    for(j=0;j<q;j++)
      printf("%d\t",a[i][j]);
    printf("\n");
  }
  getch();
  printf("enter  the  dimension  of  the  input
vector: ");
  scanf("%d",&p1);
  if(p1==p)
  {
    for(i=0;i<p;i++)
    {
      for(j=0;j<q;j++)
    e[i][j]=a[i][j];
    }
  }
```

```
else
{
  for(i=0;i<p;i++)
   {
     for(j=0;j<q;j++)
  e[j][i]=a[i][j];
   }
q=p;
p=p1;
}
printf("\nenter the input vector:\n");
for(i=0;i<p;i++)
  scanf("%d",&x[i]);
for(i=0;i<p;i++)
{
  x1[i]=x[i];
}
for(j=0;j<q;j++)
{
  c2=0;
  for(i=0;i<p;i++)
    c2=c2+(x1[i]*e[i][j]);
  y[j]=c2;
}
for(i=0;i<q;i++)
{
  if(y[i]>0)
    y[i]=1;
  if(y[i]<0)
    y[i]=-1;
  if(y[i]==0)
    y[i]=0;
}
for(i=0;i<q;i++)
  y1[i]=y[i];
for(k=1;(t==1)&&(t1==1);k++)
{
  for(i=0;i<p;i++)
    b[s1][i]=x[i];
  s1++;
  for(i=0;i<p;i++)
  {
    c1=0;
    for(j=0;j<q;j++)
  c1=c1+(y[j]*e[i][j]);
    x[i]=c1;
  }
  for(i=0;i<p;i++)
  {
    if(x[i]>0)
  x[i]=1;
    if(x[i]<0)
  x[i]=-1;
    if(x[i]==0)
  x[i]=0;
  }
```

```c
 printf("\n\n");
 for(i=0;i<p;i++)
 {
   if(abs(x[i]+x1[i])<=1)
 x[i]=x[i]+x1[i];
 }
 for(j=0;j<s1;j++)
 {
   t1=0;
   for(i=0;i<p;i++)
   {
 if(b[j][i]!=x[i])
   t1=1;
   }
   if(t1==0)
   break;
 }
 for(i=0;i<q;i++)
   c[s][i]=y[i];
 s++;
 for(j=0;j<q;j++)
 {
   d=0;
   for(i=0;i<p;i++)
 d=d+(x[i]*e[i][j]);
   y[j]=d;
 }
 for(i=0;i<q;i++)
 {
   if(y[i]>0)
 y[i]=1;
   if(y[i]<0)
 y[i]=-1;
   if(y[i]==0)
 y[i]=0;
 }
 printf("\n\n");
 for(i=0;i<q;i++)
 {
   if(abs(y[i]+y1[i])<=1)
 y[i]=y[i]+y1[i];
 }
 for(j=0;j<s;j++)
 {
   t=0;
   for(i=0;i<q;i++)
   {
 if(c[j][i]!=y[i])
   t=1;
   }
   if(t==0)
   break;
 }
}
printf("\nthe input vector is:\n");
for(i=0;i<p;i++)
```

```
    printf("%d\t",x1[i]);
  for(i=0;i<p;i++)
    b[s1][i]=x[i];
  for(i=0;i<q;i++)
    c[s][i]=y[i];
  if(t1==0)
  {
    for(j=0;j<=s1;j++)
    {
      v1=0;
      for(i=0;i<p;i++)
      {
    if(b[j][i]!=x[i])
      v1=1;
      }
      if(v1==0)
      break;
    }
    u1=j;
    printf("\n");
    for(j=u1;j<=s1;j++)
    {
      printf("\nthe  input  vector  for  the  %d
iteration:",j+1);
      for(i=0;i<p;i++)
    printf("%d\t",b[j][i]);
      printf("\n");
    }
  }
  else
  {
    for(j=0;j<=s;j++)
    {
      v=0;
      for(i=0;i<q;i++)
      {
    if(c[j][i]!=y[i])
      v=1;
      }
      if(v==0)
      break;
    }
    u=j;
    printf("\n");
    for(j=u;j<=s;j++)
    {
      printf("\nthe  input  vector  for  the  %d
iteration:",j+1);
      for(i=0;i<q;i++)
    printf("%d\t",c[j][i]);
      printf("\n");
    }
  }
}
```

Appendix Eight

# C-PROGRAM FOR CFRM

**PROGRAM IN C LANGUAGE TO FIND THE FIXED POINT FOR A GIVEN INPUT VECTOR FOR CFRM**

```c
#include<stdio.h>
#include<math.h>
void main()
{
  int
a[8][8],b[8][8],c[8][8],y[8],y1[8],e[8][8],m[8][8
];
  int
x[8],x1[8],i,j,k,t=1,t1=1,d,c1,c2,s=0,u,v,u1,v1,s
1=0;
  int p,q,p1,l;
  clrscr();
  printf("enter the number of rows and columns of
the matrix: ");
  scanf("%d %d",&p,&q);
  printf("enter the number of experts: ");
  scanf("%d",&l);
  for(i=0;i<p;i++)
  {
    for(j=0;j<q;j++)
      a[i][j]=0;
  }
  for(k=0;k<l;k++)
  {
    printf("enter the %d expert's matrix\n",k+1);
    for(i=0;i<p;i++)
    {
      for(j=0;j<q;j++)
      {
      scanf("%d",&m[i][j]);
      a[i][j]=a[i][j]+m[i][j];
      }
    }
  }
  printf("\nThe given matrix is:\n");
  for(i=0;i<p;i++)
  {
    for(j=0;j<q;j++)
      printf("%d\t",a[i][j]);
    printf("\n");
  }
```

```
   getch();
   printf("enter    the    dimension    of    the    input
vector: ");
   scanf("%d",&p1);
   if(p1==p)
   {
     for(i=0;i<p;i++)
     {
       for(j=0;j<q;j++)
    e[i][j]=a[i][j];
     }
   }
   else
   {
     for(i=0;i<p;i++)
     {
       for(j=0;j<q;j++)
    e[j][i]=a[i][j];
     }
   q=p;
   p=p1;
   }
   printf("\nenter the input vector:\n");
   for(i=0;i<p;i++)
     scanf("%d",&x[i]);
   for(i=0;i<p;i++)
   {
     x1[i]=x[i];
   }
   for(j=0;j<q;j++)
   {
     c2=0;
     for(i=0;i<p;i++)
       c2=c2+(x1[i]*e[i][j]);
     y[j]=c2;
   }
   for(i=0;i<q;i++)
   {
     if(y[i]>0)
       y[i]=1;
     if(y[i]<0)
       y[i]=-1;
     if(y[i]==0)
       y[i]=0;
   }
   for(i=0;i<q;i++)
     y1[i]=y[i];
   for(k=1;(t==1)&&(t1==1);k++)
   {
     for(i=0;i<p;i++)
       b[s1][i]=x[i];
     s1++;
     for(i=0;i<p;i++)
     {
       c1=0;
       for(j=0;j<q;j++)
```

```
c1=c1+(y[j]*e[i][j]);
   x[i]=c1;
 }
 for(i=0;i<p;i++)
 {
   if(x[i]>0)
x[i]=1;
   if(x[i]<0)
x[i]=-1;
   if(x[i]==0)
x[i]=0;
 }
 printf("\n\n");
 for(i=0;i<p;i++)
 {
   if(abs(x[i]+x1[i])<=1)
x[i]=x[i]+x1[i];
 }
 for(j=0;j<s1;j++)
 {
   t1=0;
   for(i=0;i<p;i++)
   {
if(b[j][i]!=x[i])
  t1=1;
   }
   if(t1==0)
   break;
 }
 for(i=0;i<q;i++)
   c[s][i]=y[i];
 s++;
 for(j=0;j<q;j++)
 {
   d=0;
   for(i=0;i<p;i++)
d=d+(x[i]*e[i][j]);
   y[j]=d;
 }
 for(i=0;i<q;i++)
 {
   if(y[i]>0)
y[i]=1;
   if(y[i]<0)
y[i]=-1;
   if(y[i]==0)
y[i]=0;
 }
 printf("\n\n");
 for(i=0;i<q;i++)
 {
   if(abs(y[i]+y1[i])<=1)
y[i]=y[i]+y1[i];
 }
 for(j=0;j<s;j++)
 {
```

```
      t=0;
      for(i=0;i<q;i++)
      {
   if(c[j][i]!=y[i])
     t=1;
      }
      if(t==0)
      break;
    }
  }
  printf("\nthe input vector is:\n");
  for(i=0;i<p;i++)
    printf("%d\t",x1[i]);
  for(i=0;i<p;i++)
    b[s1][i]=x[i];
  for(i=0;i<q;i++)
    c[s][i]=y[i];
  if(t1==0)
  {
    for(j=0;j<=s1;j++)
    {
      v1=0;
      for(i=0;i<p;i++)
      {
   if(b[j][i]!=x[i])
     v1=1;
      }
      if(v1==0)
      break;
    }
    u1=j;
    printf("\n");
    for(j=u1;j<=s1;j++)
    {printf("\nthe  input  vector  for  the  %d
iteration:",j+1);
      for(i=0;i<p;i++)
    printf("%d\t",b[j][i]);
      printf("\n");
    }
  }
  else
  {for(j=0;j<=s;j++)
    {v=0;
      for(i=0;i<q;i++)
      {
   if(c[j][i]!=y[i])
     v=1;
      }
      if(v==0)
      break;
    }
    u=j;
    printf("\n");
    for(j=u;j<=s;j++)
    {
```

```
      printf("\nthe  input  vector  for  the  %d
iteration:",j+1);
      for(i=0;i<q;i++)
   printf("%d\t",c[j][i]);
      printf("\n");
     }}
}
```

Appendix Nine

# NEUTROSOPHY

The neutrosophics are based on neutrosophy, which is an extension of dialectics. Neutrosophy is a theory developed by Florentin Smarandache as a generalization of dialectics. For more on neutrosophy, please refer:
http://arxiv.org/ftp/math/papers/0010/0010099.pdf

This theory considers every notion or idea <A> together with its opposite or negation <Anti-A> and the spectrum of "neutralities" <Neut-A> (i.e. notions or ideas located between the two extremes, supporting neither <A> nor <Anti-A>). The <Neut-A> and <Anti-A> ideas together are referred to as <Non-A>. The theory proves that every idea <A> tends to be neutralized and balanced by <Anti-A> and <Non-A> ideas—as a state of equilibrium. It was K.Atanassov (1986) who first introduced the Intuitionistic Fuzzy Set(IFS), describing a degree of membership of an element to a set , amd a degree of non-membership of that element to the set, what's left was considered as indeterminacy. This was generalized to Neutrosophic Set (NS) in 1995 by Florentin Smarandache.

**Distinctions between Intutionistic Fuzzy Logic (IFL) and Neutrosophic Logic (NL)**

The distinctions between IFL and NL {plus the corresponding intuitionistic fuzzy set (IFS) and neutrosophic set (NS)} are the following.

a) Neutrosophic Logic can distinguish between *absolute truth* (truth in all possible worlds, according to Leibniz) and *relative truth* (truth in at least one world), because NL(absolute truth)=$1^+$ while NL(relative truth)=1.  This has application in philosophy (see the neutrosophy).  That's why the unitary standard interval [0, 1] used in IFL has been extended to the unitary non-standard interval $]^-0, 1^+[$ in NL. Similar distinctions for absolute or relative falsehood, and absolute or relative indeterminacy are allowed in NL.

b) In NL there is no restriction on T, I, F other than they are subsets of $]^-0, 1^+[$, thus: $^-0 \leq \inf T + \inf I + \inf F \leq \sup T + \sup I + \sup F \leq 3^+$. This non-restriction allows paraconsistent, dialetheist, and incomplete information to be characterized in NL {i.e. the sum of all three components if they are defined as points, or sum of superior limits of all three components if they are defined as subsets can be >1 (for paraconsistent information coming from different sources) or < 1 for incomplete information}, while that information can not be described in IFL because in IFL the components T (truth), I (indeterminacy), F (falsehood) are restricted either to t+i+f=1 or to $t^2 + f^2 \leq 1$, if T, I, F are all reduced to the points t, i, f respectively, or to sup T + sup I + sup F = 1 if T, I, F are subsets of [0, 1].

c) In NL the components T, I, F can also be *non-standard* subsets included in the unitary non-standard interval $]^-0, 1^+[$, not only *standard* subsets included in the unitary standard interval [0, 1] as in IFL.

d) NL, like dialetheism, can describe paradoxes, NL(paradox) = (1, I, 1), while IFL can not describe a paradox because the sum of components should be 1 in IFL

**Distinctions between Neutrosophic Set (NS) and Intuitionistic Fuzzy Set (IFS).**

a) Neutrosophic Set can distinguish between *absolute membership* (i.e. membership in all possible worlds; we have extended Leibniz's absolute truth to absolute membership) and *relative membership* (membership in at least one world but not in all), because NS(absolute membership element)=$1^+$ while NS(relative membership element)=1. This has application in philosophy (see the neutrosophy). That's why the unitary standard interval [0, 1] used in IFS has been extended to the unitary non-standard interval $]^-0, 1^+[$ in NS.

Similar distinctions for *absolute or relative non-membership*, and *absolute or relative indeterminant appurtenance* are allowed in NS.

b) In NS there is no restriction on T, I, F other than they are subsets of $]^-0, 1^+[$, thus: $^-0 [ \inf T + \inf I + \inf F [ \sup T + \sup I + \sup F [ 3^+$.
The inequalities (2.1) and (2.4) of IFS are relaxed.

This non-restriction allows paraconsistent, dialetheist, and incomplete information to be characterized in NS {i.e. the sum of all three components if they are defined as points, or sum of superior limits of all three components if they are defined as subsets can be >1 (for paraconsistent information coming from different sources), or < 1 for incomplete information}, while that information can not be described in IFS because in IFS the components T (membership), I (indeterminacy), F (non-membership) are restricted either to t+i+f=1 or to $t^2 + f^2 \leq 1$, if T, I, F are all reduced to the points t, i, f respectively, or to sup T + sup I + sup F = 1 if T, I, F are subsets of [0, 1].

Of course, there are cases when paraconsistent and incomplete informations can be normalized to 1, but this procedure is not always suitable.

c) Relation (2.3) from interval-valued intuitionistic fuzzy set is relaxed in NS, i.e. the intervals do not necessarily belong to Int[0,1] but to [0,1], even more general to ]-0, 1+[.

d) In NS the components T, I, F can also be *non-standard* subsets included in the unitary non-standard interval $]^{-}0, 1^{+}[$, not only *standard* subsets included in the unitary standard interval [0, 1] as in IFS.

e) NS, like dialetheism, can describe paradoxist elements, NS(paradoxist element) = (1, I, 1), while IFL can not describe a paradox because the sum of components should be 1 in IFS.

f) The connectors in IFS are defined with respect to T and F, i.e. membership and non-membership only (hence the Indeterminacy is what's left from 1), while in NS they can be defined with respect to any of them (no restriction).

g) Component "I", indeterminacy, can be split into more subcomponents in order to better catch the vague information we work with, and such, for example, one can get more accurate answers to the *Question-Answering Systems* initiated by Zadeh (2003). {In Belnap's four-valued logic (1977) indeterminacy is split into Uncertainty (U) and Contradiction (C), but they were inter-related.}

# INDEX

# ABOUT THE AUTHORS

**Dr.W.B.Vasantha Kandasamy** is an Associate Professor in the Department of Mathematics, Indian Institute of Technology Madras, Chennai, where she lives with her husband Dr.K.Kandasamy and daughters Meena and Kama. Her current interests include Smarandache algebraic structures, fuzzy theory, coding/ communication theory. In the past decade she has guided eight Ph.D. scholars in the different fields of non-associative algebras, algebraic coding theory, transportation theory, fuzzy groups, and applications of fuzzy theory of the problems faced in chemical industries and cement industries. Currently, six Ph.D. scholars are working under her guidance. She has to her credit 255 research papers of which 203 are individually authored. Apart from this, she and her students have presented around 294 papers in national and international conferences. She teaches both undergraduate and post-graduate students and has guided over 41 M.Sc. and M.Tech. projects. She has worked in collaboration projects with the Indian Space Research Organization and with the Tamil Nadu State AIDS Control Society. She has authored a Book Series, consisting of ten research books on the topic of Smarandache Algebraic Structures which were published by the American Research Press.

She can be contacted at vasantha@itm.ac.in
You can visit her work on the web at: http://mat.iitm.ac.in/~wbv

**Dr.Florentin Smarandache** is an Associate Professor of Mathematics at the University of New Mexico, Gallup Campus, USA. He published over 60 books and 80 papers and notes in mathematics, philosophy, literature, rebus. In mathematics his research papers are in number theory, non-Euclidean geometry, synthetic geometry, algebraic structures, statistics, and multiple valued logic (fuzzy logic and fuzzy set, neutrosophic logic and neutrosophic set, neutrosophic probability). He contributed with proposed problems and solutions to the Student Mathematical Competitions. His latest interest is in information fusion were he works with Dr.Jean Dezert from ONERA (French National Establishment for Aerospace Research in Paris) in creasing a new theory and plausible and paradoxical reasoning (DSmT).

He can be contacted at smarand@unm.edu